
\ifx\shlhetal\undefinedcontrolsequence\let\shlhetal\relax\fi

\input amstex
\expandafter\ifx\csname mathdefs.tex\endcsname\relax
  \expandafter\gdef\csname mathdefs.tex\endcsname{}
\else \message{Hey!  Apparently you were trying to
  \string\input{mathdefs.tex} twice.   This does not make sense.} 
\errmessage{Please edit your file (probably \jobname.tex) and remove
any duplicate ``\string\input'' lines}\endinput\fi




\catcode`\X=12\catcode`\@=11

\def\n@wcount{\alloc@0\count\countdef\insc@unt}
\def\n@wwrite{\alloc@7\write\chardef\sixt@@n}
\def\n@wread{\alloc@6\read\chardef\sixt@@n}
\def\r@s@t{\relax}\def\v@idline{\par}\def\@mputate#1/{#1}
\def\l@c@l#1X{\firstpart.#1}\def\gl@b@l#1X{#1}\def\t@d@l#1X{{}}

\def\RED#1{\pfeilsw
            #1%
				              \pfeilso}

\def\crossrefs#1{\ifx\all#1\let\tr@ce=\all\else\def\tr@ce{#1,}\fi
   \n@wwrite\cit@tionsout\openout\cit@tionsout=\jobname.cit 
   \write\cit@tionsout{\tr@ce}\expandafter\setfl@gs\tr@ce,}
\def\setfl@gs#1,{\def\@{#1}\ifx\@\empty\let\next=\relax
   \else\let\next=\setfl@gs\expandafter\xdef
   \csname#1tr@cetrue\endcsname{}\fi\next}
\def\m@ketag#1#2{\expandafter\n@wcount\csname#2tagno\endcsname
     \csname#2tagno\endcsname=0\let\tail=\all\xdef\all{\tail#2,}
   \ifx#1\l@c@l\let\tail=\r@s@t\xdef\r@s@t{\csname#2tagno\endcsname=0\tail}\fi
   \expandafter\gdef\csname#2cite\endcsname##1{\expandafter
     \ifx\csname#2tag##1\endcsname\relax\RED{##1}\else\csname#2tag##1\endcsname\fi
     \expandafter\ifx\csname#2tr@cetrue\endcsname\relax\else
     \write\cit@tionsout{#2tag ##1 cited on page \folio.}\fi}
   \expandafter\gdef\csname#2page\endcsname##1{\expandafter
     \ifx\csname#2page##1\endcsname\relax\RED{##1}\else\csname#2page##1\endcsname\fi
     \expandafter\ifx\csname#2tr@cetrue\endcsname\relax\else
     \write\cit@tionsout{#2tag ##1 cited on page \folio.}\fi}
   \expandafter\gdef\csname#2tag\endcsname##1{\expandafter
      \ifx\csname#2check##1\endcsname\relax
      \expandafter\xdef\csname#2check##1\endcsname{}%
      \else\immediate\write16{Warning: #2tag ##1 used more than once.}\fi
      \multit@g{#1}{#2}##1/X%
      \write\t@gsout{#2tag ##1 assigned number \csname#2tag##1\endcsname\space
      on page \number\count0.}%
   \csname#2tag##1\endcsname}}

\def\multit@g#1#2#3/#4X{\def\t@mp{#4}\ifx\t@mp\empty%
      \global\advance\csname#2tagno\endcsname by 1 
      \expandafter\xdef\csname#2tag#3\endcsname
      {#1\number\csname#2tagno\endcsnameX}%
   \else\expandafter\ifx\csname#2last#3\endcsname\relax
      \expandafter\n@wcount\csname#2last#3\endcsname
      \global\advance\csname#2tagno\endcsname by 1 
      \expandafter\xdef\csname#2tag#3\endcsname
      {#1\number\csname#2tagno\endcsnameX}
      \write\t@gsout{#2tag #3 assigned number \csname#2tag#3\endcsname\space
      on page \number\count0.}\fi
   \global\advance\csname#2last#3\endcsname by 1
   \def\t@mp{\expandafter\xdef\csname#2tag#3/}%
   \expandafter\t@mp\@mputate#4\endcsname
   {\csname#2tag#3\endcsname\lastpart{\csname#2last#3\endcsname}}\fi}
\def\t@gs#1{\def\all{}\m@ketag#1e\m@ketag#1s\m@ketag\t@d@l p
\let\realscite\scite
\let\realstag\stag
   \m@ketag\gl@b@l r \n@wread\t@gsin
   \openin\t@gsin=\jobname.tgs \re@der \closein\t@gsin
   \n@wwrite\t@gsout\openout\t@gsout=\jobname.tgs }
\outer\def\localtags{\t@gs\l@c@l}
\outer\def\globaltags{\t@gs\gl@b@l}
\outer\def\newlocaltag#1{\m@ketag\l@c@l{#1}}
\outer\def\newglobaltag#1{\m@ketag\gl@b@l{#1}}

\newif\ifpr@ 
\def\m@kecs #1tag #2 assigned number #3 on page #4.%
   {\expandafter\gdef\csname#1tag#2\endcsname{#3}
   \expandafter\gdef\csname#1page#2\endcsname{#4}
   \ifpr@\expandafter\xdef\csname#1check#2\endcsname{}\fi}
\def\re@der{\ifeof\t@gsin\let\next=\relax\else
   \read\t@gsin to\t@gline\ifx\t@gline\v@idline\else
   \expandafter\m@kecs \t@gline\fi\let \next=\re@der\fi\next}
\def\pretags#1{\pr@true\pret@gs#1,,}
\def\pret@gs#1,{\def\@{#1}\ifx\@\empty\let\n@xtfile=\relax
   \else\let\n@xtfile=\pret@gs \openin\t@gsin=#1.tgs \message{#1} \re@der 
   \closein\t@gsin\fi \n@xtfile}

\newcount\sectno\sectno=0\newcount\subsectno\subsectno=0
\newif\ifultr@local \def\ultralocal{\ultr@localtrue}
\def\firstpart{\number\sectno}
\def\lastpart#1{\ifcase#1 \or a\or b\or c\or d\or e\or f\or g\or h\or 
   i\or k\or l\or m\or n\or o\or p\or q\or r\or s\or t\or u\or v\or w\or 
   x\or y\or z \fi}

\def\resetall{\global\advance\sectno by 1\subsectno=0
   \gdef\firstpart{\number\sectno}\r@s@t}
\def\resetsub{\global\advance\subsectno by 1
   \gdef\firstpart{\number\sectno.\number\subsectno}\r@s@t}
\def\newsection#1\par{\resetall\vskip0pt plus.3\vsize\penalty-250
   \vskip0pt plus-.3\vsize\bigskip\bigskip
   \message{#1}\leftline{\bf#1}\nobreak\bigskip}
\def\subsection#1\par{\ifultr@local\resetsub\fi
   \vskip0pt plus.2\vsize\penalty-250\vskip0pt plus-.2\vsize
   \bigskip\smallskip\message{#1}\leftline{\bf#1}\nobreak\medskip}


\newdimen\marginshift

\newdimen\margindelta
\newdimen\marginmax
\newdimen\marginmin

\def\margininit{       
\marginmax=3 true cm                  
				      
\margindelta=0.1 true cm              
\marginmin=0.1true cm                 
\marginshift=\marginmin
}    

\def\t@gsjj#1,{\def\@{#1}\ifx\@\empty\let\next=\relax\else\let\next=\t@gsjj
   \def\@@{p}\ifx\@\@@\else
   \expandafter\gdef\csname#1cite\endcsname##1{\citejj{##1}}
   \expandafter\gdef\csname#1page\endcsname##1{\RED{##1}}
   \expandafter\gdef\csname#1tag\endcsname##1{\tagjj{##1}}\fi\fi\next}
\newif\ifshowstuffinmargin
\showstuffinmarginfalse
\def\jjtags{\ifx\shlhetal\relax 
  \else
\ifx\shlhetal\undefinedcontrolseq
\else
\showstuffinmargintrue
\ifx\all\relax\else\expandafter\t@gsjj\all,\fi\fi \fi
}

\def\tagjj#1{\realstag{#1}\oldmginpar{\zeigen{#1}}}
\def\citejj#1{\rechnen{#1}\mginpar{\zeigen{#1}}}     

\def\rechnen#1{\expandafter\ifx\csname stag#1\endcsname\relax\RED{#1}\else
                           \csname stag#1\endcsname\fi}

\newdimen\theight

\def\marginfont{\sevenrm}

\def\trymarginbox#1{\setbox0=\hbox{\marginfont\hskip\marginshift #1}%
		\global\marginshift\wd0 
		\global\advance\marginshift\margindelta}

\def \oldmginpar#1{%
\ifvmode\setbox0\hbox to \hsize{\hfill\rlap{\marginfont\quad#1}}%
\ht0 0cm
\dp0 0cm
\box0\vskip-\baselineskip
\else 
             \vadjust{\trymarginbox{#1}%
		\ifdim\marginshift>\marginmax \global\marginshift\marginmin
			\trymarginbox{#1}%
                \fi
             \theight=\ht0
             \advance\theight by \dp0    \advance\theight by \lineskip
             \kern -\theight \vbox to \theight{\rightline{\rlap{\box0}}%
\vss}}\fi}

\newdimen\upordown
\global\upordown=8pt
\font\tinyfont=cmtt8 
\def\mginpar#1{\smash{\hbox to 0cm{\kern-10pt\raise7pt\hbox{\tinyfont #1}\hss}}}
\def\mginpar#1{{\hbox to 0cm{\kern-10pt\raise\upordown\hbox{\tinyfont #1}\hss}}\global\upordown-\upordown}


\def\t@gsoff#1,{\def\@{#1}\ifx\@\empty\let\next=\relax\else\let\next=\t@gsoff
   \def\@@{p}\ifx\@\@@\else
   \expandafter\gdef\csname#1cite\endcsname##1{\zeigen{##1}}
   \expandafter\gdef\csname#1page\endcsname##1{?}
   \expandafter\gdef\csname#1tag\endcsname##1{\zeigen{##1}}\fi\fi\next}
\def\verbatimtags{\showstuffinmarginfalse
\ifx\all\relax\else\expandafter\t@gsoff\all,\fi}
\def\zeigen#1{\hbox{$\scriptstyle\langle$}#1\hbox{$\scriptstyle\rangle$}}


\def\margintag#1{\ifshowstuffinmargin\oldmginpar{\zeigen{#1}}\fi}

\def\(#1){\edef\dot@g{\ifmmode\ifinner(\hbox{\noexpand\etag{#1}})
   \else\noexpand\eqno(\hbox{\noexpand\etag{#1}})\fi
   \else(\noexpand\ecite{#1})\fi}\dot@g}

\newif\ifbr@ck
\def\eat#1{}
\def\[#1]{\br@cktrue[\br@cket#1'X]}
\def\br@cket#1'#2X{\def\temp{#2}\ifx\temp\empty\let\next\eat
   \else\let\next\br@cket\fi
   \ifbr@ck\br@ckfalse\br@ck@t#1,X\else\br@cktrue#1\fi\next#2X}
\def\br@ck@t#1,#2X{\def\temp{#2}\ifx\temp\empty\let\neext\eat
   \else\let\neext\br@ck@t\def\temp{,}\fi
   \def\teemp{#1}\ifx\teemp\empty\else\rcite{#1}\fi\temp\neext#2X}
\def\resetbr@cket{\gdef\[##1]{[\rtag{##1}]}}
\def\references{\resetbr@cket\newsection References\par}

\newtoks\symb@ls\newtoks\s@mb@ls\newtoks\p@gelist\n@wcount\ftn@mber
    \ftn@mber=1\newif\ifftn@mbers\ftn@mbersfalse\newif\ifbyp@ge\byp@gefalse
\def\defm@rk{\ifftn@mbers\n@mberm@rk\else\symb@lm@rk\fi}
\def\n@mberm@rk{\xdef\m@rk{{\the\ftn@mber}}%
    \global\advance\ftn@mber by 1 }
\def\rot@te#1{\let\temp=#1\global#1=\expandafter\r@t@te\the\temp,X}
\def\r@t@te#1,#2X{{#2#1}\xdef\m@rk{{#1}}}
\def\b@@st#1{{$^{#1}$}}\def\str@p#1{#1}
\def\symb@lm@rk{\ifbyp@ge\rot@te\p@gelist\ifnum\expandafter\str@p\m@rk=1 
    \s@mb@ls=\symb@ls\fi\write\f@nsout{\number\count0}\fi \rot@te\s@mb@ls}
\def\byp@ge{\byp@getrue\n@wwrite\f@nsin\openin\f@nsin=\jobname.fns 
    \n@wcount\currentp@ge\currentp@ge=0\p@gelist={0}
    \re@dfns\closein\f@nsin\rot@te\p@gelist
    \n@wread\f@nsout\openout\f@nsout=\jobname.fns }
\def\m@kelist#1X#2{{#1,#2}}
\def\re@dfns{\ifeof\f@nsin\let\next=\relax\else\read\f@nsin to \f@nline
    \ifx\f@nline\v@idline\else\let\t@mplist=\p@gelist
    \ifnum\currentp@ge=\f@nline
    \global\p@gelist=\expandafter\m@kelist\the\t@mplistX0
    \else\currentp@ge=\f@nline
    \global\p@gelist=\expandafter\m@kelist\the\t@mplistX1\fi\fi
    \let\next=\re@dfns\fi\next}
\def\symbols#1{\symb@ls={#1}\s@mb@ls=\symb@ls} 
\def\bigsymbol{\textstyle}
\symbols{\bigsymbol\ast,\dagger,\ddagger,\sharp,\flat,\natural,\star}
\def\ftnumbers{\ftn@mberstrue} \def\ftsymbols{\ftn@mbersfalse}
\def\paginal{\byp@ge} \def\resetftnumbers{\ftn@mber=1}
\def\ftnote#1{\defm@rk\expandafter\expandafter\expandafter\footnote
    \expandafter\b@@st\m@rk{#1}}

\long\def\jump#1\endjump{}
\def\ssum{\mathop{\lower .1em\hbox{$\textstyle\Sigma$}}\nolimits}

\def\qed{\nobreak\kern 1em \vrule height .5em width .5em depth 0em}
\def\newneq{\hbox{\rlap{\hbox to 1\wd9{\hss$=$\hss}}\raise .1em 
   \hbox to 1\wd9{\hss$\scriptscriptstyle/$\hss}}}
\def\subsetne{\setbox9 = \hbox{$\subset$}\mathrel{\hbox{\rlap
   {\lower .4em \newneq}\raise .13em \hbox{$\subset$}}}}
\def\supsetne{\setbox9 = \hbox{$\subset$}\mathrel{\hbox{\rlap
   {\lower .4em \newneq}\raise .13em \hbox{$\supset$}}}}

\def\vbar{\mathchoice{\vrule height6.3ptdepth-.5ptwidth.8pt\kern-.8pt}
   {\vrule height6.3ptdepth-.5ptwidth.8pt\kern-.8pt}
   {\vrule height4.1ptdepth-.35ptwidth.6pt\kern-.6pt}
   {\vrule height3.1ptdepth-.25ptwidth.5pt\kern-.5pt}}
\def\f@dge{\mathchoice{}{}{\mkern.5mu}{\mkern.8mu}}
\def\b@c#1#2{{\rm \mkern#2mu\vbar\mkern-#2mu#1}}
\def\b@b#1{{\rm I\mkern-3.5mu #1}}
\def\b@a#1#2{{\rm #1\mkern-#2mu\f@dge #1}}
\def\bb#1{{\count4=`#1 \advance\count4by-64 \ifcase\count4\or\b@a A{11.5}\or
   \b@b B\or\b@c C{5}\or\b@b D\or\b@b E\or\b@b F \or\b@c G{5}\or\b@b H\or
   \b@b I\or\b@c J{3}\or\b@b K\or\b@b L \or\b@b M\or\b@b N\or\b@c O{5} \or
   \b@b P\or\b@c Q{5}\or\b@b R\or\b@a S{8}\or\b@a T{10.5}\or\b@c U{5}\or
   \b@a V{12}\or\b@a W{16.5}\or\b@a X{11}\or\b@a Y{11.7}\or\b@a Z{7.5}\fi}}

\catcode`\X=11 \catcode`\@=12




\let\thischap\jobname

\def\partof#1{\csname returnthe#1part\endcsname}
\def\chapof#1{\csname returnthe#1chap\endcsname}

\def\setchapter#1,#2,#3;{%
  \expandafter\def\csname returnthe#1part\endcsname{#2}%
  \expandafter\def\csname returnthe#1chap\endcsname{#3}%
}

\def\setprevious#1 #2 {%
  \expandafter\def\csname set#1page\endcsname{\input page-#2}
}

\setchapter  E53,B,N;       \setprevious E53 null
\setchapter  300z,B,A;       \setprevious 300z E53
\setchapter  88r,B,I;       \setprevious 88r 300z
\setchapter 300a,A,II.A;      \setprevious 300a 88r
\setchapter 300b,A,II.B;       \setprevious 300b 300a
\setchapter 300c,A,II.C;       \setprevious 300c 300b
\setchapter 300d,A,II.D;       \setprevious 300d 300c
\setchapter 300e,A,II.E;       \setprevious 300e 300d
\setchapter 300f,A,II.F;       \setprevious 300f 300e
\setchapter 300g,A,II.G;       \setprevious 300g 300f
\setchapter  600,B,III;       \setprevious  600 300g
\setchapter  705,B,IV;       \setprevious   705 600
\setchapter  734,B,V;        \setprevious   734 705
\setchapter  E46,B,VI;      \setprevious    E46 734
\setchapter  838,B,VII;      \setprevious   838 E46
\setchapter  300x,B,;      \setprevious   300x 838

\newwrite\pageout
\def\rememberpagenumber{\let\setpage\relax
\openout\pageout page-\jobname  \relax \write\pageout{\setpage\the\pageno.}}

\def\recallpagenumber{\csname set\jobname page\endcsname
\def\headmark##1{\rightheadtext{\chapof{\jobname}.##1}}\WRITETOC}
\def\setupchapter#1{\leftheadtext{\chapof{\jobname}. #1}}

\def\setpage#1.{\pageno=#1\relax\advance\pageno1\relax}

\def\cprefix#1{
\edef\theotherpart{\partof{#1}}\edef\theotherchap{\chapof{#1}}%
\ifx\theotherpart\thispart
   \ifx\theotherchap\thischap 
    \else 
     \theotherchap%
    \fi
   \else 
     \theotherchap\fi}

\def\sectioncite[#1]#2{%
     \cprefix{#2}#1}

\def\chaptercite#1{Chapter \cprefix{#1}}

\edef\thispart{\partof{\thischap}}
\edef\thischap{\chapof{\thischap}}

\def\lastpage of '#1' is #2.{\expandafter\def\csname lastpage#1\endcsname{#2}}

\def\xCITE#1{\chaptercite{#1}}
\def\yCITE[#1]#2{\cprefix{#2}.\scite{#2-#1}}

\newwrite\writetoc
\immediate\openout\writetoc \jobname.toc
\def\addcontents#1{\def\WRITETOC{\immediate\write\writetoc{\noexpand\tocentry{\chapof{\jobname}}{#1}{\number\pageno}}}}


\def\spuriousreset{}


\expandafter\ifx\csname citeadd.tex\endcsname\relax
\expandafter\gdef\csname citeadd.tex\endcsname{}
\else \message{Hey!  Apparently you were trying to
\string\input{citeadd.tex} twice.   This does not make sense.} 
\errmessage{Please edit your file (probably \jobname.tex) and remove
any duplicate ``\string\input'' lines}\endinput\fi

\sectno=-1   
\localtags
\jjtags
\NoBlackBoxes
\define\mr{\medskip\roster}
\define\sn{\smallskip\noindent}
\define\mn{\medskip\noindent}
\define\bn{\bigskip\noindent}
\define\ub{\underbar}
\define\wilog{\text{without loss of generality}}
\define\ermn{\endroster\medskip\noindent}

\define\dbcu{\dsize\bigcup}
\define \nl{\newline}
\define\ortp{\text{\bf tp}}
\define\rest{\restriction}

\newbox\noforkbox \newdimen\forklinewidth
\forklinewidth=0.3pt   
\setbox0\hbox{$\textstyle\bigcup$}
\setbox1\hbox to \wd0{\hfil\vrule width \forklinewidth depth \dp0
                        height \ht0 \hfil}
\wd1=0 cm
\setbox\noforkbox\hbox{\box1\box0\relax}
\def\unionstick{\mathop{\copy\noforkbox}\limits}
\def\nonfork#1#2_#3{#1\unionstick_{\textstyle #3}#2}
\def\nonforkin#1#2_#3^#4{#1\unionstick_{\textstyle #3}^{\textstyle #4}#2}     
%
\setbox0\hbox{$\textstyle\bigcup$}
\setbox1\hbox to \wd0{\hfil{\sl /\/}\hfil}
\setbox2\hbox to \wd0{\hfil\vrule height \ht0 depth \dp0 width
                                \forklinewidth\hfil}
\wd1=0cm
\wd2=0cm
\newbox\doesforkbox
\setbox\doesforkbox\hbox{\box1\box0\relax}
\def\nunionstick{\mathop{\copy\doesforkbox}\limits}

\def\fork#1#2_#3{#1\nunionstick_{\textstyle #3}#2}
\def\forkin#1#2_#3^#4{#1\nunionstick_{\textstyle #3}^{\textstyle #4}#2}     

\magnification=\magstep 1
\documentstyle{amsppt}

{    
\catcode`@11

\ifx\alicetwothousandloaded@\relax
  \endinput\else\global\let\alicetwothousandloaded@\relax\fi

\gdef\subjclass{\let\savedef@\subjclass
 \def\subjclass##1\endsubjclass{\let\subjclass\savedef@
   \toks@{\def\usualspace{{\rm\enspace}}\eightpoint}%
   \toks@@{##1\unskip.}%
   \edef\thesubjclass@{\the\toks@
     \frills@{{\noexpand\rm2000 {\noexpand\it Mathematics Subject
       Classification}.\noexpand\enspace}}%
     \the\toks@@}}%
  \nofrillscheck\subjclass}
} 


\expandafter\ifx\csname alice2jlem.tex\endcsname\relax
  \expandafter\xdef\csname alice2jlem.tex\endcsname{\the\catcode`@}
\else \message{Hey!  Apparently you were trying to
\string\input{alice2jlem.tex}  twice.   This does not make sense.}
\errmessage{Please edit your file (probably \jobname.tex) and remove
any duplicate ``\string\input'' lines}\endinput\fi

\expandafter\ifx\csname bib4plain.tex\endcsname\relax
  \expandafter\gdef\csname bib4plain.tex\endcsname{}
\else \message{Hey!  Apparently you were trying to \string\input
  bib4plain.tex twice.   This does not make sense.}
\errmessage{Please edit your file (probably \jobname.tex) and remove
any duplicate ``\string\input'' lines}\endinput\fi

\def\renewcommand{\newcommand}	       
\edef\cite{\the\catcode`@}%
\catcode`@ = 11
\let\@oldatcatcode = \cite
\chardef\@letter = 11
\chardef\@other = 12
%
%
%
%
\def\@innerdef#1#2{\edef#1{\expandafter\noexpand\csname #2\endcsname}}%
%
%
\@innerdef\@innernewcount{newcount}%
\@innerdef\@innernewdimen{newdimen}%
\@innerdef\@innernewif{newif}%
\@innerdef\@innernewwrite{newwrite}%
%
%
%
\def\@gobble#1{}%
%
%
%
\ifx\inputlineno\@undefined
   \let\@linenumber = \empty 
\else
   \def\@linenumber{\the\inputlineno:\space}%
\fi
%
%
%
\def\@futurenonspacelet#1{\def\cs{#1}%
   \afterassignment\@stepone\let\@nexttoken=
}%
\begingroup 
\def\\{\global\let\@stoken= }%
\\ 
\endgroup
\def\@stepone{\expandafter\futurelet\cs\@steptwo}%
\def\@steptwo{\expandafter\ifx\cs\@stoken\let\@@next=\@stepthree
   \else\let\@@next=\@nexttoken\fi \@@next}%
\def\@stepthree{\afterassignment\@stepone\let\@@next= }%
%
%
%
\def\@getoptionalarg#1{%
   \let\@optionaltemp = #1%
   \let\@optionalnext = \relax
   \@futurenonspacelet\@optionalnext\@bracketcheck
}%
%
%
\def\@bracketcheck{%
   \ifx [\@optionalnext
      \expandafter\@@getoptionalarg
   \else
      \let\@optionalarg = \empty
      \expandafter\@optionaltemp
   \fi
}%
\def\@@getoptionalarg[#1]{%
   \def\@optionalarg{#1}%
   \@optionaltemp
}%
%
%
%
\def\@nnil{\@nil}%
\def\@fornoop#1\@@#2#3{}%
\def\@for#1:=#2\do#3{%
   \edef\@fortmp{#2}%
   \ifx\@fortmp\empty \else
      \expandafter\@forloop#2,\@nil,\@nil\@@#1{#3}%
   \fi
}%
\def\@forloop#1,#2,#3\@@#4#5{\def#4{#1}\ifx #4\@nnil \else
       #5\def#4{#2}\ifx #4\@nnil \else#5\@iforloop #3\@@#4{#5}\fi\fi
}%
\def\@iforloop#1,#2\@@#3#4{\def#3{#1}\ifx #3\@nnil
       \let\@nextwhile=\@fornoop \else
      #4\relax\let\@nextwhile=\@iforloop\fi\@nextwhile#2\@@#3{#4}%
}%
%
%
%
\@innernewif\if@fileexists
\def\@testfileexistence{\@getoptionalarg\@finishtestfileexistence}%
\def\@finishtestfileexistence#1{%
   \begingroup
      \def\extension{#1}%
      \immediate\openin0 =
         \ifx\@optionalarg\empty\jobname\else\@optionalarg\fi
         \ifx\extension\empty \else .#1\fi
         \space
      \ifeof 0
         \global\@fileexistsfalse
      \else
         \global\@fileexiststrue
      \fi
      \immediate\closein0
   \endgroup
}%
%
%
%
%
\def\bibliographystyle#1{%
   \@readauxfile
   \@writeaux{\string\bibstyle{#1}}%
}%
\let\bibstyle = \@gobble
%
%
\let\bblfilebasename = \jobname
\def\bibliography#1{%
   \@readauxfile
   \@writeaux{\string\bibdata{#1}}%
   \@testfileexistence[\bblfilebasename]{bbl}%
   \if@fileexists
      \nobreak
      \@readbblfile
   \fi
}%
\let\bibdata = \@gobble
%
%
\def\nocite#1{%
   \@readauxfile
   \@writeaux{\string\citation{#1}}%
}%
\@innernewif\if@notfirstcitation
%
%
\def\cite{\@getoptionalarg\@cite}%
%
%
\def\@cite#1{%
   \let\@citenotetext = \@optionalarg
   \printcitestart
   \nocite{#1}%
   \@notfirstcitationfalse
   \@for \@citation :=#1\do
   {%
      \expandafter\@onecitation\@citation\@@
   }%
   \ifx\empty\@citenotetext\else
      \printcitenote{\@citenotetext}%
   \fi
   \printcitefinish
}%
\newif\ifweareinprivate
\weareinprivatetrue
\ifx\shlhetal\undefinedcontrolseq\weareinprivatefalse\fi
\ifx\shlhetal\relax\weareinprivatefalse\fi
\def\@onecitation#1\@@{%
   \if@notfirstcitation
      \printbetweencitations
   \fi
   \expandafter \ifx \csname\@citelabel{#1}\endcsname \relax
      \if@citewarning
         \message{\@linenumber Undefined citation `#1'.}%
      \fi
     \iftrue 
      \expandafter\gdef\csname\@citelabel{#1}\endcsname{%
\strut 
\vadjust{\vskip-\dp\strutbox
\vbox to 0pt{\vss\parindent0cm \leftskip=\hsize 
\advance\leftskip3mm
\advance\hsize 4cm\strut\openup-4pt 
\rightskip 0cm plus 1cm minus 0.5cm ?  #1 ?\strut}}
         {\tt
            \escapechar = -1
            \nobreak\hskip0pt\pfeilsw
            \expandafter\string\csname#1\endcsname
             \pfeilso
            \nobreak\hskip0pt
         }%
      }%
     \else  
      \expandafter\gdef\csname\@citelabel{#1}\endcsname{%
            {\tt\expandafter\string\csname#1\endcsname}
      }%
     \fi  
   \fi
   \csname\@citelabel{#1}\endcsname
   \@notfirstcitationtrue
}%
%
%
\def\@citelabel#1{b@#1}%
%
%
\def\@citedef#1#2{\expandafter\gdef\csname\@citelabel{#1}\endcsname{#2}}%
%
%
%
\def\@readbblfile{%
   \ifx\@itemnum\@undefined
      \@innernewcount\@itemnum
   \fi
   \begingroup
      \def\begin##1##2{%
         \setbox0 = \hbox{\biblabelcontents{##2}}%
         \biblabelwidth = \wd0
      }%
      \def\end##1{}
      %
      %
      \@itemnum = 0
      \def\bibitem{\@getoptionalarg\@bibitem}%
      \def\@bibitem{%
         \ifx\@optionalarg\empty
            \expandafter\@numberedbibitem
         \else
            \expandafter\@alphabibitem
         \fi
      }%
      \def\@alphabibitem##1{%
         \expandafter \xdef\csname\@citelabel{##1}\endcsname {\@optionalarg}%
         \ifx\biblabelprecontents\@undefined
            \let\biblabelprecontents = \relax
         \fi
         \ifx\biblabelpostcontents\@undefined
            \let\biblabelpostcontents = \hss
         \fi
         \@finishbibitem{##1}%
      }%
      \def\@numberedbibitem##1{%
         \advance\@itemnum by 1
         \expandafter \xdef\csname\@citelabel{##1}\endcsname{\number\@itemnum}%
         \ifx\biblabelprecontents\@undefined
            \let\biblabelprecontents = \hss
         \fi
         \ifx\biblabelpostcontents\@undefined
            \let\biblabelpostcontents = \relax
         \fi
         \@finishbibitem{##1}%
      }%
      \def\@finishbibitem##1{%
         \biblabelprint{\csname\@citelabel{##1}\endcsname}%
         \@writeaux{\string\@citedef{##1}{\csname\@citelabel{##1}\endcsname}}%
         \ignorespaces
      }%
      %
      %
      \let\em = \bblem
      \let\newblock = \bblnewblock
      \let\sc = \bblsc
      \frenchspacing
      \clubpenalty = 4000 \widowpenalty = 4000
      \tolerance = 10000 \hfuzz = .5pt
      \everypar = {\hangindent = \biblabelwidth
                      \advance\hangindent by \biblabelextraspace}%
      \bblrm
      \parskip = 1.5ex plus .5ex minus .5ex
      \biblabelextraspace = .5em
      \bblhook
      \input \bblfilebasename.bbl
   \endgroup
}%
%
%
\@innernewdimen\biblabelwidth
\@innernewdimen\biblabelextraspace
%
%
%
\def\biblabelprint#1{%
   \noindent
   \hbox to \biblabelwidth{%
      \biblabelprecontents
      \biblabelcontents{#1}%
      \biblabelpostcontents
   }%
   \kern\biblabelextraspace
}%
%
%
%
\def\biblabelcontents#1{{\bblrm [#1]}}%
%
%
\def\bblrm{\rm}%
%
%
\def\bblem{\it}%
%
%
\def\bblsc{\ifx\@scfont\@undefined
              \font\@scfont = cmcsc10
           \fi
           \@scfont
}%
%
%
\def\bblnewblock{\hskip .11em plus .33em minus .07em }%
%
%
\let\bblhook = \empty
%
%
%
\def\printcitestart{[}
\def\printcitefinish{]}
\def\printbetweencitations{, }
\def\printcitenote#1{, #1}
%
%
%
\let\citation = \@gobble
%
%
%
\@innernewcount\@numparams
%
%
\def\newcommand#1{%
   \def\@commandname{#1}%
   \@getoptionalarg\@continuenewcommand
}%
%
%
\def\@continuenewcommand{%
   \@numparams = \ifx\@optionalarg\empty 0\else\@optionalarg \fi \relax
   \@newcommand
}%
%
%
\def\@newcommand#1{%
   \def\@startdef{\expandafter\edef\@commandname}%
   \ifnum\@numparams=0
      \let\@paramdef = \empty
   \else
      \ifnum\@numparams>9
         \errmessage{\the\@numparams\space is too many parameters}%
      \else
         \ifnum\@numparams<0
            \errmessage{\the\@numparams\space is too few parameters}%
         \else
            \edef\@paramdef{%
               \ifcase\@numparams
                  \empty  No arguments.
               \or ####1%
               \or ####1####2%
               \or ####1####2####3%
               \or ####1####2####3####4%
               \or ####1####2####3####4####5%
               \or ####1####2####3####4####5####6%
               \or ####1####2####3####4####5####6####7%
               \or ####1####2####3####4####5####6####7####8%
               \or ####1####2####3####4####5####6####7####8####9%
               \fi
            }%
         \fi
      \fi
   \fi
   \expandafter\@startdef\@paramdef{#1}%
}%
%
%
%
%
\def\@readauxfile{%
   \if@auxfiledone \else 
      \global\@auxfiledonetrue
      \@testfileexistence{aux}%
         \begingroup
            \endlinechar = -1
            \catcode`@ = 11
\citation{Sh:c}
\citation{Sh:1}
\citation{Eh57}
\citation{BoNe94}
\citation{Mo65}
\citation{Sh:a}
\citation{Sh:a}
\citation{Sh:c}
\citation{HHL00}
\citation{Sh:10}
\citation{Sh:a}
\citation{Sh:93}
\citation{KiPi98}
\citation{GIL02}
\citation{Sh:225}
\citation{Sh:225a}
\citation{BLSh:464}
\citation{LwSh:489}
\citation{HyTu91}
\citation{HShT:428}
\citation{HySh:474}
\citation{HySh:529}
\citation{HySh:602}
\citation{CoSh:919}
\citation{BlSh:156}
\citation{BlSh:156}
\citation{BlSh:156}
\citation{Bl85}
\citation{Sh:197}
\citation{Sh:205}
\citation{Sh:284c}
\citation{LwSh:560}
\citation{LwSh:687}
\citation{LwSh:871}
\citation{Sh:840}
\citation{Sh:e}
\citation{GbTl06}
\citation{Sh:c}
\citation{Sh:800}
\citation{Sh:e}
\citation{Sh:800}
\citation{Ke71}
\citation{BaFe85}
\citation{Bal0x}
\citation{Mw85}
\citation{KM67}
\citation{Ke71}
\citation{Sh:3}
\citation{Sh:54}
\citation{Hy98}
\citation{HySh:629}
\citation{HySh:632}
\citation{HySh:629}
\citation{GrLe02}
\citation{GrLe0x}
\citation{GrLe00a}
\citation{Le0x}
\citation{Le0y}
\citation{ChKe62}
\citation{ChKe66}
\citation{He74}
\citation{Str76}
\citation{HeIo02}
\citation{Sh:54}
\citation{BY0y}
\citation{BeUs0x}
\citation{Pi0x}
\citation{ShUs:837}
\citation{Ke71}
\citation{Fr75}
\citation{Sh:48}
\citation{Sh:48}
\citation{Sh:87a}
\citation{Sh:87b}
\citation{GrHa89}
\citation{MkSh:366}
\citation{HaSh:323}
\citation{Zi0xa}
\citation{Zi0xb}
\citation{Sh:300}
\citation{ShVi:635}
\citation{Va02}
\citation{GrSh:222}
\citation{GrSh:238}
\citation{GrSh:259}
\citation{Sh:394}
\citation{Gr91}
\citation{BlSh:330}
\citation{BlSh:360}
\citation{BlSh:393}
\citation{GrVa0xa}
\citation{GrVa0xb}
\citation{BKV0x}
\citation{MaSh:285}
\citation{KlSh:362}
\citation{Sh:472}
\citation{Bal0x}
\citation{Sh:88}
\citation{Jn56}
\citation{Jn60}
\citation{Sh:88}
\citation{Jn56}
\citation{Jn60}
\citation{Jn56}
\citation{Jn60}
\citation{Jn56}
\citation{Jn60}
\citation{MoVa62}
\citation{MoVa62}
\citation{Sh:88}
\citation{Sh:87b}
\citation{Sh:87a}
\citation{Sh:87a}
\citation{Sh:87b}
\citation{Sh:300}
\citation{Sh:155}
\citation{Sh:300}
\citation{Sh:576}
\citation{Sh:88}
\citation{Sh:e}
\citation{Sh:576}
\citation{Sh:576}
\citation{Sh:576}
\citation{Sh:c}
\citation{Sh:200}
\citation{Bal88}
\citation{Sh:3}
\citation{Sh:394}
\citation{Sh:394}
\citation{Sh:56}
\citation{Sh:1}
\citation{Sh:93}
\citation{Sh:576}
\citation{Sh:576}
\citation{Sh:868}
\citation{Sh:576}
\citation{Sh:576}
\citation{Sh:87a}
\citation{Sh:87b}
\citation{Sh:87b}
\citation{Sh:87a}
\citation{Sh:87b}
\citation{Sh:87a}
\citation{Sh:87b}
\citation{Sh:F735}
\citation{Sh:576}
\citation{Sh:603}
\citation{BlSh:862}
\citation{Sh:603}
\citation{Sh:576}
\citation{Sh:576}
\citation{Sh:576}
\citation{Sh:322}
\citation{Sh:87a}
\citation{Sh:87b}
\citation{MaSh:285}
\citation{KlSh:362}
\citation{Sh:472}
\citation{Sh:576}
\citation{Sh:576}
\citation{HaSh:323}
\citation{ShVi:648}
\citation{MaSh:285}
\citation{KlSh:362}
\citation{Sh:472}
\citation{Sh:394}
\citation{Bal0x}
\citation{Sh:842}
\citation{Sh:842}
\citation{Sh:842}
\citation{Sh:c}
\citation{Sh:c}
\citation{Sh:a}
\citation{Sh:c}
\citation{Sh:c}
\citation{Sh:429}
\citation{He92}
\citation{Sh:839}
\citation{Sh:576}
\citation{JrSh:875}
\citation{BlSh:862}
\citation{Zi0xa}
\citation{Zi0xb}
\citation{Sh:F709}
\citation{Sh:E54}
\citation{Sh:F735}
\citation{Sh:F782}
\citation{Sh:F888}
\citation{Sh:E56}
\citation{Sh:F841}
\citation{Ke70}
\citation{Sh:c}
\citation{Sh:842}
\citation{Sh:3}
\citation{HySh:676}
\citation{GrLe0x}
\citation{Sh:839}
\citation{Sh:300}
\citation{Sh:87a}
\citation{Sh:87b}
\citation{Sh:842}
\citation{Sh:300}
\citation{Sh:c}
\citation{Sh:576}
\citation{Sh:F888}
\citation{Sh:842}
\citation{Sh:576}
\citation{Sh:842}
\citation{ignore-this-bibtex-warning}
\citation{Sh:48}
\citation{Sh:c}
\citation{Sh:839}
\citation{Sh:576}
\citation{Sh:576}
\citation{Sh:87a}
\citation{Sh:87b}
\citation{Sh:48}
\citation{Sh:576}
\citation{Sh:E36}
\citation{Sh:56}
\citation{Sh:603}
\citation{Sh:603}
\citation{ignore-this-bibtex-warning}
\citation{Sh:48}
\citation{Sh:87a}
\citation{Sh:87b}
\citation{Sh:87a}
\citation{Sh:87b}
\citation{Sh:48}
\citation{Sh:c}
\citation{Sh:87a}
\citation{Sh:87b}
\citation{Sh:48}
\citation{Sh:87a}
\citation{Sh:87b}
\citation{Sh:868}
\citation{Sh:868}
\citation{Sh:87a}
\citation{Sh:87b}
\citation{Sh:87a}
\citation{Sh:87b}
\citation{Sh:87a}
\citation{Sh:87b}
\citation{Sh:576}
\citation{Sh:46}
\citation{Sh:576}
\citation{Sh:576}
\citation{Sh:88}
\citation{Sh:87a}
\citation{Sh:87b}
\citation{Sh:88}
\citation{Sh:48}
\citation{Sh:88}
\citation{Mw85a}
\citation{Sh:88}
\citation{Sh:F709}
\citation{Sh:88}
\citation{Ke71}
\citation{KM67}
\citation{Sh:87a}
\citation{Sh:87b}
\citation{Sh:3}
\citation{Sh:c}
\citation{Sh:48}
\citation{Ke70}
\citation{Mo70}
\citation{Sh:a}
\citation{Ke70}
\citation{Sh:43}
\citation{Sc76}
\citation{DvSh:65}
\citation{Sh:f}
\citation{Sh:E45}
\citation{Sh:420}
\citation{Sh:E12}
\citation{Sh:c}
\citation{Sh:c}
\citation{Sh:c}
\citation{Sh:394}
\citation{Jo56}
\citation{Jo60}
\citation{KM67}
\citation{Sh:c}
\citation{Sh:c}
\citation{Sh:c}
\citation{Sh:868}
\citation{Sh:c}
\citation{Sh:f}
\citation{Ha61}
\citation{Sh:48}
\citation{Sh:c}
\citation{GrSh:259}
\citation{Sh:48}
\citation{Mo70}
\citation{Mo70}
\citation{Bg}
\citation{Sh:202}
\citation{Sh:43}
\citation{BKM78}
\citation{Sh:48}
\citation{Sh:87a}
\citation{Sh:87b}
\citation{Ke71}
\citation{Sh:48}
\citation{Sh:87a}
\citation{Sh:87b}
\citation{Sh:c}
\citation{Sh:48}
\citation{Sh:87a}
\citation{Sh:87b}
\citation{GrSh:174}
\citation{Sh:87a}
\citation{Sh:87b}
\citation{Sh:87a}
\citation{Sh:87b}
\citation{Sh:87a}
\citation{Sh:87b}
\citation{ignore-this-bibtex-warning}
\citation{Sh:155}
\citation{Sh:c}
\citation{Sh:c}
\citation{Sh:c}
\citation{Sh:c}
\citation{Sh:E54}
\citation{Sh:300}
\citation{Sh:a}
\citation{Sh:16}
\citation{GrSh:222}
\citation{GrSh:238}
\citation{GrSh:259}
\citation{McSh:55}
\citation{GrSh:174}
\citation{Sh:a}
\citation{Sh:c}
\citation{Sh:715}
\citation{Sh:c}
\citation{Sh:c}
\citation{Sh:c}
\citation{Sh:3}
\citation{Sh:e}
\citation{Sh:c}
\citation{MaSh:285}
\citation{Sh:300}
\citation{Sh:e}
\citation{Sh:220}
\citation{Sh:e}
\citation{Sh:16}
\citation{GrSh:222}
\citation{GrSh:259}
\citation{Sh:c}
\citation{Sh:300}
\citation{Sh:e}
\citation{Sh:e}
\citation{Sh:220}
\citation{Sh:e}
\citation{Sh:220}
\citation{Sh:e}
\citation{Sh:e}
\citation{Sh:c}
\citation{GrSh:222}
\citation{GrSh:259}
\citation{Sh:c}
\citation{Sh:c}
\citation{Sh:c}
\citation{Sh:c}
\citation{Sh:c}
\citation{Sh:715}
\citation{Mo65}
\citation{Sh:c}
\citation{ignore-this-bibtex-warning}
\citation{BlSh:330}
\citation{BlSh:360}
\citation{BlSh:393}
\citation{Sh:839}
\citation{Sh:c}
\citation{Sh:c}
\citation{Sh:e}
\citation{Sh:300}
\citation{Sh:3}
\citation{Sh:3}
\citation{Sh:54}
\citation{Sh:3}
\citation{Sh:54}
\citation{Sh:3}
\citation{Sh:3}
\citation{Sh:3}
\citation{Sh:3}
\citation{Sh:c}
\citation{Sh:3}
\citation{Sh:3}
\citation{Sh:3}
\citation{Sh:54}
\citation{Sh:3}
\citation{Sh:3}
\citation{Sh:c}
\citation{HySh:629}
\citation{HySh:676}
\citation{Sh:3}
\citation{ignore-this-bibtex-warning}
\citation{Sh:c}
\citation{Sh:132}
\citation{Sh:c}
\citation{Sh:48}
\citation{Sh:e}
\citation{Sh:420}
\citation{Sh:108}
\citation{Sh:88a}
\citation{Sh:420}
\citation{Sh:351}
\citation{Sh:420}
\citation{Sh:420}
\citation{Sh:351}
\citation{Sh:E54}
\citation{Sh:e}
\citation{DvSh:65}
\citation{Sh:b}
\citation{Sh:87b}
\citation{Sh:c}
\citation{Sh:300}
\citation{Sh:e}
\citation{Sh:300}
\citation{Sh:e}
\citation{Sh:e}
\citation{Sh:300}
\citation{Sh:e}
\citation{ignore-this-bibtex-warning}
\citation{Sh:c}
\citation{Sh:c}
\citation{Sh:a}
\citation{Sh:c}
\citation{Sh:c}
\citation{Sh:a}
\citation{Sh:c}
\citation{Sh:c}
\citation{Sh:c}
\citation{Sh:c}
\citation{Sh:c}
\citation{Sh:c}
\citation{Sh:c}
\citation{Sh:c}
\citation{Sh:c}
\citation{Sh:a}
\citation{Sh:c}
\citation{Sh:c}
\citation{Sh:E54}
\citation{ignore-this-bibtex-warning}
\citation{Sh:c}
\citation{Sh:87a}
\citation{Sh:87b}
\citation{Sh:c}
\citation{Sh:E54}
\citation{Sh:300}
\citation{Sh:e}
\citation{Sh:e}
\citation{Sh:E54}
\citation{Sh:E54}
\citation{Sh:c}
\citation{Sh:c}
\citation{Sh:c}
\citation{Sh:a}
\citation{Sh:c}
\citation{Sh:c}
\citation{Sh:c}
\citation{Sh:c}
\citation{ignore-this-bibtex-warning}
\citation{RuSh:117}
\citation{Sh:f}
\citation{Sh:E54}
\citation{Sh:e}
\citation{Sh:300}
\citation{Sh:e}
\citation{Sh:e}
\citation{Sh:300}
\citation{Sh:E54}
\citation{Sh:E54}
\citation{Sh:E54}
\citation{ignore-this-bibtex-warning}
\citation{Sh:E54}
\citation{Sh:E54}
\citation{KoSh:796}
\citation{ignore-this-bibtex-warning}
\citation{Sh:576}
\citation{Sh:576}
\citation{Sh:576}
\citation{Sh:48}
\citation{Sh:576}
\citation{JrSh:875}
\citation{Sh:F888}
\citation{Sh:576}
\citation{Sh:48}
\citation{Sh:576}
\citation{Sh:87a}
\citation{Sh:87b}
\citation{Sh:87a}
\citation{Sh:87b}
\citation{Sh:48}
\citation{Sh:88}
\citation{Sh:576}
\citation{Sh:576}
\citation{Sh:c}
\citation{Sh:87a}
\citation{Sh:87b}
\citation{Sh:576}
\citation{Sh:87a}
\citation{Sh:576}
\citation{Sh:576}
\citation{Sh:576}
\citation{Sh:576}
\citation{Sh:576}
\citation{Sh:576}
\citation{Sh:576}
\citation{Sh:48}
\citation{Mw85a}
\citation{Sh:88}
\citation{Sh:300}
\citation{Sh:576}
\citation{Sh:842}
\citation{Sh:576}
\citation{Sh:576}
\citation{Sh:87a}
\citation{Sh:87b}
\citation{Sh:87a}
\citation{Sh:87b}
\citation{Sh:87b}
\citation{JrSh:875}
\citation{Sh:842}
\citation{Sh:88}
\citation{Sh:576}
\citation{Sh:48}
\citation{Sh:87a}
\citation{Sh:87b}
\citation{Sh:c}
\citation{Sh:88}
\citation{Sh:87a}
\citation{Sh:87b}
\citation{Sh:48}
\citation{Sh:576}
\citation{Sh:c}
\citation{Sh:c}
\citation{HuSh:342}
\citation{Sh:c}
\citation{Sh:88}
\citation{Sh:88}
\citation{Sh:48}
\citation{Sh:48}
\citation{Sh:48}
\citation{Sh:48}
\citation{Sh:48}
\citation{Ke71}
\citation{Sh:48}
\citation{Sh:48}
\citation{Sh:48}
\citation{Sh:48}
\citation{Sh:48}
\citation{Sh:48}
\citation{Sh:48}
\citation{Sh:48}
\citation{Sh:48}
\citation{Sh:48}
\citation{Sh:48}
\citation{Sh:48}
\citation{Sh:48}
\citation{Sh:48}
\citation{Sh:576}
\citation{Sh:576}
\citation{Sh:576}
\citation{Sh:576}
\citation{Sh:576}
\citation{Sh:576}
\citation{Sh:576}
\citation{Sh:576}
\citation{Sh:576}
\citation{Sh:576}
\citation{Sh:576}
\citation{Sh:576}
\citation{Sh:576}
\citation{Sh:576}
\citation{Sh:576}
\citation{Sh:576}
\citation{Sh:576}
\citation{Sh:576}
\citation{Sh:576}
\citation{Sh:576}
\citation{Sh:576}
\citation{Sh:576}
\citation{Sh:576}
\citation{Sh:482}
\citation{Sh:576}
\citation{Sh:576}
\citation{Sh:832}
\citation{Sh:576}
\citation{Sh:576}
\citation{Sh:576}
\citation{Sh:576}
\citation{Sh:576}
\citation{Sh:576}
\citation{Sh:576}
\citation{Sh:576}
\citation{Sh:576}
\citation{Sh:842}
\citation{Sh:576}
\citation{Sh:48}
\citation{ignore-this-bibtex-warning}
\citation{Sh:F735}
\citation{Sh:c}
\citation{Sh:F735}
\citation{Sh:576}
\citation{Sh:576}
\citation{Sh:48}
\citation{Sh:c}
\citation{Sh:c}
\citation{Sh:225a}
\citation{Sh:429}
\citation{Sh:c}
\citation{ShHM:158}
\citation{Las88}
\citation{Sh:F735}
\citation{Sh:c}
\citation{Sh:576}
\citation{Sh:576}
\citation{Sh:31}
\citation{Sh:576}
\citation{Sh:c}
\citation{Sh:576}
\citation{Sh:F735}
\citation{Sh:c}
\citation{Sh:F735}
\citation{Sh:F735}
\citation{Sh:842}
\citation{Sh:842}
\citation{Sh:F735}
\citation{Sh:c}
\citation{Sh:839}
\citation{Sh:839}
\citation{Sh:300}
\citation{Sh:e}
\citation{Sh:576}
\citation{Sh:87b}
\citation{Sh:c}
\citation{Sh:c}
\citation{Sh:87b}
\citation{Sh:87b}
\citation{Sh:c}
\citation{Sh:87a}
\citation{Sh:87b}
\citation{Sh:87a}
\citation{Sh:87b}
\citation{Sh:842}
\citation{Sh:394}
\citation{Sh:842}
\citation{MaSh:285}
\citation{KlSh:362}
\citation{Sh:472}
\citation{Sh:394}
\citation{Sh:842}
\citation{ignore-this-bibtex-warning}
\citation{Sh:F782}
\citation{Sh:842}
\citation{Sh:842}
\citation{Sh:F820}
\citation{Sh:F782}
\citation{Sh:394}
\citation{Di}
\citation{Sh:g}
\citation{Sh:394}
\citation{Di}
\citation{Sh:300}
\citation{Sh:e}
\citation{Sh:88}
\citation{Sh:F782}
\citation{Kin66}
\citation{Sh:F820}
\citation{Sh:F782}
\citation{Lv71}
\citation{Sh:e}
\citation{Sh:e}
\citation{MiRa65}
\citation{Sh:620}
\citation{Sh:F782}
\citation{Sh:394}
\citation{Sh:394}
\citation{Sh:394}
\citation{Sh:394}
\citation{Sh:394}
\citation{Sh:394}
\citation{Sh:589}
\citation{Sh:g}
\citation{Sh:g}
\citation{Sh:394}
\citation{Sh:394}
\citation{Sh:394}
\citation{Sh:394}
\citation{Sh:394}
\citation{Sh:F782}
\citation{KlSh:362}
\citation{Sh:472}
\citation{Sh:394}
\citation{Sh:394}
\citation{Sh:394}
\citation{ignore-this-bibtex-warning}
\citation{Sh:576}
\citation{Sh:603}
\citation{Sh:576}
\citation{Sh:603}
\citation{Sh:705}
\citation{Sh:c}
\citation{MaSh:285}
\citation{KlSh:362}
\citation{Sh:472}
\citation{Sh:576}
\citation{Sh:603}
\citation{Sh:576}
\citation{Sh:576}
\citation{Sh:603}
\citation{Sh:576}
\citation{Sh:576}
\citation{Sh:576}
\citation{Sh:603}
\citation{Sh:576}
\citation{Sh:F888}
\citation{Sh:300}
\citation{BlSh:862}
\citation{Sh:351}
\citation{Sh:F888}
\citation{Sh:E45}
\citation{Sh:g}
\citation{Sh:g}
\citation{KjSh:409}
\citation{Sh:603}
\citation{Sh:430}
\citation{Sh:F888}
\citation{Sh:922}
\citation{Sh:576}
\citation{ignore-this-bibtex-warning}
\citation{Sh:87b}
\citation{Sh:88r}
\citation{Sh:88}
\citation{Sh:576}
\citation{Sh:603}
\citation{Sh:88r}
\citation{Sh:705}
\citation{Sh:576}
\citation{Sh:87b}
\citation{Sh:88}
\citation{Sh:576}
\citation{Sh:603}
\citation{Sh:87b}
\citation{Sh:842}
\citation{Sh:576}
\citation{Sh:576}
\citation{Sh:576}
\citation{Sh:576}
\citation{Sh:576}
\citation{Sh:576}
\citation{Sh:603}
\citation{Sh:576}
\citation{Sh:f}
\citation{Sh:b}
\citation{Sh:f}
\citation{Sh:f}
\citation{DvSh:65}
\citation{Sh:b}
\citation{Sh:f}
\citation{Sh:E45}
\citation{Sh:576}
\citation{Sh:576}
\citation{Sh:576}
\citation{Sh:603}
\citation{Sh:603}
\citation{Sh:576}
\citation{Sh:576}
\citation{Sh:576}
\citation{Sh:E45}
\citation{Sh:576}
\citation{Sh:576}
\citation{Sh:87b}
\citation{Sh:48}
\citation{Sh:F888}
\citation{Sh:F888}
\citation{Sh:576}
\citation{Sh:F841}
\citation{Sh:F841}
\citation{GiSh:577}
\citation{Sh:576}
\citation{Sh:842}
\citation{Sh:E45}
\citation{Sh:f}
\citation{Sh:460}
\citation{Sh:829}
\citation{Sh:576}
\citation{Sh:603}
\citation{Sh:g}
\citation{Sh:g}
\citation{Sh:430}
\citation{Sh:g}
\citation{Sh:g}
\citation{Sh:f}
\citation{Sh:g}
\citation{Sh:g}
\citation{Sh:g}
\citation{Sh:g}
\citation{Sh:g}
\citation{Sh:g}
\citation{Sh:g}
\citation{Sh:g}
\citation{Sh:g}
\citation{J}
\citation{Sh:460}
\citation{Sh:829}
\citation{Ha61}
\citation{Sh:g}
\citation{Sh:420}
\citation{Sh:g}
\citation{Sh:g}
\citation{Ha61}
\bibstyle{lit-plain}
\bibdata{lista,listb,listx,listf,liste}
\@citedef{Bal88}{Bal88}
\@citedef{Bal0x}{Bal0x}
\@citedef{BlSh:862}{BlSh 862}
\@citedef{Bl85}{Bl85}
\@citedef{BKV0x}{BKV0x}
\@citedef{BLSh:464}{BLSh 464}
\@citedef{BlSh:156}{BlSh 156}
\@citedef{BlSh:330}{BlSh 330}
\@citedef{BlSh:360}{BlSh 360}
\@citedef{BlSh:393}{BlSh 393}
\@citedef{BaFe85}{BaFe85}
\@citedef{BKM78}{BKM78}
\@citedef{BY0y}{BY0y}
\@citedef{BeUs0x}{BeUs0x}
\@citedef{BoNe94}{BoNe94}
\@citedef{Bg}{Bg}
\@citedef{ChKe66}{ChKe66}
\@citedef{ChKe62}{ChKe62}
\@citedef{CoSh:919}{CoSh:919}
\@citedef{DvSh:65}{DvSh 65}
\@citedef{Di}{Di}
\@citedef{Eh57}{Eh57}
\@citedef{Fr75}{Fr75}
\@citedef{GiSh:577}{GiSh 577}
\@citedef{GbTl06}{GbTl06}
\@citedef{Gr91}{Gr91}
\@citedef{GrHa89}{GrHa89}
\@citedef{GIL02}{GIL02}
\@citedef{GrLe0x}{GrLe0x}
\@citedef{GrLe00a}{GrLe00a}
\@citedef{GrLe02}{GrLe02}
\@citedef{GrSh:259}{GrSh 259}
\@citedef{GrSh:174}{GrSh 174}
\@citedef{GrSh:238}{GrSh 238}
\@citedef{GrSh:222}{GrSh 222}
\@citedef{GrVa0xa}{GrVa0xa}
\@citedef{GrVa0xb}{GrVa0xb}
\@citedef{Ha61}{Ha61}
\@citedef{HHL00}{HHL00}
\@citedef{HaSh:323}{HaSh 323}
\@citedef{He74}{He74}
\@citedef{HeIo02}{HeIo02}
\@citedef{He92}{He92}
\@citedef{HuSh:342}{HuSh 342}
\@citedef{Hy98}{Hy98}
\@citedef{HySh:474}{HySh 474}
\@citedef{HySh:529}{HySh 529}
\@citedef{HySh:632}{HySh 632}
\@citedef{HySh:602}{HySh 602}
\@citedef{HySh:629}{HySh 629}
\@citedef{HySh:676}{HySh 676}
\@citedef{HShT:428}{HShT 428}
\@citedef{HyTu91}{HyTu91}
\@citedef{JrSh:875}{JrSh 875}
\@citedef{J}{J}
\@citedef{Jn56}{Jn56}
\@citedef{Jo56}{Jo56}
\@citedef{Jn60}{Jn60}
\@citedef{Jo60}{Jo60}
\@citedef{KM67}{KM67}
\@citedef{Ke70}{Ke70}
\@citedef{Ke71}{Ke71}
\@citedef{KiPi98}{KiPi98}
\@citedef{Kin66}{Kin66}
\@citedef{KjSh:409}{KjSh 409}
\@citedef{KlSh:362}{KlSh 362}
\@citedef{KoSh:796}{KoSh 796}
\@citedef{Las88}{Las88}
\@citedef{LwSh:871}{LwSh 871}
\@citedef{LwSh:489}{LwSh 489}
\@citedef{LwSh:560}{LwSh 560}
\@citedef{LwSh:687}{LwSh 687}
\@citedef{Lv71}{Lv71}
\@citedef{Le0x}{Le0x}
\@citedef{Le0y}{Le0y}
\@citedef{McSh:55}{McSh 55}
\@citedef{MaSh:285}{MaSh 285}
\@citedef{Mw85a}{Mw85a}
\@citedef{Mw85}{Mw85}
\@citedef{MkSh:366}{MkSh 366}
\@citedef{MiRa65}{MiRa65}
\@citedef{MoVa62}{MoVa62}
\@citedef{Mo65}{Mo65}
\@citedef{Mo70}{Mo70}
\@citedef{Pi0x}{Pi0x}
\@citedef{RuSh:117}{RuSh 117}
\@citedef{Sc76}{Sc76}
\@citedef{Sh:88r}{Sh 88r}
\@citedef{Sh:E12}{Sh:E12}
\citation{Sh:g}
\@citedef{Sh:F888}{Sh:F888}
\@citedef{Sh:322}{Sh 322}
\@citedef{Sh:482}{Sh 482}
\@citedef{Sh:922}{Sh:922}
\@citedef{Sh:F820}{Sh:F820}
\@citedef{Sh:832}{Sh 832}
\@citedef{Sh:840}{Sh 840}
\@citedef{Sh:e}{Sh:e}
\@citedef{Sh:F782}{Sh:F782}
\@citedef{Sh:800}{Sh 800}
\@citedef{Sh:F841}{Sh:F841}
\@citedef{Sh:F735}{Sh:F735}
\@citedef{Sh:842}{Sh 842}
\@citedef{Sh:839}{Sh 839}
\@citedef{Sh:705}{Sh 705}
\@citedef{Sh:868}{Sh 868}
\@citedef{Sh:1}{Sh 1}
\@citedef{Sh:3}{Sh 3}
\@citedef{Sh:10}{Sh 10}
\@citedef{Sh:16}{Sh 16}
\@citedef{Sh:31}{Sh 31}
\@citedef{Sh:48}{Sh 48}
\@citedef{Sh:46}{Sh 46}
\@citedef{Sh:43}{Sh 43}
\@citedef{Sh:54}{Sh 54}
\@citedef{Sh:56}{Sh 56}
\@citedef{Sh:a}{Sh:a}
\@citedef{Sh:108}{Sh 108}
\@citedef{Sh:93}{Sh:93}
\@citedef{Sh:b}{Sh:b}
\@citedef{Sh:132}{Sh 132}
\@citedef{Sh:87a}{Sh 87a}
\@citedef{Sh:87b}{Sh 87b}
\@citedef{Sh:202}{Sh 202}
\@citedef{Sh:200}{Sh 200}
\@citedef{Sh:205}{Sh 205}
\@citedef{Sh:197}{Sh 197}
\@citedef{Sh:155}{Sh 155}
\@citedef{Sh:88a}{Sh 88a}
\@citedef{Sh:88}{Sh 88}
\@citedef{Sh:220}{Sh 220}
\@citedef{Sh:225}{Sh 225}
\@citedef{Sh:300}{Sh 300}
\@citedef{Sh:225a}{Sh 225a}
\@citedef{Sh:c}{Sh:c}
\@citedef{Sh:284c}{Sh 284c}
\@citedef{Sh:429}{Sh 429}
\@citedef{Sh:351}{Sh 351}
\@citedef{Sh:420}{Sh 420}
\@citedef{Sh:g}{Sh:g}
\@citedef{Sh:430}{Sh 430}
\@citedef{Sh:f}{Sh:f}
\@citedef{Sh:394}{Sh 394}
\@citedef{Sh:620}{Sh 620}
\@citedef{Sh:589}{Sh 589}
\@citedef{Sh:460}{Sh 460}
\@citedef{Sh:576}{Sh 576}
\@citedef{Sh:472}{Sh 472}
\@citedef{Sh:603}{Sh 603}
\@citedef{Sh:715}{Sh 715}
\@citedef{Sh:829}{Sh 829}
\@citedef{ShHM:158}{ShHM 158}
\@citedef{ShUs:837}{ShUs 837}
\@citedef{ShVi:648}{ShVi 648}
\@citedef{ShVi:635}{ShVi 635}
\@citedef{Sh:E45}{Sh:E45}
\@citedef{Sh:E54}{Sh:E54}
\@citedef{Sh:E56}{Sh:E56}
\@citedef{Sh:F709}{Sh:F709}
\@citedef{Sh:E36}{Sh:E36}
\@citedef{Str76}{Str76}
\@citedef{Va02}{Va02}
\@citedef{Zi0xa}{Zi0xa}
\@citedef{Zi0xb}{Zi0xb}

\endgroup
      \immediate\openout\@auxfile = \jobname.aux
   \fi
}%
%
%
\newif\if@auxfiledone
\ifx\noauxfile\@undefined \else \@auxfiledonetrue\fi
%
%
%
%
\@innernewwrite\@auxfile
\def\@writeaux#1{\ifx\noauxfile\@undefined \write\@auxfile{#1}\fi}%
%
%
%
\ifx\@undefinedmessage\@undefined
   \def\@undefinedmessage{No .aux file; I won't give you warnings about
                          undefined citations.}%
\fi
%
%
\@innernewif\if@citewarning
\ifx\noauxfile\@undefined \@citewarningtrue\fi
%
%
%
\catcode`@ = \@oldatcatcode

\def\pfeilso{\leavevmode
            \vrule width 1pt height9pt depth 0pt\relax
           \vrule width 1pt height8.7pt depth 0pt\relax
           \vrule width 1pt height8.3pt depth 0pt\relax
           \vrule width 1pt height8.0pt depth 0pt\relax
           \vrule width 1pt height7.7pt depth 0pt\relax
            \vrule width 1pt height7.3pt depth 0pt\relax
            \vrule width 1pt height7.0pt depth 0pt\relax
            \vrule width 1pt height6.7pt depth 0pt\relax
            \vrule width 1pt height6.3pt depth 0pt\relax
            \vrule width 1pt height6.0pt depth 0pt\relax
            \vrule width 1pt height5.7pt depth 0pt\relax
            \vrule width 1pt height5.3pt depth 0pt\relax
            \vrule width 1pt height5.0pt depth 0pt\relax
            \vrule width 1pt height4.7pt depth 0pt\relax
            \vrule width 1pt height4.3pt depth 0pt\relax
            \vrule width 1pt height4.0pt depth 0pt\relax
            \vrule width 1pt height3.7pt depth 0pt\relax
            \vrule width 1pt height3.3pt depth 0pt\relax
            \vrule width 1pt height3.0pt depth 0pt\relax
            \vrule width 1pt height2.7pt depth 0pt\relax
            \vrule width 1pt height2.3pt depth 0pt\relax
            \vrule width 1pt height2.0pt depth 0pt\relax
            \vrule width 1pt height1.7pt depth 0pt\relax
            \vrule width 1pt height1.3pt depth 0pt\relax
            \vrule width 1pt height1.0pt depth 0pt\relax
            \vrule width 1pt height0.7pt depth 0pt\relax
            \vrule width 1pt height0.3pt depth 0pt\relax}

\def\pfeilsw{ \leavevmode 
            \vrule width 1pt height0.3pt depth 0pt\relax
            \vrule width 1pt height0.7pt depth 0pt\relax
            \vrule width 1pt height1.0pt depth 0pt\relax
            \vrule width 1pt height1.3pt depth 0pt\relax
            \vrule width 1pt height1.7pt depth 0pt\relax
            \vrule width 1pt height2.0pt depth 0pt\relax
            \vrule width 1pt height2.3pt depth 0pt\relax
            \vrule width 1pt height2.7pt depth 0pt\relax
            \vrule width 1pt height3.0pt depth 0pt\relax
            \vrule width 1pt height3.3pt depth 0pt\relax
            \vrule width 1pt height3.7pt depth 0pt\relax
            \vrule width 1pt height4.0pt depth 0pt\relax
            \vrule width 1pt height4.3pt depth 0pt\relax
            \vrule width 1pt height4.7pt depth 0pt\relax
            \vrule width 1pt height5.0pt depth 0pt\relax
            \vrule width 1pt height5.3pt depth 0pt\relax
            \vrule width 1pt height5.7pt depth 0pt\relax
            \vrule width 1pt height6.0pt depth 0pt\relax
            \vrule width 1pt height6.3pt depth 0pt\relax
            \vrule width 1pt height6.7pt depth 0pt\relax
            \vrule width 1pt height7.0pt depth 0pt\relax
            \vrule width 1pt height7.3pt depth 0pt\relax
            \vrule width 1pt height7.7pt depth 0pt\relax
            \vrule width 1pt height8.0pt depth 0pt\relax
            \vrule width 1pt height8.3pt depth 0pt\relax
            \vrule width 1pt height8.7pt depth 0pt\relax
            \vrule width 1pt height9pt depth 0pt\relax
      }


\def\widestnumber#1#2{}

\def\citewarning#1{\ifx\shlhetal\relax 
    \else
    \par{#1}\par
    \fi
}

\def\rm{\fam0 \tenrm}

\def\fakesubhead#1\endsubhead{\bigskip\noindent{\bf#1}\par}



%
%
%

%

\font\textrsfs=rsfs10
\font\scriptrsfs=rsfs7
\font\scriptscriptrsfs=rsfs5

\newfam\rsfsfam
\textfont\rsfsfam=\textrsfs
\scriptfont\rsfsfam=\scriptrsfs
\scriptscriptfont\rsfsfam=\scriptscriptrsfs

\edef\oldcatcodeofat{\the\catcode`\@}
\catcode`\@11

\def\Cal@@#1{\noaccents@ \fam \rsfsfam #1}

\catcode`\@\oldcatcodeofat


\expandafter\ifx \csname margininit\endcsname \relax\else\margininit\fi

\long\def\red#1\endred{}
\long\def\green#1\endgreen{}
\long\def\blue#1\endblue{}
\long\def\private#1\endprivate{}

\def\endred{ \unmatched endred! }
\def\endgreen{ \unmatched endgreen! }
\def\endblue{ \unmatched endblue! }
\def\endprivate{ \unmatched endprivate! }

\ifx\latexcolors\undefinedcs\def\latexcolors{}\fi

\def\emptycs{}
\def\evaluatelatexcolors{%
        \ifx\latexcolors\emptycs\else
        \expandafter\xxevaluate\latexcolors\xxfertig\evaluatelatexcolors\fi}
\def\xxevaluate#1,#2\xxfertig{\setupthiscolor{#1}%
        \def\latexcolors{#2}}


\font\smallfont=cmsl7
\def\rutgerscolor{\ifmmode\else\endgraf\fi\smallfont
\advance\leftskip0.5cm\relax}
\def\setupthiscolor#1{\edef\tmptmpcs{\noexpand\bgroup\noexpand\rutgerscolor
\noexpand\def\noexpand\currentcolor{#1}%
\noexpand}%
\expandafter\let\csname#1\endcsname\tmptmpcs
\def\tmptmpcs{\checkColorUnmatched{#1}\popthecolor}
\expandafter\let\csname end#1\endcsname\tmptmpcs}

\def\checkColorUnmatched#1{\def\expectcolor{#1}%
    \ifx\expectcolor\currentcolor   
    \else \edef\failhere{\noexpand\tryingToClose '\currentcolor' with end\expectcolor}\failhere\fi}

\def\currentcolor{???}

\def\popthecolor{\ifmmode\else\endgraf\fi\egroup}

\expandafter\def\csname#1\endcsname{}

\evaluatelatexcolors

 \let\outerhead\head
 \def\head{\innerhead}
 \let\innerhead\outerhead

 \let\outersubhead\subhead
 \def\subhead{\innersubhead}
 \let\innersubhead\outersubhead

 \let\outersubsubhead\subsubhead
 \def\subsubhead{\innersubsubhead}
 \let\innersubsubhead\outersubsubhead

 \let\outerproclaim\proclaim
 \def\proclaim{\innerproclaim}
 \let\innerproclaim\outerproclaim

 %
 %
 %
 %

\def\demo#1{\medskip\noindent{\it #1.\/}}
\def\enddemo{\smallskip}

\def\remark#1{\medskip\noindent{\it #1.\/}}
\def\endremark{\smallskip}

%
%

\pageheight{8.5truein}
 \topmatter
 \title{Non-structure in $\lambda^{++}$ using instances of WGCH} \endtitle

 \author {Saharon Shelah \thanks {\null\newline I would like to thank 
 Alice Leonhardt for the beautiful typing. \null\newline
 The author would like to thank the Israel Science Foundation for
 partial support of this research (Grant No. 242/03). Publication
 838. } \endthanks} \endauthor  

 \affil{The Hebrew University of Jerusalem \\
 Einstein Institute of Mathematics \\
 Edmond J. Safra Campus, Givat Ram \\
 Jerusalem 91904, Israel
  \medskip
  Department of Mathematics \\
  Hill Center-Busch Campus \\
   Rutgers, The State University of New Jersey \\
  110 Frelinghuysen Road \\
  Piscataway, NJ 08854-8019 USA} \endaffil

 \endtopmatter
%
\document

\head {\S0 Introduction} \endhead  \resetall \sectno=0
 \spuriousreset
\bigskip

Our aim is to prove the results of the form ``build
complicated/many models of cardinality $\lambda^{++}$ by approximation
of cardinality $\lambda$" assuming only $2^\lambda <
2^{\lambda^+} < 2^{\lambda^{++}}$, which are needed in developing
classification of a.e.c., i.e. in this book and the related works
(this covers \cite{Sh:87b}, \cite{Sh:88} redone in \xCITE{88r} 
and \cite{Sh:576}, \cite{Sh:603} which are redone in
\chaptercite{E46} + \xCITE{838} and \cite{Sh:87b} which is redone by
\xCITE{600}, \xCITE{705}, \cite{Sh:842}, so we ignore, e.g. 
\cite{Sh:576} now), 
fulfilling  promises, uniformizing and correcting inaccuracies there
and doing more.  But en-route we spend time on the structure side.

As in \cite[\S3]{Sh:576} we consider a version of construction framework, 
trying to give sufficient conditions for constructing many models 
of cardinality $\partial^+$ by approximations of cardinality $<
\partial$ so $\lambda^+$ above correspond to $\partial$.  Compared
to \cite[\S3]{Sh:576}, the present version is hopefully more
transparent.
\nl
We start in \S1,\S2 (and also \S3) by giving several sufficient conditions for
non-structure, in a framework closer to the applications we have in
mind than \cite[\S3]{Sh:576}.  The price is delaying the actual
proofs and losing some generality.  
Later (mainly in \S4, but also in \S6 and \S8) we do the 
applications, usually each is quoted
(in some way) elsewhere.  Of course, it is a delicate question
how much should we repeat the background which exists when the quote
was made.

The ``many" is interpreted as $\ge
\mu_{\text{unif}}(\lambda^{++},2^{\lambda^+})$, see \scite{838-7f.14} why
this is almost equal to $2^{\lambda^{++}}$.  Unfortunately, there is
here no one theorem covering all cases.  But if a ``lean" version
suffice for us, which means that we assume the very
weak set theoretic 
assumption ``WDmId$_{\lambda^+}$ (a normal ideal on $\lambda^+$) is not
$\lambda^{++}$-saturated", (and, of course, we are content with getting $\ge
\mu_{\text{unif}}(\lambda^{++},2^{\lambda^+}))$ \ub{then} all the
results can be deduced from weak coding, i.e. Theorem \scite{838-2b.3}.  In
this case, some parts are redundant and the paper is neatly divided to two:
structure part and non-structure part
and we do now describe this.

First, in \S1 we define a 
(so called) nice framework to deal with such theorems, and in the
beginning of \S2 state the theorem but we replace $\lambda^+$ by a regular
uncountable cardinal $\partial$.  This is done in a way closed to
the applications we have in mind as deduced in \S4.
In Theorem 
\scite{838-2b.3} the model in $\partial^+$ is approximated triples by $(\bar
M^\eta,\bar{\bold J}^\eta,\bold f^\eta)$ for $\eta \in {}^{\partial^+
>}2$ increasing with $\eta$ where $\bar M^\eta$ is an 
increasing chain of length $\partial$ of models of cardinality $<
\partial$ and for each $\eta \in {}^{\partial^+}2$ the sequence
$\langle \cup\{M^{\eta \rest \varepsilon}_\alpha:\alpha
< \partial\}:\varepsilon < \partial^+\rangle$ is increasing; 
similarly in the other such theorems.

Theorem \scite{838-2b.3} is not proved in \S2.  It is proved in \S9, \S10,
specifically in \scite{838-10k.17}.  Why?  In the proof we apply relevant set
theoretic results (see in the end of \S0 and more in
\S9 on weak diamond and failure of strong
uniformization), for this it is helpful to decide that
the universe of each model (approximating the desired one) 
is $\subseteq \partial^+$ and to add commitments $\bar{\Bbb F}$
on the amalgamations
used in the construction called amalgamation choice function.

So model theoretically they look artificial though the theorems are stronger.

Second, we deal with the applications in \S4, actually in
\S4(A),(C),(D),(E), so we have in 
each case to choose ${\frak u}$, the
construction framework and prove the required properties.  By our choice this
goes naturally.  

But we would like to eliminate the extra assumption ``WDmId$_{\lambda^+}$
is $\lambda^{++}$-saturated".  So in \S2 and \S3 there are additional
``coding" theorems.  Some still need the ``amalgamation choice
function", others, as we have a stronger model theoretic assumption do
not need such function so their proof is not delayed to \S10.

Probably the most interesting case is proving the density of
$K^{3,\text{uq}}_{\frak s}$, i.e. of uniqueness triples $(M,N,a) \in
K^{3,\text{bs}}_{\frak s}$ for ${\frak s}$ an (almost) good
$\lambda$-frame, a somewhat weaker version of the (central) notion of
\xCITE{600}.  Ignoring for a minute the ``almost" this is an
important step in (and is promised in) \sectioncite[\S5]{600}.  
The proof is done in two stages.  In the first stage we consider, in
\S6,  a wider class 
$K^{3,\text{up}}_{\frak s} \subseteq K^{3,\text{bs}}_{\frak s}$ than
$K^{3,\text{uq}}_{\frak s}$ and prove that its failure to be dense implies
non-structure.  This is done in \S6, the proof is easier when 
$\lambda > \aleph_0$ or
at least ``${\Cal D}_{\lambda^+}$ is not $\lambda^{++}$-saturated";
(but this  is unfortunate for an application to \xCITE{88r}).
But for the proofs in \S6 we
need before this in \S5 to prove some
``structure positive theory" claims even if ${\frak s}$ is a good
$\lambda$-frame; we need more in the almost good case.

So naturally we assume (categoricity in $\lambda$ and) density of
$K^{3,\text{up}}_{\frak s}$ and prove (in \S7) that WNF$_{\frak s}$, a
weaker relative of NF$_{\frak s}$, is a weak ${\frak s}$-non-forking relation
on ${\frak K}_{\frak s}$ respecting ${\frak s}$ and that ${\frak s}$
is actually a  good $\lambda$-frame (see \scite{838-7v.37}(1)); 
both results are helpful.
 
The second stage (in \S8) is done in two substages.  In the first
substage we deal with a delayed version of uniqueness, proving that its
failure implies non-structure.  In
the second substage we assume delayed uniqueness but
$K^{3,\text{uq}}_{\frak s}$ is not dense and we get another non-structure
\ub{but} relying on a positive consequence of density of
$K^{3,\text{up}}_{\frak s}$ (that is, on a weak form of NF, see \S7).

Why do we deal with almost good $\lambda$-frames?  
By \sectioncite[\S3]{600}(E) from an a.e.c. ${\frak K}$
categorical in $\lambda,\lambda^+$ which have $\in
[1,\mu_{\text{unif}}(\lambda^{++},2^{\lambda^+}))$ non-isomorphic
models in $\lambda^{++}$ we construct a good $\lambda^+$-frame ${\frak
s}$ with ${\frak K}^{\frak s} \subseteq {\frak K}$.
The non-structure theorem stated (and used) there is fully proven in \S4.
However, not only do we
use the $\lambda^{++}$-saturation of the weak diamond ideal on
$\lambda^+$, but it is a good
$\lambda^+$-frame rather than a good $\lambda$-frame.

This does not hamper us in \xCITE{705} but still is regretable.
Now we ``correct" this but the price is getting an \ub{almost} good
$\lambda$-frame, noting that such ${\frak s}$ is 
proved to exist in \sectioncite[\S8]{E46}, the revised version of \cite{Sh:576}.
However, to arrive to those points in \chaptercite{E46}, \cite{Sh:576} we
have to prove the density of minimal types under the weaker
assumptions, i.e. without the saturation of the ieal
WDmId$_{\lambda^+}$ together
\xCITE{838} + \sectioncite[\S3,\S4]{E46} gives a full proof.  This
requires again on developing some positive theory, so in \S5 we do here some
positive theory.  Recall that \cite{Sh:603}, \cite{Sh:576} are
subsumed by them.

We can note that in building models $M \in K^{\frak s}_{\lambda^{++}}$ for
${\frak s}$ an almost good $\lambda$-frame, for convenience we use
disjoint amalgamation.  This may seem harmless but proving the 
density of the minimal
triples this is not obvious; without assuming this we have to use
$\langle M_\alpha,h_{\alpha,\beta}:\alpha < \beta < \partial\rangle$
instead of increasing $\langle M_\alpha:\alpha < \partial\rangle$; a
notationally 
cumbersome choice.  So we use a congruence relation $=_\tau$ but the
models we construct are not what we need.  We have to take their
quotient by $=_\tau$, which has to, e.g. have the right cardinality.
But we can take care that $|M|/ \equiv_\tau$ has cardinality $\lambda^{++}$ 
and $a \in M \Rightarrow
|a/\equiv_\tau| = \lambda^{++}$.  For the almost good $\lambda$-frame case
this follows if we use not just models $M$ which
are $\lambda^+$-saturated above $\lambda$ but if $M_0
\le_{{\frak K}[{\frak s}]} M_1 \le_{{\frak K}[{\frak s}]} M,M_0 \in
K_{\frak s},M_1 \in K^{\frak s}_{\le \lambda}$ and $p \in {\Cal
S}^{\text{bs}}_{\frak s}(M_0)$ then for some $a \in M$, for every
$M'_1 \le_{{\frak K}[{\frak s}]} M_1$ of cardinality $\lambda$
including $M_0$, the type $\ortp_{\frak s}(a,M'_1,M)$ is the non-forking
extension of $p$, so not a real problem.
\bn
\ub{Reading Plans}:  \ub{The miser model - theorist Plan A}:

If you like to see only the results \ub{quoted} elsewhere 
(in this book), willing to
assume an extra weak set theoretic assumption, this is the plan for you.

The results are all in \S4, more exactly \S4(A),(C),(D),(E).  They all
need only \scite{838-2b.3} relying on \scite{838-2b.1}, but the
rest of \S2 and \S3 are irrelevant as well as \S5 - \S8.

To understand what \scite{838-2b.3} say you have to read \S1 (what is ${\frak
u}$; what are ${\frak u}$-free rectangles; assuming ${\frak K}$ is
categorical in $\lambda^+$ you can ignore the ``almost").  You may
take \scite{838-2b.3} on belief, so you are done; otherwise you have to
see \S9 and \scite{838-10k.1} - \scite{838-10k.17}.
\bn
\ub{The pure model - theorist Plan B}:

Suitable if you like to know about the relatives of ``good $\lambda$-frames".
Generally see \S5 - \S8.

In particular on ``almost good $\lambda$-frames" see \S5; but better
first read \scite{838-1a.1} - \scite{838-1a.27}, which deal with a related
framework called ``nice construction framework" and in \S6 learn of the class
$K^{3,\text{up}}_{\frak s} \subseteq K^{3,\text{bs}}_{\frak s}$ with a
weak version of uniqueness.  By quoting we get non-structure if they fail
density.  Then in \S7 learn on weak non-forking relations WNF on
${\frak K}_\lambda$ which respects ${\frak s}$, it is interesting when
we assume $K^{3,\text{up}}_{\frak s}$ has density
or reasonably weak existence assumption, because then
we can prove that the definition given such existence, 
and this implies that ${\frak s}$ is a good $\lambda$-frame
(not just almost).  In \S8 we prove density of uinqueness triples
$(K^{3,\text{uq}}_{\frak s})$ in $K^{3,\text{bs}}_{\frak s}$, so quote
non-structure theorems.
\bn
\ub{The set theorist Plan C}:

Read \S1,\S2,\S3,\S9,\S10,\S11 this presents construction in $\partial^+$
by approximation of cardinality $\le \lambda$.
\bn
\margintag{0z.1}\ub{\stag{838-0z.1} Notation}:
\mn
1) ${\frak u}$, a construction framework, see \S1, in particular
Definition \scite{838-1a.3}
\sn
2) Triples $(M,N,\bold J) \in \text{ FR}^\ell_{\frak u}$, see Definition
\scite{838-1a.3}
\sn
2A) $\bold J$ and $\bold I$ are $\subseteq M \in K_{\frak u}$
\sn
3) $\bold d$, and also $\bold e$, a 
${\frak u}$-free rectangle (or triangle), 
see Definitions \scite{838-1a.7}, \scite{838-1a.11}
\sn
4) $K^{\text{qt}}_{\frak u}$, the set of $(\bar M,\bar{\bold J},\bold f)$,
see Definition \scite{838-1a.29}, where, in particular:

\quad 4A) $\bold f$ a function from $\partial$ to $\partial$

\quad 4B) $\bar{\bold J} = 
\langle \bold J_\alpha:\alpha < \partial\rangle,\bold
J_\alpha \subseteq M_{\alpha +1} \backslash M_\alpha$

\quad 4C) $\bar M = \langle M_\alpha:
\alpha < \partial\rangle$ where $M_\alpha \in
{\frak K}_{< \partial}$ is $\le_{\frak K}$-increasing continuous
\sn
5) Orders (or relations) on $K^{\text{qt}}_{\frak u}:\le^{\text{at}}_{\frak u},
\le^{\text{qt}}_{\frak u},\le^{\text{qs}}_{\frak u}$
\sn
6) $\bold c$, a colouring (for use in weak diamond)
\sn
7) $\Bbb F$ (usually $\bar{\Bbb F}$),
for amalgamation choice functions, see Definition \scite{838-10k.5}
\sn
8) ${\frak g}$, a function from $K^{\text{qt}}_{\frak u}$ to itself,
etc., see Definition \scite{838-1a.43} 
\nl

\quad (for defining ``almost every ...")
\nl
9) Cardinals $\lambda,\mu,\chi,\kappa,\theta,\partial$, \ub{but} here 
$\partial = \text{ cf}(\partial) >
\aleph_0$, see \scite{838-1a.15}(2) in \scite{838-1a.15}(1B), 
${\Cal D}_\partial$ is the club filter on $\partial$
\sn
10) $\bold F$ in the definition of limit model, see Definition
\xCITE{88r}, marginal.
\bigskip

\definition{\stag{838-0z.3} Definition}  1) For $K$ a set or a class of
 models let $\dot I(K) = \{M/\cong:M \in K\}$, so it is a cardinality
 or $\infty$.
\nl
2) For a class $K$ of models let $\dot I(\lambda,K) = \dot
I(K_\lambda)$ where $K_\lambda = \{M \in K:\|M\| = \lambda\}$.
\nl
3) For ${\frak K} = (K,\le_{\frak K})$ let $\dot I({\frak K}) = I(K)$ and
$\dot I(\lambda,{\frak K}) = \dot I(\lambda,K)$.
\enddefinition
\bigskip

\remark{Remark}  We shall use in particular
$\dot I({\frak K}^{{\frak u},{\frak h}}_{\partial^+})$, see Definition
\scite{838-1a.45}. 
\endremark
\bn
We now define some set theoretic notions
(we use mainly the ideal WDmTId$_\partial$ and the cardinals 
$\mu_{\text{wd}}(\partial),\mu_{\text{unif}}(\partial^+,2^\partial)$).
\bigskip

\definition{\stag{838-0z.5} Definition}  Fix $\partial$ regular and
uncountable.
\nl
1) For $\partial$ regular uncountable, $S \subseteq \partial$ and
$\bar \chi = \langle \chi_\alpha:\alpha < \partial \rangle$ but only
$\bar \chi \rest S$ matters so we can use any $\bar \chi 
= \langle \chi_\alpha:\alpha \in S'\rangle$ where $S \subseteq S'$
let

$$
\align
\text{WDmTId}(\partial,S,\bar \chi) = \biggl\{ A:&A \subseteq 
\dsize \prod_{\alpha \in S} \chi_\alpha,\text{ and for some 
function (= colouring)} \\  
  &\bold c \text{ with domain } \dsize \bigcup_{\alpha < \partial} {}^\alpha
(2^{< \partial}) \text{ mapping }
{}^\alpha(2^{< \partial}) \text{ into } \chi_\alpha, \\
  &\text{ for every } \eta \in A, \text{ for some } f \in 
{}^\partial(2^{< \partial}) \text{ the set} \\
  &\{ \delta \in S:\eta (\delta) = \bold c(f \restriction \delta)\} 
\text{ is not stationary (in } \partial) \biggr\}.
\endalign
$$
\mn
(Note: WDmTId stands for weak diamond target ideal; 
of course, if we increase the $\chi_\alpha$ we get a bigger
ideal); the main case is
when $\alpha \in S \Rightarrow \chi_\alpha = 2$ this is the
weak diamond, see below.
\nl
1A) Here we can 
replace $2^{< \partial}$ by any set of this cardinality, and so
we can replace $f \in {}^\partial(2^{< \partial})$ by $f_1,\dotsc,f_n \in
{}^\partial(2^{< \partial})$ and $f \restriction \delta$ by $\langle
f_1 \restriction \delta,\dotsc,f_n \restriction \delta \rangle$ and
$\bold c(f \restriction \delta)$ by $\bold c'(f_1 \restriction
\delta,\dotsc,f_n \restriction \delta)$ so with $\bold c'$ being an 
$n$-place function; justified in \cite[AP,\S1]{Sh:f}.
\nl
2)\footnote{in \cite[AP,\S1]{Sh:b}, \cite[AP,\S1]{Sh:f} we express 
cov$_{\text{wdmt}}(\partial,S) > \mu^*$ by allowing $f(0) \in \mu^* <
\mu$}

$$
\text{cov}_{\text{wdmt}}(\partial,S,\bar \chi) = 
\text{ Min}\biggl\{ |{\Cal P}|:{\Cal P} \subseteq \text{
WDmTId}(\partial,S,\bar \chi) 
\text{ and } \dsize \prod_{\alpha \in S} \chi_\alpha  
\subseteq \dsize \bigcup_{A \in {\Cal P}} A \biggr\} 
$$

$$
\align
\text{WDmTId}_{< \mu}(\partial,S,\bar \chi) = 
\biggl\{ A \subseteq \dsize \prod_{\alpha \in S} 
\chi_\alpha:&\text{ for some } 
i^* < \mu \text{ and} \tag"{$3)$}" \\
  &A_i \in \text{ WDmTId}(\partial,S,\bar \chi) 
\text{ for} \\
  &\,i < i^* \text{ we have } A \subseteq \dsize \bigcup_{i < i^*} A_i
\biggr\}
\endalign
$$

$$
\text{WDmTId}(\partial) = \text{WDmTId}(\partial,\partial,2) \tag
"{$4)\,\,(a)$}" 
$$

$$
\text{WDmId}_{< \mu}(\partial,\bar \chi) = \{S \subseteq \partial:
\text{cov}_{\text{wdmt}}(\partial,S,\bar \chi) < \mu\}. 
\tag "{$4)\,\,(b)$}"
$$
\mn
5) Instead of ``$< \mu^+$" we may write $\le \mu$ or just $\mu$;
if we omit $\mu$ we mean $(2^{< \partial})$.  
If $\bar \chi$ is constantly $2$ we may omit it, see below, if
$\chi_\alpha = 2^{|\alpha|}$ we may write pow instead of $\bar \chi$;
all this in the parts above and below.
\mn
6) Let $\mu_{\text{wd}}(\partial,\bar \chi) = 
\text{ cov}_{\text{wdmt}}(\partial,\partial,\bar \chi)$.
\mn
7)  We say that the weak diamond holds on $\lambda$ \ub{if} 
$\partial \notin \text{ WDmId}(\partial)$.
\enddefinition
\bigskip

\remark{Remark}  This is used in \scite{E46-3c.22} \yCITE[3c.22]{E46}, \yCITE[6f.13]{E46}.  Note that by
\scite{838-0z.7}(1A) that $\mu_{\text{wd}}(\lambda^+)$ is large (but $\le
2^{\lambda^+}$, of course).
\endremark
\bn
A relative is
\definition{\stag{838-0z.6} Definition}  Fix $\partial$ regular and
uncountable.
\nl
1) For $\partial$ regular uncountable, $S \subseteq \partial$ and
$\bar \chi = \langle \chi_\alpha:\alpha < \partial \rangle$ let

$$
\align
\text{UnfTId}(\partial,S,\bar \chi) = \biggl\{ A:&A \subseteq 
\dsize \prod_{\alpha \in S} \chi_\alpha,\text{ and for some 
function (= colouring)} \\  
  &\bold c \text{ with domain } \dsize \bigcup_{\alpha < \partial} {}^\alpha
(2^{< \partial}) \text{ mapping }
{}^\alpha(2^{< \partial}) \text{ into } \chi_\alpha, \\
  &\text{ for every } \eta \in A, \text{ for some } f \in 
{}^\partial(2^{< \partial}) \text{ the set} \\
  &\{ \delta \in S:\eta (\delta) \ne \bold c(f \restriction \delta)\} 
\text{ is not stationary (in } \partial) \biggr\}.
\endalign
$$
\mn
(Note: UnfTId stands for uniformization target ideal; of course, if 
we increase the $\chi_\alpha$ we get a smaller
ideal); when $\alpha \in S \Rightarrow \chi_\alpha = 2$ this is the
weak diamond, i.e. as in \scite{838-0z.5}(1), similarly below.
\nl
1A) Also here we can 
replace $2^{< \partial}$ by any set of this cardinality, and so
we can replace $f \in {}^\partial(2^{< \partial})$ by $f_1,\dotsc,f_n \in
{}^\partial(2^{< \partial})$ and $f \restriction \delta$ by $\langle
f_1 \restriction \delta,\dotsc,f_n \restriction \delta \rangle$ and
$\bold c(f \restriction \delta)$ by $\bold c'(f_1 \restriction
\delta,\dotsc,f_n \restriction \delta)$ so with $\bold c'$ being an 
$n$-place function, justified in \cite[AP,\S1]{Sh:f}.

$$
\text{cov}_{\text{unf}}(\partial,S,\bar \chi) = 
\text{ Min}\biggl\{ |{\Cal P}|:{\Cal P} \subseteq 
\text{ UnfTId}(\partial,S,\bar \chi) 
\text{ and } \dsize \prod_{\alpha \in S} \chi_\alpha  
\subseteq \dsize \bigcup_{A \in {\Cal P}} A \biggr\} \tag"{$2)$}"
$$

$$
\align
\text{UnfTId}_{< \mu}(\partial,S,\bar \chi) = 
\biggl\{ A \subseteq \dsize \prod_{\alpha \in S} 
\chi_\alpha:&\text{ for some } 
i^* < \mu \text{ and} \tag"{$3)$}" \\
  &A_i \in \text{ UnfTId}(\partial,S,\bar \chi) 
\text{ for} \\
  &\,i < i^* \text{ we have } A \subseteq \dsize \bigcup_{i < i^*} A_i
\biggr\}
\endalign
$$

$$
\text{UnfId}_{< \mu}(\partial,\bar \chi) = \biggl\{ S \subseteq \partial:
\text{cov}_{\text{unf}}(\partial,S,\bar \chi) < \mu \biggr\}. \tag "{$4)$}"
$$
\mn
5) Instead of ``$< \mu^+$" we may write $\le \mu$ or just $\mu$, 
if we omit $\mu$ we mean $(2^{< \partial})$.  
If $\bar \chi$ is constantly $2$ we may omit it, if
$\chi_\alpha = 2^{|\alpha|}$ we may write pow instead of $\bar \chi$;
all this in the parts above and below.
\nl
6) $\mu_{\text{unif}}(\partial,\bar \chi)$ where $\bar \chi = \langle
   \chi_\alpha:\alpha < \partial\rangle$ is Min$\{|{\Cal P}|:{\Cal P}$
   is a family of subsets of $\dsize \prod_{\alpha < \partial}
   \chi_\alpha$ with union $\dsize \prod_{\alpha < \partial}
   \chi_\alpha$ and for each $A \in {\Cal P}$ there is a function
   $\bold c$ with domain $\dbcu_{\alpha < \partial} \dsize
   \prod_{\beta < \alpha} \chi_\beta$ such that $f \in A \Rightarrow
   \{\delta \in S:\bold c(f \rest \delta) = f(\delta)\}$ is not
   stationary$\}$.
\nl
7) $\mu_{\text{unif}}(\partial,\chi) = \mu_{\text{unif}}(\partial,\bar
   \chi)$ where $\bar \chi = \langle \chi:\alpha < \partial\rangle$
   and $\mu_{\text{unif}}(\partial,< \chi)$ means
   sup$\{\mu_{\text{unif}}(\partial,\chi_1):\chi_1 < \chi\}$;
   similarly in the other definitions above.  If $\chi = 2$ we may
   omit it.
\enddefinition
\bn
By Devlin Shelah \cite{DvSh:65}, \cite[XIV,1.5,1.10]{Sh:b}(2);1.18(2),1.9(2) 
(presented better in \cite[AP,\S1]{Sh:f} we have:
\proclaim{\stag{838-0z.7} Theorem} 
\nl
1)  If $\partial = \aleph_1,
2^{\aleph_0} < 2^{\aleph_1},\mu \le (2^{\aleph_0})^+$
\ub{then} $\partial \notin { \text{\rm WDmId\/}}_{<
\mu}(\partial)$. 
\nl
2) If $2^\theta = 2^{< \partial} < 2^\partial,\mu = (2^\theta)^+$, or just:
for some $\theta,2^\theta = 2^{< \partial} < 2^\partial,\mu \le
2^\partial$, and $\chi^\theta < \mu$ for $\chi < \mu$, \ub{then} 
$\partial \notin { \text{\rm WDmId\/}}_\mu(\partial)$ equivalently
${}^\partial 2 \notin \text{\rm WDmTId}_\mu(\partial)$.  So
$(\mu_{\text{wd}}(\partial))^\theta = 2^\partial$.
\nl
3) Assume $2^\theta = 2^{< \partial} < 2^\partial,\mu \le 2^\partial$ and
\mr
\item "{$(a)$}"   $\mu \le \partial^+$ or ${\text{\rm cf\/}}([\mu_1]
^{\le \partial},\subseteq) < \mu$ for $\mu_1 < \mu$ and
\sn
\item "{$(b)$}"   ${\text{\rm cf\/}}(\mu) > \partial \,
{ \text{\rm or\/}} \, \mu \le (2^{< \partial})^+$.
\ermn
\ub{Then} ${\text{\rm WDmId\/}}_{< \mu}(\partial,\bar \chi)$ is a normal 
ideal on $\partial$ and {\rm WDmTId}$_{< \mu}(\partial,\chi)$ is a
{\rm cf}$(\mu)$-complete ideal on ${}^\partial 2$.   
[If this ideal is not trivial, \ub{then} 
$\partial = { \text{\rm cf\/}}(\partial) > \aleph_0,
2^{< \partial} < 2^\partial$.] 
\nl
4) {\rm WDmTId}$_{< \mu}(\partial,S,\bar \chi)$ is 
{\rm cf}$(\mu)$-complete ideal on $\dsize \prod_{\alpha \in S} \chi_\alpha$.
\endproclaim
\bigskip

\remark{\stag{838-0z.9} Remark}  0) Compare to \S9,\S10, mainly \scite{838-7f.15}.
\nl
1) So if cf$(2^\partial) < \mu$ (which 
holds if $2^\partial$ is singular and $\mu = 2^\partial$) 
then \scite{838-0z.7}(3) implies that there is 
$A \subseteq {}^\partial 2,|A| < 2^\partial,A \notin
\text{ WDmTId}(\partial)$. 
\newline
2) Some related definitions appear in \cite[\S1]{Sh:E45}, 
mainly DfWD$_{< \mu}(\partial)$, but presently we ignore them.
\nl
3) We did not look again at the case $(\forall \sigma < \lambda)(2^\sigma <
2^{< \partial} < 2^\partial)$. 
\nl
4) Recall that for an a.e.c. ${\frak K}$:
\mr
\item "{$(a)$}"   if $K_\lambda \ne \emptyset$ but ${\frak K}$ has no
$\le_{\frak K}$-maximal model in $K_\lambda$ then $K_{\lambda^+} \ne
\emptyset$
\sn
\item "{$(b)$}"   if ${\frak K}$ is categorical in $\lambda$ and
LS$({\frak K}) \le \lambda$ \ub{then} $K_{\lambda^+} \ne \emptyset$
iff $K$ has no $\le_{\frak K}$-maximal model in $K_\lambda$.
\ermn
5) About $\mu_{\text{wd}}(\partial)$ see 
\yCITE[1a.16]{E46}, \yCITE[2b.13]{E46}, \yCITE[6f.4]{E46}.
\endremark
\bigskip

\definition{\stag{838-0z.13} Definition}  1) We say that a normal ideal
$\Bbb I$ on a regular uncountable cardinal $\lambda$ is
$\mu$-saturated \ub{when} we cannot find a sequence $\bar A = \langle
A_i:i < \mu\rangle$ such that $A_i \subseteq \lambda,A_i \notin \Bbb
I$ for $i < \mu$ and $A_i \cap A_j \in \Bbb I$ for $i \ne j < \mu$; if
$\mu \le \lambda^+$ without loss of generality  $A_i \cap A_j \in [\lambda]^{<
\lambda}$.
\nl
2) Similarly for a normal filter on a regular uncountable cardinal $\lambda$.
\enddefinition
\goodbreak

\head {\S1 Nice contruction framework} \endhead  \resetall \sectno=1
 \spuriousreset
\bigskip

We define here when ${\frak u}$ is a nice construction framework.  Now
${\frak u}$ consists of an a.e.c. ${\frak K}$ with LS$({\frak K}) <
\partial_{\frak u} = \text{ cf}(\partial_{\frak u})$, and enables us
to build a model in $K_{\partial^+}$ by approximations of cardinality
$< \partial := \partial_{\frak u}$.

Now for notational reasons we
prefer to use increasing sequences of models rather than directed
systems, i.e., 
sequences like $\langle M_\alpha,f_{\beta,\alpha}:\alpha \le \beta <
\alpha^* \rangle$ with $f_{\beta,\alpha}:M_\alpha \rightarrow M_\beta$
satisfying $f_{\gamma,\beta} \circ f_{\beta,\alpha} = f_{\alpha,\gamma}$
for $\alpha \le \beta \le \gamma < \alpha^*$.
For this it is very desirable to have disjoint amalgamation; however,
in one of the major applications (the density of minimal types,
see here in \S4(A),(B) or in \cite[\S3]{Sh:576} used in
\sectioncite[\S3,\S4]{E46}) we do not
have this.  In \cite[\S3]{Sh:576} the solution was to allow non-standard
interpretation of the equality (see Definition \scite{838-1a.19} here).
Here we choose another formulation: we have
$\tau \subseteq \tau({\frak u})$ such that we are interested in the
non-isomorphism of the $\tau$-reducts $M^{[\tau]}$ of the $M$'s
constructed, see Definition \scite{838-1a.15}.  Of course, this is only a
notational problem.

The main results on such ${\frak u}$ appear later; a major
 theorem is \scite{838-2b.3}, deducing non-structure results
assuming the weak coding property.  This and similar theorems, assuming other
variant of the coding property, are dealt with in \S2,\S3.  They all
have (actually lead to) the form ``if most 
triples $(\bar M,\bar{\bold J},\bold f) \in
K^{\text{qt}}_{\frak u}$ has, in some sense 2 (or many, say 
$2^{< \partial}$) extensions which are (pairwise) incompatible in suitable
sense, \ub{then} we build a suitable tree $\langle(\bar M_\eta,\bar{\bold
J}_\eta,f_\eta):\eta \in {}^{\partial^+ >} 2 \rangle$ and letting
$M_\eta = \cup\{M^\eta_\alpha:\alpha < \lambda\}$ for $\eta \in
{}^{\partial^+ >} 2$ and $M_\nu := \cup\{M_{\nu \restriction
\alpha}:\alpha < \partial^+\}$ for $\nu \in {}^{\partial^+}2$ we
have: among $\langle M_\nu:\nu \in {}^{\partial^+}2\rangle$ many are
non-isomorphic (and in $K_{\partial^+}$).  
Really, usually the indexes are $\eta \in {}^{\partial^+>}(2^\partial)$ and
 the conditions speak on amalgamation in
${\frak K}_{\frak u}$, i.e. on models of cardinality $< \partial$ but
using FR$_1$, FR$_2$, see below.

As said earlier, in the framework defined below we (relatively) prefer
transparency and simplicity on generality, e.g. we can weaken
``${\frak K}_{\frak u}$ is an a.e.c." and/or make FR$^+_\ell$ is
axiomatic and/or use more than atomic successors (see \scite{838-10k.23} +
\scite{838-10k.27}). 

In \scite{838-1a.1} - \scite{838-1a.11} we introduce our frameworks ${\frak
u}$ and ${\frak u}$-free rectangles/triangles; in \scite{838-1a.13} - the
dual of ${\frak u}$, and in \scite{838-1a.15} - \scite{838-1a.21} we justify
the disjoint amalgamation through ``$\tau$ is a ${\frak
u}$-sub-vocabulary", so a reader not bothered by this point can
ignore it, then in \scite{838-1a.24} we consider another property of
${\frak u}$, monotonicity and in \scite{838-1a.27} deal with variants of
${\frak u}$.

In \scite{838-1a.29} - \scite{838-1a.51} we introduce a class
$K^{\text{qt}}_{\frak u}$ of triples $(\bar M,\bar{\bold J},\bold f)$
serving as approximations of size $\partial$, some relations and orders
on it and variants, and define what it means ``for almost every such
triple" (if $K_{\frak u}$ is categorical in $\partial$ this is
usually easy and in many of our applications for most $(\bar
M,\bar{\bold J},\bold f)$ the model 
$\cup\{M_\alpha:\alpha < \partial\}$ is saturated (of cardinality $\partial$).
\bn
\margintag{1a.1}\ub{\stag{838-1a.1} Convention}:  If not said otherwise, ${\frak u}$ 
is as in Definition \scite{838-1a.3}.
\bigskip

\definition{\stag{838-1a.3} Definition}  We say that ${\frak u}$ is a nice 
construction framework \ub{when} (the demands are for $\ell=1,2$ and
later (D) means (D)$_1$ and (D)$_2$ and (E) means (E)$_1$ and (E)$_2$):
\mr
\item "{$(A)$}"  ${\frak u}$ consists of $\partial,{\frak K} =
(K,\le_{\frak K}),\text{FR}_1,\text{FR}_2,\le_1,\le_2$ 
(also denoted by $\partial_{\frak u},{\frak K}^{\frak u} =
{\frak K}^{\text{up}}_{\frak u} = (K^{\text{up}}_{\frak u},\le_{\frak u}),
\text{FR}_1^{\frak u},\text{FR}_2^{\frak u},
\le^1_{\frak u},\le^2_{\frak u}$) and let $\tau_{\frak u} =
\tau_{\frak K}$.   The indexes 1 and 2 can be replaced by ver
(vertical \footnote{Hard but immaterial choice.  We construct a model
of cardinality $\partial^+$ by a sequence of length $\partial^+$
approximations, each of the form $\langle M_\alpha,\bold
J_\alpha:\alpha < \partial\rangle,M_\alpha \in K_{< \partial}$ is
$\le_{{\frak K}_{< \partial}}$-increasing and $(M_\alpha,M_{\alpha
+1},\bold J_\alpha) \in \text{ FR}_2$.  If $\langle M'_\alpha,\bold
J'_\alpha:\alpha < \partial\rangle$ is an immediate successor in the
$\partial^+$-direction of $\langle M_\alpha,\bold J_\alpha:\alpha <
\partial\rangle$ then for most $\alpha,M_\alpha \le_{\frak u}
M'_\alpha$ and $(M_\alpha,M'_\alpha,\bold I_\alpha) \in \text{ FR}_1$
for suitable $\bold I_\alpha$, increasing with $\alpha$ 
and $(M'_\alpha,M'_{\alpha +1},\bold J'_\alpha) \in \text{ FR}_2$ is
$\le^2_{\frak u}$-above
$(M_\alpha,M_{\alpha + 1},\bold J_\alpha)$.  Now the natural order on
FR$_2$ leads in the horizontal direction.}, 
direction of $\partial$) and hor (horizontal, direction of
$\partial^+$) respectively
\sn
\item "{$(B)$}"   $\partial$ is regular uncountable
\sn
\item "{$(C)$}"  ${\frak K} = {\frak K}^{\text{up}}_{\frak u} = 
(K,\le_{\frak K})$ is an a.e.c., $K \ne \emptyset$ of course with
LS$({\frak K}) < \partial$ (or just $(\forall M \in K)(\forall A \in
[M]^{< \partial})(\exists N)(A \subseteq N \le_{\frak K} M \wedge
\|N\| < \partial))$.  Let ${\frak K}_{\frak u} = {\frak K}_{< \partial}
= ({\frak K}_{< \partial},\le_{\frak K} \restriction K_{< \partial})$
where $K_{< \partial} = K_{\frak u} \restriction \{M:M \in
K^{\text{up}}_{\frak u}$ has cardinality $< \partial\}$ and ${\frak
K}[{\frak u}] = K^{\text{up}}_{\frak u}$.
\nl
(To prepare for weaker versions we can start with ${\frak K}_{\frak u}$, a 
$(< \partial)$-a.e.c. 
\footnote{less is used, but natural for our applications, see \S9};
this means $K$ is a class of models of cardinality $< \partial$, and
in AxIII, the existence of union we add the assumption that
the length of the union is $<
\partial$ (here equivalently the union has cardinality $< \partial$)
and we replace ``LS$({\frak K}_{< \partial})$ exists" by $K_{< \partial} \ne
\emptyset$ and let ${\frak K}^{\text{up}}_{\frak u}$ be its 
lifting up, as in \yCITE[0.31]{600} and 
we assume $K \ne \emptyset$ so ${\frak K} = 
{\frak K}^{\text{up}}_{\frak u}$ and we write 
$\tau_{\frak u} = \tau_{\frak K}$ and $\le_{\frak
K}$ for $\le_{{\frak K}^{\text{up}}}$ and $\le_{\frak u}$ for
$\le_{{\frak K}_{< \partial}}$)
\sn
\item "{$(D)_\ell$}"  $(a) \quad \text{ FR}_\ell$ is a class of triples of the
form $(M,N,\bold J)$, closed under \nl

\hskip25pt isomorphisms, let FR$^+_\ell = \text{ FR}^{{\frak
u},+}_\ell$ be the family of $(M,N,\bold J) \in \text{ FR}_\ell$ \nl

\hskip25pt such that $\bold J \ne \emptyset$; 
\sn
\item "{${{}}$}"  $(b) \quad$ if $(M,N,\bold J) \in \text{ FR}_\ell$
\ub{then} $M \le_{\frak u} N$ hence both are from $K_{\frak u}$
\nl

\hskip25pt  so of cardinality $< \partial$
\sn
\item "{${{}}$}"  $(c) \quad$ if $(M,N,\bold J) \in \text{ FR}_\ell$ 
\ub{then}\footnote{If we use the a.e.c. ${\frak K}'$ defined in
\scite{838-1a.19} we, in fact, weaken this demand to ``$\bold J \subseteq N$".
This is done, e.g. in 
the proof of \scite{838-e.1} that is in Definition \scite{838-e.1A}.}
$\bold J$ is a set of elements of $N \backslash M$
\sn
\item "{${{}}$}"  $(d) \quad$ if $M \in K_{\frak u}$ \ub{then}\footnote{We can
weaken this and in some natural example we have less, but we
circumvent this, via \scite{838-1a.19}, see \scite{838-e.1}(c); this applies
to (E)(b)$(\beta)$, too; see Example \scite{838-2b.8}} for some
$N,\bold J$ we have $(M,N,\bold J) \in \text{ FR}^+_\ell$  
\sn 
\item "{${{}}$}"  $(e) \quad$ if $M \le_{\frak u} N \in K_{\frak u}$
\ub{then} \footnote{not a great loss if we demand $M = N$; but then we
have to strengthen the amalgamation demand (clause (F)); this is 
really needed only for $\ell=2$}
$(M,N,\emptyset) \in \text{ FR}^{\frak u}_\ell$
\sn
\item "{$(E)_\ell$}"  $(a) \quad \le_\ell = \le^\ell_{\frak u}$ is 
a partial order on FR$_\ell$, closed under isomorphisms
\sn
\item "{${{}}$}"  $(b)(\alpha) \quad$ if $(M_1,N_1,\bold J_1) \le_\ell
(M_2,N_2,\bold J_2)$ then $M_1 \le_{\frak u} M_2,
N_1 \le_{\frak u} N_2$ and 
\nl

\hskip30pt $\bold J_1 \subseteq \bold J_2$
\sn
\item "{${{}}$}"  $\quad (\beta) \quad$ moreover $N_1 \cap M_2 = M_1$
(disjointness) 
\sn
\item "{${{}}$}"  $(c) \quad$ if 
$\langle (M_i,N_i,\bold J_i):i < \delta \rangle$ is $\le_\ell$-increasing
continuous \nl

\hskip30pt (i.e. in limit we take unions)
and $\delta < \partial$ \ub{then} the union \nl

\hskip30pt $\dbcu_{i < \delta} (M_i,N_i,\bold J_i) := (\dbcu_{j < \delta}
M_i,\dbcu_{i < \delta} N_i,\dbcu_{i < \delta} \bold J_i)$ belongs to
FR$_\ell$ and \nl

\hskip30pt  $j < \delta \Rightarrow (M_j,N_j,\bold J_j) \le_\ell (\dbcu_{i <
\delta} M_i,\dbcu_{i < \delta} N_i,\dbcu_{i < \delta} \bold J_i)$
\sn
\item "{${{}}$}"  $(d) \quad$ if $M_1 \le_{\frak u} M_2
\le_{\frak u} N_2$ and $M_1 \le_{\frak u} N_1
\le_{\frak u} N_2$ \ub{then}
\nl

\hskip30pt  $(M_1,N_1,\emptyset) \le_\ell (M_2,N_2,\emptyset)$ 
\sn
\item "{$(F)$}"  (amalgamation) if $(M_0,M_1,\bold I_1) \in
\text{ FR}_1,(M_0,M_2,\bold J_1) \in \text{ FR}_2$ 
and $M_1 \cap M_2 = M_0$ 
\ub{then} we can find \footnote{we can ask for $M'_3 \le_{\frak K}
M''_3$ and demand $(M_0,M_1,\bold J_1) \le_2 (M_2,M'_3,\bold J'_1),
(M_0,M_2,\bold I_1) \le (M_1,M''_3,\bold I'_2)$, no real harm here but
also no clear gain} $M_3,\bold I_2,\bold J_2$ such that 
$(M_0,M_1,\bold I_1) \le_1 
(M_2,M_3,\bold I_2)$ and $(M_0,M_2,\bold J_1) \le_2 (M_1,M_3,\bold
J_2)$ hence $M_\ell \le_{\frak u} M_3$ for $\ell=0,1,2$.
\endroster
\enddefinition
\bigskip

\proclaim{\stag{838-1a.5} Claim}  1) ${\frak K}_{\frak u}$ has disjoint
amalgamation.
\nl
2) If $\ell=1,2$ and $(M_0,M_1,\bold I_1) \in \text{\rm FR}_\ell$ and
$M_0 \le_{\frak u} M_2$ and $M_1 \cap M_2 = M_0$ \ub{then} we
can find a pair $(M_3,\bold I^*_2)$ such that: $(M_0,M_1,\bold I)
\le^\ell_{\frak u} (M_2,M_3,\bold I^*_2) \in \text{\rm FR}^{\frak u}_\ell$. 
\endproclaim
\bigskip

\demo{Proof}  1) Let $M_0 \le_{{\frak K}_{\frak u}} M_\ell$ for
$\ell=1,2$ and for simplicity $M_1 \cap M_2 = M_0$.  Let $\bold I_1 =
\emptyset$, so by condition (D)$_1$(e) of Definition \scite{838-1a.3} we
have $(M_0,M_1,\bold I_1) \in \text{\rm FR}_1$.  Now apply part (2),
(for $\ell=1$).
\nl
2) By symmetry without loss of generality  $\ell=1$.  Let $\bold J_1 := \emptyset$, so
by Condition (D)$_2$(e) of Definition \scite{838-1a.3} we have
$(M_0,M_2,\bold J_1) \in \text{\rm FR}^{\frak u}_2$.  So
$M_0,M_1,\bold I_1,M_2,\bold J_1$ satisfies the assumptions of
condition (F) of Definition \scite{838-1a.3} hence there are $M_3,\bold
I_2,\bold J_2$ as guaranteed there so in particular $(M_0,M_1,\bold
I_1) \le^2_{\frak u} (M_2,M_3,\bold I_2)$ so the pair $(M_3,\bold I_2)$ is as
required.  \hfill$\square_{\scite{838-1a.5}}$
\enddemo
\bigskip

\definition{\stag{838-1a.7} Definition}  1) We say that $\bold d$ is a
${\frak u}$-free $(\alpha,\beta)$-rectangle or
is ${\frak u}$-non-forking $(\alpha,\beta)$-rectangle 
(we may omit ${\frak u}$ when clear from the context) \ub{when}:
\mr
\item "{$(a)$}"  $\bold d = (\langle M_{i,j}:i \le
\alpha,j \le \beta \rangle,\langle \bold J_{i,j}:i < \alpha,j \le
\beta \rangle,\langle \bold I_{i,j}:i \le \alpha,j < \beta \rangle)$,
\nl
(we may add superscript $\bold d$, the ``$i < \alpha$", ``$j < \beta$"
are not misprints)
\sn
\item "{$(b)$}"  $\langle (M_{i,j},M_{i,j+1},\bold I_{i,j}):i \le \alpha
\rangle$ is $\le_1$-increasing continuous for each $j < \beta$
\sn
\item "{$(c)$}"  $\langle (M_{i,j},M_{i+1,j},\bold J_{i,j}):j \le
\beta \rangle$ is $\le_2$-increasing continuous for each $i < \alpha$
\sn
\item "{$(d)$}"  $M_{\alpha,\beta} = \cup\{M_{i,\beta}:i < \alpha\}$
if $\alpha,\beta$ are limit ordinals.
\ermn
2) For $\bold d^1$ a ${\frak u}$-free $(\alpha,\beta)$-rectangle
and $\alpha_1 \le \alpha,\beta_1 \le \beta$ let 
$\bold d^2 = \bold d^1 \restriction (\alpha_1,\beta_1)$ means:
\mr
\item "{$(a)$}"  $\bold d^2$ is a ${\frak u}$-free
$(\alpha_1,\beta_1)$-rectangle, see \scite{838-1a.9} below 
\sn
\item "{$(b)$}"  the natural equalities: $M^{{\bold d}^1}_{i,j} =
M^{{\bold d}^2}_{i,j},\bold J^{{\bold d}^1}_{i,j} = \bold J^{{\bold
d}^2}_{i,j},\bold I^{{\bold d}_1}_{i,j} = \bold I^{{\bold d}_2}_{i,j}$
when both sides are well defined.
\ermn
3) $\bold d^2 = \bold d^1 \restriction 
([\alpha_1,\alpha_2],[\beta_1,\beta_2])$ when $\alpha_1 \le \alpha_2
\le \alpha,\beta_1 \le \beta_2 \le \beta$ is defined similarly.
\nl
4) For $\bold d$ as above we may also write $\alpha_{\bold d},
\alpha(\bold d)$ for $\alpha$ and 
$\beta_{\bold d},\beta(\bold d)$ for $\beta$ 
and if $\bold I^{\bold d}_{i,j}$ is a singleton we
may write $\bold I^{\bold d}_{i,j} = \{a^{\bold d}_{i,j}\}$ and may
just write $(M^{\bold d}_{i,j},M^{\bold d}_{i,j+1},a^{\bold d}_{i,j})$ and if
$\bold J^{\bold d}_{i,j}$ is a singleton we may write 
$\bold J^{\bold d}_{i,j} = \{b^{\bold d}_{i,j}\}$ and may write
$(M^{\bold d}_{i,j},M^{\bold d}_{i+1,j},b^{\bold d}_{i,j})$.  
Similarly in Definition \scite{838-1a.11} below.
\nl
5) We may allow $\alpha = \partial$ and or $\beta \le \partial$, but
we shall say this.
\enddefinition
\bigskip

\demo{\stag{838-1a.9} Observation}  1) The restriction in Definition
\scite{838-1a.7}(2) always gives a ${\frak u}$-free 
$(\alpha_1,\beta_1)$-rectangle.
\nl
2) The restriction in Definition \scite{838-1a.7}(3) always gives a
${\frak u}$-free $(\alpha_2-\alpha_1,\beta_2-\beta_1)$-rectangle.
\nl
3) If $\bold d$ is a ${\frak u}$-free $(\alpha,\beta)$-rectangle, \ub{then}
\mr
\item "{$(e)$}"  $\langle M^{\bold d}_{i,j}:i \le \alpha\rangle$ is $\le_{\frak
u}$-increasing continuous for $j \le \beta$
\sn
\item "{$(f)$}"  $\langle M^{\bold d}_{i,j}:j \le \beta\rangle$ is $\le_{\frak
u}$-increasing continuous for each $i \le \alpha$.
\ermn
4) In Definition \scite{838-1a.11} below, clause (d), when $j < \beta$
or $j = \beta \wedge (\beta$ successor) follows from (b).  Similarly
for the pair of clauses (e),(c).
\nl
5) Assume that $\alpha_1 \le \alpha_2,\beta_1 \ge \beta_2$ and $\bold
   d_\ell$ is ${\frak u}$-free $(\alpha_\ell,\beta_\ell)$-rectangle
   for $\ell=1,2$ and $\bold d_1 \rest (\alpha_1,\beta_2) = \bold d_2
   \rest (\alpha_1,\beta_2)$ and $M^{\bold d_1}_{\alpha_1,\beta_1}
   \cap M^{\bold d_2}_{\alpha_2,\beta_2} = M^{\bold
   d_\ell}_{\alpha_1,\beta_2}$.  \ub{Then} we can find a ${\frak
   u}$-free $(\alpha_2,\beta_1)$-rectangle $\bold d$ such that $\bold
   d \rest (\alpha_\ell,\beta_\ell) = \bold d_\ell$.
\enddemo
\bigskip

\demo{Proof}  Immediate, e.g. in (5) we use clause (F) of Definition
\scite{838-1a.3} for each $\alpha \in [\alpha_1,\alpha_2),\beta \in
[\beta_2,\beta_1)$ in a suitable induciton.  \hfill$\square_{\scite{838-1a.9}}$
\enddemo
\bigskip

\definition{\stag{838-1a.11} Definition}  We say that $\bold d$ is a 
${\frak u}$-free $(\bar \alpha,\beta)$-triangle
or ${\frak u}$-non-forking $(\bar \alpha,\beta)$-triangle  
\ub{when} $\bar \alpha = \langle \alpha_i:i \le \beta \rangle$ 
is a non-decreasing\footnote{not unreasonable to demand $\bar\alpha$
 to be increasing continuous} sequence
of ordinals and (letting $\alpha := \alpha_\beta$):
\mr
\item "{$(a)$}"  $\bold d = (\langle M^{\bold d}_{i,j}:i \le \alpha_j,j
\le \beta \rangle,\langle \bold J_{i,j}:i < \alpha_j,j \le \beta
\rangle,\langle \bold I_{i,j}:i \le \alpha_j,j < \beta \rangle)$
\sn
\item "{$(b)$}"  $\langle (M^{\bold d}_{i,j},M^{\bold d}_{i,j+1},
\bold I_{i,j}):i \le \alpha_j\rangle$ is $\le_1$-increasing 
continuous for each $j<\beta$
\sn
\item "{$(c)$}"  $\langle (M^{\bold d}_{i,j},M^{\bold d}_{i+1,j},
\bold J_{i,j}):j \le \beta$
and $j$ is such that $i +1 \le \alpha_j\rangle$ 
is $\le_2$-increasing continuous for each $i < \alpha$ 
\sn
\item "{$(d)$}"  for each $j \le \beta$ the sequence $\langle
M_{i,j}:i \le \alpha_j\rangle$ is $\le_{\frak u}$-increasing
continuous
\sn
\item "{$(e)$}"  for each $i_* \in [\alpha_{j_*},\alpha),j_* \le \beta$ the
sequence $\langle M_{i_*,j}:j \in [j_*,\beta]\rangle$ 
is $\le_{\frak K}$-increasing continuous.
\endroster
\enddefinition
\bigskip

\definition{\stag{838-1a.13} Definition/Claim}  1) For nice construction
framework ${\frak u}_1$ let ${\frak u}_2 = \text{\rm dual}({\frak u}_1)$ be
the unique nice construction framework ${\frak u}_2$ such that:
$\partial_{{\frak u}_2} = \partial_{{\frak u}_1},{\frak K}_{{\frak
u}_2} = {\frak K}_{{\frak u}_1}$ (hence ${\frak K}^{\text{up}}_{{\frak
u}_2} = K^{\text{up}}_{{\frak u}_1}$, etc) and $(\text{FR}^{{\frak
u}_2}_\ell,\le^\ell_{{\frak u}_2}) = (\text{FR}^{{\frak u}_1}_{3 -
\ell},\le^{3 - \ell}_{{\frak u}_1})$ for $\ell=1,2$.
\nl
2) We call ${\frak u}_1$ self-{\rm dual} when {\rm dual}$({\frak u}_1) 
= {\frak u}_1$.
\nl
3) In part (1), if addition if $\bold d_1$ is ${\frak u}_1$-free rectangle
\ub{then} there is a unique $\bold d_2 = \text{\rm dual}(\bold d_1)$
which is a ${\frak u}_2$-free rectangle such that:

$$
\alpha_{\bold d_2} = \beta_{\bold d_1},\beta_{\bold d_2} =
\alpha_{\bold d_1},M^{\bold d_2}_{i,j} = M^{\bold d_1}_{j,i}
$$

$$
\bold I^{\bold d_2}_{i,j} = \bold J^{\bold d_1}_{j,i} \text{ and }
\bold J^{\bold d_2}_{i,j} = \bold I^{\bold d_1}_{i,j}.
$$
\enddefinition
\bigskip

\definition{\stag{838-1a.15} Definition}  1) We say $\tau$ is a weak
${\frak u}$-sub-vocabulary \ub{when}:
\mr
\item "{$(a)$}"  $\tau \subseteq \tau_{\frak u} = \tau_{{\frak
K}_{\frak u}}$ except that $=_\tau$ is, in $\tau_{\frak u}$, a two-place
predicate such that for every $M \in {\frak K}_{\frak u}$ hence even $M
\in {\frak K}^{\text{up}}_{\frak u}$, the relation 
$=^M_\tau$ is an equivalence relation on Dom$(=^M_\tau) = \{a:a =_\tau
b \vee b =_\tau a$ for some
$b \in M\}$ and is a congruence relation for all $R^M 
\restriction \text{ Dom}(=^M_\tau),F^M \restriction 
\text{ Dom}(=^M_\tau)$ for $R,F \in \tau$ and $F^M$ maps
\footnote{we may better ask less: for $F \in \tau$ a function symbol
letting $n = \text{ arity}_\tau(F)$, so $F^{M^{[\tau]}}$ is a
function with $n$-place from Dom$(=^M_\tau)/=^M_\tau$ to itself and
$F^M$ is $\{(a_0,a_1,\dotsc,a_n):(a_0/=^M =
F(a_1/=^M),\dotsc,a_n/=^M)\}$, i.e. the graph of $M^{[\tau]}$, so we
treat $F$ as an (arity$_\tau(F)+1$)-place predicate; neither real
change nor a real gain}
the set Dom$(=^M_\tau)$ into itself for any function symbol $F \in
\tau$.
\ermn
So
\nl
1A) For $M \in K^{\text{up}}_{\frak u}$, the model 
$M^{[\tau]}$ is defined naturally, e.g. with universe
\nl
Dom$(=^M_\tau)/=^M_\tau$ and $K^\tau,K^\tau_{\mu}$ are defined
accordingly.  Let $M_1 \cong_\tau M_2$ means $M^{[\tau]}_1 \cong
M^{[\tau]}_2$. 
\nl
1B) Let $\dot I_\tau(\lambda,K^{\text{up}}_{\frak u}) 
= \{M^{[\tau]}/\cong:M \in
K^{\frak u}_\lambda$ and $M^{[\tau]}$ has cardinality $\lambda\}$.
\nl
1C) We say that $\tau$ is a strong ${\frak u}$-sub-vocabulary \ub{when}
we have clause (a) from above 
and \footnote{note that it is important for us that
the model we shall construct will be of cardinality $\partial^+$; this clause
will ensure that the approximations will be of cardinality $\partial$
for $\alpha < \partial^+$ large enough and the final model (i.e. for
$\alpha = \partial^+$) will be of cardinality $\partial^+$.  
This is the reason for a
preference to $\le_1$, however there is no 
real harm in demanding clauses (b) + (c) for $\ell=2$, too.  But see
\scite{838-1a.17}, i.e. if $|\tau| \le \partial$ and we get $\mu >
2^\partial$ pairwise non-isomorphic models of cardinality $\le
\partial^+$, clearly only few (i.e. $\le 2^\partial$)  of them have
cardinal $< \partial^+$; so this problem is not serious to begin
with.}
\mr
\item "{$(b)$}"  if $(M,N,\bold I) \in \text{ FR}^+_1$ \ub{then}
for some $c \in \bold I \cap \text{Dom}(=^N_\tau)$ we have $N \models ``\neg 
(c =_\tau d)"$ for every $d \in M$
\sn
\item "{$(c)$}"  if  $(M_1,N_1,\bold I_1) \le_1 (M_2,N_2,\bold I_2)$
and $c \in \bold I_1$ is as in clause (b) for $(M_1,N_1,\bold I_1)$ 
\ub{then} $c \in \bold I_2$ is as in (b) for $(M_2,N_2,\bold I_2)$.
\ermn
2) We say that $N_1,N_2$ are $\tau$-isomorphic over $\langle M_i:i <
\alpha \rangle$ \ub{when}: $\ell \in \{1,2\} \wedge i < \alpha 
\Rightarrow M_i \le_{\frak u} N_\ell$
and there is a $\tau$-isomorphism $f$ of $N_1$ onto $N_2$ over
$\dbcu_{i < \alpha} M_i$ which means: $f$ is an isomorphism from
$N^{[\tau]}_1$ onto $N^{[\tau]}_2$ which is the identity on the
universe of $M^{[\tau]}_i$ for each $i < \alpha$. 
\nl
2A) In part (2), if 
$\alpha = 1$ we may write $M_0$ instead $\langle M_i:i
<1\rangle$ and we can replace $M_0$ by a set $\subseteq N_1 \cap N_2$.
If $\alpha=0$ we may omit ``over $\bar M$".
\nl
3) We say that $N_1,N_2$ are $\tau$-incompatible extensions of
$\langle M_i:i < \alpha \rangle$ \ub{when}:
\mr
\item "{$(a)$}"  $M_i \le_{\frak u} N_\ell$ for $i < \alpha,\ell =
1,2$
\sn
\item "{$(b)$}"  if $N_\ell \le_{\frak u} N'_\ell$ for $\ell=1,2$ then
$N'_1,N'_2$ are not $\tau$-isomorphic over $\langle M_i:i < \alpha
\rangle$.
\ermn
4) We say that $N^1_2,N^2_2$ are $\tau$-incompatible (disjoint)
amalgamations of $N_1,M_2$ over $M_1$ \ub{when} $(N_1 \cap M_2 = M_1$ and):
\mr
\item "{$(a)$}"  $M_1 \le_{\frak u} N_1 \le_{\frak u} N^\ell_2$ and 
$M_1 \le_{\frak u} M_2 \le_{\frak u} N^\ell_2$ for $\ell=1,2$
(equivalently $M_1 \le_{\frak u} N_1 \le_{\frak u} N^\ell_2,M_1
\le_{\frak u} M_2 \le_{\frak u} N^\ell_2$)
\sn
\item "{$(b)$}"  if $N^\ell_2 \le_{\frak u} N^{\ell,*}_2$ for
$\ell=1,2$ \ub{then} $(N^{1,*}_2)^{[\tau]},(N^{2,*}_2)^{[\tau]}$ are 
not $\tau$-isomorphic over $M_2 \cup N_1$, i.e. over $M^{[\tau]}_2
\cup N^{[\tau]}_1$.
\ermn
5) We say $\tau$ is a ${\frak K}$-sub-vocabulary or $K$-sub-vocabulary
\ub{when} clause (a) of part (1) holds replacing $K_{\frak u}$ by $K$;
similarly in parts (1A),(1B).
\enddefinition
\bigskip

\demo{\stag{838-1a.17} Observation}  Concerning \scite{838-1a.15}(1B) we may be
careless in checking the last condition, $= \lambda$, i.e. $\le\lambda$
usually suffice, because if $|\{M^{[\tau]}/\cong:M
\in K_\lambda,\|M^{[\tau]}\| < \lambda\}| < \mu$ then in
proving $\dot I_\tau(\lambda,K^{\text{up}}_{\frak u}) \ge \mu$ we may
omit it.
\enddemo
\bigskip

\remark{Remark}  1) But we give also remedies by 
FR$^+_\ell$, i.e., clause (c) of \scite{838-1a.15}(1).
\nl
2) We also give reminders in the phrasing of the coding
properties.
\nl
3) If $|\tau| < \lambda$ and $2^{< \lambda} < \lambda$ the demand in
\scite{838-1a.17} holds. 
\endremark
\bigskip

\demo{Proof}  Should be clear.  \hfill$\square_{\scite{838-1a.17}}$
\enddemo
\bigskip

\definition{\stag{838-1a.19} Definition}  1) For any  
$(< \partial)$-a.e.c. ${\frak K}$ let ${\frak K}'$ be the $(<
\partial$)-a.e.c. defined like ${\frak K}$ only adding the two-place
predicate $=_\tau$, demanding it to be a
congruence relation, i.e.
\mr
\item "{$(a)$}"  $\tau' = \tau({\frak K}') = \tau \cup \{=_\tau\}$
where $\tau = \tau({\frak K})$
\sn
\item "{$(b)$}"  $K' = \{M:M$ is a $\tau'$-model, $=^M_\tau$ is a
congruence relation and $M/=^M_\tau$ belongs to $\kappa$ and $\|M\| <
\partial\}$
\sn
\item "{$(c)$}" $M \le_{{\frak K}'} N$ \ub{iff} $M \subseteq N$ and
the following function is a $\le_{\frak K}$-embedding of $M/=^M_\tau$
into $N/=^N_\tau:f(a/=^M_\tau) = a/=^N_\tau$
\ermn
(see Definition \scite{838-1a.15}(1)).
\nl
1A) Similarly for ${\frak K}$ an a.e.c. or a $\lambda$-a.e.c.
\nl
2) This is a special case of Definition \scite{838-1a.15}.
\nl
3) We can interpret $M \in K$ as $M' \in K'$ just letting $M' \rest
\tau = M,=^{M'}_\tau$ is equality on $|M|$.
\nl
4) A model $M' \in {\frak K}$ is called $=_\tau$-full when $a \in M'
\Rightarrow \|M'\| = |\{b \in M':M' \models a =_\tau b\}|$.
\nl
5) A model $M' \in {\frak K}$ is called $(\lambda,=_\tau)$-full
\ub{when} $a \in M' \Rightarrow \lambda \le |\{b \in M':M' \models ``a
=_\tau b"\}|$.
\nl
6) A model $M'$ is called $=_\tau$-fuller \ub{when} it is
$=_\tau$-full and $\|M'\|$ is the cardinality of $M'/=^{M'}_\tau$.
\enddefinition
\bigskip

\proclaim{\stag{838-1a.20} Claim}  Assume ${\frak K}$
is\footnote{We can use ${\frak K}$ is an a.e.c. and 
have similar results.} a $(< \partial)$-a.e.c.  
and\footnote{now pedantically ${\frak K}$ may be both a $(<
\partial_1)$-a.e.c. and a $(< \partial_2)$-a.e.c., e.g. if $\partial_2
= \partial^+_1,K_{\partial_1} = \emptyset$, so really $\partial$
should be given} ${\frak K}'$ is from \scite{838-1a.19} and $\lambda < \partial$.
\nl
0) ${\frak K}'_\lambda$ is a $\lambda$-a.e.c.
\nl
1) If $M',N' \in K'_\lambda$ \ub{then} $(M'/=^{M'}_\tau) \in
K_{\le\lambda}$ and $(N'/=^{N'}_\tau) \in K_{\le\lambda}$ and if in
addition $M' \le_{{\frak K}'_\lambda} N'$ \ub{then} (up to identifying
$a/=^{M'}_\tau$ with $q/=^{N'}_\tau$) we have $(M'/=^{M'}_\tau) 
\le_{{\frak K}_{\le\lambda}} (N'/=^{N'}_\tau)$, i.e. 
$(M')^{[\tau]} \le_{\frak K} (N')^{[\tau]}$.
\nl
2) If $M' \subseteq N'$ are $\tau'_{\frak K}$-models of cardinality
$\lambda$ and $=^{M'}_\tau,=^{N'}_\tau$ are congruence relation on $M'
\restriction \tau_{\frak K},N' \restriction \tau_{\frak K}$
respectively, \ub{then}
\mr
\item "{$(a)$}"  $M' \in K'_{< \partial}$ \ub{iff} $(M'/=^{M'}_\tau) \in
K_{< \partial}$
\sn
\item "{$(b)$}"  $M' \le_{{\frak K}'_{< \partial}} N'$ \ub{iff}
$(M'/=^{M'}_\tau) \le_{{\frak K}_{< \partial}} (N'/=^{N'}_\tau)$;
pedantically\footnote{or define when $f$ is a $\le_{{\frak K}'_{<
\partial}}$-embedding of $M^1$ into $N'$} $(M'/=^{N'}_\tau)
\le_{{\frak K}_\lambda}(N'/=^{N'}_\tau)$ or make the natural identification
\sn
\item "{$(c)$}"  if $M',N'$ are $=_\tau$-fuller then in clauses
(a),(b) we can replace ``$\le \lambda"$ by ``$= \lambda$"
\sn
\item "{$(d)$}"  if $K_{< \lambda} = \emptyset$ then in clauses
(a),(b) we can replace ``$\le \lambda"$ by ``$= \lambda$".
\ermn
3) ${\frak K}'$ has disjoint amalgamation \ub{if} ${\frak K}$ has
   amalgamation.
\nl
4) $K_\lambda \subseteq K'_\lambda$ and $\le_{{\frak K}_\lambda} =
\le_{{\frak K}'_\lambda} \restriction K_\lambda$.
\nl
5) $({\frak K}^{\text{up}}_{< \mu})' = ({\frak K}'_{<
   \mu})^{\text{up}}$ for any $\mu$ so we call it $K'_{< \mu}$ and
   ${\frak K}'_\mu = ({\frak K}'_{< \mu^+})_\mu$.
\nl
6) For every $\mu$, 
\mr
\item "{$(a)$}"  $\dot I(\mu,K) = |\{M'/\cong:M' \in K'_\mu$ 
is $=_\tau$-fuller$\}|$
\sn
\item "{$(b)$}"  $K_\mu = \{M'/=^{M'}_\tau:M' \in K'_\mu$ is
$\mu$-fuller$\}$ under the natural identification
\sn
\item "{$(c)$}"  $K_{\le \mu} = \{M/ =^M_\tau:M \in K'_\mu\}$ under
the natural identification
\sn
\item "{$(d)$}"  if $M',N^1 \in K'_\mu$ are $\mu$-full \ub{then} $M'
\cong N' \Leftrightarrow (M'/=^{M'}_\tau) \cong (N'/=^{N'}_\tau)$.
\endroster
\endproclaim
\bigskip

\demo{Proof}  Straight.  \hfill$\square_{\scite{838-1a.20}}$
\enddemo
\bn
\margintag{1a.21}\ub{\stag{838-1a.21} Exercise}:  Assume ${\frak K},{\frak K}'$ are as in
\scite{838-1a.19}.
\nl
1) If $\lambda \ge |\tau_{\frak K}|$ and $2^\lambda < 2^{\lambda^+}$
\ub{then} $\dot I(\lambda^+,{\frak K}) + 2^\lambda = \dot
   I(\lambda^+,{\frak K}') + 2^\lambda$, so 
if $\dot I(\lambda^+,{\frak K}) > 2^\lambda$ 
or $\dot I(\lambda^+,{\frak K}') > 2^\lambda$
then they are equal.
\nl
2) If $\lambda > |\tau_{\frak K}|$ and $2^{< \lambda} < 2^\lambda$
   \ub{then} $\dot I(\lambda,{\frak K}) + 2^{< \lambda} = \dot
   I(\lambda^+,{\frak K}') + 2^{< \lambda}$ (and as above).
\bigskip

\remark{Remark}  Most of our examples satisfies montonicity, see
below.
\nl
But not so FR$_1,\le_1$ in \S4(C).
\endremark
\bn
\margintag{1a.24}\ub{\stag{838-1a.24} Exercise}:  Let ${\frak u}$ be a nice construction
framework, as usual.
\nl
1) [Definition] We say ${\frak u}$ satisfies (E)$_\ell$(e),
monotonicity, when:
\mr
\item "{$(E)_\ell(e)$}"  if $(M,N,\bold J) \in \text{ FR}^\ell_{\frak
u}$ and $N \le_{\frak u} N'$ 
then $(M,N,\bold J) \le^\ell_{\frak u} (M,N,\bold J') \in 
\text{ FR}^\ell_{\frak u}$.
\ermn
1A)  Let (E)(e) mean (E)$_1$(e) + (E)$_2$(e).
\nl
2) [Claim] Assume ${\frak u}$ has monotonicity.

Assume $\bold d$ is a ${\frak u}$-free
 $(\alpha_2,\beta_2)$-rectangle, $h_1,h_2$ is an increasing
continuous function from $\alpha_1 + 1,\beta_1 +1$ into $\alpha_2 + 1,
\beta_2 + 1$ respectively.  Then $\bold d'$
is a  ${\frak u}$-free rectangle where we define $\bold d'$ by:
\mr
\item "{$(a)$}"  $\alpha(\bold d') = \alpha_1,\beta(\bold d') =
\beta_1$
\sn
\item "{$(b)$}"  $M^{\bold d'}_{i,j} = M^{\bold d}_{h_1(i),h_2(1)}$ if
$i \le \alpha_1,j \le \beta_1$
\sn
\item "{$(c)$}"  $\bold J^{\bold d'}_{i,j} = \bold J^{\bold
d}_{h_1(i),h_2(j)}$ for $i \le \alpha_1,j < \beta_1$
\sn
\item "{$(d)$}"  $\bold J^{\bold d'}_{i,j} = \bold J^{\bold
d}_{h_1(i),h_2(j)}$ for $i < \alpha_1,j \le \beta_1$.
\ermn
3)  [Claim] Phrase and prove the parallel of part (2) 
for ${\frak u}$-free triangles.
\bigskip

\demo{\stag{838-1a.27} Observation}  Assume ${\frak u}$ is a nice
construction framework except that we omit clauses (D)$_\ell$(e) +
(E)$_\ell$(d) for $\ell=1,2$ but satisfying Claim \scite{838-1a.5}.
We can show that ${\frak u}'$ is a nice construction framework where we
define ${\frak u}'$ like ${\frak u}$ but, for $\ell=1,2$, we replace 
FR$_\ell,\le_\ell$ by FR$'_\ell,\le'_\ell$ defined as follows:
\mr
\item "{$(a)$}"   FR$'_\ell = \{(M_1,M_2,\bold J):(M_1,M_2,\bold J) 
\in \text{ FR}_\ell$ or $M_1 \le_{\frak u} 
M_2 \in K_{< \partial}$ and $\bold J = \emptyset\}$
\sn
\item "{$(b)$}"  $\le'_\ell = \{((M_1,N_1,\bold J'),(M_2,N_2,\bold J'')):
(M_1,N_1,\bold J') \le^\ell_{\frak u} (M_2,N_2,\bold J'')$ \ub{or}
$M_1 \le_{\frak u} N_1,\bold J' = \emptyset,
M_1 \le_{\frak u} M_2,N_1 \le_{\frak u} N_2,N_1 \cap M_2 = M_1$ and 
$(M_2,N_2,\bold J'') \in \text{ FR}_\ell$ \ub{or} $M_1 
\le_{\frak u} N_1 \le_{\frak u} N_2,M_1
\le_{\frak u} M_2 \le_{\frak u} N_2,N_1 \cap M_2 = M_1$ 
and $\bold J' = \emptyset = \bold J''\}$.
\endroster
\enddemo
\bigskip

\demo{Proof}  Clauses (A),(B),(C) does not change, most subclauses of
(D)$_\ell$(a),(b),(d),(E)$_\ell$(a),(b) hold by the parallel 
for ${\frak u}$ and the choice of FR$'_\ell,\le'_\ell$; clauses
(D)$_\ell$(e) and (E)$_\ell$(d) holds by the choice of
(FR$'_\ell,\le'_\ell$); and clause (F) holds by clause (F) for 
${\frak u}$ and Claim \scite{838-1a.5}.  Lastly
\mn
\ub{Condition (E)$_\ell$(c)}:

So assume $\langle (M_i,N_i,\bold J_i):i < \delta\rangle$ be
increasing continuous, where $\delta$ is a limit ordinal; and let
$(M_\delta,N_\delta,\bold J_\delta) = (\cup\{M_i:i < \delta\},
\cup\{N_i:i < \delta\}, \cup \{\bold J_i:i < \delta\})$.

First, assume $i < \delta \Rightarrow \bold J_i = \emptyset$ hence
$\bold J_\delta = \emptyset$ and the desired conclusion holds
trivially (by the properties of a.e.c. and our definition of ${\frak u}'$).

Second, assume $i < \delta \nRightarrow \bold J_i \ne \emptyset$ hence
 $j := \text{ min}\{i:\bold J_i \ne \emptyset\}$
 is well defined and let $\delta' = \delta - j$, it is a limit ordinal.  Now
use the ``${\frak u}$ satisfies the Condition (E)$_\ell$(c)" for the
sequence $\langle (M_{j+i},N_{j+i},\bold J_{j+i}:i < \delta'\rangle$
and $\le^\ell_{\frak u}$ being transitive.  
\nl
${{}}$  \hfill$\square_{\scite{838-1a.27}}$
\enddemo
\bn
\centerline {$* \qquad * \qquad *$}
\bn
Now we define the approximations of size $\partial$; note that the
notation $\le_{\text{qt}}$ and the others below hint that they are
quasi orders, this will be justified later in \scite{838-1a.37}(2),
\ub{but} not concerning $\le^{\text{at}}_{\frak u}$.  On the
existence of canonical limits see \scite{838-1a.37}(4).
\definition{\stag{838-1a.29} Definition}  1) We let $K^{\text{qt}}_\partial =
K^{\text{qt}}_{\frak u}$ be the class of triples $(\bar M,\bar{\bold
J},\bold f)$ such that
\mr
\item "{$(a)$}"  $\bar M = \langle M_\alpha:\alpha < \partial \rangle$
is $\le_{\frak u}$-increasing continuous, so $M_\alpha \in K_{\frak u} \,
(= K^{\frak u}_{< \partial})$
\sn
\item "{$(b)$}"  $\bar{\bold J} = \langle \bold J_\alpha:\alpha < \partial
\rangle$
\sn
\item "{$(c)$}"  $\bold f \in {}^\partial \partial$
\sn
\item "{$(d)$}"  $(M_\alpha,M_{\alpha +1},\bold J_\alpha) \in \text{
FR}_2$ for $\alpha < \partial$.
\ermn
1A) We call $(\bar M,\bar{\bold J},\bold f) \in K^{\text{qt}}_{\frak
u}$ non-trivial \ub{if} for stationarily many $\delta < \partial$ for some
$i < \bold f(\delta)$ we have $\bold J_{\delta +i} \ne
\emptyset$ that is $(M_{\delta +i},M_{\delta +i+1},\bold J_{\delta +i}) 
\in \text{ FR}^+_2$.
\nl
1B) If ${\Cal D}$ is a normal filter on $\partial$ let
$K^{\text{qt}}_{\Cal D} = K^{\text{qt}}_{{\frak u},{\Cal D}}$ be the
class of triples $(\bar M,\bar J,\bold f) \in K^{\text{qt}}_{\frak u}$
such that
\mr
\item "{$(e)$}"  $\{\delta < \partial:\bold f(\delta)=0\} \in {\Cal D}$.
\ermn
1C) When we have $(\bar M^x,\bar{\bold J}^x,\bold f^x)$ then
$M^x_\alpha,\bold J^x_\alpha$ for $\alpha < \partial$ 
has the obvious meaning and
$M^x_\partial$  or just $M^x$ is $\cup\{M^x_\alpha:\alpha < \partial\}$
\nl
2) We define the two-place relation $\le_{\text{qt}} = 
\le^{\text{qt}}_{\frak u}$ on 
$K^{\text{qt}}_{\frak u}$ as follows: $(\bar M^1,\bar{\bold J}^1,\bold
f^1) \le^{\text{qt}}_{\frak u} 
(\bar M^2,\bar{\bold J}^2,\bold f^2)$ \ub{if} they are equal (and $\in
K^{\text{qt}}_{\frak u}$) or 
for some club $E$ of $\partial$ (a witness) we have:
\mr
\item "{$(a)$}"  $(\bar M^k,\bar{\bold J}^k,\bold f^k) \in
K^{\text{qt}}_{\frak u}$ for $k=1,2$
\sn
\item "{$(b)$}"  $\delta \in E \Rightarrow \bold f^1(\delta) \le \bold
f^2(\delta)$
\sn
\item "{$(c)$}"  $\delta \in E \and i \le \bold f^1(\delta)
\Rightarrow M^1_{\delta +i} \le_{\frak u} M^2_{\delta +i}$
\sn
\item "{$(d)$}"  $\delta \in E \and i < \bold f^1(\delta) \Rightarrow
(M^1_{\delta +i},M^1_{\delta +i+1},\bold J^1_{\delta +i}) \le_2
(M^2_{\delta +i},M^2_{\delta +i+1},\bold J^2_{\delta +i})$
\sn
\item "{$(e)$}"  $\delta \in E \and i \le \bold f^1(\delta) \Rightarrow
M^2_{\delta+i} \cap \dbcu_{\alpha < \partial} M^1_\alpha = M^1_{\delta
+i}$, disjointness.
\ermn
3) We define the two place relation $\le_{\text{at}} = 
\le^{\text{at}}_{\frak u}$ on $K^{\text{qt}}_{\frak u}$: 
$(\bar M^1,\bar{\bold J}^1,\bold f^1) \le^{\text{at}}_{\frak u} 
(\bar M^2,\bar{\bold J}^2,\bold f^2)$ 
\ub{if} for some club $E$ of $\partial$ and $\bar{\bold I}$ (the
witnesses) we have (a)-(e) as in part (2) and
\mr
\item "{$(f)$}"  $\bar{\bold I} = \langle \bold I_\alpha:\alpha < \partial
\rangle$ and $\langle (M^1_\alpha,M^2_\alpha,\bold I_\alpha):\alpha
\in \cup\{[\delta,\delta + \bold f^1(\delta)]:\delta \in E\}\rangle$ is
$\le_1$-increasing continuous, so we may use $\langle \bold
I_\alpha:\alpha \in \cup\{[\delta,\delta + \bold f^1(\delta)]:\delta \in
E\rangle$ only.
\ermn
3A) We say $(\bar M^1,\bar{\bold J}^1,\bold f^1),(\bar M^2,
\bar{\bold J}^2,\bold f^2)$ are equivalent when for a club of $\delta <
\partial$ we have $\bold f^1(\delta) = \bold f^2(\delta)$ and $i \le
\bold f^1(\delta) \Rightarrow M^1_{\delta +i} = M^2_{\delta +i}$ and $i <
\bold f^1(\delta) \Rightarrow \bold J^1_{\delta +i} = \bold
J^2_{\delta +i}$. \nl
3B) Let $(\bar M^1,\bar{\bold J}^1,\bold f^1) <^{\text{at}}_{\frak u}(\bar
M^2,\bar{\bold J}^2,\bold f^2)$ mean that in part (3) in addition
\mr
\item "{$(g)$}"   for some $\alpha \in
\cup\{[\delta,\delta + \bold f(\delta)]:\delta \in E\}$ the triple
$(M^1_\alpha,M^2_\alpha,\bold I_\alpha)$ belongs to FR$^+_1$.
\ermn
4) We say $(\bar M^\delta,\bar{\bold J}^\delta,\bold f^\delta)$ is a
canonical limit of $\langle (\bar M^\alpha,\bar{\bold J}^\alpha,
\bold f^\alpha):\alpha < \delta \rangle$ \ub{when}:
\mr
\item "{$(a)$}"  $\delta < \partial^+$
\sn
\item "{$(b)$}"  $\alpha < \beta < \delta \Rightarrow (\bar
M^\alpha,\bar{\bold J}^\alpha,\bold f^\alpha) \le^{\text{qt}}_{\frak u} (\bar
M^\beta,\bar{\bold J}^\beta,\bold f^\beta)$
\sn
\item "{$(c)$}"  for some increasing continuous sequence
$\langle \alpha_\varepsilon:\varepsilon < \text{ cf}(\delta)\rangle$
of ordinals with limit $\delta$ we have:
\ermn
\ub{Case 1}:  cf$(\delta) < \partial$.

For some club $E$ of $\partial$ we have:
\mr
\item "{$(\alpha)$}"  $\zeta \in E \Rightarrow \bold f^\delta(\zeta) = 
\sup\{\bold f^{\alpha_\varepsilon}(\zeta):\varepsilon < \text{
cf}(\delta)\}$; 
\sn
\item "{$(\beta)$}"  $\zeta \in E 
\and \xi < \text{ cf}(\delta) \and i \le \bold f^{\alpha_\xi}(\zeta)$ implies 
that $M^\delta_{\zeta +i} = \cup\{M^{\alpha_\varepsilon}_{\zeta
+i}:\varepsilon < \text{ cf}(\delta)$ satisfies $\varepsilon \ge \xi\}$; 
\sn
\item "{$(\gamma)$}"  $\zeta \in E \and \xi < \text{ cf}(\delta) \and
i < \bold f^\xi(\zeta)$ implies that 
$\bold J^\delta_{\zeta +i} = \cup\{\bold
J^{\alpha_\varepsilon}_{\zeta +i}:\varepsilon < \text{ cf}(\delta)$ satisfies
$\varepsilon \ge \xi\}$ 
\sn
\item "{$(\delta)$}"   if $\zeta \in E$ and $j = \bold f^\delta(\zeta) > 
\bold f^{\alpha_\varepsilon}(\zeta)$ for every
$\varepsilon < \text{ cf}(\delta)$ then $M^\delta_{\zeta +j} =
\cup\{M^\delta_{\zeta +i}:i<j\}$.
\endroster
\bn
\ub{Case 2}:  cf$(\delta) = \partial$.

Similarly, using diagonal unions. 
\nl
4A) We say $\langle(\bar M^\alpha,\bar{\bold J}^\alpha,\bold
f^\alpha):\alpha < \alpha(*)\rangle$ is 
$\le^{\text{qt}}_{\frak u}$-increasing 
continuous \ub{when} it is $\le^{\text{qt}}_{\frak u}$-increasing and
for every limit ordinal $\delta < \alpha(*)$, the triple $(\bar
M^\delta,\bar{\bold J}^\delta,\bold f^\delta)$ is a canonical limit of
$\langle (\bar M^\alpha,\bar{\bold J}^\alpha,\bold f^\alpha):\alpha <
\delta\rangle$.
\nl
5) We define the relation $\le_{\text{qs}} = 
\le^{\text{qs}}_{\frak u}$ on $K^{\text{qt}}_{\frak u}$
by: $(\bar M',\bar{\bold J}',\bold f') \le_{\text{qs}} (\bar M'',\bar{\bold
J}'',\bold f'')$ \ub{if} there is a $\le^{\text{at}}_{\frak u}$-tower
$\langle (\bar M^\alpha,\bar{\bold J}^\alpha,
\bold f^\alpha):\alpha \le \alpha(*)\rangle$ witnessing it, meaning
that is is a sequence such that:
\mr
\item "{$(a)$}"  the sequence is $\le^{\text{qt}}_{\frak u}$-increasing of
length $\alpha(*) + 1 < \partial^+$
\sn
\item "{$(b)$}"  $(\bar M^0,\bar{\bold J}^0,\bold f^0) = (\bar
M',\bar{\bold J}',\bold f')$
\sn
\item "{$(c)$}"  $(\bar M^{\alpha(*)},\bar{\bold J}^{\alpha(*)},
\bold f^{\alpha(*)}) = (\bar M'',\bar{\bold J}'',\bold f'')$
\sn
\item "{$(d)$}"  $(\bar M^\alpha,\bar{\bold J}^\alpha,\bold f^\alpha)
\le^{\text{at}}_{\frak u} (\bar M^{\alpha +1},\bar{\bold J}^{\alpha
+1},\bold f^{\alpha +1})$ for $\alpha < \alpha(*)$
\sn
\item "{$(e)$}"  if $\delta \le \alpha(*)$ is a limit ordinal then
$(\bar M^\delta,\bar{\bold J}^\delta,\bold f^\delta)$ is a canonical
limit of $\langle (\bar M^\alpha,\bar{\bold J}^\alpha,\bold
f^\alpha):\alpha < \delta \rangle$.
\ermn
5A) Let $<_{\text{qs}} = <^{\text{qs}}_{\frak u}$ be defined 
similarly but for at least one
$\alpha < \alpha(*)$ we have $(\bar M^\alpha,\bar{\bold J}^\alpha,
\bold f^\alpha) <^{\text{at}}_{\frak u} (\bar
M^{\alpha +1},\bar{\bold J}^{\alpha +1},\bold f^{\alpha +1})$.
\nl
5B) Let $\le_{\text{qr}} = \le^{\text{qr}}_{\frak u}$ be defined as in
part (5) but in clause (d) we use $<^{\text{at}}_{\frak u}$.
Similarly for $<_{\text{qr}} = <^{\text{qr}}_{\frak u}$, i.e. when
$\alpha(*) > 0$.
\nl
6) We say that
$\langle (\bar M^\alpha,\bar{\bold J}^\alpha,\bold f^\alpha):\alpha
< \alpha(*)\rangle$ is $\le^{\text{qs}}_{\frak u}$-increasing
continuous \ub{when} it is $\le^{\text{qs}}_{\frak u}$-increasing and
clause (e) of part (5) holds.  Similarly for $\le^{\text{qr}}_{\frak u}$.
\enddefinition
\bn
Some obvious properties are (see more in Observation \scite{838-1a.37}).
\demo{\stag{838-1a.31} Observation}  1) $(\bar M^1,\bar{\bold J}^1,\bold f^1)
\le^{\text{at}}_{\frak u} (\bar M^2,\bar{\bold J}^2,\bold f^2)$
\ub{iff} for some club $E$ of $\partial$ and sequence $\bar{\bold I} =
\langle \bold I_\alpha:\alpha \in \cup\{[\delta,\delta + \bold
f^1(\delta)]:\alpha \in E\}\rangle$ we have clause (a),(b),(c) of
Definition \scite{838-1a.29}(2) and
\mr
\item "{$(f)'$}"  $(\langle M^1_{\delta +i}:i \le \bold f^1(\delta)
\rangle,\langle M^2_{\delta +i}:i \le \bold f^1(\delta)
\rangle,\langle \bold J^1_{\delta +i}:i < \bold f^1(\delta)\rangle,
\langle\bold J^2_{\delta +i}:i < \bold f^1(\delta)\rangle,
\langle \bold I_{\delta +i}:i \le \bold f^1(\delta)\rangle)$ is a 
${\frak u}$-free $(\bold f^1(\delta),1)$-rectangle
\sn
\item "{$(f)'_1$}"  if $\delta_1 < \delta_2$ are from $E$ then
$(M^1_{\delta_1 + \bold f^1(\delta_1)},M^2_{\delta_1+ \bold f^1(\delta_1)},
\bold I_{\delta_1 + \bold f^1(\delta_1)}) \le_2
(M^1_{\delta_2},M^2_{\delta_2},\bold I_{\delta_2})$.
\ermn
2) The relation $\le^{\text{qt}}_{\frak u},\le^{\text{at}}_{\frak u},
\le^{\text{qs}}_{\frak u},\le^{\text{qr}}_{\frak u}$ and
$\le^{\text{at}}_{\frak u}$ are preserved by equivalence, see
Definition \scite{838-1a.29}(3A) (and equivalence is an equivalence
relation) and so are $<_{\text{at}},<_{\text{qt}},<_{\text{qs}},<_{\text{qr}}$.
\enddemo
\bigskip

\demo{Proof}  Straightforward.  \hfill$\square_{\scite{838-1a.31}}$
\enddemo
\bigskip

\remark{\stag{838-1a.33} Remark}  1) In some of our applications 
it is natural to redefine the partial order 
$\le^{\text{qs}}_{\frak u}$ we use on
$K^{\text{qt}}_{\frak u}$ as the closure of
a more demanding relation.
\nl
2) If we demand FR$_1 = \text{ FR}^+_1$ hence we omit 
clause $(D)_1(e),(E)_1(d)$ of Definition \scite{838-1a.3}, really
$<^{\text{at}}_{\frak u}$ is the same as $\le^{\text{at}}_{\frak
u}$.  In Definition \scite{838-1a.29}(3) we can choose $\bold I_\alpha =
\emptyset$, then we get $\le^{\text{qt}}_{\frak u}$.
But even so we would like to be able to  say ``repeat \S1
with the following modifications".  If in Definition \scite{838-1a.29}(5)
clause (d) we use $<^{\text{at}}_{\frak u}$, i.e. use
$<^{\text{qr}}_{\frak u}$,  the difference below is small. 
\nl
3) Note that below $K^{{\frak u},*}_\partial \subseteq
K^{\text{up}}_{{\frak u},\partial}$ and $K^{{\frak u},*}_{\partial^+} \subseteq
K^{\text{up}}_{{\frak u},\partial^+}$, see the definition below.  
\nl
4) Should we use $\le^{\text{qs}}_{\frak u}$ or $\le^{\text{qr}}_{\frak
   u}$ (see Definition \scite{838-1a.29}(5A),(5B))?  So far it does not matter.
\endremark
\bigskip

\definition{\stag{838-1a.35} Definition}  1) $K^{{\frak u},*}_\partial = \{M:M =
\cup\{M_\alpha:\alpha < \partial\}$ for some non-trivial $(\bar M,\bar{\bold
J},\bold f) \in K^{\text{qt}}_{\frak u}\}$, recalling Definition
\scite{838-1a.29}(1A). 
\nl
2) $K^{{\frak u},*}_{\partial^+} = \{\bigcup\{M^\gamma:\gamma <
\partial^+\}:\langle(\bar M^\gamma,\bar{\bold J}^\gamma,\bold
f^\gamma):\gamma < \partial^+ \rangle$ is 
$\le^{\text{qs}}_{\frak u}$-increasing
continuous and for unboundedly many $\gamma < \partial^+$ we have
$(\bar M^\alpha,\bar{\bold J}^\alpha,\bold f^\alpha) 
<^{\text{qs}}_{\frak u} (\bar M^{\alpha +1},
\bar{\bold J}^{\alpha +1},\bold f^{\alpha +1})$ and as usual $M^\gamma =
\cup\{M^\gamma_\alpha:\alpha < \partial\}\}$. 
\enddefinition
\bigskip

\demo{\stag{838-1a.37} Observation}  1) $K^{\text{qt}}_{\frak u} \ne
\emptyset$; moreover it has non-trivial members.
\nl
2) The two-place relations $\le^{\text{qt}}_{\frak u}$ and 
$\le^{\text{qs}}_{\frak u},\le^{\text{qr}}_{\frak u}$ are quasi
orders and so are $<_{\text{qt}},<_{\text{qs}},<_{\text{qr}}$
\ub{but} not necessarily $\le_{\text{at}},<_{\text{at}}$. 
\nl
3)  Assume $(\bar M^1,\bar{\bold J}^1,\bold f^1) \in
K^{\text{qt}}_{\frak u}$ and $\alpha < \partial,
(M^1_\alpha,M^2_\alpha,\bold I^*) \in \text{ FR}_1$ and
$\bold f^1 \le_{{\Cal D}_\partial} \bold f^2 \in {}^\partial \partial$
and $M^2_\alpha \cap M^1 = M^1_\alpha$.  
\ub{Then} we can find $\bar M^2,\bar{\bold J}^2,E$ and $\bar{\bold I} =
\langle \bold I_\alpha:\alpha \in \cup\{[\delta,\bold
f^1(\delta)]:\delta \in E\}\rangle$ such that
\mr
\item "{$(a)$}"  $(\bar M^1,\bold J^1,\bold f^1) 
\le^{\text{at}}_{\frak u} (\bar
M^2,\bar{\bold J}^2,\bold f^2)$ as witnessed by $E,\bar{\bold I}$
\sn
\item "{$(b)$}"  if $\beta \in \cup\{[\delta,\bold
f^2(\delta)]:\delta \in E\}$ then 
$(M^1_\alpha,M^2_\alpha,\bold I^*) \le_1 (M^1_\beta,M^2_\beta,\bold I_\beta)$
\sn
\item "{$(c)$}"  if $(M^1_\alpha,M^2_\alpha,\bold I^*) \in 
\text{ FR}^+_1$ then $(\bar M^1,\bar{\bold J}^1,\bold f^1) 
<^{\text{qs}}_{\frak u} (\bar M^2,\bar{\bold J}^2,\bold f^2)$ moreover
$(\bar M^1,\bar{\bold J}^1,\bold f^1) 
<^{\text{at}}_{\frak u} (\bar M^2,\bar{\bold J}^2,\bold f^2)$.
\ermn
4) If $\langle(\bar M^\alpha,\bar{\bold J}^\alpha,\bold
f^\alpha):\alpha < \delta \rangle$ is $\le^{\text{qt}}_{\frak u}$-increasing
continuous (i.e. we use canonical limits) and $\delta$ is a limit
ordinal $< \partial^+$,
\ub{then} it has a canonical limit $(\bar M^\delta,\bar{\bold
J}^\delta,\bold f^\delta)$ which is unique up to equivalence (see
\scite{838-1a.29}(3A)).  Similarly for $\le^{\text{qs}}_{\frak u}$ and
$\le^{\text{qr}}_{\frak u}$.
\nl
5) $K^{{\frak u},*}_\partial,K^{{\frak u},*}_{\partial^+}$ are non-empty and
 included in $K^{\frak u}$, in fact in $K^{\frak u}_\partial,K^{\frak
u}_{\partial^+}$, respectively.  Also if $\tau$ is a strong ${\frak
 u}$-sub-vocabulary and $M \in K^{{\frak u},*}_\partial$ or $M \in
 K^{{\frak u},*}_{\partial^+}$ then $M^{[\tau]}$ has cardinality
 $\partial$ or $\partial^+$ respectively.  If $\tau$ is a weak ${\frak
 u}$-subvocabulary we get only $\le \partial,\le \partial^+$ respectively.
\enddemo
\bigskip

\demo{Proof}  1) We choose $M_i \in K_{\frak u} = 
{\frak K}_{< \partial},\le_{\frak u}$-increasing continuous with 
$i$ as follows.  For $i=0$ use ${\frak
K}_{< \partial} \ne \emptyset$ by clause (C) of Definition \scite{838-1a.3}.  
For $i$ limit note that $\langle M_j:j < i \rangle$ is
$\le_{\frak K}$-increasing continuous and $j < i \Rightarrow
M_j \in {\frak K}_{\frak u}$ 
hence $j<i \Rightarrow \|M_j\| < \partial$ but $i <
\partial$ and $\partial$ is regular (by Definition \scite{838-1a.3}, clause (B))
so $M_i := \cup\{M_j:j<i\}$ has cardinality $< \partial$ hence (by
clause (C) of Definition \scite{838-1a.3}) $M_i \in {\frak K}_{\frak u}$
and $j<i \Rightarrow M_j \le_{\frak u} M_i$.  For $i=j+1$ by clause
$(D)_2(d)$ of Definition \scite{838-1a.3} there are $M_i,\bold J_j$
such that $(M_j,M_i,\bold J_j) \in \text{\rm FR}^+_2$.  Choose $\bold f \in
{}^\partial \partial$, e.g., $\bold f(\alpha)=1$.
\nl
Now $(\langle M_i:i < \partial \rangle,\langle \bold J_i:i < \partial
\rangle,\bold f) \in K^{\text{qt}}_{\frak u}$ is as required; moreover
is non-trivial, see Definition \scite{838-1a.29}(1A). 
\nl
2) We first deal with $\le^{\text{qt}}_{\frak u}$.

Trivially $(\bar M,\bar{\bold J},\bold f) \le^{\text{qt}}_{\frak u} (\bar
M,\bar{\bold J},\bold f)$ for $(M,\bar{\bold J},\bold f) \in
K^{\text{qt}}_{\frak u}$. \nl
[Why?  It is witnessed by $E = \partial$ (as 
$(M_i,M_i,\emptyset) \in \text{ FR}_1$ by clause $(D)_1(e)$ of
Definition \scite{838-1a.3} and $i < j < \partial \Rightarrow
(M_i,M_i,\emptyset) \le_1 (M_j,M_j,\emptyset)$ by clause $(E)_1(d)$ of
Definition \scite{838-1a.3}.] 
\nl
So assume $(\bar M^\ell,\bar{\bold J}^\ell,\bold f) \le^{\text{qt}}_{\frak u}
(\bar M^{\ell +1},\bar{\bold J}^{\ell +1},\bold f^\ell)$ and let it
be witnessed by $E_\ell$ for $\ell = 1,2$.  Let $E =
E_1 \cap E_2$, it is a club of $\partial$.   For every $\delta \in E$
by clause (b) of Definition \scite{838-1a.29}(2) we 
have $\bold f^1(\delta) \le \bold f^2(\delta) \le \bold
f^3(\delta)$ hence $\bold f^1(\delta) \le \bold f^3(\delta)$
and for $i \le \bold f^1(\delta)$ by clause (c) of Definition
\scite{838-1a.29}(2), clearly $M^1_{\bold f^1(\delta)+i} \le_{\frak u} M^2_{\bold
f^1(\delta)+i} \le_{\frak u} M^3_{\bold f^1(\delta)+i}$ so as 
${\frak K}_{\frak u}$ is a $(< \partial)$-a.c.e. we have $M^1_{\delta +i}
\le_{\frak u} M^3_{\delta +i}$.  Similarly as 
$\le_{\frak u}^2$ is a quasi order by clause (E)$_2$(a) of Definition
\scite{838-1a.3} we have $\delta \in E \and 
i < \bold f^1(\delta) \Rightarrow
(M^1_{\delta +i},M^1_{\delta +i+1},\bold J^1_{\delta +i}) 
\le_{\frak u}^2 (M^3_{\delta +i},M^3_{\delta +i+1},
\bold J^3_{\delta +i})$ so clause (d) of Definition \scite{838-1a.29}(2) holds.  

Also if $\delta \in E$ and 
$i \le \bold f^1(\delta)$ then $M^2_{\delta +i} \cap
(\cup\{M^1_\gamma:\gamma < \partial\}) = M^1_{\delta +i}$ and
$M^3_{\delta  +i} \cap (\cup\{M^2_\gamma:\gamma < \partial\}) =
M^2_{\delta +i}$ hence $M^3_{\delta +i} \cap(\cup\{M^1_\gamma:\gamma <
\partial\}) = M^1_{\delta +i}$, i.e. clause (e) there holds and clause
(a) is trivial.

Together really $(\bar M^1,\bar{\bold J}^1,\bold
f^1) \le^{\text{qt}}_{\frak u} (\bar M^3,\bar{\bold J}^3,\bold f^3)$.
\nl
So $\le^{\text{qt}}_{\frak u}$ is actually a quasi order.  As for 
$\le^{\text{qs}}_{\frak u}$ and $\le^{\text{qr}}_{\frak u}$, this
follows by the result on 
$\le^{\text{qt}}_{\frak u}$ and the definitions.  Similarly for
$<^{\text{qt}}_{\frak u},<^{\text{qs}}_{\frak u},<^{\text{qr}}_{\frak u}$.
\nl
3) By induction on $\beta \in [\alpha,\partial)$ we choose
$(g_\beta,M^2_\beta,\bold I_\beta)$ and $\bold J_\beta$ (but $\bold
J_\beta$ is chosen in the $(\beta +1)$-th step) such that
\mr
\item "{$(a)$}"  $M^2_\beta \in {\frak K}_{\frak u}$ is 
$\le_{\frak u}$-increasing continuous
\sn
\item "{$(b)$}"  $g_\beta$ is a $\le_{\frak u}$-embedding of
$M^1_\beta$ into $M^2_\beta$, increasing and continuous with $\beta$
\sn
\item "{$(c)$}"  $(g_\beta(M^1_\beta),M^2_\beta,\bold I_\beta) \in
\text{ FR}^{\frak u}_1$ is $\le^1_{\frak u}$-increasing and
continuous 
\sn
\item "{$(d)$}"  if $\beta = \gamma +1$ then
$(g_\gamma(M^1_\gamma),g_\beta(M^1_\beta),g_\beta(\bold J^1_\beta)) 
\le_{\frak u}^2 (M^2_\gamma,M^2_\beta,\bold J^2_\beta)$.
\ermn
For $\beta = \alpha$ let $g_\beta = \text{ id}_{M^1_\alpha},M^2_\beta$ 
as given, $\bold I_\beta = \bold I^*$, so by
the assumptions all is O.K.

For $\beta$ limit use clause (E)$_1$(c) of Definition \scite{838-1a.3}.

For $\beta = \gamma +1$ use clause $(F)$ of Definition \scite{838-1a.3}.
Having carried the induction, by renaming without loss of generality  $g_\beta = \text{
id}_{M^1_\beta}$ for $\beta < \partial$.  So clearly we are done. 
\nl
4),5)  Easy, too.    \hfill$\square_{\scite{838-1a.37}}$
\enddemo
\bigskip

\definition{\stag{838-1a.39} Definition}  1) $K^{\text{qt}}_\alpha =
K^{\text{qt}}_{{\frak u},\alpha}$ is defined as in Definition
\scite{838-1a.29}(1) above but $\alpha =
\text{ Dom}(\bold f) = \ell g(\bar{\bold J}) = \ell g(\bar M)-1$, where
$\alpha \le \partial$.
\nl
2) $K^{\text{qt}}_{< \alpha} = K^{\text{qt}}_{{\frak u},< \alpha},
K^{\text{qt}}_{\le \alpha} = K^{\text{qt}}_{{\frak u},\le \alpha}$ are defined
similarly.  
\enddefinition
\bn
\centerline {$* \qquad * \qquad *$}
\bn
\margintag{1a.41}\ub{\stag{838-1a.41} Discussion}:  1) Central here in \xCITE{838} are 
``$\tau$-coding properties" meaning that they will
help us in building $M \in K^{\frak u}_{\partial^+}$, moreover in
$K^{{\frak u},*}_{\partial^+}$ (or $K^{{\frak u},{\frak
h}}_{\partial^+}$, see below) such that 
we can code some subset of $\partial^+$
by the isomorphism type of $M^{[\tau]}$; that is during the
construction, choosing $(\bar M^\alpha,\bar{\bold J}^\alpha,\bold
f^\alpha) \in K^{\text{qt}}_{\frak u}$ which are $\le^{\text{qs}}_{\frak
u}$-increasing with $\alpha < \partial^+$, we shall have
enough free decisions.  This means that, arriving to the $\alpha$-th triple we
have continuations which are incompatible in some sense.  This will be
done in \S2,\S3.
\nl
2) The following definition 
will help phrase coding properties which
holds just for ``almost all" triples from $K^{\text{qt}}_{\frak u}$.  Note that
in the weak version of coding we have to preserve $\bold f(\delta) = 0$ for
enough $\delta$'s.
\nl
3) In the applications we have in mind,
$\partial = \lambda^+$, the set of $(\bar M,\bar{\bold J},\bold f) \in
K^{\text{qt}}_{\frak u}$ for which $M_{\lambda^+} = 
\cup\{M_i:i < \lambda^+\}$ is
saturated above $\lambda$, is dense enough which for our purpose means that for
almost every $(\bar M,\bar{\bold J},\bold f)$ this holds.
\nl
4) Central in our proof will be having ``for almost all $(\bar M,\bar{\bold
J},\bold f) \in K^{\text{qt}}_{\frak u}$ in some sense, satisfies
....".  The first version (almost$_3$, in \scite{838-1a.43}(0)), is related
to Definition \scite{838-1a.47}.
\nl
5) The version of Definition \scite{838-1a.43} we shall use mostly in
   \scite{838-1a.43}(3C), ``almost$_2$...", which means that for some
   stationary $S \subseteq \partial$, we demand the sequences to
   ``strictly $S$-obey ${\frak g}$"; and from Definition \scite{838-1a.47}
   is \scite{838-1a.47}(7), ``$\{0,2\}$-almost".
\bigskip

\definition{\stag{838-1a.43} Definition}  0) We say that ``almost$_3$ every
$(\bar M,\bar{\bold J},\bold f) \in K^{\text{qt}}_{\frak u}$ satisfies
Pr" \ub{when} there is a function ${\frak h}$ witnessing it which
means:
\mr
\item "{$(*)_1$}"  $(\bar M^\delta,\bar{\bold J}^\delta,\bold
f^\delta)$ satisfies Pr \ub{when}:  the sequence $\bold x = 
\langle(\bar M^\alpha,\bar{\bold J}^\alpha,\bold f^\alpha):
\alpha \le \delta \rangle$ is 
$\le^{\text{qs}}_{\frak u}$-increasing continuous 
($\delta < \partial^+$ a limit ordinal, of course) and obeys ${\frak h}$
which means that for some
unbounded subset $u$ of $\delta$ for every $\alpha \in u$ the sequence
$\bold x$ does
obey$_3$ or 3-obeys ${\frak h}$ which means $(\bar M^{\alpha +1},
\bar{\bold J}^{\alpha +1},\bold f^{\alpha +1}) = {\frak h}((\bar
M^\alpha,\bar{\bold J}^\alpha,\bold f^\alpha)$ (and for notational
simplicity the universes of $M^\alpha_\partial,M^{\alpha +1}_\partial$
are sets of ordinals\footnote{we may alternatively restrict yourself to models
with universe $\subseteq \partial^+$ or use a universal choice
function.  Also if we use ${\frak h}
(\langle \bar M^\beta,\bar{\bold J}^\beta,\bold f^\beta):\beta \le
\alpha\rangle)$ the difference is minor: make the statement a little
cumbersome and the checking a little easier.  Presently we do not
distinguish the two versions.}); we may write obeys instead 3-obeys
when this is clear from the context; also below
\sn
\item "{$(*)_2$}"  above $\bold f^\alpha(i) = 0 \Rightarrow 
\bold f^{\alpha +1}(i) = 0$ for a club of $i < \partial$ and
$\{i:\bold f^\alpha(i) > 0\}$ is stationary for every $\alpha \le \delta$
\sn
\item "{$(*)_3$}"  ${\frak h}$ has the domain and range implicit in
$(*)_1 + (*)_2$
\sn
\item "{$(*)_4$}"  we shall restrict ourselves to a case where each
of the models
$M^\alpha_i$ above have universe $\subseteq \partial^+$, (or just be a
set of ordinals) thus avoiding the problem of global choice; similarly
below (e.g. in part (3)).
\ermn
1) We say that the pair $((\bar
M^1,\bar{\bold J}^1,\bold f^1),(\bar M^2,\bar{\bold J}^2,\bold f^2))$
does $S$-obey or $S$-obeys$_1$ the function 
${\frak g}$ (or $(\bar M^2,\bar{\bold J}^2,\bold f^2)$ does $S$-obeys
or $S$-obeys$_1 \,{\frak g}$ above $(\bar M^1,\bar{\bold J},
\bold f^1))$, \ub{when} for some $\bar{\bold I}$ and $E$
we have
\mr
\item "{$(a)$}"  $S$ is a stationary subset of $\partial$ and $E$ is a
club of $\partial$
\sn
\item "{$(b)$}"  $(\bar M^1,\bar{\bold J}^1,\bold f^1)
<^{\text{at}}_{\frak u} (\bar M^2,\bar{\bold J}^2,\bold f^2)$ as
witnessed by $E$ and $\bar{\bold I}$
\sn
\item "{$(c)$}"  for stationarily many $\delta \in S$
{\roster
\itemitem{ $\odot$ }   the triple
$(\bar M^2 \restriction (\delta + \bold f^2(\delta)+1),\bar{\bold J}^2
\restriction (\delta + \bold f^2(\delta)),\bar{\bold I} \restriction
(\delta + \bold f^2(\delta)) +1)$ is equal to, (in particular
\footnote{alternatively we can demand (as in \S9,\S10) that: 
the universe of $M^1_\partial$ and of
$M^2_\partial$ is an ordinal $< \partial^+$} ${\frak g}$ is well defined in
this case) ${\frak g}(\bar M^1,\bar{\bold J}^1,\bold f^1,
\bar M^2 \restriction (\delta + \bold f^1(\delta) +1),
\bar{\bold J}^2 \restriction (\delta + \bold f^1(\delta),
\bar{\bold I} \restriction (\delta + \bold f^1(\delta) +1),S)$
\nl
or at least
\sn
\itemitem{ $\odot'$ }  for some $\gamma_1 \le \gamma_2$ from the interval
$[\bold f^1(\delta),\bold f^2(\delta)]$, the triple $(\bar M^2 \rest
(\gamma_2 +1),\bar{\bold J}_2 \rest \gamma_2,\bar{\bold I} \rest
(\gamma_2 +1))$ is equal to ${\frak g}(\bar M^1,\bar{\bold J}^1,\bold
f^1,\bar M^2 \rest (\gamma_1 +1)),\bar{\bold J}^2 \rest
\gamma_1,\bar{\bold I} \rest (\gamma_1 +1),S)$. 
\endroster}
\ermn
1A) Saying ``strictly $S$-obeys$_1$" mean that in clause (1)(c) we
replace ``stationarily many $\delta \in S$" by ``every $\delta \in E
\cap S$ (we can add the ``strictly" in other places, too).
 Omitting $S$ means for some stationary $S \subseteq \partial$; we  
may assume ${\frak g}$ codes $S$ and in this case we 
write $S = S_{\frak g}$ and can omit $S$.  In
the end of clause (1)(c), if the resulting value does not depend on some of
the objects written as arguments we may omit them.
We may use $\bar{\frak g} = \langle {\frak g}_S:S \subseteq
\partial$ stationary$\rangle$  and obeying $\bar{\frak g}$ means obeying
${\frak g}_S$ for some $S$ (where ${\frak g} = {\frak g}_S \Rightarrow
S_{\frak g} = S$).
\nl
2) A $\le^{\text{qs}}_{\frak u}$-increasing continuous sequence 
$\langle(\bar M^\alpha,\bar{\bold J}^\alpha,\bold f^\alpha):\alpha \le \delta
\rangle$ obeys$_1$ or 1-obeys $\bar{\frak g}$ when 
$\delta$ is a limit ordinal $<
\partial^+$ and for some unbounded $u \subseteq \delta$ there is a sequence
$\langle S_\alpha:\alpha \in u \cup \{\delta\}\rangle$ of stationary subsets of
$\partial$ decreasing modulo ${\Cal D}_\partial$ such that
for each $\alpha \in u$, the pair $((\bar M^\alpha,\bar{\bold J}^\alpha,
\bold f^\alpha),(\bar M^{\alpha +1},\bar{\bold J}^{\alpha +1},
\bold f^{\alpha +1}))$ strictly $S_\alpha$-obeys ${\frak g}$.
\nl
2A) In part (2) we say  $S$-obeys$_1$ when $S_\alpha \subseteq S$ mod
${\Cal D}_\partial$ for $\alpha \in u \cup \{\delta\}$.  Similarly for
$\bar S'$-obey$_1$ when $\bar S' = \langle S'_\alpha:\alpha \in
u'\rangle$ and $\alpha \in u \cup \{\delta\} \Rightarrow S^*_\alpha
\subseteq S_\alpha$ and $u \cup \{\delta\} \subseteq u'$.
\nl
2B) In part (2) we say strictly $S$-obeys$_1$ when this holds in each
case.
\nl
3) We say ``almost$_1$ every $(\bar M,\bar{\bold J},\bold f) \in
K^{\text{qt}}_{\frak u}$ satisfies Pr" \ub{when} 
there is a function ${\frak g}$ witnessing it, which means (note: the
use of ``obey" guarantees ${\frak g}$ is as in part (2) and not as
implicitly required on ${\frak h}$ in part (0)):
\mr
\item "{$(c)$}"  if $\langle (\bar M^\alpha,\bar{\bold J}^\alpha,\bold
f^\alpha):\alpha \le \delta \rangle$ is $\le^{\text{qs}}_{\frak u}$-increasing
continuous obeying ${\frak g}$  and $\delta < \partial^+$ a limit
ordinal \ub{then} $(\bar M^\delta,\bar{\bold
J}^\delta,\bold f^\delta)$ satisfies the property Pr.
\ermn
3A) We add ``above $(\bar M,\bar{\bold J},\bold f')$" \ub{when} we demand in
clause (c) that $(\bar M^0,\bar{\bold J}^0,\bold f^0) = (\bar
M',\bar{\bold J}',\bold f')$.
\nl
3B) We replace \footnote{again assume that all elements are
ordinals $< \partial^+$} almost$_1$  by $S$-almost$_2$ \ub{when} we 
require that the sequence ``strictly $S$-obeys ${\frak g}$".
\nl
3C) We replace\footnote{if we replaced it by ``for a set
of $\delta$'s which belongs to ${\Cal D}$", ${\Cal D}$ a normal filter
on $\partial$, the difference is minor.}
 almost$_1$ by almost$_2$ \ub{when} for every stationary
$S \subseteq \partial,S$-almost$_1$ every triple 
$(\bar M,\bar{\bold J},\bold f) \in
K^{\text{qt}}_{\frak u}$ satisfies Pr; and ``$S$-almost$_2$" we ?.
\enddefinition
\bigskip

\definition{\stag{838-1a.45} Definition}  1) For ${\frak h}$ as \footnote{we
shall assume that no ${\frak h}$ is both as required in Definition
\scite{838-1a.43} and as required in Definition \scite{838-1a.45}(0).} in
Definition \scite{838-1a.43}(0) we define 
$K^{{\frak u},{\frak h}}_{\partial^+}$ 
as the class of models $M$ such that for some
$\le^{\text{qs}}_{\frak u}$-increasing continuous sequence $\bold x =
\langle(\bar M^\alpha,\bar{\bold J}^\alpha,\bold f^\alpha):\alpha <
\partial^+\rangle$ of members of $K^{\text{qt}}_{\frak u}$ such that a
club of $\delta < \partial^+,\bold x \restriction (\delta +1)$ obeys
${\frak h}$ is the sense of part (0) of Definition \scite{838-1a.43}
respectively, we have $M = \cup\{M^\alpha_\partial:\alpha <
\partial^+\}$.
\nl
2) For ${\frak g}$ as in Definition \scite{838-1a.43}(1),(2) we define
$K^{{\frak u},{\frak g}}_{\partial^+}$ similarly.
\nl
3) We call ${\frak h}$ as in \scite{838-1a.43}(0) appropriate$_3$ or
3-appropriate and ${\frak g}$ as in Definition \scite{838-1a.43}(1),(2) 
we call appropriate$_\ell$ or $\ell$-appropriate for $\ell=1,2$; we 
may add ``${\frak u}$-"if not clear from the context. 
\nl
4) As in parts (1),(2) for ${\frak h}$ as in (any relevant part of)
   Definition \scite{838-1a.47} below.
\nl
5) Also $K^{{\frak u},\bar{\frak h}}_{\partial^+} = \cap\{K^{{\frak
   u},{\frak h}_\varepsilon}_{\partial^+}:\varepsilon < \ell
   g(\bar{\frak h})\}$ where each $K^{{\frak u},{\frak
   h}_\varepsilon}_{\partial^+}$ is well defined.
\enddefinition
\bigskip

\definition{\stag{838-1a.47} Definition}   1) We say
$\langle(\bar M^\zeta,\bar{\bold J}^\zeta,\bold f^\zeta):\alpha 
< \alpha_*\rangle$ does obey$_0$ (or $0$-obey) the function ${\frak h}$ in
$\zeta$ \ub{when} $\xi + 1 < \alpha_*$ and 
(if $\alpha_* =2$ we can omit $\zeta$):
\mr
\item "{$(a)$}"  $(\bar M^\varepsilon,\bar{\bold J}^\varepsilon,
\bold f^\varepsilon) \in K^{\text{qt}}_{\frak u}$ is
$\le^{\text{qt}}_{\frak u}$-increasing continuous (with $\varepsilon$)
\sn
\item "{$(b)$}"  $M^\zeta_\partial$ and even $M^{\zeta +1}_\partial$ 
has universe an ordinal $< \partial^+$
\sn
\item "{$(c)$}"  there is a club $E$ of $\partial$ and sequence
$\langle \bold I_\alpha:\alpha < \partial\rangle$ witnessing $(\bar
M^\zeta,\bar{\bold J}^\zeta,\bold f^\zeta) \le^{\text{at}}_{\frak u}
(\bar M^{\zeta +1},\bar{\bold J}^{\zeta +1},\bar{\bold J},\bold
f^{\zeta +1})$ such that $(M^\zeta_{\text{min}(E)},
M^{\zeta +1}_{\text{min}(E)},\bold I_{\text{min}(E)})
= {\frak h}(\langle(\bar M^\xi,\bar{\bold J}^\xi,\bold f^\xi):\xi \le
\zeta\rangle) \in \text{ FR}^+_1$.
\ermn
2) We say that ${\frak h}$ is ${\frak u}$-appropriate$_0$ or 
${\frak u}-0$-appropriate \ub{when}: 
${\frak h}$ has domain and range as required in part (1), particularly
clause (c).
We may say 0-appropriate or appropriate$_0$ when ${\frak u}$ 
is clear from the context and we say ``$(\bar M^\zeta,\ldots),(\bar
M^{\zeta +1},\ldots)$ does $0$-obeys ${\frak h}$".
\nl
2A) We say the function ${\frak h}$ is ${\frak u}$-1-appropriate \ub{when} its
domain and range are as required in Definition \scite{838-1a.43}(3); in
this case $S_{\frak h} = S$.
\nl
2B) We say the function ${\frak h}$ is 
${\frak u}$-2-appropriate for $S$ \ub{when}
$S \subseteq \partial$ is stationary and its domain and range are as
required in Definition \scite{838-1a.43}(3B), i.e. \scite{838-1a.43}(3).
\nl
2C) If in (2B) we omit $S$ this means that $\bar{\frak h} = \langle
{\frak h}_S:S \subseteq \partial$ is stationary$\rangle$, each ${\frak
h}_S$ as above.
\nl
3) For 0-appropriate ${\frak h}$ we define ${\frak K}^{{\frak
u},{\frak h}}_{\partial^+}$ to be the family of models $M$, with
universe $\partial^+$ for simplicity, as the set of models of the form
$\cup\{M^\zeta_{\partial^+}:\zeta < \partial^+\}$ where 
$\langle(\bar M^\zeta,\bar{\bold J}^\zeta,\bold f^\zeta):\zeta <
\partial^+\rangle$ is $\le_{\text{qt}}$-increasing continuous and
0-obeys ${\frak h}$ in $\zeta$ for unboundedly many $\zeta <
\partial^+$.  Similarly for the other ${\frak h}$, see below.
\nl
4) We say ${\frak h}$ is ${\frak u}-\{0,2\}$-appropriate or ${\frak
u}$-appropriate for $\{0,2\}$ if ${\frak h} = {\frak h}_0 \cup {\frak
h}_2$ and ${\frak h}_\ell$ is $\ell$-appropriate for $\ell=0,2$; we
may omit ${\frak h}$ when clear from the context.
\nl
5) For a $\{0,1\}$-appropriate ${\frak h}$ letting ${\frak h}_0,{\frak
h}_1$ be as above we say $\langle (\bar M^\alpha,\bar{\bold
J}^\alpha,\bold f^\alpha):\alpha < \alpha(*)\rangle$ does
$\{0,1\}$-obeys ${\frak h}$ in $\zeta < \alpha(*)$ when $((\bar
M^\zeta,\bar{\bold J}^\zeta,\bold f^\zeta),(\bar M^{\zeta
+1},\bar{\bold J}^{\zeta +1},\bold f^\zeta))$ does $\ell$-obey ${\frak
h}_\ell$ for $\ell=0,2$.  We say strictly $\{0,2\}-S$-obeys ${\frak
h}$ in $\zeta$ when for stationary $S \subseteq \partial$, for
unboundedly many $\zeta < \alpha(*)$ the
pair $0$-obeys ${\frak h}_0$ and strictly $1-S$-obeys ${\frak h}_1$.
\nl
6) For a $\{0,2\}$-appropriate ${\frak h}$, we say $\langle(\bar
M,\bar{\bold J}^\alpha,\bold f^\alpha):\alpha < \delta \le
\partial^+\rangle$ does $\{0,2\}$-obey ${\frak h}$ when this holds for
some stationary $S \subseteq \partial$ for unboundedly many 
$\zeta < \delta$ the sequence strictly $\{0,1\}-S$-obey ${\frak h}$.
Similarly we define ``the sequence $\{0,2\}-S$-obeys ${\frak h}$". 
\nl
7) ``$\{0,2\}$-almost every $(\bar M,\bar{\bold J},\bold f)$ (or every
$(\bar M,\bar{\bold J},\bold f)$ above $(\bar M^*,\bar{\bold
J}^*,\bold f^*$))" is defined similarly to Definition \scite{838-1a.43}.
\enddefinition
\bigskip

\demo{\stag{838-1a.49} Observation}  1) For any $\varepsilon^* < \partial^+$
and sequence $\langle {\frak h}_\varepsilon:
\varepsilon < \varepsilon^*\rangle$ of 3-appropriate ${\frak h}$, 
there is an 3-appropriate
${\frak h}$ such that ${\frak K}^{{\frak u},{\frak h}}_{\partial^+}
\subseteq \cap\{K^{{\frak u},
{\frak h}_\varepsilon}_{\partial^+}:\varepsilon < \varepsilon^*\}$ and
similarly for $<^{\text{qs}}_{\frak u}$-increasing sequences of
$K^{\text{qt}}_{\frak u}$ length $< \partial^+$.
\nl
2) $K^{{\frak u},{\frak h}}_{\partial^+} \subseteq 
K^{{\frak u},*}_{\partial^+}$ for any 3-appropriate function ${\frak h}$.
\nl
3) Similarly to parts (1)+(2) for ${\frak g}$ as in Definition
\scite{838-1a.43}(2). 
\nl
4) Similarly to parts (1) + (2) for $\{0,2\}$-appropriate ${\frak h}$,
see Definition \scite{838-1a.47}(4),(5),(6).
\enddemo
\bigskip

\remark{\stag{838-1a.51} Remark}  1) Concerning \scite{838-1a.49}, if in
Definition \scite{838-1a.43}(1)(c) we do not allow $\odot'$, then we
better\footnote{in the cases we would like to
apply \scite{838-1a.49} there is no additional price for this.} in
\scite{838-1a.43}(2) add $\bar S = \langle S^\varepsilon:\varepsilon <
\partial\rangle$ such that: $S^\varepsilon \subseteq \partial$ is
stationary, $\varepsilon < \zeta < \lambda \Rightarrow S^\varepsilon
\cap S^\zeta = \emptyset$ and $\varepsilon < \partial \wedge \alpha
\in u \cup \{\delta\} \Rightarrow S^\varepsilon \cap S_\alpha$ is
stationary.
\nl
2) A priori ``almost$_3$" look the most
natural, but we shall use as our main case ``$\{0,2\}$-almost".  We try to
explain below.
\nl
3) Note that
\mr
\item "{$(a)$}"  in the proof of e.g. \scite{838-10k.17} we use
$K^{\text{rt}}_{\frak u}$ not $K^{\text{qt}}_{\frak u}$, i.e. carry
$\bar{\Bbb F}$; this does not allow us the freedom which ``almost$_3$"
require
\sn
\item "{$(b)$}"  model theoretically here usually there is a special model in
$K^{\frak u}_\partial$, normally the superlimit or
saturated one, and we try to take care 
building the tree $\langle(\bar M_\eta,\bar{\bold J}_\eta,\bold
f_\eta,(\bar{\Bbb F}_\eta)):\eta \in {}^{\partial >}
(2^\partial)\rangle$ that, e.g. $\eta \in {}^\gamma(2^\partial)
\wedge \partial|\gamma \Rightarrow M^\eta_\partial$ is saturated.
\ermn
In the `almost$_3$" case this looks straight; in successor of
successor cases we can take care. 
\nl
4) We like to guarantee that for ``almost" all
$(\bar M^\eta,\bar{\bold J}^\eta,\bold f^\eta)$ the model 
$M^\eta_\partial \in K^{{\frak u},*}_\partial$ is saturated so that 
we have essentially one case.
If we allow in the ``almost", for, e.g. $\gamma +2$, to
choose some initial segment in $(\bar M^\eta,\bar{\bold J}^\eta,\bold
f^\eta)$ for $\eta$ of length $\gamma +1$, 
this guarantees saturation of $M^\eta_\partial$ if cf$(\ell g(\eta)) =
\partial$, but
\mr
\item "{$(c)$}"  set theoretically we do not know that
$S^{\partial^+}_\partial = \{\delta < \partial^+:\text{cf}(\delta) =
\partial\}$ is not in the relevant ideal (in fact, even under GCH,
$\diamondsuit_{S^{\partial^+}_\partial}$ may fail)
\sn
\item "{$(d)$}"  if $K_\partial$ is categorical, 
there is no problem.  However, if we know less, e.g. that there is a
superlimit one, or approximation, using the almost$_2$, in $\gamma =
\gamma' + 2$, we can guarantee that $M_\eta$ for $\eta \in
{}^\gamma(2^\partial)$ is up to isomorphism the superlimit one
\sn
\item "{$(e)$}"  we may conclude that it is better to work with
$K^{\text{rt}}_{\frak u}$ rather than $K^{\text{qt}}_{\frak u}$, see
Definition \scite{838-10k.5}(1).  This
is true from the point of view of the construction but it is model
theoretically less natural.
\ermn
5) We may in Definition \scite{838-1a.47} demand on $\langle(\bar
   M^\alpha,\bar{\bold J}^\alpha,\bold f^\alpha):\alpha <
   \delta_*\rangle$ satisfies several ${\frak h}$'s of different kinds
   say of $\{0,2\}$ and of 3; make little difference.
\nl
6) In the usual application here for ${\frak u} = {\frak
   u}^\ell_{\frak s}$ for some ${\frak g}$, if $\langle
   (M^\alpha,\bar{\bold J}^\alpha,\bold f^\alpha):\alpha \le
   \delta\rangle$ is $\le^{\text{qt}}_{\frak u}$-increasing continuous
   and $\delta = u,u := \{\alpha < \delta:((\bar M^\alpha,\bar{\bold
   J}^\alpha,\bold f^\alpha),(\bar M^{\alpha +1},\bar{\bold
   J}^{\alpha+1},\bold f^{\alpha_1}))$ does strictly 
$S$-obey ${\frak g}\}$, \ub{then}
   $M^\delta_\partial$ is saturated.  But without this extra
   knowledge, the fact that for $\alpha \in u$ we may have $S_\alpha$
   disjoint to other may be hurdle.  But using ``strictly obey$_1$"
   seems more general and the definition of ``almost$_2$" fits this feeling.
\endremark
\goodbreak

\head {\S2 Coding properties and non-structure} \endhead  \resetall \sectno=2
 \spuriousreset
\bigskip

We now come to the definition of the properties we shall use as
sufficient conditions for non-structure starting with Definition
\scite{838-2b.1}; in this section and \S3 we shall define also some
relatives needed for sharper results, those properties have parallel
cases as in Definition  \scite{838-2b.1}.
\demo{\stag{838-2b.0} Hypothesis}  We assume
${\frak u}$ be a nice construction framework and, 
$\tau$ a weak ${\frak u}$-sub-vocabulary, see Definition \scite{838-1a.15}(1). 
\enddemo
\bigskip

\remark{Remark}  The default 
value is $\tau_{\frak u} = \tau({\frak K}_{\frak u})$ or better the
pair $(\tau,\tau_{\frak u})$ such that
$\tau_{\frak u} = \tau'$, as in Definition \scite{838-1a.15}(1) and
\scite{838-2b.8}(1),(2); see also ${\frak u}$ has faked equality, see
\scite{838-3r.84} later.

Among the variants of weak $\tau$-coding in Definition \scite{838-2b.1}
the one we shall use most is \scite{838-2b.1}(5), ``${\frak u}$ has the
weak $\tau$-coding$_1$ above $(\bar M^*,\bar{\bold J}^*,\bold f^*)$".
\endremark
\bigskip

\definition{\stag{838-2b.1} Definition}  1) We say that $M \in K_{\frak u}$ 
has the weak $\tau$-coding$_0$-property (in ${\frak u}$) \ub{when}:
\mr
\item "{$(A)$}"   if $N,\bold I$ are such that $(M,N,\bold I) \in \text{
FR}^+_1$ \ub{then} $(M,N,\bold I)$ has the weak $\tau$-coding$_0$ property, 
\nl
where:
\sn
\item "{$(B)$}"  $(M,N,\bold I)$ has the weak $\tau$-coding$_0$
property \ub{when} we can find
$(M_*,N_\ell,\bold I_\ell) \in \text{ FR}_1$ for $\ell=1,2$ satisfying
{\roster
\itemitem{ $(a)$ }  $(M,N,\bold I) \le^1_{\frak u} 
(M_*,N_\ell,\bold I_\ell)$ for $\ell=1,2$
\sn
\itemitem{ $(b)$ }  $M_* \cap N = M$ (follows)
\sn
\itemitem{ $(c)$ }  $N_1,N_2$ are $\tau$-incompatible amalgamations of
$M_*,N$ over $M$ in $K_{\frak u}$, (see Definition \scite{838-1a.15}(4)).
\endroster}
\ermn
1A) We say that $(M,N,\bold I) \in \text{ FR}^+_1$ has the true weak
$\tau$-coding$_0$ \ub{when}: if $(M,N,\bold I) \le^1_{\frak u} (M',N',\bold
I)$ then $(M',N',\bold I')$ has the weak $\tau$-coding$_0$ property,
i.e. satisfies the requirement in 
clause (B) of part (1).
\nl
1B) ${\frak u}$ has the explicit weak $\tau$-coding$_0$ property \ub{when}
every $(M,N,\bold I) \in \text{ FR}^+_{\frak u}$ has the weak $\tau$-coding
property.
\nl
2) $(\bar M^*,\bar{\bold J}^*,\bold f^*) \in K^{\text{qt}}_{\frak u}$ has
the weak $\tau$-coding$_0$ property \ub{when}: for a 
club of $\delta < \partial$, not only $M = M^*_\delta$ has the true 
weak $\tau$-coding$_0$-property \ub{but} in
clause (B) of part (1) above we demand $M_* \le_{\frak u}
M^*_\gamma$ for any $\gamma < \partial$ large enough. 
\nl
3) We say that $(\bar M,\bar{\bold J},\bold f) \in
K^{\text{qt}}_{\frak u}$ has the weak
$\tau$-coding$_1$ property \footnote{the difference between coding$_0$
and coding$_1$ may seem negligible but it is crucial, e.g. in \scite{838-e.1}}
\ub{when}: (we may omit the superscript 1): recalling 
$M_\partial = \cup\{M_\alpha:\alpha < \partial\}$,
there are $\alpha(0) < \partial$ and $N_0,\bold I_0$ such that
$(M_{\alpha(0)},N_0,\bold I_0) \in \text{ FR}_1,N_0 \cap M_\partial =
M_{\alpha(0)}$ and for a club of $\alpha(1) < \partial$, if
$(M_{\alpha(0)},N_0,\bold I_1) \le^1_{\frak u} (M_{\alpha(1)},N_1,\bold I_1)$
satisfies $N_1 \cap M_\partial = 
M_{\alpha(1)}$ \ub{then} there are $\alpha(2) \in
(\alpha(1),\partial)$ and $N^\ell_2,\bold I^\ell_2$ for $\ell=1,2$ such that
$(M_{\alpha(1)},N_1,\bold I_1) \le^1_{\frak u} (M_{\alpha(2)},N^\ell_2,\bold
I^\ell_2)$ for $\ell=1,2$ and $N^1_2,N^2_2$ are $\tau$-incompatible
amalgamations of $M_{\alpha(2)},N_1$ over $M_{\alpha(1)}$ recalling
Definition \scite{838-1a.15}(4).
\nl
4) We say that $(\bar M,\bar{\bold J},\bold f)$ has the $S$-weak
$\tau$-coding$_1$ property \ub{when}: $S$ is a stationary subset of
$\partial$ and for some club $E$ of $\partial$ the demand in (3)
holds restricting ourselves to $\alpha(1) \in S \cap E$.
\nl
5) We say that ${\frak u}$ has the weak $\tau$-coding$_k$ property
\ub{when}: $\{0,2\}$-almost every $(\bar M,\bar{\bold J},\bold f) \in 
K^{\text{qt}}_u$ has the weak $\tau$-coding$_k$ property; omitting $k$
means $k=1$.  Similarly for ``above
$(\bar M,\bar{\bold J},\bold f) \in K^{\text{qt}}_{\frak u}$".
Similarly for ``$S$-weak".
\enddefinition
\bigskip

The following theorem uses a weak model theoretic assumption, but the
price is a very weak but still undesirable, additional set theoretic 
assumption (i.e. clause (c)), recall that
$\mu_{\text{unif}}(\partial^+,2^\partial)$ is defined in
\scite{838-0z.6}(7), see \scite{838-7f.14}.
\proclaim{\stag{838-2b.3} Theorem}   We have
$\dot I_\tau(\partial^+,K^{\frak u}_{\partial^+}) \ge
\mu_{\text{unif}}(\partial^+,2^\partial)$, moreover for any 
${\frak u}$-0-appropriate ${\frak h}$ (see Definition \scite{838-1a.47})
and even $\{0,2\}$-appropriate ${\frak h}$ (see Definition
\scite{838-1a.47}(3),(7) and Definition \scite{838-1a.45}) we have
$\dot I(K^{{\frak u},{\frak h}}_{\partial^+}) \ge
\mu_{\text{unif}}(\partial^+,2^\partial)$, \ub{when}:
\mr
\item "{$(a)$}"  $2^\theta = 2^{< \partial} < 2^\partial$ 
\sn
\item "{$(b)$}"  $2^\partial < 2^{\partial^+}$
\sn
\item "{$(c)$}"  the ideal
{\rm WDmId}$(\partial)$ is not $\partial^+$-saturated
\sn
\item "{$(d)$}"  ${\frak u}$ has the weak $\tau$-coding (or just the
$S$-weak $\tau$-coding property above some 
triple $(\bar M,\bar{\bold J},\bold f) \in
K^{\text{qt}}_{\frak u}$ with {\rm WDmId}$(\partial) \restriction S$
not $\partial^+$-saturated and $S \subseteq \bold f^{-1}\{0\}$).
\endroster
\endproclaim
\bigskip

\demo{Proof}  This is proved in \scite{838-10k.17}.  \hfill$\square_{\scite{838-2b.3}}$
\enddemo
\bigskip

\remark{\stag{838-2b.5} Remark}  1) Theorem \scite{838-2b.3} is 
used in \scite{838-e.1}, \scite{838-e.1Z}(1),
\scite{838-e.4}, \scite{838-e.3} and \scite{838-e.5}, for \scite{838-e.1Z}(2) we use
the variant \scite{838-2b.6}.  We
could use Theorem \scite{838-2b.7} below to get a somewhat stronger result.
\nl
In other words, e.g. it is used for 
``the minimal types are not dense in ${\Cal S}(M)$ 
for $M \in {\frak K}_\lambda$" for suitable ${\frak K}$, see
\scite{838-e.1} (and \chaptercite{E46} or the older \cite{Sh:576}, \cite{Sh:603}).  
\nl
2) We may think that here at a minor set theoretic 
price (clause (c)), we get the strongest model theoretic version.
\nl
3) We can in \scite{838-2b.3} replace ${\frak h}$ by $\bar{\frak h}$, a
   sequence of $\{0,2\}$-appropriate ${\frak h}$'s of length $\le \partial^+$.
\nl
4) In part (3), we can fix a stationary $S \subseteq \partial$ such
that WDmId$(\partial) +S$ is not $\partial^+$-saturated (so
$\partial$ is not in it) and restrict ourselves to strict $S$-obeying.
\nl
5)  We can replace assumption (d) of \scite{838-2b.3} by 
\mr
\item "{$(d)'$}"  for some ${\Cal D}$
{\roster
\itemitem{ $(i)$ }   ${\Cal D}$ is a normal filter on $\partial$ disjoint
to WDmId$(\partial)$; moreover $(\forall A \in {\Cal D})(\exists B)[B
\subseteq A \wedge \partial \backslash B \in D \wedge B \in
(\text{WDmId}(\partial))^+]$
\sn
\itemitem{ $(ii)$ }  almost$_2$ every $(\bar M,\bar{\bold J},\bold f) \in
K^{\text{qt}}_{\Cal D}$ has the weak $\tau$-coding property 
(even just above some member of $K^{\text{qt}}_{\Cal D}$).
\endroster}
\endroster
\endremark
\bn
A variant is
\proclaim{\stag{838-2b.6} Claim}   In Theorem \scite{838-2b.3} we can weaken
the assumption to ``${\frak u}$ has a weak $\tau$-coding$_2$", see below.
\endproclaim
\bigskip

\demo{Proof}  As in \scite{838-10k.17}.  \hfill$\square_{\scite{838-2b.6}}$
\enddemo
\bigskip

\definition{\stag{838-2b.6F} Definition}  1) We say that ${\frak u}$ has
the $S$-weak $\tau$-coding$_2$ property (or the $S$-weak game
$\tau$-coding property) [above $(\bar M^*,\bold J^*,\bold f^*) \in
K^{\text{qt}}_{\frak u}$] \ub{when} $\{0,2\}$-almost every $(\bar
M,\bar{\bold J},\bold f) \in K^{\text{qt}}_{\frak s}$ [above $(\bar
M^*,\bar{\bold J}^*,\bold f^*)$] has it.
\nl
2) We say $(\bar M,\bar{\bold J},\bold f) \in
K^{\text{qt}}_{\frak u}$ has the $S$-weak game $\tau$-coding$_2$ property
\ub{or} $S$-weak $\tau$-coding$_2$ property for a stationary set 
$S \subseteq \partial$ (omitting $S$
means for every such $S$) \ub{when}, recalling 
$M_\partial = \cup\{M_\alpha:\alpha <
\partial\}$, in the following game $\Game_{{\frak u},S}(\bar
M,\bar{\bold J},\bold f)$, the Coder player has a winning strategy where:
\mr
\item "{$(*)_1$}"  a play of $\Game_{{\frak u},S}$ last $\partial$
moves after the $\varepsilon$-th move a tuple
$(\alpha_\varepsilon,e_\varepsilon,\bar N^\varepsilon,\bar{\bold
J}^\varepsilon,\bold f^\varepsilon,\bar I^\varepsilon)$ is chosen such
that:
{\roster
\itemitem{ $(a)$ }  $\alpha_\varepsilon < \partial$ is increasing
continuous
\sn
\itemitem{ $(b)$ }  $e_\varepsilon$ is a closed subset of
$\alpha_\varepsilon$ such that $\zeta < \varepsilon \Rightarrow
\alpha_\zeta \in e \wedge e_\zeta = e_\varepsilon \cap \alpha_\zeta$
\sn
\itemitem{ $(c)$ }  $\bold f^\varepsilon$ is a function with domain
$e_\varepsilon$ such that $\alpha + \bold f^\varepsilon(\alpha) < \text{
min}(e_\varepsilon \cup \{\alpha_\varepsilon\} \backslash \alpha)$ and
$\bold f^\varepsilon(\alpha) \ge \bold f(\alpha)$
\sn
\itemitem{ $(d)$ }  $u_\varepsilon = \cup\{[\alpha,\alpha + \bold
f^\varepsilon(\alpha)]:\alpha \in e_\varepsilon\} \cup
\{\alpha_\varepsilon\}$ and $u^-_\varepsilon = \cup\{[\alpha,\alpha +
\bold f^\varepsilon(\alpha))):\alpha \in e_\varepsilon\}$
\sn
\itemitem{ $(e)$ }  $\bar N^\varepsilon = \langle N_\alpha:\alpha \in
u_\varepsilon\rangle$ and $\bold J^\varepsilon = \langle \bold
J_\alpha:\alpha \in u^-_\varepsilon\rangle$ and
$\bar{\bold I}^\varepsilon= \langle \bold I_\alpha:\alpha \in
u_\varepsilon\rangle$
\sn
\itemitem{ $(f)$ }  $\langle(M_\alpha,N_\alpha,\bold I_\alpha):
\alpha \in u_\varepsilon\rangle$ is $\le^1_{\frak u}$-increasing
\sn
\itemitem{ $(g)$ }  $N^\varepsilon_\alpha \cap M_\partial = M_\alpha$ for
$\alpha \in u_\varepsilon$
\sn
\itemitem{ $(h)$ }  $\bar{\bold J}^\varepsilon = \langle
\bold J^*_\alpha:\alpha \in u^-_\varepsilon\rangle$
\sn
\itemitem{ $(i)$ }  $(M_\alpha,M_{\alpha +1},\bold J_\alpha)
\le^2_{\frak u} (N_\alpha,N_{\alpha +1},\bold
J^*_\alpha)$ for $\alpha \in u^-_\varepsilon$
\sn
\itemitem{ $(j)$ }  the coder chooses
$(\alpha_\varepsilon,e_\varepsilon,\bar N^\varepsilon,\bar{\bold
J}^\varepsilon,\bar{\bold I}^\varepsilon,\bold f^\varepsilon)$ if
$\varepsilon=0$ or $\varepsilon = \zeta +1,\zeta$ a limit ordinal
$\notin S$, and otherwise the anti-coder chooses
\endroster}
\item "{$(*)_2$}"  in the end the Coder wins the play \ub{when} for a club of
$\varepsilon < \partial$, if $\varepsilon \in S$, then the triple
$(M_{\alpha_\varepsilon},N_\varepsilon,\bold I_\varepsilon)$ has the
weak $\tau$-coding$_0$ property, i.e. satisfies clause (B) of
Definition \scite{838-2b.1}(1), moreover such that $M \le_{\frak K}
M_\partial$.
\endroster
\enddefinition
\bn
We can also get ``no universal" over $M_\partial \in 
{\frak K}^{\frak u}_\partial$ (suitable for applying \scite{838-7f.9}).
\proclaim{\stag{838-2b.7} Claim}  If $(\bar M,\bold J,\bold f) \in
K^{\text{qt}}_{\frak u},M = \cup\{M_\alpha:\alpha < \partial\}$ and $M
\le_{\frak u} N_\varepsilon \in {\frak K}_{\le \mu}$ for $\varepsilon
< \varepsilon^* < \mu^+$ \ub{then} there is $(\bar M',\bar{\bold
J}',\bold f)$ satisfying $(\bar M,\bar{\bold J},\bold f) 
\le^{\frak u}_{\text{at}} (\bar M',\bar{\bold J}',\bold f')$ such that
$M'_\partial =
\cup\{M_\alpha:\alpha < \partial\} \in K_\partial$ cannot be
$\le_{{\frak K}[{\frak u}]}$-embedded into $N_\varepsilon$ for $\varepsilon <
\partial$ over $M_\partial$ \ub{provided that}:

\hskip10pt (a),(b),(d) \quad as in \scite{838-2b.3}

\hskip40pt (e) $\quad {}^\partial 2$ is not the union of
{\rm cov}$(\mu,\partial^+,\partial^+,2)$ sets from 
\nl

\hskip62pt {\rm WDmTId}$(\partial,2^{<\partial})$.
\endproclaim
\bigskip

\demo{Proof}  As in the proof of \scite{838-10k.17}, anyhow not used.
\hfill$\square_{\scite{838-2b.7}}$ 
\enddemo
\bn
\margintag{2b.8}\ub{\stag{838-2b.8} Exercise}:  1) [Definition] Call ${\frak u}$ a
semi-nice construction framework \ub{when} in Definition \scite{838-1a.3} we
omit clause (D)$_\ell$(d) and the disjointness demands (E)$_\ell$(b)$(\beta)$ 
\nl
2) For ${\frak u}$ as above we define ${\frak u}'$ as follows:
\mr
\item "{$(a)$}"  ${\frak K}_{{\frak u}'}$ is as in Definition \scite{838-1a.19}(1)
\sn
\item "{$(b)$}"  FR$^\ell_{{\frak u}'} = \{(M,N,\bold J):M
\le_{{\frak K}'_{\frak u}} N$ so both of cardinality $< \partial$ and $\bold J
\subseteq N \backslash M$ and letting $M^* = M/=_\tau,N^* = N/=_\tau$
and $\bold J^* = \{c/=^N_\tau:c \in \bold J$ and 
$(c/=_\tau) \notin M/=_\tau\}$ we have
$(M^*,N^*,\bold J^*) \in \text{ FR}^\ell_{\frak u}\}$
\nl
pedantically $=_\tau$ means $=^N_\tau$ (even $M/=^N_\tau$)
\sn
\item "{$(c)$}"  $(M_1,N_1,\bold J_1) \le^\ell_{{\frak u}'} 
(M_2,N_2,\bold J_1)$
\ub{iff} $M_1 \le_{{\frak u}'} M_2 \le_{{\frak u}'} N_2,M_1
\le_{{\frak u}'} N_1 \le_{{\frak u}'} N_2,\bold J_1 \subseteq 
\bold J_2,M_2 \cap N_1 = M_1$ and 
$(M^*_1,N^*_1,\bold J_1) \le^\ell_{\frak u} (M^*_2,N^*_2,\bold J^*)$ when we
 define them as in clause (b).
\ermn
3) (Claim) If ${\frak u}$ is a semi-nice construction framework
   \ub{then} ${\frak u}'$ is a nice construction framework.
\nl
4) [Definition] For ${\frak u}$ a semi-nice construction framework
   we define ${\frak u}''$ as in part (2) except that in clause (b) we
demand $c \in \bold J \wedge a \in M \Rightarrow N \models 
\neg(a =_\tau c)$.
\nl
5) [Claim] If ${\frak u}$ is a nice, [semi-nice], [semi-nice
satisfying (D)(d)] construction framework
\ub{then} ${\frak u}''$ is a nice, [semi-nice], [nice] construction framework.
\bn
\ub{Discussion}:  We now 
phrase further properties which are enough for the desired
conclusions under weaker set theoretic conditions.  
The main case is vertical coding (part (4) but it relies on part (1)
in Definition \scite{838-2b.9}).  On additional such properties, see later.

In the ``vertical coding" version (see Definition \scite{838-2b.9} below), we
strengthen the ``density of $\tau$-incompatibility" such that during
the proof we do not need to preserve 
``$\bold f^{-1}_\eta\{0\}$ is large" even allowing $\bold
f^{-1}_\eta\{0\} = \emptyset$.

We may say that ``vertically" means that given 
$(\bar M^1,\bar{\bold J}^1,\bold f^1) \in K^{\text{qt}}_{\frak u}$ 
building $M^2_\alpha,\bold J^2_\alpha,\bold I^2_\alpha$ by
induction on $\alpha < \partial$, arriving to some limit $\delta$, we
are committed to $M^2_{\alpha +i}$ for $i \le \bold f^1(\delta)$, but
still like to have freedom in determining the type of $M^2_\delta$
over $\cup\{M^1_\beta:\beta <\partial\}$ (see more in the proof of
Theorem \scite{838-8h.12} and Definition \scite{838-8h.15} on delayed
uniqueness, which express failure of this freedom).  In other
words the property we have is a delayed version of the weak coding.
\nl
As usual, always ${\frak u}$ is a nice construction framework.
\bigskip

\definition{\stag{838-2b.9} Definition}  1) We say that 
$(\bar M^1,\bar{\bold J}^1) = 
(\langle M^1_i:i \le \beta \rangle,\langle
\bold J^1_i:i < \beta \rangle)$ has the vertical $\tau$-coding$_0$ property 
(in ${\frak u}$) \ub{when}:
\mr
\item "{$(A)(a)$}"  $\beta < \partial$
\sn
\item "{$(b)_1$}"  $M^1_i$ is $\le_{\frak u}$-increasing continuous
for $i \le \beta$
\sn
\item "{$(c)_1$}"  $(M^1_i,M^1_{i+1},\bold J^1_i) \in \text{ FR}_2$
for $i < \beta$
\sn
\item "{$(B)\,\,\,\,\,\,\,$}"  if
$(\langle M^2_i:i \le \beta \rangle,\langle \bold J^2_i:i < \beta
\rangle,\langle \bold I_i:i \le \beta \rangle)$ satisfies 
$\circledast_1$ below, \ub{then} we can find $\gamma_\ell,
M^1_*,\bold I^\ell_*$ and $M^{2,\ell}_i$  
(for $i \in (\beta,\gamma_\ell])$ and $\bold J^{2,\ell}_i$ (for $i \in
[\beta,\gamma_\ell))$, for $\ell = 1,2$ satisfying $\circledast_2$ 
 below \ub{where}, letting 
$M^{2,\ell}_i = M^2_i$ for $i \le \beta$ and
$\bold J^{2,\ell}_i = \bold J^2_i$ for $i < \beta$, we have:
{\roster
\itemitem{ $\circledast_1$ }  $(d) \quad M^2_i \, (i \le \beta)$ is
$\le_{\frak u}$-increasing continuous
\sn
\itemitem{ ${{}}$ }  $(e) \quad M^2_i \cap M^1_\beta = M^1_i$
\sn
\itemitem{ ${{}}$ }  $(f) \quad \langle (M^1_i,M^2_i,\bold I_i):i \le
\beta \rangle$ is $\le_{\frak u}^1$-increasing continuous and
\nl

\hskip35pt $(M^1_i,M^2_i,\bold I_i) \in \text{ FR}^+_1$
\sn
\itemitem{ ${{}}$ }  $(g) \quad (M^1_i,M^1_{i+1},\bold J^1_i)
\le_{\frak u}^2
(M^2_i,M^2_{i+1},\bold J^2_i) \in \text{ FR}_2$ for $i < \beta$
\sn
\itemitem{ ${{}}$ }  $(h) \quad (M^1_0,M^2_0,\bold I_0) = (M,N,\bold I)$
\sn
\itemitem{ $\circledast_2$ }  $M^{2,1}_{\gamma_1},M^{2,2}_{\gamma_2}$ 
are $\tau$-incompatible amalgamation of $M^1_*,M^2_0$ over 
$M^1_0$ in ${\frak K}_{< \partial}$ and for $\ell=1,2$ we have
\sn
\itemitem{ $\circledast_{2,\ell}$ }  $(a)' \quad 
\beta < \gamma_\ell < \partial$
\sn
\itemitem{ ${{}}$ }  $(b)'_1 \quad M^1_\beta \le_{\frak u} M^1_*$
\sn
\itemitem{ ${{}}$ }  $(c)'_1 \quad (M^{2,\ell}_i,M^{2,\ell}_{i+1},
\bold J^{2,\ell}_i) \in \text{ FR}_2$ for $i < \gamma_\ell$
\sn
\itemitem{ ${{}}$ }  $(d)' \quad M^{2,\ell}_i$ (for $i \le \gamma_\ell$)
is $\le_{\frak u}$-increasing continuous
\sn
\itemitem{ ${{}}$ }  $(e)' \quad M^1_* \le_{\frak u} M^{2,\ell}_{\gamma_\ell}$
\sn
\itemitem{ ${{}}$ }  $(f)' \quad (M^1_\beta,M^2_\beta,\bold I_\beta)
\le^1_{\frak u} (M^1_*,M^{2,\ell}_{\gamma_\ell},\bold I^\ell_*)$.
\endroster}
\ermn
1A)  We say that $(M,N,\bold I)$ has the vertical $\tau$-coding$_0$ 
property when:  
\ub{if} $(\bar M^1,\bar{\bold J}^1) = (\langle M^1_i:i \le \beta
\rangle,\langle \bold J^1_i:i < \beta \rangle)$ satisfies clause (A) of part
(1) and $(\langle M^2_i:i \le \beta \rangle,\langle \bold J^2_i:i <
\beta \rangle,\langle \bold I_i:i \le \beta \rangle)$ satisfies
$\circledast_1$ of clause (B) of part (1) and
$(M^1_0,M^2_0,\bold I_0) = (M,N,\bold I)$ \ub{then} we can find objects
satisfying $\circledast_2$ of clause (B) of part (1).
\nl
1B)  We say that $(M,N,\bold I)$ has the true vertical $\tau$-coding$_0$
property \ub{when} it belongs to FR$^+_1$ and 
every $(M',N',\bold J')$ satisfying $(M,N,\bold I)
\le^1_{\frak u} (M',N',\bold I')$ has the vertical $\tau$-coding$_0$ property.
\nl
1C) We say that ${\frak u}$ has the explicit vertical
$\tau$-coding$_0$ property \ub{when} for every $M$ for some $N,\bold
I$ the triple $(M,N,\bold I) \in \text{ FR}^+_1$ has the true vertical $\tau$-coding$_0$ property.
\nl
2) $(\bar M^*,\bar{\bold J}^*,\bold f^*) \in K^{\text{qt}}_{\frak u}$
has the vertical $\tau$-coding$_0$ property \ub{when} for a club of
$\delta < \partial$, the pair $(\bar M^1,\bar{\bold J}^1) = (\langle
M^*_{\delta +i}:i \le \bold f^*(\delta)\rangle,\langle \bold
J^*_{\delta +1}:i < \bold f^*(\delta)\rangle)$ satisfies  part (1) even
demanding $M^1_* \le_{\frak K} M^*_\partial$.
\nl
3) We say that $(\bar M,\bar{\bold J},\bold f) \in
K^{\text{qt}}_{\frak u}$ has the vertical $\tau$-coding$_1$ property \ub{when}
(we may omit the subscript 1) we can find
$\alpha(0) < \partial$ and $(M_{\alpha(0)},N_*,\bold I_*) \in 
\text{ FR}_1$ satisfying $N_* \cap M_\partial = 
M_{\alpha(0)}$ such that: for a club of $\delta < \partial$ the pair 
$(\bar M^1,\bar{\bold J}^1) = (\langle M_{\delta +i}:i \le
\bold f(\delta) \rangle,\langle \bold J_{\delta +i}:i < \bold
f(\delta) \rangle)$ satisfies part (1) when in clause (B) where we
\mr
\item "{$(i)$}"   restrict ourselves to
the case $(M_{\alpha(0)},N_*,\bold I_*) \le^1_{\frak u} 
(M^1_0,M^2_0,\bold I_0)$
\sn
\item "{$(ii)$}"  demand that $M^1_* <_{{\frak K}[{\frak u}]} M_\partial$.
\ermn
4) We say that $(\bar M,\bar{\bold J},\bold f) \in
K^{\text{qt}}_{\frak u}$ has the $S$-vertical $\tau$-coding$_1$
property \ub{when}: $S$ is a stationary subset of $\partial$ and for
club $E$ of $\partial$ the requirement in part (3) holds when we restrict
ourselves to $\delta \in S \cap E$.
\nl
4A) We say that $(\bar M,\bar{\bold J},\bold f) \in
K^{\text{qt}}_{\frak u}$ has the $S$-vertical $\tau$-coding$_2$
property as in Definition \scite{838-2b.6F}.
\nl
5) For $k=0,1,2$ we say ${\frak u}$ has the vertical $\tau$-coding$_k$ 
property \ub{when} $\{0,2\}$-almost
every $(\bar M,\bar{\bold J},\bold f) \in K^{\text{qt}}_{\frak u}$ has
it.  If $k = 1$ we may omit it.  
Similarly adding ``above $(\bar M^*,\bar{\bold J}^*,\bold f^*)$" and/or
$S$-vertical, for stationary $S \subseteq \partial$.
\enddefinition
\bigskip

The following observation is easy but very useful.
\demo{\stag{838-2b.11} Observation}   1) Assume that some $(M,N,\bold I) \in
\text{ FR}^+_1$ has the true vertical $\tau$-coding$_0$ property 
(from Definition \scite{838-2b.9}(1B)).   If $(\bar M,\bar{\bold J},\bold f) 
\in K^{\text{qt}}_{\frak u}$ and $M_\partial := \cup\{M_\alpha:\alpha <
\partial\}$ is saturated (above $\lambda$, for ${\frak K} = 
{\frak K}^{\frak u}$) \ub{then} $(\bar
M,\bar{\bold J},\bold f)$ has the vertical $\tau$-coding property.  
See \scite{838-2b.9}(2),\scite{838-2b.9}(2A) used in \scite{838-e.2C}(5).
\nl
2) If $(\bar M,\bar{\bold J},\bold f)$ has the vertical $\tau$-coding$_0$ 
property \ub{then} it has the vertical $\tau$-coding$_1$ property.
\nl
3) Similarly (to part (2)) for weak $\tau$-coding.
\nl
4) Recalling ${\frak K}_{\frak u}$ has amalgamation (by claim \scite{838-1a.5}(1))
\mr
\item "{$(a)$}"  if $M \in {\frak K}_{\frak u} \Rightarrow |{\Cal
S}_{\frak K}(M)| \le \partial$ \ub{then} there is a saturated $M \in
K^{\frak u}_\partial$
\sn
\item "{$(b)$}"  if every $M \in K^{{\frak u},*}_\partial$ is
saturated and every $(M,N,\bold I) \in \text{ FR}^+_1$ has weak 
$\tau$-coding$_0$, \ub{then} every $(\bar M,\bar{\bold J},\bold f) \in
K^{\text{qt}}_{\frak u}$ has weak $\tau$-coding
\sn
\item "{$(c)$}"  similarly for vertical $\tau$-coding
\sn
\item "{$(d)$}"  similarly replacing ``every 
$M \in K^{{\frak u},*}_\partial$" by ``$M_\partial$ is saturated 
for $\{0,2\}$-almost every 
$(\bar M,\bar{\bold J},\bold f) \in K^{\text{qt}}_{\frak
u}$" [or just above $(\bar M^*,\bar{\bold J}^*,\bold f^*) \in
K^{\text{qt}}_{\frak u}$.]
\endroster     
\enddemo
\bigskip

\demo{Proof}  Should be clear.  \hfill$\square_{\scite{838-2b.8}}$
\enddemo
\bigskip

\proclaim{\stag{838-2b.13} Theorem}  We have
$\dot I_\tau(\partial^+,K^{\frak u}_{\partial^+}) \ge
\mu_{\text{unif}}(\partial^+,2^\partial)$; moreover 
$\dot I(K^{{\frak u},{\frak h}}_{\partial^+}) \ge
\mu_{\text{unif}}(\partial^+,2^\partial)$ for 
any $\{0,2\}$-appropriate ${\frak h}$ (see Definitions
\scite{838-1a.47}(2),(7), \scite{838-1a.45}) \ub{when}:
\mr
\item "{$(a)$}"  $2^\theta = 2^{< \partial} < 2^\partial$ 
\sn
\item "{$(b)$}"  $2^\partial < 2^{\partial^+}$ 
\sn
\item "{$(c)$}"  ${\frak u}$ has the vertical $\tau$-coding$_1$
property (or at least ${\frak u}$ has the $S$-vertical
$\tau$-coding$_\tau$ property above
some triple from $K^{\text{qt}}_{\frak u}$ for stationary $S \in
(\text{\rm WDmId}(\partial))^+$ (recall $\tau$ is a weak
${\frak u}$-sub-vocabulary, of course, by \scite{838-2b.0})).
\endroster
\endproclaim
\bigskip

\remark{Remark}  Theorem \scite{838-2b.13} is used in \scite{838-e.2} and in
\scite{838-8h.12}.
\endremark
\bigskip

\demo{Proof}  Proved in \scite{838-10k.19}.  \hfill$\square_{\scite{838-2b.13}}$
\enddemo
\bn
\centerline{$* \qquad * \qquad *$}
\bn
\margintag{2b.14}\ub{\stag{838-2b.14} Discussion}:  1) In a 
sense the following property ``horizontal 
$\tau$-coding" is dual to the previous one ``vertical $\tau$-coding",
it is ``horizontal", i.e. in the $\partial^+$-direction.  This will
result in building $(\bar M^\eta,\bar{\bold J}^n,\bar{\bold f}^\eta)$
for $\eta \in {}^{\partial^+ >}(2^{< \partial})$ such that letting
$M^\eta_\partial = \cup\{M^{\eta \restriction \alpha}_\partial:
\alpha < \partial^+\}$ we have
$\eta \ne \nu \in {}^{\partial^+}(2^\partial) \Rightarrow 
M^\eta_\partial,M^\nu_\partial$ are not isomorphic over 
$M^{<>}_\partial$, so the set theory is simpler.
\nl
2) Note that in \scite{838-2b.15}(4) below we could ask less than ``for a club",
e.g. having a winning strategy is the natural game; similarly in other
definitions of coding properties, as in Exercise \scite{838-2b.6}.
\bigskip

\definition{\stag{838-2b.15} Definition}  1) We say that $(M_0,M_1,\bold J_2)
\in \text{ FR}_2$ has the horizontal $\tau$-coding$_0$ property
\ub{when}: if $\circledast_1$ holds
\ub{then} we can find $N^{+,\ell}_1,\bold I^\ell_5,\bold J^\ell_5$ 
for $\ell=1,2$ such that $\circledast_2$ holds \ub{when}:
\mr
\item "{$\circledast_1$}"  $(a) \quad 
M_0 \le_{\frak u} N_0 \le_{\frak u} N_1,M_0 \le_{\frak u} M_1
\le_{\frak u} N_1$
\sn
\item "{${{}}$}"  $(b) \quad (M_0,M_1,\bold J_2) 
\le_{\frak u}^2 (N_0,N_1,\bold J_3)$ so both are from FR$_2$
\sn
\item "{${{}}$}"  $(c) \quad (N_0,N^+_0,\bold I_4) \in \text{ FR}^+_1$ and $N^+
\cap N_1 = N_0$
\sn
\item "{$\circledast_2$}"  $(\alpha) \quad (N_0,N^+_0,\bold I_4) 
\le^1_{\frak u} (N_1,N^{+,\ell}_1,\bold I^\ell_5)$
\sn
\item "{${{}}$}"  $(\beta) \quad (N_0,N_1,\bold J_3) 
\le^2_{\frak u} (N^+_0,N^{+,\ell}_1,\bold J^\ell_5)$
\sn
\item "{${{}}$}"  $(\gamma) \quad N^{+,1}_1,N^{+,2}_1$ are 
$\tau$-incompatible amalgamations of $M_1,N^+_0$ over 
\nl

\hskip25pt $M_0$ in $K_{\frak u}$.
\ermn
2) We say that $(\bar M,\bar{\bold J},\bold f) \in 
K^{\text{qt}}_{\frak u}$ has the $S$-horizontal 
$\tau$-coding$_0$ property \ub{when} $S$ is a stationary subset of
$\partial$ and for a club of $\delta \in S$, the triple
$(M_\delta,M_{\delta +1},\bold J_\delta)$ has it and $\bold f(\delta)
> 0$.  If $S = \partial$ we may omit it.
\nl
3) We say that $(\bar M,\bar{\bold J},\bold f) \in
K^{\text{qt}}_{\frak u}$ has the horizontal 
$\tau$-coding$_1$ property \ub{when} (we may omit the 1): 
\mr
\item "{$\circledast$}"  for $\{0,2\}$-almost every $(\bar
M',\bar{\bold J}',\bold f') \in K^{\text{qt}}_{\frak u}$ satisfying 
$(\bar M,\bar{\bold J},\bold f) 
\le^{\text{qs}}_{\frak u} (\bar M',\bar{\bold J}',\bold f')$ 
we can find $\alpha < \partial$ and $N,\bold I^*$ satisfying
$(M'_\alpha,N,\bold I^*) \in
\text{ FR}^+_2$ and $N \cap M'_\partial = M'_\alpha$ such that
{\roster
\itemitem{ $\boxdot$ }  for a club of $\delta < \partial$ if
$(M'_\alpha,N,\bold I^*) \le_1 (M'_\delta,N',\bold I') \in \text{
FR}^+_1$ and $N' \cap M'_\partial = M'_\alpha$, \ub{then} the
conclusion in \scite{838-2b.15}(1) above holds with $M_{\alpha,i(*)},
M'_{i(*)}$,
\nl

\hskip20pt $M_\delta,M'_\delta,N',\bold I',\bold J_\delta,\bold
J'_\delta$ here standing for $M_0,N_0,M_1,N_1,N^+_0,\bold I_4,\bold
J_2,\bold J_3$
\nl

\hskip20pt  there.
\endroster}
\ermn
4) We replace coding$_1$ by coding$_2$ when in $\boxdot$ we use the
game version, as in \scite{838-2b.6F}.
\nl
5) We say ${\frak u}$ has the horizontal $\tau$-coding$_k$ property
\ub{when} some $(\bar M,\bar{\bold J},\bold f) \in
K^{\text{qt}}_{\frak u}$ has it.
\enddefinition
\bigskip

\proclaim{\stag{838-2b.17} Claim}  The coding$_0$ implies the coding$_1$
versions in Definition \scite{838-2b.15} for $(\bar M,\bar{\bold J},\bold f)
\in K^{\text{qt}}_{\frak u}$ in Definition \scite{838-2b.15}(4) and for 
${\frak u}$ in Definition \scite{838-2b.15}(5). 
\endproclaim
\bigskip

\demo{Proof}  Should be clear.  \hfill$\square_{\scite{838-2b.17}}$
\enddemo
\bigskip

\proclaim{\stag{838-2b.19} Theorem}  We have
$\dot I_\tau(\partial^+,K^{\frak u}_{\partial^+}) \ge 2^{\partial^+}$;
moreover $\dot I(K^{{\frak u},{\frak h}}_{\partial^+}) \ge
2^{\partial^+}$ for any $\{0,2\}$-appropriate ${\frak h}$ (see Definition
\scite{838-1a.47}(4),(5),(6)), \ub{when}: 
\mr
\item "{$(a)$}"  $2^\theta = 2^{< \partial} < 2^\partial$ 
\sn
\item "{$(b)$}"  $2^\partial < 2^{\partial^+}$
\sn
\item "{$(c)$}"  the ideal {\rm WDmId}$(\partial)$ is not 
$\partial^+$-saturated
\sn
\item "{$(d)$}"  ${\frak u}$ has the horizontal $\tau$-coding property
(or just the $S$-horizontal $\tau$-coding$_2$ property for some
stationary $S \subseteq \partial$).
\endroster
\endproclaim
\bigskip

\remark{\stag{838-c.10.2} Remark}  1) Actually not used here.
\nl
2) What does this add compared to
\scite{838-2b.3}? getting $\ge 2^{\partial^+}$ rather than
$\ge \mu_{\text{unif}}(\partial^+,2^\partial)$. 
\endremark
\bigskip

\demo{Proof}  Proved in \scite{838-10k.21}.  \hfill$\square_{\scite{838-2b.19}}$
\enddemo
\goodbreak

\head {\S3 Invariant codings} \endhead  \resetall \sectno=3
 \spuriousreset
\bigskip

The major notion of this section is (variants of) uq-invariant coding
properties.  In our context, the point of coding properties is in
esssence that their  failure gives that there are many of uniqueness
triples, $(M,N,\bold J)$ ones, i.e. such that: if $(M,N,\bold I) \le_1
(M',N'_\ell,\bold I'')$ for $\ell=1,2$ then $N'_1,N'_2$ are compatible
over $N \cup M'$.  For uq-invariant we ask for less: if $(M,N,\bold I)
\le (M',N',\bold I')$ and $\bold d$ is a ${\frak u}$-free
$(\xi,0)$-rectanle with $(M^{\bold d}_{0,0},M^{\bold d}_{\alpha(\bold
d),0}) = (M,M')$, \ub{then} we can ``lift" $\bold d''$, i.e. find a
${\frak u}$-free $(\xi +1,1)$-rectangle $\bold d^+$ such that $\bold
d^+ \rest (\xi,0) = \bold d,(M^{\bold d}_{0,0},M^{\bold d}_{0,1},\bold
I^{\bold d}_{0,0}) = (M,N,\bold I)$ and $N' \le_{\frak u} M^{\bold
d}_{\xi +1,1}$.

So we look at the simplest version, the weak $\xi$-uq-invariant coding,
Definition \scite{838-3r.1}, we can consider a ``candidate" $(M,N,\bold I)
\in \text{ FR}^1_{\frak u}$ and challenge $\bold d$ (so $M^{\bold
d}_{0,0} = M$) and looks for a pair of amalgamation which are
incompatile in a specific way, \ub{but} unlike in \S2, they are not
symmetric.  One is really not an amalgamation \ub{but} a family of
those exhibiting $\bold d$ is ``liftable", and ``promise to continue
to do so in the future", in the $\partial^+$-direction. The real one
just has to contradict it.  

Another feature is that instead of considering isomorphisms over
$M^{\bold d}_{\alpha(\bold d),0} \cup N$ we consider isomorphisms over
$M^{\bold d}_{\alpha(\bold d),0}$ with some remnants of preserving
$N$; more specifically we consider two ${\frak u}$-free rectangles
${\frak d}_1,{\frak d}_2$ which continues the construction in those two
ways and demands $M^{\bold d_1}_{\alpha(\bold d_1)}$ is mapped onto
$M^{\bold d_2}_{0,\alpha(\bold d_2)}$.

There are more complicating factors: we have for a candidate
$(M,N,\bold I)$, for every $M',M \le_{\frak u} M'$ to find a ${\frak
u}$-free $(\xi,0)$-rectangle $\bold d$ with $M' = M^{\bold d}_{0,0}$
such that it will serve against $(M',N',\bold I')$ whenever
$(M,N,\bold I) \le_1 (M',N',\bold I')$, rather than choosing $\bold d$
after $(N',\bold I')$ is chosen, i.e. this stronger version is needed.
 The case $\xi < \partial$ should be clear but still we allow $\xi = 
\partial$, however 
then given $(N',\bold I')$ we take $\bold d \rest (\xi',0)$ for some
$\xi' < \xi'$.

We can use only Definition \scite{838-3r.1}, Claim \scite{838-3r.3}, and
Conclusion \scite{838-3r.7}, for which ``$2^\theta = 2^{< \partial} <
2^\partial < 2^{\partial^+} +$ the extra WdmId$(\partial)$ is not
$\partial^+$-saturated" is needed, (if $\xi < \partial$ less is
needed; however $\xi < \partial$ shall not be enough.  
But to get the sharp results (with the extra
assumption) for almost good $\lambda$-frames we need a more elaborate
approach - using \ub{vertical} $\xi$-uq-invariant coding, see
Definition \scite{838-3r.19}.

Actually we shall use an apparently weaker version, the so called semi
$\xi$-uq-invariant.  However, we can derive from it the vertical
version under reasonable demands on ${\frak u}$; this last proof is of
purely model theoretic characters.  We also consider other variants.
\bn
In this section we usually do not use the $\tau$ from \scite{838-3r.0},
i.e. use $\tau = \tau_{\frak u}$ as it is
not required presently.
\demo{\stag{838-3r.0} Hypothesis}  We assume ${\frak u}$ is a nice
construction framework and $\tau$ is a weak ${\frak u}$-sub-vocabulary.
\enddemo
\bigskip

\remark{Remark}  In Definition \scite{838-3r.1}(1) below ``$\bold d_*$
witnesses not being able to lift $\bold d$", of course we can ensure
it can be lifted.
\endremark
\bigskip

\definition{\stag{838-3r.1} Definition}  Let $\xi \le \partial +1$, if we
omit it we mean $\xi = \partial +1$.
\nl
1) We say that $(M,N,\bold I) 
\in \text{ FR}^1_{\frak u}$ has weak $\xi$-uq-invariant coding$_0$ \ub{when}:
\mr
\item "{$\circledast$}"  if $M \le_{\frak u} M'$ and $M' \cap N = M$ \ub{then}
there are an ordinal $\alpha < \xi$ and a
${\frak u}$-free $(\alpha,0)$-rectangle $\bold d$ so $\alpha 
= \partial$ is O.K., such that:
{\roster
\itemitem{ $(a)$ }  $M^{\bold d}_{0,0} = M'$ and $M^{\bold
d}_{\alpha,0} \cap N = M$
\sn
\itemitem{ $(b)$ }  for every $N',\bold I'$ such that $(M',N',\bold
I')$ is $\le^1_{\frak u}$-above $(M,N,\bold I)$ and $N' \cap 
M^{\bold d}_{\alpha,0} = M'$ we can find $\alpha',\bold I^1,\alpha_*$ 
and $\bold d_*$ such that $\alpha' \le \alpha,\alpha' < \partial$
(no harm\footnote{by natural monotonicity, similarly in \scite{838-3r.11},
\scite{838-3r.19}, \scite{838-3r.71}} in $\alpha < \partial \Rightarrow
\alpha' = \alpha$) and:
\endroster}
\item "{$(\alpha)$}"  $\bold d_*$ is a ${\frak u}$-free
$(\alpha_*,0)$-rectangle and $\alpha_* < \partial$
\sn
\item "{$(\beta)$}"  $M^{\bold d_*}_{0,0} = N'$ and 
$M^{\bold d}_{\alpha',0} \le_{\frak u} M^{\bold d_*}_{\alpha_*,0}$
\sn
\item "{$(\gamma)$}"  $(M',N',\bold I') \le^1_{\frak u} (M^{\bold
d}_{\alpha,0},M^{\bold d_*}_{\alpha_*,0},\bold I^1)$
\sn
\item "{$(\delta)$}"  there are no $\bold d_1,\bold d_2$ such that
{\roster
\itemitem{ $\bullet_1$ }  $\bold d_\ell$ is a ${\frak u}$-free rectangle
for $\ell=1,2$
\sn
\itemitem{ $\bullet_2$ }  $\alpha(\bold d_1) = \alpha_*$ and $\bold d_1
\restriction (\alpha_*,0) = \bold d_*$
\sn
\itemitem{ $\bullet_3$ }  $\alpha(\bold d_2) \ge \alpha'$ and $\bold d_2
\restriction (\alpha',0) = \bold d \rest (\alpha',0)$  
\sn
\itemitem{ $\bullet_4$ }  $(M^{\bold d_2}_{0,0},M^{\bold d_2}_{0,1},
\bold I^{\bold d_2}_{0,0}) = (M',N',\bold I')$ or just $(M,N,\bold I)
\le^1_{\frak u} (M^{\bold d_2}_{0,0},M^{\bold d_2}_{0,1},I^{\bold d_2}_{0,0})$
\sn
\itemitem{ $\bullet_5$ }   there are $N \in {\frak K}_{\frak u}$ and
$\le_{\frak u}$-embeddings $f_\ell$ of 
$N^{\bold d_\ell}_{\alpha(\bold d_\ell),\beta(\bold d_\ell)}$ into $N$ for
$\ell=1,2$ such that $f_1 \restriction M^{\bold d}_{\alpha',0} 
= f_2 \restriction M^{\bold d}_{\alpha',0}$ and
$f_1,f_2$ maps $M^{\bold d_1}_{0,\beta(\bold d_1)}$ onto 
$M^{\bold d_2}_{0,\beta(\bold d_2)}$.
\endroster}
\ermn
2) We say that $(\bar M^*,\bar{\bold J}^*,\bold f^*) \in
K^{\text{qt}}_{\frak u}$ has the weak $\xi$-uq-invariant 
coding$_0$ property \ub{when}: if $\alpha(0) < \partial$ and
$(M_*,N_0,\bold I_0) \in \text{ FR}^+_1,M_* \le_{\frak u}
M^*_{\alpha(0)}$ and $N_0 \cap M^*_\partial = \emptyset$
\ub{then} for some club $E$ of $\partial$, for every $\delta \in E$
the statement $\circledast$ of part (1) holds, with $(M_*,N_0,\bold
I_0),M_\delta$ here standing for $(M,N,\bold I),M'$ there
but with some changes:
\mr
\item "{$(*)_1$}"   $\bold d$ is such that $M^{\bold d}_{\alpha',0} 
\le_{\frak u} M^*_\beta$ for each $\alpha' < \alpha$
 for any $\beta < \partial$ large enough and
\sn
\item "{$(*)_2$}"  in clause (b) we demand $N' \cap M^*_\partial=
M_\delta$ and $M^{\bold d_*}_{\alpha_*} \le_{\frak u} M^*_\beta$ for
every $\beta < \partial$ large enough.
\ermn
3)  We say that $(\bar M^*,\bar{\bold J}^*,\bold f) \in
K^{\text{qt}}_{\frak u}$ has the weak $\xi$-uq-invariant coding$_1$
property as in part (2) but require only that there are such
$\alpha(0),(M_*,N_0,\bold I_0)$ so without loss of generality 
$M_* = M^*_{\alpha(0)}$.
\nl
3A) We say that $(\bar M,\bar{\bold J},\bold f)$ has the $S$-weak 
$\xi$-uq-invariant coding$_2$ property \ub{when} we combine the above
with Definition \scite{838-2b.6F}.
\nl
4) For $k=0,1$ we say that $(\bar M,\bar{\bold J},\bold f)$ has the $S$-weak 
$\xi$-uq-invariant
coding$_k$ property \ub{when}: $S$ is a stationary subset of
$\partial$ and for some club $E$ of $\partial$ the demand in part (2)
if $k=0$, part (3) if $k=1$ holds restricting ourselves to $\delta \in 
S \cap E$.
\nl
5) We say that ${\frak u}$ has the $S$-weak $\xi$-uq-invariant 
coding$_k$ property when: $\{0,2\}$-almost 
every $(\bar M,\bar{\bold J},\bold f) \in 
K^{\text{qt}}_{\frak u}$ has the $S$-weak $\xi$-uq-invariant coding$_k$ property.  
Similarly for ``above $(\bar M^*,\bar{\bold J}^*,\bold f^*) 
\in K^{\text{qt}}_{\frak u}$".  If $S
= \partial$ we may omit it; if $k=1$ we may omit it.
\enddefinition
\bigskip

\proclaim{\stag{838-3r.3} Claim}  Assume ($\xi \le \partial +1$ and):
\mr
\item "{$(a)$}"  $S \subseteq \partial$ is stationary
\sn
\item "{$(b)$}"  $(\bar M,\bar{\bold J},\bold f) 
\in K^{\text{qt}}_{\frak u}$ has the $S$-weak $\xi$-uq-invariant 
coding property 
\sn
\item "{$(c)$}"  $\bold f \restriction S$ is constantly zero
\sn
\item "{$(d)$}"  $\xi \le \partial +1$ and if $\xi = \partial$ then
the ideal {\rm WDmId}$(\partial) + (\lambda \backslash S)$ is not
$\partial$-saturated.
\ermn
\ub{Then} we can find $\langle (\bar N^\eta,\bar{\bold J}^\eta,
\bold f^\eta):\eta \in {}^\partial 2 \rangle$ such that
\mr
\item "{$(\alpha)$}"  $(\bar M,\bar{\bold J},\bold f) 
\le^{\text{qs}}_{\frak u} (\bar N^\eta,\bar{\bold J}^\eta,\bold f^\eta)$
\sn
\item "{$(\beta)$}"  $\bold f^\eta (\partial \backslash S) = \bold f
\restriction (\partial \backslash S)$
\sn
\item "{$(\gamma)$}"  if $\eta^1 \ne \eta^2 \in {}^\partial 2$ and
$(\bar N^{\eta^\ell},\bar{\bold J}^{\eta^\ell},\bold f^{\eta^\ell}) 
\le^{\text{qt}}_{\frak u} (\bar N^\ell,\bar{\bold J}^\ell,\bold f^\ell)$ for 
$\ell=1,2$ \ub{then} $N^1_\partial,N^2_\partial$ are not isomorphic
over $M_\partial$.
\endroster
\endproclaim
\bigskip

\remark{\stag{838-3r.5} Remark}  Note that in Definition \scite{838-3r.1}(1) we 
choose the ${\frak u}$-free
$(\alpha_\delta,0)$-rectangle $\bold d_\delta$ for every $\delta \in S$
before we have arrived to choosing $N^\eta_\delta$.  This will be a
burden in applying this.
\endremark
\bigskip

\demo{Proof}  For simplicity we first assume $\xi \le \partial$.
Let $\langle S_\varepsilon:\varepsilon \le \partial
\rangle$ be a sequence of pairwise disjoint stationary subsets of
$\partial$ with union $S \backslash \{0\}$ such that $\varepsilon < \text{
min}(S_\varepsilon)$ (exists; if $\xi = \partial$ by assumption (d), otherwise 
if $\partial$ successor by applying Ulam matrixes, in
general by a theorem of Solovay).  Without loss of generality $0
\notin S$.

By assumption (b) we can find $\alpha(0) < \partial$ and $N_0,\bold I_0$
such that $N_0 \cap M_\partial = M_{\alpha(0)}$ and
$(M_{\alpha(0)},N_0,\bold I) \in \text{ FR}_1$ and a club
$E_0$ of $\partial$ such that and for every $\delta \in S \cap E_0$ we can
choose $\alpha_\delta,\bold d_\delta$ as in 
$\circledast$ of Definition \scite{838-3r.1}(1) with
$M_{\alpha(0)},N_0,\bold I_0,M_\delta$ here standing for $M,N,\bold
I,M'$ there but demanding $M^{\bold d_\delta}_{\alpha_\delta,0} 
\le_{\frak u} M_{\beta_\delta}$
for some $\beta_\delta \in (\delta,\partial)$, see $(*)_1$ of
Definition \scite{838-3r.1}(2).  Without loss of generality $\alpha(0)$ is
a successor ordinal.

Let

$$
\align
E_1 = \{\delta \in E_0:&\delta > \alpha(0) \text{ and if } \delta(1)
\in \delta \cap S \cap E_0 \\
  &\text{then } \beta_{\delta(1)} \le \delta, \text{ i.e. }
M^{\bold d_{\delta(1)}}_{\alpha_{\delta(1)},0} \le_{\frak u} M_\delta \\
  &\text{ and } (\forall \beta <  
\delta)(\beta \times \beta < \delta \wedge \bold f(\beta) < \delta)\}.
\endalign
$$
\mn
Clearly $E_1$ is a club of $\partial$.

We now choose $\langle (\delta^\rho,u^\rho,\bar N^\rho,\bar M^\rho,
\bar{\bold J}^{0,\rho},\bar{\bold J}^{1,\rho},\bold f^\rho,
\bold I^\rho,e^\rho):
\rho \in {}^i 2 \rangle$ by induction on $i < \partial$ such that
\mr
\item "{$\circledast$}"  $(a) \quad (\alpha) \quad 
\delta^\rho < \partial$ belongs to
$E_1 \cup \{\alpha(0)\}$
\sn
\item "{${{}}$}"  $\qquad (\beta) \quad e^\rho$ is a closed 
subset of $E_1 \cap \delta^\rho$ 
\sn
\item "{${{}}$}"  $\qquad (\gamma) \quad \text{ min}(e^\rho) = \alpha(0)$,
also $N^\rho_{\alpha(0)} = N_0,\bold I^\rho_{\alpha(0)} = \bold I$
\sn
\item "{${{}}$}"  $(b) \quad \bold f^\rho:e^\rho \rightarrow 
\delta^\rho$
\sn 
\item "{${{}}$}"  $(c) \quad \bold f^\rho(\alpha) < \beta$ if $\alpha
< \beta$ are from $e^\rho$ 
\sn
\item "{${{}}$}"  $(d) \quad$ if $\alpha \in e^\rho \backslash S$ then
$\bold f^\rho(\alpha) = \bold f(\alpha)$
\sn
\item "{${{}}$}"  $(e) \quad (\alpha) \quad 
\bar N^\rho = \langle N^\rho_i:i \le \delta^\rho \rangle$ and $\bar
M^\rho = \langle M^\rho_i:i \le \delta^\rho\rangle$ are
\nl

\hskip25pt  $\le_{\frak u}$-increasing continuous
\sn
\item "{${{}}$}"  $\quad \quad (\beta) \quad
 \langle(M^\rho_i,N^\rho_i,\bold I^{\rho \rest i}):i \le \delta\rangle$ is
 $\le^1_{\frak u}$-increasing continuous
\sn
\item "{${{}}$}"  $(f) \quad (\alpha) \quad 
\bar{\bold J}^{\ell,\rho} = \langle \bold J^{\ell,\rho}_i:
i < \delta^\rho \rangle$ 
\sn
\item "{${{}}$}"  $\quad \quad (\beta) \quad 
(M^\rho_i,M^\rho_{i+1},\bold J^{0,\rho}_i) \in 
\text{FR}_2$ for $i < \delta^\rho$
\sn
\item "{${{}}$}"  $\quad \quad (\gamma) \quad
(M^\rho_i,M^\rho_{i+1,\rho},\bold J^{0,\rho}_i) \le^2_{\frak u} 
(N^\rho_i,N^\rho_{i+1},\bold J^{1,\rho}_i)$ when $i < \delta^\rho
\and$
\nl

\hskip35pt $(\exists \alpha)(\alpha \in e^\rho \cap S \and 
\alpha \le i \le \alpha + \bold f^\rho(\alpha))$
\sn
\item "{${{}}$}"  $(g) \quad$ if $\alpha \in e^\rho$ \ub{then}
$M^\rho_\alpha = M_\alpha$ and $N_\alpha \cap M_\partial = M_\alpha$,
\sn
\item "{${{}}$}"  $(h) \quad$ if $\alpha \in e^\rho \backslash S$ 
then ($\bold f^\rho(\alpha) = \bold f(\alpha)$ and) $i \le \bold
f(\alpha) \Rightarrow M^\rho_{\alpha +i} =$
\nl

\hskip25pt $M_{\alpha +i} \wedge N^\rho_{\alpha +i} \cap M_\partial =
M_{\alpha +i} = M^\rho_{\alpha +i}$ and $i < \bold f(\alpha) 
\Rightarrow \bold J^{0,\rho}_{\alpha +i} = \bold J_{\alpha +i}$
\nl

\hskip25pt   and $i < \bold f(\alpha) \Rightarrow 
(M^\rho_{\alpha +i},M^\rho_{\alpha +i+1},\bold J^{0,\rho}_{\alpha +i}) 
\le^2_{\frak u} (N^\rho_{\alpha +i},N^\rho_{\alpha +i+1},
\bold J^{1,\rho}_{\alpha +i})$
\sn
\item "{${{}}$}"  $(i) \quad$ if $\varrho \triangleleft \rho$ then
$\delta_\varrho < \delta_\rho,e^\varrho = \delta_\varrho \cap
e^\rho,\bar M^\varrho \triangleleft \bar M^\rho,
\bar N^\varrho \triangleleft \bar N^\rho,\bar{\bold J}^{\ell,\varrho} 
\triangleleft \bar{\bold J}^{\ell,\rho}$ 
\nl

\hskip25pt for $\ell=0,1$ and $\bar{\bold I}^\varrho \triangleleft
\bar{\bold I}^\rho$
\sn
\item "{${{}}$}"  $(j) \quad$ if $\varepsilon < \partial$ and
$\alpha \in e^\rho \cap S$ and $\rho(\alpha) = 1$ then 
\nl

\hskip25pt $(\langle M^\rho_{\alpha +i}:i \le 
\bold f^\rho(\alpha)\rangle,\langle \bold J^{0,\rho}_{\alpha +i}:i <
\bold f^\rho(\alpha)\rangle)$ is equal to $\bold d_\alpha$ 
\sn
\item "{${{}}$}"  $(k) \quad$ if $\varepsilon < \partial,\alpha \in
e^\rho \cap S$ and $\rho(\alpha) = 0$ \ub{then}
\nl

\hskip25pt there is $\bold d_*$ as in clause (b) of 
$\circledast$ of Definition \scite{838-3r.1}(1), for transparency
\nl

\hskip25pt  of successor length $\beta$ 
with $(M_{\alpha(0)},N_0,\bold I_0),(M^\rho_\alpha,N^\rho_\alpha,
\bold I^\rho_\alpha),\bold d_\alpha,\bold d_*$
\nl

\hskip25pt  here standing for $(M,N,\bold I),(M',N',\bold I'),
\bold d,\bold d_*$ there, such that
\nl

\hskip25pt $\bold f^\rho(\alpha) = \alpha',N^\rho_{\alpha +i,0} =
M^{\bold d_*}_i$ for $i \le \bold f^\rho(\alpha)$,
\nl

\hskip25pt $\bold J^{1,\rho}_{\alpha +i} = \bold J^{\bold d_*}_{i,0}$
for $i < \bold f^\rho(\alpha)$, and (for transparency) 
\nl

\hskip25pt $M_{\alpha + i} = M_\alpha,\bold
J^{0,\rho}_{\alpha +i} = \emptyset$ for $i \le \bold f^\rho(\alpha)$.
\ermn
Why can we construct?
\mn
\ub{Case 1}:  \ub{$i=0$}

For $\rho \in {}^i 2$ let $\delta_\rho  = \alpha(0)$ and by
clause $(a)(\gamma)$ we define the rest (well for $i < \alpha(0)$ and
$\ell=0,1$ let
$M^\ell_i = M_{\alpha(0)},N^\rho_i = N_0,\bold J^{\ell,\rho}_i =
\emptyset,\bold I^\rho_i = \bold I$).
\mn
\ub{Case 2}:  \ub{$i$ is a limit ordinal}

For $\rho \in {}^i2$, let $\delta_\rho = \cup\{\delta_{\rho \rest j}:j
< i\}$ so $\delta_\rho \in E_1,\delta_\rho = \sup(E_1 \cap \alpha)$.

By continuity we can define also the others.
\mn
\ub{Case 3}:  \ub{$i = j+1$}

Let $\rho \in {}^j 2$ and we define for $\rho \char 94 \langle \ell
\rangle$ for $\ell=0,1$ and first we deal only with 
$i \le \delta + \bold f^{\rho \char 94 <\ell>}(\delta)$.
\bn
\ub{Subcase 3A}:  $\delta_\rho \notin S$

We use clause (h) of $\circledast$ and \scite{838-1a.9}(5).
\mn
\ub{Subcase 3B}:  $\delta \in S$

If $\ell=1$ then we 
use $\bold d_\alpha$ for $\rho \char 94 \langle \ell\rangle$ as in
clause (j) of $\circledast$ so the proof is as in subcase 3A.  
If $\ell=0$ clause (b) of
$\circledast$ of Definition \scite{838-3r.1} can be applied with
$(M,N,\bold I),(M',N',\bold I'),\bold d$ there standing for
$(M_{\alpha(0)},N_0,\bold I_0),(M^\rho_{\delta_\rho},N^\rho_{\delta_\rho},\bold
I^\rho_{\delta_\rho}),\bold d_{\delta_\rho}$ here; so we can find $\alpha'_\rho,
\bold I^1_\rho,\alpha^\rho_*,\bold d^\rho_*$ as there (presently
$\alpha'_\rho$ there can be $\alpha_\delta$); and without loss of generality  $\ell
g(\bold d^\rho_*)$ is a successor ordinal.

Now we choose:
\mr
\item "{$(*)$}"  $(a) \quad \bold f^{\rho \char 94 <0>}(\delta_\rho) = \ell
g(\bold d^\rho_*)$
\sn
\item "{${{}}$}"  $(b) \quad M^{\rho \char 94 <0>}_{\delta_\rho +i} =
M^\rho_{\delta_\rho}$ for $i < \bold f^{\rho \char 94 <0>}(\delta)$
\sn
\item "{${{}}$}"  $(c) \quad M^{\rho \char 94 <0>}_{\delta_\rho +i} = 
M^{\bold d_\delta}_{\alpha'_\rho,0}$ for $i = \bold f^{\rho \char 94
<0>}(\delta_\rho)$
\sn
\item "{${{}}$}"  $(d) \quad N^{\rho \char 94 <0>}_{\delta_\rho +i} = 
N^{\bold d^\rho_*}_{i,0}$ for $i \le \bold f^{\rho \char 94 <0>}
(\delta_\rho)$
\sn
\item "{${{}}$}"  $(e) \quad \bold J^{0,\rho \char 94 <\ell>}_{\delta_\rho
+i} = \emptyset,\bold J^{1,\rho \char 94 <\ell>}_{\delta_\rho +i} = 
\bold J^{\bold d^\rho_*}_{i,0}$ for $i < \bold f^{\rho \char 94
<0>}(\delta_\rho)$.
\ermn
Clearly clause (k) holds.  This ends the division to cases 3A,3B.

Lastly, choose $\delta_{\rho \char 94 <\ell>} \in E_1$ large enough; 
we still have to choose $(M^{\rho \char 94 <\ell>}_i,N^{\rho \char 94
<\ell>}_i,\bold I^{\rho \char 94 <\ell>})$ for 
$i \in (\delta_\rho + \bold f^\rho(\delta_\rho),
\delta_{\rho \char 94<\ell>}]$; we choose them all equal, 
$M^{\rho \char 94 <\ell>}_i =
M_{\delta_{\rho \char 94 <\ell>}}$ and use \scite{838-1a.5}(2) to choose 
$N^{\rho \char 94 <\ell>}_i,\bold I^{\rho \char 94 <\ell>}_i$.  Then
let $\bold J^{m,\rho \char 94 <\ell>}_i = \emptyset$ when $m < 2 \and
i \in [\delta_\rho + \bold f^\rho(\delta_\rho),\delta_{\rho \char 94 <\ell>})$.

So we have carried the induction.  For $\rho \in {}^\partial 2$ we
define $(\bar M^\rho,\bar{\bold J}^\rho,\bold f^\rho) \in
K^{\text{qt}}_{\frak u}$ by $M^\rho_\alpha = N^{\rho \restriction
i}_\alpha,\bold J^\rho_\alpha = \bar{\bold J}^{1,\rho \restriction
i}_\alpha,\bold f^\rho(\alpha) = \bold f^{\rho \restriction i}(\alpha)$ for
every $i < \partial$ large enough.  Easily
\mr
\item "{$\odot_1$}"  $(\bar M,\bar{\bold J},\bold f)
\le^{\text{qt}}_{\frak u}
(\bar M^\rho,\bar{\bold J}^\rho,\bold f^\rho)$.
\ermn
Next
\mr
\item "{$\odot_2$}"   for $\nu \in {}^\partial 2$ let $\rho_\nu = \rho[\nu]
\in {}^\partial 2$ be defined by $\rho_\nu(i) = \nu(\varepsilon)$ if $i \in
S_\varepsilon \wedge \varepsilon < \partial$ and zero otherwise.  
\ermn
So it is enough to prove
\mr
\item "{$\odot_3$}"  if $\nu_1 \ne \nu_2 \in {}^\partial 2$ and 
$(\bar M^{\rho[\nu_\ell]},\bar{\bold J}^{\rho[\nu_\ell]},
\bold f^{\rho[\nu_\ell]}) \le^{\text{qt}}_{\frak u} 
(\bar M^\ell,\bar{\bold J}^\ell,\bold f^\ell)$ for $\ell =1,2$, 
\ub{then} $M^1_\partial,M^2_\partial$ are not isomorphic 
over $M_\partial$.
\ermn
Why this holds?  As $\nu_1 \ne \nu_2$, by symmetry without loss of generality  for some
$\varepsilon <\partial$ we have $\nu_1(\varepsilon)
=1,\nu_2(\varepsilon)=0$, and let $f$ be an isomorphism from
$M^1_\partial$ onto $M^2_\partial$ over $M_\partial$ and let $E^\ell$
witness $(\bar M^{\rho[\nu_\ell]},\bar{\bold J}^{\rho[\nu_\ell]},\bold
f^{\rho[\nu_\ell]}) \le^{\text{qt}}_{\frak u} (\bar M^\ell,\bar{\bold
J}^\ell,\bold f^\ell)$ for $\ell=1,2$.  Let $ E := E^1 \cap E^2 \cap
\{\delta < \partial:\delta \in \cup\{e^{\rho[\nu_\ell] \restriction
i}:i < \partial\}$ and $\delta = \text{ otp}(\delta \cap
 e_{\rho[\nu_\ell] \rest \delta})$ for $\ell=1,2$ and $f$ maps
$M^1_\delta$ onto $M^2_\delta\}$, clearly it
is a club of $\partial$.

Hence there is $\delta \in S_\varepsilon \cap E$, so
$\rho[\nu_1](\delta) = 1,\rho[\nu_2](\delta) = 0$, now the
contradiction is easy, recalling:
\mr
\item "{$(*)_1$}"   clause $(\delta)$ of Defintion \scite{838-3r.1}(1)
\sn
\item "{$(*)_2$}"  $\bold d_\delta$ does not depend on $\rho[\nu_\ell]
\restriction \delta$.
\ermn
We still owe the proof in the case $\xi = \partial +1$, it is similar
with two changes.  The first is in the choice of $\bold d_\delta$, as
now $\alpha_{\bold d_\delta}$ may be $\partial$, so $\beta_\delta$ may be
$\partial$, hence we have to omit ``$\beta_{\delta(1)} \le \delta$" in the
definition of $E_1$.  In clause $\circledast(j),(k)$ we should replace
$\bold d_\alpha$ by $\bold d_\alpha \rest (0,\gamma_{\rho \rest
\alpha})$ so $\alpha + \gamma_{\rho \rest \delta} < \text{
min}(e_\rho \backslash (\alpha +1))$ where $\gamma_{\rho \rest
\alpha} < \text{ min}\{\ell g(\bold d_\alpha) +1,\partial\}$ which is
the minimal $\alpha$ when we apply $\alpha$ in Definition
\scite{838-3r.1}.

Second, the choice of $\langle \rho_\nu:\nu \in {}^\partial 2\rangle$
with $\rho_\nu \in {}^\partial 2$ is more involved.  

For each $\varepsilon < \partial$ we choose $\varrho_\varepsilon \in
{}^{(S_\varepsilon)}2$ such that:
\mr
\item "{$\boxdot$}"  if $\rho_1,\rho_2 \in {}^\partial 2$ then for
stationarily many $\delta \in S_\varepsilon$ we have: $\bold
f^{\rho_1}(\delta) \le \bold f^{\rho_2}(\delta) \Leftrightarrow
\varrho_\varepsilon(\delta) = 1$
\ermn
(noting that $\bold f^{\rho_\alpha}(\delta)$ depends only on $\rho
\rest \delta$).
\sn
[Why possible?  As $S_\varepsilon$ is not in the weak diamond ideal.]

Then we replace $\odot_2$ by
\mr
\item "{$\odot'_2$}"  for $\nu \in {}^\partial 2$ let $\rho_\nu =
\rho[\nu] \in {}^\partial 2$ be defined by $\rho_\nu(i) =
\varrho_\varepsilon(i) + \nu(\varepsilon)$ mod $2$ when $i \in
S_\varepsilon \wedge \varepsilon < \partial$ and $\rho_\nu(i)$ is zero
otherwise.
\ermn
Why this is O.K.?  I.e. we have to prove $\odot_3$ in this case.  Why
this holds?  As $\nu_1 \ne \nu_2$ by symmetry without loss of generality 
$\nu_1(\varepsilon) =1,\nu_2(\varepsilon)=0$ and let $f,E$ be as
before.

Let $\rho_\ell := \rho[\nu_\ell] \in {}^\partial 2$, so by the choice
of $\varrho_\varepsilon$ there is $\delta \in E \cap S_\varepsilon$
such that

$$
\bold f^{\rho_1}(\delta) \le \bold f^{\rho_2}(\delta) \Leftrightarrow
\varrho_\varepsilon(\delta)=1.
$$
\mn
First assume $\varrho_\varepsilon(\delta) =1$ hence
$\bold f^{\rho_1}(\delta) \le \bold f^{\rho_2}(\delta)$, so $\rho_1(\delta)
= \varrho_\varepsilon(\delta) + \nu_1(\varepsilon) = 1 + 1 = 0$ mod
$2$ so $\rho_1(\delta)=0$ and $\rho_2(\delta) =
\varrho_\varepsilon(\delta) + \nu_2(\varepsilon) = 1 +0 =1$ mod $2$,
so $\nu_2(\delta) =1$.

Now we continue as before, because what we need there for $(\bar
M^{\rho[\nu_2]},\bar{\bold J}^{\rho[\nu_2]},\bold f^{\rho[\nu_2]})$ in
$\delta$ is satisfied for $\delta + \bold f^{\nu[\rho_2]}(\delta)$
hence also for $\delta + \bold f^{\nu[\rho_1]}$.

The other case, $\varrho_\varepsilon(\delta)=0$, is similar;
exchanging the roles.    \hfill$\square_{\scite{838-3r.3}}$ 
\enddemo
\bn
The following conclusion will be used in \scite{838-6u.24}, \scite{838-6u.25}.
\demo{\stag{838-3r.7} Conclusion}  We have
$\dot I(\partial^+,{\frak K}) \ge
\mu_{\text{unif}}(\partial^+,2^\partial)$ and moreover 
$\dot I(K^{{\frak u},{\frak h}}_{\partial^+}) \ge
\mu_{\text{unif}}(\partial^+,2^\partial)$ for any
$\{0,2\}$-appropriate ${\frak h}$ \ub{when} ($\xi \le \partial +1$ and):
\mr
\item "{$(a)$}"  $2^\partial < 2^{\partial^+}$ and $\partial > \aleph_0$
\sn
\item "{$(b)$}"  $(\alpha) \quad {\Cal D}_\partial$, the club filter
on $\partial$, is not $\partial^+$-saturated
\sn
\item "{${{}}$}"  $(\beta) \quad$ if $\xi = \partial +1$ then
$\text{ WDmId}(\partial) + (\partial \backslash S)$ is not
$\partial^+$-saturated
\sn
\item "{$(c)$}"  $\{0,2\}-S$-almost every $(\bar M,\bar{\bold J},\bold f)$ has
the weak $\xi$-uq-invariant coding property even just above $(\bar
M^*,\bar{\bold J},\bold f)$, so $S \subseteq \partial$ is stationary.
\endroster
\enddemo
\bigskip

\demo{Proof}  By \scite{838-3r.3} we can apply \scite{838-3r.9} below using
\mr
\item "{$\odot$}"  if ${\Cal D}$ is a normal filter on the regular
uncountable $\partial,{\Cal D}$ not $\partial^+$-saturated then also
the normal ideal generated by ${\Cal A}$ is not $\partial^+$-saturated
where ${\Cal A} = \{A \subseteq \partial:A \in {\Cal D}$ or $\partial
\backslash A \in {\Cal A}^+$ but ${\Cal D} + (\partial \backslash A)$
is $\partial$-saturated$\}$.
\endroster
\nl ${{}}$ \hfill$\square_{\scite{838-3r.7}}$ 
\enddemo
\bigskip

\proclaim{\stag{838-3r.9} Theorem}  1) We have 
$\dot I(\partial^+,K^{\frak u}_{\partial^+}) \ge
\mu_{\text{unif}}(\partial^+,2^\partial)$; moreover 
$\dot I(K^{{\frak u},{\frak h}}_{\partial^+}) \ge
\mu_{\text{unif}}(\partial^+,2^\partial)$ for any 
${\frak u}-\{0,2\}$-appropriate ${\frak h}$ (see Definition
\scite{838-1a.47}) \ub{when}:
\mr
\item "{$(a)$}"   $2^\partial < 2^{\partial^+}$
\sn
\item "{$(b)$}"  ${\Cal D}$ is a non-$\partial^+$-saturated normal
filter on $\partial$
\sn
\item "{$(c)$}"   for $\{0,2\}$-almost every $(\bar M,\bar{\bold
J},\bold f) \in K^{\text{qt}}_{\frak u}$ (maybe above some such
triple $(\bar M^*,\bar{\bold J}^*,\bold f^*)$ satisfying ${\Cal D}_\partial +
\bold f^{-1}\{0\}$ is not $\partial^+$-saturated), if $S \subseteq
\partial$ belongs to ${\Cal D}^+$ and $\bold f \rest S$ is constantly zero then
we can find a sequence $\langle(\bar M^\alpha,\bar{\bold
J}^\alpha,\bold f^\alpha):\alpha < 2^\partial\rangle$ such that
{\roster
\itemitem{ $(\alpha)$ }  $(\bar M,\bar{\bold J},\bold f)
\le^{\text{qt}}_{\frak u} (\bar M^\alpha,\bar{\bold J}^\alpha,\bold
f^\alpha)$ and $\bold f^\alpha \rest (\partial \backslash S) = \bold f
\rest (\partial \backslash S)$
\sn   
\itemitem{ $(\beta)$ }  if $\alpha(1) \ne \alpha(2) < 2^\partial$ and
$(\bar M^{\alpha(i)},\bar{\bold J}^{\alpha(i)},\bold f^{\alpha(i)})
\le^{\text{qt}}_{\frak u} 
(\bar M^{\ell,*},\bar{\bold J}^{\ell,*},\bold f^{\ell,*})$ for
$\ell=1,2$ then $M^{1,*}_\partial,M^{2,*}_\partial$ are not isomorphic
over $M_\partial$.
\endroster}
\ermn
2) Similarly omitting the ``$\partial^+$-saturation" demands in
clauses (b),(c) and omitting $\bold f \rest S$ is constantly zero in
clause (c).
\endproclaim
\bigskip

\demo{Proof}  1) By Observation \scite{838-1a.49}(4) 
without loss of generality  ${\frak h} = {\frak h}_0
\cup {\frak h}_2$ witness clause (c) of the assumption; we shall use
${\frak h}_2$ for $S^*_0$ so without loss of generality  ${\frak h}_2$ is a $2-S^*_0$-appropriate.
By clause (b) of the assumption let $\bar S^*$ be such that
\mr
\item "{$\odot$}"  $(a) \quad \bar S^* = \langle S^*_\alpha:\alpha <
\partial^+\rangle$
\sn
\item "{${{}}$}"  $(b) \quad S^*_\alpha \subseteq \partial$ for
$\alpha < \partial^+$
\sn
\item "{${{}}$}"  $(c) \quad S^*_\alpha \backslash S^*_\beta \in
[\partial]^{< \partial}$ for $\alpha < \beta < \partial^+$
\sn 
\item "{${{}}$}"  $(d) \quad S^*_{\alpha +1} \backslash S^*_\alpha$ is
a stationary subset of $\partial$; moreover $\in {\Cal D}^+$
\sn
\item "{${{}}$}"  $(e) \quad S = S^*_0$ is stationary; it includes 
$\{\delta < \partial:\bold f^*(\delta) > 0\}$ when $(\bar M^*,\bar{\bold
J}^*,\bold f^*)$ is given.
\ermn
Let $(\bar M^*,\bar{\bold J}^*,\bold f^*)$ be as in clause (c) of the
assumption; so $\{0,2\}$-almost every
$(\bar M^{**},\bar{\bold J}^{**},\bold f^{**}) \in K^{\text{qt}}_{\frak u}$
which is $\le_{\text{qs}}$-above $(\bar M^*,\bold J^*,\bold f^*)$ is
as there witnessed by ${\frak h}$.

Now we
choose $\langle (\bar M^\eta,\bar{\bold J}^\eta,\bold f^\eta):\eta \in
{}^\alpha(2^\partial)\rangle$ by induction on $\alpha < \partial^+$ such that
\mr
\item "{$\circledast$}"  $(a) \quad (\bar M^\eta,\bar{\bold J}^\eta,
\bold f^\eta) \in K^{\text{qt}}_{\frak u}$, and is equal to
$(\bar M^*,\bold J^*,\bold f^*)$ if $\eta = <>$
\sn
\item "{${{}}$}"  $(b) \quad \langle(\bar M^{\eta \restriction
\beta},\bar{\bold J}^{\eta \restriction \beta},\bold f^{\eta
\restriction \beta}):\beta \le \ell g(\eta) \rangle$ is 
$\le^{\text{qs}}_{\frak u}$-increasing continuous 
\sn
\item "{${{}}$}"  $(c) \quad \bold f^\eta \restriction (\partial \backslash
S^*_{\ell g(\eta)+1})$ is constantly zero
\sn
\item "{${{}}$}"  $(d) \quad$ if $\ell g(\eta) = \beta +2 \le \alpha$
and $\nu = \eta \restriction (\beta +1))$ then \nl

\hskip25pt  $(\bar M^\eta,\bar{\bold J}^\eta,
\bold f^\eta) \le^{\text{at}}_{\frak u}(\bar M^\nu,\bar{\bold J}^\nu,
\bold f^\nu)$ and this pair strictly $S$-obeys ${\frak h}$
\sn
\item "{${{}}$}"  $(e) \quad$ if $\ell g(\eta) = \delta <
\alpha,\delta$ limit or zero, $\varepsilon^0 < \varepsilon^1 <
2^\partial$ and
\nl 

\hskip25pt
$(\bar M^{\eta \char 94
<\varepsilon^\ell>},\bar{\bold J}^{\eta \char 94
<\varepsilon^\ell>},\bold f^{\eta \char 94 <\varepsilon^\ell>})$
$ \le_{\frak u}^{\text{qt}}$
$(\bar M^\ell,\bar{\bold J}^\ell,\bold f^\ell)$ for
$\ell=0,1$
\nl 

\hskip15pt
\ub{then}
$M^1_\partial,M^2_\partial$ are not isomorphic over $M^\eta_\partial$.
\ermn
The inductive construction is straightforward:

if $\alpha = 0$ let $(\bar M^\eta,\bar{\bold J}^\eta,\bold f^\eta) =
(\bar M^*,\bar{\bold J}^*,\bold f^*)$

if $\alpha$ is limit use claim \scite{838-1a.37}(4)

if $\alpha = \beta +2$ use clause $\circledast(d)$

if $\alpha = \delta +1,\delta$ limit or zero use clause (c) of the assumption to
satisfy clause $\circledast(e)$.

Having carried the induction, let $M_\eta = \cup\{M^{\eta \restriction
\alpha}:\alpha < \partial^+\}$ for $\eta \in
{}^{\partial^+}(2^\partial)$.  By \scite{838-7f.7} we get that
$|\{M^\eta/\cong:\eta \in {}^{\partial^+}(2^\partial)\}| \ge
\mu_{\text{unif}}(\partial^+,2^\partial)$ so we are done.
\nl
2) Similarly.  \hfill$\square_{\scite{838-3r.7}}$ 
\enddemo
\bn
\centerline {$* \qquad * \qquad *$}
\bn
We now note how we can replace the $\xi$-uq-invariant by $\xi$-up-invariant, a
relative, not used.
\definition{\stag{838-3r.11} Definition}  Let $\xi \le \partial +1$.
\nl
1) We say that $(M,N,\bold I)
\in \text{ FR}^1_{\frak u}$ has the weak $\xi$-up-invariant coding
property \ub{when}:
\mr
\item "{$\circledast$}"  if $M \le_{\frak u} M'$ and $M' \cap N =
\emptyset$ \ub{then} there are $\alpha_\ell < \xi$ and ${\frak
u}$-free $(\alpha_\ell,0)$-rectangle $\bold d_\ell$ for $\ell=1,2$
such that:
{\roster
\itemitem{ $(a)$ }  $M^{\bold d_1}_{0,0} = M' = M^{\bold d_2}_{0,0}$
\sn
\itemitem{ $(b)$ }  $M^{\bold d_1}_{\alpha_1,0} = M^{\bold d_2}_{\alpha_2,0}$
\sn
\itemitem{ $(c)$ }  $M^{\bold d_\ell}_{\alpha_\ell,0} \cap N = M$
\sn
\itemitem{ $(d)$ }  if $(M,N,\bold I) \le_1 (M',N',\bold I')$ and 
$M^{\bold d_\ell}_{\alpha_\ell,0} \cap N' = M'$ \ub{then} there are no
$\alpha'_\ell \le \alpha_\ell,\alpha'_\ell < \partial,
\beta_\ell < \partial$ and ${\frak u}$-free
$(\alpha'_\ell,\beta_\ell)$-rectangles $\bold d^\ell$ for $\ell=1,2$
such that
\sn
\itemitem{ ${{}}$ }  $\quad \bullet_1 \quad \bold d^\ell \restriction 
(\alpha_\ell,0) = \bold d_\ell \rest (\alpha'_\ell,0)$
\sn
\itemitem{ ${{}}$ }  $\quad \bullet_2 \quad (M',N',\bold I') \le_1 
(M^{\bold d^\ell}_{0,0},M^{\bold d^\ell}_{0,1},\bold I^{\bold
d^\ell}_{\alpha,0})$
\sn
\itemitem{ ${{}}$ }  $\quad \bullet_3 \quad$ there are $N'',f$ such
that $M^{\bold d^2}_{\alpha'_2,\beta_2} \le_{\frak u} N$ and $f$ is a
$\le_{\frak u}$-embedding 
\nl

\hskip40pt of $N^{\bold d^1}_{\alpha'_1,\beta_1}$ into
$N$ over $M^{\bold d_1}_{\alpha'_1,0} = M^{\bold d_2}_{\alpha_2,0}$
mapping $M^{\bold d^1}_{0,\beta_1}$ 
\nl

\hskip40pt onto $M^{\bold d^2}_{0,\beta_2}$.
\endroster}
\ermn
2)-5) As in Definition \scite{838-3r.1} replacing uq by up.
\enddefinition
\bigskip

\proclaim{\stag{838-3r.15} Claim}  Like \scite{838-3r.3} replacing
uq-invariant by up-invariant.
\endproclaim
\bigskip

\demo{Proof}  Similar.   \hfill$\square_{\scite{838-3r.15}}$
\enddemo
\bigskip

\demo{\stag{838-3r.17} Conclusion}  Like \scite{838-3r.7} replacing
uq-invariant by up-invariant (in clause (c)).
\enddemo
\bigskip

\demo{Proof}  Similar. \hfill$\square_{\scite{838-3r.21}}$
\enddemo
\bn
\centerline {$* \qquad * \qquad *$}
\bn
Another relative is the vertical one.
\definition{\stag{838-3r.19} Definition}  Let $\xi \le \partial +1$,
omitting $\xi$ means $\partial +1$.  We say that $(M,N,\bold I) \in
\text{ FR}^1_{\frak u}$ has the vertical $\xi$-uq-invariant coding$_1$
property \ub{when}:
\mr
\item "{$\otimes$}"  if $\alpha_0 < \partial$ and $\bold d_0$ is an
${\frak u}$-free $(\alpha_0,0)$-rectangle satisfying $M \le_{\frak u}
M^{\bold d_0}_{0,0}$ and $M^{\bold d_0}_{\alpha_0,0} \cap N = M$ \ub{then}
there are $\alpha,\bold d$ such that:
{\roster
\itemitem{ $(a)$ }  $\alpha_0 < \alpha < \xi$
\sn
\itemitem{ $(b)$ }  $\bold d$ is a ${\frak u}$-free
$(\alpha,0)$-rectangle, though $\alpha$ is possibly $\partial$ this is O.K.
\sn
\itemitem{ $(c)$ }  $\bold d \restriction (\alpha_0,0) = \bold d_0$
\sn
\itemitem{ $(d)$ }  $M^{\bold d}_{\alpha,0} \cap N = M$
\sn
\itemitem{ $(e)$ }  for every $N',\bold I'$ such that $(M^{\bold
d}_{0,0},N',\bold I') \in \text{ FR}_1$ 
is $\le^1_{\frak u}$-above $(M,N,\bold I)$ and
$N' \cap M^{\bold d}_{\alpha,0} = M^{\bold d}_{0,0}$ 
we can find $\alpha',\alpha_*,\bold d_*,\bold I'',M''$ such that
\sn
\itemitem{ ${{}}$ }  $(\alpha) \quad \alpha' \le \alpha,\alpha' <
\partial,\alpha_0 \le \alpha_*$ 
\sn
\itemitem{ ${{}}$ }  $(\beta) \quad \bold d_*$ is a ${\frak u}$-free
$(\alpha_*,0)$-rectangle
\sn
\itemitem{ ${{}}$ }  $(\gamma) \quad M^{\bold d_*}_{0,0} = N'$
\sn
\itemitem{ ${{}}$ }  $(\delta) \quad$ there is an ${\frak u}$-free
$(\alpha_0,1)$-rectangle $\bold d'$ such that 
\nl

\hskip35pt $\bold d' \restriction (\alpha_0,0) = 
\bold d_0,\bold d' \restriction
([0,\alpha_0),[1,1]) = \bold d_* \restriction (\alpha_0,0)$ and
\nl

\hskip35pt $\bold I^{\bold d'}_{0,0} = \bold I'$ and 
$(M^{\bold d'}_{\alpha_0,0},M^{\bold d'}_{\alpha_0,1},\bold I') \le^1_{\frak u}
(M'',M^{\bold d_*}_{\alpha_*,0},\bold I'')$
\sn
\itemitem{ ${{}}$ }  $(\varepsilon) \quad$ there are no $\bold d_1,
\bold d_2$ such that
\sn
\itemitem{ ${{}}$ }  $\qquad \bullet_1 \quad \bold d_\ell$ is a
${\frak u}$-free rectangle for $\ell=1,2$
\sn
\itemitem{ ${{}}$ }  $\qquad \bullet_2 \quad \alpha(\bold d_1) =
\alpha_*$ and $\bold d_1 \restriction (\alpha_*,0) = \bold d_*$
\sn
\itemitem{ ${{}}$ }  $\qquad \bullet_3 \quad \alpha(\bold d_2) \ge
\alpha'$ and $\bold d_2 \restriction (\alpha',0) = \bold d \rest (\alpha',0)$
\sn
\itemitem{ ${{}}$ }  $\qquad \bullet_4 \quad (M^{\bold
d_2}_{0,0},M^{\bold d_2}_{0,1},\bold I^{\bold d_2}_{0,0}) =
(M^{\bold d}_{0,0},N',\bold I')$ or just 
\nl

\hskip48pt $(M^{\bold d}_{0,0},N',\bold I') \le^1_{\frak u} (M^{\bold
d_2}_{0,0},M^{\bold d_2}_{0,1},\bold I^{\bold d_2}_{0,0})$
\sn
\itemitem{ ${{}}$ }  $\qquad \bullet_5 \quad$ there are $N_1,N_2,f$
such that $M^{\bold d_\ell}_{0,\beta(\bold d_\ell)} \le_{\frak u}
N_\ell$ for $\ell=1,2$
\nl

\hskip35pt  and $f$ is an isomorphism from $N_1$ onto $N_2$
over $M^{\bold d}_{\alpha,0}$
\nl

\hskip35pt  mapping $M^{\bold d_1}_{0,\beta(\bold
d_1)}$ onto $M^{\bold d_2}_{0,\beta(\bold d_2)}$.
\endroster}
\ermn
2) We say that $(\bar M^*,\bar{\bold J}^*,\bold f^*) \in
K^{\text{qt}}_{\frak u}$ has the vertical uq-invariant coding$_1$
propety as in Definition \scite{838-3r.1}(2) only $(\langle M_{\delta +i}:i \le
\bold f(\delta)\rangle,\langle\bold J_{\delta +i}:i < \bold
f(\delta)\rangle)$ play the role of $\bold d_0$ in part (1).  In all
parts coding means coding$_1$.
\nl
3),4),5) Parallely to Definition \scite{838-3r.1}. 
\enddefinition
\bigskip

\proclaim{\stag{838-3r.21} Theorem}  Like \scite{838-3r.3} 
using vertical $\xi$-uq-invariant coding in clause (b) and omitting
clause (c) of the assumption and omit clause $(\beta)$ in the
conclusion.  \hfill$\square_{\scite{838-3r.21}}$
\endproclaim
\bigskip

\demo{Proof}  Similarly.
\enddemo
\bigskip

\demo{\stag{838-3r.23} Conclusion}  Like \scite{838-3r.7} replacing
clause (b)$(\beta)$ of the assumption (by $\xi = \partial +1
\Rightarrow (\exists \theta)2^\theta = 2^{< \partial} < 2^\partial$)
and with vertical $\xi$-uq-invariant 
coding instead of the $\xi$-uq-invariant one (in clause (c), can use
$S = \partial$).
\enddemo
\bigskip

\demo{Proof}  Similar to the proof of \scite{838-3r.7}.
\hfill$\square_{\scite{838-3r.23}}$ 
\enddemo
\bn
\centerline {$* \qquad * \qquad *$}
\bn
\margintag{3r.67}\ub{\stag{838-3r.67} Discussion}:  The intention below is to help in \S6 to 
eliminate the assumption
``{\rm WDmId}$(\lambda^+)$ is not $\lambda^{++}$-saturated" when ${\frak
s}$ fails existence for $K^{3,\text{up}}_{{\frak s},\lambda^+}$.
We do using the following relatives, semi and vertical,
from Definition \scite{838-3r.71}, \scite{838-3r.19} are interesting because
\mr
\item "{$(a)$}"   under reasonable conditions (see 
Definition \scite{838-3r.84}) the first implies the second
\sn
\item "{$(b)$}"   the second, as in Theorem \scite{838-2b.13} is enough
for non-structure without the demand on saturation of WDmId$(\partial)$
\sn
\item "{$(c)$}"   the first needs a weak version of a model theoretic
assumption (in the application)
\sn
\item "{$(d)$}"  (not used) the semi-version implies the weak version
(from \scite{838-3r.1}).
\endroster
\bigskip

\definition{\stag{838-3r.71} Definition}  Let $\xi \le \partial +1$.
\nl
1) We say that ${\frak u}$ has
the semi $\xi$-uq-invariant coding$_1$ property [above some $(\bar
M^*,\bar{\bold J}^*,\bold f^*)$]
 \ub{when} for
$\{0,2\}$-almost every $(\bar M,\bar{\bold J},\bold f) \in
K^{\text{qt}}_{\frak u}$ [above $(\bar M^*,\bar{\bold J}^*,\bold
f^*)$] for some $(\alpha,N,\bold I)$ we have $\alpha < \partial,M \cap
N = M_\alpha$ and
$(M_\alpha,N,\bold I) \in \text{ FR}^+_1$ has the semi $\xi$-uq-invariant
coding$_1$ property, see below but restricting ourselves to
$M',M^{\bold d}_{\alpha_{\bold d},0} \in \{M_\beta:\beta \in
(\alpha,\partial)\}$.  Here and in part (2) we may write coding
instead of coding$_1$.
\nl
2) We say that $(M,N,\bold I) \in
\text{ FR}^1_{\frak u}$ has the semi $\xi$-uq-invariant coding$_1$ property
(we may omit the 1) \ub{when}: if $M \le_{\frak u} M'$ 
and $M' \cap N = M$ \ub{then} we can find $\bold d$ such that:
\mr
\item "{$\circledast$}"  $(a) \quad \bold d$ is a ${\frak u}$-free 
$(\alpha_{\bold d},0)$-rectangle with $\alpha_{\bold d} < \xi$ so
$\alpha_{\bold d} \le \partial$
\sn
\item "{${{}}$}"  $(b) \quad M^{\bold d}_{0,0} = M'$ and $M^{\bold
d}_{\alpha(\bold d),0} \cap N = M$
\sn
\item "{${{}}$}"  $(c) \quad$ for any $N',\bold I'$; if $(M,N,\bold I) \le_1
(M',N',\bold I')$ and $M^{\bold d}_{\alpha(\bold d),0} \cap N' = M'$
\nl

\hskip25pt  \ub{then} we can find $\alpha',N'',\bold I''$ satisfying
$\alpha' \le \alpha,\alpha' < \partial$ and
$(M',N',\bold I') \le_1$ 
\nl

\hskip25pt  $(M^{\bold d}_{\alpha',0},N'',\bold I'')
\in \text{ FR}^1_{\frak u}$ such that for no triple $(\bold e,f,N_*)$
\nl

\hskip25pt  do we have:
{\roster
\itemitem{ $(\alpha)$ }  $\bold e$ is a ${\frak u}$-free rectangle
\sn
\itemitem{ $(\beta)$ }  $\alpha_{\bold e} = \alpha_{\bold d}$ and
$\bold e \rest (\alpha',0) = \bold d \rest (\alpha',0)$
\sn
\itemitem{ $(\gamma)$ }  $(M^{\bold e}_{0,0},M^{\bold e}_{0,1},\bold
I^{\bold e}_{0,0}) = (M^{\bold d}_{0,0},N'',\bold I'')$
\sn
\itemitem{ $(\delta)$ }  $M^{\bold e}_{\alpha',\beta(\bold e)}
\le_{\frak u} N_*$
\sn
\itemitem{ $(\varepsilon)$ }  $f$ is a $\le_{\frak u}$-embedding of
$N''$ into $N_*$
\sn
\itemitem{ $(\zeta)$ }  $f \rest M^{\bold d}_{\alpha',0}$ is
the identity
\sn
\itemitem{ $(\eta)$ }  $f$ maps $N'$ into $M^{\bold e}_{0,\beta(\bold e)}$.
\endroster}
\endroster
\enddefinition
\bigskip

\remark{\stag{838-3r.72} Remark}  1) This is close to 
Definition \scite{838-3r.1} but simpler,
cover the applications here and fit Claim \scite{838-3r.88}.
\nl
2) We could have phrased the other coding properties similarly.
\endremark
\bigskip

\proclaim{\stag{838-3r.74} Claim}  If $(M,N,\bold I) \in 
\text{\rm FR}^1_{\frak u}$ has the semi $\xi$-uq-invariant 
coding property \ub{then}
it has the weak $\xi$-uq-invariant coding property, see 
Definition \scite{838-3r.1}.
\endproclaim
\bigskip

\demo{Proof}  Should be clear.  \hfill$\square_{\scite{838-3r.74}}$
\enddemo
\bn
The following holds in our natural examples when we add the fake, i.e.
artificial equality and it is natural to demand ${\frak u} = {\frak
u}^{[*]}$, see Definition \scite{838-3r.85}.
\definition{\stag{838-3r.84} Definition}  1) We say that ${\frak u}$ has the fake
equality $=_*$ \ub{when}:
\mr
\item "{$(a)$}"  $\tau_{\frak K}$ has only predicates and
some two-place relation $=_* \in \tau_{\frak K}$ is,
for every $M \in K$, interpreted as an equivalence relation which is a
congruence relation on $M$
\sn
\item "{$(b)$}"  $M \in K$ iff $M/=^M_*$ belongs to $K$
\sn
\item "{$(c)$}"  for $M \subseteq N$ both from $K$ we have 
$M \le_{\frak u} N$ iff $(M/=^N_*) \le_{\frak s} (N/=^N_*)$
\sn
\item "{$(d)$}"  assume $M \le_{\frak u} N$ and $\bold I \subseteq N
\backslash M$ and $\bold I' = \{d \in \bold I:(\forall c \in M)(\neg c
=^N d)\},\ell \in \{1,2\}$.  If $(M,N,\bold I) \in \text{ FR}_\ell$
then $(M/=^N,N/=^N,\bold I'/=^N) \in \text{ FR}_\ell$ which implies
$(M,N,\bold I') \in \text{ FR}_\ell$
\sn
\item "{$(e)$}"  if $M \subseteq N$ are from $K$ and $\bold I
\subseteq \{d \in N:(\forall c \in M)(\neg c =^N d)\}$ and $\ell \in
\{1,2\}$ \ub{then} $(M,N,\bold I) \in \text{ FR}_\ell$ iff
$(M/=^N,\bold I/=^N) \in \text{ FR}_\ell$.
\ermn
1A) In part (1) we may say that ${\frak u}$ has the fake equality
$=_*$ or $=_*$ is a fake equality for ${\frak u}$. 
\nl
2) We say ${\frak u}$ is hereditary \ub{when} every 
$(M,N,\bold I) \in \text{ FR}^+_1$ is hereditary, see below.
\nl
3) We say $(M,N,\bold I) \in \text{ FR}^1_{\frak u}$ is hereditary \ub{when}:
\mr
\item "{$(a)$}"  if $\bold d$ is ${\frak u}$-free $(1,2)$-rectangle,
$M \le_{\frak u} M^{\bold d}_{0,0}$ and $(M,N,\bold I) \le^1_{\frak u}
(M^{\bold d}_{0,1},M^{\bold d}_{0,2},\bold I^{\bold d}_{0,1})$ \ub{then}
$(M,N,\bold I) \le^1_{\frak u} (M^{\bold d}_{0,0},M^{\bold  d}_{0,2},
\bold I^{\bold d}_{0,1}) \le^1_{\frak u} (M^{\bold
 d}_{1,0},M^{\bold d}_{1,2},\bold I^{\bold d}_{1,1})$.
\ermn
4) We say ${\frak u}$ is hereditary for the fake equality $=_*$ when
   every $(M,N,\bold I) \in \text{ FR}^+_1$ is hereditary for $=_*$
   which means that clause (a) of part (3) above holds, 
$=_*$ is a fake equality for ${\frak u}$ and:
\mr
\item "{$(b)$}"  if $\bold d$ is a ${\frak u}$-free $(0,2)$-rectangle,
 $M \le_{\frak u} M^{\bold d}_{0,0}$ and $(M,N,\bold I) \le^1_{\frak
 u} (M^{\bold d}_{0,1},M^{\bold d}_{0,2},\bold I^{\bold d}_{0,1})$
 \ub{then} we can find $M_1,f,M_2$ such that:
{\roster
\itemitem{ $(\alpha)$ }  $f$ is an isomorphism from $M^{\bold
 d}_{0,1}$ onto $M_1$ over $M^{\bold d}_{0,0}$
\sn
\itemitem{ $(\beta)$ }  $M^{\bold d}_{0,0} \le_{\frak u} 
M_1 \le_{\frak u} M_2$ and $M_{0,2} \le_{\frak u} M_2$
\sn
\itemitem{ $(\gamma)$ }  $|M_2| = |M_1| \cup |M^{\bold d}_{0,2}|$
\sn
\itemitem{ $(\delta)$ }  $M_2 \models ``c =_* f(c)"$ if $c \in
M_{0,1}$
\sn
\itemitem{ $(\varepsilon)$ }  $(M_{0,0},M_{0,2},\bold I^{\bold d}_{0,1})
\le^1_{\frak u} (M_1,M_2,\bold I^{\bold d}_{0,1})$.
\endroster}
\ermn
5)  In parts (2),(3),(4) we can replace hereditary by weakly
hereditary \ub{when}: in clause (a) we assume $\bold I =  
\bold I^{\bold d}_{0,1} = \bold I^{\bold d}_{1,1}$ and in clause (b) we assume
$\bold I = \bold I^{\bold d}_{0,1}$.
\enddefinition
\bigskip

\definition{\stag{838-3r.85} Definition}  For ${\frak u}$ is a nice
    construction framework we define ${\frak u}^{[*]} = {\frak
    u}^{[*]}$ like ${\frak u}$ except that, for $\ell=1,2$ we have
    $(M_1,N_1,\bold I_1) \le^\ell_{{\frak u}[*]} (M_2,N_2,\bold I_2)$
    \ub{iff} $(M_1,N_1,\bold I_1) \le^\ell_{\frak u} (M_2,N_2,\bold
    I_2)$ and $\bold I_1 \ne \emptyset \Rightarrow \bold I_1 = \bold I_2$.
\enddefinition
\bigskip

\demo{\stag{838-3r.86} Observation}  1) ${\frak u}$ has the fake equality =
(i.e. the standard equality is also a fake equality).
\nl
2) ${\frak u}'$ as defined in \scite{838-2b.8}(2) has the fake equality
   $=_\tau$ (and is a nice construction framework, see
   \scite{838-2b.8}(3)).
\nl
3) If ${\frak u}$ is hereditary \ub{then} ${\frak u}''$ as defined in
   \scite{838-2b.8}(4) is hereditary and even hereditary for the fake
   equality $=_\tau$ and is a nice construction framework, see
   \scite{838-2b.8}(5). 
\nl
4) ${\frak u}^{[*]}$ is a nice construction framework, and if ${\frak
   u}$ is weakly hereditary [for the fake equality $=_*$] \ub{then}
   ${\frak u}^{[*]}$ is hereditary [for the fake equality $=_*$]. 
\nl
5) If ${\frak u}$ is hereditarily (for the fake equality $=_*$)
\ub{then} ${\frak u}^{[*]}$ is hereditarily (for the fake equality $=_*$).
\enddemo
\bigskip

\demo{Proof}  Check (really \scite{838-3r.0}).  \hfill$\square_{\scite{838-3r.86}}$
\enddemo
\bigskip

\proclaim{\stag{838-3r.88} Claim}  Let $\xi \in \{\partial,\partial +1\}$
or just $\xi \le \partial +1$.
Assume ${\frak u} =$ 
{\rm dual}$({\frak u})$ has fake equality $=_*$ and is hereditary for $=_*$.

If $(M,N,\bold I)$ has the semi $\xi$-uq-invariant coding property,
\ub{then} $(M,N,\bold I)$ has the vertical $\xi$-uq-invariant coding property.
\endproclaim
\bigskip

\remark{Remark}  So no ``$(\bar M,\bar{\bold J},\bold f) \in
K^{\text{qt}}_{\frak u}$" here, but applying it we use ${\frak
K}^{\frak u}$-universal homogeneous $M_\partial$.
\endremark
\bigskip

\demo{Proof}  So let $\bold d_0$ be a ${\frak u}$-free $(\alpha_{\bold
d},0)$-rectangle satisfying $M \le_{\frak u} M^{\bold d_0}_{0,\delta}$
and $M^{\bold d_0}_{\alpha(\bold d_0),0} \cap N=M$ so $\alpha_{\bold
d_0} < \partial$ and we should find
$\bold d_1$ satisfying the demand in Definition \scite{838-3r.19}(1), 
this suffice.  As we are assuming
that ``$(M,N,\bold I)$ has the semi uq-invariant coding property", 
there is $\bold e_0$ satisfying the
demands on $\bold d$ in \scite{838-3r.71}(1)$\circledast$ with $(M,N,\bold
I,M^{\bold d_0}_{0,0})$ here standing for $(M,N,\bold I,M')$ there.

Without loss of generality
\mr
\item "{$(*)_1$}"  $M^{\bold e_0}_{\alpha(\bold e_0),0} \cap M^{\bold
d_0}_{\alpha(\bold d),0} = M^{\bold d_0}_{0,0}$.
\ermn
Let
\mr
\item "{$(*)_2$}"  $\bold e_1 = \text{\rm dual}(\bold e_0)$ so as ${\frak
u} =$ {\rm dual}$({\frak u})$ clearly $\bold e_1$ is a ${\frak u}$-free
$(0,\alpha(\bold e_0))$-rectangle, so $\beta(\bold e_1) = \alpha(\bold e_0)$.
\ermn
Now by \scite{838-1a.9}(5) for some $\bold e_2$ (note: even the case
$\beta(\bold e_1) = \alpha(\bold e_0) = \partial$ is O.K.)
\mr
\item "{$(*)_3$}"  $(a) \quad \bold e_2$ is a ${\frak u}$-free
$(\alpha(\bold d_0),\beta(\bold e_1))$-rectangle
\sn
\item "{${{}}$}"  $(b) \quad \bold e_2 \rest (\alpha(\bold d_0),0) 
= \bold d_0$
\sn
\item "{${{}}$}"  $(c) \quad \bold e_2 \rest (0,\beta(\bold e_1))) =
\bold e_1$
\ermn
and \wilog
\mr
\item "{${{}}$}"  $(d) \quad M^{\bold e_2}_{\alpha(\bold e_2),\beta(\bold e_2)}
\cap N = M$.
\ermn
Now we choose $\bold d$ by
\mr
\item "{$(*)_4$}"  $\bold d$ is the ${\frak u}$-free $(\alpha(\bold d_0)
+ \alpha(\bold e_0),0)$-rectangle such that:
\sn
\item "{${{}}$}"  $(a) \quad M^{\bold d}_{\alpha,0}$ is 
$M^{\bold d_0}_{\alpha,0}$ if $\alpha \le \alpha(\bold d_0)$
\sn
\item "{${{}}$}"  $(b) \quad M^{\bold d}_{\alpha,0}$ is 
$M^{\bold e_2}_{\alpha(\bold d_0),\alpha-\alpha(\bold d_0)}$ for $\alpha \in
[\alpha(\bold d_0),\alpha(\bold d_0) + \alpha(\bold e_0)]$
\sn
\item "{${{}}$}"  $(c) \quad \bold J^{\bold d}_{\alpha,0} = 
\bold J^{\bold d_0}_{\alpha,0}$ if $\alpha < \alpha(\bold d_0)$
\sn
\item "{${{}}$}"  $(d) \quad \bold J^{\bold d}_{\alpha,0} = 
\bold I^{\bold e_2}_{\alpha(\bold d_0),\alpha - \alpha(\bold d_0)}$ 
if $\alpha = [\alpha(\bold d_0),\alpha(\bold d_0) + \alpha(\bold e_0))$.
\ermn
[Why is this O.K.?  Check; the point is that ${\frak u} =$
{\rm dual}$({\frak u})$.]

And we choose $\bold d_{2,\gamma}$ for $\gamma < \text{
min}\{\alpha(\bold d_0) +1,\partial\}$
\mr
\item "{$(*)_5$}"  $\bold d_{2,\gamma}$ is the ${\frak u}$-free 
$(\gamma +1,0)$-rectangle such that:
\sn
\item "{${{}}$}"  $(a) \quad M^{\bold d_{2,\gamma}}_{\alpha,0}$ is 
$M^{\bold d_0}_{\alpha,0}$ if $\alpha \le \gamma$
\sn
\item "{${{}}$}"  $(b) \quad M^{\bold d_{2,\gamma}}_{\alpha,0}$ is $M^{\bold
e_2}_{\alpha(\bold d),\alpha(\bold e_0)}$ if $\alpha = \gamma +1$
\sn
\item "{${{}}$}"  $(c) \quad \bold J^{\bold d_{2,\gamma}}_{\alpha,0} = 
\bold J^{\bold d_0}_{\alpha,0}$ if $\alpha < \gamma$
\sn
\item "{${{}}$}"  $(d) \quad \bold J^{\bold d_{2,\gamma}}_{\alpha,0}$
is $\emptyset$ if $\alpha = \gamma$.
\ermn
[Why is this O.K.?  Check.]

So it is enough to show
\mr
\item "{$(*)_6$}"  $\bold d$ is as required on $\bold d$ in
\scite{838-3r.19} for our given $\bold d_0$ and $(M,N,\bold I)$.
\ermn
But first note that
\mr
\item "{$(*)_7$}"  $\bold d = \bold d_1$ is as required in clauses
(b),(c),(d) of Definition \scite{838-3r.19}(1).
\ermn
Now, modulo $(*)_7$, clearly $(*)_6$ means that we have to show that
\mr
\item "{$\boxtimes$}"  if $(M,N,\bold I) \le_1 (M^{\bold
d}_{0,0},N',\bold I')$ and $M^{\bold d_1}_{\alpha(\bold d_1),0} \cap N
= M$, i.e. $N',\bold I'$ are as in $\otimes(e)$ of \scite{838-3r.19}(1),
\ub{then} we shall find $\alpha',\alpha_*,\bold d_*,\bold I'',M''$
as required in clauses $(\alpha)-(\varepsilon)$ of (e) 
from Definition \scite{838-3r.19}(1).
\ermn
By the choice of $\bold e_0$ to be as in Definition \scite{838-3r.71}
before $(*)_1$, for $(N',\bold J')$ from $\boxtimes$ 
there are $\alpha',\alpha_*,\bold e_*,N'',\bold I''_*$ such that
\mr
\item "{$\odot_1$}"  $(\alpha) \quad \alpha' \le \alpha(\bold e_0)$
and $\alpha' < \partial$
\sn
\item "{${{}}$}"  $(\beta) \quad \bold e_*$ is a ${\frak u}$-free
$(\alpha_*,0)$-rectangle
\sn
\item "{${{}}$}"  $(\gamma) \quad M^{\bold e_*}_{0,0} = N'$
\sn
\item "{${{}}$}"  $(\delta) \quad$ there is a ${\frak u}$-free
$(\alpha',1)$-rectangle $\bold e'$ such that $\bold e' \rest
(\alpha',0) = \bold e_0 \rest (\alpha',0)$
\nl

\hskip25pt  and $\bold e' \rest
(\alpha',[1,1]) = \bold e_* \rest (\alpha',0)$
\sn
\item "{${{}}$}"  $(\varepsilon) \quad$ there are no $\bold d_1,\bold
d_2$ such that
{\roster
\itemitem{ ${{}}$ }  $\bullet_1 \quad \bold d_\ell$ is a ${\frak
u}$-free rectangle for $\ell=1,2$
\sn
\itemitem{ ${{}}$ }  $\bullet_2 \quad \alpha(\bold d_1) =
\alpha_*,\bold d_2 \rest (\alpha_*,0) = \bold e_*$
\sn
\itemitem{ ${{}}$ }  $\bullet_3 \quad \alpha(\bold d_2,0) \ge \alpha'$
and $\bold d_2 \rest (\alpha',0) = \bold d \rest (\alpha',0)$
\sn
\itemitem{ ${{}}$ }  $\bullet_4 \quad (M^{\bold d_2}_{0,0},M^{\bold
d_2}_{0,1},\bold I^{\bold d_2}_{0,0})$ is $(M',N',\bold I)$ or just
$\le_1$-above it
\sn
\itemitem{ ${{}}$ }  $\bullet_5 \quad$ there are $N^*_1,N^*_2,f$ such
that $M^{\bold d_\ell}_{\alpha(\bold d_\ell),\beta(\bold d_\ell)}
\le_{\frak u} N^*_\ell$ for $\ell=1,2$
\nl

\hskip25pt  and $f$ is an isomorphism from $N^*_1$ onto $N^*_2$ over
\nl

\hskip25pt  $M^{\bold d}_{\alpha(\bold d),0}$ which maps
$M^{\bold d_1}_{0,\beta(\bold d_1)}$ onto $M^{\bold
d_2}_{0,\beta(\bold d_2)}$.
\endroster}
\ermn
Without loss of generality
\mr
\item "{$\odot_2$}"   $N''\cap M^{\bold e_2}_{\alpha',
\beta(\bold e_2)} =  M^{\bold e_0}_{\alpha',0}$.
\ermn
Now by induction on $\alpha \le \alpha(\bold d)$ we choose
$(M^*_\alpha,N^*_\alpha,\bold I^*_\alpha),\bold J^*_\alpha$ such that
\mr
\item "{$\odot_3$}"  $(a) \quad (M^*_\alpha,N^*_\alpha,\bold I^*_\alpha)
\in \text{ FR}^1_{\frak u}$
\sn
\item "{${{}}$}"  $(b) \quad \langle(M^*_\gamma,N^*_\gamma,\bold
I^*_\gamma):\gamma \le \alpha\rangle$ is $\le^1_{\frak u}$-increasing
continuous
\sn
\item "{${{}}$}"  $(c) \quad (M^*_\alpha,N^*_\alpha,\bold I^*_\alpha)
= (M^{\bold e_0}_{\alpha',0},N'',\bold I'')$ for $\alpha =0$
\sn
\item "{${{}}$}"  $(d) \quad M^*_\alpha = 
M^{\bold e_2}_{\alpha,\alpha'}$
\sn
\item "{${{}}$}"  $(e) \quad N^*_\alpha \cap M^{\bold
e_2}_{\alpha(\bold e_2),\alpha'} = M^*_\alpha$
\sn
\item "{${{}}$}"  $(f) \quad$ if $\alpha = \alpha_1 +1$ then
$(M^{\bold e_2}_{\alpha_1,\alpha'},M^{\bold e_2}_{\alpha_1
+1,\alpha'},\bold J^{\bold e_2}_{\alpha_1,\alpha(\bold e_0)}) \le_2
(N^*_{\alpha_1},N^*_{\alpha_1 +1},\bold J^*_{\alpha_1})$.
\ermn
Now we shall use the assumption ``${\frak u}$ is hereditary for $=_*$" to
finish.

Choose
\mr
\item "{$\odot_4$}"  $f$ is an isomorphism from $M^*_{\alpha(\bold
d)}$ onto a model $M^*$ such that $M^{\bold e_2}_{\alpha',
\beta(\bold e_2)} = M^* \cap M^*_{\alpha(\bold d)}$ and $f \rest 
M^{\bold e_2}_{\alpha',\beta(\bold e_2)}$ is
the identity as well as $f \rest N''$
\sn
\item "{$\odot_5$}"  $M^{**}$ is the unique model $\in K_{\frak u}$
such that $M^*_{\alpha(\bold d)} \subseteq M^{**}$ and $M^* \subseteq
M^{**}$ and $c \in M^*_{\alpha(\bold d)} \Rightarrow M^{**} \models
``c = f(c)"$.
\ermn
Lastly
\mr
\item "{$\odot_6$}"  $\bold d_*$ is the following ${\frak u}$-free
$(\alpha(\bold d) +1,1)$-rectangle:
{\roster
\itemitem{ $(a)$ }  $\bold d_* \rest (\alpha(\bold d),0) = \bold d$
\sn
\itemitem{ $(b)$ }  $M^{\bold d_*}_{\alpha(\bold d)+1} =
M^{\bold e_2}_{\alpha(\bold e_2),\alpha}$ and $\bold
J^{\bold d_*}_{\alpha(\bold d),0} = \emptyset$
\sn
\itemitem{ $(c)$ }  $M^{\bold d_*}_{\alpha(\bold d)+1,1}$ is $M^{**}$
\sn
\itemitem{ $(d)$ }  if $\alpha \le \alpha(\bold d)$ the $M^{\bold
d_4}_{\alpha,1}$ is $M^{\bold d}_{\alpha,0} \cup f(M^{\bold
e_2}_{\alpha,\alpha'})$, i.e. the submodel of $M^{**}$ with
this universe
\sn
\itemitem{ $(e)$ }  if $\alpha \le \alpha(\bold d)$ then $\bold
I^{\bold d_*}_{\alpha,0} = f(\bold I^*_\alpha)$ and 
$\bold I^{\bold d_*}_{\alpha(\bold d)} = f(\bold I^*_{\alpha(\bold d)})$
\sn
\itemitem{ $(f)$ }  if $\alpha < \alpha(\bold d)$ then $\bold J^{\bold
d_*}_{\alpha,1} = f(\bold J^{\bold e_2}_{\alpha,\alpha'})$.
\endroster}
\ermn
\hfill$\square_{\scite{838-3r.88}}$
\enddemo
\bn
To phrase a relative of \scite{838-3r.88}, we need:
\definition{\stag{838-3r.89} Definition}  1) We say ${\frak u}$ satisfies
(E)$_\ell$(f), is interpolative for $\ell$ or has interpolation for
$\ell$ \ub{when}: 

if $(M_1,N_1,\bold I_1) \le^\ell_{\frak u} (M_2,N_2,\bold I_2)$ then
$(M_1,N_1,\bold I_1) \le^\ell_{\frak u} (M_2,N_2,\bold I_1)
\le^\ell_{\frak u} (M_2,N_2,\bold I_2)$.
\nl
2) (E)(f) means (E)$_1$(f) + (E)$_2$(f). 
\enddefinition
\bigskip

\remark{Remark}  This is related to but is different from monotonicity,
see \scite{838-1a.24}(1).
\endremark
\bigskip

\proclaim{\stag{838-3r.89.4} Claim}  1) Assume that for $\ell \in \{1,2\},{\frak u}$
satisfies (E)$_\ell$(f).  For every ${\frak u}$-free
$(\alpha,\beta)$-rectangle $\bold d$, also $\bold e$ is a ${\frak
u}$-free $(\alpha,\beta)$-rectangle where $M^{\bold e}_{i,j} =
M^{\bold d}_{i,j},\ell =1 \Rightarrow \bold J^{\bold e}_{i,j} = \bold
J^{\bold d}_{i,j},\ell = 1 \Rightarrow \bold I^{\bold e}_{i,j} = \bold
I^{\bold d}_{0,j}$ and $\ell =2 \Rightarrow \bold J^{\bold e}_{i,j} = \bold
J^{\bold d}_{i,0},\ell = 2 \Rightarrow \bold I^{\bold e}_{i,j} = 
\bold I^{\bold d}_{i,j}$.
\nl
2) Similarly for ${\frak u}$-free triangle.
\nl
3) If ${\frak s}$ satisfies (E)(f) then in part (1) we can let
$\bold J^{\bold e}_{i,j} = \bold J^{\bold d}_{i,0},\bold I^{\bold
e}_{i,j} = \bold I^{\bold d}_{0,j}$.
\endproclaim
\bigskip

\demo{Proof}  Easy.  \hfill$\square_{\scite{838-3r.89.4}}$
\enddemo
\bn
Now we can state the variant of \scite{838-3r.88}.
\proclaim{\stag{838-3r.90} Claim}  Assume $\xi \le \partial +1$ and 
${\frak u} =$ {\rm dual}$({\frak u})$ has fake equality $=_*$, is weakly
hereditary for $=_*$ and has interpolation.

If $(M,N,\bold I) \in K^{\text{qt}}_{\frak u}$ has the semi 
$\xi$-uq-invariant coding property \ub{then} $(M,N,\bold I) \in
K^{\text{qt}}_{\frak u}$ has the vertical $\xi$-uq-invariant coding property.
\endproclaim
\bigskip

\demo{Proof}  Similar to the proof of \scite{838-3r.88}, except that
\mr
\item "{$(A)$}"  in $\odot_3$ we add $\bold I'' = \bold I'$, justified
by monotonicity; note that not only is it a legal choice but still
exemplify \scite{838-3r.19}
\sn
\item "{$(B)$}"  in $\odot_4$ we add $\bold I^*_\alpha = \bold I^*_0 =
\bold I'' = \bold I'$.
\endroster
\enddemo
\bigskip

\proclaim{\stag{838-3r.91} Theorem}  We have $\dot I(K^{{\frak u},{\frak
h}}_{\partial^+}) \ge \mu_{\text{unif}}(\partial^+,2^\partial)$
\ub{when}:
\mr
\item "{$(a)$}"  $2^\partial < 2^{\partial^+}$
\sn
\item "{$(b)$}"  some $(M,N,\bold I) \in \text{\rm FR}_2$ has
the vertical $\xi$-uq-invariant 
coding property, see Definition \scite{838-3r.19}
\sn
\item "{$(c)$}"  ${\frak h}$ is a $\{0,2\}$-appropriate witness that
for $\{0,2\}$-almost every $(\bar M,\bar{\bold J},\bold f) \in
K^{\text{qt}}_{\frak u}$, the model $M_\partial$
is $K_{\frak u}$-model homogeneous
\sn
\item "{$(d)$}"  $\xi = \partial$ \ub{or} $\xi = \partial^+$ and
$2^\theta = 2^{< \partial} < 2^\partial$.
\endroster
\endproclaim
\bigskip

\remark{\stag{838-3r.91.3} Remark}  1) We can phrase other theorems in this way.
\nl
2) So if we change (b) to semi uq-invariant by \scite{838-3r.88} it
suffices to add, e.g.
\mr
\item "{$(d)$}"   ${\frak u}$ is hereditary for the faked equality
$=_*$.
\endroster
\endremark
\bigskip

\demo{Proof}  Easily by clause (b) we know ${\frak u}$ has the vertical
$\xi$-uq-invariant coding property.  Now we apply \scite{838-3r.21},
i.e. as in the proof of theorems \scite{838-3r.7} using
\scite{838-3r.9}(2) rather than \scite{838-3r.9}(1) and immitating the proof
of \scite{838-3r.3}.  \hfill$\square_{\scite{838-3r.91}}$
\enddemo
\bn
\margintag{3r.91A}\ub{\stag{838-3r.91A} Exercise}:  Prove the parallel of the first
sentence of the proof of \scite{838-3r.91} to other coding properties.
\bn
\centerline {$* \qquad * \qquad *$}
\bn
\margintag{3r.92}\ub{\stag{838-3r.92} Discussion}:  We can repeat \scite{838-3r.67} -
\scite{838-3r.91} with a game version.  That is we replace Definition
\scite{838-3r.71} by \scite{838-3r.92C} and Definition \scite{838-3r.19} 
and then can immitate \scite{838-3r.88}, \scite{838-3r.91} in
\scite{838-3r.92K}, \scite{838-3r.92P}.
\bigskip

\definition{\stag{838-3r.92C} Definition}  Let $\xi \le \partial +1$.
\nl
1) We say that ${\frak u}$ has the
$S$-semi $\xi$-uq-invariant coding$_2$ property, [above $(\bar
M^*,\bar{\bold J}^*,\bold f^*) \in K^{\text{qt}}_{\frak u}$] \ub{when}
$\{0,2\}$-almost every $(\bar M,\bar{\bold J},\bold f) \in 
K^{\text{qt}}_{\frak u}$ [above $(\bar M^*,\bar{\bold J}^*,\bold
f^*)$] has it, see below; if $S = \partial$ we may omit it.
\nl
2) We say that $(\bar M,\bar{\bold J},\bold f) \in
 K^{\text{qt}}_{\frak u}$ has the $S$-semi $\xi$-invariant
coding$_2$ property \ub{when} we can choose $\langle \bold
 d_\delta:\delta \in S\rangle$ such that
\mr
\item "{$\circledast$}"  $(a) \quad \bold d_\delta$ is a ${\frak
 u}$-free $(\alpha(\bold d_\delta),0)$-rectangle
\sn
\item "{${{}}$}"  $(b) \quad M^{\bold d_\delta}_{0,0} = M_\delta$
\sn
\item "{${{}}$}"  $(c) \quad M^{\bold d_\delta}_{\alpha(\bold
d_\delta),0} \le_{\frak u} M_\partial$
\sn
\item "{${{}}$}"  $(d) \quad$ in the following game the player Coder
has a winning strategy; the
\nl

\hskip25pt  game is defined as in Definition \scite{838-2b.6F}(2) 
except that the deciding
\nl

\hskip25pt  who wins a play, i.e. we replace $(*)_2$ by
{\roster
\itemitem{ $(*)''_2$ }  in the end of the play the player Coder wins
the play \ub{when}:
\nl

\quad for a club of $\delta \in \partial$ if $\delta \in S$ \ub{then} there
are $N'',\bold I''$ such that
\nl

\quad $(M_\delta,N_\delta,\bold I_\alpha)
\le_{\frak u} (M^{\bold d_\delta}_{\alpha(\bold d_\delta),0},N'',\bold
I'')$ and for no $(\bold e,f,N_*)$ do we have
\nl

\quad the parallel of $(\alpha) - (\eta)$ of 
clause (c) of Definition \scite{838-3r.71}.
\endroster}
\ermn
3) We define when ${\frak u}$ or $(\bar M,\bar{\bold J},\bold f)$ has
   the $S$-vertical $\xi$-uq-invariant coding$_2$ property as in parts
   (1),(2) replacing \scite{838-3r.71} by \scite{838-3r.19}.
\enddefinition
\bigskip

\proclaim{\stag{838-3r.92K} Claim}  Like \scite{838-3r.88} using Definition
\scite{838-3r.92C}.
\endproclaim
\bigskip

\proclaim{\stag{838-3r.92P} Theorem}  Like \scite{838-3r.91} using
\scite{838-3r.92C}.
\endproclaim
\bn
\centerline{$* \qquad * \qquad *$}
\bn
We now point out some variants of the construction framework, 
here amalgamation may fail (unlike
\scite{838-1a.5}(1) used in \scite{838-2b.11}(4) but not usually).  This
relates to semi a.e.c.
\definition{\stag{838-3r.94} Definition}  We define when ${\frak u}$ is a
weak nice construction framework as in Definition \scite{838-1a.3} but we
considerably weaken the demands of ${\frak K}_{\frak u}$ being an
a.e.c.
\mr
\item "{$(A)$}"  ${\frak u}$ consists of $\partial,\tau_{\frak
u},{\frak K}_{\frak u} = (K_{\frak u},\le_{\frak u})$, FR$_1$,
FR$_2,\le_1,\le_2$ (also denoted by FR$^{\frak u}_1$, FR$^{\frak
u}_2,\le^1_{\frak u},\le^2_{\frak u}$) 
\sn
\item "{$(B)$}"  $\partial$ is regular uncountable
\sn
\item "{$(C)$}"  $(a) \quad \tau_{\frak u}$ is a vocabulary
\sn
\item "{${{}}$}"  $(b) \quad K_{\frak u}$ is a non-empty class of
$\tau_{\frak u}$-models of cardinality $< \partial$ closed under
\nl

\hskip25pt  isomorphisms (but $K^{{\frak u},*}_\partial,K^{{\frak
u},*}_{\partial^+}$ are defined below)
\sn
\item "{${{}}$}"  $(c) \quad$ (restricted union) if $\ell \in
\{1,2\}$ then
\sn
\item "{${{}}$}"  $\qquad (\alpha) \quad \le_{\frak u}$ is a 
partial order on $K_{\frak u}$,
\sn
\item "{${{}}$}"  $\qquad (\beta) \quad \le_{\frak u}$ is closed
under isomorphism
\sn
\item "{${{}}$}"  $\qquad (\gamma) \quad$ restricted union existence:
if $\ell=1,2$ and
\nl

\hskip35pt $\langle M_\alpha:\alpha \le \delta\rangle$ is 
$\le_{\frak K}$-increasing continuous, 
$\delta$ a limit ordinal $< \partial$
\nl

\hskip35pt  and $(M_\alpha,M_{\alpha +1},\bold I_\alpha) \in 
\text{ FR}^\ell_{\frak u}$ for $\alpha < \delta$ and $\delta =
\sup\{\alpha < \delta$:
\nl

\hskip35pt $(M_\alpha,M_{\alpha +1},\bold I_\alpha) \in 
\text{ FR}^{\ell,+}_{\frak u}$ \ub{then} $M_\delta := \{M_\alpha:\alpha <
\delta\}$ belongs to 
\nl

\hskip35pt  $K_{\frak u}$ and $\alpha < \delta \Rightarrow
M_\alpha \le_{\frak u} M_\delta$
\sn
\item "{${{}}$}"  $(d) \quad$ restricted smoothness: in clause (c)$(\gamma)$ if
$\alpha < \delta \Rightarrow M_\alpha \le_{\frak u} N$ \ub{then}
\nl

\hskip25pt $M_\delta \le_{\frak u} N$
\sn
\item "{$(D)_\ell$}"  as in Definition \scite{838-1a.3}
\sn
\item "{$(E)_\ell$}"  as in Definition \scite{838-1a.3} but we replace
clause (c) by
\sn
\item "{${{}}$}"  $(c)' \quad ((M_\delta,N_\delta,\bold J_\delta) \in
\text{ FR}^\ell_{\frak u}$ and $\alpha < \delta \Rightarrow
(M_\alpha,N_\alpha,\bold J_\alpha) \le^\ell_{\frak u}
(M_\delta,N_\delta,\bold J_\delta)$
\nl

\hskip25pt  \ub{when}:
{\roster
\itemitem{ ${{}}$ }  $(\alpha) \quad \delta < \partial$ is a limit ordinal
\sn
\itemitem{ ${{}}$ }   $(\beta) \quad \langle(M_\alpha,N_\alpha,\bold
J_\alpha):\alpha < \delta\rangle$ is $\le^\ell_{\frak u}$-increasing
continuous
\sn
\itemitem{ ${{}}$ }  $(\gamma) \quad 
(M_\delta,N_\delta,\bold J_\delta) = (\cup\{M_\alpha:
\alpha < \delta\},\cup\{N_\alpha:\alpha < \delta\},
\cup\{\bold J_\alpha:\alpha < \delta\})$
\sn
\itemitem{ ${{}}$ }  $(\delta) \quad M_\delta$ is a $\le_{\frak u}$-upper bound
of $\langle M_\alpha:\alpha \ge \delta\rangle$
\sn
\itemitem{ ${{}}$ }  $(\varepsilon) \quad N_\delta$ is a $\le_{\frak u}$-upper
bound of $\langle N_\alpha:\alpha < \delta\rangle$
\sn
\itemitem{ ${{}}$ }  $(\zeta) \quad M_\delta \le_{\frak u} N_\delta$
\endroster}
\item "{$(F)$}"  As in Definition \scite{838-1a.3}.
\endroster
\enddefinition
\bigskip

\remark{\stag{838-3r.96} Remark}  1) We may in condition (E)$_\ell$(c) use
essentially a ${\frak u}$-free $(\delta,2)$-rectangle (or
$(2,\delta)$-rectangle).
\nl
2) A stronger version of (E)$_\ell$(c) is: (E)$_2$(c)$^+$ as in
   (E)$_1$(c)$'$ adding:
\mr
\item "{$(\eta)$}"  $(M_\alpha,M_{\alpha +1},\bold I_\alpha)
\le^1_{\frak u} (N_\alpha,N_{\alpha +1},\bold I^1_\alpha)$
\sn
\item "{$(\theta)$}"  $(M_\alpha,N_\alpha,\bold I^m_\alpha) \in \text{
FR}^{1,+}_{\frak u}$ for $m=0,1$ for unboundedly many $\alpha <
\delta$ (so clause $(\zeta)$ follows).
\ermn
(E)$_1$(c)$^+$ means clause (E)$_2$(c)$^+$ is satisfied by 
{\rm dual}$({\frak u})$.
\nl
3) We may demand for $(\bar M,\bar{\bold J},\bold f) \in
K^{\text{qt}}_{\frak u}$ that for a club of $\delta$, if $\bold
f(\delta) > 0$ then:
\mr
\item "{$(a)$}"  $\bold f(\delta)$ is a limit ordinal
\sn
\item "{$(b)$}"  $\bold f(\delta) = \sup\{i < \bold
f(\delta):(M_{\delta +i},M_{\delta +i+1},\bold J_{\delta +i}) \in
\text{ FR}^{2,+}_{\frak u}\}$
\sn
\item "{$(c)$}"  in the examples coming for an almost good
$\lambda$-frame ${\frak s}$, see \S5: if $i < \delta$ and $p \in {\Cal
S}^{\text{bs}}_{\frak s}(M_{\delta +i})$ then $\delta = \sup\{j:j \in
(i,\delta)$ and $\bold J_{\delta +j} = \{a_{\delta +j}\}$ and
$\ortp_{\frak s}(a_{\delta +j},M_{\delta +j},M_{\delta +j+1})$ is a
non-forking extension of $p$; see more in \S5 on this.
\ermn
4) We may in part (3)(b) and in \scite{838-3r.98}(A)(a) below restrict
   ourselves to sucessor $i$.
\endremark
\bigskip

\proclaim{\stag{838-3r.98} Lemma}  We can repeat \S1 + \S2 (and \S3) with Definition
\scite{838-3r.94} instead of Definition \scite{838-1a.3} with the following
changes:
\mr
\item "{$(A)$}"  from Definition \scite{838-1a.29}:
{\roster
\itemitem{ $(a)$ }  in the Definition of $(\bar M,\bar{\bold J},\bold
f) \in K^{\text{qt}}_{\frak u}$ we demand $S = \{\delta <
\partial:\bold f(\delta) > 0\}$ is stationary and for a club of
$\delta \in S,\bold f(\delta)$ is a limit ordinal and $i < \bold
f(\delta) \Rightarrow (M_{\delta +i},M_{\delta +i+1},\bold J_{\delta
+i}) \in \text{\rm FR}^{2,+}_{\frak u}$
\sn
\itemitem{ $(b)$ }  we redefine $\le^{\text{qr}}_{\frak u},
\le^{\text{qs}}_{\frak u}$ as $\le^{\text{qt}}_{\frak u}$ was redefined
\endroster}
\item "{$(B)$}"  proving \scite{838-1a.37}:
{\roster
\itemitem{ $(a)$ }  in part (1), we should be given a stationary $S
\subseteq \partial$ and for $\alpha < \partial$ let $\bold f(\alpha) =
\omega$
\sn
\itemitem{ $(b)$ }  in part (2), we use the restricted version of
union existence  and smoothness
\sn
\itemitem{ $(c)$ }  in part (3), we demand
$(M^1_\alpha,M^2_\alpha,\bold I^*) \in \text{\rm FR}^+_1$ and let $E \subseteq
\lambda \backslash \alpha$ witness $(\bar M,\bar{\bold J},\bold f) \in
K^{\text{qt}}_{\frak u}$ and $S = \{\delta \in E:\bold f(\delta) >
0\}$ and inthe induction we use just $\langle M_\beta:\beta\in
\{\alpha\} \cup \bigcup\{[\delta,\delta,\bold f(\delta)]:\delta \in
S\}\rangle$. 
\endroster}
\endroster
\endproclaim
\bn
\margintag{3r.99}\ub{\stag{838-3r.99} Exercise}:  Rephrase this section with
$\tau$-coding$_k$ instead coding$_k$.
\sn
[Hint: 1) Of course, we replace ``coding" by ``$\tau$-coding" and
isomorphic by $\tau$-isomorphic.
\nl
2) We replace $\bullet_5$ of \scite{838-3r.1}(1)$(\delta)$ by:
\mr
\item "{$\bullet'_5$}"  there are $N_1,N_2$ such that $N^{\bold
d_\ell}_{\alpha(\bold d_\ell),\beta(\bold d_\ell)} \le_{\frak u}
N_\ell$ for $\ell=1,2$ and there is a $\tau$-isomorphism $f$ from $N_1$
onto $N_2$ extending id$_{M^{\bold d}_{\alpha',0}}$ and mapping
$M^{\bold d_1}_{0,\beta(\bold d_1)}$ onto $M^{\bold
d_2}_{0,\beta(\bold d_2)}$.
\ermn
3) In Claim \scite{838-3r.3} in the end of clause $(\gamma)$
``$N^1_\partial,N^2_\partial$ are not $\tau$-isomorphic over
$M_\gamma$", of course.
\nl
4) In \scite{838-3r.7} replace $\dot I$ by $\dot I_\tau$.
\nl
5) In \scite{838-3r.9}, in the conclusion replace $\dot I$ by $\dot
   I_\tau$, in clause (c)$(\beta)$ use ``not $\tau$-isomorphic".
\nl
6) In Definition  \scite{838-3r.11} like (2), in Definition
   \scite{838-3r.19}(e)$(\delta)$ as in (2).
\nl
7) Change \scite{838-3r.71}(c)$(\varepsilon)$ as in (2), i.e.
\mr
\item "{$(\varepsilon)'$}"  $f$ is a $\tau$-isomorphism from $N^*$
onto $N_*$ for some $N^*$ such that $N'' \le_{\frak u} N^*$.]
\endroster
\goodbreak

\head{\S4 Straight Applications of (weak) coding} \endhead  \resetall \sectno=4
 \spuriousreset
\bigskip

Here, to try to exemplify the usefulness of Theorem \scite{838-2b.3}, the
``lean" version, i.e. using weak coding, we revisit 
older non-structure results.  First, recall
that the aim of
\cite{Sh:603} or better \sectioncite[\S3,\S4]{E46} is 
to show that the set of minimal types in ${\Cal S}^{\text{na}}_{\frak K}(M),
M \in {\frak K}_\lambda$ is dense, \ub{when}:
\mr
\item "{$\boxtimes^\lambda_{\frak K}$}"  $2^\lambda < 2^{\lambda^+} <
2^{\lambda^{++}},{\frak K}$ is an a.e.c. with LS$({\frak K}) \le
\lambda$, categorical in $\lambda,\lambda^+$ and have a medium number
of models in $\lambda^{++}$ (hence ${\frak K}$ has amalgamation in
$\lambda$ and in $\lambda^+$).
\ermn
More specifically we have to justify Claim \yCITE[3c.34]{E46} 
when the weak diamond ideal on $\lambda^+$ is not
$\lambda^{++}$-saturated and we have to justify claim 
\yCITE[4d.17]{E46} when some $M \in K_{\lambda^+}$ is saturated; in 
both cases inside the proof there we quote results from here.

We interpret medium as $1 \le \dot I(\lambda^{++},K) <
\mu_{\text{unif}}(\lambda^{++},2^{\lambda^+})$ (where the latter is
usually $2^{\lambda^{++}}$).  This is done in \scite{838-e.1}-\scite{838-e.2T}, 
i.e., where we prove the non-structure parts relying on the one hand on
the pure model theoretic part done in \chaptercite{E46} and on the other hand
on coding theorems from \S2.  
More elaborately, as we are relying on Theorem \scite{838-2b.3}, in \scite{838-e.1} -
\scite{838-e.1.7}, i.e. \S4(A) we assume that the normal ideal
WDmId$(\lambda^+)$ is not $\lambda^{++}$-saturated and prove for appropriate
${\frak u}$ that (it is a nice construction framework and) it has
 the weak coding property.  Then in \scite{838-e.2} - \scite{838-e.2T},
i.e. \S4(B) relying on Theorem \scite{838-2b.13},  
we assume more model theory and (for the appropriate ${\frak u}$) 
prove the vertical coding property, hence eliminate
the extra set theoretic assumption (but retaining the relevant cases
of the WGCH, i.e.  $2^\lambda < 2^{\lambda^+} < 2^{\lambda^{++}}$).

Second, we relook at the results in \sectioncite[\S6]{E46},
i.e. \cite[\S6]{Sh:576} which were originally 
proved relying on \cite[\S3]{Sh:576}.  That is, our aim
is to prove the density of uniqueness triples $(M,N,a)$ in
$K^{3,\text{na}}_\lambda$, assuming medium
number of models in $\lambda^{++}$, and set theoretically $2^\lambda <
2^{\lambda^+} < 2^{\lambda^{++}}$ and in addition assume (for now) the
non-$\lambda^{++}$-saturation of the
weak diamond ideal on $\lambda^+$.  So we use the ``weak coding"
from Definition \scite{838-2b.1}, Theorem \scite{838-2b.3} 
(see \scite{838-e.4}, i.e. \S4(D)).  The elimination of the
extra assumption is delayed as it is more involved (similarly for
\S4(E)).

Third, we fulfill the (``lean" version of the) promise from
\sectioncite[\S5]{600}, proving density of uniqueness triples in
$K^{3,\text{bs}}_{\frak s}$, for ${\frak s}$ a good $\lambda$-frame,
also originally relying on \cite[\S3]{Sh:576}, see \scite{838-e.5}, i.e. \S4(E).

Fourth, we deal with the promises from
\sectioncite[\S5]{88r} by Theorem \scite{838-2b.3} in \scite{838-e.3} -
\scite{838-e.3.12}, i.e. \S4(C).

But still we owe the ``full version", this is \S4(F) in which 
we eliminate the extra set theoretic result relying
on the model theory from \S5-\S8.
\bn
\centerline{$* \qquad * \qquad *$}
\bn
\ub{(A) \quad Density of the minimal types for ${\frak K}_\lambda$}

\proclaim{\stag{838-e.1} Theorem}    We have $\dot I(\lambda^{++},K) \ge
\mu_{\text{unif}}(\lambda^{++},2^{\lambda^+})$ \ub{when}:
\mr
\item "{$\odot$}"   $(a) \quad 2^\lambda < 2^{\lambda^+} <
2^{\lambda^{++}}$
\sn
\item "{${{}}$}"  $(b) \quad$ the ideal {\rm WDmId}$(\lambda^+)$, a
normal ideal on $\lambda^+$, is not $\lambda^{++}$-saturated 
\sn
\item "{${{}}$}"  $(c) \quad {\frak K}$, an a.e.c. with 
{\rm LS}$({\frak K}) \le \lambda$, has amalgamation in $\lambda$, the
{\rm JEP} 
\nl

\hskip25pt  in $\lambda$, for simplicity and $K_{\lambda^+} \ne
\emptyset$;
\sn
\item "{${{}}$}"  $(d) \quad$  for every $M \in K_{\lambda^+}$ and
$\le_{\frak K}$-representation $\bar M = \langle M_\alpha:\alpha <
\lambda^+\rangle$ of 
\nl

\hskip25pt $M$ we can find $(\alpha_0,N_0,a)$, i.e., a triple
$(M_{\alpha_0},N_0,a)$  such that: 
{\roster
\itemitem{ ${{}}$ }   $(\alpha) \quad 
M_{\alpha_0} \le _{{\frak K}_\lambda} N_0$
\sn
\itemitem{ ${{}}$ }  $(\beta) \quad a \in N_0 \backslash M_{\alpha_0}$ and
$\ortp_{\frak K}(a,M_{\alpha_0},N_0)$ is not realized in $M$ 
\sn
\itemitem{ ${{}}$ }  $(\gamma) \quad$ if $\alpha_0 < \alpha_1 <
\lambda^+,M_{\alpha_1} \le_{\frak K} N_1$ and $f$ is a $\le_{\frak
K}$-embedding of 
\nl

\hskip40pt $N_0$ into $N_1$ over $M_{\alpha_0}$ \ub{then} we can find
$\alpha_2  \in (\alpha_1,\lambda^+)$ 
\nl

\hskip40pt such that $M_{\alpha_2},N_1$ are not uniquely amalgamated over
\nl

\hskip40pt  $M_{\alpha_1}$ (in ${\frak K}_\lambda$), i.e.
$\text{\rm NUQ}_\lambda(M_{\alpha_1},M_{\alpha_2},N_1)$, see
\yCITE[3c.4]{E46}(2). 
\endroster}
\endroster
\endproclaim
\bigskip

\remark{\stag{838-e.1H} Remark}  1) Used in Claim \yCITE[3c.34]{E46}, more
exactly the relative \scite{838-e.1Z} is used.
\nl
2) A further question, mentioned in \yCITE[2b.33]{E46}(3) concern $\dot I
\dot E(\lambda^{++},K)$ but we do not deal with it here.
\nl
3) Recall that for $M \in K_\lambda$ the type $p \in 
{\Cal S}^{\text{na}}_{{\frak K}_\lambda}(M)$ is minimal \ub{when} 
there is no $M_2 \in K_\lambda$ which $\le_{\frak K}$-extends 
$M$ and $p$ has at least 2
extensions in ${\Cal S}^{\text{na}}_{{\frak K}_\lambda}(M_2)$; see Definition
\yCITE[1a.34]{E46}.
\nl
3A) We say that in ${\frak K}_\lambda$ the minimal types are dense
when: for any $M \in K_\lambda$ and $p \in 
{\Cal S}^{\text{na}}_{{\frak K}_\lambda}(M)$ there is a pair $(N,q)$
such that $M \le_{{\frak K}_\lambda} N$ and $q \in 
{\Cal S}^{\text{na}}_{{\frak K}_\lambda}(N)$ is minimal and extend $p$
(see \yCITE[1a.34]{E46}(1A)).
\nl
4) A weaker version of Clause (d) of \scite{838-e.1} holds when 
any $M \in K_{\lambda^+}$ is saturated (above $\lambda$) and the
minimal types are not dense (i.e. omit subclause $(\beta)$ and in
subclause $(\gamma)$
add $f(a) \notin M_{\alpha_1}$; the proof is similar (but using
\scite{838-e.2K}).  Actually, \scite{838-e.1} as phrased is useful normally
only when $2^\lambda > \lambda^+$, but otherwise we use \scite{838-e.1Z}.
\nl
5) In \yCITE[3c.34]{E46} we work more to justify
a weaker version (d)$''$ of Lemma \scite{838-e.1Z} below which suffice.
\endremark
\bn
Similarly
\proclaim{\stag{838-e.1Z} Lemma}  1) Like \scite{838-e.1} but we replace clause (d)
by:
\mr
\item "{$(d)'$}"  there is a superlimit $M \in K_{\lambda^+}$ and for
it clause (d) of \scite{838-e.1} holds.
\ermn
2) For $\tau$ a $K$-sub-vocabulary, see \scite{838-1a.15}(5), we have
$\dot I_\tau(\lambda^{++},K) \ge
\mu_{\text{unif}}(\lambda^{++},2^{\lambda^+})$ when (a),(b),(c) of
\scite{838-e.1} holds and
\mr
\item "{$(d)''$}"   there\footnote{Why not $K'_\lambda$?  Just because we
use this notation in \scite{838-1a.19}.} is $K''_{\lambda^+}$ such that
{\roster
\itemitem{ $(\alpha)$ }  $K''_{\lambda^+} \subseteq
K_{\lambda^+},K''_{\lambda^+} \ne \emptyset$ and $K''_{\lambda^+}$ is
closed under unions of $\le_{\frak K}$-increasing 
continuous sequences of length $\le \lambda^+$
\sn
\itemitem{ $(\beta)$ }   there is $M_* \in K''_{\lambda^+}$ such that for
any $\le_{{\frak K}_\lambda}$-increasing continuous sequence
$\langle M^1_\alpha:\alpha < \lambda^+\rangle$ with union $M^1 \in
K''_{\lambda^+}$ satisfying $M_* \le_{\frak K} M^1$, in the
following game the even player has a winning strategy.
In the $\alpha$-th move a triple $(\beta_\alpha,M_\alpha,f_\alpha)$ is
chosen such that $\beta_\alpha < \lambda^+,
M_\alpha \in {\frak K}_\lambda$ has universe $\subseteq
\lambda^+$ and $f_\alpha$ is a $\le_{\frak K}$-embedding of
$M^1_{\beta_\alpha}$ into $M_\alpha$, all three are
increasing continuous with $\alpha$.
Of course, the $\alpha$-th move is done by the even/odd iff
$\alpha$ is even/odd.
Lastly, in the end the even player wins iff $\cup\{M_\alpha:\alpha <
\lambda^+\}$ belongs to $K''_{\lambda^+}$
 and for a club of $\alpha < \lambda^+$ for some $\gamma \in
(\beta_\alpha,\lambda^+)$ and some $N \in K_\lambda$ such that
there is an isomorphism from $N$ onto $M_\alpha$ extending $f_\alpha$
we have {\rm NUQ}$_\tau(M^1_{\beta_\alpha},N,M^1_\gamma)$,
i.e. $N,M_\gamma$ can be amalgamated over $M_{\beta_\alpha}$ in
${\frak K}_\lambda$ in at least two $\tau$-incompatible ways.
\endroster}
\endroster
\endproclaim
\bigskip

\remark{\stag{838-e.1R} Remark}  0) We may replace $\beta_\alpha$ in clause
(d)$''(\beta)$ above by $\beta = \alpha +1$, many times it does not matter.
\nl
1) We may weaken the model theoretic assumption $(d)''$ of
   \scite{838-e.1Z}(2) so weaken (d) in \scite{838-e.1} and $(d)'$ in \scite{838-e.1Z} 
if we strengthen the set theoretic assumptions, e.g.
\mr
\item "{$(*)_1$}"  for some stationary $S \subseteq
S^{\lambda^{++}}_{\lambda^+}$ we have $S \in \check I[\lambda^{++}]$ 
but $S \notin \text{ WDmId}(\lambda^{++})$
\sn
\item "{$(*)_2$}"  in \scite{838-e.1Z}, in subclause $(\alpha)$ of clause
(d)$''$ we weaken the closure under union of $K''_{\lambda^+}$ to: for
a $\le_{\frak K}$-increasing sequence in $K''_{\lambda^+}$ of length
$\lambda^+$, its union belongs to $K''_{\lambda^+}$.
\ermn
2) If, e.g. $\lambda = \lambda^{<\lambda}$ and $\bold V = \bold V^{\Bbb Q}$
where $\Bbb Q$ is the forcing notion of 
adding $\lambda^+$-Cohen subsets of $\lambda$ 
and the minimal types are not dense \ub{then} 
$\dot I(\lambda^+,K) = 2^{\lambda^+}$, (hopefully see more in \cite{Sh:E45}).
\nl
3) If in part (2) of \scite{838-e.1Z}, we may consider omitting 
the amalgamation, but demand ``no maximal model in ${\frak K}_\lambda$".  
However, the minimality may hold for uninteresting reasons.
\nl
4)  This is used in \yCITE[3c.34]{E46}.
\nl
5) We may assume $2^\lambda \ne \lambda^+$ as essentially the case
$2^\lambda = \lambda^+$ is covered 
\footnote{by the assumption ${\frak K}_\lambda$ has amalgamation so
togther with $2^\lambda = \lambda^+$, if $|\tau_K| \le \lambda$ or
 just $M \in {\frak K}_\lambda \Rightarrow |{\Cal S}_{{\frak
 K}_\lambda}(M)| \le \lambda^+$ \ub{then} there is a model $M^* \in {\frak
K}_{\lambda^+}$ which is saturated above $\lambda$; in general this is
 how it is used in \sectioncite[\S4]{E46}; but if there is $M \in 
{\frak K}_\lambda$ with $|{\Cal S}_{{\frak K}_\lambda}(M)| >
 \lambda^+$; we can use \scite{838-e.1}.} by Lemma \scite{838-e.2} below. 
\nl
6) In clause (d)$''(\beta)$ of Lemma \scite{838-e.1Z}, we can let the
 even player choose also for $\alpha = \delta +1$ for $\delta \in S$
   \ub{when} $S \subseteq \lambda^+$ but WDmId$(\lambda^+) +S$ is not
   $\lambda^{++}$-saturated. 
\endremark
\bigskip

\demo{Proof of \scite{838-e.1}}  We shall apply Theorem \scite{838-2b.3}.
So (model theoretically) we have an a.e.c. ${\frak K}$  
with LS$({\frak K}) \le \lambda$, and we are interested in proving
$\dot I(\lambda^{++},
{\frak K}) \ge \mu_{\text{unif}}(\lambda^{++},2^{\lambda^+})$.

We shall define (in Definition \scite{838-e.1A} below) a nice construction 
framework ${\frak u}$ such that $\partial_{\frak u} = \lambda^+$; 
the set theoretic assumptions of \scite{838-2b.3} hold; i.e.
\mr
\item "{$(a)$}"  $\lambda < \partial$ and $2^\lambda = 2^{< \partial} <
2^\partial$; i.e. we choose $\theta := \lambda$ and this holds 
by clause $\odot(a)$ of Theorem \scite{838-e.1}
\sn
\item "{$(b)$}"  $2^\partial < 2^{\partial^+}$; holds by clause $\odot(a)$ of
\scite{838-e.1} 
\sn
\item "{$(c)$}"  the ideal WDmID$(\partial)$ is not
$\partial^+$-saturated; holds by clause $\odot(b)$ of \scite{838-e.1}
\ermn
We still have to find ${\frak u}$ (and $\tau$) as required in clause
(d) of Theorem \scite{838-2b.3}.  We define it in Definition \scite{838-e.1A}
below, in particular we let ${\frak K}_{\frak u} = {\frak K}'_\lambda,\tau =
\tau_{\frak K}$, see
Definition \scite{838-1a.19} where ${\frak K}$ is the a.e.c. from \scite{838-e.1Z} 
hence the conclusion ``$\dot I_\tau(\partial^+,K^{{\frak u},{\frak
h}}_{\partial^+}) \ge \mu_{\text{unif}}(\partial^+,2^\partial) > 
2^\partial$ for any $\{0,2\}$-appropriate function ${\frak h}$" of
Theorem \scite{838-2b.3} implies that $\dot I(\lambda^{++},K^{\frak u}) \ge
\mu_{\text{unif}}(\lambda^{++},2^{\lambda^+})$ as required using
\scite{838-e.1C}(3); we can use Exercise \scite{838-1a.21}.  
So what we should actually prove is that we can find such nice
construction frameworks ${\frak u}$ with the weak coding 
property which follows from ${\frak
u}$ having the weak coding property (by \scite{838-2b.11}(2),(3)) and
${\frak h}$ such that every $M \in K^{{\frak u},{\frak
h}}_{\lambda^{++}}$ is $\tau$-fuller.  
This is done in \scite{838-e.1B}, \scite{838-e.1C} below.
\hfill$\square_{\scite{838-e.1}}$ 
\enddemo
\bigskip

\demo{Proof of \scite{838-e.1Z}}  1) By part (2), in particular letting
$K''_{\lambda^+} = \{M \in K_{\lambda^+}:M$ is superlimit$\}$.

But why does subclause $(\beta)$ of clause (d)$''$ hold?  Let $M_*
\in K''_{\lambda^+}$ be superlimit, let $\langle M^1_\alpha:\alpha <
\lambda^+\rangle$ be $\le_{{\frak K}_\lambda}$-increasing continuous
with union $M^1 \in K''_{\lambda^+}$ 
and assume $M_* \le_{\frak K} M^1$.  Without loss of
generality $M^1$ has universe $\{\beta < \lambda^+:\beta$ odd$\}$.

We shall prove that the even player has a winning strategy.  We
describe it as follows: the even player in the $\alpha$-th move also
choose $(N_\alpha,g_\alpha)$ for $\alpha$ even and also for $\alpha$
odd (after the odd's move) such that
\mr
\item "{$(*)$}"  $(a) \quad N_\alpha \in K_{\lambda^+}$ is superlimit
and $M^1 \le_{\frak K} N_0$
\sn
\item "{${{}}$}"  $(b) \quad$ the universe of $N_\alpha$ is $\{\gamma
< \lambda^+:\gamma$ is not divisible by $\lambda^\alpha\}$
\sn
\item "{${{}}$}"  $(c) \quad N_\beta \le_{\frak K} N_\alpha$ for
$\beta < \alpha$
\sn
\item "{${{}}$}"  $(d) \quad g_\alpha$ is a $\le_{\frak K}$-embedding
of $M_\alpha$ into $N_\alpha$
\sn
\item "{${{}}$}"  $(e) \quad g_\beta \subseteq g_\alpha$ for $\beta <
\alpha$
\sn
\item "{${{}}$}"  $(f) \quad$ if $\alpha$ is odd then the universe of
$g_{\alpha +1}(M_{\alpha +1})$ includes $\alpha = N_\alpha \cap \alpha$.
\ermn
It should be clear that the even player can do this.  Also for any such
play $\langle(M_\alpha,f_\alpha,N_\alpha,g_\alpha):\alpha <
\lambda^+\rangle$ we have $\lambda^+ = \cup\{N_\alpha \cap
\alpha:\alpha < \lambda^+\} \subseteq \cup\{g_\alpha(M_\alpha):\alpha
< \lambda^+\} \subseteq \cup\{N_\alpha:\alpha < \lambda^+\} \subseteq
\lambda^+$, so $g = \cup\{g_\alpha:\alpha < \lambda^+\}$ is an
isomorphism from $\cup\{M_\alpha:\alpha < \lambda^+\}$ onto
$\cup\{N_\alpha:\alpha < \lambda^+\}$.  As the latter is superlimit we
are done (so for this part being $(\lambda^+,\lambda^+)$-superlimit suffice).
\nl
2) The proof is like the proof of \scite{838-e.1} but we have to use a
variant of \scite{838-2b.3}, i.e. we use the variant of weak coding where we use a
game, see \scite{838-2b.6}.  \hfill$\square_{\scite{838-e.1Z}}$
\enddemo
\bigskip

\definition{\stag{838-e.1A} Definition}   [Assume clause (c) of
\scite{838-e.1}.]

We define ${\frak u} = {\frak u}^1_{{\frak K}_\lambda}$ as follows 
(with $\tau({\frak u}) = (\tau_{\frak K})'$ so $=_\tau$ is
a congruence relation in $\tau({\frak u})$, O.K. by \scite{838-1a.19};
this is a fake equality \scite{838-3r.84}(1), \scite{838-3r.86})
\mr
\item "{$(a)$}"  $\partial_{\frak u} = \lambda^+$
\sn
\item "{$(b)$}"  essentially ${\frak K}_{\frak u} = {\frak
K}_\lambda$; really ${\frak K}'_\lambda$ (i.e. $=_\tau$ is a 
congruence relation!)
\sn
\item "{$(c)$}" $\text{\rm FR}_1 = \{(M,N,\bold I):
M \le_{{\frak K}_{\frak u}} N,\bold I \subseteq N \backslash M$ 
empty or a singleton $\{a\}\}$
\sn
\item "{$(d)$}"  $\text{FR}_2 = \text{ FR}_1$
\sn
\item "{$(e)$}"  $(M_1,N_1,\bold J_1) \le_\ell (M_2,N_2,\bold J_2)$ \ub{when}
{\roster
\itemitem{ $(i)$ }  both triples are from FR$_\ell$ 
\sn
\itemitem{ $(ii)$ }  $M_1 \le_{\frak K} M_2,N_1 \le_{\frak K} N_2$ 
and $\bold J_1 \subseteq \bold J_2$
\sn
\itemitem{ $(iii)$ }  $M_2 \cap N_1 = M_1$.
\endroster}
\endroster
\enddefinition
\bigskip

\demo{\stag{838-e.1B} Observation}   ${\frak u}$ is a nice construction
framework which is self-dual.
\enddemo
\bigskip

\demo{Proof}  Easy and the proof of \scite{838-e.2M} can
serve when we note that (D)$_1$(d) and (F) are obvious in our context,
 recalling we have fake equality.   
 \hfill$\square_{\scite{838-e.1B}}$
\enddemo
\bigskip

\demo{\stag{838-e.1C} Observation}  [Assume $\lambda,{\frak K}$ are as in
 \scite{838-e.1} or \scite{838-e.1Z}(1) or \scite{838-e.1Z}(2).]
\nl
1) ${\frak u}$ has the weak coding property (see Definition \scite{838-2b.1}).
\nl
2) If $(\bar M,\bar{\bold J},\bold f) \in K^{\text{qt}}_{\frak u}$ and so
$M_\partial \in K_{\lambda^+}$ but for \scite{838-e.1Z}(2) the model
$M_\partial$ is superlimit \ub{then} $(\bar M,\bar{\bold J},\bold f)$
has weak coding property.
\nl
3) For some ${\frak u}-0$-appropriate function ${\frak h}$,  every
$M \in K^{{\frak u},{\frak h}}_{\lambda^{++}}$ is 
$\tau_{\frak s}$-fuller, i.e. the model $M/=^M$ has cardinality $\lambda^{++}$
and the set $a/=^M$ has cardinality $\lambda^{++}$ for every $a \in M$.
\enddemo
\bigskip

\demo{Proof}  Easy.
\nl
1) By part (2) and (3).
\nl
2) For proving \scite{838-e.1} by clause (d) there 
we choose $(\alpha(0),N_0,\bold I_0)$, it is as
   required in Definition \scite{838-2b.1}(3), 
noting that we can get the necessary disjointness
because ${\frak K}'_\lambda$ has fake equality.  Similarly for
\scite{838-e.1Z}(1). 

For proving \scite{838-e.1Z}(2) we fix a winning strategy {\bf st} for the
even player in the game from clause (d) of \scite{838-e.1Z}.  Again by the
fake equality during the game we can demand $\alpha_1 < \alpha
\Rightarrow f_\alpha(M_{\beta_\alpha}) \cap M_\alpha =
f_{\alpha_1}(M_{\beta_1})$. 
\nl
3) By \scite{838-1a.49}(3) is suffice to deal separately with each aspect
   of being $\tau_{\frak s}$-fuller.

First, we choose a ${\frak u}$-0-appropriate function
${\frak h}_0$ such that if $((\bar M^1,\bar{\bold J}^1,
\bold f^1),(\bar M^2,\bar{\bold J}^2,\bold f^2))$ does 0-obeys
${\frak h}_0$ as witnessed by $(E,\bar{\bold I})$ \ub{then} for any
$\delta \in E,(M^1_\delta,M^2_\delta,\bold I_\delta) \in 
\text{ FR}^+_1$ and 
 there is $a \in \bold I_\delta$ such that $c \in
 M^1_\partial \Rightarrow M^2_\partial \models \neg(a =_\tau c)$;
this is possible as in every $(\bar M^1,\bar{\bold J}^1,\bold f^1) \in
K^{\text{qt}}_{\frak u}$ there are $\alpha < \lambda^+$ and $p \in
{\Cal S}^{\text{na}}(M^1_\alpha)$ not realized in $M^1 =
\cup\{M^1_\beta:\beta < \alpha\}$.  Why?  For \scite{838-e.1} by
\scite{838-e.1}(d)$(\beta)$, similarly for \scite{838-e.1Z}(1) and for
\scite{838-e.1Z}(2), if it fails, then the even player cannot win, because
\mr
\item "{$(*)$}"  if $M_0 \le_{\frak u} M_\ell$ for $\ell=1,2$ and
$(\forall b \in M_1)(\exists b \in M_1)(\exists a \in M_0)(M_1 \models
a =_\tau b)$ \ub{then} $M_1,M_2$ can be uniquely disjointly amalgamated in
${\frak K}_{\frak u}$.
\ermn
Second, we choose a ${\frak u}-0$-suitable ${\frak h}_2$ such that if
$\langle(\bar M^\alpha,\bar{\bold J}^\alpha,\bold f^\alpha):\alpha <
\partial^+\rangle$ does $0$-obeys ${\frak h}_2$, then for every
$\alpha < \partial^+$ and $a \in M_\alpha$ for $\lambda^{++}$ many
$\delta \in (\alpha,\lambda^{++})$, we have $(a/=^{M_\delta}_\tau) \subset
(a/=^{M_{\delta +1}}_\tau)$.    \hfill$\square_{\scite{838-e.1C}}$
\enddemo
\bigskip

\margintag{e.1F}\ub{\stag{838-e.1F} Example}  For ${\frak K},{\frak u}$, {\bf st} as in the
proof of \scite{838-e.1C}(2).
\nl
1) For some initial segment $\bold x =
\langle(\beta_\alpha,M_\alpha,f_\alpha):\alpha \le
\alpha_*\rangle$, of a play of the game of length $\alpha_* <
   \lambda^+$ in which the even player uses the strategy {\bf st}, for
any longer such initial segment $\langle (\beta_\alpha,M_\alpha,
f_\alpha):\alpha \le \alpha_{**}\rangle$ of such a play we have
   $M_{\alpha_*} \cap f_{\alpha_{**}}(M^1_{\beta_{\alpha_{**}}}) =
   f_{\alpha_*}(M^1_{\beta_{\alpha_*}})$ and
   $f_{\alpha_*}(M^1_{\beta_{\alpha_*}}) <_{\frak K} M_\alpha$.
\nl
[Why?  As in the proof of the density of reduced triples; just think.]
\nl
2)  Moreover if $c \in M_\partial \backslash M_{\beta_{\alpha_*}}$
    then $f_{\alpha_*}(\ortp_{\frak
    K}(c,M_{\beta_{\alpha_*}},M_\partial))$ is not realized in
    $M_{\alpha_*}$ (recall the definition of NUQ).
\hfill$\square_{\scite{838-e.1}}$ 
\bigskip

\remark{\stag{838-e.1.7} Remark}  So we have finished proving \scite{838-e.1}, 
\scite{838-e.1Z}.
\endremark
\bn
\centerline {$* \qquad * \qquad *$}
\bn
\ub{(B) \quad Density of minimal types: without
$\lambda^{++}$-saturation of the ideal}
\bigskip

The following takes care of
\yCITE[4d.17]{E46}, of its assumptions, (a)-(g) are listed 
in \yCITE[4d.14]{E46}. 
\proclaim{\stag{838-e.2} Theorem}   We have $\dot I(\lambda^{++},K) \ge
\mu_{\text{unif}} (\lambda^{++},2^{\lambda^+})$ \ub{when}:
\mr
\item "{$(a)$}"  $2^\lambda < 2^{\lambda^+} < 2^{\lambda^{+2}}$
\sn
\item "{$(b)$}"  ${\frak K}$ is an abstract elementary class,
${\text{\rm LS\/}}({\frak K})
\le \lambda$ 
\sn
\item "{$(c)$}"  $K_{\lambda^{++}} \ne \emptyset$,
\sn
\item "{$(d)$}"  ${\frak K}$ has amalgamation in $\lambda$
\sn
\item "{$(e)$}"  the minimal types, for ${\frak K}_\lambda$ are not
dense, see \scite{838-e.1H}(3A)
\sn
\item "{$(f)$}"  ${\frak K}$ is categorical in $\lambda^+$
or at least has a superlimit model in $\lambda^+$
\sn
\item "{$(g)$}"  there is $M \in K_{\lambda^+}$ which is saturated (in
${\frak K}$) above $\lambda$. 
\endroster
\endproclaim
\bigskip

\demo{Proof}   The proof is broken as in the other cases.
\enddemo
\bigskip

\definition{\stag{838-e.2K} Definition}  We define ${\frak u} = {\frak
u}^2_{{\frak K}_\lambda}$ as in \scite{838-e.1A} so ${\frak K}_{\frak u} =
{\frak K}'_\lambda$ but replacing clauses (c),(e) by (we shall use the fake
equality only for having disjoint amalgamation):
\mr
\item "{$(c)'$}"  $\text{\rm FR}_1 = \{(M,N,\bold I):
M \le_{{\frak K}_{\frak u}} N,\bold I \subseteq N \backslash M$ 
empty or a singleton $\{a\}$ and if $(a/=^M) \notin M/=^M$ \ub{then}
$\ortp_{{\frak K}_\lambda}(a,M,N)$ has no
minimal extension, i.e. and $\ortp_{{\frak
K}_\lambda}((a/=^M_\tau,(M/=^M_\tau),(N/=^N_\tau))$ 
has no minimal extension$\}$
\sn
\item "{$(e)$}"  $(M_1,N_1,\bold J_1) \le_\ell (M_2,N_2,\bold J_2)$
\ub{when}: clauses (i),(ii),(iii) there and
{\roster
\itemitem{ $(iv)$ }  if $b \in \bold J_1$
and $(\forall a \in M_1)(\neg a =^{N_1}_\tau b)$ \ub{then} $(\forall a
\in M_2)(\neg a =^{N_2}_\tau b)$.
\endroster}
\endroster
\enddefinition
\bigskip

\demo{\stag{838-e.2E} Observation}  Without loss of generality ${\frak K}$
has (jep)$_\lambda$ and for every $M \in K_\lambda$ there is 
$p \in {\Cal S}^{\text{na}}_{{\frak K}_\lambda}(M)$ with no minimal extension.
\enddemo
\bigskip

\demo{Proof}  Why?

Because if $(p_*,M_*)$ are as required we can replace ${\frak K}$ by
${\frak K}^* = {\frak K} \restriction \{M \in K$: there is a
$\le_{\frak K}$-embedding of $M_*$ into $M\}$.  Clearly
${\frak K}^*$ satisfies the older requirements and if $h$ is a
$\le_{\frak K}$-embedding of $M_*$ into $M \in {\frak K}^*_\lambda$
then $h(p_*)$ can be extended to some $p \in {\Cal S}_{{\frak
K}^*_\lambda}(M) = {\Cal S}_{{\frak K}_\lambda}(M)$ as required.  Why
it can be extended?  As any triples $(M,N,a) \in
K^{3,\text{na}}_\lambda$ with no minimal extension has the extension
property, see \yCITE[2b.7]{E46}(1).  \hfill$\square_{\scite{838-e.2E}}$

The first step is to prove that (and its proof 
includes a proof of \scite{838-e.1B}).
\enddemo
\bigskip

\proclaim{\stag{838-e.2M} Claim}   ${\frak u}$ is a nice construction
framework (if it is as in \scite{838-e.2E}).
\endproclaim
\bigskip

\demo{Proof}   Clauses (A),(B),(C) of Definition \scite{838-1a.3} are
obvious.  Also (D)$_1 =$ (D)$_2$ and (E)$_1 =$ (E)$_2$ 
as FR$_1 = \text{ FR}_1$ and $\le_1 = \le_2$.  Now (D)$_1$(a), (b), (c), (e) and
(E)$_1$(a), (b)$(\alpha)$, (c), (d) holds by the definition of ${\frak
u}$.  Concerning (D)$_1$(d), by assumption (e) of Lemma \scite{838-e.2}
clearly FR$^+_1 \ne \emptyset$ and by Observation \scite{838-e.2E} we have
$[M \in K_\lambda \Rightarrow (M,N,a) \in \text{ FR}_1$ for some pair
$(N,\bold I)]$, i.e. (D)$_1$(d) holds.
As ${\frak K}_\lambda$ has amalgamation and (D)$_1$(d) holds, clearly
\mr
\item "{$(*)$}"  if $M \le_{{\frak K}_\lambda} N$ then for some pair
$(N',\bold J)$ we have $N \le_{{\frak K}_\lambda} N'$ and $(M,N',\bold
J) \in \text{ FR}^+_1$.
\ermn
Concerning (E)$_1$(b)$(\beta)$, just remember that $=_\tau$ is a fake
 equality, recalling subclause (iii) of clause (e) of Definition
\scite{838-e.1A}.  Now the main point is amalgamation = clause (F) of
Definition \scite{838-1a.3}.  We first ignore the $=_\tau$ and
disjointness, that is, we work in ${\frak K}_\lambda$ not 
${\frak K}'_\lambda$; easily this suffices.   So we assume
that $(M_0,M_\ell,\bold J_\ell) \in \text{ FR}_\ell$ for $\ell=1,2$
and by $(*)$ above without loss of generality  $\ell=1,2 \Rightarrow \bold J_\ell \ne
\emptyset$ so let $\bold J_\ell = \{a_\ell\}$.  Let $p_\ell = 
\ortp_{{\frak K}_\lambda}(a_\ell,M_0,M_\ell)$ and by 
\yCITE[1a.43]{E46}(1) we can find
a reduced $(M'_0,M'_1,a_1) \in K^{3,\text{na}}_\lambda$ which is
$\le_{\text{na}}$-above $(M_0,M_1,a_1)$.  We can apply Claim
\yCITE[2b.7]{E46}(1) because: first Hypothesis \yCITE[2b.0]{E46}
holds (as ${\frak K}$ is an a.e.c., LS$({\frak K}) \le \lambda$ and
$K_\lambda \ne \emptyset$) and, second, (amg)$_\lambda$ holds by assumption
(d) of \scite{838-e.2}.  So by \yCITE[2b.7]{E46}(1), the extension property for
such types (i.e. ones with no minimal extensions) holds, so there are
$M'_2$ such that $M'_0 \le_{{\frak K}_\lambda}
M'_2$ and $\le_{\frak K}$-embedding $g$ of $M_2$ into $M'_2$ over
$M_0$ such that $g(a_2) \notin M'_0$.

Again by \yCITE[2b.7]{E46}(1) we can find 
$M''_1 \in K_\lambda$ such that
$M'_2 \le_{\frak K} M''_1$ and $f$ which is a $\le_{\frak K}$-embedding of
$M'_1$ into $M''_1$ such that $f(a_1) \notin M'_2$.

By the definition of ``$(M'_0,M'_1,a)$ is reduced", see Definition
\yCITE[1a.34]{E46} it follows that $f(M'_1) \cap M'_2 
= M'_0$, so $M'_1 \cap M'_2 =
M'_0$.  In particular $f(a_1) \in M'_2 \wedge g(a_2) \notin f(M'_1)$
so we are done.  Now the result with disjointness follows because
$=_\tau$ is a fake equality.     \hfill$\square_{\scite{838-e.2M}}$ 
\enddemo
\bigskip

\proclaim{\stag{838-e.2C} Claim}  1) ${\frak u}$ has the 
vertical coding property, see Definition \scite{838-2b.9}(5).
\nl
2) If $(M,N,\bold I) \in \text{\rm FR}^+_2$ and $a \in \bold I \and b
 \in M \Rightarrow \neg(a =^{N_1}_\tau b)$ \ub{then} this triple
 has the true vertical coding$_0$ property (see
Definition \scite{838-2b.9}(1B).
\nl
3) ${\frak K}_{\lambda^+}$ has a superlimit model which is saturated.
\nl
4) For almost$_2$ every $(\bar M,\bar{\bold J},\bold f) \in
   K^{\text{qt}}_{\frak u}$ the model $M_\partial$ is saturated.
\nl
5) Every $(\bar M,\bar{\bold J},\bold f) \in K^{\text{qt}}_{\frak u}$
has the vertical coding property (see Definition
\scite{838-2b.9}(3)) \ub{when} $M_{\lambda^+} 
\in {\frak K}_{\lambda^+}$ is saturated.
\nl
6) For some ${\frak u}-0$-appropriate function ${\frak h}$, for every $M \in
K^{{\frak u},{\frak h}}_{\lambda^{++}}$ the model $M/=^M$ have
cardinality $\lambda^{++}$ and the set $a/=^M_\tau$ has cardinality
$\lambda^{++}$ for every $a \in M$.  
\endproclaim
\bigskip

\demo{Proof}  1) Follows by part (3),(4),(5).
\nl
2) By the choice of FR$^{\frak u}_1$ for some $a,\bold I 
= \{a\}$ and the type $\ortp_{{\frak K}_\lambda}(a,M,N)$ has no minimal 
extension.  To prove the true vertical coding
property assume that $(\langle M^\ell_i:i \le \beta \rangle,\langle
\bold J^\ell_i:i < \beta \rangle)$ for $\ell=1,2$ and $\langle \bold
I_i:i \le \beta \rangle$ are as in Definition \scite{838-2b.9}(1), i.e.,
they form a ${\frak u}$-free $(\beta,1)$-rectangle with 
$(M,N,\{a\}) \le^1_{\frak u} (M^1_0,M^2_0,\bold I_0)$; i.e. there is
$\bold d$, a ${\frak u}$-free $(\beta,1)$-rectangle such that $M^{\bold
d}_{i,0} = M^1_i,M^{\bold d}_{i,1} = M^2_i,\bold J^{\bold d}_{i,\ell} 
= \bold J^\ell_i,\bold I^{\bold d}_{i,\ell} = \bold I_i$).

So $M \le_{{\frak K}_\lambda} M^1_\beta,
N \le_{{\frak K}_\lambda} M^2_\beta,a \in N \backslash M^1_\beta$ and
$\bold I_i = \{a\}$.  As
$\ortp_{{\frak K}_\lambda}(a,M,N)$ has no minimal extension 
we can find $M^1_{\beta +1}$ such
that $M^1_\beta \le_{{\frak K}_\lambda} M^1_{\beta +1}$ and
$\ortp_{{\frak K}_\lambda}(a,M^1_\beta,M^2_\beta)$ has at least 
two non-algebraic extensions
in ${\Cal S}_{{\frak K}_\lambda}(M^1_\beta)$, hence we can 
choose $p_1 \ne p_2 \in {\Cal S}^{\text{na}}_{{\frak K}_\lambda}
(M^1_{\beta +1})$ extending $\ortp_{{\frak K}_\lambda}
(a,M^1_\beta,M^2_\beta)$.  
Now treating equality as congruence
without loss of generality  $M^1_{\beta +1} \cap M^2_\beta = M^1_\beta$ and there are
$N^1,N^2 \in K_\lambda$ such that $M^1_{\beta +1} \le_{\frak K}
N^\ell,M^2_\beta \le_{\frak K} N^\ell$ and 
$\ortp_{{\frak K}_\lambda}(a,M^1_{\beta +1},N_\ell) = p_\ell$ for $\ell=1,2$.  
\nl
Letting $M^{2,\ell}_{\beta +1} := N_\ell$ we are done. 
\nl
3)  If ${\frak K}$ is categorical in $\lambda^+$ then the desired
 conclusion holds as 
every $M \in {\frak K}_{\lambda^+}$ is saturated above $\lambda$ by
clause (g) of the assumption of \scite{838-e.2}.  If
${\frak K}$ only has a superlimit model in ${\frak K}_{\lambda^+}$
as there is a $M' \in {\frak K}_{\lambda^+}$ saturated above
$\lambda$, necessarily the superlimit $M' \in {\frak K}_{\lambda^+}$
is saturated above $\lambda$ by \yCITE[2b.13]{E46}(4).
\nl
4) We prove the existence of ${\frak g}$ for the ``almost$_2$" (or use the
proof of \scite{838-e.1Z}(1)).  Now recalling $(*)_4$ of \scite{838-1a.43}(1) 
for each $M \in K_\lambda$ with universe $\in [\lambda^{++}]^\lambda$, 
we can choose a sequence $\langle p_{M,\alpha}:
\alpha < \lambda^+\rangle$ listing
${\Cal S}_{{\frak K}_\lambda}(M)$.  When defining the value 
${\frak g}(\bar M^1,\bar{\bold J}^1,\bold f^1,\bar M^2 \rest (\delta +
\bold f^1(\delta)+1)),\bar{\bold J}^2 \rest (\delta + \bold
f^1(\delta),\bar{\bold I} \rest (\delta + \bold f^1(\delta) +1),S)$,
see Definition \scite{838-1a.43}(1)(c) we just realize all 
$p_{M^2_i,j}$ with $i,j < \delta$.  Recalling that by part (3) the
union of a $\le_{\frak K}$-increasing sequence of length $<
\lambda^{++}$ of saturated members of $K_{\lambda^+}$ is saturated, we
are done.
\nl
5) Holds by Observation \scite{838-2b.11}(1).
\nl
6) Easy, as in \scite{838-e.1C}(3).  \hfill$\square_{\scite{838-e.2C}}$
\enddemo
\bigskip

\demo{Continuation of the Proof of \scite{838-e.2}}  By the 
\scite{838-1a.37}, \scite{838-e.2M}, \scite{838-e.2C}(1) we can apply 
Theorem \scite{838-2b.13}.  \hfill$\square_{\scite{838-e.2}}$
\enddemo
\bigskip

\remark{\stag{838-e.2T} Remark}  So we have finished proving \scite{838-e.2}.
\endremark
\bn
\centerline {$* \qquad * \qquad *$}
\bn
\ub{(C) \quad The symmetry property of PC$_{\aleph_0}$ classes}
\bigskip

Here we pay a debt from Theorem \yCITE[5.23]{88r}(1), so naturally we
assume knowledge of \sectioncite[\S5]{88r}; of course later results
supercede this.  Also we can avoid this subsection altogether, dealing with the
derived good $\aleph_0$-frame in \yCITE[Ex.1]{600}.
\bigskip

\proclaim{\stag{838-e.3} Theorem}  $\dot I(\aleph_2,{\frak K}) \ge
\mu_{\text{unif}}(\aleph_2,2^{\aleph_1})$ moreover $\dot
I(\aleph_2,{\frak K}(\aleph_1\text{-saturated})) \ge
\mu_{\text{unif}}(\aleph_2,2^{\aleph_1})$ \ub{when}:
\mr
\item "{$\circledast$}"  $(a) \quad$ (set theory)
{\roster
\itemitem{ ${{}}$ }  $(\alpha) \quad 2^{\aleph_0} < 2^{\aleph_1} < 
2^{\aleph_0}$ and 
\sn
\itemitem{ ${{}}$ }  $(\beta) \quad$ {\rm WDmId}$(\aleph_1)$ 
is not $\aleph_2$-saturated
\endroster}
\item "{${{}}$}"  $(b) \quad {\frak K}$, an a.e.c., is 
$\aleph_0$-representable,
i.e., is {\rm PC}$_{\aleph_0}$-a.e.c., see Definition 
\yCITE[1.4]{88r}(4),(5)
\sn
\item "{${{}}$}"  $(c) \quad {\frak K}$ is as in
\yCITE[4.5]{88r},
\yCITE[5.0]{88r}
\sn
\item "{${{}}$}"  $(d) \quad {\frak K}$ fails the symmetry property or 
the uniqueness of two sided 
\nl

\hskip25pt stable amalgamation, see Definition \yCITE[5.21]{88r},
equivalently
\sn
\item "{${{}}$}"  $(d)' \quad {\frak K}$ fails the uniqueness of
one-sided amalgamation
\sn
\item "{${{}}$}"  $(e) \quad \bold D$ is countable, see Definition
\yCITE[5.1]{88r},
\yCITE[5.6.1]{88r}
and \sectioncite[\S3]{600}(B).
\endroster
\endproclaim
\bigskip

\remark{Remark}  1) On omitting ``WDmId$(\aleph_1)$  
is not $\aleph_2$-saturated", see Conclusion \scite{838-e.6L}.
\nl
2) Clause (e) of \scite{838-e.3} is reasonable as we can without loss of
generality assume it by Observation \yCITE[5.23.8]{88r}.
\endremark
\bigskip

\demo{Proof}  Let $\lambda = \aleph_0$,  \wilog
\mr
\item "{$(f)$}"  for $M \in K$, any finite sequence is coded by an
element.
\ermn

Now by \sectioncite[\S3 ]{600}(B), i.e. \yCITE[Ex.1]{600} we 
have ${\frak s}_{\aleph_0}$, which we call ${\frak s} = {\frak
s}_{\frak K} = {\frak s}^1_{\frak K}$, is a good
$\aleph_0$-frame as defined (and proved) there, moreover 
${\frak s}_\lambda$ is type-full and ${\frak K}^{\frak s} = {\frak K}$.

The proof is broken, as in other cases, i.e., we prove it by Theorem
\scite{838-2b.3} 
which is O.K. as by \scite{838-e.3.8} + \scite{838-e.3.12} below its assumptions holds.
\enddemo
\bigskip

\remark{Remark}  Note in the following definition FR$_1$, FR$_2$
are quite different even though $\le^1_{\frak u},\le^2_{\frak u}$
 are the same, except the domain. 
\endremark
\bigskip

\definition{\stag{838-e.3.4} Definition}  We define ${\frak u} = 
{\frak u}^3_{{\frak K}_{\aleph_0}}$ by
\mr
\item "{$(a)$}"  $\partial = \partial_{\frak u} =\aleph_1$
\sn
\item "{$(b)$}"  ${\frak K}_{\frak u} = {\frak K}_{\aleph_0}$ so
$K^{\text{up}}_{\frak u} = {\frak K}$
\sn
\item "{$(c)$}"  FR$_2$ is the family of triples $(M,N,\bold J)$ such
that:
{\roster
\itemitem{ $(\alpha)$ }  $M \le_{\frak K} N \in {\frak K}_{\aleph_0}$
\sn  
\itemitem{ $(\beta)$ }  $\bold J \subseteq N \backslash M$ and 
$|\bold J| \le 1$
\endroster}
\item "{$(d)$}"  $(M_1,N_1,\bold J_1) \le_2 (M_2,N_2,\bold J_2)$ \ub{iff}
{\roster
\itemitem{ $(\alpha)$ }  both are from FR$_1$
\sn  
\itemitem{ $(\beta)$ }  $M_1 \le_{\frak K} M_2$ and $N_1 \le_{\frak K}
N_2$ 
\sn
\itemitem{ $(\gamma)$ }   $\bold J_1 \subseteq \bold J_2$
\sn
\itemitem{ $(\delta)$ }   if $\bar c \in \bold J$ then gtp$(\bar
c,M_2,N_2)$ is the stationarization of gtp$(\bar c,M_1,N_1)$
\endroster}
\item "{$(e)$}"  FR$_1$ is the class of triples $(M,N,\bold J)$ such
that
{\roster
\itemitem{ $(\alpha)$ }  $M \le_{\frak K} N$ are countable
\sn
\itemitem{ $(\beta)$ }  $\bold J \subseteq N \backslash M$ or, less
pedantically, $\bold J \subseteq {}^{\omega >} N \backslash {}^{\omega>} M$ 
\sn
\itemitem{ $(\gamma)$ }  if $|\bold J| > 1$ then $\bold J =
{}^{\omega >} N \backslash {}^{\omega>} M$  and $N$ is
$(\bold D(M),\aleph_0)^*$-homogeneous
\endroster} 
\sn
\item "{$(f)$}"  $\le_1$ is defined as in clause (d) \ub{but} on
FR$_2$.
\endroster
\enddefinition
\bigskip

\proclaim{\stag{838-e.3.8} Claim}  1) ${\frak u}$ is a nice construction
framework.
\nl
2) For almost$_2$ all triples $(\bar M,\bar{\bold J},\bold f) \in
K^{\text{qt}}_{\frak u}$ the model $M$ is saturated (for ${\frak K}$).  
\endproclaim  
\bigskip

\demo{Proof}  1) The main points are:
\mr
\item "{$\boxtimes_1$}"  $(M,N,\bold J) \in \text{ FR}_2$ and
$n < \omega \Rightarrow (M_n,N_n,\bold J_n) \le_2 (M,N,\bold J)$ \ub{when}:
{\roster
\itemitem{ $(a)$ }  $(M_n,N_n,\bold J_n) \in  \text{ FR}_2$
\sn
\itemitem{ $(b)$ }  $(M_n,N_n,\bold J_n) \le_2 (M_{n+1},N_{n+1},
\bold J_{n+1})$ for $n < \omega$
\sn
\itemitem{ $(c)$ }   $M = \cup\{M_n:n < \omega\}$
\sn
\itemitem{ $(d)$ }   $N = \cup\{N_n:n < \omega\}$
\sn
\itemitem{ $(e)$ }   $\bold J = \cup\{\bold J_n:n < \omega\}$.
\endroster}
\ermn
[\ub{Why $\boxtimes_1$ holds}?  See \yCITE[5.13]{88r}(9).]
\mr
\item "{$\boxtimes_2$}"  $(M,N,\bold J) \in \text{ FR}_1$ and $n <
\omega \Rightarrow (M_n,N_n,\bold J_n) \le_1 (M,N,\bold J)$ \ub{when}
{\roster
\itemitem{ $(a)$ }  $(M_n,N_n,\bold J_n) \in \text{ FR}_1$
\sn
\itemitem{ $(b)$ }  $(M_n,N_n,\bold J_n) \le_1 (M_{n+1},N_{n+1},\bold
J_{n+1})$ for $n < \omega$
\sn
\itemitem{ $(c)$ }  $M = \cup\{M_n:n < \omega\}$
\sn
\itemitem{ $(d)$ }  $N = \cup\{N_n:n < \omega\}$
\sn
\itemitem{ $(e)$ }  $\bold J = \cup\{\bold J_n:n < \omega\}$.
\endroster}
\ermn
[Why does $\boxtimes_2$ holds?  If $|\bold J| \le 1$ then the proof is
similar to the one in $\boxtimes_1$, so assume that $|\bold J| > 1$,
so as $n < \omega \Rightarrow \bold J_n \subseteq \bold J_{n+1}$ by
clause (b) of $\boxtimes_2$ and $\bold J = \cup\{\bold J_n:n <
\omega\}$ by clause (e) of the $\boxtimes_2$ 
necessarily for some $n,|\bold J_n| > 1$,
so without loss of generality  $|\bold J_n| > 1$ for every $n < \omega$.  So by the
definition of FR$_1$, we have:
\mr
\item "{$(*)_1$}"  $\bold J_n = {}^{\omega >}(N_n) \backslash
{}^{\omega >}(M_n)$
\sn
\item "{$(*)_2$}"  $N_n$ is $(\bold D(M_n),\aleph_0)^*$-homogeneous.
\ermn
Hence easily
\mr
\item "{$(*)_3$}"  $\bold J = {}^{\omega >}N \backslash {}^{\omega >} M$
\ermn
and as in the proof of $\boxtimes_1$ clearly
\mr
\item "{$(*)_4$}"  if $n < \omega$ and $\bar c \in \bold J_n$ then
gtp$(\bar c,M,N)$ is the stationarization of gtp$(\bar c,M_n,N_n)$.
\ermn
So the demands for ``$(M_n,N_n,\bold J_n) \le_1 (M,N,\bold J)$" holds except
that we have to verify
\mr
\item "{$(*)_5$}"  $N$ is $(\bold D(M),\aleph_0)^*$-homogeneous.
\ermn
[Why this holds?  Assume $N \le_{{\frak K}_{\aleph_0}} N^+,\bar a \in
{}^{\omega >}N,\bar b \in {}^{\omega >}(N^+)$.  So gtp$(\bar a \char
94 \bar b,M,N^+) \in \bold D(M)$, hence by \yCITE[5.13]{88r}(9)
 it is the stationarization of
gtp$(\bar a \char 94 \bar b,M_{n_0},N^+)$ for some $n_0 < \omega$.
Also for some $n_1 < \omega$, we have $\bar a \in {}^{\omega
>}(N_{n_1})$.  Now some $n < \omega$ is $\ge n_0,n_1$, so by $(*)_1 + (*)_2$
 for some $\bar b' \in {}^{\omega >}(N_n)$ the type
gtp$(\bar a \char 94 \bar b',M_n,N_n)$ is equal to
gtp$(\bar a \char 94 \bar b,M_n,N^+)$ so $\bar b' \in \bold J_n$.  But
by $(*)_4$, also the 
type gtp$(\bar a \char 94 \bar b',M,N)$ is a stationarization
of gtp$(\bar a \char 94 \bar b',M_n,N) =$ gtp$(\bar a \char 94 \bar b,
M_n,N^+)$ hence gtp$(\bar a \char 94 \bar
b',M,N) = \text{ gtp}(\bar a \char 94 \bar b,M,N^+)$ so we are done.] 
\nl
Thus we have finished proving $\boxtimes_2$]
\mr
\item "{$\boxtimes_3$}"   clause (F) of Definition \scite{838-1a.3} holds.
\ermn
[Why $\boxtimes_3$?  By \yCITE[5.20]{88r}.]
\nl
Together we have finished proving ${\frak u}$ is a nice construction framework.
\nl
2) Note that
\mr
\item "{$\boxtimes_4$}"   $M_{\omega_1} \in K_{\aleph_1}$ is saturated
\ub{when}:
{\roster
\itemitem{ $(a)$ }  $\langle M_\alpha:\alpha \le \omega_1\rangle$ is
$\le_{\frak K}$-increasing continuous
\sn
\itemitem{ $(b)$ }  $\alpha < \omega_1 \Rightarrow M_\alpha \in
K_{\aleph_0}$
\sn
\itemitem{ $(c)$ }  $M_{\alpha +1}$ is $(\bold
D(M_\alpha),\aleph_0)^*$-homogeneous for $\alpha < \omega_1$
\nl
or just
\sn
\itemitem{ $(c)'$ }  if $\alpha < \omega_1$ and $p \in \bold
D(M_\alpha)$ is a 1-type then for some $\beta \in [\alpha,\omega_1)$
and  some $c \in M_{\beta +1}$, 
gtp$(c,M_\beta,M_{\beta +1})$ is the stationarization of $p$.
\endroster}
\ermn
[Why?  Obvious.]

Let $S \subseteq \omega_1$ be stationary, so clearly it suffices to prove:
\mr
\item "{$\boxtimes_5$}"  there is ${\frak g}$ as in Definition
\scite{838-1a.43}(1), \scite{838-1a.45}(3) for $S$ and our ${\frak u}$ such that: 
\nl

if $\langle
(\bar M^\alpha,\bar{\bold J}^\alpha,\bold f^\alpha):\alpha \le
\delta\rangle$ is $\le^{\text{qs}}_{\frak u}$-increasing continuous
1-obeying ${\frak g}$ 
\nl

(and $\delta < \partial^+$ is a limit ordinal)
\ub{then} $M^\delta \in {\frak K}_{\aleph_1}$ is saturated.
\ermn
[Why?  Choose ${\frak g}$ such that if the pair $((\bar M',\bar{\bold
J}',\bold f'),(\bar M'',\bar{\bold J}'',\bold f''))$ does 1-obey ${\frak g}$
then for every $\alpha < \partial$ and $p \in \bold D(M''_\alpha)$
we have
\mr
\item "{$(*)$}"  for stationarily many $\delta \in S$ for some $i <
\bold f''(\delta)$ the type gtp$(a,M''_\alpha,M''_{\delta +i+1})$ is
the stationarization of $p$ where $\bold J''_{\delta +i} = \{a\}$.
\ermn
Now assume that $\langle (\bar M^\zeta,\bar{\bold J}^\zeta,\bold
f^\zeta):\zeta < \delta\rangle$ is $\le_{\text{qs}}$-increasing (for
our present ${\frak u}$) and $\delta = \sup(u)$ where $u =
\{\zeta:((\bar M^\zeta,\bold J^\zeta,\bold f^\zeta),(\bar M^{\zeta
+1},\bar{\bold J}^{\zeta +1},\bold f^{\zeta +1}))$ does 1-obeys
${\frak g}\}$, and we should prove that $M^\delta :=
\cup\{M^\zeta_\partial:\zeta < \partial\}$ is saturated.  Without loss
of generality $u$ contains all odd ordinals $< \delta$ and $\delta =
\text{ cf}(\delta)$.  If $\delta = \aleph_1$ this is obvious, and if
$\delta = \aleph_0$ just use non-forking of types, and the criterion
in $\boxtimes_4$ using $(*)$.  So $\boxtimes_5$ is proved.]
\hfill$\square_{\scite{838-e.3.8}}$ 
\enddemo
\bigskip

\proclaim{\stag{838-e.3.12} Claim}  ${\frak u}$ has the weak coding
property.
\endproclaim
\bigskip

\demo{Proof}  Clearly by clause $(d)'$ of the assumption, i.e.
by Definition \yCITE[5.21]{88r}(2),(3) there are $N_\ell \, (\ell \le
2)),N'_3,N''_3$ such that:
\mr
\item "{$(*)_1$}"  $(a) \quad N_0 \le_{\frak u} N_1$ and $N_0
\le_{\frak u} N_2$
\sn
\item "{${{}}$}"   $(b) \quad N_\ell \le_{\frak u} N'_3$ and $N_\ell
\le_{\frak u} N''_3$ for $\ell=1,2$
\sn
\item "{${{}}$}"   $(c) \quad N_0,N_1,N_2,N'_3$ is in one-sided
amalgamation, i.e. 
\nl

\hskip25pt $\bar a \in {}^{\omega >}(N_1) \Rightarrow
(N_0,N_1,\{\bar a\}) \le^2_{\frak u} (N_2,N'_3,\{\bar a\})$ 
\nl

\hskip25pt (hence $N_1 \cap N_2 = N_0$)
\sn
\item "{${{}}$}"   $(d) \quad N_0,N_1,N_2,N''_3$ is in one sided
amalgamation
\sn
\item "{${{}}$}"   $(e) \quad$ there are no $(N_3,f)$ such that $N''_3
\le_{\frak u} N_3$ and $f$ is a $\le_{\frak u}$-embedding
\nl

\hskip25pt  of $N'_3$ into $N_3$ over $N_1 \cup N_2$.
\ermn
Now \wilog
\mr
\item "{$(*)_2$}"   $(N_0,N_1,\bold J) \in \text{ FR}_1$ where $\bold
J \in {}^{\omega >}(N_1) \backslash {}^{\omega >}(N_0)$ such that
$|\bold J| > 1$.
\ermn
[Why?  We can find $N^+_1 \in K_{\frak u}$ which is $(\bold
D(N_1),\aleph_0)^*$-homogeneous over $N_1$ and without loss of generality  $N^+_1 \cap
N'_3 = N_1 = N^+_1 \cap N''_3$.  Now we can find $N^*_3$ and $N^{**}_3
\in K_{\frak u}$ such that $(N_1,N^+_1,N'_3,N^*_3)$ as well as
$(N_1,N^+_1,N''_3,N^{**}_3)$ is in one sided stable amalgamation.  It
follows that $(N_0,N^+_1,N_2,N^*_3,N^{**}_3)$ satisfies all the
requirements on $(N_0,N_1,N_2,N'_3,N''_3)$ and in addition the demand
in $(*)_2$ so we are done.]

Also \wilog 
\mr
\item "{$(*)_3$}"   $(N_0,N_2,\bold I) \in \text{ FR}_2$ and $|\bold
I|=1$ and $N_2$ is $(\bold D(N_0),\aleph_0)^*$-homogeneous.
\ermn
[Why?  Similarly to $(*)_1$.]

To prove ${\frak u}$ has the weak coding we can assume (the saturation
is justified by \scite{838-e.3.8}(2))
\mr
\item "{$(*)_4$}"  $(\bar M,\bar{\bold J},\bold f) \in
K^{\text{qt}}_{\frak u}$ and $M := \cup\{M_\alpha:\alpha < \omega_1\}$
is saturated for ${\frak K}$.
\ermn
Now by renaming without loss of generality  
\mr  
\item "{$(*)_5$}"  $N_0 \le_{\frak K} M_{\alpha(0)}$ and $N_1 \cap M =
N_0$ and $(N_0,N_1,\bold J) \le_1 (M_{\alpha(0)},N'_1,\bold J')$ and
$N'_1 \cap M = M_{\alpha(0)}$.
\ermn
It suffices to prove that $(\alpha(0),N'_1,\bold J')$ is as required
in \scite{838-2b.1}(3).
Next by the definition of ``having the weak coding property", for our
purpose we can assume we are given $(N'',\bold J'')$ such that
\mr
\item "{$(*)_6$}"  $\alpha(0) \le \delta < \omega_1$ and
$(N_0,N'_1,\bold J') \le_1 (M_\delta,N'',\bold J'')$.
\ermn
By the definition of $\le_1$ we know that
\mr
\item "{$(*)_7$}"  $N''$ is $(\bold
D(M_\delta),\aleph_0)^*$-homogeneous over $M_\delta$.
\ermn
As $\cup\{M_\alpha:\alpha < \omega_1\}$ is saturated (for $K^{\frak u}$)
we can find $\beta \in (\delta,\omega_1)$ such that $M_\beta$ is
$(\bold D(M_\delta),\aleph_0)^*$-homogeneous over $M_\delta$. 

As ${\frak K}$ is categorical in $\aleph_0$
\mr
\item "{$(*)_8$}"  there is an isomorphism $f_0$ from $N_0$ onto
$M_\delta$.
\ermn
Similarly using the uniqueness over $N_0$ of a countable $(\bold
D(M_0),\aleph_0)^*$-homogeneous model over $N_0$
\mr
\item "{$(*)_9$}"  there are isomorphisms $f_1,f_2$ from $N_1,N_2$
onto $N'',M_\beta$ respectively extending $f_0$.
\ermn
Lastly, $M_\beta,N$ can be amalgamated over $M_\delta$
in the following two ways:
\mr
\item "{$\odot_1$}"  there are $f',M' \in K_{\frak u}$ such that $f'$
is an isomorphism from $N'_3$ onto $M'$ extending $f_1 \cup f_2$ 
\sn
\item "{$\odot_2$}"  there are $f'',M'' \in K_{\frak u}$ such that $f''$
is an isomorphism from $N^1_3$ onto $M''$ extedning $f_1 \cup f_2$.
\ermn
This is clearly enough.  The rest should be clear.]  
\hfill$\square_{\scite{838-e.3.12}}$
\enddemo
\bigskip

\demo{Proof of \scite{838-e.3}}  By the claims above.
\hfill$\square_{\scite{838-e.3}}$ 
\enddemo
\bn
\centerline {$* \qquad * \qquad *$}
\bn
\ub{(D) \quad Density of $K^{3,\text{uq}}_\lambda$ when minimal triples
are dense}
\bigskip

Having taken care of \sectioncite[\S3,\S4]{E46} and of
\sectioncite[\S5]{88r}, we now deal with proving the
non-structure results of \sectioncite[\S6]{E46}, i.e. \cite[\S6]{Sh:576}, 
relying on \scite{838-2b.3} instead of \cite[\S3]{Sh:576}.  Of course, later
we prove stronger results but have to work harder, both model
theoretically (including ``${\frak s}$ is almost a good 
$\lambda$-frame") and set theoretically (using (vertical coding in)
Theorem \scite{838-2b.13} and \S3 rather than (weak coding in) Theorem
\scite{838-2b.3}).  
\bigskip

This is used in \yCITE[6f.13]{E46}.
\proclaim{\stag{838-e.4} Theorem}  The non-structure results of
\yCITE[6f.13]{E46}, Case 1 holds.

It details: $\dot I(\lambda^{++},{\frak K}) \ge
\mu_{\text{unif}}(\lambda^{++},2^{\lambda^+})$ \ub{when} we are assuming:
\mr
\item "{$(A)$}"  set theoretically: 
{\roster
\itemitem{ $(a)$ }  $2^\lambda < 2^{\lambda^+} < 2^{\lambda^{++}}$ and
\sn
\itemitem{ $(b)$ }   the weak diamond ideal on $\partial :=
\lambda^+$ is not $\partial^+$-saturated
\endroster}
\item "{$(B)$}"  model theoretically:
{\roster
\itemitem{ $(a)$ }  ${\frak K}$ is an a.e.c., $\lambda \ge { \text{\rm
LS\/}}({\frak K})$
\sn
\itemitem{ $(b)$ }   $(\alpha) \quad {\frak K}$ is categorical in $\lambda$
\sn
\itemitem{ ${{}}$ }  $(\beta) \quad {\frak K}$ is categorical in
$\lambda^+$ or just has a superlimit model in $\lambda^+$
\sn
\itemitem{ $(c)$ }  $(\alpha) \quad {\frak K}$ has amalgamation in $\lambda$
\sn
\itemitem{ ${{}}$ }  $(\beta) \quad {\frak K}$ is stable in $\lambda$
\ub{or} just $M \in K_\lambda \Rightarrow 
|{\Cal S}^{\text{min}}_{\frak K}(M)| \le \lambda$
\sn
\itemitem{ $(d)$ }   $(\alpha) \quad$ the minimal types are dense (for $M \in 
{\frak K}_\lambda$) 
\sn
\itemitem{ ${{}}$ }  $(\beta) \quad$ for $M \in K_\lambda$ the set
${\Cal S}^{\text{min}}_{{\frak K}_\lambda}(M) = \{p \in {\Cal
S}_{{\frak K}_\lambda}(M):p$ minimal$\}$ 
\nl

\hskip35pt is inevitable
\sn
\itemitem{ ${{}}$ }  $(\gamma) \quad$ the $M \in 
K^{\text{slm}}_{\lambda^+}$ is saturated above $\lambda$ 
\sn
\itemitem{ $(e)$ }   above (by $\le_{\text{na}}$) some $(M^*,N^*,a) \in
K^{3,\text{na}}_\lambda$ there is no triple with the uniqueness
property, i.e. from $K^{3,\text{uq}}_\lambda$, see \yCITE[6f.2]{E46}.
\endroster}
\endroster
\endproclaim
\bigskip

\remark{\stag{838-e.4.0} Remark}  1) Note: every $M \in K_{\lambda^+}$ is saturated
above $\lambda$ when the first, stronger version of (B)(b)$(\beta)$
holds noting (B)(c)$(\beta) +$ (B)(d)$(\beta)$.
\nl
2) When we use the weaker version of clause (b)$(\beta)$, i.e. ``there is 
superlimit $M \in K_{\lambda^+}$" then we have to prove that for
almost$_2$ every $(\bar M,\bar{\bold J},\bold f)$, the model $M_{\lambda^+}$ is
saturated above $\lambda$ which, as in earlier cases, can be done; see
\yCITE[2b.13]{E46}(4).
\nl
3) Concerning clause (B)(d): ``the minimal types are dense", it
follows from (amg)$_\lambda +$ (stb)$_\lambda$, i.e. from clause (c)
recalling \yCITE[2b.4]{E46}(4).
\nl
4) Note that \scite{838-e.4.2}, \scite{838-e.4.3} does not depend on clause
   (A)(b) of \scite{838-e.4}.
\endremark
\bigskip

\definition{\stag{838-e.4.1} Definition}  We define ${\frak u} = {\frak
u}_4 = {\frak u}^4_{{\frak K}_\lambda}$ as follows:
\mr
\item "{$(a)$}"  $\partial_{\frak u} = \lambda^+$
\sn
\item "{$(b)$}"  ${\frak K}_{\frak u} = {\frak K}_\lambda$ or
pedantically ${\frak K}'_\lambda$, see Definition \scite{838-1a.19}
\sn
\item "{$(c)_1$}"  $\text{FR}^{\frak u}_1$ is the set of triples 
$(M,N,\bold I)$ satisfying $M \le_{\frak K} N \in K_\lambda,\bold I =
\emptyset$ or $\bold I = \{a\}$ and the type
$\ortp_{{\frak K}_\lambda}(a,M,N)$ is minimal, pedantically, if $a/=^N
\notin M/=^N$ then $\ortp(a/=^N,M/=^N,N/=^N)$ is minimal
\sn
\item "{$(c)_2$}"  $(M_1,N_1,\bold I_1) 
\le_1 (M_2,N_2,\bold I_2)$ iff (both are
FR$^{\frak u}_1$ and) $M_1 \le_{\frak K} M_2,N_1 \le_{\frak K} N_2$
and $\bold I_1 \subseteq \bold I_2$ (in the non-trivial cases,
equivalently, $\bold I_1 = \bold I_2$), pedantically, if $(a/=^N)
\notin M/=^N$ then $\ortp(a/=^N,M/=^N,N/=^N)$ is minimal
\sn
\item "{$(d)$}"  FR$^{\frak u}_2 = { \text{\rm FR\/}}^{\frak u}_1$ and
$\le^2_{\frak u} = \le_{\frak u}^1$.
\endroster
\enddefinition
\bigskip

\proclaim{\stag{838-e.4.2} Claim}  ${\frak u}$ is a nice construction
framework which is self-dual.
\endproclaim
\bigskip

\demo{Proof}  Easy (amalgamation, i.e. clause (F) of Definition
\scite{838-1a.3} holds by the proof of symmetry in Axiom
(E)(f) in proof of Theorem \yCITE[8h.1]{E46}).  
\hfill$\square_{\scite{838-e.4.2}}$ 
\enddemo
\bigskip

\proclaim{\stag{838-e.4.3} Claim}  1) Every $(M,N,\bold I) \in 
{ \text{\rm FR\/}}_1$ such that $[a \in \bold I \and b \in M
 \Rightarrow \neg a =^N_\tau b]$ has the true weak coding property (see 
Definition \scite{838-2b.1}(1A)).
\nl
2) ${\frak u}$ has the weak coding property.
\nl
3) For almost$_2$ every $(\bar M,\bar{\bold J},\bold f) \in
   K^{\text{qt}}_{\frak u}$ the model $M_\partial \in K_{\lambda^+}$
   is satuarted above $\lambda$.
\nl
4) For some ${\frak u}-\{0,2\}$-appropriate function 
${\frak h}$, for every $M \in
K^{{\frak u},{\frak h}}_{\lambda^{++}}$ the model $M/=^M$ has
cardinality $\lambda^{++}$ and is saturated above $\lambda$. 
\endproclaim
\bigskip

\demo{Proof}  1) Straight.
\nl
2) By part (1) above and part (3) below.
\nl
3) By clauses (B)(c)$(\alpha),(\beta)$ 
of \scite{838-e.4}, clearly there is a $M \in
K_{\lambda^+}$ which is saturated above $\lambda$.  If in (B)(b)$(\beta)$
we assume categoricity in $\lambda^+$ then 
every $M \in K_{\lambda^+}$ is saturated above $\lambda$, but then it
is obvious that part (1) implies part (2) by \scite{838-2b.11}(4)(b).
For any stationary $S \subseteq \partial$, we choose
${\frak h}$ such that
\mr
\item "{$(*)$}"  if $((\bar M^1,\bar{\bold J}^1,\bold f^1),(\bar
M^2,\bar{\bold J}^2,\bold f^2))$ does 2-obeys ${\frak h}$ then: for
stationarily many $\delta \in S$ there is successor ordinal $i <
\bold f^2(\delta)$ such that $M^2_{\delta +i+1}$ is $<_{{\frak
K}_\lambda}$-universal over $M^2_{\delta +i}$ (hence $M^2_\partial$ is
saturated above $\lambda$ and is the superlimit model in 
${\frak K}_{\lambda^+}$). 
\ermn
Alternatively do as in \scite{838-e.2C}(4), using \yCITE[8h.1]{E46}.
\nl
4) As in \scite{838-e.1C}(3).
\enddemo
\bigskip

\demo{Proof of \scite{838-e.4}}  By \scite{838-e.4.2}, \scite{838-e.4.3} and
 Theorem \scite{838-2b.3}.  \hfill$\square_{\scite{838-e.4}}$
\enddemo
\bn
\centerline {$* \qquad * \qquad *$}
\bn
\ub{(E) Density of $K^{3,\text{uq}}_{\frak s}$ for good $\lambda$-frames}
\bigskip
We now deal with the non-structure proof in \yCITE[nu.6]{600}, that is
justifying why the density of $K^{3,\text{uq}}_{\frak s}$ holds.
\bigskip

Before we state the theorem, in order to get rid of the problem of
disjoint amalgamation, one of the ways is to note:
\definition{\stag{838-e.5D} Definition}  Assume that ${\frak s}$ is a good
$\lambda$-frame (or just an almost good $\lambda$-frame see Definition
in \scite{838-5t.3} below or just a pre-$\lambda$-frame, see \yCITE[8h.3]{E46}).
\nl
1) We say that ${\frak s}$ has fake equality $=_*$ when 
${\frak K}_{\frak s}$ has the fake equality $=_*$, see Definition
\scite{838-3r.84}(1) and $\ortp_{\frak s}(a,M_1,M_2)$ does not 
fork over $M_0$ iff $M_0 \le_{\frak s} M_1
\le_{\frak s} M_2,a \in M_2$ and letting $M'_\ell =
M_\ell/=^{M_2}_*$ we have $(\forall b \in M_1)(\neg(a =^{M_2}_* b))
\Rightarrow \ortp_{\frak s}(a/=^{M_2}_*,M'_1,M'_2)$ does not fork
over $M'_0$.
\nl
2) We define ${\frak s}' = (K_{{\frak s}'},
{\Cal S}^{\text{bs}}_{{\frak s}'},\nonfork{}{}_{{\frak s}'})$ as follows:
\mr
\item "{$(a)$}"  ${\frak K}_{{\frak s}'} = K'_{\frak s}$, see
\scite{838-1a.19} as in \scite{838-e.1A}
\nl
so $\tau_{{\frak s}'} = \tau'_{\frak s} = 
\tau_{\frak s} \cup\{=_\tau\}$ and a $\tau'_{\frak s}$-model $M$
belongs to $K_{{\frak s}'}$ iff $=^M_*$ is a congruence relation and
the model $M/=^{M'}_\tau$ belongs to ${\frak K}_{\frak s}$
\sn
\item "{$(b)$}"  for $M' \in K_{{\frak s}'}$ we let ${\Cal
S}^{\text{bs}}_{{\frak s}'}(M') = \{\ortp_{{\frak K}_{{\frak
s}'}}(a,M',N'):M' \le_{{\frak s}'} N'$ and $\ortp_{\frak s}
(a/=^{N'}_\tau,M'/=^{M'}_\tau,N'/=^{N'}_\tau) \in {\Cal S}^{\text{bs}}_{\frak
s}(M'/=^{M'}_\tau)$ or $a \in N' \backslash M'$ but $(\exists b \in
M')(a =_\tau b)\}$
\sn
\item "{$(c)$}"  $\ortp_{{\frak s}'}(a,M'_1,M'_2)$ does not fork over
$M'_0$ when $M'_0 \le_{{\frak s}'} M'_1 \le_{{\frak s}'} M'_2$ and either
$\ortp_{\frak s}(a/=^{M'_2}_\tau,M'_1
=^{M'_2}_\tau,M'_2/=^{M'_2}_\tau)$ does not fork over $M'_0/=^{M'_2}$
or for some $b \in M'_0$ we have $M'_2 \models ``a =_\tau b"$ but $a
\notin M'_1$.
\endroster
\enddefinition
\bigskip

\proclaim{\stag{838-e.5E} Claim}  Let ${\frak s},{\frak s}'$ be as in
\scite{838-e.5D}(2).
\nl
1) If ${\frak s}$ is a good $\lambda$-frame \ub{then}
${\frak s}'$ is a good $\lambda$-frame and if ${\frak s}$ is an
almost good $\lambda$-frame then ${\frak s}'$ is an almost good
$\lambda$-frame; and if ${\frak s}$ is a pre-$\lambda$-frame then
${\frak s}'$ is a pre-$\lambda$-frame.  In all cases ${\frak s}'$
has the fake equality $=_\tau$.
\nl
2) For $\mu \ge \lambda,\dot I(\mu,K^{\frak s}) =|\{M'/\cong:M' \in
K^{{\frak s}'}_\mu$ and is $=_\tau$-fuller, that is 
$a \in M' \Rightarrow \|M'/=^{M'}_\tau\| = \mu = |\{b \in M':a 
=^{M'}_\tau b\}|\}$.
\nl
3) If $M' \in K^{{\frak s}'}$ \ub{then} $M'$ is $\lambda^+$-saturated above
$\lambda$ for ${\frak s}'$ iff $M'/=^{M'}_\tau$ is
$\lambda^+$-saturated above $\lambda$ for ${\frak s}$ and $M'$ is
$(\lambda^+,=_\tau)$-full (recalling \scite{838-1a.19}(5A)).
\endproclaim
\bigskip

\remark{\stag{838-e.5F} Remark}  1) By \scite{838-e.5E}(2), the proof of ``$\dot
I(\mu,K^{{\frak s}'})$ is $\ge \chi$" here usually gives ``$\dot
I(\mu,K^{\frak s})$ is $\ge \chi$".  
\nl
2) We define ${\frak s}'$ such that for some $0$-appropriate ${\frak
h}$, if $\langle (\bar M^\alpha,\bar{\bold J}^\alpha,\bold
f^\alpha):\alpha < \partial^+ \rangle$ is $\le_{\text{qt}}$-increasing
continuous 0-obeying ${\frak h}$, then $M=
\cup\{M^\alpha_\partial:\alpha < \partial^+\}$ satisfies the
condition in \scite{838-e.5E}(2); it does not really matter if we need
$\{0,2\}$-appropriate ${\frak h}$.
\nl
3) Recall Example \scite{838-1a.21} as an alternative to \scite{838-e.5E}(2).
\nl
4) Another way to deal with disjointness is by \scite{838-5t.43},
\scite{838-5t.45} below.
\endremark
\bigskip

\demo{Proof}  Easy and see \scite{838-1a.20}.  \hfill$\square_{\scite{838-e.5E}}$
\enddemo
\bigskip

\proclaim{\stag{838-e.5} Theorem}  Like \scite{838-e.4} but dealing with
${\frak s}$, i.e. replacing clause (B) by clause (B)$'$ stated below;
that is, $\dot I(\lambda^{++},K^{\frak s}) \ge
\mu_{\text{unif}}(\lambda^{++},2^{\lambda})$ \ub{when}: 
\mr
\item "{$(A)$}"  set theoretically: 
{\roster
\itemitem{ $(a)$ }  $2^\lambda < 2^{\lambda^+} < 2^{\lambda^{++}}$ and
\sn
\itemitem{ $(b)$ }   the weak diamond ideal on $\partial :=
\lambda^+$ is not $\partial^+$-saturated
\endroster}
\sn
\item "{$(B)'$}"  model theoretic 
{\roster
\itemitem{ $(a)$ }  ${\frak s}$ is a good $\lambda$-frame (or just an
almost good $\lambda$-frame, see \scite{838-5t.3}) with 
${\frak K}_{\frak s} = {\frak K}_\lambda$ 
\sn
\itemitem{ $(b)$ }   density of $K^{3,\text{uq}}_{\frak s}$ fail,
i.e. for some $(M,N,a) \in K^{3,\text{bs}}_{\frak s}$ we have
\nl
\hskip25pt
$(M,N,a) \le_{\text{bs}}
(M',N',a) \Rightarrow (M',N',a) \notin
K^{3,\text{uq}}_{\frak s}$, see Definition \scite{838-5t.5}.
\endroster}
\endroster
\endproclaim
\bigskip

\demo{Proof}  We apply \scite{838-2b.3}, its assumption holds by 
Definition \scite{838-e.5I} and Claim \scite{838-e.5L} below applied to 
${\frak s}'$ from \scite{838-e.5D} by using \scite{838-e.5E}.
\enddemo
\bigskip

\definition{\stag{838-e.5I} Definition}  Let ${\frak s}$ be as in
\scite{838-e.5} (or just a pre-$\lambda$-frame).  We let 
${\frak u} = {\frak u}_{\frak s} = {\frak u}^1_{\frak s}$  be defined as
\mr
\item "{$(a)$}"  $\partial_{\frak u} = \lambda^+_{\frak s}$
\sn
\item "{$(b)$}"   ${\frak K}_{\frak u} = {\frak K}_{\frak s}$
\sn
\item "{$(c)$}"  FR$^{\frak u}_1 = \{(M,N,\bold I):M 
\le_{{\frak K}_{\frak u}} N,\bold I = \emptyset$ or $\bold I = \{a\}$
where $a \in N$ and $\ortp_{\frak s}(a,M,N) 
\in {\Cal S}^{\text{bs}}_{\frak s}(M)\}$
\sn
\item "{$(d)$}"  $\le^1_{\frak u}$ is defined by
$(M_1,N_1,\bold I_1) 
\le_{\frak u}^1 (M_2,N_2,\bold I_2)$ \ub{when} both are from
FR$^{\frak u}_1,M_1 \le_{\frak s} M_2,N_1 \le_{\frak s} N_2,
\bold I_1 \subseteq \bold I_2,M_1 = M_2 \cap N_1$ 
and if $\bold I_1 = \{a\}$ then
$\ortp_{\frak s}(a,M_2,N_2)$ does not fork (for ${\frak s}$) 
over $M_1$ (so if $\bold I_2 = \{a_\ell\}$ for $\ell=1,2$ this means
$(M_1,N_1,a_1) \le_{\text{bs}} (M_2,N_2,a_1)$)
\sn
\item "{$(e)$}"  FR$^{\frak u}_2 = \text{ FR}^{\frak u}_1$ and 
$\le_{\frak u}^2 = \le_{\frak u}^1$.
\endroster
\enddefinition
\bigskip

\proclaim{\stag{838-e.5L} Claim}  Let 
${\frak u} = {\frak u}_{{\frak s}'}$ where ${\frak s}'$ is from
Definition \scite{838-e.5D} or ${\frak u} = {\frak s}$ except when we
mention equality (or $=_\tau$-fuller).
\nl
1) ${\frak u}$ is a nice construction framework which is self dual.
\nl
2) For almost$_2$ every $(\bar M,\bar{\bold J},\bold f) \in
K^{\text{qt}}_{\frak u}$ the model $M := \cup\{M_\alpha:\alpha < \lambda^+\}$
is saturated, see Definition \scite{838-1a.43}(3C), see \scite{838-e.1Z}.
\nl
3) ${\frak u}$ has the weak coding property.
\nl
4) There is a ${\frak u}$-0-appropriate function ${\frak h}$ such that
every $M \in K^{{\frak u},{\frak h}}_{\lambda^{++}}$ is
$\lambda^+$-saturated above $\lambda$ and is $=_\tau$-fuller (hence
$M/=^M_\tau$ has cardinality $\lambda^{++}$).
\nl
5) Moreover, there is a ${\frak u}-\{0,2\}$-appropriate function
${\frak h}$ such that if $\langle(\bar M^\alpha,\bar{\bold
J}^\alpha,\bold f^\alpha):\alpha < \lambda^{++}\rangle$ obeys ${\frak
h}$ \ub{then} for some club $E$ of $\lambda^{++}$ the model
$M^\delta_\partial$ is saturated above $\lambda$ for $\delta \in E$
and $\cup\{M^\zeta_\partial:\zeta < \lambda^{++}\}$ is $=_\tau$-fuller.
\nl
6) Also ${\frak u}$ 
\mr
\item "{$(\alpha)$}"  satisfies (E)$_\ell$(e), 
monotonicity (see \scite{838-1a.24}(1))
\sn
\item "{$(\beta)$}"    is hereditary (see Definition \scite{838-3r.84}(2),(3))
\sn
\item "{$(\gamma)$}"   if ${\frak u} = {\frak u}_{\frak s},{\frak s}$
from \scite{838-e.5D}(2) then $=_\tau$ is a fake equality for ${\frak u}$,
(see Definition \scite{838-3r.74}(1))
\sn
\item "{$(\delta)$}"   ${\frak u}$ is hereditary for the 
fake equality $=_\tau$, (see Definition \scite{838-3r.84}(4))
\sn
\item "{$(\varepsilon)$}"   ${\frak u}$ is interpolative, see
Definition \scite{838-3r.89}.
\endroster
\endproclaim
\bigskip

\remark{\stag{838-e.5N} Remark}  1) In claim \scite{838-5t.21} we 
shall deal with the almost
good case, (see Definition \scite{838-5t.3}), the proof below serves
there too.
\nl
2) In \scite{838-e.5L}, only clause $(B)'(a)$ from the assumptions 
of Theorem \scite{838-e.5} is used except in part (3) which uses
also  clause (B)$'$(b).
\nl
3) Part (6) of \scite{838-e.5L} is used only in \scite{838-6u.34}.
\nl
4) Most parts of \scite{838-e.5L} holds also for ${\frak u} = 
{\frak u}_{\frak s}$, i.e. we have to omit the statements on
$=_\tau$-fuller, fake equality.
\endremark
\bigskip

\demo{Proof}  1) Note
\mn
\ub{Clause (D)$_\ell$(d)}:  Given $M \in K_{\frak s}$, it is not 
$<_{\frak s}$-maximal hence there is $N$ such that $M <_{\frak s} N$ hence by
density (Ax(D)(c) of (almost) good $\lambda$-frames) there is $c
\in N$ such that $\ortp_{\frak s}(c,M,N) \in 
{\Cal S}^{\text{bs}}_{\frak s}(M)$, so $(M,N,a) \in 
\text{ FR}^+_1$, as required.
\mn
\ub{Clause $(E)_\ell(c)$}:  Preservation under increasing union.

Holds by axiom $(E)(h)$ of Definition of \yCITE[1.1]{600} of ${\frak
s}$ being a good $\lambda$-frame (and similarly for being an almost good
$\lambda$-frame). 
\mn
\ub{Clause (F), amalgamation}:

This holds by symmetry axiom $(E)(i)$ of Definition \yCITE[1.1]{600}
of ${\frak s}$ being a good $\lambda$-frame (and similarly for ${\frak
s}$ being an almost good $\lambda$-frame).
The disjointness is not problematic in proving clause (F) of Definition
\scite{838-1a.3} because 
\mr
\item "{$(*)_1$}"  for ${\frak u} = {\frak u}_{\frak s}$
we can prove it (when ${\frak K}$ is categorical in $\lambda$, see 
\scite{838-5t.45}, \scite{838-5t.23} below) and 
it follows by our allowing the use of $=_\tau$ or use ${\frak s}'$
(see \scite{838-e.5L}).
\ermn
2) We just use
\mr
\item "{$(*)$}"   $M$ is saturated $(\in {\frak K}_{\lambda^+})$ when
{\roster
\itemitem{ $(a)$ }  $M = \cup\{M_\alpha:\alpha < \lambda^+\}$
\sn  
\itemitem{ $(b)$ }  $M_\alpha \in K_{\frak s}$ is $\le_{\frak
s}$-increasing continuous
\sn
\itemitem{ $(c)$ }   if $p \in {\Cal S}^{\text{bs}}_{\frak
s}(M_\alpha)$ then for some $\beta \in [\alpha,\lambda^+)$ the
non-forking extension $q \in {\Cal S}^{\text{bs}}_{\frak s}(M_\beta)$
of $p$ is realized in $M_{\beta +1}$ (or just in some
$M_\gamma,\gamma \in (\beta,\lambda^+))$.
\endroster}
\ermn
See \sectioncite[\S4]{600}; more fully see the proof of part (5).
\nl
3) Let $(M,N,a) \in K^{3,\text{bs}}_{\frak s}$ be such that there is
   no triple $(M',N',a) \in K^{3,\text{uq}}_{\frak s}$ which is
   $<_{\text{bs}}$-above it, exists by clause (B)$'$(b) from Theorem
   \scite{838-e.5}.  Let $\bold I = \{a\}$, so if $(M,N,\bold I)
   \le^1_{\frak u} (M',N',\bold I')$ then ($\bold I' = \bold I =
   \{a\}$ and) $(M',N',a) \in K^{3,\text{bs}}_{\frak s} \backslash
   K^{3,\text{uq}}_{\frak s}$ hence there are $M'',N_1,N_2$ such that
   $(M',N',\bold I) \le^1_{\frak u} (M'',N_\ell,\bold I)$ and $N_1,N_2$
   are $\le_{\frak s}$-incompatible amalgamations of $M'',N'$ over
   $M'$.  This shows that $(M',N',\bold I)$ has the true weak coding
   property.  As for almost$_2$ every triple 
$(\bar M,\bar{\bold J},\bold f) \in K^{\text{qt}}_{\frak s},M_\partial =
   M_{\lambda^+}$ is saturated, by \scite{838-2b.11}(4) and part (5) we
   get that ${\frak u}$ has the weak coding property. 
\nl
4) Easy to check by \scite{838-1a.20} or as in (5).  
\nl
5) We choose ${\frak h}$ such that:
\mr
\item "{$\boxtimes$}"   if $\bold x = \langle(\bar
M^\zeta,\bar{\bold J}^\zeta,\bold f^\zeta):\zeta \le \zeta(*)\rangle$
is $\le_{\text{qt}}$-increasing continuous and obey ${\frak h}$ is
$\xi < \zeta(*)$ \ub{then}
{\roster
\itemitem{ $(\alpha)$ }   $(\bar M^\xi,\bar{\bold J}^\xi,\bold
f^\xi) <^{\text{at}}_{\frak u} (\bar M^{\xi +1},\bar{\bold J}^{\xi +1},
\bold f^{\alpha +1})$ and let it be witnessed by $E,\bar{\bold I}$
\sn
\itemitem{ $(\beta)$ }  $M^{\xi +1}_{\delta +1}$ is brimmed over
$M^\xi_\delta$ for a club of $\delta < \lambda^+$
\sn
\itemitem{ $(\gamma)$ }  if $\zeta \le \xi$ is minimal
such that one of the cases occurs, then the demand in the first of the
cases below holds:
\endroster}
\endroster
\bn
\ub{Case A}:  There is $a \in M^\xi_\partial$ such that
$a/=^{M^\xi}_\tau$ is $\subseteq M^\zeta_\partial$ and $\zeta
< \xi$.  

\ub{Then} for some such $a',M^{\xi +1}_\partial \models ``a' =_\tau
b"$ (but $b \notin M^\zeta_\partial$), in fact $b' \in \bold
I_\alpha$ for some $\alpha < \partial$ large enough.
\bn
\ub{Case B}:  $\zeta < \xi$, not Case A (for $\zeta$) but for some
$\alpha < \partial$ and $p \in {\Cal S}^{\text{bs}}_{\frak
s}(M^\zeta_\alpha)$ for no $\varepsilon \in [\zeta,\xi)$ are there $a
\in M^{\varepsilon +1}_\partial$ such that $\ortp_{\frak
s}(a,M^\varepsilon_\beta,M^{\varepsilon +1}_\beta)$ is a non-forking
extension of $p$ for every $\beta < \partial$ large enough.  

\ub{Then} for some such
$(p,\alpha)$ we have $\ortp_{\frak s}(b,M^\zeta_\beta,M^{\zeta
+1}_\partial)$ is a non-forking extension of $p$ for every $\beta <
\partial$ large enough.
\bn
\ub{Case C}:  $\zeta < \xi$, Cases A,B fail for $\zeta$ 
and there is a pair $(\alpha,p)$ such
that $\alpha < \partial,p \in {\Cal S}^{\text{bs}}_{\frak
s}(M^\zeta_\alpha)$ such that for no $\varepsilon \in [\zeta,\xi)$ is
the set $S := \{\delta < \partial$: there is $i < \bold f^{\varepsilon
+1}(\delta)$ such that $\ortp_{\frak s}(a_{\bold J^{\varepsilon
+1}_{\delta +i}},M^{\varepsilon +1}_{\delta +i},M^{\varepsilon
+1}_{\delta +i+1})$ is a non-forking extension of $p\}$ stationary.

\ub{Then} for some such pair $(\alpha,p)$, the condition above holds for $\xi$.
\bn
\ub{Case D}:  $\zeta = \xi$.

Does not matter.
\nl
6) Easy, too.  \hfill$\square_{\scite{838-e.5L}}$ 
\enddemo
\bn
\centerline {$* \qquad * \qquad *$}
\bn
\ub{(F) \quad The better versions of the results}:

Here we prove the better versions of the results, i.e. without using
on ``WDmId$(\lambda^+)$ is $\lambda^{++}$-saturated" but relying on
later sections.

Of course, the major point is reproving the results of \S4(E),
i.e. ``non-structure for a good $\lambda$-frame ${\frak s}$ failing the
density of $K^{3,\text{uq}}_{\frak s}$", we have to rely on \S5-\S8.

We also deal with \S4(D); here we rely on \sectioncite[\S8]{E46}, so we get an
almost good $\lambda$-frame ${\frak s}$ (rather than good
$\lambda$-frames).  But in \S5-\S8 we deal also with this more general
case (and in \S7, when we discard a non-structure case, we prove
that ${\frak s}$ is really a good $\lambda$-frame).

Lastly, we revisit \S4(C).
\bigskip

\proclaim{\stag{838-e.6A} Theorem}  1) In Theorem \scite{838-e.5} we can omit
the assumption (A)(b).
\nl
2) $\dot I(\lambda^{++},{\frak K}) \ge
   \mu_{\text{unif}}(\lambda^{++},2^{\lambda^+})$ and moreover $\dot
   I(\lambda^{++},K^{{\frak u},{\frak h}}) \ge
   \mu_{\text{unif}}(\lambda^{++},2^{\lambda^+})$ for ${\frak u} =
{\frak u}^1_{\frak s}$ from Definition \scite{838-e.5I} or ${\frak u} =
   {\frak u}^3_{\frak s}$ from Definition \scite{838-8h.4} and 
any ${\frak u}-\{0,2\}$-appropriate function ${\frak h}$, \ub{when}:
\mr
\item "{$(A)$}"  (set theoretic assumption), $2^\lambda < 2^{\lambda^+} <
2^{\lambda^{++}}$
\sn
\item "{$(B)$}"  (model theoretic assumptions), 
{\roster
\itemitem{ $(a)$ }  ${\frak s}$ is an almost good $\lambda$-frame
\sn
\itemitem{ $(b)$ }  ${\frak s}$ is categorical in $\lambda$
\sn
\itemitem{ $(c)$ }  ${\frak s}$ is not a good $\lambda$-frame \ub{or}
${\frak s}$ (is a good $\lambda$-frame which) fail density for
$K^{3,\text{uq}}_{\frak s}$.
\endroster}
\endroster
\endproclaim
\bigskip

\remark{\stag{838-e.6C} Remark}  1) This proves \yCITE[0z.1]{E46} from
 \yCITE[8h.9]{E46}, proving the main theorem \yCITE[0z.1]{E46}.
\nl
2) We can phrase the theorem also as: if (A),(B)(a),(B)(b) holds and
the $\dot I(K^{{\frak u},{\frak h}}_{\lambda^{++}}) <
   \mu_{\text{unif}} (\lambda^{++},2^\lambda)$ \ub{then} ${\frak s}$
   is a good $\lambda$-frame for which $K^{3,\text{uq}}_{\frak s}$ is
dense in $(K^{3,\text{bs}}_{\frak s},\le_{\text{bs}})$ (and so has
   existence for $K^{3,\text{uq}}_{\frak s}$). 
\endremark
\bigskip

\demo{Proof}  1) This is a special case of part (2).
\nl
2) Toward contradiction assume that the desired conclusion fail.

First, the Hypothesis \scite{838-5t.1} of \S5 holds for ${\frak s}$ hence
its results.  Second, the Hypothesis \scite{838-6u.1} of \S6 apply hence
its results.  So consider conclusion \scite{838-6u.27}(2); its
assumption ``$2^{\lambda^+} < 2^{\lambda^{++}}$" holds by assumption
(A) here, and its assumption ``$\dot I(K^{{\frak s},
{\frak h}}_{\lambda^{++}}) <
\mu_{\text{unif}}(\lambda^{++},2^{\lambda^+})$ for ${\frak u} =
{\frak u}^1_{\frak s}$ for some ${\frak u}-\{0,2\}$-appropriate ${\frak h}$"
holds by our present assumption toward contradiction and its
assumption ``${\frak K}_{\frak s}$ is categorical" holds by clause
(B)(b) of the assumption of \scite{838-e.6A}.

Hence the conclusion of \scite{838-6u.27} holds which says that
\mr
\item "{$(*)$}"  ${\frak s}$ has existence for
$K^{3,\text{up}}_{{\frak s},\xi}$ for every $\xi \le \lambda^+$, see
Definition \scite{838-6u.7}.
\ermn
Now consider Hypothesis \scite{838-7v.1}; now part (1) there (${\frak s}$
is an almost good $\lambda$-frame) holds by the present assumption
(B)(a), part (2) there was just proven; part (3) there (${\frak s}$ is
categorical in $\lambda$) holds by the present assumption (B)(b), and
lastly, part (4) there (disjointness) is proved in \scite{838-5t.45}.  So
Hypothesis \scite{838-7v.1} of \S7 holds hence the results of that section
up to \scite{838-7v.38} apply.

In particular, WNF$_{\frak s}$ defined in \scite{838-7v.5}(1),(2) is well
defined and by \scite{838-7v.33}(1)
is a weak non-forking relation on ${}^4(K_{\frak s})$ respecting
${\frak s}$.  Also ${\frak s}$ is a good $\lambda$-frame by Lemma
\scite{838-7v.37}(1) so the first possibility in clause (B)(c) of
\scite{838-e.6A} does not hold.  By inspection all parts of Hypothesis \scite{838-8h.1}
of \S8 holds hence the results of that section apply.

Now in Claim \scite{838-8h.61}, its conclusion fails as this means our
assumption toward contradiction and among its assumptions, clause (a),
saying ``$2^\lambda < 2^{\lambda^+} < 2^{\lambda^{++}}$" holds by
clause (A) of \scite{838-e.6A}, clause (c) saying ``$K$ is categorical in
$\lambda$" holds by clause (B)(b) of \scite{838-e.6A} and clause (d)
saying ``${\frak u} = {\frak u}_{\frak s}^3$ from \scite{838-8h.4} has
existence for $K^{3,\text{up}}_{{\frak s},\lambda^+}$" was proved
above.  So clause (b) of \scite{838-8h.61} fails, i.e. ${\frak s}$ fails
the non-uniqueness for WNF$_{\frak s}$, but by \scite{838-8h.34}(1) this
implies that we have uniqueness for WNF.

Lastly, we apply Observation \scite{838-8h.34}(2), it has two assumptions,
the first ``${\frak s}$ has existence for $K^{3,\text{up}}_{{\frak
s},\lambda^+}"$, was proved above, and second ``${\frak s}$ has
uniqueness for WNF" has just been proved; so the conclusion of
\scite{838-8h.34} holds.  This means ``${\frak s}$ has existence for
$K^{3,\text{uq}}_{\frak s}"$, so also the second possibility of clause
(B)(c) of \scite{838-e.6A} fails; a contradiction.
 \hfill$\square_{\scite{838-e.6A}}$
\enddemo
\bigskip

\proclaim{\stag{838-e.6H} Theorem}  1) In Theorem \scite{838-e.4} we can omit the
assumption (A)(b) at least if $K$ is categorical in $\lambda^+$.
\nl
2) $\dot I(\lambda^{++},K^{{\frak u},{\frak h}}) \ge
   \mu_{\text{unif}}(\lambda^{++},2^{\lambda^+})$ \ub{when}:
\mr
\item "{$(A)$}"  (set theoretic) $2^\lambda < 2^{\lambda^+} <
   2^{\lambda^+}$
\sn
\item "{$(B)$}"  (model theoretic) as in \scite{838-e.4}, but ${\frak K}$
categorical in $\lambda^+$
\sn
\item "{$(C)$}"  ${\frak u} = u^4_{{\frak K}_\lambda}$, see Definition
\scite{838-e.4.1}, ${\frak h}$ is a ${\frak u}-\{0,2\}$-appropriate function.
\endroster
\endproclaim
\bigskip

\remark{Remark}  This theorem is funny as \sectioncite[\S6]{E46} and in
particular \yCITE[6f.13]{E46} is a shortcut, but we prove this by a
detour (using \sectioncite[\S8]{E46}) so in a sense \scite{838-e.6H} is less
natural than \scite{838-e.4}; but no harm done.
\endremark
\bigskip

\demo{Proof}  1) By part (2) and \scite{838-e.4.3}(4) recalling
\scite{838-e.4.0}(4).
\nl
2) Toward contradiction, assume that the desired conclusion
fails.  By \yCITE[8h.1]{E46} there is an almost good
$\lambda$-frame ${\frak s}$ such that ${\frak K}_{\frak s} = {\frak
K}_\lambda$ and ${\Cal S}^{\text{bs}}_{\frak s}(M)$ is the set of
minimal $p \in {\Cal S}^{\text{bs}}_{\frak s}(M)$.

Note that categoricity in $\lambda^+$ is used in \chaptercite{E46} to
deduce the stability in $\lambda$ for minimal types and the set of minimal
types in ${\Cal S}_{\frak K}(M)$ being inevitable,
but this is assumed in clause (B)(d) of the assumption of \scite{838-e.4}, 
so natural to conjecture that it is not needed, see \chaptercite{E46}.

Now using the meaning of the assumption (B)(e) of Theorem \scite{838-e.4}
is that ``$K^{3,\text{uq}}_{\frak s}$ is not dense in
$(K^{3,\text{bs}}_{\frak s},\le_{\text{bs}})$" so we can apply Theorem
\scite{838-e.6A} to get the desired result.
\nl
${{}}$  \hfill$\square_{\scite{838-e.6H}}$
\enddemo
\bigskip

\proclaim{\stag{838-e.6L} Theorem}  1) In Theorem \scite{838-e.3} we can weaken
the set theoretic assumption, omitting the extra assumption
(a)$(\beta)$.
\nl
2) $\dot I(\aleph_2,K^{{\frak u},{\frak h}}) \ge
\mu_{\text{unif}}(\aleph_2,2^{\aleph_1})$ \ub{when}:
\mr
\item "{$(a)$}"  (set theory) $2^{\aleph_0} < 2^{\aleph_1} < 2^{\aleph_2}$
\sn
\item "{$(b)-(e)$}"  as in \scite{838-e.4}
\sn
\item "{$(f)$}"  ${\frak u} = {\frak u}^3_{{\frak K}_{\aleph_0}}$ from
Definition \scite{838-e.3.4} and 
${\frak h}$ is a ${\frak u}-\{0,2\}$-appropriate function.
\endroster
\endproclaim
\bn
\margintag{e.6L.1}\ub{\stag{838-e.6L.1} Discussion}:  1) This completes a promise from
\sectioncite[\S5]{88r}.  You may say that once we prove in
\sectioncite[\S3]{600}(B) that ${\frak s} = {\frak s}^1_{\frak K}$ is a
good $\aleph_0$-frame we do not need to deal with ${\frak K}$ any more, so
no need of \scite{838-e.6L}.  In addition to keeping promises this is only
partially true because of the following.
\nl
2) First, arriving to ${\frak s}^+$, see \sectioncite[\S1]{705}, we do not
know that $\le_{{\frak s}(+)} = \le_{\frak K} \rest 
{\frak K}_{{\frak s}(+)}$, because this is proved only if ${\frak s}$
is good$^+$ (see \sectioncite[\S1]{705}).  Now by looking at the
definitions (and \yCITE[Ex.1]{600}, equality of the various types), we
know that ${\frak s}$ being good$^+$ is equivalent to the symmetry
property, i.e. every one sided stable
amalgamation.  We prove that its failure implies non-structure in 
\scite{838-e.6L.3}, \scite{838-e.6L.7}, \scite{838-e.6L.12} below.
\nl
3) Another point is that even if ${\frak s}$ is weakly successful 
(i.e. we have existence for $K^{3,\text{uq}}_{\frak s}$), we
 can define NF = NF$_{\frak s}$ and so 
we have unique non-forking amalgamation, it is
not clear that this is equal to the one/two sided stable
amalgamation from \xCITE{88r}. 
\nl
4) Also defining ${\frak s}^{+n}$ as in \xCITE{600} we may hope not
to shrink $K^{{\frak s}(+n)}$, i.e. to get all the 
$(\aleph_0,n)$-properties (as in \cite{Sh:87b}).  If we start with
$\psi \in \Bbb L_{\omega_1,\omega}(\bold Q)$ as in \cite{Sh:48} this
seems straight, in general, this is a priori not clear, hopefully see
\cite{Sh:F888}.
\nl
5) Concerning (2) above, we like to use \scite{838-3r.71} - \scite{838-3r.88}
   in the proof as in the proof of \scite{838-6u.23}.  If we have used ${\frak u} =
   {\frak u}^3_{\frak R}$ from Definition \scite{838-e.3.4}, this fails,
   e.g. it is not self dual.  We can change (FR$_2,\le_2$) to make it
   symmetric but still it will fail ``hereditary", so it is natural to
   use ${\frak u}_2$ defined in \scite{838-e.6L.3} below, but then we
   still need $(M_\delta,N_\delta,\bold I_\delta) \in \text{
   FR}^1_{{\frak u}_2}$ to ensure $N_\delta$ is $(\bold
   D(M_\delta),\aleph_0)^*$-homogeneous over $M_\delta$.  This can be
   done by using the game version of the coding property.  This is
   fine but was not our ``main road" so rather we use the theorem on
   ${\frak u}_3$ but use ${\frak u}_5$ to apply \S3.

A price of using \S3 is having to use fake equality.  Also together
with symmetry, we deal with lifting free $(\alpha,0)$-rectangles.
\nl
6) To complete the proof of \scite{838-e.6L}, by \scite{838-e.6L.12} it
   suffices to prove the uniqueness of two-sided stable amalgamation.
   We use \S8 and toward this we define WNF$_*$, prove that it is a
   weak non-forking relation of ${\frak K}_{\aleph_0}$ respecting
   ${\frak s}$, using the ``lifting" from \S5.  Then we can apply \S8.
\nl
7) A drawback of \scite{838-e.6L} as well as \scite{838-e.3} and
\sectioncite[\S3]{600}(B) is that we restrict ourselves to a countable
$\bold D$.  Now in \xCITE{88r} this is justified as it is proved
that for some increasing continuous sequence $\langle \bold
D_\alpha:\alpha < \omega_1\rangle$ with each $\bold D_\alpha$
countable, $\bold D = \cup\{\bold D_\alpha:\alpha < \omega_1\}$,
i.e. for every $M \in K_{\aleph_0}$, the sequence $(\bold
D_\alpha(M):\alpha < \omega_1\rangle$ is an increasing sequence of
sets of types with union $\bold D(M)$.  However, from the positive
results on every $\bold D_\alpha$ we can deduce positive results on
$\bold D$.  See, hopefully, in \cite{Sh:F888}.
\bigskip

\remark{\stag{838-e.6L.2} Remark}  1) Assumption (d) of \scite{838-e.3} gives: usually
$(M,N,\bold I) \in \text{ FR}^1_{\frak u}$ has non-uniqueness,
i.e. when $\bold I = ({}^{\omega >} N) \backslash ({}^{\omega >} M)$.  
We like to work as in \scite{838-e.6A}.
\nl
2) So as indirectly there we would like to use \scite{838-3r.91}; for this we need
the vertical uq-invariant whereas we naturally get failure of
   the semi uq-invariant coding property.  So we would like to quote
   \scite{838-3r.88} but this requires ${\frak u}$ to be self dual.
\nl
3) Hence use also a relative of ${\frak u}$ from \scite{838-e.6L.17}, 
for it we prove the implication and from this deduce what we need for the old.
\nl
4) Our problem is to prove that ${\frak s} = {\frak s}_{\aleph_0}$ is
   good$^+$, equivalently prove the symmetry property, this is done in
   Claim \scite{838-e.6L.12}.  It is natural to apply \scite{838-3r.67} -
   \scite{838-3r.91}.
\nl
5) The proof of \scite{838-e.6L} will come later.
\endremark
\bigskip

\definition{\stag{838-e.6L.3} Definition}  In \scite{838-e.6L} we define
${\frak u}_5 = {\frak u}^5_{{\frak K}_{\aleph_0}}$ as follows ($\ell$ is $1,2$)
\mr
\item "{$(a)$}"  $\partial = \partial_{\frak u} = \aleph_1$
\sn
\item "{$(b)$}"  ${\frak K}_{\frak u} = {\frak K}_{\aleph_0}$ more
pedantically ${\frak K}_{\frak u} = {\frak K}'_{\aleph_0}$
\sn
\item "{$(c)_1$}"  FR$^{\frak u}_1$ is the class of triples
$(M,N,\bold I)$ such that $\bold I \subseteq {}^{\omega >} N
\backslash {}^{\omega >} M$
\sn
\item "{$(c)_2$}"  $(M_1,N_1,\bold I_1) \le_1 (M_1,N_2,\bold I_2)$ iff
both are from FR$^{\frak u}_1,M_1 \le_{\frak K} M_2,N_1 \le_{\frak K}
N_2$ and $\bar c \in \bold I \Rightarrow \text{ gtp}(\bar c,M_2,N_2)$
is the stationarization of gtp$(\bar c,M_1,N_1)$
\sn
\item "{$(d)$}"  FR$_2 = \text{ FR}_1$ and $\le^2_{\frak u} =
\le^1_{\frak u}$.
\ermn
Now we have to repeat various things.
\enddefinition
\bigskip

\proclaim{\stag{838-e.6L.7} Claim}  1) ${\frak u}_5$ is a nice construction
framework.
\nl
2) For almost$_2$ every triples $(\bar M,\bar{\bold J},\bold f) \in
 K^{\text{qt}}_{{\frak u}_5}$ the model $M_\partial = M_{\lambda^+}$ belongs to
   $K_{\lambda^+}$ and is saturated.
\nl
3) ${\frak u}_5$ has fake equality $=_\tau$ and is monotonic, see
   Definition (\scite{838-1a.24}(1)), and weakly hereditary for the fake
   equality $=_\tau$, see Definition \scite{838-3r.84}(5) and
   interpolative (see Definition \scite{838-3r.89}).
\endproclaim
\bigskip

\demo{Proof}  Should be clear (e.g. part (2) as in \scite{838-6u.3}).  
\hfill$\square_{\scite{838-e.6L.7}}$
\enddemo
\bigskip

\proclaim{\stag{838-e.6L.12} Claim}  $\dot I(\lambda^{++},{\frak K}) \ge
\mu_{\text{unif}}(\aleph_2,2^{\aleph_0})$ and moreover $\dot
   I(\aleph_2,{\frak K}(\aleph_1-\text{saturated})) \ge
   \mu_{\text{unif}}(\aleph_2,2^{\aleph_2})$ \ub{when}:
\mr
\item "{$\circledast$}"  $(a)(\alpha),(b),(c),(e)$ from \scite{838-e.3} and
\sn
\item "{${{}}$}"  $(d)''(\alpha) \quad {\frak K}$ fails the symmetry
   property \ub{or}
\sn
\item "{${{}}$}"  $\quad\,\,\, (\beta) \quad {\frak K}$ fails the lifting
property, see Definition \scite{838-e.6L.17} below.
\endroster
\endproclaim
\bigskip

\definition{\stag{838-e.6L.17} Definition}  We define a 4-place relation
WNF$_*$ on $K_{\aleph_0}$ as follows:\nl
WNF$_*(M_0,M_1,M_2,M_3)$
\ub{when}
\mr
\item "{$(a)$}"  $M_\ell \in K_{\aleph_0}$ for $\ell \le 3$
\sn
\item "{$(b)$}"   $M_0 \le_{\frak K} M_\ell \le_{\frak K} M_3$ for
$\ell=1,2$
\sn
\item "{$(c)$}"   for $\ell = 1,2$ if $\bar a \in {}^{\omega
>}(M_\ell)$ then gtp$(\bar a,M_{3-\ell},M_3)$ is the stationarization
 of gtp$(\bar a,M_0,M_3) = \text{ gtp}(\bar a,M_0,M_\ell)$.
\endroster
\enddefinition
\bigskip

\definition{\stag{838-e.6L.19} Definition}  We say that $({\frak
 K},\text{WNF}_*)$ has the lifting property when WNF$_*$ satisfies
 clause (g) of Definition \scite{838-7v.35}, see the proof of
 \scite{838-7v.35}, i.e. if 
{\rm WNF}$_*(M_0,N_0,M_1,N_1)$ and $\alpha < \lambda^+$ and
$\langle M_{0,i}:i \le \alpha\rangle$ is $\le_{\frak s}$-increasing
continuous, $M_{0,0} = M_0$ and $N_0 \le_{{\frak K}_\lambda},
M_{0,\alpha}$ \ub{then} we can find a
$\le_{{\frak K}_\lambda}$-increasing continuous sequence $\langle M_{1,i}:i \le
\alpha +1 \rangle$ such that $M_{1,0} = M_1,N_1 \le_{\frak K}
M_{1,\alpha +1}$ and for each $i < \alpha$ we have {\rm WNF}$_*
(M_{0,i},M_{0,i+1},M_{1,i},M_{1,i+1})$ for $i < \alpha$.
\enddefinition
\bigskip

\demo{Proof of \scite{838-e.6L.12}}  We start as in the proof 
of \scite{838-e.3}, choosing the 
good $\aleph_0$-frame ${\frak s} = {\frak s}_{\aleph_0}$ and define
   ${\frak u} = {\frak u}^3_{\frak K}$ as there, (except having the fake
inequality which causes no problem), so it is a nice
   construction framework by \scite{838-e.3.8}(1) and for almost$_2$ all
   triples $(\bar M,\bar{\bold J},\bold f) \in K^{\text{qt}}_{\frak
   u}$ the model $M \in {\frak K}_{\aleph_1}$ is saturated (by
   \scite{838-e.3.8}(2)).

Now Theorem \scite{838-3r.91} gives the right conclusion, so to suffice to
verify its assumptions.  Of course, ${\frak u}$ is as required in
Hypothesis \scite{838-3r.0}.

Clause (a) there means $2^{\aleph_0} < 2^{\aleph_0} < 2^{\aleph_2}$ (as
$\partial_{\frak u} = \aleph_1$ and we choose $\theta = \aleph_0$),
which holds by clause $(a)(\alpha)$ of the present claim.

Clause (c) there says that for $\{0,2\}$-almost every $(\bar M,\bar{\bold
J},\bold f) \in K^{\text{qt}}_{\frak u}$ the model $M_\partial \in
K_{\aleph_2}$ is $K_{\frak u}$-model homogeneous; this holds and can
be proved as in \scite{838-6u.3}.

We are left with clause (b), i.e. we have to prove that some
$(M,N,\bold I) \in \text{ FR}^{\frak u}_1$ has the vertical
uq-invariant coding property, see Definition \scite{838-3r.19}.  Choose
$(M,N,\bold I) \in \text{ FR}^{\frak u}_1$ such that $|\bold I| > 1$,
hence $N$ is $(\bold D(M),N)^*$-homogeneous and $\bold I = ({}^{\omega
>}N) \backslash ({}^{\omega >}M)$ and we shall prove that it has the
vertical uq-invariant coding, so assume
\mr
\item "{$(*)$}"  $\bold d_0$ is a ${\frak u}$-free $(\alpha_{\bold
d},0)$-rectangle satisfying $M \le_{\frak s} M' = M^{\bold d}_{0,0}$
and $M^{\bold d}_{\alpha(\bold d),0} \cap N = M$.
\ermn
We have to find $\bold d$ as required in Definition
\scite{838-3r.19}.

Note that by the choice of $(M,N,\bold I)$, and the assumption
``${\frak K}$ fails the symmetry property" we can find $(M_*,N_*)$ and
then $\bar c$
\mr
\item "{$(*)_2$}"  $(a) \quad M \le_{\frak s} M_* \le_{\frak s} N_*$
and $N \le_{\frak s} N_*$
\sn
\item "{${{}}$}"  $(b) \quad M,N,M_*,N_*$ is in one-sided stable
amalgamation, i.e. if $\bar b \in {}^{\omega >} N$
\nl

\hskip25pt  then gtp$(\bar b,M_*,N_*)$ is the 
stationarization of gtp$(\bar b,M,N)$
\sn
\item "{${{}}$}"  $(c) \quad M,M_*,N,N_*$ is not in one sided stable
amalgamation, so
\sn
\item "{${{}}$}"  $(c)^+ \quad \bar c \in {}^{\omega >}(M_*)$ and
gtp$(\bar c,N,N_*)$ is not the stationarization of gtp$(\bar
c,M,M_*)$.
\ermn
We like to apply the semi version, i.e. Definition \scite{838-3r.71} and
Claim \scite{838-3r.88}.  There are technical difficulties so we apply it
to ${\frak u}_5$, see  \scite{838-e.6L.3}, \scite{838-e.6L.7} above and in the
end increase the models to have the triples in FR$^1_{\frak u}$ and
use \scite{838-3r.90} instead of \scite{838-3r.88}, so
all should be clear.  

Alternatively, works only with ${\frak u}_5$ but use the game version
of the coding theorem.  
\nl
${{}}$   \hfill$\square_{\scite{838-e.6L.12}}$
\enddemo
\bigskip

\proclaim{\stag{838-e.6L.21} Claim}  1) If (${\frak K}$,{\rm WNF}$_*$) has lifting,
see Definition \scite{838-e.6L.17}, \ub{then} {\rm WNF}$_*$ is a weak non-forking
relationg of ${\frak K}_{\aleph_0}$ respecting ${\frak s}$ with
disjointness (\scite{838-7v.35}(3)).
\nl
2) {\rm WNF}$_*$ is a pseudo non-forking relation of 
${\frak K}_{\aleph_0}$ respecting ${\frak s}$ meaning clauses (a)-(f)
with disjointness, see the proof \ub{or} see Definition \scite{838-7v.35}(4),(3).
\nl
3) If ${\frak K}_{\aleph_0}$ satisfies symmetry then in clause (c) of
   Definition \scite{838-e.6L.17}, it is enough if it holds for one $\ell$. 
\endproclaim
\bigskip

\demo{Proof}  1) We should check all the clauses of Definition
   \scite{838-7v.35}, so see \yCITE[nf.0X]{600} or the proof of
   \scite{838-7v.33}(1).
\mn
\ub{Clause (a)}:  WNF$_*$ is a 4-place relation on ${\frak
   K}_{\aleph_0}$.
\nl
[Why?  By Definition \scite{838-e.6L.17}, in particular clause (a).]
\mn
\ub{Clause (b)}:  WNF$_*(M_0,M_1,M_2,M_3)$ implies $M_0 \le_{\frak K}
M_\ell \le_{\frak K} M_3$ for $\ell=1,2$ and is preserved by
isomorphisms.
\nl
[Why?  By Definition \scite{838-e.6L.17}, in particular clause (c).]
\mn
\ub{Clause (c)}:  Monotonicity
\nl
[Why?  By properties of gtp, see \yCITE[5.13]{88r}(1).]
\mn
\ub{Clause (d)}:  Symmetry
\nl
[Why?  Read Definition \scite{838-e.6L.17}.]
\mn
\ub{Clause (e)}:  Long Transitivity

As in the proof of \yCITE[5.19]{88r}.
\mn
\ub{Clause (f)}:  Existence

This is proved in \yCITE[5.22]{88r}.
\mn
\ub{Clause (g)}:  Lifting, see Definition \scite{838-e.6L.17}.

This holds by an assumption.

WNF$_*$ \ub{respects} ${\frak s}$ and has
\ub{disjointness}.  

Clear by the definition (in particular of gtp).
\nl
2) The proof included in the proof of part (1).
\nl
3) Should be clear.  \hfill$\square_{\scite{838-e.6L.21}}$
\enddemo
\bigskip

\demo{Proof of \scite{838-e.6L}}  Let $\lambda = \aleph_0$ and toward
contradiction assume that $\dot I(\lambda^{++},{\frak
K}(\lambda^+$-saturated)) $<
\mu_{\text{unif}}(\lambda^{++},2^\lambda)$.

As in the proof of \scite{838-e.3}, by \yCITE[Ex.1]{600} ${\frak s} :=
{\frak s}_{\aleph_0}$ is a good $\lambda$-frame categorical in
$\lambda$.  By Theorem \scite{838-e.6A}, recalling our assumption toward
contradiction, $K^{3,\text{uq}}_{\frak s}$ is dense in
$K^{3,\text{bs}}_{\frak s}$ hence ${\frak s}$ has existence for
$K^{3,\text{uq}}_{\frak s}$; i.e. is weakly successful, but we shall
not use this. 

By \scite{838-e.6L.12} we know that ${\frak K}$ has the symmetry property,
hence the two-sided stable amalgamation fails uniqueness.  Also by
\scite{838-e.6L.12} we know that it has the lifting property, so by
\scite{838-e.6L.21}, \scite{838-e.6L.12} we know that WNF$_*$ is a 
weak non-forking relation on
${\frak K}_{\aleph_0}$ which respects ${\frak s}$, so Hypothesis
\scite{838-8h.1} holds.

Let ${\frak u}$ be defined as in \scite{838-8h.4} (for our given ${\frak
s}$ and WNF$_*$).  Now we try to apply Theorem \scite{838-8h.61}.  Its
conclusion fails by our assumption toward contradiction and clause
(a),(b),(c) there holds.  So clause (b) there fails so by \scite{838-8h.34}(1).
So we can conclude that we have uniqueness for WNF$_*$ by
\scite{838-8h.34}(2) clearly ${\frak s}$ has existence for
$K^{3,\text{bs}}_{\frak s}$, i.e. is weakly successful.

So ${\frak s}^+$ is a well defined good $\lambda^+$-frame, see
\xCITE{705}.  By \sectioncite[\S8]{600}, \xCITE{705} and our assumption
toward contradiction,  we know that ${\frak s}$ is
successful.  Now if ${\frak s}$ is not good$^+$ then ${\frak K}$ fails
the symmetry property hence by \scite{838-e.6L.12} we get contradiction, so
ncessarily ${\frak s}$ is good$^+$ hence we have $\le_{{\frak s}(+)} =
\le_{\frak s} \rest K_{{\frak s}(+)}$.  This proves that the saturated
$M \in {\frak K}_{\lambda^+}$ is super limit (see \xCITE{705} also 
this is \yCITE[5.24]{88r}).  \hfill$\square_{\scite{838-e.6L}}$
\enddemo
\goodbreak

\head {\S5 On almost good $\lambda$-frames} \endhead  \resetall \sectno=5
 \spuriousreset
\bigskip

Accepting ``WDmId$(\partial)$ is not $\partial^+$-saturated" where
$\partial = \lambda^+$ we have accomplished in \S4 the applications we
promised.  Otherwise for \sectioncite[\S5]{600} we have to prove for a good
$\lambda$-frame ${\frak s}$ that among the triples $(M,N,a) \in
K^{3,\text{bs}}_{\frak s}$, the ones with uniqueness are not dense,
as otherwise non-structure in $\lambda^{++}$ follows.  
Toward this (in this section) 
we have to do some (positive structure side) work  
which may be of self interest.  Now in
the case we get ${\frak s}$ in \sectioncite[\S3]{600} starting  from 
\cite{Sh:576} or better from \sectioncite[\S8]{E46}, but with 
$\lambda_{\frak s} = \lambda$ rather than
$\lambda_{\frak s} = \lambda^+$, we have to start with an 
almost good $\lambda$-frames ${\frak s}$ rather than with 
good $\lambda$-frames.  However, there is a price: for 
eliminating the non-$\partial^+$-saturation of
the weak diamond ideal and for using the ``\ub{almost} $\lambda$-good"
version, we have to work more.

First, we shall not directly try to prove density 
of uniqueness triples $(M,N,\bold J)$ but just the density of 
poor relatives like $K^{3,\text{up}}_{\frak s}$.

Second, we have to prove some positive results, particularly in the
almost good $\lambda$-frame case.  This is done here in \S5 and more is done
in \S7 assuming existence for $K^{3,\text{up}}_{{\frak s},\lambda^+}$
justified by 
the non-structure result in \S6 and the complimentary
full non-structure result is proved in \S8.
\bigskip

\demo{\stag{838-5t.1} Hypothesis}  ${\frak s}$ is an almost good
$\lambda$-frame (usually categorical in $\lambda$) and for
transparency ${\frak s}$ has disjointness,
see Definitions \scite{838-5t.3}, \scite{838-5t.9} below; the disjointness 
is justified in the Discussion
\scite{838-5t.11} and not used in \scite{838-5t.41} - \scite{838-5t.49} which in
fact prove it and let $\partial = \lambda^+$.
\enddemo
\bigskip

\definition{\stag{838-5t.3} Definition}  ``${\frak s}$ is an almost good
$\lambda$-frame" is defined as in \yCITE[1.1]{600} except that we
weaken (E)(c) to (E)(c)$^-$ and strengthen (D)(d) to (D)(d)$^+$ where
(recall $\ortp_{\frak s} = \ortp_{{\frak K}_{\frak s}})$:
\mn
\ub{Ax(E)(c)$^-$}:  the local character 
\sn
\ub{if} $\langle M_i:i \le \delta +1\rangle$ is
$\le_{\frak s}$-increasing continuous and the set
$\{i < \delta:N_i <^*_{\frak s} N_{i+1}$, i.e. $N_{i+1}$ is 
universal over $N_i\}$ is unbounded in $\delta$ \ub{then} for some 
$a \in M_{\delta +1}$ the type 
$\ortp_{\frak s}(a,M_\delta,M_{\delta +1})$ belongs to ${\Cal
S}^{\text{bs}}_{\frak s}(M_\delta)$ and 
does not fork over $M_i$ for some $i < \delta$
\sn
Ax$(D)(d)^+ \quad$ if $M \in K_{\frak s}$ then ${\Cal S}_{{\frak
K}_{\frak s}}(M)$ has cardinality $\le \lambda$ 
\nl
(for good $\lambda$-frame this holds by \yCITE[4a.1]{600}).
\enddefinition
\bn
As in \xCITE{600}
\definition{\stag{838-5t.5} Definition}  1) $K^{3,\text{bs}}_{\frak s}$
is the class of triples $(M,N,a)$ such that $M \le_{\frak s} N$ and $a
\in N \backslash M$.
\nl
2) $\le_{\text{bs}} = \le^{\text{bs}}_{\frak s}$ is the following
two-place relation (really partial order) 
on $K^{3,\text{bs}}_{\frak s}$.  We let 
$(M_1,N_1,a_1) \le_{\text{bs}} (M_2,N_2,a_2)$ iff $a_1=a_2,M_1
\le_{\frak s} M_2,N_1 \le_{\frak s} N_2$ and $\ortp_{\frak
s}(a_1,N_1,N_2)$ does not fork over $M_1$.
\enddefinition
\bigskip

\proclaim{\stag{838-5t.7} Claim}  1) $K^{3,\text{bs}}_{\frak s}$ and
$\le_{\text{bs}}$ are preserved by isomorphisms.
\nl
2) $\le_{\text{bs}}$ is a partial order on $K^{3,\text{bs}}_{\frak s}$.
\nl
3) If $\langle(M_\alpha,N_\alpha,a):\alpha < \delta\rangle$ is
$\le_{\text{bs}}$-increasing and $\delta < \lambda^+$ is a limit
ordinal and $M_\delta := \cup\{M_\alpha:\alpha < \delta\},N_\delta :=
\cup\{N_\alpha:\alpha < \delta\}$ \ub{then} $\alpha < \delta
\Rightarrow (M_\alpha,N_\alpha,a) \le_{\text{bs}}
(M_\delta,N_\delta,a) \in K^{3,\text{bs}}_{\frak s}$ (using Ax(E)(h)).
\endproclaim
\bigskip

\demo{Proof}  Easy.  \hfill$\square_{\scite{838-5t.7}}$
\enddemo
\bigskip

\definition{\stag{838-5t.9} Definition}   We say ${\frak s}$ has
disjointness or is disjoint \ub{when}:
\mr
\item "{$(a)$}"  strengthen Ax(C), i.e. ${\frak K}_{\frak s}$ 
has disjoint amalgamation
which means that: if $M_0 \le_{\frak s} M_\ell$ for $\ell=1,2$ and $M_1 \cap
M_2 = M_0$ \ub{then} for some $M_{\frak s} \in K_{\frak s}$ we have
$M_\ell \le_{\frak s} M_{\frak s}$ for $\ell=0,1,2$
\sn
\item "{$(b)$}"  strengthen Ax(E)(i) by disjointness: if above we
assume in addition that $(M_0,M_\ell,a_\ell) \in
K^{3,\text{bs}}_{\frak s}$ for $\ell=1,2$ \ub{then} we can add 
$(M_0,M_\ell,a_\ell) \le_{\text{bs}} 
(M_{3 - \ell},M_{\frak s},a_\ell)$ for $\ell=1,2$.
\endroster
\enddefinition
\bn
\margintag{5t.11}\ub{\stag{838-5t.11} Discussion}:  How ``expensive" is the (assumption of)
disjoint amalgamation (in Ax(C) and Ax(E)(i))?
\nl
1)   We can ``get it for free" by using 
${\frak K}'$ and ${\frak s}'$, see Definition
\scite{838-1a.19} and \scite{838-e.5D}(2) so we assume it.
\nl
2) Alternatively we can prove it assuming categoricity in $\lambda$
(see \scite{838-5t.45} which relies on \scite{838-5t.43}).
\nl
3)  So usually we shall ignore this point.
\bigskip

\demo{\stag{838-5e.23} Exercise}  There is a good $\lambda$-frame without
disjoint amalgamation.
\enddemo
\bn
[Hint:  Let
\mr
\item "{$\circledast_1$}"  $(a) \quad \tau = \{F\},F$ a unary function
\sn
\item "{${{}}$}"   $(b) \quad \psi$ the first order sentence $(\forall
x,y)[F(x)\ne x \wedge F(y) \ne y \rightarrow F(x) =$
\nl

\hskip25pt $F(y)] \wedge ((\forall x)[F(F(x)) = F(x)]$
\sn
\item "{${{}}$}"  $(c) \quad K = \{M:M$ is a 
$\tau$-model of $\psi\}$, so $M \in
K \Rightarrow |\varphi(M)| \le 1$ where we
\nl

\hskip25pt  let $\varphi(x) = (\exists y)(F(y) = x \wedge y \ne x)$
\sn
\item "{${{}}$}"  $(d) \quad M \le_{\frak K} N'$ iff 
$M \subseteq N$ are from $K$
and $\varphi(M) = \varphi(N) \cap M$.
\ermn
Now note
\mr
\item "{$(*)_1$}"  ${\frak K} := (K,\le_{\frak K})$ is an a.e.c. with
LS$({\frak K}) = \aleph_0$
\sn
\item "{$(*)_2$}"  ${\frak K}$ has amalgamation.
\ermn
[Why?  If $M_0 \le_{\frak K} M_\ell$ for $\ell=1,2$, then separate the
proof to three cases: the first when $\varphi(M_0) = \varphi(M_1) =
\varphi(M_2) = \emptyset$ the second when $|\varphi(M_1)| +
|\varphi(M_2)| = 1$ and $\varphi(M_0) = \emptyset$; 
in the third case $\varphi(M_0) =
\varphi(M_1) = \varphi(M_2)$ is a singleton.]

So
\mr
\item "{$(*)_3$}"   ${\frak K}_\lambda = (K_\lambda,\le_{\frak K}
\restriction K_\lambda)$.
\ermn
Now we define ${\frak s}$ by letting
\mr
\item "{$(*)_4$}"  $(a) \quad {\frak K}_{\frak s} = {\frak
K}_\lambda$
\sn
\item "{${{}}$}"  $(b) \quad K^{3,\text{bs}}_{\frak s} := 
\{(M,N,a):M \le_{{\frak K}_\lambda} N,a \in N \backslash M 
\backslash \varphi(N)\}$
\sn
\item "{${{}}$}"  $(c) \quad$ for $M_1 \le_{\frak K} M_2 
\le_{\frak K} M_3$ and $a \in M_3$ we say 
$\ortp_{\frak K}(a,M_2,M_3)$ does not
\nl

\hskip25pt  fork over $M_1$ iff $a \notin M_2 \and F^{M_3}(a) \notin
M_2 \backslash M_1$.
\ermn
Lastly
\mr
\item "{$(*)_5$}"  ${\frak s}$ is a good $\lambda$-frame.
\ermn
[Why?  Check.  E.g.
\bn
\ub{Ax(D)(c)}:  Density

So assume $M <_{\frak s} N$ now if there is $a \in N \backslash M
\backslash \varphi(M)$ then a $\ortp(a,M,N) \in {\Cal
S}^{\text{bs}}_{\frak s}(M)$ so $a$ is as required.  Otherwise, as $M
\ne M$ necessary $\varphi(N)$ is non-empty and $\subseteq N \backslash
M$, let it be $\{b\}$.  By the definition of $\varphi$ there is $a \in
N$ such that $F^N(a) = b \wedge a \ne b$ so necessarily $a \notin M$ and is as
required.]
\bn
\ub{Ax(E)(e)}:  Uniqueness

The point is that:
\mr
\item "{$(*)_6$}"  if $\varphi(M) = \emptyset$ then
${\Cal S}^{\text{bs}}_{\frak s}(M)$ contains just two types $p_1,p_2$ such
that if $p_\ell = \ortp(a,M,N) \Rightarrow$ then $\ell = 1 \Rightarrow 
F^N(a) = a \in N \backslash M$ and $\ell = 2 \Rightarrow F^N(a) \in N
\backslash M \backslash \{a\}$
\sn
\item "{$(*)_7$}"  if $\varphi(M) = \{b\}$ then
${\Cal S}^{\text{bs}}_{\frak s}(M)$ contains just two types $p_1,p_2$ such
that $p_1$ is as above and $p_2 = \ortp(a,M,N) \Rightarrow F^N(a) = b$
\sn
\item "{$(*)_8$}"  there are $M_0 \le_{\frak K} M_1 = M_2$ such that
$\varphi(M_1) \ne \emptyset = \varphi(M_0)$ and let $\varphi(M_\ell) =
\{b\}$ for $\ell=1,2$ so $b \in M_\ell \backslash M_0$.  
So we cannot disjointly amalgamate $M_1,M_2$ over $M_0$.
\ermn
[Why?  Think.]
\nl
So we are done with Example \scite{838-5e.23}.]
\bn
Recalling \yCITE[0.21]{600},
          \yCITE[0.22]{600}:
\proclaim{\stag{838-5t.15} Claim}  1) For $\kappa = \text{\rm cf}(\kappa)
\le \lambda$
\mr
\item "{$(a)$}"   there is a $(\lambda,\kappa)$-brimmed $M \in {\frak
K}_{\frak s}$, in fact $(\lambda,\kappa)$-brimmed over $M_0$ for any
pregiven $M_0 \in {\frak K}_{\frak s}$
\sn
\item "{$(b)$}"  $M$ is unique up to isomorphism over $M_0$ (but we fix
$\kappa$)
\sn
\item "{$(c)$}"  if $M \in {\frak K}_{\frak s}$ is
$(\lambda,\kappa)$-brimmed over $M_0$ \ub{then} it is $\le_{\frak
s}$-universal over $M_0$.
\ermn
2) So the superlimit $M \in {\frak K}_{\frak s}$ is
$(\lambda,\kappa)$-brimmed for every $\kappa \le \text{\rm cf}(\kappa)
\le \lambda$ hence is brimmed.
\nl
3) If $\kappa = \text{\rm cf}(\kappa) \le \lambda$ and $M_1
\le_{\frak s} M_2$ are both $(\lambda,\kappa)$-brimmed and $\Gamma
\subseteq {\Cal S}^{\text{bs}}_{\frak s}(M_2)$ has cardinality $<
\kappa$ and every $p \in \Gamma$ does not fork over $M_1$ \ub{then}
there is an isomorphism $f$ from $M_2$ onto $M_1$ such that $p \in
\Gamma \Rightarrow f(p) = p \restriction M_1$.
\endproclaim
\bigskip

\demo{Proof}  1) By Definition \yCITE[0.21]{600} and Claim
\yCITE[0.22]{600} because ${\frak K}_\lambda$ has amalgamation, the
JEP recalling \yCITE[1.1]{600} and has no 
$<_{{\frak K}_\lambda}$-maximal member (having a superlimit model).
\nl
2) By the definition of being brimmed and ``superlimit in ${\frak
K}_{\frak s}$" which exists by ``${\frak s}$ is almost good
$\lambda$-frame".
\nl
3) Exactly as in the proof of \yCITE[stg.9]{705}.
\hfill$\square_{\scite{838-5t.15}}$  
\enddemo
\bigskip

\remark{\stag{838-5t.16} Remark}  1) It seems that there is no great harm in 
weakening (E)(h) to (E)(h)$^-$ as in (E)(c)$^-$, but also no urgent
need, where:.

Ax(E)(h)$^-$: assume
$\langle M_i:i \le \delta\rangle$ is $\le_{\frak s}$-increasing
continuous and $\delta = \sup\{i:M_{i+1}$ is $\le_{\frak s}$-universal
over $M_i\}$.  If $p \in {\Cal S}_{{\frak K}_{\frak s}}(M_\delta)$ and
$i < \delta \Rightarrow p \restriction M_i \in {\Cal
S}^{\text{bs}}_{\frak s}(M_i)$ \ub{then} $p \in {\Cal
S}^{\text{bs}}_{\frak s}(M_\delta)$.
\nl
2) That is, if we weaken Ax(E)(h) then we are drawn to further
problems.  After defining ${\frak u}$, does 
$\le_\ell$-increasing sequence $\langle(M_i,N_i,\bold J_i):i <
\delta\rangle$ has the union as a $\le_\ell$-upper bound?  If
$M_{i+1}$ is universal over $M_i$ for $i < \delta$ this is O.K.,
but using triangles in the limit we have a problem; see part (4) below.
\nl
3) Why ``no urgent need"?  The case which draws us to consider
Ax(E)(c)$^-$ is \sectioncite[\S8]{E46}, i.e. by the ${\frak s}$ derived
there satisfies Ax(E)(h).  So we may deal
with it elsewhere, \cite{Sh:F841}.
\nl
4) When we deal with ${\frak u}$ derived from such ${\frak s}$,
i.e. as in part (1) we may demand:
\mr
\item "{$(A)$}"  First, dealing with ${\frak u}$-free rectangles 
and triangles we add
{\roster
\itemitem{ $(a)$ }  $\bold I^{\bold d}_{i,j} = \emptyset$ when $j$ is
a limit ordinal
\sn
\itemitem{ $(b)$ }  $\bold J^{\bold d}_{i,j} = \emptyset$ when $i$ is
a limit ordinal
\sn
\itemitem{ $(c)$ }  it is everywhere universal (each $M_{i+1,j+j}$ is
$\le_{\frak s}$-universal over 
$M^{\bold d}_{i+1,j} \cup M^{\bold d}_{i,j+1}$ (or at least
each)
\endroster}
\item "{$(B)$}"  similarly with $K^{\text{qt}}_{\frak s}$,
i.e. defining $(\bar M,\bar{\bold J},\bold f) \in K^{\text{qt}}_{\frak
s}$ we add the demands
{\roster
\itemitem{ $(a)$ }  $\bold J_\delta = \emptyset$ for limit $\delta$
\sn
\itemitem{ $(b)$ }   if $\delta \in S \cap E$, then $\bold d$, 
the $({\frak u}_{\frak s})$-free $(\bold f(\delta),0)$-rectangle
$(\langle M_{\delta + i}:i \le \bold f(\delta)\rangle,\langle \bold
J_{\delta +i}:i < \bold f(\delta)\rangle)$ is strongly full (defined
as below) and so for any $\theta \le \lambda$, if $\lambda|j$ and 
$i < j \le \bold f(\delta)$ and cf$(\delta) = \theta$ then $M_{\delta,j}$ is
$(\lambda,\theta)$-brimmed over $M_{\delta +i}$.
\endroster}
\item "{$(C)$}"  similarly for $\le^{\text{at}}_{{\frak u}_{\frak
s}},\le^{\text{qr}}_{{\frak u}_{\frak s}}$.
\endroster
\endremark
\bigskip

\definition{\stag{838-5t.19} Definition}  1) We define ${\frak u} = {\frak
u}_{\frak s} = {\frak u}^1_{\frak s}$ as in Definition 
\scite{838-e.5I}, it is denoted by ${\frak u}$ for
this section so may be omitted and we may write ${\frak s}$-free
instead of ${\frak u}_{\frak s}$-free.
\nl
2) We say $(M,N,\bold I) \in \text{ FR}_\ell$ realizes $p$ \ub{when}
$\ortp_{\frak s}(a_{\bold I},M,N) = p$, recalling $\bold I =
\{a_{\bold I}\}$.
\enddefinition
\bigskip

\proclaim{\stag{838-5t.21} Claim}  1) ${\frak u}_{\frak s}$ is a nice
construction framework which is self-dual.
\nl
2) Also ${\frak u}_{\frak s}$ is monotonic and hereditary and interpolative.
\endproclaim
\bigskip

\remark{\stag{838-5t.23} Remark}  1) Here we use ``${\frak s}$ has
disjointness" proved in \scite{838-5t.45}.
\nl
2) Even without \scite{838-5t.45}, if ${\frak s} = {\frak s}'_1$ for 
some almost good $\lambda$-frame ${\frak s}_1$ 
then ${\frak s}$ has disjointness.
\nl
3) Mostly it does not matter if we use ${\frak u}^1_{\frak s},
{\frak s}'$ from \scite{838-e.5D} (see \scite{838-6u.31}) but in proving
\scite{838-6u.21}(1), the use of ${\frak u}^1_{{\frak s}'}$ is preferable;
alternatively in defining nice construction framework we waive the 
disjointness, which is a cumbersome but not serious change.
\endremark
\bigskip

\demo{Proof}  As in \scite{838-e.5L} except disjointness which holds by
Hypothesis \scite{838-5t.3} and is justified by \scite{838-5t.11} above, (or see 
\scite{838-5t.45} below).   \hfill$\square_{\scite{838-5t.21}}$
\enddemo
\bigskip

\remark{\stag{838-5t.25} Remark}  Because we assume on ${\frak s}$ only
that it is an almost good $\lambda$-frame we
have to be more careful as (E)(c) may fail, in particular in
proving brimmness in triangles of the right kind.
I.e. for a ${\frak u}$-free $(\bar \alpha,\beta)$- triangle $\bold d$, 
we need that in the ``vertical sequence", $\langle M^{\bold
d}_{i,\beta}:i \le \alpha_\beta\rangle$ the highest model
$M^{\bold d}_{\alpha_\beta,\beta}$ is brimmed over the lowest
$M^{\bold d}_{0,\beta}$.  This motivates the following.
\endremark
\bigskip

\definition{\stag{838-5t.27} Definition}  1) We say $\bar M$ is a
$(\Gamma,\delta)$-correct sequence \ub{when}: $\Gamma \subseteq {\Cal
S}^{\text{bs}}_{\frak s}(M_\delta)$, the sequence $\bar M =
\langle M_\alpha:\alpha \le \alpha(*)\rangle$ is $\le_{\frak
s}$-increasing continuous, $\delta \le \alpha(*)$ and: \ub{if} $N \in
K_{\frak s}$ satisfies $M_\delta <_{\frak s} N$ 
\ub{then} for some $c \in N \backslash M_\delta$ and
$\alpha < \delta$ the type $\ortp_{\frak s}(c,M_\delta,N)$ does not fork over
$M_\alpha$ and belongs to $\Gamma$.
\nl
2) We omit $\delta$ when this holds for every limit $\delta \le \alpha(*)$.
\nl
3) We say $\Gamma$ is $M$-inevitable when $\Gamma \subseteq {\Cal
S}^{\text{bs}}_{\frak s}(M)$ and: if $M <_{\frak s} N$
\ub{then} some $p \in {\Cal S}^{\text{bs}}_{\frak s}(M) \cap \Gamma$ is
realized in $N$.
\nl
4) Using a function ${\Cal S}^*$ instead of $\Gamma$ we mean we use
$\Gamma = {\Cal S}^*(M_\delta)$.
\nl
5) We may omit $\Gamma$ (and write $\delta$-correct) above \ub{when}
$\Gamma$ is ${\Cal S}^{\text{bs}}_{\frak s}$.
\nl
6) For $\bar M = \langle M_\alpha:\alpha \le \alpha(*)\rangle$ let
correct$_\Gamma(\bar M) = \{\delta \le \alpha(*):\bar M$ is 
$(\delta,\Gamma)$-correct so $\delta$ is a limit ordinal$\}$ and we
may omit $\Gamma$ if $\Gamma = {\Cal S}^{\text{bs}}_{\frak
s}(M_{\alpha(*)})$.
\nl
7) If $\bold d$ is a ${\frak u}_{\frak s}$-free $(\bar
   \alpha,\beta)$-triangle let $\Gamma_{\bold d} = \{p \in {\Cal
   S}^{\text{bs}}_{\frak s}(M_{\alpha_\beta,\beta}):p$ does not fork
   over $M_{i,j}$ for some $j < \beta,i < \alpha_j\}$.
\enddefinition
\bigskip

\definition{\stag{838-5t.29} Definition}  1) We say that 
$\bold d$ is a brimmed (or universal) 
${\frak u}_{\frak s}$-free or ${\frak s}$-free triangle \ub{when}:
\mr
\item "{$(a)$}"  $\bold d$ is a ${\frak u}_{\frak s}$-free triangle
\sn
\item "{$(b)$}"  if $i < \alpha_j(\bold d)$ and $j < \beta(\bold d)$
\ub{then} $M^{\bold d}_{i+1,j+1}$ is brimmed (or universal) over
$M^{\bold d}_{i+1,j}$.
\ermn
2) We say strictly brimmed (universal) \ub{when} also 
$M^{\bold d}_{i+1,j+1}$ is
brimmed/universal over $M^{\bold d}_{i+1,j} \cup M^{\bold d}_{i,j+1}$ when 
$j < \beta,i < \alpha_j(\bold d)$.  
\nl
2A) We say that $\bold d$ is a weakly brimmed (or weakly universal)
${\frak u}_{\frak s}$-free or ${\frak s}$-free triangle \ub{when}:
\mr
\item "{$(a)$}"  ${\bold d}$ is a ${\frak u}_{\frak s}$-free triangle
\sn
\item "{$(b)$}"  if $j_1 < \beta_{\bold d},i_1 < \alpha_j(\bold d)$ then we
can find a pair $(i_2,j_2)$ such that $j_1 \le j_2 < \beta_{\bold d},i_1 \le
i_2 < \alpha_j({\bold d})$ and $M^{\bold d}_{i_2+1,j_2+1}$ is brimmed (or
is $\le_{\frak s}$-universal) over $M^{\bold d}_{i_2,j_2}$ or just over
$M^{\bold d}_{i_1,j_1}$.
\ermn
2B) We say that $\bold d$ is a weakly brimmed (weakly universal)
${\frak u}_{\frak s}$-free rectangle \ub{when} it and its dual (see Definition
\scite{838-1a.13}(3)) are weakly brimmed ${\frak u}_{\frak s}$-free
triangles.  Similarly for brimmed, strictly brimmed, universal,
strictly universal.
\nl
3) We say that a ${\frak u}_{\frak s}$-free triangle 
$\bold d$ is full \ub{when}:
\mr
\item "{$(a)$}"  $\lambda_{\frak s}$ divides 
$\alpha_{\beta(\bold d)}(\bold d)$ and $\beta(\bold d)$ is a limit ordinal and
$\bar \alpha$ is continuous or just $\alpha_{\beta(\bold d)}(\bold d) 
= \cup\{\alpha_j(\bold d):j<\beta\}$
\sn
\item "{$(b)$}"  if $i < \alpha_j(\bold d),j < \beta := \beta(\bold d)$ and $p
\in {\Cal S}^{\text{bs}}_{\frak s}(M_{i,j})$ \ub{then}:
{\roster
\itemitem{ $(*)$ }  the following subset $S^{\bold d}_p$ of
$\alpha_\beta(\bold d)$ has order type $\ge \lambda_{\frak s}$ where
$S^{\bold d}_p := \{i_1 < \alpha_\beta(\bold d)$: for some $j_1 \in
(j,\beta)$ we have $i \le i_1 < \alpha_{j_1}(\bold d)$ and for some $c \in 
\bold J_{i_1,j_1}$ the type 
$\ortp_{\frak s}(c,M_{i_1,j_1},M_{i_1+1,j_1})$ is a 
non-forking extension of $p\}$.
\endroster}
\ermn
3A) We say that the ${\frak u}_{\frak s}$-free triangle 
is strongly full \ub{when}:
\mr
\item "{$(a)$}"  as in part (3)
\sn
\item "{$(b)$}"  if $i < \alpha_j(\bold d)$ and $j < \beta(\bold d)$
and $p \in {\Cal S}^{\text{bs}}_{\frak s}(M^{\bold d}_{i,j})$ 
then for $\lambda$ ordinals $i_1 \in [i,i+\lambda)$ for some 
$j_1 \in [j_2,\beta)$ the type $\ortp_{\frak s}
(b^{\bold d}_{i_1,j_1},M^{\bold d}_{i_1,j_1}
M^{\bold d}_{i_1 + 1,j_1})$ is a non-forking extension of $p$ and $i_1
< \alpha_{j_1}(\bold d)$, of course.
\ermn
3B) We say a ${\frak u}_{\frak s}$-free rectangle 
$\bold d$ is full [strongly full] \ub{when} both
$\bold d$ and its dual are full [strongly full] 
${\frak u}_{\frak s}$-free triangles.
\enddefinition
\bigskip

\demo{\stag{838-5t.31} Observation}  1) If $\bar M = \langle M_\alpha:\alpha
\le \delta\rangle$ is $\le_{\frak s}$-increasing continuous, $\delta$
a limit ordinal and $\delta = \sup\{\alpha < \delta:M_{\alpha +1}$ or
just $M_\beta$ for some $\beta \in (\alpha,\delta)$, is $\le_{\frak
s}$-universal over $M_\alpha\}$ \ub{then} $\delta \in 
\text{ correct}(\bar M)$.
\nl
2) If cf$(\delta_\ell) = \kappa$ and $\bar M^\ell = \langle
M^\ell_\alpha:\alpha \le \delta_\ell\rangle$ is $\le_{\frak
s}$-increasing continuous and $h_\ell:\kappa \rightarrow \delta_\ell$
is increasing with $\delta_\ell = \sup(\text{Rang}(h_\ell))$ for
$\ell=0,1$ and $\varepsilon < \kappa \Rightarrow
M^1_{h_1(\varepsilon)} = M^2_{h_2(\varepsilon)}$ \ub{then} $\bar M^1$ is
$\delta_1$-correct iff $\bar M^2$ is $\delta_2$-correct.
\nl
3) Instead of the $h_1,h_2$ is part (2) it suffices that $(\forall
\alpha < \delta_\ell)(\exists \beta < \delta_{3 - \ell})(M^\ell_\alpha
\le_{\frak s} M^{3-\ell}_\beta)$ for $\ell=1,2$.  Also instead of
``$\delta_\ell$-correct" we can use $(\Gamma,\delta_\ell)$-correct.
\nl
4) If $\bar M = \langle M_\alpha:\alpha \le \delta\rangle$ is
   $\le_{\frak s}$-increasing continuous and $\Gamma_{\bar M} = \{p
   \in {\Cal S}^{\text{bs}}_{\frak s}(M_\delta):p$ does not fork over
   $M_\alpha$ for some $\alpha < \delta\}$ \ub{then} $\bar M$ is
   $\delta$-correct iff $\bar M$ is $(\delta,\Gamma)$-correct.
\enddemo
\bigskip

\demo{Proof}  1) By Ax(E)(c)$^-$ and the definition of correct$(\bar M)$.
\nl
2),3),4)  Read the definitions.    \hfill$\square_{\scite{838-5t.31}}$ 
\enddemo
\bigskip

\demo{\stag{838-5t.33} Observation}   Assume $\bold d^{\text{ver}}$ 
is a ${\frak u}_{\frak s}$-free $(\alpha,0)$-rectangle,
$\bold d_{\text{hor}}$ is a ${\frak u}_{\frak s}$-free 
$(0,\beta)$-rectangle and
$M^{\bold d^{\text{ver}}}_{0,0} = M^{\bold d_{\text{hor}}}_{0,0}$.
\ub{Then} there is a pair $(\bold d,f)$ such that:
\mr
\item "{$(a)$}"  $\bold d$ is a strictly brimmed 
${\frak u}$-free $(\alpha,\beta)$-rectangle
\sn
\item "{$(b)$}"  $\bold d \restriction (0,\beta) = \bold
d_{\text{hor}}$
\sn
\item "{$(c)$}"  $f$ is an isomorphism from $\bold d^{\text{ver}}$
onto $\bold d \restriction (\alpha,0)$ over $M^{\bold d}_{0,0}$
\sn
\item "{$(d)$}"  if $\alpha$ is divisible  by $\lambda$ and
$S^{\text{ver}} = \{\alpha' < \alpha:
\bold J^{\bold d^{\text{ver}}}_{\alpha'} = \emptyset\}$ is an
unbounded subset of $\alpha$ of order-type divisible by $\lambda$ (so
in particular $\alpha$ is a limit ordinal)
\ub{then} we can add $\bold d$ as a triangle is full; 
we can add strongly full if
$\alpha' < \alpha \Rightarrow \lambda = |S^{\text{ver}} \cap [\alpha',\alpha' +
\lambda]|$
\sn
\item "{$(e)$}"  if $S_{\text{hor}} := \{\beta' < \beta:
\bold I^{\text{hor}}_\beta = 0\}$ is an
unbounded subset of $\beta$ of order-type divisible by $\lambda$
\ub{then} we can add ``{\rm dual}$(\bold d)$ as a 
triangle is full"; and we can add strongly full if
$\beta' < \beta \Rightarrow \lambda = |S_{\text{hor}} \cap 
[\beta',\beta' + \lambda]|$.
\endroster
\enddemo
\bigskip

\demo{Proof}  Easy.  \hfill$\square_{\scite{838-5t.33}}$
\enddemo
\bn
\margintag{5t.39}\ub{\stag{838-5t.39} Exercise}:  Show the obvious implication concerning
the notions from Definition \scite{838-5t.29}.  Let $\bold d$ be a ${\frak
u}$-free $(\bar \alpha,\beta)$-triangle and $\bold e$ be a ${\frak
u}$-free $(\alpha,\beta)$-rectangle.
\nl
1) $\bold d$ strictly brimmed implies $\bold d$ is brimmed which
   implies $\bold d$ is weakly brimmed.
\nl
2) Like (1) replacing brimmed by universal.
\nl
3) If $\bold d$ is strictly brimmed/brimmed/weakly brimmed \ub{then}
   $\bold d$ is strictly universal/universal/weakly universal.
\nl
4) If $\bold d$ is strongly full \ub{then} $\bold d$ is full.
\nl
5) Similarly for the rectangle $\bold e$.
\nl
6) If $\bold e$ is strictly brimmed/brimmed/weakly brimmed \ub{then}
so is {\rm dual}$(\bold e)$.
\nl
7) If $\bold e$ is strictly universal/universal/weakly universal 
\ub{then} so is {\rm dual}$(\bold e)$.
\bigskip

\proclaim{\stag{838-5t.36} The Correctness Claim}  1) Assume 
$\delta < \lambda^+$ is a
limit ordinal, $\bar M^\ell = \langle M^\ell_\alpha:\alpha \le
\delta\rangle$ is $\le_{\frak s}$-increasing continuous sequence for
$\ell=1,2$ and $\alpha < \delta \Rightarrow M^1_\alpha \le_{\frak s}
M^2_\alpha$ and $M^1_\delta = M^2_\delta$.  If $\bar M^1$ is
$\delta$-correct \ub{then} $\bar M^2$ is $\delta$-correct.
\nl
2) $M_\delta$ is $(\lambda,\text{\rm cf}(\delta))$-brimmed over $M_0$;
   moreover over $M_i$ for any $i < \delta$ \ub{when}:
\mr
\item "{$(a)$}"  $\delta$ is a limit ordinal divisible by $\lambda$
(the divisibility follows by clause (d))
\sn
\item "{$(b)$}"  $\bar M = \langle M_\alpha:\alpha \le \delta\rangle$
is $\le_{\frak s}$-increasing continuous
\sn
\item "{$(c)$}"  $\bar M$ is $(\delta,\Gamma)$-correct, so $\Gamma
\subseteq \{p \in {\Cal S}^{\text{bs}}_{\frak s}(M_\delta):p$ does not
fork over $M_\alpha$ for some $\alpha < \delta\}$, if $\Gamma$ is
equal this means
$\delta \in \text{\rm correct}(\bar M)$, recalling
Definition \scite{838-5t.27}(1),(6)
\sn 
\item "{$(d)$}"  if $\alpha < \delta$ and $p \in \{q \rest M_\alpha:q
\in \Gamma\} \subseteq
{\Cal S}^{\text{bs}}_{\frak s}(M_\alpha)$ \ub{then}: for $\ge \lambda$
ordinals $\beta \in (\alpha,\delta)$ \ub{there} is
$c \in M_{\beta +1}$ such that $\ortp_{\frak s}(c,M_{\beta},M_{\beta +1})$ is a
non-forking extension of $p$.
\ermn
2A) $M_\delta$ is $(\lambda,\text{\rm cf}(\delta))$-brimmed 
over $M_0$ \ub{when} clauses (a),(b) of part (2) holds and
$\delta  = \sup(S)$ where $S = \{\delta':\delta' < \delta$ 
and $\bar M \restriction (\delta' +1)$ satisfies clauses (a)-(d) from
part (2)$\}$.
\nl
3) Assume $\bold d$ is a ${\frak u}_{\frak s}$-free $(\bar
\alpha,\beta)$-triangle, $\beta$ is a limit ordinal, $\bar\alpha$ is
continuous (or just $\alpha_\beta = \alpha_\beta(\bold d) 
= \cup\{\alpha_{j+1}:j <  \beta\})$ 
and $\bar\alpha \restriction \beta$ is not eventually constant,
\mr
\item "{$(a)$}"  if $\bold d$ is brimmed or just weakly universal
\ub{then} $\langle M^{\bold d}_{\alpha,\beta}:\alpha 
\le \alpha_\beta\rangle$ is $\alpha_\beta$-correct; moreover is
correct for $(\alpha_\beta,\Gamma_{\bold d})$ recalling $\Gamma_{\bold
d} = \{p \in {\Cal S}^{\text{bs}}_{\frak s}(M^{\bold
d}_{\alpha_\beta,\beta}):p$ does not fork over $M^{\bold d}_{i,j}$ for
some $j < \beta,i < \alpha_j\}$
\sn
\item "{$(b)$}"  if $\bold d$ is weakly universal and full \ub{then} 
$M_{\alpha_\beta,\beta}$ is brimmed over $M_{i,\beta}$ for every
$i < \alpha_\beta$.
\ermn
3A) Assume $\bold d$ is an ${\frak u}_{\frak s}$-free 
$(\alpha,\beta)$-rectangle
\mr
\item "{$(a)$}"  if $\bold d$ is  brimmed or just weakly universal
(see Definition \scite{838-5t.29}(2A)) and {\rm cf}$(\alpha) = 
\text{\rm cf}(\beta) \ge \aleph_0$ 
\ub{then} $\langle M^{\bold d}_{i,\beta}:i \le
\alpha\rangle$ is $\alpha$-correct and even $(\alpha,\Gamma_{\bold
d})$-correct 
\sn
\item "{$(b)$}"  if in clause (a), $\bold d$ is full \ub{then} $M^{\bold
d}_{\alpha,\beta}$ is $(\lambda,\text{\rm cf}(\alpha))$-brimmed over
$M_{i,\beta}$ for $i < \alpha$
\sn
\item "{$(c)$}"  if $\bold d$ is strictly brimmed (or just strictly
universal, see Definition \scite{838-5t.29}(2))
and strongly full, (see Definition \scite{838-5t.29}(3A),(3B)) 
and $\lambda^2 \omega$ divides $\alpha$ (but no requirement on the
cofinalities) \ub{then} 
$M^{\bold d}_{\alpha,\beta}$ is $(\lambda,\text{\rm cf}(\alpha))$-brimmed over
$M^{\bold d}_{i,\beta}$ for every $i < \alpha$.
\ermn
4) For $M \in K_{\frak s}$ there is $N \in K_{\frak s}$ which is
brimmed over $M$ and is unique up to isomorphism over $M$ (so in other
words, if $M_\ell$ is $(\lambda,\kappa_\ell)$-brimmed over $M$ for
$\ell=1,2$ \ub{then} $N_1,N_2$ are isomorphic over $M$).
\endproclaim
\bigskip

\demo{Proof}  1) Assume $M^2_\delta <_{\frak s} N$, hence $M^1_\delta 
<_{\frak s} N$ so as $\delta \in \text{ correct}(\bar M^1)$
necessarily for some pair $(p,\alpha)$ we have: $\alpha < \delta$ and
$p \in {\Cal S}^{\text{bs}}_{\frak s}(M^1_\delta)$ is realized in $N$
and does not fork over $M^1_\alpha$.  As $M^1_\alpha \le_{\frak s}
M^2_\alpha \le_{\frak s} M^2_\delta = M^1_\delta$ and montonicity of
non-forking it follows that $``p \in {\Cal S}^{\text{bs}}_{\frak
s}(M^2_\delta)$ does not fork over $M^2_\alpha"$, and, of course, $p$
is realized in $N$.  So $(p,\alpha)$ are as required in the definition
of ``$\delta \in \text{ correct}(\bar M^2)$".
\nl
2) Similar to \sectioncite[\S4]{600} but we give a full self-contained proof.
   The ``moreover" can be proved by renaming.

Let $\langle {\Cal U}_\alpha:\alpha < \delta\rangle$ be an
increasing continuous sequence of subsets of $\lambda$ such that
$|{\Cal U}_0| = \lambda,|{\Cal U}_{\alpha +1} \backslash 
{\Cal U}_\alpha| = \lambda$.  We choose a triple $(\bar{\bold
a}^\alpha,N_\alpha,f_\alpha)$ by induction on $\alpha \le \delta$ such that:
\mr
\item "{$\oplus$}"  $(a) \quad N_\alpha \in {\frak K}_{\frak s}$ and
$N_0 = M_0$
\sn
\item "{${{}}$}"  $(b) \quad f_\alpha$ is a $\le_{\frak s}$-embedding
of $M_\alpha$ into $N_\alpha$
\sn
\item "{${{}}$}"  $(c) \quad \langle N_\beta:\beta \le \alpha\rangle$
is $\le_{\frak s}$-increasing continuous
\sn
\item "{${{}}$}"  $(d) \quad \langle f_\beta:\beta \le \alpha\rangle$
is $\subseteq$-increasing continuous and $f_0 = \text{ id}_{M_0}$
\sn
\item "{${{}}$}"  $(e) \quad \bar{\bold a}^\alpha = \langle a_i:i
\in{\Cal U}_\alpha\rangle$, so $\beta < \alpha \Rightarrow 
\bar{\bold a}^\beta = \bar{\bold a}^\alpha \restriction {\Cal
U}_\beta$
\sn
\item "{${{}}$}"  $(f) \quad \bar{\bold a}^\alpha$ lists the elements of
$N_\alpha$ each appearing $\lambda$ times
\sn
\item "{${{}}$}"  $(g) \quad$ if $\alpha = \beta +1$ then $N_\alpha$
is $\le_{\frak s}$-universal over $N_\beta$
\sn
\item "{${{}}$}"  $(h) \quad$ if $\alpha = \beta +1$ and ${\Cal
W}_\beta : = \{i \in {\Cal U}_\beta$: for some $c \in M_{\alpha +1}
\backslash M_\alpha$ we 
\nl

\hskip25pt have $f_\beta(\ortp_{\frak s}(c,M_\beta,M_\alpha)) = \ortp_{\frak
s}(a_i,f_\beta(M_\beta),N_\beta)\}$ is not empty
\nl

\hskip25pt  and $i_\beta :=
\text{ min}({\Cal W}_\beta)$ then $a_{i_\beta} \in \text{
Rang}(f_\alpha)$.
\ermn
There is no problem to carry the definition; and by clauses (c),(g) of
$\oplus$ obviously
\mr
\item "{$\odot$}"  $N_\delta$ is $(\lambda,\text{cf}(\delta))$-brimmed
over $N_0$ (hence over $f_0(M_0))$.
\ermn
Also by renaming without loss of generality  $f_\alpha = \text{ id}_{M_\alpha}$ for
$\alpha \le \delta$ hence $M_\delta \le_{\frak s} N_\delta$.

Now if $M_\delta = N_\delta$ then by $\odot$ we are done.  Otherwise by clause
(c) of the assumption $\bar M$ is $(\delta,\Gamma)$-correct
hence by Definition \scite{838-5t.27}(1), for some $c \in N_\delta
\backslash f_\delta(M_\delta)$ and $\alpha_0 < \delta$ the type $\ortp_{\frak
s}(c,M_\delta,N_\delta)$ belongs to $\Gamma \subseteq
{\Cal S}^{\text{bs}}_{\frak s}(M_\delta)$ 
and does not fork over $M_{\alpha_0}$.  As
$\langle N_\beta:\beta \le \delta\rangle$ is $\le_{\frak
s}$-increasing continuous, for some $\alpha_1 < \delta$ we have $c \in
N_{\alpha_1}$, hence for some $i_* \in {\Cal U}_{\alpha_1}$ we have $c =
a_{i_*}$.  So $\alpha_2 := \text{ max}\{\alpha_0,\alpha_1\} < \delta$
and by clause (d) of the assumption the set ${\Cal W} := \{\alpha <
\delta:\alpha \ge \alpha_2$ and for some $c' \in M_{\alpha +1}$ the type 
$\ortp_{\frak s}(c',M_\alpha,M_{\alpha +1})$ is a non-forking
extension of $\ortp_{\frak s}(c,M_{\alpha_0},N_\delta)\}$
has at least $\lambda$ members and is $\subseteq [\alpha_2,\delta)$
and by the monotonicity and uniqueness
properties of non-forking we have ${\Cal W} = 
\{\alpha < \delta:\alpha \ge \alpha_2$ and some 
$c' \in M_{\alpha +1}$ realizes $\ortp_{\frak s}
(c,M_\alpha,N_\delta)$ in $M_{\alpha +1}\}$.  Now for every
$\alpha \in {\Cal W} \subseteq [\alpha_2,\delta)$ the set ${\Cal
W}_\alpha$ defined in clause (h) of $\oplus$ above
is not empty, in fact, $i_* \in {\Cal W}_\alpha$ hence $\beta \in
{\Cal W} \subseteq [\alpha_2,\delta) \Rightarrow i_\beta = 
\text{ min}({\Cal W}_\beta) \le i_*$ but $|{\Cal W}| = \lambda$, so 
by cardinality consideration for some $\beta_1 < \beta_2$ from 
${\Cal W}$ we have $i_{\beta_1} = i_{\beta_2}$ but $a_{i_{\beta_1}} \in 
\text{ Rang}(f_{\beta_1+1}) \subseteq \text{ Rang}(f_{\beta_2})$
whereas $a_{i_{\beta_2}} \notin \text{ Rang}(f_{\beta_2})$,
contradiction.
\nl
2A) If $\alpha < \beta \in S$ then by part (2) applied to the sequence
$\langle M_{\alpha + \gamma}:\gamma \le \beta - \alpha\rangle$, 
the model $M_\beta$ is
$(\lambda,\text{cf}(\beta))$-brimmed over $M_\alpha$ hence $M_\beta$
is $\le_{\frak s}$-universal over $M_\alpha$ by \scite{838-5t.15}(1)(c).  
Choose an increasing continuous sequence $\langle
\alpha_\varepsilon:\varepsilon < \text{ cf}(\delta)\rangle$ with limit
$\delta$ such that $\varepsilon < \text{ cf}(\delta) \Rightarrow
\alpha_{\varepsilon +1} \in S$ and $\alpha_0 = 0$, so clearly $\langle
M_{\alpha_\varepsilon}:\varepsilon < \text{ cf}(\delta)\rangle$
exemplifies that $M_\delta$ is $(\lambda,\text{ cf}(\delta))$-brimmed
over $M_0$.
\nl
3) \ub{Clause (a)}:

Note that necessarily cf$(\alpha_\beta) = \text{ cf}(\beta)$ as $\bar
\alpha = \langle \alpha_j:j \le \beta\rangle$ is non-decreasing 
and $\bar \alpha \restriction \beta$ is not eventually constant.

Let $\langle \beta_\varepsilon:\varepsilon < \text{ cf}
(\alpha_\beta)\rangle,\langle \gamma_\varepsilon:\varepsilon < 
\text{ cf}(\alpha_\beta)\rangle$ be increasing continuous sequences of
ordinals with limit $\beta,\alpha_\beta$ respectively 
such that $\gamma_\varepsilon \le \alpha_{\beta_\varepsilon}$ 
for every $\varepsilon < \text{ cf}(\alpha_\beta)$.

We now choose a pair $(i_\varepsilon,j_\varepsilon)$ by induction on
$\varepsilon < \text{ cf}(\alpha_\beta)$ such that:
\mr
\item "{$\odot$}"  $(a) \quad j_\varepsilon < \beta$ is
increasing continuous with $\varepsilon$
\sn
\item "{${{}}$}"  $(b) \quad i_\varepsilon \le \alpha_{j_\varepsilon}$
is increasing continuous with $\varepsilon$
\sn
\item "{${{}}$}"  $(c) \quad$ if $\varepsilon = \zeta +1$ then
$M^{\bold d}_{i_\varepsilon,j_\varepsilon}$ is $\le_{\frak
s}$-universal over $M^{\bold d}_{i_\zeta,j_\zeta}$
\sn
\item "{${{}}$}"  $(d) \quad$ if $\varepsilon = \zeta +1$ then
$i_\varepsilon > \gamma_\varepsilon,j_\varepsilon > \beta_\varepsilon$.
\ermn
There is no problem to carry the definition as $\bold d$ is weakly
universal, see Definition \scite{838-5t.29}(2A) and $\bar \alpha$ not
eventually constant.  Now the sequence
$\langle i_\varepsilon:\varepsilon < \text{ cf}(\alpha_\beta)\rangle$
is increasing with limit $\alpha_\beta$ (by clause $\odot$(d)), and
$\langle j_\varepsilon:\varepsilon < \beta\rangle$ is an increasing
continuous sequence and has
limit $\beta$ (as $\langle \alpha_j:j < \beta\rangle$ is not
eventually constant), hence
\mr
\item "{$(*)$}"  $\langle M^{\bold
d}_{i_\varepsilon,j_\varepsilon}:\varepsilon < \text{
cf}(\beta)\rangle$ is $\le_{\frak s}$-increasing continuous with union
$M^{\bold d}_{\alpha_\beta,\beta}$.
\ermn
Let $(i_{\text{cf}(\beta)},j_{\text{cf}(\beta)}) :=
(\alpha_\beta,\beta)$.

So by $\odot(c) + (*)$ it follows that $\langle M^{\bold
d}_{i_\varepsilon,j_\varepsilon}:\varepsilon \le \text{
cf}(\beta)\rangle$ satisfies the assumptions of claim
\scite{838-5t.31}(1), hence its conclusion, i.e. the sequence $(\langle M^{\bold
d}_{i_\varepsilon,j_\varepsilon}:\varepsilon \le \text{
cf}(\beta)\rangle$ is cf$(\beta)$-correct.

We shall apply part (1) of the present claim \scite{838-5t.36}.  Now the
pair $(\langle M^{\bold d}_{i_\varepsilon,j_\varepsilon}:\varepsilon
\le \text{ cf}(\beta)\rangle,\langle M^{\bold
d}_{i_\varepsilon,\beta}:\varepsilon \le \text{ cf}(\beta)\rangle)$
satisfies its assumptions hence its conclusion holds and it says that
 $\langle M^{\bold d}_{i_\varepsilon,\beta}:
\varepsilon \le \text{ cf}(\beta)\rangle$ is
cf$(\beta)$-correct.  As $\langle i_\varepsilon:\varepsilon \le \text{
cf}(\beta)\rangle$ is increasing continuous with last element $i_\varepsilon =
\alpha_\beta$ and also $\langle M^{\bold d}_{i,\beta}:i \le
\alpha_\beta\rangle$ is $\le_{\frak s}$-increasing continuous also
$\langle M^{\bold d}_{i,\beta}:i \le \alpha_\beta\rangle$ is
$\alpha_\beta$-correct, by Observation \scite{838-5t.31}(2), as required.
\bn
\ub{Clause (b)}:

We shall apply part (2) of the present claim on the sequence $\langle M^{\bold
d}_{i,\beta}:i \le \alpha_\beta\rangle$ and $\Gamma = \Gamma_{\bold d}
:= \{p \in {\Cal S}^{\text{bs}}_{\frak s}(M^{\bold
d}_{\alpha_\beta,\beta}):p$ does not fork over $M^{\bold d}_{i,j}$ for
some $i < \alpha_j,j < \beta\}$.  By the definition of an
${\frak s}$-free triangle it is $\le_{\frak s}$-increasing continuous
hence clause (b) of part (2) is satisfied.  As $\bold d$ is full by
clause (a) of Definition \scite{838-5t.29}(3) the ordinal $\alpha_\beta
= \alpha_\beta(d)$ is divisible by $\lambda$, i.e. clause (a) of part
(2) holds.  Clause (c) of the assumption of
part (2) is satisfied because we have proved clause (a) here.

As for clause (d) of part (2) let $i_1 < \alpha_\beta$ and $p_1 \in 
{\Cal S}^{\text{bs}}_{\frak s}(M^{\bold d}_{i_1,\beta})$ be given; let $p_2 \in
\Gamma_{\bold d} = {\Cal S}^{\text{bs}}_{\frak s}
(M^{\bold d}_{\alpha_\beta,\beta})$ be a
non-forking extension of $p_1$.  By the definition of 
$\Gamma_{\bold d}$,  see part (2) we
can find $j_2 < \beta$ and $i_2 \le \alpha_{j_2}$ such that $p_2$ does not
fork over $M^{\bold d}_{i_2,j_2}$.
By monotonicity, without loss of generality  $i_2 \ge i_1$ and $i_2 < \alpha_{j_2}$.  
As $\bold d$ is full (see clause (b) of Definition \scite{838-5t.29}(3)) 
we can find $S \subseteq [i_2,\alpha_\beta)$ of
order type $\ge \lambda_{\frak s}$ such that for each $i \in S$ 
an ordinal $j_*(i) < \beta$ satisfying $j_*(i) > j_2$  and an 
element $c \in \bold J^{\bold d}_{i,j_*(i)}$
 such that $i < \alpha_{j_*(i)}$ and $\ortp_{\frak s}
(c,M^{\bold d}_{i,j_*(i)},M^{\bold d}_{i+1,j_*(i)})$ is a 
non-forking extension of
$p_2 \restriction M^{\bold d}_{i_2,j_2}$.  So by the definition of
${\frak u}_{\frak s}$ we have $\bold J^{\bold d}_{i,j_*(i)} = \{c\}$
and so by the definition of ``$\bold d$ is a ${\frak u}_{\frak
s}$-free triangle" we have 
$(M^{\bold d}_{i,j_*(i)},M^{\bold d}_{i+1,j_*(i)},c) \le^1_{\frak
u} (M^{\bold d}_{i,\beta},M^{\bold d}_{i+1,\beta},c)$.
Hence $(M^{\bold d}_{i,\beta},M^{\bold d}_{i+1,\beta},c)$ 
realizes a non-forking extension of $p_2 \restriction
M_{i_2,j_2}$, hence, by the uniqueness of non-forking extensions the
element $c$ realizes $p_2 \restriction M^{\bold d}_{i,\beta}$.  So
$S$ is as required in clause (d) of the assumption of part (2).

So all the assumptions of part (2) applied to the sequence
 $\langle M^{\bold d}_{\alpha,\beta}:\alpha \le \alpha_\beta\rangle$
 and the set $\Gamma_{\bold d}$ are satisfied hence its conclusion which
says that $M^{\bold d}_{\alpha_\beta,\beta}$ is
$(\lambda,\text{cf}(\alpha_\beta))$-brimmed over $M_{i,\beta}$ for
every $i < \alpha_\beta$.  So we are done proving clause (b) hence part (3).
\nl
3A) We prove each clause.
\bn
\ub{Clause (a)}:

So $\theta := \text{ cf}(\alpha) = \text{ cf}(\beta)$.  Let $\langle
\alpha_\varepsilon:\varepsilon < \theta\rangle$ be an increasing sequence
of ordinals with limit $\alpha$ and $\langle
\beta_\varepsilon:\varepsilon < \theta\rangle$ be an increasing
sequence of ordinals with limit $\beta$.  
Now for each $\varepsilon < \theta$ we can find $i \in
(\alpha_\varepsilon,\alpha)$ and $j \in (\beta_\varepsilon,\beta)$
such that $M^{\bold d}_{i,j}$ is $\le_{\frak s}$-universal over
$M^{\bold d}_{\alpha_\varepsilon,\beta_\varepsilon}$; this holds as we
are assuming $\bold d$ is weakly universal, see Definition \scite{838-5t.29}(2A).

By monotonicity without loss of generality  $i \in \{\alpha_\zeta:\zeta \in
(\varepsilon,\theta)\}$ and $j \in \{\beta_\zeta:\zeta \in
(\varepsilon,\theta)\}$.  So without loss of generality  $M^{\bold d}_{\alpha_{\varepsilon
+1},\beta_{\varepsilon +1}}$ is $\le_{\frak s}$-universal over
$M^{\bold d}_{\alpha_\varepsilon,\beta_\varepsilon}$ for $\varepsilon
< \theta$.  Let $\alpha_\theta : = \alpha,\beta_\theta = \beta$.  

Hence by Observation \scite{838-5t.31}(1) we have $\theta \in 
\text{ correct}(\langle M^{\bold d}_{\alpha_\varepsilon,\beta_\varepsilon}:
\varepsilon \le \theta\rangle)$ which means $\theta \in \text{
correct}_{\Gamma_{\bold d}}(\langle M^{\bold
d}_{\alpha_\varepsilon,\beta_\varepsilon}:\varepsilon \le
\theta\rangle)$, see \scite{838-5t.31}(4).  
So by part (1) also $\theta \in \text{ correct}_{\Gamma_{\bold d}}
(\langle M^{\bold d}_{\alpha_\varepsilon,\beta}:\varepsilon \le
\theta\rangle)$, hence by Observation \scite{838-5t.31}(2) also $\alpha
\in \text{ correct}_{\Gamma_{\bold d}}(\langle M^{\bold d}_{i,\beta}:
i \le \alpha\rangle)$, as required.
\bn
\ub{Clause (b)}:

We can apply part (2) of the present claim to the sequence $\langle
M^{\bold d}_{i,\beta}:i\le \alpha\rangle$.  This is similar to the
proof of clause (b) of part (3), alternatively letting $\alpha'_j =
\sup\{\alpha_\varepsilon:\varepsilon < \theta$ and $\beta_\varepsilon
\le j\}$ for $j < \beta$ and $\bar \alpha' = \langle \alpha'_j:j \le
\beta\rangle$  use part (3) for the
${\frak u}$-free triangle $\bold d \rest (\bar \alpha,\beta)$,
i.e. $\bar M^{\bold d} = \langle M_{i,j}:j \le \beta$ and $i \le
\alpha'_j\rangle$, etc.; this applies to clause (a), too.
\bn
\ub{Clause (c)}:

Note that ``cf$(\alpha) = \text{ cf}(\beta)$" is not assumed.

We use part (2A) of the present claim, but we elaborate. Let $S :=\{\alpha' <
\alpha:\alpha'$ is divisible by $\lambda$ and has cofinality
cf$(\beta)\}$.

Now $S$ is a subset of $\alpha$, unbounded (as for every $i < \alpha$
we have $i + \lambda(\text{cf}(\beta)) \in S \cup \{\alpha\}$ and $i +
\lambda(\text{cf}(\beta)) \le i + \lambda^2 < i + \lambda^2 \omega \le
\alpha)$, hence it is enough to show that $i_1 < i_2 \in S \Rightarrow
M^{\bold d}_{i_2,\beta}$ is brimmed over $M^{\bold d}_{i_1,\beta}$.

Now this follows by clause (b) of part (3A) which we have
just proved applied to $\bold d' =
\bold d \restriction (i_2,\beta)$, it is a ${\frak u}$-free
$(i_2,\beta)$-rectangle, it is strongly full hence full and cf$(i_2) = 
\text{ cf}(\beta) \ge \aleph_0$.  So the assumptions of part (2A)
holds hence its conclusion so we are done.
\nl
4) Let $\kappa_1,\kappa_2$ be regular $\le \lambda$ and choose
$\alpha_\ell = \lambda^2 \times \kappa_\ell$, for $\ell=1,2$.  Let $M
\in K_{\frak s}$ and define a ${\frak u}$-free
$(\alpha_1,0)$-rectangle by $M^{\bold d^{\text{ver}}_0}_{(i,0)} = M$
for $i \le \alpha_1$ and 
$\bold J^{\bold d^{\text{ver}}_0}_{(i,0)} = \emptyset$ for $i < \alpha_1$.

Define a ${\frak u}$-free $(0,\alpha_2)$-rectangle $\bold
d_{\text{hor}}$ by $M^{\bold d_{\text{hor}}}_{(0,j)} = M$ for $j \le
\alpha_2$ and $\bold I^{\bold d}_{(0,j)} = \emptyset$ for $j <
\alpha_2$.
By Observation \scite{838-5t.33} there is 
a strongly full strictly brimmed ${\frak u}$-free 
$(\alpha_1,\alpha_2)$-rectangle $\bold d$ such that its dual is
strongly full too (and automatically strictly brimmed recalling
\scite{838-5t.39}(6)).

We can apply clause (c) of part (3A) with $(\bold d,\alpha_1,\alpha_2)$ here
standing for $(\bold d,\alpha,\beta)$ there; so we can conclude in
particular that $M_{\alpha_1,\alpha_2}$ is
$(\lambda,\text{cf}(\alpha_1))$-brimmed over $M^{\bold d}_{0,\alpha_2}$ hence
over $M_{0,0} = M$.  But ${\frak u}_{\frak s}$ is self-{\rm dual} 
so {\rm dual}$(\bold d)$
is a ${\frak u}_{\frak s}$-free $(\alpha_2,\alpha_1)$-rectangle, and
by the choice of $\bold d$ (recalling \scite{838-5t.33}(d)) it is still
strongly full and, e.g. by \scite{838-5t.39}(7) is 
universal.  So applying clause (c) of part (3A) we
get that $M^{\text{\rm dual}(\bold d)}_{\alpha_2,\alpha_1}$ is
$(\lambda,\text{cf}(\alpha_2))$-brimmed over $M^{\text{\rm dual}(\bold
d)}_{0,\alpha_1}$ hence over $M^{\text{\rm dual}(\bold d)}_{0,0} =
M^{\bold d}_{0,0} = M$.  However $M^{\text{\rm dual}(\bold
d)}_{\alpha_2,\alpha_1} = M^{\bold d}_{\alpha_1,\alpha_2}$ so this
model $\le_{\frak s}$-extend $M$ and is
$(\lambda,\text{cf}(\alpha_\ell))$-brimmed over it for $\ell=1,2$; 
this means $(\lambda,\kappa_\ell)$-brimmed over $M$. So as for each
regular $\kappa \le \lambda$, the
``$(\lambda,\kappa)$-brimmed model over $M$ for some regular $\kappa
\le \lambda$" is unique up to
isomorphism over $M$ we conclude that the brimmed model over $M$ is
unique, so we are done.
\nl
${{}}$  \hfill$\square_{\scite{838-5t.36}}$
\enddemo
\bigskip

\proclaim{\stag{838-5t.37} Claim}   If $M_1 \le_{\frak s} M_2$ are
brimmed and $p_i \in {\Cal S}^{\text{bs}}_{\frak s}(M_2)$ does not
fork over $M_1$ for $i <i_* < \lambda_{\frak s}$ \ub{then} for some
isomorphism $\pi$ from $M_2$ onto $M_1$ we have $i < i_* \Rightarrow
\pi(p_i) = p_i \restriction M_1$.
\endproclaim
\bigskip

\demo{Proof}   Easy, by \scite{838-5t.36}(4) and \scite{838-5t.15}(3), i.e. as in 
\yCITE[stg.9]{705}.  \hfill$\square_{\scite{838-5t.37}}$
\enddemo
\bn
\centerline {$* \qquad * \qquad *$}
\bigskip

Another way to deal with disjointness is through reduced triples
(earlier we have used \scite{838-1a.19}, \scite{838-1a.20}, \scite{838-e.5D} (with
some repetitions).
\definition{\stag{838-5t.41} Definition}  1) $K^{3,\text{rd}}_{\frak s}$
is the class of triples $(M,N,a) \in K^{3,\text{bs}}_{\frak s}$ which
are reduced\footnote{This is different from our choice in Definition
\yCITE[1a.34]{E46}(2), \ub{but} here this is for a given almost good
$\lambda$-frame, there it is for a $\lambda$-a.e.c. ${\frak K}$.}
 which means: if $(M,N,a) \le_{\text{bs}} (M_1,N_1,a) \in
K^{3,\text{bs}}_{\frak s}$ then $N \cap M_1 = M$.
\nl
2) We say that ${\frak s}$ has existence for $K^{3,\text{rd}}_{\frak
   s}$ \ub{when} for every $M \in K_{\frak s}$ and $p \in {\Cal
   S}^{\text{bs}}_{\frak s}(M)$ there is a pair $(N,a)$ such that the
   triple $(M,N,a) \in K^{3,\text{rd}}_{\frak s}$ realizes $p$,
   i.e. $p = \ortp_{\frak s}(a,M,N)$. 
\nl
3) Let $\xi^{\text{rd}}_{\frak s}$ be the minimal
$\xi$ from Claim \scite{838-5t.43}(4) below for $M \in 
K^{3,\text{bs}}_{\frak s}$ which is superlimit.
\nl
3A) For $M \in K_{\frak s}$, let $\xi^{\text{rd}}_{{\frak s},M} =
\xi^{\text{rd}}_M$ be the minimal $\xi < \lambda^+_{\frak s}$ in
\scite{838-5t.43}(4) below for $M$ \ub{when} it exists, $\infty$ otherwise
(well defined, i.e. $< \infty$ if ${\frak s}$ has existence for
$K^{3,\text{rd}}_{\frak s}$). 
\enddefinition
\bigskip

\proclaim{\stag{838-5t.43} Claim}  1) For every $(M,N,a) \in
K^{3,\text{bs}}_{\frak s}$ there is $(M_1,N_1,a) \in
K^{3,\text{rd}}_{\frak s}$ such that $(M,N,a) \le_{\text{bs}}
(M_1,N_1,a)$ and morever $M_1,N_1$ are brimmed over $M,N$ respectively.
\nl
2) $K^{3,\text{rd}}_{\frak s}$ is closed under increasing unions of
length $< \lambda^+$, i.e. if $\delta < \lambda^+$ is a limit ordinal
and $(M_\alpha,N_\alpha,a) \in K^{3,\text{rd}}_{\frak s}$ is 
$\le^{\text{bs}}_{\frak s}$-increasing with $\alpha(< \delta)$ and $M_\delta :=
\cup\{M_\alpha:\alpha < \delta\}$ and $N_\delta := \cup\{N_\alpha:\alpha <
\delta\}$ \ub{then} $(M_\delta,N_\delta,a) \in K^{3,\text{rd}}_{\frak s}$
and $\alpha < \delta \Rightarrow (M_\alpha,N_\alpha,a) 
\le^{\text{bs}}_{\frak s} (M_\delta,N_\delta,a)$.
\nl
3) If ${\frak K}_{\frak s}$ is categorical (in $\lambda$) \ub{then}
${\frak s}$ has existence for $K^{3,\text{rd}}_{\frak s}$.
\nl
4) For every $M \in K_{\frak s}$ and $p \in {\Cal
S}^{\text{bs}}_{\frak s}(M)$ there are $\xi < \lambda^+$, a
$\le_{\frak s}$-increasing continuous $\bar M = \langle
M_\alpha:\alpha \le \xi\rangle$ and $\bar a = \langle a_\alpha:\alpha
< \xi\rangle$ such that $M_0 = M,M_\xi$ is brimmed over $M_0$ and
each $(M_\alpha,M_{\alpha +1},a_\alpha)$ is a reduced member of 
$K^{3,\text{bs}}_{\frak s}$ and $a_0$ realizes $p$ in $M_1$, 
\ub{provided that}\footnote{Why not $\xi \le \lambda$?  The
bookkeeping is O.K. but then we have to use Ax(E)(c) in the end, but
see Exercise \scite{838-5t.52}.} 
${\frak s}$ is categorical or $M$ is brimmed (equivalently superlimit)
or ${\frak s}$ has existence for $K^{3,\text{rd}}_{\frak s}$.
\nl
5) If $(M_1,N_1,a) \in K^{3,\text{rd}}_{\frak s}$ and $(M_1,N_1,a)
   \le_{\text{bs}} (M_2,N_2,a)$ then $M_2 \cap N_1 = M_1$.
\endproclaim
\bigskip

\demo{Proof}  1),2),3)  Easy, or see details in \cite{Sh:F841}.
\nl
4) As in the proof of \scite{838-7v.19}(2B) using Fodor lemma; for the
``$M$ brimmed" case, use the moreover from part (1).
\nl
5) By the definitions.  \hfill$\square_{\scite{838-5t.43}}$ 
\enddemo
\bigskip

\demo{\stag{838-5t.45} Conclusion}  (Disjoint amalgamation) 
Assume ${\frak K}_{\frak s}$ is categorical or just has existence for
$K^{3,\text{rd}}_{\frak s}$ recalling \scite{838-5t.43}(3).  
If $(M,N_\ell,a_\ell) \in
K^{3,\text{bs}}_{\frak s}$ for $\ell=1,2$ and $N_1 \cap N_2 = M$
\ub{then} there is $N_3 \in K_{\frak s}$ such that $(M,N_\ell,a_3)
\le_{\text{bs}} (N_{3 - \ell},N_3,a_\ell)$ for $\ell=1,2$.

Hence ${\frak s}$ has disjointness, see definition \scite{838-5t.9}.
\enddemo
\bigskip

\demo{Proof}  Straight by \scite{838-5t.43}(4) similarly to Observation
\scite{838-5t.33} using \scite{838-5t.43}(5), of course.  
\hfill$\square_{\scite{838-5t.45}}$
\enddemo
\bn
\margintag{5t.47}\ub{\stag{838-5t.47} Question}:   Is \scite{838-5t.45} true without 
categoricity (and without assuming existence for $K^{3,\text{rd}}_{\frak s}$)?
\bigskip

\remark{\stag{838-5t.49} Remark}  So we can redefine ${\frak u}$ such that
the amalgamation is disjoint by restricting ourselves to ${\frak
s}_{[M]},M \in K_{\frak s}$ superlimit or assuming ${\frak s}$ has
existence for $K^{3,\text{rd}}_{\frak s}$.  

We can work with ${\frak u}$ which includes
disjointness so this enters the definition of 
$K^{3,\text{up}}_{\frak s}$ defined in \scite{838-6u.7}, 
so this is a somewhat different property
or, as we prefer, we ignore this using $=_\tau$ as in \S1.
\endremark
\bn
\margintag{5t.52}\ub{\stag{838-5t.52} Exercise}:   If $M$ is superlimit and 
$(M,N,a) \in K^{3,\text{bs}}_{\frak s}$ \ub{then} for some $\le_{\frak
s}$-increasing continuous sequence $\bar M = \langle M_\alpha:\alpha <
\lambda\rangle$ and $\bar a = \langle a_\alpha:\alpha < \lambda\rangle,
\bold d$ we have $M =  M_0 \le_{\frak s} N \le_{\frak s} M_\lambda$ 
and each $(M_\alpha,M_{\alpha +1},a_\alpha) \in K^{3,\text{bs}}_{\frak
s}$ is reduced and $a = a_0$ and $M_\lambda$ is 
brimmed over $M_0$.
\bn
[\ub{Hint}:  Let $\bar{\Cal U} = \langle {\Cal U}_\alpha:\alpha \le
\lambda\rangle$ be an increasing continuous sequence of subsets of
$\lambda$ such that $|{\Cal U}_0| = \lambda = |{\Cal U}_{\alpha +1}
\backslash {\Cal U}_0|$ and min$({\Cal U}_\alpha) \ge \alpha$ for $\alpha <
\lambda$ and ${\Cal U}_\lambda = \lambda$.

Let $\alpha_i = i$ for $i \le \lambda$ and $\bar\alpha^i= \langle
\alpha_j:j \le i\rangle$.  We choose pairs $(\bold d_\beta,\bar
p^\beta)$ by induction on $\beta \le \lambda$ such that
\mr
\item "{$\circledast$}"  $(a) \quad \bold d_\beta$ is a ${\frak
u}$-free $(\bar \alpha^\beta,\beta)$-triangle
\sn
\item "{${{}}$}"  $(b) \quad M^{\bold d_\beta}_{i+1,j+1}$ is
$\le_{\frak s}$-universal over $M^{\bold d_\beta}_{i+1,j} \cup
M^{\bold d_\beta}_{i,j+1}$ when $j < \beta$ and $i < \alpha_j$
\nl

\hskip25pt and $M^{\bold d_\beta}_{0,j+1}$ is $\le_{\frak
s}$-universal over $M^{\bold d_\beta}_{0,j}$ when $j < \beta$
\sn
\item "{${{}}$}"  $(c) \quad \bar p^\beta = \langle p_i:i \in {\Cal
U}_\alpha\rangle$ list $\cup\{{\Cal S}^{\text{bs}}_{\frak s}(M^{\bold
d_\beta}_{i,j}):j \le \beta,i \le \alpha_j\}$ each appearing 
\nl

\hskip25pt $\lambda$ times
\sn
\item "{${{}}$}"  $(d) \quad$ if $\beta = 2\alpha +1,\alpha \in {\Cal
U}_\varepsilon$ \ub{then} $\bold J^{\bold d_\beta}_{2 \alpha,\beta}
\ne \emptyset$ and letting $a$ be the unique 
\nl

\hskip25pt  member of $\bold J^{\bold d_\beta}_{2 \alpha,\beta}$, the type
$\ortp_{\frak s}(a,M^{\bold d_\beta}_{2 \alpha,\beta},
M^{\bold d_\beta}_{p,\beta})$ is a
non-forking 
\nl

\hskip25pt  extension of $p_\alpha$
\sn
\item "{${{}}$}"  $(e) \quad$ if $\beta = 2 \alpha +2$ and $\alpha \in
{\Cal U}_\varepsilon$ then there is $(M'_\beta,N'_\beta,\bold J_\beta)
\in K^{3,\text{rd}}_{\frak s}$ 
\nl

\hskip25pt  such that  $(M^{\bold d_\beta}_{\varepsilon,2 \alpha +1},
M^{\bold d_\beta}_{\varepsilon +1,2 \alpha +1},
\bold J^{\bold d_\beta}_{\varepsilon,2 \alpha +1})
\le^{\text{bs}}_{\frak s} (M'_\beta,N'_\beta,\bold J_\beta)
\le^{\text{bs}}_{\frak s}$
\nl

\hskip25pt $(M^{\bold d_\beta}_{\varepsilon,\beta},
M^{\bold d_\beta}_{\varepsilon +1,\beta},
\bold J^{\bold d_\beta}_{\varepsilon,\beta})$]
\sn
\item "{${{}}$}"  $(f) \quad M^{\bold d_0}_{0,0} = M$ and $(M^{\bold
d_1}_{0,1},M^{\bold d_1}_{1,1},\bold J^{\bold d_1}_{0,1}) =
(M,N,\{a\})$.
\ermn
Now $\langle M^{\bold d_\lambda}_{i,\lambda}:i \le \lambda\rangle$ is
as required except ``$M = M_{0,\lambda} \le_{\frak s} N \le_{\frak s}
M_{\lambda,\lambda}$".  But $(M,N,\{a\}) = (M^{\bold d_1}_{0,1},
M^{\bold d_1}_{1,1},\{a\}) \le^1_{\frak u} (M^{\bold
d_\lambda}_{0,\lambda},M^{\bold d_\lambda}_{1,\lambda},\bold J^{\bold
d_\lambda}_{0,\lambda})$, that is $(M,N,a) \le_{\text{bs}} (M^{\bold
d_\lambda}_{0,\lambda},a)$ and both $M$ and $M^{\bold
d_\lambda}_{0,\lambda}$ are brimmed equivalently superlimit, hence by
\scite{838-5t.15}(3) there is an isomorphism $\pi$ from $M$ onto $M$ onto
$M^{\bold d_\lambda}_{0,\lambda}$ mapping $\ortp(a,M,N)$ to
$\ortp(a,M^{\bold d_\lambda}_{0,\lambda},M^{\bold
d_\lambda}_{1,\lambda})$.  Recalling $M^{\bold
d_\lambda}_{\lambda,\lambda}$is brimmed over $M^{\bold
d_\lambda}_{0,\lambda}$ we can extend $\pi$ to a $\le_{\frak
s}$-embedding $\pi^+$ of $N$ into $M^{\bold
d_\lambda}_{\lambda,\lambda}$ mapping $a$ to itself, so renaming we
are done.]
\goodbreak

\head {\S6 Density of weak version of uniqueness} \endhead  \resetall \sectno=6
 \spuriousreset
\bigskip

We would like to return to ``density of 
$K^{3,\text{uq}}_{\frak s}$", where ${\frak
s}$ is a good $\lambda$-frame or just almost good $\lambda$-frames,
i.e. to eliminate the (weak) extra set theoretic assumption in the
non-structure results from failure of density of
$K^{3,\text{uq}}_{\frak s}$.
But we start with a  notion $K^{3,\text{up}}_{\frak s}$, weaker then
$K^{3,\text{uq}}_{\frak s}$ related
to weak non-forking relations defined later in Definition
\scite{838-7v.35}.  Defining weak non-forking we shall waive 
uniqueness, but still can
``lift ${\frak u}_{\frak s}$-free $(\alpha,0)$-rectangles".  We now
look for a dichotomy - either a non-structure results applying the
theorems of \S2 or actually \S3, \ub{or} density of
$K^{3,\text{up}}_{\frak s}$.  Using the last possibility in the subsequent
sections \S7, \S8 we get a similar dichotomy with $K^{3,\text{uq}}_{\frak s}$.

It turns out here that what we prove is somewhat weaker than density
for $K^{3,\text{up}}_{\frak s}$ in some ways.
Mainly we prove it for the $K^{3,\text{up}}_{{\frak s},\xi}$ version
for each $\xi \le \lambda^+$.  Actually for every $M \in K_{\frak s}$ and $p \in
{\Cal S}^{\text{bs}}_{\frak s}(M)$ what we find is a triple $(M_1,N_1,a) \in
K^{3,\text{up}}_{\frak s}$ such that 
$M \le_{\frak s} M_2$ and the type $\ortp_{\frak s}(a,M_1,N_1)$ is a
non-forking extension of $p$; not a serious difference when 
${\frak s}$ is categorical which is reasonable for our purposes.
Eventually in this section we have to use ${\frak K}_{\frak s}$ with
fake equality (to apply \scite{838-3r.74}), this is justified in \scite{838-6u.31}.
\bn
\ub{Discussion}:  Why do we deal with $K^{3,\text{up}}_{{\frak
s},\xi}$ for $\xi \le \lambda^+$ rather than with
$K^{3,\text{up}}_{\frak s}$?  The point is that, e.g. in the
weak/semi/vertical uq-invariant coding property in \S3
(see Definitions \scite{838-3r.1}, \scite{838-3r.71}, \scite{838-3r.19}), 
given $(\bar M,\bar{\bold J},\bold
f) \in K^{\text{qt}}_{\frak u}$ and $(M_{\alpha(0)},N_0,\bold I) \in
\text{ FR}^1_{\frak u}$, for a club of $\delta < \partial$ we promise
the existence of a ${\frak u}$-free $(\alpha,0)$-rectangle $\bold
d_\delta$, which is O.K. for \ub{every} $N_\delta$ such that
$(M_{\alpha(*)},N_0,\bold I) \le^1_{\frak u} (M_\delta,N_\delta,\bold
I)$.  So the failure gives (not much more than) that for every 
${\frak u}$-free $(\alpha,0)$-rectangle $\bold d$ and $p \in {\Cal
S}^{\text{bs}}_{\frak s}(M^{\bold d}_{\alpha,0})$ there is a pair
$(N,\bold I)$ such that $(M,N,\bold I)$ realizes $p$ and $\bold d$ is
what we call uq-orthogonal to $(M,N,\bold I)$.  
We like to invert the quantifiers,
i.e. ``for $p \in {\Cal S}^{\text{bs}}_{\frak s}(M)$ there is $(N,\bold I)$
such that for every $\bold d$....".  Of course, we assume categoricity
(of $K_{\frak s}$, that is in $\lambda$), but we need to use a ``universal
$\bold d$".  This is guaranteed by \scite{838-6u.19} for
$K^{3,\text{up}}_{{\frak s},\xi}$ for any $\xi \le \lambda^+$ (but we
have to work more for $\xi = \lambda^+$, i.e. for all $\xi < \lambda^+$ at
once, i.e. for $K^{3,\text{up}}_{\frak s}$). 
\bigskip

\demo{\stag{838-6u.1} Hypothesis}  1) As in \scite{838-5t.1} and for transparency
${\frak s}$ has disjointness.
\nl
2) ${\frak u} = {\frak u}_{\frak s}$, see Definition \scite{838-5t.19},
Claim \scite{838-5t.21}, so $\partial = \lambda^+$.
\enddemo
\bigskip

\proclaim{\stag{838-6u.3} Claim}  1) For almost$_2$ all 
$(M,\bar{\bold J},\bold f) \in K^{\text{qt}}_{\frak u}$ 
the model $M_\partial$ belongs to $K^{\frak s}_{\lambda^+}$ and is
saturated (above $\lambda$).
\nl
2) If ${\frak s}$ has the fake equality $=_\tau$ (e.g. ${\frak s} =
   {\frak t}'$ where ${\frak t}$ is an almost good $\lambda$-frame and
   ${\frak t}'$ is defined as in \scite{838-e.5D}(1), 
${\frak u} = {\frak u}^1_{\frak s}$, see \scite{838-5t.19},\scite{838-5t.21},
   then for some ${\frak u}-0$-appropriate 
${\frak h}$, if $\langle (\bar M^\alpha,
\bar{\bold J}^\alpha,\bold f^\alpha):\alpha <
\partial^+_{\frak u}\rangle$ is $\le^{\text{qt}}_u$-increasing
 continuous and obeys ${\frak h}$, \ub{then} $M =
   \cup\{M^\alpha_j:\alpha < \partial^+_{\frak u}\rangle$ is $=_\tau$-fuller.
\endproclaim
\bigskip

\remark{Remark}  1) See Definition \scite{838-1a.43}.
\nl
2) Compare with \scite{838-e.5L}.
\endremark
\bigskip

\demo{Proof}  1) We choose ${\frak g}$ as in 
Definition \scite{838-1a.43}(2) such that:
\mr
\item "{$(*)_1$}"  if $S \subseteq \partial$ is a stationary subset of
$\partial$ and the pair $((\bar M^1,\bar{\bold J}^1,\bold f^1),(\bar
M^2,\bar{\bold J}^2,\bold f^2))$ strictly $S$-obeys ${\frak g}$ \ub{then}:
{\roster
\itemitem{ $\odot$ }    for some club $E$ of $\partial$ for every $\delta
\in S \cap E$, we have
\sn
\itemitem{ ${{}}$ }  $(a) \quad$ if $i < \bold f^1(\delta)$ then $M^2_{\delta
+i+1}$ is $\le_{\frak s}$-universal over $M^1_{\delta +i+1} \cup
M^2_{\delta +i}$
\sn
\itemitem{ ${{}}$ }  $(b) \quad \bold f^2(\delta)$ is $> \bold f^1(\delta)$ and
is divisible by $\lambda^2$
\sn
\itemitem{ ${{}}$ }  $(c) \quad$ if 
$i \in [\bold f^1(\delta),\bold f^2(\delta))$ then
\sn
\itemitem{ ${{}}$ }  $\quad (\alpha) \quad M^2_{\delta +i+1}$ 
is $\le_{\frak s}$-universal over $M^2_{\delta +i}$
\sn
\itemitem{ ${{}}$ }  $\quad (\beta) \quad$ if $p \in 
{\Cal S}^{\text{bs}}_{\frak s}(M^2_{\delta +i})$ \ub{then} for
$\lambda$ ordinals $i_1 \in [i,i + \lambda)$ the type
\nl

\hskip45pt $\ortp_{\frak s}(a_{\bold J^2_{\delta +i}},M^2_{\delta +i},
M^2_{\delta +i+1})$ is a non-forking extension of $p$
\nl

\hskip45pt  where $b_{\bold J^2_{\delta +i}}$ is the unique member of 
$\bold J^2_{\delta +i}$.
\endroster}
\ermn
We can find such ${\frak g}$.  Now
\mr
\item "{$(*)_2$}"  assume $\langle(\bar M^\alpha,\bar{\bold
J}^\alpha,\bold f^\alpha):\alpha \le \delta\rangle$ is
$\le^{\text{qt}}_{\frak u}$-increasing continuous (see Definition
\scite{838-1a.29}(4A)) and obey
${\frak g}$ (i.e. for some\footnote{we can use a decreasing sequence
 of $S$'s but then we really use the last one only,
the point is that ${\frak g}$ treat each $\delta \in S$ in
 the same way (rather than dividing it according to tasks, a
 reasonable approach, but not needed here}
 stationary $S \subseteq \partial$ for unboundedly many $\alpha < \delta$).
\ub{Then} $M^\delta_\partial$ is saturated above $\lambda$.
\ermn
[Why?  Let $\kappa = \text{ cf}(\delta)$, of course $\aleph_0 \le
\kappa \le \lambda^+$.

We can find an increasing
continuous sequence $\langle \alpha_\varepsilon:\varepsilon < 
\text{ cf}(\delta) = \kappa\rangle$ of ordinals with limit $\delta$ such
that: 
\mr
\item "{$(*)_3$}"  if $\varepsilon = \zeta +1 < \kappa$ and
$\varepsilon$ is an even ordinal then
$\alpha_{\varepsilon +1} = \alpha_\varepsilon +1$ and letting
$\alpha_{\text{cf}(\delta)} = \delta$ the pair
$\langle(\bar M^{\alpha_\varepsilon},\bar{\bold J}^{\alpha_\varepsilon},\bold
f^{\alpha_\varepsilon}),(\bar M^{\alpha_{\varepsilon +1}},\bar{\bold
J}^{\alpha_{\varepsilon +1}},\bold f^{\alpha_{\varepsilon
+1}})\rangle$ does $S$-obey ${\frak g}$.
\ermn
Now clearly
\mr
\item "{$(*)_4$}"  for every $\varepsilon < \lambda^+$ satisfying 
$\varepsilon \le
\kappa$ for some club $E = E_\varepsilon$ of $\lambda^+$, for every
$\delta \in S \cap E$ we have:
{\roster
\itemitem{ $(a)$ }  $\langle \bold f^{\alpha_\zeta}(\delta):\zeta \le
\varepsilon\rangle$ is non-decreasing and is continuous
\sn
\itemitem{ $(b)$ }  $\bold d_\delta$ is a ${\frak u}$-free $(\bar
\alpha^\delta,\beta)$-triangle where $\bar \alpha^\delta = \langle
\bold f^{\alpha_\zeta}(\delta):\zeta \le \varepsilon\rangle,\beta
= \varepsilon$ such that $M^{\bold d_\delta}_{i,\zeta} =
M^{\alpha_\zeta}_{\delta +i}$ for $\zeta \le \varepsilon,i \le 
\bold f^{\alpha_\zeta}(\delta)$ and $\bold J^{\bold
d_\delta}_{i,\zeta} = \bold J^{\alpha_\zeta}_{\delta +i}$ for $\zeta \le
\varepsilon,i < \bold f^{\alpha_\zeta}(\delta)$ and $\bold I^{\bold
d_\delta}_{i,\zeta} = \emptyset$ for $\zeta < \varepsilon,i \le
\bold f^{\alpha_\zeta}(\delta)$.  
\endroster}
\ermn
Sorting out the definition by \scite{838-5t.27}(3A)(a) and the
correctness claim \scite{838-5t.36}(3)(b), clearly:
\mr
\item "{$(*)$}"  for a club of $\gamma < \partial$, if $\gamma \in S$
\ub{then} $M^{\alpha_\kappa}_{\gamma + \bold f^{\alpha_\kappa}(\gamma)}$ is
brimmed over $M^{\alpha_\kappa}_\gamma$ (and even
$M^{\alpha_\kappa}_{\gamma +i}$ for $i < \bold f^{\alpha_\kappa}(\gamma))$
\ermn
which means that $M^\delta_{\gamma + \bold f^\delta(\gamma)}$ is brimmed over
$M^\delta_\gamma$.

This is clearly enough.
\nl
2) Should be clear.    \hfill$\square_{\scite{838-6u.3}}$
\enddemo
\bigskip

\demo{\stag{838-6u.5} Conclusion}  For any $(\bar M,\bar{\bold J},\bold f)
\in K^{\text{qt}}_{\frak u}$ and for stationary $S \subseteq \partial$
there is ${\frak u}-2$-appropriate ${\frak g}$ with $S_{\frak g} = S$ 
(see \scite{838-1a.43} and \scite{838-1a.45}) such that:
\mr
\item "{$\circledast$}"  if $\langle(\bar M^\alpha,\bar{\bold J}^\alpha,\bold
f^\alpha):\alpha \le \partial\rangle$ is $\le^{\text{qs}}_{\frak
u}$-increasing continuous and 2-obeys ${\frak g}$ such that $(\bar
M^0,\bar{\bold J}^0,\bold f^0) = (\bar M,\bar{\bold J},\bold f)$ (so
$\bold f^\partial(\delta) = \sup\{\bold f^\alpha(\delta):\alpha < \delta\}$
for a club of $\delta < \partial$)
\ub{then} for a club of $\delta < \partial$ the model
$M^\partial_\delta = \cup\{M^\alpha_\delta:
\alpha < \delta\} \in K_{\frak s}$ is brimmed over $M_\delta = M^0_\delta$.
\endroster
\enddemo
\bigskip

\demo{Proof}  Similar to the proof of \scite{838-6u.3}; only the
${\frak u}$-free triangle is flipped, i.e. it is a {\rm dual}$({\frak
u})$-free triangle but {\rm dual}$({\frak u}) = {\frak u}$.  
\hfill$\square_{\scite{838-6u.5}}$
\enddemo
\bigskip

\definition{\stag{838-6u.7} Definition}  Let $1 \le \xi \le \lambda^+$, if we
omit it we mean $\xi = \lambda^+$. 
\nl
1) $K^{3,\text{up}}_{{\frak s},\xi}$ is 
the class of $(M,N,a) \in K^{3,\text{bs}}_{\frak s}$ which has 
up-$\xi$-uniqueness, which means:
\mr
\item "{$\circledast$}"  if $(M,N,a) \le^1_{\frak u} 
(M',N',a)$ and ${\bold d}$ is a ${\frak u}$-free 
$(\alpha,0)$-rectangle with $\alpha \le \xi,\alpha < \lambda^+$ 
satisfying $(M^{\bold d}_{0,0},M^{\bold d}_{\alpha,0}) = (M,M')$ \ub{then}
$\bold d$ can be lifted for $((M,N,a),N')$ which means:
{\roster
\itemitem{ $\boxdot$ }    we can find a ${\frak u}$-free 
$(\alpha +1,1)$-rectangle $\bold d^*$ and $f$ such that
$\bold d^* \restriction (\alpha,0) = {\bold d},f(N) \le_{\frak s} 
M^{\bold d^*}_{0,1},f(a) = a^{\bold d^*}_{0,0}$, i.e. $\bold I^{\bold
d^*}_{0,0} = \{f(a)\}$  and $f$ is a 
$\le_{\frak s}$-embedding of $N'$ into $M^{\bold d^*}_{\alpha +1,1}$
over $M'$ hence also $f \rest M = \text{ id}_M$.
\endroster}
\ermn
2) We say that $K^{3,\text{up}}_{{\frak s},\xi}$ is dense (in
$K^{3,\text{bs}}_{\frak s}$) or ${\frak s}$ has density for
$K^{3,\text{up}}_{{\frak s},\xi}$ \ub{when} for every $(M_0,N_0,a) \in
K^{3,\text{bs}}_{\frak s}$ there is $(M_1,N_1,a) \in
K^{3,\text{up}}_{{\frak s},\xi}$ such that $(M_0,N_0,a) \le_{\text{bs}}
(M_1,N_1,a)$.
\nl
3) We say that ${\frak s}$ has (or satisfies)
existence for $K^{3,\text{up}}_{{\frak s},\xi}$ or
$K^{3,\text{up}}_{{\frak s},\xi}$ has (or satisfies) existence \ub{when}:  
if $M \in K_{\frak s}$ and $p \in {\Cal S}^{\text{bs}}_{\frak s}(M)$ 
\ub{then} for some pair $(N,a)$ we have
$(M,N,a) \in K^{3,\text{up}}_{{\frak s},\xi}$ and 
$p = \ortp_{\frak s}(a,M,N)$.
\nl
3A) We say that ${\frak s}$ has existence for 
$K^{3,\text{up}}_{{\frak s},< \xi}$ 
\ub{when} ${\frak s}$ has existence for 
$K^{3,\text{up}}_{{\frak s},\zeta}$ for every $\zeta \in [1,\xi)$;
similarly in the cases below.
\nl
4) Let $K^{3,\text{up}}_{{\frak s},\xi} = 
\cap\{K^{3,\text{up}}_{{\frak s},\zeta}:\zeta < \xi\}$.
\nl
5) Let $K^{3,\text{up+rd}}_{{\frak s},\xi}$ be defined as
$K^{3,\text{up}}_{{\frak s},\xi} \cap K^{3,\text{rd}}_{\frak s}$
recalling Definition \scite{838-5t.41}, and we repeat parts (2),(3),(4) for it.
\enddefinition
\bigskip

\demo{\stag{838-6u.9} Observation}  1) $K^{3,\text{up}}_{{\frak s},\xi} \subseteq
K^{3,\text{up}}_{{\frak s},\zeta} \subseteq K^{3,\text{rd}}_{\frak s}$
recalling Definition \scite{838-5t.41} when $1 \le \zeta \le \xi \le
\lambda^+$.
\nl
2) $\xi = \lambda^+ \Rightarrow K^{3,\text{up}}_{{\frak
s},\xi} = \cap\{K^{3,\text{up}}_{{\frak s},\zeta}:\zeta \in
[1,\xi)\}$.
\nl
3) The triple $(M,N,a) \in K^{3,\text{bs}}_{\frak s}$ does not belong
to $K^{3,\text{up}}_{{\frak s},\xi}$ \ub{iff} we can find $\bold d_1,\bold
d_2$ such that for $\ell=1,2$
\mr
\item "{$\square$}"  $(a) \quad \bold d_\ell$ is a ${\frak
u}$-free $(\alpha_\ell,1)$-rectangle and $\alpha_\ell < \text{ min}\{\xi
+1,\lambda^+_{\frak s}\}$
\sn
\item "{${{}}$}"  $(b) \quad (M^{\bold d_\ell}_{0,0},M^{\bold d_\ell}_{0,1},
\bold I^{\bold d_\ell}_{0,0}) = (M,N,\{a\})$
\sn
\item "{${{}}$}"  $(c) \quad M^{\bold d_1}_{\alpha_1,0} = M^{\bold
d_2}_{\alpha_2,0}$ and $M^{\bold d_1}_{0,0} = M = M^{\bold d_2}_{0,0}$ 
\sn
\item "{${{}}$}"  $(d) \quad$ there is no triple $(\bold d,f)$ such
that
{\roster
\itemitem{ ${{}}$ }  $(\alpha) \quad \bold d$ is a ${\frak u}$-free
$(\alpha_{\bold d},\beta_{\bold d})$-rectangle
\sn
\itemitem{ ${{}}$ }  $(\beta) \quad 
\bold d \restriction (\alpha_1,1) = \bold d_1$ 
so $\beta_{\bold d} \ge 1,\alpha_{\bold d} \ge \alpha_1$
\sn
\itemitem{ ${{}}$ }   $(\gamma) \quad f$ is a $\le_{\frak s}$-embedding of
$M^{\bold d_2}_{\alpha_2,1}$ into 
$M^{\bold d}_{\alpha_{\bold d},\beta_{\bold d}}$ over 
\nl

\hskip30pt $M^{\bold d}_{\alpha_1,0} =
M^{\bold d_1}_{\alpha(\bold d_1),0}$ mapping $M^{\bold d_2}_{0,1}$
into $M^{\bold d}_{0,\beta(\bold d)}$ and $a^{\bold d_2}_{0,0}$
to itself.
\endroster}
\endroster
\enddemo
\bigskip

\demo{Proof}  0) By the definition; as for a ${\frak u}$-free
$(\alpha,\beta)$-rectangle or $(\bar \alpha,\beta)$-triangle $\bold d$
we have: if $i = \text{ min}\{i_1,i_2\},j = \text{ min}\{j_1,j_2\}$
and $M^{\bold d}_{i_1,j_1},M^{\bold d}_{i_2,j_2}$ are well defined
then $M^{\bold d}_{i,j} = M^{\bold d}_{i_1,j_1} \cap M^{\bold
d}_{i_2,j_2}$.
\nl
1) By the definition (no need of categoricity).
\nl
2) By \scite{838-5t.15}(3).
\nl
3) Straight, recalling that ${\frak u}$ satisfies monotonicity,
(E)(e), see \scite{838-1a.24} but we elaborate.
\enddemo
\bn
\ub{The Direction $\Rightarrow$}:

So $(M,N,a)$ belongs to $K^{3,\text{bs}}_{\frak s}$ but not to
$K^{3,\text{up}}_{{\frak s},\xi}$.  So by Definition \scite{838-6u.7}(1)
there is $(M',N',\bold d)$ exemplifying the failure of $\circledast$
from \scite{838-6u.7}(1), which means
\mr
\item "{$\odot$}"  $(a) \quad (M,N,a) \le^1_{\frak u} (M',N',a)$,
i.e. see Definition \scite{838-5t.19}, i.e. \scite{838-e.5I}, i.e. 
\nl

\hskip25pt $M = N \cap M'$ and $(M,N,a) \le_{\text{bs}} (M',N',a)$
\sn
\item "{${{}}$}"  $(b) \quad \bold d$ is a ${\frak u}$-free
$(\alpha,0)$-rectangle with $\alpha \le \xi,\alpha < \lambda^+$
\sn
\item "{${{}}$}"  $(c) \quad (M^{\bold d}_{0,0},M^{\bold
d}_{\alpha,0}) = (M,M')$
\sn
\item "{${{}}$}"  $(d) \quad \bold d$ cannot be lifted for
$(M,N,a,N')$, i.e. there is no pair $(\bold d^*,f)$
\nl

\hskip25pt  such that
{\roster
\itemitem{ $(\alpha)$ }   $\bold d^*$ is a ${\frak u}$-free $(\alpha
+1,1)$-rectangle 
\sn
\itemitem{ $(\beta)$ }  $\bold d^* \rest (\alpha,0) = \bold d$
\sn
\itemitem{ $(\gamma)$ }    $f$ is a $\le_{\frak s}$-embedding of $N'$
into $M^{\bold d^*}_{\alpha +1,1}$ over $M' = M^{\bold d}_{\alpha,0}$
\sn
\itemitem{ $(\delta)$ }   $f(N) \le_{\frak s} M^{\bold d^*}_{0,1}$ and
$f(a) = a^{\bold d^*}_{0,0}$.
\endroster}
\ermn
We define $\bold d_2$ by
\mr
\item "{$\odot_2$}"  $\bold d_2$ is the ${\frak u}$-free
$(1,1)$-rectangle if
{\roster
\itemitem{ $(a)$ }  $(M^{\bold d_2}_{0,0},M^{\bold d_2}_{1,0},\bold
J^{\bold d_2}_{0,0}) = (M,M',\emptyset)$
\sn
\itemitem{ $(b)$ }   $(M^{\bold d_2}_{0,0},M^{\bold d_2}_{1,0},\bold
I^{\bold d_2}_{0,0}) = (M,N,\{a\})$
\sn
\itemitem{ $(c)$ }    $M^{\bold d_2}_{1,1} = N,\bold J^{\bold
d_2}_{0,1} = \emptyset,\bold I^{\bold d_2}_{1,0} = \{a\}$.
\endroster}
\ermn
We choose $\bold d_1$ such that
\mr
\item "{$\odot_2$}"  $(a) \quad \bold d_2$ is a ${\frak u}$-free
$(\alpha,1)$-rectangle
\sn
\item "{${{}}$}"  $(b) \quad \bold d_2 \rest (\alpha,0) = \bold d^*$
\sn
\item "{${{}}$}"  $(c) \quad (M^{\bold d_2}_{0,0},M^{\bold
d_2}_{0,1},\bold I^{\bold d_2}_{0,0}) = (M,N,\{a\})$.
\endroster
\bn
\ub{The Direction $\Leftarrow$}:

Choose $\bold d = \bold d_1$ and use Exercise \scite{838-1a.24}.
\hfill$\square_{\scite{838-6u.9}}$ 
\bigskip

\definition{\stag{838-6u.11} Definition}  We say $(M,N,a)$ is a
up-orthogonal to $\bold d$ \ub{when}:
\mr
\item "{$\circledast$}"  $(a) \quad (M,N,a) \in K^{3,\text{bs}}_{\frak s}$
\sn
\item "{${{}}$}"  $(b) \quad \bold d$ is a ${\frak u}$-free $(\alpha_{\bold
d},0)$-rectangle
\sn
\item "{${{}}$}"  $(c) \quad M^{\bold d}_{0,0} = M$
\sn
\item "{${{}}$}"  $(d) \quad$ \ub{Case 1}:  
$N \cap M^{\bold d}_{\alpha(\bold d),0} = M$.
\nl
If $N_1$ satisfies
 $(M,N,a) \le^1_{\frak u} (M^{\bold d}_{\alpha(\bold d),0},N_1,a)$,
\nl

\hskip25pt  so $N_1$ does $\le_{\frak s}$-extends 
$M^{\bold d}_{\alpha(\bold d),0}$ and
$N$, \ub{then} the rectangle $\bold d$
\nl

\hskip25pt  can be lifted for $((M,N,a),N_1)$; 
\nl

\hskip25pt \ub{Case 2}:  possibly $N \cap M^{\bold d}_{(\alpha(\bold
d),0)} \ne M$.
\nl
We replace $(N,a)$ by $(N',a')$
\nl

\hskip25pt  such that $N' \cap
M^{\bold d}_{\alpha(\bold d),0} = M$ and there 
\nl

\hskip25pt is an isomorphism from $N$ onto $N'$ over $M$ mapping $a$ to $a'$.
\endroster
\enddefinition
\bn
We now consider a relative of Definition \scite{838-6u.7}.
\definition{\stag{838-6u.13} Definition}  Let $\xi \le \lambda^+$ but $\xi
\ge 1$, if $\xi = \lambda^+$ we may omit it.
\nl
1) We say that ${\frak s}$ has almost-existence for 
$K^{3,\text{up}}_{{\frak s},\xi}$ \ub{when}: if 
$M \in K_{\frak s}$ and $p \in 
{\Cal S}^{\text{bs}}_{\frak s}(M)$ we have
\mr
\item "{$\odot_{M,p}$}"  if $\alpha \le \xi$ and
$\bold d$ is a ${\frak u}$-free $(\alpha,0)$-rectangle with 
$M^{\bold d}_{0,0} = M$ (yes, we allow $\alpha = \xi = \lambda^+$)
\ub{then} there is
a triple $(M,N,a) \in K^{3,\text{bs}}_{\frak s}$ such that $p =
\ortp_{\frak s}(a,M,N)$ and $(M,N,a)$ is up-orthogonal to $\bold d$.
\ermn
2) We say that ${\frak s}$ has the weak density for
   $K^{3,\text{up}}_{{\frak s},\xi}$ \ub{when}: if 
$M \in K_{\frak s}$ and $p \in {\Cal S}^{\text{bs}}_{\frak s}(M)$ 
\ub{then} for some $(M_1,p_1)$ the demand
$\odot_{M_1,p_1}$ in part (1) holds and $M \le_{\frak s} M_1$ and $p_1 \in
{\Cal S}^{\text{bs}}_{\frak s}(M_1)$ is a non-forking extension of
$p$.  
\nl
3) We write ``almost-existence/weak density for
   $K^{3,\text{up}}_{{\frak s},< \xi}$" when this holds for every
   $\xi' < \xi$.
\enddefinition
\bigskip

\demo{\stag{838-6u.14} Observation}   Assume ${\frak s}$ is categorical (in
$\lambda_{\frak s}$) and $\xi \le \lambda^+$. 
\nl
1) Then ${\frak s}$ has weak
density for $K^{3,\text{up}}_{{\frak s},\xi}$ \ub{iff}
${\frak s}$ has almost existence for $K^{3,\text{up}}_{{\frak s},\xi}$.
\nl
2) If ${\frak s}$ has existence for
$K^{3,\text{up}}_{{\frak s},\xi}$ \ub{then} ${\frak s}$ has almost
existence for $K^{3,\text{up}}_{{\frak s},\xi}$.
\enddemo
\bigskip

\demo{Proof}  1) The weak density version implies the existence
version, i.e. the first implies the second because if 
$M \le_{\frak s} M_1$ and $p_1 \in {\Cal
S}^{\text{bs}}_{\frak s}(M_1)$ does not fork over $M$ then there is an
isomorphism $f$ from $M_1$ onto $M$ mapping $p$ to $p \rest M$, see
\scite{838-5t.37}.  The inverse is obvious.
\nl
2) Read the definitions.    \hfill$\square_{\scite{838-6u.14}}$
\enddemo
\bn
\margintag{6u.15}\ub{\stag{838-6u.15} Discussion}:  Below we fix $p_* \in {\Cal
S}^{\text{bs}}_{\frak s}(M_*)$ and look only at stationarization of
$p_*$.  We shall use the failure of almost-existence for
$K^{3,\text{up}}_{{\frak s},\xi}$ to get non-structure.

We first present a proof in the case 
${\Cal D}_\partial$ is not $\partial^+$-saturated (see
\scite{838-6u.24}) but by a more complicated proof this is not necessary, see
\scite{838-6u.23}.  As it happens, we  do not assume 
$2^\lambda < 2^{\lambda^+}$, but still assume $2^{\lambda^+} <
2^{\lambda^{++}}$ using the failure of weak density for
$K^{3,\text{up}}_{{\frak s},< \lambda^+}$ to get an up-invariant
coding property (avoiding the problem we encounter when we try to use
$\bold d_\delta$ depending on $N^\eta_\delta$).

So in \scite{838-6u.24} for each $\alpha < \lambda^{++} = \partial^+$ we ``give" a
stationary $S_\alpha \subseteq \partial$ almost disjoint to $S_\beta$
for $\beta < \alpha^+$.

Well, we have for the time being decided to deal only with
up-uniqueness arguing that it will help to deal with ``true" uniqueness.
Also, in the non-structure we use failure of weak density for
$K^{3,\text{up}}_{{\frak s},\xi}$, whereas for the positive side, in
\S7, we use existence.  The difference is that for existence we have
$(N,a)$ for given $(M,p,\xi)$ whereas for almost existence we are
given $(M,p,\bold d)$.   However, we now prove their equivalence. To
get the full theorem \scite{838-6u.23} we use \scite{838-3r.19} -
\scite{838-3r.91}.
\bigskip

\proclaim{\stag{838-6u.16} Claim}  Assume ${\frak s}$ is categorical; if
$\xi \le \lambda^+$ and ${\frak s}$ has almost-existence for
$K^{3,\text{up}}_{{\frak s},\xi}$ \ub{then} ${\frak s}$ has
existence for $K^{3,\text{up}}_{{\frak s},\xi}$.
\endproclaim
\bigskip

\demo{Proof}  This follows from the following two subclaims,
\scite{838-6u.17}, \scite{838-6u.19}.
\enddemo
\bigskip

\proclaim{\stag{838-6u.17} Subclaim}  If $(M,N,a) \in
K^{3,\text{bs}}_{\frak s}$ is up-orthogonal to $\bold d_2$, \ub{then}
it is up-orthogonal to $\bold d_1$ \ub{when}:
\mr
\item "{$\circledast$}"  $(a) \quad \bold d_\ell$ is a ${\frak u}$-free
$(\alpha_\ell,0)$-rectangle for $\ell=1,2$ where $\alpha_\ell \le \lambda^+$
\sn
\item "{${{}}$}"  $(b) \quad M^{\bold d_1}_{0,0} = M = M^{\bold d_2}_{0,0}$
\sn
\item "{${{}}$}"  $(c) \quad h$ is an increasing function from $\alpha_1$ to
$\alpha_2$
\sn
\item "{${{}}$}"  $(d) \quad f$ is an $\le_{\frak s}$-embedding of $M^{\bold
d_1}_{\alpha_1,0}$ into $M^{\bold d_2}_{\alpha_2,0}$
\sn
\item "{${{}}$}"  $(e) \quad f \restriction M^{\bold d_1}_{0,0} 
= \text{\rm id}_M$
\sn
\item "{${{}}$}"  $(f) \quad$ if $\beta < \alpha_1$ then
{\roster
\itemitem{ ${{}}$ }   $(\alpha) \quad f(b^{\bold d_1}_{\beta,0}) = 
b^{\bold d_2}_{h(\beta),0}$
\sn
\itemitem{ ${{}}$ }   $(\beta) \quad f$ maps $M^{\bold d_1}_{\beta,0}$ 
into $M^{\bold d_2}_{h(\beta),0}$
\sn
\itemitem{ ${{}}$ }   $(\gamma) \quad \ortp_{\frak s}(b^{\bold
d_2}_{h(\beta),0},M^{\bold d_2}_{h(\beta),0},M^{\bold
d_2}_{h(\beta)+1,0})$ does not fork over $f(M^{\bold d_1}_{\beta,0})$.
\endroster}
\endroster
\endproclaim
\bigskip

\demo{Proof}  Without loss of generality $f$ is the identity and
$M^{\bold d_2}_{\alpha_2,0} \cap N = M$.

So assume $N_1 \in K_{\frak s}$ is a $\le_{\frak s}$-extension of $N$
and of $M^{\bold d_1}_{\alpha_1,0}$ such that $(M,N,a) \le_{\text{bs}}
(M^{\bold d_1}_{\alpha_1,0},N_1,a)$ and we should prove the existence
of a suitable lifting.  Without loss of generality $N_1 \cap 
M^{\bold d_2}_{\alpha_2,0} = M^{\bold d_1}_{\alpha_1,0}$.  Hence there
is $N_2$ which does $\le_{\frak s}$-extend $M^{\bold
d_2}_{\alpha_2,0}$ and $N_1$ and $(M^{\bold d_1}_{\alpha_1,0},N_1,a)
\le^1_{\frak u} (M^{\bold d_2}_{\alpha_2,0},N_2,a)$; but
$\le^1_{\frak u}$ is a partial order hence $(M,N,a) \le^1_{\frak u}
(M^{\bold d_2}_{\alpha_{2,0}},N_2,a)$.

Recall that we are assuming $(M,N,a)$ is up-orthogonal to $\bold d_2$
hence we can find $\bold d^2,f$ as in Definition \scite{838-6u.11}, i.e.
as in $\odot$ inside Definition \scite{838-6u.7}(1), so $\bold d^2$ is a
${\frak u}$-free $(\alpha_2 +1,1)$-rectangle, $\bold d^2 \rest (\alpha_2,0)
= \bold d_2,f$ is a $\le_{\frak s}$-embedding of $N_2$ into $M^{\bold
d_2}_{\alpha_2+1,1}$ over $M^{\bold d_2}_{\alpha_2,0}$
 mapping $N$ into $M^{\bold d}_{0,1}$ satisfying
$f(a) = a^{\bold d^2}_{0,0}$, note: $\bold I^{\bold d}_{\alpha,0} = \{a\}$
for $\alpha \le \alpha_2 +1$.  Now we define
$\bold d^1$, a ${\frak u}$-free $(\alpha_1+1,1)$-rectangle by
\mr
\item "{$\boxtimes$}"  $(a) \quad \bold d^1 \restriction (\alpha_1,0) 
= \bold d_1$
\sn
\item "{${{}}$}"  $(b) \quad M^{\bold d^1}_{\alpha,1}$ is
$M^{\bold d^2}_{h(\alpha),1}$ if $\alpha \le \alpha_1$ is a non-limit
ordinal and is 
\nl

\hskip25pt  $\cup\{M^{\bold d^2}_{\beta,1}:\beta < h(\alpha)\}$ if $\alpha
\le \alpha_1$ is a limit ordinal and is $M^{\bold d^2}_{\alpha_1
+1,1}$
\nl

\hskip25pt  if $\alpha = \alpha_1 +1$
\sn
\item "{${{}}$}"  $(c) \quad M^{\bold d^1}_{\alpha_1 + 1,0} = M^{\bold
d^2}_{\alpha_1,0}$
\sn
\item "{${{}}$}"  $(d) \quad \bold I^{\bold d_1}_{\alpha,0} = 
\bold I^{\bold d_2}_{h(\alpha),0}$ for $\alpha \le \alpha_1$
\sn
\item "{${{}}$}"  $(e) \quad \bold J^{\bold d_1}_{\alpha,1} = \bold
J^{\bold d_2}_{h(\alpha),0}$ for $\alpha < \alpha_1$
\sn
\item "{${{}}$}"  $(f) \quad \bold I^{\bold d^1}_{\alpha_1 + 1,0} =
\bold I^{\bold d^2}_{\alpha_2 +1,0}$
\sn
\item "{${{}}$}"  $(g) \quad \bold J^{\bold d^1}_{\alpha_1,0} =
\emptyset = \bold J^{\bold d^1}_{\alpha_1,1}$.
\ermn
Now check.    \hfill$\square_{\scite{838-6u.17}}$
\enddemo
\bigskip

\proclaim{\stag{838-6u.19} Claim}  For every $\alpha_1 \le \lambda^+$ there
is $\alpha_2 \le \lambda^+$ (in fact $\alpha_2 = \lambda \alpha_1$ is
O.K.) such that: for every $M \in K_{\frak s}$ there is a ${\frak
u}$-free $(\alpha_2,0)$-rectangle $\bold d_2$ with $M^{\bold
d_2}_{\alpha_2,0} = M$ such that
\mr
\item "{$(*)$}"  if $\bold d_1$ is a ${\frak u}$-free
$(\alpha_1,0)$-rectangle with $M^{\bold d_1}_{0,0} = M$ \ub{then} there
are $h,f$ as in $\circledast$ of \scite{838-6u.17}.
\endroster
\endproclaim
\bigskip

\demo{Proof}  Let $\langle {\Cal U}_i:i \le \lambda\rangle$ be a
$\subseteq$-increasing continuous sequence of subsets of $\lambda$
such that ${\Cal U}_\lambda = \lambda$,  
min$({\Cal U}_i) \ge i,\lambda = |{\Cal U}_0|$ and $\lambda
= |{\Cal U}_{\alpha +1} \backslash {\Cal U}_\alpha|$ for $\alpha <
\lambda$.  We now choose $(M_i,\bar p_i,a_i)$ by induction on $i \le
\lambda \alpha_1$ such that
\mr
\item "{$\oplus$}"  $(a) \quad \langle M_i:j \le i\rangle$ is
$\le_{\frak s}$-increasing continuous
\sn
\item "{${{}}$}"  $(b) \quad M_0 = M$
\sn
\item "{${{}}$}"  $(c) \quad$ if $i=j+1$ then $M_i$ is $\le_{\frak
s}$-universal over $M_j$
\sn
\item "{${{}}$}"  $(d) \quad$ if $i = \lambda i_1 + i_2$ and 
$i_1 < \alpha_1,i_2 < \lambda$ then $\bar p_i 
= \langle p^{i_1}_\varepsilon:\varepsilon \in
{\Cal U}_{i_2}\rangle$ where
\nl

\hskip25pt $p^{i_1}_\varepsilon \in \cup{\Cal S}^{\text{bs}}_{\frak s}
(M_{\lambda i_1 + i}):i \le i_2\}$
\sn
\item "{${{}}$}"   $(e) \quad$ if $i = \lambda i_1 + i_2,i_1 < \alpha_1,
j \le i_2 < \lambda$
and $q \in {\Cal S}^{\text{bs}}_{\frak s}(M_{\lambda i_1 + i_2})$ then
\nl

\hskip25pt $(\exists^\lambda \zeta \in {\Cal U}_{i_2})(p^i_\zeta =q)$
\sn
\item "{${{}}$}"  $(f) \quad$ if $i = \lambda i_1 +i_2,i_2 = j_2 +1 <
\lambda$, and $j_2 \in {\Cal U}_\varepsilon$
then $a_{i-1} \in M_i$ and 
\nl

\hskip25pt  the type $\ortp_{\frak s} (a_{i-1},M_{i-1},M_i)$ 
is a non-forking extension of $p^{i_1}_{j_2}$.
\ermn
Now choose $\bold d_2 = (\langle M_i:i \le \lambda
\alpha\rangle,\langle a_i:i < \lambda\alpha\rangle)$.
So assume $\bold d_1$ is a ${\frak u}$-free $(\alpha_1,0)$-rectangle
with $M^{\bold d_1}_{0,0} = M = M^{\bold d_2}_{0,0}$.  We now choose a
pair $(f_i,h_i)$ by induction on $i \le \alpha_1$ such that: $f_i$ is
a $\le_{\frak s}$-embedding of $M^{\bold d_1}_i$ into $M^{\bold
d_2}_{\lambda,i},h_i:i \rightarrow \lambda i$ such that $h(j) \in
(\lambda j,\lambda + \lambda),f_i$ is $\subseteq$-increasing, $h_i$ is
$\subseteq$-increasing, $\ortp_{\frak s}(a^{\bold
d_2}_{h(j),0},M^{\bold d_2}_{h(j)+1,0})$ is a non-forking extension of
$f_{j+1}(\ortp_{\frak s}(a^{\bold d_1}_j,M^{\bold d_1}_j,M^{\bold
d_1}_{j+1}))$ and check; in
fact $\oplus$ gives more than necessary.  
\hfill$\square_{\scite{838-6u.19}} \,\,\, \square_{\scite{838-6u.16}}$
\enddemo
\bigskip 

\proclaim{\stag{838-6u.21} Claim}   Assume $\xi \le \lambda^+$ and
$M_* \in K_{\frak s}$ and $p_* \in {\Cal S}^{\text{bs}}_{\frak s}(M_*)$ 
witness that ${\frak s}$ fails the weak density for 
$K^{3,\text{up}}_{{\frak s},\xi}$, see
Definition \scite{838-6u.13}.
\nl
1) If $(M,N,\{a\}) \in K^{3,\text{bs}}_{\frak s}$ and $M_* \le_{\frak
s} M$ and $\ortp_{\frak s}(a,M,N)$ is a non-forking extension of $p_*$
\ub{then} $(M,N,\{a\})$ has the weak $\xi$-uq-invariant coding property for
${\frak u}$; pedantically assuming ${\frak s}$ has fake equality see
Definition \scite{838-e.5D}, see \scite{838-6u.31},
similarly in part (2); on this coding property, 
(see Definition \scite{838-3r.1}(1).
\nl
2) Moreover the triple $(M,N,\{a\}) \in \text{\rm FR}^1_{\frak u}$ has
the semi $\xi$-uq-invariant coding property, (see Definition \scite{838-3r.71}). 
\endproclaim
\bigskip

\remark{Remark}  If below \scite{838-6u.24} suffices for us then part (2)
of \scite{838-6u.21} is irrelevant.
\endremark
\bigskip

\demo{Proof}  1) Read the definitions, i.e. Definition
\scite{838-6u.13}(2) on the one hand and Definition \scite{838-3r.1}(1) on the
other hand.  Pedantically one may worry that in \scite{838-6u.13}(2) we
use $\le^{\text{bs}}_{\frak s}$, where disjointness is not required
whereas in $\le^\ell_{\frak u}$ it is, however as we allow using fake
equality in ${\frak K}_{\frak u}$ this is not problematic.
\nl
2) Similar.  \hfill$\square_{\scite{838-6u.21}}$
\enddemo
\bn
Now we arrive to the main result of the section.
\proclaim{\stag{838-6u.23} Theorem}  $\dot I(\lambda^{++},{\frak K}^{\frak s}) \ge
\mu_{\text{unif}} (\lambda^{++},2^{\lambda^+})$ and even $\dot
I(\lambda^{++},K^{\frak s}(\lambda^+$-saturated above $\lambda)) \ge
\mu_{\text{unif}}(\lambda^{++},2^{\lambda^+})$ and even $\dot
I(K^{{\frak s},{\frak h}}_{\lambda^{++}}) \ge
\mu_{\text{unif}}(\lambda^{++},2^{\lambda^+})$ for any ${\frak
u}_{\frak s}-\{0,2\}$-appropriate ${\frak h}$ \ub{when}:
\mr
\item "{$(a)$}"  $2^{\lambda^+} < 2^{\lambda^{++}}$
\sn
\item "{$(b)$}"   $(\alpha) \quad {\frak s}$ fail the weak density for 
$K^{3,\text{up}}_{{\frak s},\xi}$ where $\xi \le \lambda^+$ or just
for $\xi = \lambda^+$
\sn
\item "{${{}}$}"  $(\beta) \quad$ if $\xi = \lambda^+$ then $2^\lambda
< 2^{\lambda^+}$.
\endroster
\endproclaim
\bn
Before we prove \scite{838-6u.23} we prove a weaker variant when we
strengthen the set theoretic assumption.
\proclaim{\stag{838-6u.24} Theorem}  Like \scite{838-6u.23} but we add to the
assumption
\mr
\item "{$(c)$}"  ${\Cal D}_{\lambda^+}$ is not
$\lambda^{++}$-saturated (see \scite{838-6u.25}(1) below).
\endroster
\endproclaim
\bigskip

\remark{\stag{838-6u.25} Remark}  0) In the section main conclusion,
\scite{838-6u.27}, if we add clause (c) of \scite{838-6u.24} to the
assumptions then we can rely there on \scite{838-6u.24} instead of on
\scite{838-6u.23}.
\nl
1) Recall that $\lambda > 
\aleph_0 \Rightarrow {\Cal D}_{\lambda^+}$ is not 
$\lambda^+$-saturated by Gitik-Shelah \cite{GiSh:577}, hence this
extra set theoretic assumption is quite weak.
\nl
2) We use \scite{838-3r.7} in proving \scite{838-6u.24}.
\nl
3) We can add the version with ${\frak h}$ to the other such theorems.
\endremark
\bigskip

\demo{Proof of \scite{838-6u.24}}  We can choose a stationary $S \subseteq
\partial = \lambda^+$ such that ${\Cal D}_{\lambda^+} + (\lambda^+
\backslash S)$ is not $\lambda^{++}$-saturated.
We shall apply Theorem \scite{838-3r.7} for the $S$ we just chose for
$\xi'$ which is $\lambda +1$ if $\xi = \lambda$ and is $\xi$ if $\xi <
\lambda^+$.

We have to verify \scite{838-3r.7}'s assumption 
(recalling $\partial = \lambda^+$):
clauses (a) + (b) of \scite{838-3r.7} holds by clauses (a) + (c) of the
assumption of \scite{838-6u.24} if $\xi < \lambda^+$ and clauses (a) +
(b)$(\beta)$ of the assumptions of \scite{838-6u.24} if $\xi = \lambda^+$.
Clause (c) of \scite{838-3r.7} holds by \scite{838-6u.3}, \scite{838-6u.21}(1), whose
assumption holds by clause (b) of the assumption of \scite{838-6u.23}.
Really we have to use \scite{838-6u.21}(1), \scite{838-1a.15}(6).
\hfill$\square_{\scite{838-6u.23}}$ 
\enddemo
\bigskip

\demo{\stag{838-6u.27} Conclusion}  Assume $2^{\lambda^+} <
2^{\lambda^{++}}$ and $\xi \le \lambda^+$ but $\xi = \lambda^+
\Rightarrow 2^\lambda < 2^{\lambda^+}$.
 \nl
1) If $\dot I(\lambda^{++},K^{\frak s}(\lambda^+$-saturated)) is $<
\mu_{\text{unif}}(\lambda^{++},2^{\lambda^+})$ and ${\frak K}_{\frak
s}$ is categorical \ub{then} ${\frak s}$ has existence for
$K^{3,\text{up}}_{{\frak s},\xi}$ for every $\xi \le \lambda^+$.
\nl
2) Similarly for $\dot I({\frak K}^{{\frak s},{\frak
h}}_{\lambda^{++}})$ for any ${\frak u}_{\frak s}-\{0,2\}$-appropriate
${\frak h}$. 
\enddemo
\bigskip

\demo{Proof}  Let $\xi \le \lambda^+$.  
We first try to apply Theorem \scite{838-6u.23}.  Its
conclusion fails, but among its assumptions clauses (a) and (b)$(\beta)$ hold
by our present assumptions.  
So necessarily the demand in clause (b)$(\alpha)$ of 
\scite{838-6u.23} fails.  So we
have deduced that ${\frak s}$ has weak density for
$K^{3,\text{up}}_{{\frak s},\xi}$.  By Observation \scite{838-6u.14}(1),
recalling we are assuming that ${\frak K}_{\frak s}$ is categorical, 
it follows that ${\frak s}$ has almost existence 
for $K^{3,\text{up}}_{{\frak s},\xi}$.
By Claim \scite{838-6u.16} again recalling we are assuming ${\frak
K}_{\frak s}$ is categorical we can deduce that ${\frak s}$ has existence for 
$K^{3,\text{up}}_{{\frak s},\xi}$.

So we have gotten the desired conclusion.  But we
still have to prove Theorem \scite{838-6u.23} in the general case.
\hfill$\square_{\scite{838-6u.27}}$ 
\enddemo
\bigskip

\proclaim{\stag{838-6u.31} Claim}  Assume ${\frak t}'$ is an almost good
$\lambda$-frame derived from
${\frak t}$ as in Definition \scite{838-e.5D} and ${\frak u}' = 
{\frak u}^1_{{\frak t}'}$ is as defined in Definition \scite{838-5t.19},
i.e. Definition \scite{838-e.5I} (and see Claims \scite{838-e.5L},
\scite{838-5t.21}).

\ub{Then}
\mr
\item "{$\boxtimes$}"  $(a) \quad {\frak t}'$ satisfies all the
assumptions on ${\frak t}$ in \scite{838-6u.1}
\sn
\item "{${{}}$}"  $(b) \quad$ so all that we have proved on $({\frak
t},{\frak u})$ in this section apply to 
\nl

\hskip25pt $({\frak t}',{\frak u}')$, too
\sn
\item "{${{}}$}"  $(c) \quad {\frak u}'$ has fake equality $=_\tau$
(see Definition \scite{838-3r.84}(1)
\sn
\item "{${{}}$}"  $(d) \quad {\frak u}'$ is hereditary for the fake
equality $=_\tau$ (see Definition \scite{838-3r.84}(4)).
\endroster
\endproclaim
\bigskip

\demo{Proof}  Clause (a) holds by Claim \scite{838-e.5E}.  Clause (b)
holds by Claim \scite{838-e.5L}, \scite{838-5t.21}.  Clause (c) holds by
direct inspection on \scite{838-e.5I}(6)$(\gamma)$.

In clause (d), ``${\frak u}'$ is hereditary for the fake equality
$=_\tau$", i.e. satisfies clause (a)
of Definition \scite{838-3r.84}(3) by Claim \scite{838-e.5L}(6)$(\beta)$,
\scite{838-5t.21} applied to ${\frak t}'$.

Lastly, to prove ``${\frak u}'$ is hereditary for the fake equality
$=_\tau$" we have to show that it satisfies clause (b) of
\scite{838-3r.84}(4), which holds by \scite{838-e.5L}(6)$(\delta)$.
\hfill$\square_{\scite{838-6u.31}}$ 
\enddemo
\bn
\margintag{6u.34}\ub{\stag{838-6u.34} Proof of \scite{838-6u.23}}:  We shall use Claim
\scite{838-6u.31} to derive $({\frak s}',{\frak u}')$, so we can 
use the results of this section to $({\frak
s},{\frak u})$ and to $({\frak s}',{\frak u}')$.  Now by \scite{838-6u.3}
for some ${\frak u}_{{\frak s}'}-2$-appropriate ${\frak h}$, every $M
\in K^{{\frak u}',{\frak h}}_{\partial^+}$ is $\tau$-fuller, see
Definition \scite{838-1a.15}(6), so by \scite{838-1a.15}(6) it is enough to prove
 Theorem \scite{838-6u.23} for
$({\frak s}',{\frak u}')$.  Now by Claim \scite{838-6u.21}(1) and clause
(b)$(\alpha)$ of the assumption we know that some $(M,N,\{a\}) \in \text{
FR}^1_{{\frak u}'}$ has the semi uq-invariant coding property 
(for ${\frak u}'$).  Also ${\frak u}'$ has the fake equality $=_\tau$
and is hereditary for it by \scite{838-6u.31} and is self dual by \scite{838-5t.21}(1).

Hence in Claim \scite{838-3r.88} all the assumptions hold for ${\frak
u}',(M,N,a)$, hence its conclusion holds, i.e. $(M,N,\{a\})$ has the
weak vertical $\xi$-uq-invariant coding property.  This means that clause
(b) from the assumptions of Theorem \scite{838-3r.91} holds.  Clause (a)
there means $2^\lambda < 2^{\lambda^+} < 2^{\lambda^{++}}$ (choosing
$\theta := \lambda,\partial = \lambda^+$) and clause (c) holds by
\scite{838-6u.3}.  So we are done.  Having shown that the assumptions of
Theorem \scite{838-3r.91} hold, we get its conclusion, which is the
conclusion of the present theorem (reclaling we show that it suffices
to prove it for ${\frak s}'$, so we are done.
\hfill$\square_{\scite{838-6u.23}}$ 
\goodbreak

\head {\S7 Pseudo uniqueness} \endhead  \resetall \sectno=7
 \spuriousreset
\bigskip

Our explicit main aim is to help in \S8 to show that under the
assumptions of \chaptercite{E46}, i.e. \cite{Sh:576} we can get a good
$\lambda$-frame not just a good $\lambda^+$-frame as done in
\sectioncite[\S3]{600}(D),3.7. 
For this we deal with almost good frames (see
\scite{838-7v.1} and Definition \scite{838-5t.3}) and assume existence for
$K^{3,\text{up}}_{{\frak s},\lambda^+}$ (see Definition 
\scite{838-6u.7}(3A) justified by \scite{838-6u.27}) and get
enough of the results of \sectioncite[\S6]{600} and few from
\xCITE{705}.  This means that 
WNF$_{\frak s}$ is defined in \scite{838-7v.5} and proved
to be so called ``a weak non-forking relation on ${\frak K}_{\frak s}$
respecting ${\frak s}$"; we also look at almost good $\lambda$-frames
with such relations and then prove that they are good
$\lambda$-frames in \scite{838-7v.37}(2).  Those results
are used in the proof of \scite{838-e.6A} (and in \S8).

But this has also interest in itself as in general we 
like to understand pre-$\lambda$-frames which are not as
good as the ones considered in \xCITE{705}, i.e. weakly
successful good $\lambda$-frames.  
We will try to comment on this, too.  Note that below even if ${\frak
s}$ is a good $\lambda$-frame satisfying Hypothesis \scite{838-7v.1} which
is weakly successful (i.e. we have existence for
$K^{3,\text{uq}}_{\frak s}$, still
WNF$_{\frak s}$ defined below is not in general equal to NF$_{\frak s}$).
We may wonder, is the assumption (3) of
\scite{838-7v.1} necessary?  The problem is in \scite{838-7v.27}, \scite{838-7v.21}.
\bn
Till \scite{838-7v.38} we use:
\demo{\stag{838-7v.1} Hypothesis}  1) ${\frak s}$ is an 
almost good $\lambda$-frame (see Definition \scite{838-5t.3}).
\nl
2) ${\frak s}$ has existence\footnote{if ${\frak s}$ is a good
 $\lambda$-frame, then actually $K^{3,\text{up}}_{{\frak s},<
 \lambda^+}$ is enough, see \scite{838-6u.7}(3A); 
the main difference is in the proof of \scite{838-7v.27}}
 for $K^{3,\text{up}}_{{\frak s},\lambda^+}$, see
Definition \scite{838-6u.7}(3A) and sufficient condition in \scite{838-6u.27}.
\nl
3) ${\frak s}$ is categorical in $\lambda$ (used only from
\scite{838-7v.27} on).
\nl
4) ${\frak s}$ has disjointness (see Definition \scite{838-5t.9};
used only from \scite{838-7v.19} on, just for transparency, in fact
follows from parts (1) + (3) by \scite{838-5t.45}). 
\enddemo
\bigskip

\definition{\stag{838-7v.3} Definition}  1) Let ${\frak u} = {\frak u}_{\frak s}$ be
as in Definition \scite{838-e.5I} and \scite{838-5t.19}. 
\nl
2) Let FR$_{\frak s}$ be FR$^\ell_{{\frak u}_{\frak s}}$ for $\ell=1,2$
(they are equal).
\enddefinition
\bigskip

\definition{\stag{838-7v.5} Definition}  1) Assume $\xi \le \lambda^+$.
Let WNF$^{\xi}_{\frak s}(M_0,N_0,M_1,N_1)$ mean that: $M_0 \le_{\frak
s} N_0 \le_{\frak s} N_1,M_0 \le_{\frak s} M_1 \le_{\frak s} N_1$ and
 if $\alpha < \xi$ so $\alpha < \lambda^+$ 
and $\bold d$ is an ${\frak u}$-free $(0,\alpha)$-rectangle
and $f$ is a $\le_{\frak s}$-embedding of $N_0$ into $M^{\bold d}_{0,\alpha}$
such that $f(M_0) = M^{\bold d}_{0,0}$ \ub{then} we can find a model
$N^*$ and a ${\frak u}$-free $(1,\alpha)$-rectangle $\bold d^+$
satisfying $\bold d^+ \restriction (0,\alpha) = \bold d$ and $M^{\bold
d^+}_{1,\alpha} \le_{\frak s} N^*$ and 
 $\le_{\frak s}$-embedding $g \supseteq f$ of
$N_1$ into $N^*$ such that $M^{\bold d^+}_{1,0} = g(M_1)$.
\nl
2) If $\xi = \lambda^+$ we may omit it.  So WNF$^\xi_{\frak s}$ is
also considered as the class of such quadruples of models.
\enddefinition
\bigskip

\remark{\stag{838-7v.7} Remark}  0) Definition \scite{838-7v.5}(1) is dull for
$\xi = 0$.
\nl
1) So this definition 
is not obviously symmetric but later we shall prove it is.
\nl
2) Similarly, it seemed that 
the value of $\xi$ is important, but we shall show that
for $\xi < \lambda^+$ large enough it is not \ub{when} ${\frak s}$ is
a good $\lambda$-frame; see e.g. \scite{838-7v.11}.
\nl
3) In Definition \scite{838-7v.5} we may ignore $\xi = 0$ as it
   essentially says nothing.
\endremark
\bigskip

\demo{\stag{838-7v.9} Observation}  1) If $1 \le \xi \le \lambda^+$ and
WNF$^\xi_{\frak s}(M_0,N_0,M_1,N_1)$ and
$(M_0,N_0,a) \in K^{3,\text{bs}}_{\frak s}$ \ub{then} $(M_0,N_0,a)
\le_{\text{bs}} (M_1,N_1,a)$ and in particular $(M_1,N_1,a) \in
K^{3,\text{bs}}_{\frak s}$.
\nl
2) WNF$^\xi_{\frak s}$ is $\subseteq$-decreasing with $\xi$.
\nl
3) WNF$_{\frak s} = \text{ WNF}^{\lambda^+}_{\frak s} =
\cap\{\text{WNF}^\xi_{\frak s}:\xi < \lambda^+\}$. 
\nl
4) In Definition \scite{838-7v.5}(1) in the end we can weaken ``$M^{\bold
   d^+}_{1,0} = g(M_1)$" to ``$g(M_1) \le_{\frak s} M^{\bold d^+}_{1,0}$".
\enddemo
\bigskip

\demo{Proof}  1) Straight, use $\bold d$, the 
${\frak u}$-free $(0,1)$-rectangle
such that $M^{\bold d}_{0,0} = M_0$ and $M^{\bold d}_{0,1} = M$
and $a^{\bold d}_{0,0} = a$, i.e. $\bold I^{\bold d}_{0,0} = \{a\}$.
\nl
2),3) Trivial.
\nl
4) Given $(N^*,\bold d^1,g)$ as in Definition \scite{838-7v.9}(4).  We
   define $\bold d'$, a ${\frak u}$-free $(1,\alpha)$-rectangle by
   $M^{\bold d'}_{i,j}$ is $g(M_1)$ if $(i,j) = (1,0)$ and is
   $M^{\bold d^+}_{i,j}$ when $i \le 1 \and j \le \alpha \and (i,j)
   \ne (1,0)$ and $\bold I^{\bold d'}_{0,j} = 
\bold I^{\bold d'}_{1,j} = \bold I^{\bold d}_{0,j}$ for $j < \alpha$ and $\bold
   J^{\bold d'}_{0,j} = \emptyset$ for $j \le \alpha$.  The only
   non-obvious point is why $(M^{\bold d'}_{0,0},
M^{\bold d'}_{0,1},\bold I^{\bold d'}_{0,0}) \le^1_{\frak u} (M^{\bold
   d'}_{1,0},M^{\bold d'}_{1,1},\bold I^{\bold d'}_{0,0})$ which means
   $(M^{\bold d}_{0,0},M^{\bold d}_{0,1},\bold I^{\bold d}_{0,0})
   \le^1_{\frak u} (g(M_1),M^{\bold d^+}_{1,1},\bold I^{\bold d^+}_{0,0})$.  
This is because ${\frak u}$ is interpolative by
\scite{838-e.5L}(6)$(\varepsilon)$, see Definition \scite{838-3r.89}.
\hfill$\square_{\scite{838-7v.9}}$
\enddemo
\bigskip

\proclaim{\stag{838-7v.11} Claim}  [Monotonicity]  Assume $1 \le \xi \le
\lambda^+$.  
If {\rm WNF}$^\xi_{\frak s}
(M_0,N_0,M_1,N_1)$ and $M_0 \le_{\frak s} N'_0 \le_{\frak s} N_0$
and $M_0 \le_{\frak s} M'_1 \le M_1$ and $N_1 \le_{\frak s} N'_1,N'_0
\cup M'_1 \subseteq N''_1 \le_{\frak s} N'_1$ 
\ub{then} {\rm WNF}$^\xi_{\frak s} (M_0,N'_0,M'_1,N''_1)$ holds.
\endproclaim
\bigskip

\demo{Proof}  It is enough to prove that for the case three of
the equalities $N'_0 = N_0,M'_1 = M_1,N''_1 = N'_1,N'_1 = N_1$ hold.  Each
follows: in the case $N'_0 \ne N_0$ by the Definition \scite{838-7v.5}, in the
case $M'_1 \ne M_1$ by \scite{838-7v.9}(4), and in the cases $N''_1 = N'_1
\vee N'_1 = N_1$ by amalgamation (in ${\frak K}_{\frak s}$) and the definition
of WNF$^\xi_{\frak s}$ and \scite{838-7v.9}(4).
\nl
${{}}$  \hfill$\square_{\scite{838-7v.11}}$
\enddemo
\bigskip

\demo{\stag{838-7v.13} Observation}  [${\frak s}$ categorical in $\lambda$
or $M_0$ is brimmed or ${\frak s}$ has existence for
$K^{3,\text{rd}}_{\frak s}$, see Definition \scite{838-5t.41}.]

If $\xi > \xi^{\text{rd}}_{M_0}$, see Definition \scite{838-5t.41} 
and WNF$^\xi_{\frak s}(M_0,N_0,M_1,N_1)$ 
\ub{then} $M_1 \cap N_0 = M_0$. 
\enddemo
\bigskip

\remark{Remark}  1) Recall that if 
${\frak s}$ is an almost good $\lambda$-frame
\ub{then} it has density of $K^{3,\text{rd}}_{\frak s}$ hence if
${\frak s}$ is categorical then it also has 
existence for $K^{3,\text{rd}}_{\frak s}$.

It is convenient to assume this but not essential.  Proving density
for $K^{3,\text{up}}_{\frak s}$ we actualy prove density for
$K^{3,\text{up}}_{\frak s} \cap K^{3,\text{rd}}_{\frak s}$; moreover
$K^{3,\text{up}}_{\frak s} \subseteq K^{3,\text{rd}}_{\frak s}$.
\nl
2) Recall $\xi^{\text{rd}}_{M_0} = \xi^{\text{rd}}_{\frak s}$ if $M_0$
   is superlimit, e.g. when $K_{\frak s}$ is categorical.
\endremark
\bigskip

\demo{Proof}  Recall that letting $\alpha = \xi^{\text{rd}}_{M_0}$
 by \scite{838-5t.43}(3) or \scite{838-5t.43}(4) there is a 
${\frak u}$-free $(0,\alpha)$-rectangle
$\bold d$ such that $M^{\bold d}_{0,0} = M_0,N_0 \le_{\frak s}
M^{\bold d}_{0,\alpha}$ and each $(M^{\bold d}_{0,\alpha},M^{\bold
d}_{0,\alpha +1},a^{\bold d}_{0,\alpha})$ is reduced (see Definition
\scite{838-5t.41}).  Now apply Definition \scite{838-7v.5} to this $\bold d$.
Alternatively, recall for ${\frak u}$-free $(\alpha,\beta)$-rectangle
 (or $(\bar\alpha,\beta)$-triangle) $\bold d$ we have $M^{\bold
 d}_{i_1,j_1} \cap M^{\bold d}_{i_2,j_2} = M^{\bold
 d}_{\text{min}\{i_1,i_2\},\text{min}\{j_1,j_2\}}$, or we can use
 \scite{838-6u.9}(0).   \hfill$\square_{\scite{838-7v.13}}$
\enddemo
\bigskip

\proclaim{\stag{838-7v.15} Claim}  [Long transitivity]

Assume $1 \le \xi \le \lambda^+$.  We have 
{\rm WNF}$^\xi_{\frak s}(M_0,N_0,M_{\alpha(*)},N_{\alpha(*)})$ \ub{when}:
\mr
\item "{$(a)$}"  $\langle M_\alpha:\alpha \le \alpha(*)\rangle$ is
$\le_{\frak s}$-increasing continuous
\sn
\item "{$(b)$}"  $\langle N_\alpha:\alpha \le \alpha(*)\rangle$ is
$\le_{\frak s}$-increasing continuous
\sn
\item "{$(c)$}"  {\rm WNF}$^\xi_{\frak s}(M_\alpha,N_\alpha,M_{\alpha
+1},N_{\alpha +1})$ for every $\alpha < \alpha(*)$.
\endroster
\endproclaim
\bigskip

\remark{Remark}  1) Recall that we do not know symmetry for
WNF$_{\frak s}$ and while this claim is easy its dual is not clear at
this point.
\endremark
\bigskip

\demo{Proof}  By chasing arrows.  \hfill$\square_{\scite{838-7v.15}}$
\enddemo
\bigskip

\proclaim{\stag{838-7v.17} Claim}  [weak existence]  Assume $1 \le \xi 
\le \lambda^+$.

If $(M_0,M_1,a) \le_{\text{bs}} (N_0,N_1,a)$ and $(M_0,M_1,a) \in
K^{3,\text{up}}_{{\frak s},\xi}$ \ub{then} 
{\rm WNF}$^\xi_{\frak s}(M_0,N_0,M_1,N_1)$.
\endproclaim
\bigskip

\demo{Proof}  Let $\alpha \le \xi$ be such that 
$\alpha < \lambda^+$ and $\bold d$ 
be a ${\frak u}$-free $(0,\alpha)$-rectangle
and let $f$ be a $\le_{\frak s}$-embedding of $N_0$ into $M^{\bold
d}_{0,\alpha}$ such that $f(M_0) = M^{\bold d}_{0,0}$.  Let $N'_1 \in
K_{\frak s}$ be $\le_{\frak s}$-universal over $M^{\bold
d}_{0,\alpha}$, exist by \scite{838-5t.15}, and 
let $N'_0 = M^{\bold d}_{0,\alpha}$.

As $(N_0,N_1,a) \in K^{3,\text{bs}}_{\frak s}$ we 
can find a $\le_{\frak s}$-embedding $g$ of $N_1$ into $N'_1$
extending $f$ such that $(g(N_0),g(N_1),g(a)) \le_{\text{bs}}
(N'_0,N'_1,g(a))$ so as $\le_{\text{bs}}$ is a partial order preserved
by isomorphisms, clearly $(g(M_0),g(M_1),g(a)) \le_{\text{bs}}
(N'_0,N'_1,g(a))$.  Now as $(M_0,M_1,a) \in K^{3,\text{up}}_{{\frak s},\xi}$
it follows that $(g(M_0),g(M_1),g(a)) \in K^{3,\text{up}}_{{\frak s},\xi}$.
Applying the definition of $K^{3,\text{up}}_{{\frak s},\xi}$, see
Definition \scite{838-6u.7}(1) with
$g(M_0),g(M_1),g(a),N'_0,N'_1$, {\rm dual}$(\bold d)$ here standing for
$M,N,a,M',N',\bold d$ there, recalling that ${\frak u}$ is self-dual
we can find a ${\frak u}$-free $(\alpha +1,1)$-rectangle $\bold d^*$
and $h$ such that: $\bold d^* \rest (\alpha,0) =$ {\rm dual}$(\bold
d),h(g(M_1) \le_{\frak s} M^{\bold d^*}_{0,1},h(g(a)) = a^{\bold
d^*}_{0,0}$ and $h$ is a $\le_{\frak s}$-embedding of $N'_1$ into
$M^{\bold d^*}_{\alpha +1,1}$ over $N'_0 = M^{\bold d}_{0,\alpha}$
hence $h \rest g(M_0) = \text{ id}_{g(M_0)}$.

Now $N'_1,\bold d^+ := \text{\rm dual}(\bold d^*) \rest (1,\alpha)$ and
$h \circ g$ are as required in Definition \scite{838-7v.5} (standing for
$N^*,\bold d^+,g$ there) recalling \scite{838-7v.9}.
  \hfill$\square_{\scite{838-7v.17}}$
\enddemo
\bigskip

\proclaim{\stag{838-7v.19} Lemma}  [Amalgamation existence]  Let $\xi <
\lambda^+$ or just $\xi \le \lambda^+$ and $\xi \ge 1$.
\nl
1) If $M_0 \le_{\frak s} M_\ell$ for $\ell=1,2$ and $M_1 \cap M_2 =
M_0$ \ub{then} for some $M_3$ we have 
{\rm WNF}$^\xi_{\frak s}(M_0,M_1,M_2,M_3)$.
\nl
2) If $M_0 \le_{\frak s} M_2$ \ub{then} we can find an 
${\frak u}$-free rectangle $\bold d$ satisfying $\beta_{\bold d} = 0$
such that
\mr
\item "{$\boxtimes$}"  $(a) \quad M^{\bold d}_{0,0} = M_0$
\sn
\item "{${{}}$}"  $(b) \quad M_2 \le_{\frak s} M^{\bold d}_{\alpha(\bold d),0}$
\sn
\item "{${{}}$}"  $(c) \quad (M^{\bold d}_{\alpha,0},M^{\bold d}_{\alpha +1,0},
b^{\bold d}_{\alpha,0})$ belongs to $K^{3,\text{up}}_{{\frak s},\xi}$
for $\alpha < \alpha_{\bold d}$
\sn
\item "{${{}}$}"  $(d) \quad$ if $M_2$ is $(\lambda,*)$-brimmed over $M_0$ then
$M^{\bold d}_2 = M^{\bold d}_{\alpha_{\bold d},0}$.
\ermn
2A) If $(M_0,M_2,b) \in K^{3,\text{bs}}_{\frak s}$ we can add
$b^{\bold d}_{0,0} = b$.
\nl
2B) Assume $K^{3,*}_{\frak s} \subseteq K^{3,\text{bs}}_{\frak s}$ is
such that ${\frak s}$ has existence for $K^{3,*}_{\frak s}$ \ub{then}
in parts 2), 2A) we can replace $K^{3,\text{up}}_{{\frak s},\xi}$ by
$K^{3,*}_{\frak s}$.
\nl
3) In part (1) if $(M_0,M_\ell,b_\ell) \in K^{3,\text{bs}}_{\frak s}$
for $\ell=1,2$ \ub{then} we can add
$(M_0,M_\ell,b_\ell) \le_{\text{bs}} (M_{3 -\ell},M_3,b_\ell)$ for $\ell=1,2$.
\endproclaim
\bigskip

\demo{Proof}  1) Follows by part (3).
\nl
2) By Ax(D)(c), density, of almost good $\lambda$-frames 
there is $b \in M_2 \backslash M_0$ such that
$\ortp_{\frak s}(b,M_0,M_2) \in {\Cal S}^{\text{bs}}_{\frak s}(M_0)$,
hence $(M_0,M_2,b) \in K^{3,\text{bs}}_{\frak s}$, by the definition
of $K^{3,\text{bs}}_{\frak s}$ it follows that so we can apply
part (2A).
\nl
2A) By part (2B). 
\nl
2B) So let $(M_0,M_2,b) \in K^{3,*}_{\frak s}$ be given.  We
try to choose $(M_{0,\alpha},M_{2,\alpha})$ and if $\alpha = \beta +1$
also $a_\beta$ by induction on $\alpha < \lambda^+$ such that:
\mr
\item "{$\circledast$}"  $(a) \quad M_{\ell,\alpha} \in K_{\frak s}$
is $\le_{\frak s}$-increasing continuous for $\ell=0,2$
\sn
\item "{${{}}$}"  $(b) \quad M_{0,\alpha} \le_{\frak s} M_{2,\alpha}$
\sn
\item "{${{}}$}"  $(c) \quad (M_{0,\alpha},M_{2,\alpha}) = (M_0,M_2)$
for $\alpha =0$
\sn
\item "{${{}}$}"  $(d) \quad$ if $\alpha = \beta +1$ then
$(M_{0,\beta},M_{0,\alpha},a_\beta) \in K^{3,*}_{\frak s}$ and
$a_\beta \in M_{2,\alpha} \backslash M_{0,\alpha}$
\sn
\item "{${{}}$}"  $(e) \quad$ if $\alpha = \beta +1$ then
$M_{2,\alpha}$ is $\le_{\frak s}$-brimmed over $M_{2,\beta}$
\sn
\item "{${{}}$}"  $(f) \quad$ if $\alpha = 0$ then $a_\alpha = b$.
\ermn
By Fodor lemma we cannot choose for every $\alpha < \lambda^+$.  For
$\alpha = 0$ and $\alpha$ limit there are no problems, hence for some
$\alpha = \beta +1$, we have defined up to $\beta$ but cannot define
for $\alpha$ clearly $\beta < \lambda^+$.  First assume
\mr
\item "{$(*)$}"  $M_{0,\beta} \ne M_{2,\beta}$.
\ermn
So by Ax(D)(c) of Definition \scite{838-5t.3} of almost good
$\lambda$-frame 
we can
choose $a_\beta \in M_{2,\beta} \backslash M_{0,\beta}$ such that
$\ortp_{\frak s}(b,M_{0,\beta},M_{2,\beta}) \in {\Cal
S}^{\text{bs}}_{\frak s}(M_{0,\beta})$ and $a_\beta = b$ if $\beta =
0$.  By the assumption on $K^{3,*}_{\frak s}$ there is $N_\beta \in
K_{\frak s}$ such that $(M_{0,\alpha},N_\beta,a_\beta) \in
K^{3,*}_{\frak s}$ and $\ortp_{\frak s}(a_\beta,M_{0,\beta},N_\beta) =
\ortp_{\frak s}(a_\beta,M_{0,\beta},M_{2,\beta})$.

By the definition of orbital type (and amalgamation of 
${\frak K}_{\frak s}$) without loss of generality  for some $M'_{2,\beta}$ we have $N_\beta
\le_{\frak s} M'_{2,\beta}$ and $M_{2,\beta} \le_{\frak s}
M'_{2,\beta}$.

Let $M_{2,\alpha} \in K_{\frak s}$ be brimmed over $M'_{2,\beta}$.  So
we can choose for $\alpha$, contradiction.

Hence $(*)$ cannot hold so
$M_{0,\beta} = M_{2,\beta}$, easily $\beta \ge 1$ (as $M_2 \ne
M_0$) and by clause (e) of $\circledast$ $M_{2,\beta}$ is brimmed over $M_{2,0} = M_2$ hence over
$M_0$.  What about clause (d) of the conclusion?  It follows because any
two brimmed extensions of $M_0$ are isomorphic over it by
\scite{838-5t.36}(4) and with a little more work even over $M_0 \cup \{b\}$.
\nl
3) So let $K^{3,*}_{\frak s} = K^{3,\text{up}}_{\frak s}$ or just
   $K^{3,*}_{\frak s} \subseteq K^{3,\text{up}}_{{\frak s},\xi}$ and ${\frak
   s}$ has existence for $K^{3,*}_{\frak s}$.

Let $\bold d$ be as guaranteed in parts (2),(2A) so $a^{\bold
d}_{0,0} = b_2$ and $M_0 = M^{\bold d}_{0,0},M_2 \le_{\frak s}
M^{\bold d}_{\alpha_{\bold d},0}$.  Without loss of generality
$M^{\bold d}_{\alpha_{\bold d},0} \cap M_1 = M_0$ and now we choose
$N_\alpha$ by induction on $\alpha \le \alpha_{\bold d}$ such that
\mr
\item "{$\boxplus$}"  $(a) \quad N_\alpha \in K_{\frak s}$ is
$\le_{\frak s}$-increasing continuous
\sn
\item "{${{}}$}"  $(b) \quad M^{\bold d}_{\alpha_{\bold d},0} \cap
N_\alpha = M^{\bold d}_{\alpha,0}$
\sn
\item "{${{}}$}"  $(c) \quad M^{\bold d}_{\alpha,0} \le_{\frak s}
N_\alpha$
\sn
\item "{${{}}$}"  $(d) \quad (M^{\bold d}_{\alpha,0},M^{\bold
d}_{\alpha +1,0},a^{\bold d}_{\alpha,0}) \le_{\text{bs}}
(N_\alpha,N_{\alpha +1},a^{\bold d}_{\alpha,0})$
\sn
\item "{${{}}$}"  $(e) \quad N_\alpha = M_1$ for $\alpha = 0$.
\ermn
There is no problem to carry the choice by Hypothesis \scite{838-7v.1}(4)
and Definition \scite{838-5t.9}.   

Now for each $\alpha < \alpha_{\bold d}$ by clause (c) of $\boxtimes$
of part (2) or (2B) we have
$$(M^{\bold d}_{0,\alpha},M^{\bold d}_{0,\alpha +1},
a^{\bold d}_{\alpha,0}) \in K^{3,*}_{\frak s}
\subseteq K^{3,\text{up}}_{{\frak s},\xi}$$
recalling the choice of $\bold d$ and by clause (d) of $\boxplus$ we have
$$(M^{\bold d}_{0,\alpha},M^{\bold d}_{0,\alpha +1},
a^{\bold d}_{\alpha,0}) \le_{\text{bs}}
(N_\alpha,N_{\alpha +1},a^{\bold d}_{\alpha,0}),$$ hence by the weak 
existence Claim \scite{838-7v.17} we have WNF$^\xi_{\frak
s}(M^{\bold d}_{0,\alpha},N_\alpha,M^{\bold d}_{0,\alpha +1},N_{\alpha +1})$.

As $\langle M^{\bold d}_{0,\alpha}:\alpha \le 
\alpha_{\bold d}\rangle$ and $\langle
N_\alpha:\alpha \le \alpha_{\bold d}\rangle$ are $\le_{\frak s}$-increasing
continuous, it follows the long transitivity claim \scite{838-7v.15} that
WNF$^\xi_{\frak s}(M^{\bold d}_{0,0},N_0,M^{\bold d}_{\alpha_{\bold
d},0},N_{\alpha(\bold d)})$ which means that
WNF$^\xi_{\frak s}(M_0,M_1,M^{\bold d}_{\alpha_{\bold d},0},
N_{\alpha(\bold d)})$.
Let $M_3 := N_{\alpha(\bold d)}$, but now 
$M_0 \le_{\frak s} M_2 \le_{\frak s} M^{\bold d}_{\alpha(\bold d),0}$ 
hence by the monotonicity Claim \scite{838-7v.11} we have 
WNF$^\xi_{\frak s}(M_0,M_1,M_2,M_3)$.

This proves the desired conclusion of part (1), but there are more demands in
part (3).  One is $(M_0,M_1,b_1) \le_{\text{bs}} (M_2,M_3,b_1)$, but
$M_1 = N_0$ and $M_3 = N_{\alpha(\bold d)}$ 
so this means $(M_0,N_0,b_1) \le_{\text{bs}}
(M_2,N_{\alpha(\bold d)},b_1)$ and by monotonicity of non-forking it
suffices to show $(M_0,N_0,b_1) \le_{\text{bs}} (M^{\bold
d}_{\alpha(\bold d),0},N_{\alpha(\bold d)})$.  But recall
WNF$^\xi_{\frak s}(M_0,N_0,M^{\bold d}_{\alpha(\bold d),0},
N_{\alpha(\bold d)})$ and this implies $(M_0,N_0,b_1)
\le_{\text{bs}} (M^{\bold d}_{\alpha(\bold d),0},N_{\alpha(\bold d)},b_1)$
 by Observation \scite{838-7v.9}(1) which as said above suffices.

Now also we have chosen $a^{\bold d}_{0,0}$ as $b_2$, so by clause (d)
of $\boxplus$ for $\alpha =0$ we have easily 
$(M_0,M^{\bold d}_{1,0},b_2) = (M^{\bold d}_{0,0},M^{\bold d}_{1,0},
a^{\bold d}_{0,0}) \le_{\text{bs}} (N_0,N_1,a^{\bold d}_{0,0}) =
(M_1,N_1,b_2) \le_{\text{bs}} (M_1,N_{\alpha(\bold d)},b_2) 
= (M_1,M_3,b_2)$ but 
$M_0 \le_{\frak s} M_1 \le_{\frak s} M^{\bold
d}_{\alpha(\bold d),0}$ so by easy monotonicity we have $(M_0,M_2,b_2)
\le_{\text{bs}} (M_1,M_3,b_2)$, as desired in part (3);
so we are done.
\nl
${{}}$    \hfill$\square_{\scite{838-7v.19}}$ 
\enddemo
\bigskip

\remark{\stag{838-7v.21} Remark}   In the proof of \scite{838-7v.19}(2B), if 
${\frak s}$ is a good $\lambda$-frame, in fact, $\lambda$ steps 
in the induction suffice by a careful choice of $a_\beta$ using
bookkeeping as in the proof of \scite{838-5t.36}(1),
so we get $\alpha_{\bold d} = \lambda$.  Without this extra
hypothesis on ${\frak s}$, this is not clear.
\endremark
\bigskip

\proclaim{\stag{838-7v.23} Claim}  Assume that $1 \le \xi \le \lambda^+,
\alpha < \lambda^+,\bold d$ is a 
${\frak u}$-free $(0,\alpha)$-rectangle and $M_0 \le_{\frak s}
M^{\bold d}_{0,0} \le_{\frak s} M^{\bold d}_{0,\alpha} \le_{\frak s}
N_1$.  \ub{Then} we can find $\alpha',\bold d',h$ such that
\mr 
\item "{$\oplus$}"  $(a) \quad \alpha' \in [\alpha,\lambda^+)$
\sn
\item "{${{}}$}"  $(b) \quad \bold d'$ is 
a ${\frak u}$-free $(0,\alpha')$-rectangle
\sn
\item "{${{}}$}"  $(c) \quad h$ is an increasing function, $h:\alpha +1
\rightarrow \alpha' +1$
\sn
\item "{${{}}$}"  $(d) \quad M^{\bold d'}_{0,0} = M_0$
\sn
\item "{${{}}$}"  $(e) \quad N_1 \le_{\frak s} M^{\bold d'}_{0,\alpha'}$
\sn
\item "{${{}}$}"  $(f) \quad i \le \alpha \Rightarrow M^{\bold d}_{0,i}
\le_{\frak s} M^{\bold d'}_{0,h(i)}$
\sn
\item "{${{}}$}"  $(g) \quad$ for $i < \alpha,a^{\bold d}_{0,i} = a^{\bold
d'}_{0,h(i)}$ and $\ortp_{\frak s}(a^{\bold d'}_{0,h(i)},M^{\bold
d'}_{0,h(i)},M^{\bold d'}_{0,h(i)+1})$ does not
\nl

\hskip25pt  fork over $M^{\bold d}_{0,i}$
\sn
\item "{${{}}$}"  $(h) \quad (M^{\bold d'}_{0,\beta},
M^{\bold d'}_{0,\beta +1},a^{\bold d'}_{0,\beta}) 
\in K^{3,\text{up}}_{{\frak s},\xi}$ for $\beta < \alpha'$.
\endroster
\endproclaim
\bigskip

\demo{Proof}  We can choose $M'_i$ by induction on $i \le 1 + \alpha
+1$ such that
\mr
\item "{$(\alpha)$}"   $\langle M'_j:j \le i\rangle$ is increasing
continuous
\sn
\item "{$(\beta)$}"   $M'_0 = M_0$
\sn
\item "{$(\gamma)$}"  $M'_i$ is brimmed over $M'_j$ if $i = j+1 \le 1
+ \alpha +1$ 
\sn
\item "{$(\delta)$}"  $M_\alpha \cap M'_{1+i} = M^{\bold d}_{0,i}$ if
$i \le \alpha$
\sn
\item "{$(\varepsilon)$}"  $M^{\bold d}_{0,i} \le_{\frak s} M'_{1+i}$
if $i \le \alpha$
\sn
\item "{$(\zeta)$}"  $\ortp_{\frak s}(a^{\bold d}_{0,i},
M'_i,M'_{i+1})$ does not fork over $M^{\bold d}_{0,i}$ if $i < \alpha$
\sn
\item "{$(\eta)$}"  $N_1 \le_{\frak s} M'_{1 + \alpha +1}$.
\ermn
This is possible because we have disjoint amalgamation (see \scite{838-5t.45}).
Now for each $i \le 1 + \alpha$ use \scite{838-7v.19}(2A) with
$M'_i,M'_{i+1},a^{\bold d}_{0,i}$ here standing for $M_0,M_2,b$ there
(so clause (d) there apply).  \hfill$\square_{\scite{838-7v.23}}$
\enddemo
\bigskip

\remark{\stag{838-7v.25} Remark}  Recall that from now on we are 
assuming that ${\frak K}_{\frak s}$ is categorical.
\endremark
\bigskip

\proclaim{\stag{838-7v.27} Claim}  [Symmetry]  There is $\xi = \xi_{\frak
s} < \lambda^+$ such that $(\xi \ge \xi^{\text{rd}}_{\frak s}$ for
simplicity, see \scite{838-5t.41} and) for every  $\zeta < \lambda^+$, 
if {\rm WNF}$^\xi_{\frak s}(M_0,N_0,M_1,N_1)$ 
\ub{then} {\rm WNF}$^\zeta_{\frak s}(M_0,M_1,N_0,N_1)$ holds.
\endproclaim
\bigskip

\remark{Remark}  1) Yes, the models $N_0,M_1$ exchange places.
\nl
2) Without categoricity, $\xi = \xi_{{\frak s},M_0}$ is O.K.
\endremark
\bigskip

\demo{Proof}  By \scite{838-7v.19}(2) there are $\xi = \xi(*) < \lambda^+$ and
$\bold d$, a ${\frak u}$-free $(0,\xi)$-rectangle with each
$(M^{\bold d}_{0,\alpha},M^{\bold d}_{0,\alpha +1},
a^{\bold d}_{0,\alpha})$ belonging to $K^{3,\text{up}}_{{\frak
s},\zeta}$ for every $\zeta < \lambda^+$
such that $M^{\bold d}_{0,\xi}$ is brimmed over $M^{\bold d}_{0,0}$
and $M_0 = M^{\bold d}_{0,0}$.  Note that
the choice of $\xi$ does not depend on $\zeta,\langle
M_0,N_0,M_1,N_1\rangle$, just on $M_0$ by \scite{838-7v.19} and it does
not depend on $M_0$ recalling $K_{\frak s}$ is categorical.  

As $M^{\bold d}_{0,\xi}$ is $\le_{\frak s}$-universal
over $M^{\bold d}_{0,0}$ without loss of generality  $N_0 \le_{\frak s} 
M^{\bold d}_{0,\xi}$.  

Now let $\zeta < \lambda^+$ and recall that we assume
WNF$^\xi_{\frak s}(M_0,N_0,M_1,N_1)$.  Let $\bold d^+,N^+,f$ be as
guaranteed by Definition \scite{838-7v.5}(1) and by renaming without loss
of generality the function
$f$ is the identity.  Now for each $\alpha < \xi$, we shall apply the weak
existence claim \scite{838-7v.17}, with $M^{\bold d^+}_{0,\alpha},
M^{\bold d^+}_{1,\alpha},M^{\bold d^+}_{0,\alpha +1},M^{\bold
d^+}_{1,\alpha +1},a^{\bold d^+}_{0,\alpha}$ here standing for
$M_0,N_0,M_1,N_1,a$ as there; this is O.K. as its assumptions mean
$(M^{\bold d^+}_{0,\alpha},M^{\bold d^+}_{0,\alpha +1},a^{\bold
d^+}_{0,\alpha}) \le_{\text{bs}} (M^{\bold d^+}_{1,\alpha},
M^{\bold d^+}_{1,\alpha +1},a^{\bold d^+}_{0,\alpha})$ and 
$(M^{\bold d^+}_{0,\alpha},M^{\bold d^+}_{0,\alpha +1},
a^{\bold d^+}_{0,\alpha}) \in K^{3,\text{up}}_{{\frak s},\zeta}$ which
hold by clause $\boxtimes(c)$ of Claim \scite{838-7v.19}(2), i.e. by the
choice of $\bold d$ as $\bold d^+ \rest (0,\alpha(\bold d)) = 
\bold d$.   Hence the
conclusion of \scite{838-7v.17} applies, which gives that
 we have WNF$^\zeta_{\frak s}(M^{\bold
d^+}_{0,\alpha},M^{\bold d^+}_{1,\alpha},M^{\bold d^+}_{0,\alpha +1},
M^{\bold d^+}_{1,\alpha +1})$.
Of course $\langle M^{\bold d^+}_{0,\alpha}:\alpha \le \xi\rangle$
and $\langle M^{\bold d^+}_{1,\alpha}:\alpha \le \xi\rangle$ are
$\le_{\frak s}$-increasing continuous.  Together by the
long transitivity, claim \scite{838-7v.15} we have WNF$^\zeta_{\frak s}(M^{\bold
d^+}_{0,0},M^{\bold d^+}_{1,0},M^{\bold d^+}_{0,\xi},
M^{\bold d^+}_{1,\xi})$.  But $M^{\bold d^+}_{0,0} = M_0,
M^{\bold d^+}_{1,0} = M_1$ and $N_0 \le_{\frak
s} M^{\bold d}_{0,\xi}$ and $N_1 \le_{\frak s} N^+,M^{\bold
d^+}_{1,\xi} 
\le_{\frak s} N^+$ so by the monotonicity claim, 
\scite{838-7v.11}, we have
WNF$^\zeta_{\frak s}(M_0,M_1,N_0,N_1)$ as required.   
\hfill$\square_{\scite{838-7v.27}}$
\enddemo
\bigskip

\demo{\stag{838-7v.29} Conclusion}  If WNF$^\xi_{\frak s}
(M_0,N_0,M_1,N_1)$ and $\xi \ge \xi_{\frak s}$, see \scite{838-7v.27} 
\ub{then} $\zeta < \lambda^+ \Rightarrow 
\text{ WNF}^\zeta_{\frak s}(M_0,N_0,M_1,N_1)$ that is WNF$_{\frak
s}(M_0,N_0,M_1,N_1)$. 
\enddemo
\bigskip

\demo{Proof}  Applying \scite{838-7v.27} twice recalling \scite{838-7v.9}(3)
in the end.  \hfill$\square_{\scite{838-7v.29}}$
\enddemo
\bigskip

\proclaim{\stag{838-7v.31} Claim}  ({\rm WNF}$_{\frak s}$ lifting \ub{or}
weak uniqueness)

If {\rm WNF}$_{\frak s}(M_0,N_0,M_1,N_1)$ and $\alpha < \lambda^+$ and
$\langle M_{0,i}:i \le \alpha\rangle$ is $\le_{\frak s}$-increasing
continuous, $M_{0,0} = M_0$ and $M_{0,\alpha} = N_0$ \ub{then} we can find a
$\le_{\frak s}$-increasing continuous sequence $\langle M_{1,i}:i \le
\alpha +1 \rangle$ such that $M_{1,0} = M_1,N_1 \le_{\frak s}
M_{1,\alpha +1}$ and for each $i < \alpha$ we have {\rm WNF}$_{\frak
s}(M_{0,i},M_{0,i+1},M_{1,i},M_{1,i+1})$ for $i < \alpha$.
\endproclaim
\bigskip

\demo{Proof}  We shall use \scite{838-7v.23}. 

By induction on $i \le \alpha$ we can find $M'_i$ which is 
$\le_{\frak s}$-increasing continuous such that $M'_i \cap 
N_1 = M_{0,i},M_{0,i}
\le_{\frak s} M'_i$ and if $i = j+1$ then $M'_i$ is brimmed over
$M'_j$ and WNF$_{\frak s}(M_{0,i},M_{0,i+1},M'_i,M'_{i+1})$.

So by \scite{838-7v.19}(2) and see \scite{838-7v.27} we can find a ${\frak
u}_{\frak s}$-free $(0,\xi_{\frak s} \alpha)$-rectangle $\bold d$
such that $M^{\bold d}_{\xi_{\frak s} i} = M'_i$ for $i \le 
\xi_{\frak s} \alpha$
and $(M^{\bold d}_{\varepsilon,0},M^{\bold d}_{\varepsilon
+1,0},a^{\bold d}_{\varepsilon,0}) \in K^{3,\text{up}}_{\frak s}$ for
$\varepsilon < \xi_{\frak s} \alpha$.  Recalling $M_{0,\alpha} = N_0$,
without loss of generality  $M^{\bold d}_{0,\xi_{\frak s} \alpha} \cap N_1 = N_0$ so by 
\scite{838-7v.19}(1) we can find $N^*_1$ such that WNF$_{\frak s}(N_0,N_1,
M^{\bold d}_{0,\xi_{\frak s} \alpha},N^*_1)$.  Recalling that we are assuming
WNF$_{\frak s}(M_0,N_0,M_1,N_1)$, by symmetry, i.e. \scite{838-7v.27} we
have WNF$_{\frak s}(M_0,M_1,N_0,N_1)$ hence  
by transitivity, i.e. Claim \scite{838-7v.15} we
can deduce that WNF$_{\frak s}(M_0,M_1,M^{\bold d}_{0,\xi_{\frak s}
\alpha},N^*_1)$ hence by \scite{838-7v.27}, i.e. 
symmetry WNF$_{\frak s}(M_0,M^{\bold
d}_{0,\xi_{\frak s} \alpha},M_1,N^*_1)$.  
By the definition of WNF$_{\frak s}$, we can find a
${\frak u}_{\frak s}$-free $(1,\xi_{\frak s} \alpha +1)$-rectangle $\bold
d^+$ such that $\bold d^+ \restriction (0,\xi_{\frak s} \alpha) =
\bold d$ and $N^*_1 \le_{\frak s} M^{\bold d^+}_{1,\xi_{\frak s}\alpha
+1}$ and $M_1 = M^{\bold d^+}_{1,0}$.

By the weak existence claim \scite{838-7v.17},
we have
WNF$_{\frak s}(M^{\bold d^+}_{0,\varepsilon},M^{\bold
d^+}_{1,\varepsilon},M^{\bold d^+}_{0,\varepsilon +1},
M^{\bold d^+}_{1,\varepsilon +1})$
for each $\varepsilon < \xi_{\frak s} \alpha$.

Let $M_{1,\alpha +1} = M^{\bold
d^+}_{1,\xi_{\frak s} \alpha +1}$ and for $i \le \alpha$ let
$M_{1,i} := M^{\bold d^+}_{1,\xi_{\frak s} i}$.  So clearly $\langle M_{1,i}:i
\le \alpha +1\rangle$ is $\le_{\frak s}$-increasing continuous, $M_{1,0}
= M_1,N_1 \le_{\frak s} M_{1,\alpha +1}$.  Now to finish the proof 
we need to show, for $i <
\alpha$ that WNF$_{\frak s}(M_{0,i},M_{0,i+1},M_{1,i},M_{1,i+1})$.

For each $i < \alpha$ by the long transitivity claim,
i.e. \scite{838-7v.15} applied to $\langle M^{\bold d^+}_{0,\xi_{\frak
s} i + \varepsilon}:\varepsilon \le \xi_{\frak s}\rangle$ and $\langle
M^{\bold d^+}_{1,\xi_{\frak s}, i +\varepsilon}:\varepsilon \le
\xi_{\frak s}\rangle$ we have
WNF$_{\frak s}(M^{\bold d^+}_{0,\xi_{\frak s} i},
M^{\bold d^+}_{1,\xi_{\frak s} i},M^{\bold d^+}_{0,\xi_{\frak s}(i+1)},
M^{\bold d^+}_{1,\xi_{\frak s}(i+1)})$, by symmetry we have
WNF$_{\frak s}(M^{\bold d^+}_{0,\xi_{\frak s} i}$,
$M^{\bold d^+}_{0,\xi_{\frak s}(i+1)},
M^{\bold d^+}_{1,\xi_{\frak s}i},M^{\bold d^+}_{1,\xi_{\frak
s}(i+1)})$ which means WNF$_{\frak s}(M'_i,M'_{i+1},M_{1,i},M_{i,1+i})$.

Recall that for each $i < \alpha$ we have WNF$_{\frak
s}(M_{0,i},M_{0,i+1},M'_i,M'_{i+1})$.  By the transitivity claim
\scite{838-7v.15}, the two previous sentences imply WNF$_{\frak s}
(M_{0,i},M_{0,i+1},M_{1,i},M_{1,i+1})$ as required.
  \hfill$\square_{\scite{838-7v.31}}$
\enddemo
\bigskip

\proclaim{\stag{838-7v.33} Theorem}  1) {\rm WNF}$_{\frak s}$ is a 
weak non-forking relation on ${\frak K}_\lambda$ respecting 
${\frak s}$ and having disjointness (see Definition \scite{838-7v.35} below).
\nl
2) If {\rm NF} is a weak non-forking relation of ${\frak K}_\lambda$ 
respecting ${\frak s}$, \ub{then} $\text{\rm NF} \subseteq$ 
{\rm WNF}$_{\frak s}$.
\endproclaim
\bn
A relative of Definition \yCITE[nf.0X]{600} is:
\definition{\stag{838-7v.35} Definition}  Let ${\frak K}_\lambda$ be a
$\lambda$-a.e.c.
\nl
1) We say that NF is a weak non-forking relation on ${\frak
K}_\lambda$ \ub{when}: it satisfied the requirements in Definition
\yCITE[nf.0X]{600}(1); 
\ub{except that} we replace uniqueness (clause (g) there) by 
weak uniqueness meaning that the conclusion of Claim
\scite{838-7v.31} holds (replacing WNF$_{\frak s}$ by {\rm NF});
or see the proof of \scite{838-7v.33} below for a list or use Definition
\yCITE[nf.20.9]{600}(1).  
\nl
2) Let ${\frak t}$ be an almost good $\lambda$-frame and ${\frak
K}_\lambda = {\frak K}_{\frak t}$ and NF be a weak non-forking relation
on ${\frak K}_\lambda$.  We say that NF respects ${\frak t}$ \ub{when}: if
NF$(M_0,N_0,M_1,N_1)$ and $(M_0,N_0,a) \in K^{3,\text{bs}}_{\frak t}$
\ub{then} $\ortp_{\frak t}(a,M_1,N_1)$ does not fork over $M_0$.  We
say NF is a weak ${\frak t}$-non-forking relation when it is a weak
${\frak t}$-non-forking relation respecting ${\frak t}$. 
\nl
3) In part (1) we say NF has disjointness when WNF$(M_0,N_0,M_1,N_1)
   \Rightarrow M_0 \cap M_1 = M_0$.
\nl
4) We say NF is a pseudo non-forking relation on ${\frak K}_\lambda$
   \ub{when} we have clauses (a)-(f) of Definition \yCITE[nf.0X]{600} or see
   the proof below.  Also here parts (2),(3) are meaningful.
\enddefinition
\bigskip

\demo{Proof of \scite{838-7v.33}}  1) Let us list the conditions on NF :=
WNF$_{\frak s}$ being a weak non-forking relation let ${\frak
K}_\lambda = {\frak K}_{\frak s}$.  We shall use \scite{838-7v.29} freely.
\mn
\ub{Condition (a)}:  NF is a 4-place relation on ${\frak K}_\lambda$.
\nl
[Why?  This holds by Definition \scite{838-7v.5}(1),(2).]
\mn
\ub{Condition (b)}:  NF$(M_0,M_1,M_2,M_3)$ implies $M_0 
\le_{{\frak K}_\lambda} M_\ell \le_{{\frak K}_\lambda} M_3$ for
$\ell=1,0$ and NF is preserved by isomorphisms.
\nl
[Why?  The preservation by isomorphisms holds by the definition, and
also the order demands.]
\mn
\ub{Condition (c)$_1$}:  [Monotonicity]  If NF$(M_0,M_1,M_2,M_3)$ and
$M_0 \le_{{\frak K}_\lambda} M'_\ell \le_{{\frak K}_\lambda} M_\ell$ for
$\ell=1,2$ then NF$(M_0,M'_1,M'_2,M_3)$.
\nl
[Why?  By \scite{838-7v.11}.]
\mn
\ub{Condition (c)$_2$}:  [Monotonicity]  If NF$(M_0,M_1,M_2,M_3)$ and
$M_3 \le_{{\frak K}_\lambda} M'_3$ and $M_1 \cup M_2
\subseteq M''_3 \le_{{\frak K}_\lambda} M_3$ \ub{then} 
NF$(M_0,M_1,M_2,M''_3)$.
\nl
[Why?  By Claim \scite{838-7v.11}.]
\mn
\ub{Condition (d)}:  [Symmetry]  If NF$(M_0,M_1,M_2,M_3)$ then
NF$(M_0,M_2,M_1,M_3)$.
\nl
[Why?  By Claim \scite{838-7v.27}.]
\mn
\ub{Condition (e)}:  [Long Transitivity]  If $\alpha < \lambda^+$, 
NF$(M_i,N_i,M_{i+1},N_{i+1})$ for each $i < \alpha$ and
$\langle M_i:i \le \alpha\rangle,\langle N_i:i \le \lambda\rangle$ are
$\le_{{\frak K}_\lambda}$-increasing continuous sequences \ub{then}
NF$(M_0,N_0,M_\alpha,N_\alpha)$.
\nl
[Why?  By Claim \scite{838-7v.15}.]
\mn
\ub{Condition (f)}:  [Existence]  Assume $M_0 \le_{{\frak K}_\lambda}
M_\ell$ for $\ell=1,2$.
Then for some $M_3,f_1,f_2$ we have $M_0 \le_{{\frak K}_\lambda} M_3
\in {\frak K}_\lambda,f_\ell$ is a $\le_{{\frak K}_\lambda}$-embedding
for $M_\ell$ into $M_3$ over $M_0$ for $\ell=1,2$ and
NF$(M_0,f_1(M_1),f_2(M_2),M_3)$.  Here we have the disjoint version,
i.e. $f_2(M_1) \cap f_2(M_2) = M_2$.
\nl
[Why?  By Lemma \scite{838-7v.19}(1).]
\mn
\ub{Condition (g)}: Lifting or weak uniqueness [a replacement for uniqueness]

This is the content of \scite{838-7v.31}.

Thus we have finished presenting the definition of ``NF is a weak non-forking
relation on ${\frak K}_{\frak s}$" and proving that WNF$_{\frak s}$
satisfies those demands.

But we still owe ``WNF$_{\frak s}$ respect ${\frak s}$" where 
\ub{NF respect ${\frak s}$} 
means that if NF$(M_0,N_0,M_1,N_1)$ and $(M_0,N_0,a) \in K^{3,{\frak
s}}_{\frak s}$ then $\ortp_{\frak s}(a,M_1,N_1)$ does not fork over $M_0$,
i.e. $(M_0,N_0,a) \le_{\text{bs}} (M_1,N_1,a) \in
K^{3,\text{bs}}_{\frak s}$.
\nl
[Why?  This holds by Observation \scite{838-7v.9}(1).]

Also the disjointness of WNF is easy; use \scite{838-7v.13} and categoricity.
\enddemo
\bigskip

\demo{Proof of \scite{838-7v.33}(2)}

So assume NF$(M_0,N_0,M_1,N_1)$ and we should prove WNF$_{\frak
s}(M_0,N_0,M_1,N_1)$ so let $\bold d$ be as in Definition
\scite{838-7v.5}(1).  As NF satisfies existence, transitivity and
monotonicity without loss of generality  it suffices to deal with the case $M^{\bold
d}_{0,\alpha(\bold d)} = N_0$.

This holds by the definition of WNF$_{\frak s}$ in \scite{838-7v.5} and
clause (g) in the definition of being weak non-forking relation and
``respecting ${\frak s}$".   \hfill$\square_{\scite{838-7v.33}}$ 
\enddemo
\bigskip

At last we can get rid of the ``almost" in ``almost good
$\lambda$-frame", of course, this is under the Hypothesis
\scite{838-7v.1}, otherwise we do not know.
\bigskip

\proclaim{\stag{838-7v.37} Lemma}  1) If ${\frak t}$ is an almost good
$\lambda$-frame and {\rm WNF} is a weak non-forking relation on ${\frak
K}_\lambda$ respecting ${\frak t}$ \ub{then} ${\frak t}$ is a good
$\lambda$-frame.
\nl
2) In \scite{838-7v.1}, ${\frak s}$ is a good $\lambda$-frame.
\endproclaim
\bigskip

\demo{Proof}  Part (2) follows from part (1) and \scite{838-7v.33}(1)
above.  So recalling Definition \scite{838-5t.3} we should just 
prove that ${\frak t}$ satisfies Ax(E)(c).  Note that the Hypothesis
\scite{838-7v.38} below holds hence we are allowed to use \scite{838-7v.61}.

Let $\langle M_i:i \le \delta\rangle$ be $\le_{\frak s}$-increasing
continuous and $p \in {\Cal S}^{\text{bs}}_{\frak s}(M_\delta)$ and we
should prove that $p$ does not fork over $M_i$ for some $i < \delta$.
By renaming without loss of generality  $\delta < \lambda^+$ and $\delta$ is divisible
by $\lambda^2 \omega$ and $\varepsilon \le \lambda \wedge i < \delta
\Rightarrow M_{i+1} = M_{i+1 + \varepsilon}$.  
Let ${\frak u}$ be as in Definition
\scite{838-7v.3}, so ${\frak u}$ is 
a nice construction framework.  Let $\alpha = \delta,\beta = \delta$.

Now
\mr
\item "{$(*)_1$}"  there is $\bold d$ such that
{\roster
\itemitem{ $(\alpha)$ }  $\bold d$ is a ${\frak u}$-free
$(\alpha,\beta)$-rectangle
\sn
\itemitem{ $(\beta)$ }  $\bold d$ is a strictly brimmed, see
Definition \scite{838-5t.29}(2)
\sn
\itemitem{ $(\gamma)$ }  if $i < \alpha$ and $j < \delta$ then
WNF$(M^{\bold d}_{i,j},M^{\bold d}_{i,j+1},M^{\bold d}_{i+1,j},
M^{\bold d}_{i+1,j+1})$
\sn
\itemitem{ $(\delta)$ }  $M^{\bold d}_{i,0} = M_i$ for $i \le \delta$
\sn
\itemitem{ $(\varepsilon)$ }  $\bold J^{\bold d}_{i,j} = \emptyset$ when $i
< \delta,j \le \delta$ and $\bold I^{\bold d}_{i,j} = \emptyset$ when $i \le
\delta,j < \delta$.
\endroster}
\ermn
This will be done as in the proof of Observation \scite{838-5t.33}.

By the properties of WNF (i.e. using twice long transitivity and symmetry)
\mr
\item "{$(*)_2$}"  if $i_1 < i_2 \le \alpha$ and $j_1 < j_2 \le \beta$
\ub{then} WNF$(M^{\bold d}_{i_1,j_1},M^{\bold d}_{i_1,j_2},M^{\bold
d}_{i_2,j_1},M^{\bold d}_{i_2,j_2})$ and
WNF$(M^{\bold d}_{i_1,j_1},M^{\bold d}_{i_2,j_1},M^{\bold d}_{i_1,j_2},
M^{\bold d}_{i_2,j_2})$.
\ermn
Now we can choose a ${\frak u}$-free $(\alpha_{\bold d},\beta_{\bold
d})$-rectangle $\bold e$ such that
\mr
\item "{$(*)_3$}"  $(a) \quad M^{\bold e}_{i,j} = M^{\bold d}_{i,j}$
for $i \le \alpha_{\bold d},j \le \beta_{\bold d}$
\sn
\item "{${{}}$}"  $(b) \quad$ if $i < \alpha,j < \beta$ and $p \in
{\Cal S}^{\text{bs}}_{\frak s}(M^{\bold d}_{i,j})$ then 
for $\lambda$ ordinals $\varepsilon < \lambda$ we have:
{\roster
\itemitem{ ${{}}$ }   $(\alpha) \quad \ortp_{\frak t}(b^{\bold d}_{i +
\varepsilon,j},M^{\bold d}_{i + \varepsilon,j+1},M^{\bold d}_{i +
\varepsilon +1,j+1})$ is a non-forking extension of 
\nl

\hskip25pt $p$, recalling 
$\bold J^{\bold e}_{i + \varepsilon,j} = \{b^{\bold e}_{i + \varepsilon,j}\}$, 
\sn
\itemitem{ ${{}}$ }  $(\beta) \quad 
\ortp_{\frak t}(a^{\bold e}_{i,j +\varepsilon},
M^{\bold d}_{i + 1,j +\varepsilon},M^{\bold d}_{i + 1,j
+ \varepsilon +1})$ is a non-forking
\nl

\hskip25pt  extension of $p$, recalling
$\bold I^{\bold e}_{i,j + \varepsilon} = \{a^{\bold e}_{i + 1,j
+\varepsilon}\}$.
\endroster}
\ermn
[Why?  We can choose $\bold e \rest (\alpha,\alpha)$ by induction on
$\alpha \le \delta$.  The non-forking condition in the definition of
${\frak u}$-free holds because WNF respects ${\frak t}$ and $(*)_2$.]

So $\bold d$ is full (see Definition \scite{838-5t.29}(3),(3A), even
strongly full) hence by Claim \scite{838-5t.36}(3A)(c)
\mr
\item "{$(*)_4$}"   $M^{\bold d}_{\alpha,\beta}$ is
brimmed over $M^{\bold d}_{\alpha,0} = M_\delta$.
\ermn
Hence $p \in {\Cal S}^{\text{bs}}_{\frak s}(M_\delta) = {\Cal
S}^{\text{bs}}_{\frak t}(M^{\bold d}_{\alpha,0})$ is realized in
$M_{\alpha,\beta}$ say by $c \in M_{\alpha,\beta}$, so for some $i <
\alpha$ we have $c \in M_{i,\beta}$.  As 
WNF$(M^{\bold d}_{i,0},M^{\bold d}_{i,\beta},M^{\bold d}_{\alpha,0},  
M^{\bold d}_{\alpha,\beta})$ holds by $(*)_2$ above  
by Claim \scite{838-7v.61} below, it follows that
$\ortp_{\frak t}(c,M^{\bold d}_{\alpha,0},M^{\bold d}_{\alpha,\beta})$
does not fork over $M^{\bold d}_{i,0}$ which means $i < \delta$ and
$p$ does not fork over $M_i$ as required.  \hfill$\square_{\scite{838-7v.37}}$
\enddemo
\bn
Now for the rest of the section we replace Hypothesis \scite{838-7v.1} by
\demo{\stag{838-7v.38} Hypothesis}  Assume ${\frak s}$ is an 
almost good $\lambda$-frame categorical in $\lambda$
and WNF is a weak non-forking 
relation on ${\frak K}_{\frak s}$ respecting ${\frak s}$.
\enddemo
\bn
The following is related to the proof of  \scite{838-7v.37}.
\definition{\stag{838-7v.41} Definition}  1) We say $\bold d = \langle
  M_{\alpha,\beta}:\alpha \le \alpha_{\bold d},\beta \le \beta_{\bold
  d}\rangle$ is a WNF-free rectangle (or $(\alpha_{\bold
  d},\beta_{\bold d})$-rectangle) when:
\mr
\item "{$(a)$}"  $\langle M_{\alpha,j}:j \le \beta_{\bold d}\rangle$
  is $\le_{\frak s}$-increasing continuous for each $\alpha \le
  \alpha_{\bold d}$
\sn
\item "{$(b)$}"  $\langle M_{i,\beta}:i \le \alpha_{\bold d}\rangle$
 is $\le_{\frak s}$-increasing continuous for each $\beta \le \beta_{\bold d}$
\sn
\item "{$(c)$}"  WNF$(M_{i,j},M_{i+1,j},M_{i,j+1},M_{i+1,j+1})$ for $i
< \alpha,j < \beta$.
\ermn
2) Let $\bar \alpha = \langle \alpha_j:j \le \beta\rangle$ be $\le$-increasing.

We say $\bold d = \langle M_{i,j}:i \le \alpha_j$ and $j \le
\beta\rangle$ is a WNF-free $(\langle \alpha_j:j \le
\beta\rangle,\beta)$-triangle \ub{when}:
\mr
\item "{$(a)$}"  $\langle M_{i,j}:i \le \alpha_j\rangle$
 is $\le_{\frak s}$-increasing continuous for each $j \le \beta$
\sn
\item "{$(b)$}"  $\langle M_{i,j}:j \le \beta$ satisfies $i \le
\alpha_j\rangle$  is $\le_{\frak s}$-increasing continuous for each $i
< \alpha_\beta$
\sn
\item "{$(c)$}"    WNF$(M_{i,j},M_{i+1,j},M_{i,j+1},M_{i+1,j+1})$ for
$j < \beta,i < \alpha_j$.
\endroster
\enddefinition
\bn
Now we may note that 
some facts proved in \xCITE{705} for weakly successful good
$\lambda$-frame can be proved under Hypothesis \scite{838-7v.1} or just
\scite{838-7v.38}.  Systematically see \cite{Sh:842}.
\proclaim{\stag{838-7v.45} Claim}  $M^{\bold d}_{\alpha_\beta,\beta}$ is
brimmed over $M_{0,\beta}$ \ub{when}:
\mr
\item "{$(a)$}"  $\bar \alpha = \langle \alpha_j:\alpha \le \beta\rangle$ is 
increasing continuous
\sn
\item "{$(b)$}"  $\bold d$ is a {\rm WNF}-free
$(\bar\alpha,\beta)$-triangle 
\sn
\item "{$(c)$}"  $M_{i+1,j+1}$ is $\le_{\frak s}$-universal
over $M_{i,j}$ when $i < \alpha_j,j < \beta$.
\endroster
\endproclaim
\bigskip

\demo{Proof}  If $\beta $ and each $\alpha_j$ is divisible by
$\lambda$ we can repeat (part of the) proof of \scite{838-7v.37} and this
suffices for proving \scite{838-7v.61} hence for proving \scite{838-7v.37} (no
vicious circle!)

In general we have to repeat the proof of \scite{838-5t.36} \ub{or} first
find a WNF-free $(\langle \lambda \alpha_j:j \le \beta
\rangle,\beta)$-rectangle $\bold d'$ which is brimmed and full and
$M^{\bold d'}_{\lambda,j} = M^{\bold d^+}_{i,j}$ for $i \le
\alpha_j,j \le \beta$ and then use the first sentence.
\hfill$\square_{\scite{838-7v.45}}$ 
\enddemo
\bigskip

It is natural to replace $\le_{\text{bs}}$ by the stronger
$\le_{\text{wnf}}$ defined below (and used later).  
\definition{\stag{838-7v.49} Definition}  1) Let $\le_{\text{wnf}}$ be the
following two-place relation on $K^{2,\text{bs}}_{\frak s} :=
\{(M,N):M \le_{\frak s} N$ are from  ${\frak K}_{\frak s}\}$, we have
$(M_0,N_0) \le_{\text{wnf}} (M_1,N_1)$ iff:
\mr
\item "{$(a)$}"  $(M_\ell,N_\ell) \in K^{2,\text{bs}}_{\frak s}$ for $\ell=0,1$
\sn
\item "{$(b)$}"  WNF$(M_0,N_0,M_1,N_1)$.
\ermn
2) Let $(M_1,N_1,a) \le_{\text{wnf}} (M_2,N_2,a)$ means
$(M_\ell,N_\ell,a) \in K^{3,\text{bs}}_{\frak s}$ for $\ell=1,2$
and $(M_1,N_1) \le_{\text{wnf}} (M_2,N_2)$.
\enddefinition
\bigskip

\proclaim{\stag{838-7v.51} Claim}  1) $\le_{\text{wnf}}$ is a partial order
on $K^{2,\text{bs}}_{\frak s}$.
\nl
2) If $\langle (M_\alpha,N_\alpha):\alpha < \delta\rangle$ is
$\le_{\text{wnf}}$-increasing continuous \ub{then} $\alpha < \delta
\Rightarrow (M_\alpha,N_\alpha) \le_{\text{wnf}} (\dbcu_{\beta <
\delta} M_\beta,\dbcu_{\beta < \delta} N_\beta) \in
K^{3,\text{bs}}_{\frak s}$.
\nl
3) If $(M_1,N_1) \le_{\text{wnf}} (M_2,N_2)$ and $(M_1,N_1,a) \in
{\frak K}^{3,\text{bs}}_{\frak s}$ \ub{then} $(M_1,N_1,a) \le_{\text{bs}}
(M_2,N_2,a) \in K^{3,\text{bs}}_{\frak s}$.
\nl
4) If $M \le_{\frak s} N' \le_{\frak s} N$ and $(M,N) \in
K^{2,\text{bs}}_{\frak s}$ \ub{then} $(M,N') \le_{\text{wnf}} (M,N)$.
\nl
5) If $(M_1,N_1) \le_{\text{wnf}} (M_2,N_2)$ and $M_1 \le_{\frak s}
   N'_1 \le_{\frak s} N_1,M_1 \le_{\frak s} M'_2 \le_{\frak s} M_2$
   and $N'_1 \cup M'_2 \subseteq N'_2 \le_{\frak s} N''_2$ and $N_2
   \le_{\frak s} N''_2$ \ub{then} $(M_1,N'_1) \le_{\text{wnf}}
   (M'_2,N'_2)$.
\nl
6) Similarly to (1),(2),(4),(5) for $(K^{3,\text{bs}}_{\frak s},
\le_{\text{wnf}})$. 
\endproclaim
\bigskip

\demo{Proof}  Easy by now.  \hfill$\square_{\scite{838-7v.51}}$
\enddemo
\bn
The following is a ``downward" version of ``WNF respect ${\frak s}$"
which was used in the proof of \scite{838-7v.37}.
\proclaim{\stag{838-7v.61} Claim}  If {\rm WNF}$(M_0,N_0,M_1,N_1)$ and $c
\in N_0$ and
$(M_1,N_1,a) \in K^{3,\text{bs}}_{\frak s}$ \ub{then} $(M_0,N_0,c) \in
K^{3,\text{bs}}_{\frak s}$ and $\ortp_{\frak s}(c,M_1,N_1)$ does not
fork over $M_0$. 
\endproclaim
\bigskip

\demo{Proof}  The second phrase in the conclusion (about
``$\ortp_{\frak s}(a,M_1,N_1)$ does not fork over $M_0$") follows from
the first by ``WNF respects ${\frak t}$"; so if ${\frak s}$ is type
full (i.e. ${\Cal S}^{\text{bs}}_{\frak s}(M) = {\Cal
S}^{\text{na}}_{\frak s}(M)$) this is easy.  (So if ${\frak s}$ is
type-full the result is easy.)

Toward contradiction assume
\mr
\item "{$(*)_0$}"  $(M_0,N_0,M_1,N_1,c)$ form a counterexample.
\ermn
Also by monotonicity properties without loss of generality :
\mr
\item "{$(*)_1$}"  $(a) \quad N_0$ is brimmed over $M_0$
\sn
\item "{${{}}$}"  $(b) \quad M_1$ is brimmed over $M_0$
\sn
\item "{${{}}$}"  $(c) \quad N_1$ is brimmed over $N_0 \cup M_1$.
\ermn
Let $\langle {\Cal U}_\varepsilon:\varepsilon < \lambda_{\frak s}\rangle$ be
an increasing sequence of subsets of $\lambda_{\frak s}$ such that
$|{\Cal U}_0| = |{\Cal U}_{\varepsilon +1} \backslash {\Cal
U}_\varepsilon| = \lambda_{\frak s}$.  We now by induction on
$\varepsilon \le \lambda$ choose $\bold d_\varepsilon,
\bar{\bold a}_\varepsilon$ such that $\bold d_\varepsilon = \langle
M_{i,j}:i \le \lambda j,j \le \varepsilon\rangle$ and:
\mr
\item "{$\circledast$}"  $(a) \quad \bold d_\varepsilon$ is a WNF-free
triangle
\sn
\item "{${{}}$}"  $(b) \quad M_{i+1,j+1}$ is brimmed over $M_{i+1,j}
\cup M_{i,j+1}$ when $i < \lambda j \wedge j < \varepsilon$
\sn
\item "{${{}}$}"  $(c) \quad M_{0,j+1}$ is brimmed over $M_{0,j}$
\sn
\item "{${{}}$}"  $(d) \quad \bar{\bold a}_\varepsilon = \langle
a_\alpha:\alpha \in {\Cal U}_\varepsilon\rangle$ list the elements of
$M_{\varepsilon \lambda,\varepsilon}$
\sn
\item "{${{}}$}"  $(e) \quad$ if $\varepsilon = j+1$ and ${\Cal W}_j
\subseteq {\Cal U}_j$ defined below is $\ne \emptyset$ and $\gamma_j =
\text{ min}({\Cal W}_j)$ 
\nl

\hskip25pt then $\ortp_{\frak s}(a_{\gamma_j},M_{0,\varepsilon}, 
M_{\lambda j,\varepsilon}) \in
{\Cal S}^{\text{bs}}_{\frak s}(M_{0,\varepsilon})$ 
and $[\ortp_{\frak s}(a_{\gamma_j},M_{0,j},M_{\lambda j,j}) \in$
\nl

\hskip25pt ${\Cal S}^{\text{bs}}_{\frak s}(M_{0,j}) \Rightarrow \ortp_{\frak
s}(a_{\gamma_j},M_{0,\varepsilon},M_{\lambda j,\varepsilon})$ 
forks over $M_{0,j}]$ where 
\nl

\hskip25pt ${\Cal W}_j = \{\alpha \in {\Cal U}_i$: we can
find $M',N'$ such that WNF$(M_{0,j},M_{\lambda j,j},M',N')$ 
\nl

\hskip25pt  and $\ortp_{\frak s}(a_\alpha,M',N') \in
{\Cal S}^{\text{bs}}_{\frak s}(M')$ but is not a
non-forking extension 
\nl

\hskip25pt  of $\ortp_{\frak s}(a,M_{0,j},M_{\lambda j,j})$,
e.g. $\ortp_{\frak s}(a,M_{0,j},M_{\lambda j,j}) \notin {\Cal
S}^{\text{bs}}_{\frak s}(M_{0,j})\}$.
\ermn
There is no problem to carry the definition.

Now by \scite{838-7v.45}
\mr
\item "{$(*)_2$}"  $M_{\lambda \lambda,\lambda}$ is brimmed over
$M_{0,\lambda}$.
\ermn
So by $(*)_1 + (*)_2$ and ${\frak s}$ being categorical without loss of generality 
\mr
\item "{$(*)_3$}"  $(M_0,N_0) = (M_{0,\lambda},M_{\lambda
\lambda,\lambda})$.
\ermn
Clearly $M_{\lambda \lambda,\lambda}$ is the union of the $\le_{\frak
s}$-increasing continuous chain $\langle M_{\lambda j,j}:j \le \lambda\rangle$
hence $j_1(*)$ is well defined where: 
\mr
\item "{$(*)_4$}"  $j_1(*) = \text{ min}\{j < \lambda:c \in M_{\lambda j,j}\}$.
\ermn
By clause (d) of $\circledast$ for some $j(*)$ we have
\mr
\item "{$(*)_5$}"  $j(*) \in [j_1(*),\lambda)$ and $c \in
\{a_\alpha:\alpha \in {\Cal U}_{j(*)}\}$.
\ermn
So by the choice of $j(*),\gamma(*)$ is well defined where
\mr
\item "{$(*)_6$}"  $\gamma(*) = \text{ min}\{\gamma \in {\Cal
U}_{j(*)}:a_\gamma = c\} < \lambda$.
\ermn
Note that
\mr
\item "{$(*)_7$}"  $(M_{0,j},M_{\lambda j,j}) \le_{\text{wnf}}
(M_{0,\lambda},M_{\lambda \lambda,\lambda}) = (M_0,N_0)
\le_{\text{wnf}} (M_1,N_1)$ for $j \in [j(*),\lambda)$.
\ermn
Also
\mr
\item "{$(*)_8$}"  if $j \in [j(*),\lambda)$ then $\ortp_{\frak
s}(c,M_{0,j},M_{\lambda j,j}) \notin {\Cal S}^{\text{bs}}_{\frak
 s}(M_{0,j})$.
\ermn
[Why?  By $(*)_7$, we have WNF$_{\frak s}(M_{0,j},M_{\lambda
j,j},M_{0,\lambda},M_{\lambda \lambda,\lambda})$ hence if
$(\ortp_{\frak s}(c,M_{0,j},M_{\lambda j,j}) \in {\Cal
S}^{\text{bs}}_{\frak s}(M_{0,j})$ then $(M_{0,j},M_{\lambda j,j},c)
\in K^{3,\text{bs}}_{\frak s}$ but by $(*)_7$ we have
WNF$(M_{0,j},M_{\lambda j,j},M_1,N_1)$ and WNF respects ${\frak s}$
hence $\ortp_{\frak s}(c,M_1,N_1)$ does not 
fork over $M_{0,j}$, so by monotonicity it
does not fork over $M_{0,\lambda} = M_0$ and this contradicts $(*)_0$.]
\mr
\item "{$(*)_9$}"  if $j \in [j(*),\lambda)$ then min$({\Cal W}_j) \le
\gamma(*)$. 
\ermn
[Why?  As $\gamma(*) \in {\Cal W}_j$, i.e. satisfies the requirement
which appear in clause (e) of $\circledast$ that is $(M_1,N_1)$ here
can stand for $(M',N')$ there.]

So by cardinality considerations for some $j_1 < j_2$ from
$[j(*),\lambda)$ we have min$({\Cal W}_{j_1}) = \text{ min}({\Cal W}_{j_2})$
but this gives a contradiction as in the proof of $(*)_8$.
\hfill$\square_{\scite{838-7v.61}}$ 
\enddemo
\bn
\margintag{7v.62}\ub{\stag{838-7v.62} Exercise}:  Show that in Hypothesis \scite{838-7v.38} we
can omit ``${\frak s}$ is categorical (in $\lambda$)".
\sn
[\ub{Hint}:  The only place it is used is in showing $(*)_3$ during
the proof of \scite{838-7v.61}.  To avoid it in $\circledast$ there waive
clause (c) and add $j \le \lambda \Rightarrow M_{0,j} = M_0$.]
\bn
\margintag{7v.64A}\ub{\stag{838-7v.64A} Exercise}:  Show that in this section we can replace
clause (c) of Hypothesis \scite{838-7v.1}, i.e. ``${\frak s}$ is categorical in
$\lambda$" by
\mr
\item "{$(c)'$}"  $\dot I(K_{\frak s}) \le \lambda$
\nl
or just
\sn
\item "{$(c)''$}"   $\xi_{\frak s} = \sup\{\xi_M:M \in K_{\frak s}\} <
\lambda^+$ where for $M \in K_{\frak s}$ we let $\xi_N = \text{
min}\{\alpha_{\bold d}$: there is a ${\frak u}$-free $(\alpha_{\bold
d},0)$-rectangle such that $M^{\bold d}_{0,0} = M,(M^{\bold
d}_{i,0},M^{\bold d}_{i,0},a^{\bold d}_{i,0}) \in
K^{3,\text{up}}_{\frak s}$ and $M^{\bold d}_{\alpha(\bold d),0}$ is
universal over $M^{\bold d}_{0,0} = M\}$.
\ermn
[Hint:  The only place we use ``${\frak s}$ is
categorical in $\lambda$" is in claim \scite{838-7v.27} more fully in
\scite{838-7v.19}(2) we would like to have a bound $< \lambda^+$ on
$\alpha_{\bold d}$ not depending on $M_0$, see \scite{838-7v.21}
(and then quoting it).  It is used to define 
$\xi_{\frak s}$.  As $(c)' \Rightarrow (c)''$ without loss of generality  we assume $(c)''$.]
\bn
\margintag{7v.67}\ub{\stag{838-7v.67} Exercise}:  (Brimmed lifting, compare with 
\yCITE[1.15K]{705}(4).)  
\nl
1) For any $\le_{\frak s}$-increasing 
continuous sequence $\langle M_\alpha:\alpha
\le \alpha(*)\rangle$ with $\alpha(*) < \lambda^+_{\frak s}$ we can
find $\bar N$ such that
\mr
\item "{$\circledast$}"  $(a) \quad \bar N = \langle N_\alpha:\alpha
\le \alpha(*)\rangle$ is $\le_{\frak s}$-increasing continuous
\sn
\item "{${{}}$}"  $(b) \quad$ {\rm WNF}$(M_\alpha,N_\alpha,
M_\beta,N_\beta)$ for $\alpha < \beta \le
\alpha(*)$
\sn
\item "{${{}}$}"  $(c) \quad N_\alpha$ is brimmed over $M_\alpha$ for
$\alpha = 0$ and moreover for every $\alpha \le \alpha(*)$
\sn
\item "{${{}}$}"  $(d) \quad N_{\alpha +1}$ is brimmed over $M_{\alpha
+1} \cup N_\alpha$ for $\alpha < \alpha(*)$.
\ermn
2) In fact the moreover in clause (c) follows from (a),(b),(d); an addition
\mr
\item "{$(c)^+$}"  $N_\beta$ is brimmed over $N_\alpha \cup M_\beta$
for $\alpha < \beta \le \alpha(*)$.
\endroster
\bn
[Hint:   Similar to \scite{838-5t.36} or \scite{838-5t.37} or note that by
\scite{838-7v.37} we can quote \xCITE{705}.]  
\goodbreak

\head {\S8 Density of $K^{3,\text{uq}}_{\frak s}$ for good
 $\lambda$-frames} \endhead  \resetall 
 \spuriousreset
\bigskip

We shall prove non-structure from failure of density for
$K^{3,\text{uq}}_{\frak s}$ in two rounds.  First, in \scite{838-8h.12} -
\scite{838-8h.26} we prove the wnf-delayed version.  Second, in \scite{838-8h.37} -
\scite{838-p.4.13} - we use its conclusion to prove the general case.  Of
course, by \S6, we can assume as in \S7:
\bigskip

\proclaim{\stag{838-8h.1} Hypothesis}  We assume (after \scite{838-8h.2})
\mr
\item "{$(a)$}"  ${\frak s}$ is an almost \footnote{by \scite{838-7v.33} 
clauses (b) + (c) implies that ${\frak s}$ 
is actually a good $\lambda$-frame but we may ignore this} good $\lambda$-frame
\sn
\item "{$(b)$}"  {\rm WNF} is a weak non-forking
relation on $K_{\frak s}$ respecting ${\frak s}$ with disjointness
\nl
(not necessarily the one from Definition \scite{838-7v.5}, but Hypothesis
\scite{838-7v.38} holds).
\endroster
\endproclaim
\bn
We can justify Hypothesis \scite{838-8h.1} by
\demo{\stag{838-8h.2} Observation}  Instead clause (b) of \scite{838-8h.1} we can
assume
\mr
\item "{$(b)'$}"  ${\frak s}$ has existence for
$K^{3,\text{up}}_{{\frak s},\lambda^+}$ 
(so we may use the consequences of Conclusion \scite{838-7v.33})
\sn
\item "{$(c)$}"  ${\frak s}$ is categorical in $\lambda$.
\endroster
\enddemo
\bigskip

\demo{Proof}  Why?  First note that Hypothesis \scite{838-7v.1} holds:
\mn
\ub{Part (1)} there, ${\frak s}$ is an almost good $\lambda$-frame, is
clause (a) of Hypothesis \scite{838-8h.1}.
\mn
\ub{Part (2)} there is clause (b)$'$ assumed above.
\mn
\ub{Part (3)}, categoricity in $\lambda$, there is clause (c) above.
\mn
\ub{Part (4)} there, disjointness of ${\frak s}$, holds by
\scite{838-5t.45}.

So the results of \S7 holds, in particular the relation WNF :=
WNF$_{\frak s}$ defined in Definition \scite{838-7v.5} is a weak
non-forking relation (on $K_{\frak s}$) respecting ${\frak s}$, by
Claim \scite{838-7v.33}.  
\nl
${{}}$  \hfill$\square_{\scite{838-8h.2}}$
\enddemo
\bn
We now use a relative of ${\frak u}^1_{\frak s}$ from Definition 
\scite{838-e.5I}; this will be
the default value of ${\frak u}$ in this section so $\partial$ will be
$\partial_{\frak u} = \lambda^+$. 
\definition{\stag{838-8h.4} Definition}  For ${\frak s}$ as in
\scite{838-8h.1} we define ${\frak u} = {\frak u}^3_{\frak s}$ as
follows:
\mr
\item "{$(a)$}"    $\partial_{\frak u} = \lambda^+ (= \lambda^+_{\frak s})$
\sn
\item "{$(b)$}"    ${\frak K}_{\frak u} = {\frak K}_{\frak s}$ (or
${\frak K}'_{\frak s}$ see \scite{838-e.5D}, \scite{838-e.5E}; but not
necessary by (c) of \scite{838-8h.1})
\sn
\item "{$(c)$}"    FR$_\ell = \{(M,N,\bold J):M \le_{\frak s} N$ and
$\bold J = \emptyset$ or $\bold J = \{a\}$ and $(M,N,a) \in
K^{3,\text{bs}}_{\frak s}\}$ 
\sn
\item "{$(d)$}"    $\le_{\frak u}^\ell$ is defined by
$(M_0,N_0,\bold J_0) \le_{\frak u}^\ell (M_1,N_1,\bold J_1)$ \ub{when}
{\roster
\itemitem{ $(\alpha)$ }  WNF$(M_0,N_0,M_1,N_1)$
\sn    
\itemitem{ $(\beta)$ }  $\bold J_1 \subseteq \bold J_0$
\sn    
\itemitem{ $(\gamma)$ }  if $\bold J_0 = \{a\}$ then $\bold J_1 =
\{a\}$ hence $(M_0,N_0,a) \le_{\frak s}^{\text{bs}} (M_1,N_1,a)$ by
``WNF respects ${\frak s}$", see Hypothesis \scite{838-8h.1} and
Definition \scite{838-7v.35}; if we use WNF$_{\frak s}$ then we can quote
\scite{838-7v.33}(2), also by \scite{838-7v.9}(1).
\endroster}
\endroster
\enddefinition
\bigskip

\remark{\stag{838-8h.7} Remark}  0) The choice in \scite{838-8h.4} gives us
symmetry, etc., i.e. ${\frak u}$ is self-dual, this sometimes helps.
\nl
1) We could define FR$_1,\le_1$ as above but
\mr   
\item "{$(e)$}"   FR$_2 = \{(M,N,\bold J):M \le_{\frak s} N$ and
$\bold J = \emptyset$ or $\bold J = \{a\},(M,N,a) \in
K^{3,\text{bs}}_{\frak s}\}$
\sn
\item "{$(f)$}"    $\le_2$ is defined by:
\nl
$(M_0,N_0,\bold J_0) \le_1 (M_1,N_1,\bold J_1)$ \ub{when} (both are from
FR$_1$ and)
{\roster
\itemitem{ $(\alpha)$ }  $M_0 \le_{\frak s} M_1$ and $N_0 \le_{\frak s} N_1$
\sn    
\itemitem{ $(\beta)'$ }  $\bold J_0 \subseteq \bold J_1$
\sn
\itemitem{ $(\beta)''$ }  if $\bold J_0 = \{a\}$ then $\bold J_1 =
\{a\}$ and $(M_0,N_0,a) \le_{\frak s}^{\text{bs}} (M_1,N_1,a)$.
\endroster}
\ermn
2) We call it ${\frak u}^{3,*}_{\frak s}$.  However, then for proving
\scite{838-8h.37} we have to use ${\frak u}^{*,3}_{\frak s}$ which is
defined similarly interchanging $(\text{FR}_1,\le_1)$ with
$(\text{FR}_2,\le_2)$.  Thus we lose ``self-dual".     
\endremark
\bigskip

\proclaim{\stag{838-8h.9} Claim}  1) ${\frak u}$ is a nice construction
framework which is self-dual.
\nl
2) For almost$_2$ every $(\bar M,\bar{\bold J},\bold f) \in
K^{\text{qt}}_{\frak s}$ the model $M_\partial$ is 
saturated above $\lambda$.
\nl
3) ${\frak u}$ is monotonic (see \scite{838-1a.24}(1)), hereditary (see
   \scite{838-3r.84}(12)), hereditary for $=_*$ if $=_*$ is a fake
   equality for ${\frak s}$ (see \scite{838-e.5D}) and has interpolation
   (see \scite{838-3r.89}).
\endproclaim
\bigskip

\demo{Proof}  1) As in earlier cases (see \scite{838-e.5L}(1)),
\scite{838-5t.21}(1)).   
\nl
2) As in \scite{838-e.5L}(2) or \scite{838-5t.21}(2).  
\nl
3) check.  \hfill$\square_{\scite{838-8h.9}}$
\enddemo
\bigskip

\proclaim{\stag{838-8h.12} Theorem}  We have $\dot I(\lambda^{++},
K^{\frak s}(\lambda^+\text{-saturated})) 
\ge \mu_{\text{unif}}(\lambda^{++},2^{\lambda^+})$ and even $\dot
I(K^{{\frak u},{\frak h}}_{\lambda^{++}}) \ge
\mu_{\text{unif}}(\lambda^{++},2^{\lambda^+})$ for any ${\frak
u}-\{0,2\}$-appropriate function ${\frak h}$  \ub{when}:
\mr
\item "{$(a)$}"   $2^\lambda < 2^{\lambda^+} < 2^{\lambda^{++}}$
\sn
\item "{$(b)$}"   ${\frak u}$ fails wnf-delayed uniqueness for {\rm WNF} see
Definition \scite{838-8h.15} below.
\endroster
\endproclaim
\bigskip

\remark{Remark}  Note that we have some versions of delayed
uniqueness: the straight one, the one with WNF and the one in \S5 and more.
\endremark
\bn
Before we prove Theorem \scite{838-8h.12}
\definition{\stag{838-8h.15} Definition}  We say that (the almost good
$\lambda$-frame) ${\frak s}$ has wnf-delayed uniqueness for WNF
\ub{when}: (if WNF is clear from the context we may omit it)
\mr
\item "{$\boxtimes$}"  for every $(M_0,N_0,a) \in
K^{3,\text{bs}}_{\frak s}$ we can find $(M_1,N_1)$ such that
{\roster
\itemitem{ $(a)$ }   $(M_0,N_0,a) \le^1_{\frak u} (M_1,N_1,a)$, 
 i.e. $(M_0,N_0,a) \le_{\text{bs}} (M_1,N_1,a)$ and
\nl

\hskip25pt WNF$(M_0,N_0,M_1,N_1)$, see clause (b) of Hypothesis 
\scite{838-8h.1} and
\sn
\itemitem{ $(b)$ }   if $(M_1,N_1,a) \le^1_{\frak u} (M_\ell,N_\ell,a)$ hence
WNF$(M_1,N_1,M_\ell,N_\ell)$ for $\ell=2,3$ and $M_2 = M_3$
\ub{then} $N_2,N_3$ are $\le_{\frak s}$-compatible over $M_2 \cup
N_0$, that is we can find a pair $(f,N')$ such that
\sn
\itemitem{ ${{}}$ }   $\quad (\alpha) \quad N_3 \le_{\frak s} N'$
\sn
\itemitem{ ${{}}$ }   $\quad (\beta) \quad f$ is a $\le_{\frak
s}$-embedding of $N_2$ ito $N'$
\sn
\itemitem{ ${{}}$ }   $\quad (\gamma) \quad f$ is the identity on
$N_0$ (not necessarily $N_1$!) and on $M_2 = M_3$.
\endroster}
\endroster
\enddefinition
\bigskip

\remark{Remark}  A point of \scite{838-8h.15} is that we look for
uniqueness among $\le^1_{\frak u}$-extensions (if \scite{838-8h.2} apply then
$\le_{\text{wnf}}$-extensions) and, of course, it is
``delayed", i.e. possibly $M_0 \ne M_1$.
\endremark
\bigskip

\demo{\stag{838-8h.19} Observation}  In Definition \scite{838-8h.15} 
we can without loss of generality  demand that $M_1$ is brimmed over $M_0$.  Hence $M_1$
can be any pregiven $\le_{\frak s}$-extension of $M_0$ brimmed over it
such that $M_1 \cap N_0 = M_0$.
\enddemo
\bigskip

\demo{Proof}  Read the definition.   \hfill$\square_{\scite{838-8h.19}}$
\enddemo
\bigskip

\proclaim{\stag{838-8h.23} Claim}  If ${\frak s}$ (satisfies \scite{838-8h.1})
and fails wnf-delayed uniqueness for {\rm WNF}
(i.e. satisfies \scite{838-8h.12}(b)) \ub{then}
${\frak u} = {\frak u}^3_{\frak s}$ has vertical coding, see Definition
\scite{838-2b.9}.  
\endproclaim
\bigskip

\demo{Proof}  Straight.  \hfill$\square_{\scite{838-8h.23}}$
\enddemo
\bigskip

\demo{Proof of \scite{838-8h.12}}  Straight by the above and Theorem \scite{838-2b.13}.
\hfill$\square_{\scite{838-8h.23}}$ 
\enddemo
\bigskip

\remark{\stag{838-8h.26} Remark}  Note that the assumption of
\scite{838-8h.23}, failure of wnf-delayed uniqueness may suffice for a
stronger version of \scite{838-8h.12} because given $\eta \in
{}^{\partial^+}2$ and $(\bar M^\eta,\bar{\bold J}^\eta,\bold f^\eta)
\in K^{\text{qt}}_{\frak s}$, we can find $2^\partial$ extensions 
$\langle (\bar M^{\eta,\rho},\bar{\bold J}^{\eta,\rho},\bold
f^{\eta,\rho}):\rho \in {}^\omega2\rangle$ and $\alpha(*) < \partial$
such that $M^{\eta,\rho}_{\alpha(0)} = N_*$ and $\langle
M^{\eta,\rho}_\partial:\rho \in {}^\partial 2\rangle$ are pairwise
non-isomorphic over $M^\eta_\partial \cup N_*$.  Does this help to 
omit the assumption $2^\lambda < 2^{\lambda^+}$?
\endremark
\bigskip

\definition{\stag{838-8h.30} Definition}  1) We say ${\frak s}$ has
uniqueness for WNF \ub{when}:

if WNF$(M^k_0,M^k_1,M^k_2,M^k_3)$ for $k=1,2$ and $f_\ell$
is an isomorphism from $M^1_\ell$ onto $M^2_\ell$ for $\ell=0,1,2$ and
$f_0 \subseteq f_1,f_0 \subseteq f_2$ \ub{then} there is a pair $(N,f)$
such that $M^2_3 \le_{\frak s} N$ and $f$ is a $\le_{\frak
s}$-embedding of $M^1_3$ into $N$ extending $f_1 \cup f_2$.
\nl
2) We say $(M_0,M_1,a) \in K^{3,\text{bs}}_{\frak s}$ has
   non-uniqueness for WNF \ub{when}: if $(M_0,M_1,a) \le^1_{\frak u} 
(M'_0,M'_1,a)$ \ub{then} we can find 
$\langle M^k_\ell:\ell \le 3,k=1,2\rangle,\langle f_\ell:\ell \le 2 \rangle$
such that
\mr
\item "{$\circledast$}"  $(a) \quad (M^k_0,M^k_1,a) \le^1_{\frak u}
   (M^k_2,M^k_3,a)$ for $k=1,2$
\sn
\item "{${{}}$}"  $(b) \quad M^k_0 = M'_0,M^k_1 = M'_1$ for $k=1,2$
\sn
\item "{${{}}$}"  $(c) \quad M^1_2 = M^2_2$ and $f_2$ is the identity
on $M^2_1$
\sn
\item "{${{}}$}"  $(d) \quad$ there is no pair $(N,f)$ such that 
$M^2_3 \le_{\frak s} N$ and $f$ is a $\le_{\frak s}$-embedding 
\nl

\hskip25pt of $M^1_3$ into $N$ extending id$_{M^k_1} \cup 
\text{ id}_{M^k_2}$ (which does not depend
\nl

\hskip25pt  on $k$).
\ermn
3) We say that ${\frak s}$ has non-uniqueness for WNF when
   some triple $(M_0,M_1,a) \in K^{3,\text{bs}}_{\frak s}$ has it.
\enddefinition
\bigskip

\demo{\stag{838-8h.34} Observation}  1) Assume ${\frak s}$ is categorical (in
$\lambda$).  \ub{Then} ${\frak s}$ has non-uniqueness for WNF 
\ub{iff} it does not have uniqueness for WNF \ub{iff}
${\frak s}$ fails existence for $K^{3,\text{ur}}_{{\frak s},{\frak u}}$,
see below.
\nl
2) If ${\frak s}$ is categorical (in $\lambda$), has 
existence for $K^{3,\text{up}}_{\frak s} = 
K^{3,\text{up}}_{{\frak s},\lambda^+}$ and has uniqueness for WNF 
and \scite{838-7v.17} holds for WNF (i.e. if $(M_1,N_1,a) \in
K^{3,\text{up}}_{\frak s}$ and $(M_1,N_1,a) \le_{\text{bs}}
(M_2,N_2,a)$ then\footnote{so really WNF = WNF$_{\frak s}$}
WNF$(M_1,M_1,N_1,N_2)$) \ub{then} ${\frak s}$ 
has existence for $K^{3,\text{uq}}_{\frak s}$ and
$K^{3,\text{up}}_{\frak s} \subseteq K^{3,\text{uq}}_{\frak s}$.
\nl
3) $K^{3,\text{uq}}_{\frak s} \subseteq K^{3,\text{up}}_{\frak s}$.
\nl
4) If ${\frak s}$ has uniqueness for WNF \ub{then} WNF is a
   non-forking relation on ${\frak K}_\lambda$ respecting ${\frak s}$.
\enddemo
\bigskip

\remark{Remark}  Note that \scite{838-8h.34}(3) is used in \S4(F).
\endremark
\bigskip

\definition{\stag{838-8h.35} Definition}  1) Let $K^{3,\text{ur}}_{{\frak
s},{\frak u}}$ be the class of 
triples $(M,N,a) \in K^{3,\text{bs}}_{\frak s}$ such
that: if $(M,N,a) \le^1_{\frak u} (M',N'_\ell,a)$ for $\ell=1,2$ then
we can find a pair $(N^*,f)$ such that $N'_2 \le_{\frak s} N^*$ and
$f$ is a $\le_{\frak s}$-embedding of $N'_1$ into $N^*$ extending
id$_N \cup \text{ id}_{M'}$.
\nl
2) ${\frak s}$ has existence for $K^{3,\text{ur}}_{\frak s}$ when: if
   $M \in K_{\frak s}$ and $p \in {\Cal S}^{\text{bs}}_{\frak s}(M)$
   then for some pair $(N,a)$ the triple $(M,N,a) \in
   K^{3,\text{ur}}_{\frak s}$ realizes $p$.
\nl
3) If WNF is WNF$_{\frak s}$ and ${\frak u}$ is defined as in
\scite{838-8h.15} above then we
may omit ${\frak u}$.
\enddefinition
\bigskip

\demo{Proof}  1) By the definition
\mr
\item "{$(*)$}"  if ${\frak s}$ has non-uniqueness for WNF
\ub{then} ${\frak s}$ does not have uniqueness for ${\frak s}$.
\ermn
Now
\bn
\ub{Case 1}:  ${\frak s}$ fails existence for $K^{3,\text{ur}}_{\frak
s}$.

We shall show that ${\frak s}$ has the non-uniqueness property; this
suffices by $(*)$.  Let $(M,p)$ exemplify it and let $(N,a)$ be such
that $(M,N,a) \in K^{3,\text{bs}}_{\frak s}$ realizes $p$ and we shall
prove that $(M,N,a)$ is as required in Definition \scite{838-8h.30}(2).

Let $(M',N',a) \in K^{3,\text{bs}}_{\frak s}$ be $\le^1_{\frak
u}$-above $(M,N,a)$.  We can find $(f_*,N_*,a_*)$ such that $(M,N_*,a_*)
\in K^{3,\text{bs}}_{\frak s}$ realizes $p$ and $f_*$ is an isomorphism
from $N'$ onto $N_*$ which maps $M'$ onto $M$ and $a$ to $a_*$.  [Why
it exists?  See \scite{838-5t.37} recalling ${\frak s}$ is categorical.]  So
$(M,N_*,a_*) \in K^{3,\text{bs}}_{\frak s}$ and $\ortp(a^*,M,N_*) = p$
hence by the choice of $M$ and $p$ clearly $(M,N_*,a_*) \notin
K^{3,\text{ur}}_{{\frak s},{\frak u}}$ so by Definition
\scite{838-8h.35}(1) we can find $M',N_1,N_2$ as there such that there are
no $(N^*,f)$ as there.  But this means that for $(M,N_*,a_*)$ 
we can find $\langle M^k_\ell:\ell \le 3,k=1,2\rangle$ 
as required in \scite{838-8h.30}(2).  By
chasing maps this holds also for $(M',N',a)$ so we are done.
\bn
\ub{Case 2}:  ${\frak s}$ has existence for $K^{3,\text{ur}}_{\frak s}$. 

We shall show that ${\frak s}$ has uniqueness for WNF.
This suffices by $(*)$.

First note 
\mr
\item "{$\boxtimes$}"  if $(*)_1$ and $(*)_2$ below then $M^1_3,M^2_3$
are isomorphic for $M_2 \cup M_1$ \ub{when}:
{\roster
\itemitem{ $(*)_1$ }  $(a) \quad \bar M = \langle M_{0,\alpha}:\alpha \le
\alpha(*)\rangle$ is $<_{\frak s}$-increasing continuous
\sn
\itemitem{ ${{}}$ }  $(b) \quad (M_{0,\alpha},M_{0,\alpha +1},a_\alpha)
\in K^{3,\text{ur}}_{{\frak s},{\frak u}}$ for $\alpha < \alpha(*)$
\sn
\itemitem{ ${{}}$ }  $(c) \quad M_{0,\alpha(*)}$ is brimmed over $M_{0,0}$
\sn
\itemitem{ $(*)_2$ }  $(\alpha) \quad$ WNF$(M_0,M_1,M_2,M^k_3)$ for
$k=1,2$
\sn
\itemitem{ ${{}}$ }  $(\beta) \quad M_0 = M_{0,0}$ and $M_1 =
M_{0,\delta}$
\sn
\itemitem{ ${{}}$ }  $(\gamma) \quad M^k_3$ is brimmed over $M_1 \cup
M_2$.
\endroster}
\ermn
[Why?  As in previous arguments in \S7, we lift $\bar M$ by
$(*)_2(\alpha)$ and clause (g)
of Definition \scite{838-7v.35} 
of ``WNF is a weak non-forking relation on
${\frak K}_{\frak s}$", i.e. being as in the proof of \scite{838-7v.33}
and then use $(*)_1(b)$.]

Next
\mr
\item "{$\boxplus$}"  we can weaken $(*)_2(\beta)$ to $M_{0,0} = M_0
\le_{\frak s} M_1 \le_{\frak s} M_{0,\alpha(*)}$.
\ermn
[By the properties of WNF.]

Checking the definitions we are done recalling \scite{838-7v.19}(2B), or
pedantically repeating its proof to get $\langle M_{0,i}:i \le
\alpha(*)\rangle$ as in $(*)_1$.
\nl
2) We are assuming that 
${\frak s}$ has existence for $K^{3,\text{up}}_{\frak s}$ 
so it suffices to prove
$K^{3,\text{up}}_{\frak s} \subseteq K^{3,\text{uq}}_{\frak s}$.  So
assume $(M_0,N_0,a) \in K^{3,\text{up}}_{\frak s}$ and $(M_0,N_0,a)
\le_{\text{bs}} (M_\ell,N_\ell,a)$ for $\ell=2,3$ and $M_2 = M_3$.  By
\scite{838-7v.17} it follows that WNF$(M_0,N_0,M_\ell,N_\ell)$ so
by Definition \scite{838-8h.4} we have $(M_0,N_0,a) \le^1_{\frak u}
(M_\ell,N_\ell,a)$.  Applying Definition \scite{838-8h.30}(1) we are done.
\nl
3),4) Clear by the definitions.    \hfill$\square_{\scite{838-8h.34}}$
\enddemo
\bn
\centerline{$* \qquad * \qquad *$}
\bigskip

\proclaim{\stag{838-8h.37} Theorem}  $\dot I(\lambda^{++},K^{\frak s}) \ge
\dot I(\lambda^{++},K^{\frak s}(\lambda^+$-saturated)) $\ge
\mu_{\text{unif}}(\lambda^{++},2^{\lambda^+})$ and even $\dot
I(K^{{\frak u},{\frak h}}_{\lambda^{++}}) \ge
\mu_{\text{unif}}(\lambda^{++},2^{\lambda^+})$ for any ${\frak
u}_{\frak s}-\{0,2\}$-appropriate ${\frak h}$ (so we can 
restrict ourselves to models $\lambda^+$-saturated above $\lambda$ and
if ${\frak s} = {\frak s}'$ also to $\tau_{\frak s}$-fuller ones) \ub{when}:
\mr
\item "{$(a)$}"  $2^\lambda < 2^{\lambda^+} < 2^{\lambda^+}$ 
\sn
\item "{$(b)$}"  ${\frak s}$ has non-uniqueness for 
{\rm WNF} (for every $M \in K_{\frak s}$)
\sn
\item "{$(c)$}"  ${\frak s}$ has {\rm wnf}-delayed uniqueness for {\rm WNF}.
\endroster
\endproclaim
\bigskip

\demo{Proof}  We first prove claim \scite{838-8h.47}.
\enddemo
\bn
Note that proving them we can use freely \scite{838-8h.15}, \scite{838-8h.4},
\scite{838-8h.19} and that wnf-delayed uniqueness replaces the use of
\scite{838-10k.11}.
\bn
\margintag{8h.41}\ub{\stag{838-8h.41} Explanation}:  As 
FR$^{\frak u}_1 = \text{ FR}^{\frak u}_2$ there
is symmetry, i.e. ${\frak u}$ is self-dual.  
The wnf-delayed uniqueness was gotten vertically, i.e. from its
failure we got a non-structure result (\scite{838-8h.12}) relying on
vertical coding, i.e. \scite{838-2b.13}.  But now
we shall use it horizontally; we shall construct over $(\bar
M,\bar{\bold J},\bold f)$ with $M_\partial \in K^{\frak
s}_{\lambda^+}$ saturated above $\lambda$, a tree $\langle
(M^\rho,\bar{\bold J}^\rho,\bold f^\ell):\rho \in {}^{\partial
\ge}2\rangle$ as in weak coding but each is not as usual but 
a sequence of length $\ell g(\rho)$ such extensions.  In fact we use the
$\lambda$-wide case of \S10, i.e. \scite{838-10k.23}, \scite{838-10k.25}
\ub{without} quoting.  
So the ``non-structure" 
is done in the ``immediate successor" of $(\bar M,\bar{\bold J},\bold f)$.   
The rest of the section is intended to make the 
rest of the construction, in the 
$\partial^+$-direction, irrelevant
(well, mod ${\Cal D}_\partial$, etc) using the wnf-delayed uniqueness
assumed in clause (c) of \scite{838-8h.37}, justified by \scite{838-8h.12}.  
The net result is that we can
find $\langle(\bar M^\rho,\bar{\bold J}^\rho,\bold f^\rho):\rho \in
{}^\partial 2)$ which are $\le^{\text{qt}}_{\frak u}$-above $(\bar
M,\bar{\bold J},\bold f)$ and for $\rho \ne \nu \in {}^\partial 2$, there is
no $\le_{\frak K}$-embedding of $M^\rho_\partial$ into
$M'_\partial$ if $(\bar M^\nu,\bold J^\nu,\bold f^\nu)
\le^{\text{qs}}_{\frak u} (\bar M',\bar{\bold J}',\bold f')$.  
That comes instead of using $\bar{\Bbb F}$, the
amalgamation choice functions in \S10.

For constructing $\langle(\bar M^\rho,\bar{\bold J}^\rho,\bold
 f^\rho):\rho \in {}^\partial 2\rangle$ as above, again we use $\langle(\bar
M^{\rho,\alpha},\bar{\bold J}^{\rho,\alpha},
\bold f^{\rho,\alpha}):\rho \in {}^i 2,\alpha < \lambda\rangle$ for $i
\le \partial$ such that for a club of $\delta < \partial$ the model
$\cup\{M^{\rho,\alpha}_\delta:\alpha < \delta\}$ is
brimmed over $M^{\rho,\beta}$ for $\beta < \gamma$.
\bigskip

\proclaim{\stag{838-8h.47} Claim}  [Under the assumptions of \scite{838-8h.37}]

If $\boxtimes$ then $\circledast$ where
\mr
\item "{$\boxtimes$}"  $(a) \quad (M,N,a) \in 
K^{3,\text{bs}}_{\frak s}$ has non-uniqueness for {\rm WNF}$_{\frak s}$
\sn
\item "{${{}}$}"  $(b) \quad \delta < \lambda^+$ is divisible by
$\lambda^3$
\sn
\item "{${{}}$}"  $(c) \quad \bold d$ is a ${\frak u}$-free
$(0,\delta)$-rectangle, let $M_\alpha = M^{\bold d}_{0,\alpha}$ for
$\alpha \le \delta,a_\alpha = a^{\bold d}_{0,\alpha}$ 
\nl

\hskip25pt for $\alpha < \delta$ and $a = a_0$
\sn
\item "{${{}}$}"  $(d) \quad (M,N,a) \le_{\text{wnf}} (M_0,M_1,a)$ 
equivalently $(M,N,\{a\}) \le^1_{\frak u} (M_0,M_1,\bold I^{\bold d}_{0,0})$
\sn
\item "{${{}}$}"  $(e) \quad M_\delta$ is brimmed over $M_\alpha$ for
$\alpha < \delta$
\sn
\item "{${{}}$}"  $(f) \quad \delta \in \text{\rm correct}
(\langle M_\alpha:\alpha \le \delta\rangle)$, i.e. if 
$M_\delta <_{\frak s} N$ \ub{then}
some $p \in {\Cal S}^{\text{bs}}_{\frak s}(M_\delta)$
\nl

\hskip25pt is realized in $N$ and does not 
fork over $M_\beta$ for some $\beta < \delta$ 
\nl

\hskip25pt (on correctness, see Definition \scite{838-5t.27})
\sn
\item "{${{}}$}"  $(g) \quad (M_0,M',b) \in \text{\rm FR}_2$ and $M'
\cap M_\delta = M_0$ and $M'$ is brimmed over $M_0$
\sn
\item "{$\circledast$}"  there are $\bold d_1,\bold d_2$ such that
{\roster
\itemitem{ $(\alpha)$ }  $\quad \bold d_\ell$ is a ${\frak u}$-free
$(1,\delta +1)$-rectangle for $\ell=1,2$
\sn
\itemitem{ $(\beta)$ }  $\quad b^{\bold d}_{0,0} = b$
\sn
\itemitem{ $(\gamma)$ } $\quad \bold d_\ell \restriction (0,\delta) =
\bold d$ for $\ell=1,2$
\sn
\itemitem{ $(\delta)$ }  $\quad \bold d_1 \restriction (1,0) = 
\bold d_2 \restriction (1,0)$
\sn
\itemitem{ $(\varepsilon)$ }  $\quad M^{\bold d_1}_{1,0} =
M^{\bold d_2}_{1,0} = M'$ and $b^{\bold d_1}_{0,0} = b = b^{\bold d_2}_{0,0}$
\sn
\itemitem{ $(\zeta)$ }  $\quad M^{\bold d_1}_{\alpha,\delta},
M^{\bold d_2}_{\alpha,\delta}$ are $\tau_{\frak s}$-incompatible over
$(M^{\bold d_1}_{1,0} = M^{\bold d_2}_{1,0}) + M_1$
\sn
\itemitem{ $(\eta)$ }  $\quad$ if $k,\bold d_1,\bold d_2,f$
satisfies $\bullet_1 - \bullet_4$ below, then we 
can find a triple
\nl

\hskip25pt  $(g,N_1,N_2)$
such that $\ell=1,2 \Rightarrow M^{\bold d_\ell}_{1,\alpha(\bold
d_\ell)} \le_{\frak s} N_\ell$ and $g$ is an
\nl

\hskip25pt  isomorphism from $N_1$
onto $N_2$ extending id$_{M'} \cup f$ where
\nl

\hskip25pt  (for $\ell=1,2$):
\sn
\itemitem{ ${{}}$ }  $\bullet_1 \quad \bold d_\ell$ 
is a ${\frak u}$-free rectangle
\sn
\itemitem{ ${{}}$ }  $\bullet_2 \quad \beta(\bold d_\ell) = 1$
\sn
\itemitem{ ${{}}$ }  $\bullet_3 \quad \alpha(\bold d_\ell) \ge \delta$ and 
$\bold d_\ell \rest (0,\delta) = \bold d$
\sn
\itemitem{ ${{}}$ }  $\bullet_4 \quad f$ is an isomorphism from $M^{\bold
d_1}_{0,\alpha(\bold d_1)}$ onto 
$M^{\bold d_2}_{0,\alpha(\bold d_2)}$ over $M^{\bold d}_{0,\delta}$.
\endroster}
\endroster
\endproclaim
\bigskip

\remark{Remark}  1) In the proof we use wnf-delayed uniqueness.
\nl
2) This claim helps.
\endremark
\bigskip

\demo{Proof}  First, letting $M_* = M'$, 
we can choose $M^1_*,M^2_*$ such that (for $\ell=1,2$)
\mr
\item "{$\circledast_1$}"  $M_* \le_{\frak s} M^\ell_*$
\sn
\item "{$\circledast_2$}"   $M_1 \le_{\frak s} M^\ell_*$
\sn
\item "{$\circledast_3$}"   $M^\ell_* \cap M_\delta = M_1$
\sn
\item "{$\circledast_4$}"  WNF$(M_0,M_1,M_*,M^\ell_*)$ hence
$(M_0,M_*,b) \le_{\text{bs}} (M_1,M^\ell_*,b)$
\sn
\item "{$\circledast_5$}"  $M^1_*,M^2_*$ are $\tau$-incompatible over
$M_* + M_1$.
\ermn
[Why?  By $\boxtimes(a)$ and $\boxtimes(g)$ 
recalling Definition \scite{838-8h.30}(2).]

Second, we choose $N^\ell_*$ for $\ell=1,2$ such that
\mr
\item "{$\circledast_6$}"  WNF$(M_1,M_\delta,M^\ell_*,N^\ell_*)$ for $\ell=1,2$
\sn
\item "{$\circledast_7$}"  wnf-delayed uniqueness: if $\ell
\in\{1,2\}$ and $M_{\delta +1},N^1_{\delta+1},N^2_{\delta+1}$ satisfies
\nl

WNF$(M_\delta,M_{\delta +1},N^\ell_*,N^k_{\delta +1})$ for $k=1,2$ \ub{then}
we can find $(f,N)$ such
\nl

 that $N^2_{\delta+1} \le_{\frak s} N$ and $f$
is a $\le_{\frak s}$-embedding of $N^1_{\delta +1}$ into $N$ over
$M^1_*$ 
\nl

(hence over $M_*$) and over $M_{\delta +1}$.
\ermn
[Why is this possible?  As $M_\delta$ is brimmed over $M_1$ by cluase (e) of
$\boxtimes$ we are assuming, and
${\frak s}$ has wnf-delayed uniqueness by clause (c) of
Theorem \scite{838-8h.37} and we apply it $(M_1,M^\ell_*,b) \le_2
(M_\delta,N^\ell_*,b)$ recalling ${\frak u}$ is self-dual and \scite{838-8h.19}.]

Note that in $\circledast_6$ we can replace $N^\ell_*$ by
$N^\ell_{**}$ if $N^\ell_* \le_{\frak s} N^\ell_{**}$ or $M^\ell_*
\cup M_\delta \subseteq N^\ell_{**} \le_{\frak s} N^\ell_*$. 

Third, by the properties of WNF, for $\ell=1,2$ we can
choose $N^\ell_{**}$ and a
${\frak u}$-free $(1,\delta)$-rectangle $\bold d'_\ell$ with 
$M^{\bold d'_\ell}_{1,0} = M^\ell_*,b^{\bold d_\ell}_{0,0} = b,\bold
d'_\ell \restriction (0,\delta) = \bold d \restriction
([0,0],[1,\delta])$ and $M^{\bold d'_\delta}_{1,\delta} \le_{\frak s}
N^\ell_{**},N^\ell_* \le_{\frak s} N^\ell_{**}$.  

Now $\bold d'_1,\bold d'_2$ are as required. \hfill$\square_{\scite{838-8h.47}}$
\enddemo
\bigskip

\demo{\stag{838-p.4.13} Proof of \scite{838-8h.37}}  In this case, for variety,
instead of using a theorem on ${\frak u}$ from \S2 or \S3, we do it
directly (except quoting \scite{838-7f.7}).  
We fix a stationary $S \subseteq \partial$ such that
$\partial \backslash S \notin \text{ WDmId}(\partial)$ and $S$ is a
set of limit ordinals.

We choose
${\frak g}$ witnessing \scite{838-8h.9}(2) for $S$ so without loss of generality  $S_{\frak g} =
S$ so ${\frak g}$ is ${\frak u}$-2-appropriate.  Let ${\frak h}$ be
any ${\frak u}-\{0,2\}$-appropriate function.  
We restrict ourselves to $K^{\text{qt},*}_{\frak u} := \{(\bar
M,\bar{\bold J},\bold f) \in K^{\text{qt}}_{\frak u}:M_\partial \in
K^{\frak s}_{\lambda^+}$ is saturated (above $\lambda$), $M_\partial$
has universe an ordinal $< \partial^+$ and $\bold f \restriction
(\partial \backslash S)$ is constantly $1$ and $(\alpha \in \partial
\backslash S \Rightarrow M_{\alpha +1}$ is brimmed over $M_\alpha\}$.  
We now choose
$\langle(\bar M^\eta,\bar{\bold J}^\eta,\bold f^\eta):\eta \in
{}^\alpha(2^\partial)\rangle$ by induction on $\alpha < \partial^+$
such that
\mr
\item "{$\oplus_\alpha$}"  $(a) \quad (\bar M^\eta,\bar{\bold
J}^\eta,\bold f^\eta) \in K^{\text{qt},*}_{\frak u}$ for 
$\eta \in {}^\alpha(2^\partial)$ 
\sn
\item "{${{}}$}"  $(b) \quad \langle(\bar M^{\eta \restriction \beta},
\bar{\bold J}^{\eta \restriction \beta},\bold f^{\eta \restriction
\beta}):\beta \le \alpha\rangle$ is $\le^{\text{qt}}_{\frak
u}$-increasing continuous
\sn
\item "{${{}}$}"  $(c) \quad$ if $\alpha = \beta +1$ and $\beta$ is
non-limit and $\eta \in {}^\alpha(2^\partial)$ then the pair
$((\bar M^{\eta \restriction \beta},\bar{\bold J}^{\eta \restriction
\beta},\bold f^{\eta \rest \beta})$,
\nl

\hskip25pt $(\bar M^\eta,\bar{\bold J}^\eta,\bold f^\eta))$ 
 obeys\footnote{we may combine} ${\frak g}$ if $\alpha$ is even 
and obeys ${\frak h}$ if $\alpha$ is odd
\sn
\item "{${{}}$}"  $(d) \quad$ if $\alpha = \beta +1,\beta$ is a limit
ordinal, $\nu \in {}^\beta(2^\partial)$ and $\varepsilon_1 \ne \varepsilon_2 <
2^\partial$ and so 
\nl

\hskip25pt $\eta^\ell = \nu \char 94 \langle \varepsilon_\ell
\rangle$ is from ${}^\alpha(2^\partial)$ for $\ell=1,2$ \ub{then}
not only $M^{\eta^1}_\partial,M^{\eta^2}_\partial$
\nl

\hskip25pt  are not isomorphic over $M^\nu_\partial$,
but if $(\bar M^{\eta^\ell},
\bar{\bold J}^{\eta^\ell},\bold f^{\eta^\ell}) \le^{\text{qt}}_{\frak u} (\bar
M^\ell,\bar{\bold J}^\ell,\bold f^\ell)$ 
\nl

\hskip25pt for $\ell=1,2$ then
$M^1_\partial,M^2_\partial$ are not isomorphic over $M^\nu_\partial$.
\ermn
By \scite{838-7f.7} this suffices. 
For $\alpha=0$ and $\alpha$ limit there are no problems (well we have
to show that the limit exists which hold by \scite{838-1a.37}(4), and belongs to
$K^{\text{qt},*}_{\frak u}$, but this is easy by \scite{838-7v.67}(2)).

So assume $\alpha = \beta +1,\eta \in {}^\beta 2$ and we should choose
$\langle(\bar M^{\eta \char 94 <\varepsilon>},
\bar{\bold J}^{\eta \char 94 <\varepsilon>},\bold f^{\eta \char 94
<\varepsilon>}):\varepsilon < 2^\partial\rangle$, let $\gamma_*$ be
the universe of $M^\eta_\partial$.  

Let $E_1$ be a club of $\partial = \lambda^+$ such that if $\alpha <
\delta \in E_1$ then $\bold f(\alpha) < \delta$ and $M^\eta_\delta$ is
brimmed over $M^\eta_\alpha$.  Let $E_2 = E_1 \cup \{[\delta,\delta +
\bold f(\delta)]:\delta \in S \cap E_1\}$, and without loss of generality  $M^\eta_0$ is
brimmed and if $\delta \in S \cap E_1$ then $M_{\delta +1}$ is brimmed
over $M_\delta$ (can use ${\frak g}$ to guarantee this, or increase it
inside $M^\eta_\gamma$ with no harm).  Let $h$ be the increasing
continuous function from $\lambda^+$ onto $E_2$ and $E = \{\delta <
\lambda^+:\delta$ a limit ordinal and $h(\delta) = \delta\}$ a club of
$\lambda^+ = \partial$.

So
\mr
\item "{$\boxplus$}"  $(a) \quad (M^\eta_\alpha,M^\eta_{\alpha
+1},\bold J^\eta_\alpha) \in \text{ FR}^2_{\frak u} = \text{ FR}^1_{\frak u}$ 
\sn
\item "{${{}}$}"  $(b) \quad M^\eta_\theta$ is brimmed
\sn
\item "{${{}}$}"  $(c) \quad M^\eta_{h(\alpha +1)}$ is brimmed over
$M^\eta_{h(\alpha)}$ if $\alpha \in E_1 \cap S$
\sn
\item "{${{}}$}"  $(d) \quad$ if $h(\alpha) \in S$ then $h(\alpha +1)
= h(\alpha) +1$ (used?).
\ermn
Now ${\frak s}$ has non-uniqueness
for WNF hence we can find $(N,a)$ such that the triple
$(M^\eta_0,N,a)$ has the non-uniqueness property for WNF;
without loss of generality  $N \backslash M^\eta_0$ is $[\gamma_*,\gamma_* 
+ i_{<>})$ for some ordinal $i_{<>} \le \lambda$.

Now we choose $\bold d_\rho$ for $\rho \in {}^\varepsilon 2$ by
induction on $\varepsilon < \lambda^+$ such that (recalling ${\frak
u}$ is self-dual; note that $\bold d$ looks inverted
letting $\bar \alpha^\varepsilon = \langle \lambda (1 + \zeta):\zeta \le
\varepsilon\rangle$
\mr
\item "{$\odot$}"  for $\rho \in {}^\varepsilon 2$
{\roster
\itemitem{ $(a)$ }  $\bold d_\rho$ is an ${\frak u}$-free
$(\bar\alpha^\varepsilon,\varepsilon)$-triangle
\sn
\itemitem{ $(b)$ }  $(M^\eta_0,N,\{a\}) = (M^{\bold
d_\rho}_{0,0},M^{\bold d_\rho}_{1,0},\bold J^{\bold d_\rho}_{0,0})$
\sn
\itemitem{ $(c)$ }  if $\zeta < \varepsilon$ then $\bold d_{\rho
\restriction \zeta} = \bold d_\rho \restriction (\bar
\alpha^\zeta,\zeta)$
\sn
\itemitem{ $(d)$ }  $M^{\bold d_\rho}_{0,\zeta} =  
M^\eta_{h(\zeta)}$ for $\zeta \le \varepsilon$ and
$\bold I^{\bold d_\rho}_{0,\zeta} = \bold J^\eta_\zeta$ for $\zeta <
\varepsilon$ 
\sn
\itemitem{ $(e)$ }  $M^{\bold d_\rho}_{i+1,\zeta +1}$ is brimmed
over $M^{\bold d_\rho}_{i,\zeta +1} \cup 
M^{\bold d_\rho}_{i,\zeta +1}$ when $\zeta < \varepsilon,i < \lambda(1
+ \zeta)$
\sn
\itemitem{ $(f)$ }  if $\varepsilon = \zeta +1$ and $i < \lambda$
then\footnote{actually can waive clause (f),(g)}
$M^{\bold d_\varepsilon}_{\lambda \varepsilon + i+1,\varepsilon}$ is
brimmed over $M^{\bold d_\rho}_{\lambda \varepsilon +i,\varepsilon}$
\sn
\itemitem{ $(g)$ }  if $\varepsilon = \zeta +1$ and $i < \lambda$ and
$p \in {\Cal S}^{\text{bs}}_{\frak s}(M^{\bold d_\varepsilon}_{\lambda
\varepsilon + i,\varepsilon})$ \ub{then} for $\lambda$ ordinals $j \in
[i,\lambda)$, the non-forking extension of $p$ in ${\Cal
S}^{\text{bs}}_{\frak s}(M^{\bold d_\varepsilon}_{\lambda \varepsilon
+j,\varepsilon})$ is realized by the $b \in \bold J^{\bold
d_\varepsilon}_{\lambda \varepsilon +1,\varepsilon}$ in $M^{\bold
d_\varepsilon}_{\lambda \varepsilon +j+1,\varepsilon}$
\sn 
\itemitem{ $(h)$ }  if $\varepsilon \in E_! \cap S$ so $1 + \varepsilon =
\varepsilon$ then clause $(\eta)$ of $\circledast$ of \scite{838-8h.47}
holds with \nl
{\rm dual}$(\bold d_\rho \restriction (\lambda \varepsilon,0))$,
$\lambda \varepsilon$, 
$(M^\eta_0,N,a)$,
$(M^\eta_\varepsilon,M^\eta_{\varepsilon +1}, \bold J^\eta_\varepsilon)$,
\nl
{\rm dual}$(\bold d_{\rho \char 94 <0>} \restriction [0,\lambda \varepsilon],[\varepsilon,\varepsilon +1])$, 
{\rm dual}$(\bold d_{\rho \char 94 <1>} \restriction [0,\lambda
\varepsilon],[\varepsilon,\varepsilon +1))$ here standing for 
 $\bold d$, $\delta$, $(M,N,a)$, $(M_0,M',\{b\})$, $\bold d_1$, $\bold d_2$
there
\sn
\itemitem{ $(i)$ }  the set $M^{\bold d_\rho}_{\lambda(1 +
\varepsilon),\varepsilon} \backslash M^\eta_\varepsilon$ is
$[\gamma_*,\gamma_* + \lambda(1 + \varepsilon))$.
\endroster}
\ermn
There is no problem to carry the definition.

Lastly for $\rho \in {}^\partial 2$ we define $(\bar
M^{\eta,\rho},\bar{\bold J}^{\eta,\rho},\bold f^{\eta,\rho})$ by: (let
$E \subseteq \partial = \lambda^+$ be a thin enough club):
\mr
\item "{$\boxtimes$}"  $(a) \quad M^{\eta,\rho}_\varepsilon = 
M^{\bold d_{\rho \restriction \varepsilon}}_{\lambda(1+ \varepsilon),
\varepsilon}$ for $\varepsilon < \lambda$
\sn
\item "{${{}}$}"  $(b) \quad \bold f^{\eta,\rho} = \bold f^\eta$
\sn
\item "{${{}}$}"  $(c) \quad \bold J^{\eta,\rho}_\varepsilon = 
\bold J^\eta_\varepsilon$ when $\varepsilon \in \cup\{[\delta,\delta +
\bold f^\eta(\delta)):\delta \in E\}$.
\ermn
Now let $\langle S_\varepsilon:\varepsilon < \partial =
\lambda^+\rangle$ be a partition of $\lambda^+ \backslash S$ 
to (pairwise disjoint)
sets from $(\text{WDmId}_{\lambda^+})^+$.

Now we define a function $\bold c$:
\mr
\item "{$(*)_1$}"  its domain is the set of $\bold x =
(\rho_1,\rho_2,f,\bold d)$ such that: for some $\varepsilon \in S \cap
E \subseteq \lambda^+$
{\roster
\itemitem{ $(a)$ }  $\rho_1,\rho_2 \in {}^\varepsilon 2$
\sn
\itemitem{ $(b)$ }  $\bold d$ is a ${\frak u}$-free $(\varepsilon,1)$-rectangle
with $\bold I^{\bold d}_{\zeta,0} = \emptyset$ for $\zeta \le \varepsilon$
\sn
\itemitem{ $(c)$ }  $M^{\bold d}_{\zeta,0} = 
M^{\bold d_{\rho_2}}_{\lambda(1 + \zeta),\zeta}$ 
for $\zeta \le \varepsilon$
\sn
\itemitem{ $(d)$ }  $\bold J^{\bold d}_{\zeta,0} = 
\bold I^{\bold d_{\rho_2}}_{0,\zeta}$ for $\zeta \le \varepsilon$
\sn
\itemitem{ $(e)$ }  $M^{\bold d}_{\varepsilon,1} \backslash
M^{\bold d_{\rho_2}}_{\lambda(1 + \varepsilon),\varepsilon} \subseteq
[\gamma_* + \lambda^+,\gamma_* + \lambda^+ + \lambda^+)$
\sn
\itemitem{ $(f)$ }  $f$ is a $\le_{\frak s}$-embedding of $M^{\bold
d_{\rho_1}}_{\lambda(1+ \varepsilon),\varepsilon}$ into $M^{\bold
d}_{\varepsilon,1}$ over $M^\eta_{h(\varepsilon)} = 
M^{\bold d_{\rho_\ell}}_{0,\varepsilon}$ for $\ell =1,2$
\endroster}
\item "{$(*)_2$}"  for $\bold x = (\rho_1,\rho_2,f,\bold d)$ as above,
say with $\varepsilon = \ell g(\rho_1) = \ell g(\rho_2)$ we have
 $\bold c(\bold x) =  1$ \ub{iff}
there are $\bold d^+,f^+,N^*$ such that letting $\nu_\ell = \rho_\ell
\char 94 \langle 0 \rangle$ for $\ell=1,2$:
{\roster
\itemitem{ $(a)$ }  $\bold d^+$ is a ${\frak u}$-free 
$(\varepsilon +1,1)$-rectangle
\sn
\itemitem{ $(b)$ }  $\bold d^+ \restriction (\varepsilon,1) = \bold d$
\sn
\itemitem{ $(c)$ }  $M^{\bold d_{\nu_2}}_{\lambda(1+ \varepsilon),
\varepsilon +1} \le_{\frak s} N^*$ and $M^{\bold d^+}_{\varepsilon
+1,1} \le_{\frak s} N^*$
\sn
\itemitem{ $(d)$ }  $f^+$ is a $\le_{\frak s}$-embedding of $M^{\bold
d_{\nu_1}}_{\lambda(1 +\varepsilon),\varepsilon +1}$ into $N^*$ extending
$f \cup \text{ id}_{M^\eta_{h(\varepsilon) +1}}$.
\endroster}
\ermn
Note
\mr
\item "{$\circledast$}"  if $(\bar M^{\eta,\rho},
\bar{\bold J}^{\eta,\ell},\bold f^{\eta,\rho}) \le_{\text{qt}} (\bar M^*,\bold
J^*,\bold f^*)$ and the universe of $M^*_\partial$ is an ordinal $<
\lambda^{++}$ and $\pi$ is a one-to-one mapping from $M^*_\partial$
onto $\gamma_* + \partial +\partial$ over $\gamma_* + \partial$, then we
can find $\langle \bold e_\varepsilon:\varepsilon < \lambda^+\rangle$
such that for some club $E$ of $\partial$
{\roster
\itemitem{ $(a)$ }  $\bold e_\varepsilon$ is a ${\frak u}$-free
$(\varepsilon,1)$-rectangle
\sn  
\itemitem{ $(b)$ }  $M^{\bold e_\varepsilon}_{\zeta,0} = 
M^{\bold d_\rho}_{\lambda(1 + \zeta),\zeta}$ for 
$\zeta \le \varepsilon$
\sn  
\itemitem{ $(c)$ }  $\bold J^{\bold e_\varepsilon}_{\zeta,0} = 
\bold I^{\bold d_\rho}_{\zeta,\zeta}$ for $\zeta < \varepsilon$
\sn
\itemitem{ $(d)$ }  $M^{\bold e_\varepsilon}_{\varepsilon,1} 
\backslash M^{\bold d_\rho}_{\lambda(1 + \varepsilon),\varepsilon} 
\subseteq [\gamma_* + \lambda^+,\gamma_* + \lambda^+ + \lambda^+)$
\sn
\itemitem{ $(e)$ }  if $\varepsilon \in E$ then there is a ${\frak
u}$-free $(\varepsilon +1,1)$-rectangle $\bold  e^+_\varepsilon$ 
such that $\bold e^+_\varepsilon \restriction
(\varepsilon +1) = \bold e_{\varepsilon +1}$
\sn
\itemitem{ $(f)$ }  $\langle M^{\bold e_\varepsilon}_{\varepsilon,1}:
\varepsilon < \partial\rangle$ is a $\le_{\frak  s}$-increasing 
continuous sequence with union $\pi(M^*_\partial)$ which has universe
$\gamma_* + \lambda^+ + \lambda^+$.
\endroster}
\ermn
For each $\zeta < \lambda^+$ as $S_\zeta \subseteq \lambda^+$ does not
belong to the weak dimaond ideal, there is a sequence $\varrho_\zeta
\in {}^{(S_\zeta)}2$ such that
\mr
\item "{$(*)_3$}"  for any $\rho_1,\rho_2 \in {}^\partial2$ such that
$\rho_1 \restriction S_\zeta = 0_{S_\zeta}$ and $N_* \in K^{\frak
s}_{\lambda^+},(\bar M^*,\bar{\bold J}^*,\bold f^*) \in
K^{\text{qs}}_{\frak s}$ which is $\le_{\text{qt}}$-above $(\bar
M^{\eta,\rho_2},\bar{\bold J}^{\eta,\rho_2},\bold f^{\eta,\rho_2})$
and $|M^*_\partial| = \gamma_{**} < \partial^+$ and $\le_{{\frak
K}[{\frak s}]}$-embedding $f$ of $M^{\eta,\rho_1}_\partial$ into
$M^*_\partial$ over $M^\eta_\partial$, letting $\langle \bold
e_\delta:\delta < \partial\rangle,\pi$ be as in $\circledast$
for the pair $((\bar M^\eta,\bar{\bold J}^\eta,\bold f^\eta),(\bar
M^*,\bar{\bold J}^*,\bold f^*)$ the set $\{\delta \in
S_\varepsilon:\bold c(\rho_1 \restriction \delta,\rho_2 \restriction
\delta,f \restriction \delta,\bold e_\delta) = \varrho(\delta)\}$ is
stationary.
\ermn
Now for any $u \subseteq \partial$ we define $\rho_u \in {}^\lambda 2$ by
\mr
\item "{$(*)_4$}"  for $\zeta < \partial,\ell < 2$ let $\rho_u
\restriction S_{2 \zeta + \ell}$ be $0_{2 \zeta + \ell}$ if $[\zeta \notin u
\leftrightarrow \ell=0]$ and $\varrho_\zeta$ otherwise
\sn
\item "{$(*)_5$}"  let $\rho _u \rest S_0$ be constantly zero.
\ermn
Let $\langle u(\alpha):\alpha < 2^\partial\rangle$ list ${\Cal
P}(\gamma)$ and for $\alpha < 2^\partial$ let $(\bar M^{\eta \char
94<\alpha>},\bold J^{\eta \char 94<\alpha>},\bold f^{\eta\char 94
<\alpha>})$ be $(\bar M^{\eta,\rho_{u(\alpha)}},
\bold J^{\eta,\rho_{u(\alpha)}},\bold f^{\eta,\rho_{u(\alpha)}})$.

Clearly they are as required.   \hfill$\square_{\scite{838-8h.37}}$
\enddemo
\bn
\centerline{$* \qquad * \qquad *$}
\bn
\margintag{8h.59}\ub{\stag{838-8h.59} Exercise}:  1) Assume ${\frak t}$ is an almost good
$\lambda$-frame, ${\frak u} = {\frak u}^1_{\frak t}$ from Definition
\scite{838-e.5I} \ub{then} for some ${\frak u}-\{0,2\}$-appropriate
${\frak h}$, for every $M \in K^{{\frak t},{\frak h}}_{\lambda^{++}}$
we have
\mr
\item "{$(a)$}"  $M$ is $\lambda^+$-saturated
\sn
\item "{$(b)$}"  if $M_0 \in K_{\frak t}$ is $\le_{{\frak K}[{\frak t}]}
M$ and $p \in {\Cal S}^{\text{bs}}_{\frak t}(M_0)$
\ermn
then {\rm dim}$(p,M) =
\lambda^{++}$ that is, there is a sequence $\langle a_\alpha:\alpha <
\lambda^{++}\rangle$ of members of $M$ realizing $p$ such that: if
$M_0 \le_{\frak t} M_1 <_{{\frak K}[{\frak t}]} M$ then $\{\alpha <
\lambda^{++}:\ortp_{\frak t}(a,M_1,M)$ does not fork over $M_0\}$ is a
co-bounded subset of $\lambda^{++}$.
\nl
2) Similarly if ${\frak t}$ has existence for $K^{3,\text{up}}_{\frak
s}$ and ${\frak u} = {\frak u}^3_{\frak t}$, see Definition
\scite{838-8h.4}. 
\bigskip

\proclaim{\stag{838-8h.61} Theorem}  $\dot I(\lambda^{++},K^{\frak s}) \ge
\dot I(\lambda^{++},K^{\frak s}(\lambda^+$-saturated)) $\ge
\mu_{\text{unif}}(\lambda^{++},2^{\lambda^+})$ and even $\dot
I(K^{{\frak u},{\frak h}}_{\lambda^{++}}) \ge
\mu_{\text{unif}}(\lambda^{++},2^{\lambda^+})$ for any ${\frak
u}_{\frak s}-\{0,2\}$-appropriate ${\frak h}$ (so we can 
restrict ourselves to models $\lambda^+$-saturated above $\lambda$ and
if ${\frak s} = {\frak s}'$ also to $\tau_{\frak s}$-fuller ones) \ub{when}:
\mr
\item "{$(a)$}"  $2^\lambda < 2^{\lambda^+} < 2^{\lambda^+}$ 
\sn
\item "{$(b)$}"  ${\frak s}$ has non-uniqueness for 
{\rm WNF} (for every $M \in K_{\frak s}$)
\sn
\item "{$(c)$}"  $K$ is categorical in $\lambda$
\sn
\item "{$(d)$}"  ${\frak u}$ has existence for
$K^{3,\text{up}}_{{\frak s},\lambda^+}$.
\endroster
\endproclaim
\bigskip

\demo{Proof}  We shall use \scite{838-8h.12}, \scite{838-8h.34},
\scite{838-8h.37}.  So assume toward contradiction that the conclusion
fails.  We try to apply Theorem \scite{838-8h.12}, now its conclusion
fails by our assumption toward contradiction, and clause (a) there
which says ``$2^\lambda < 2^{\lambda^+} < 2^{\lambda^{++}}$" holds by
clause (a) of the present theorem.  So necessarily clause (b) of Theorem
\scite{838-8h.12} fails which means that ${\frak u}$ has wnf-delayed uniqueness,
see Definition \scite{838-8h.15}.

Next we try to apply Theorem \scite{838-8h.37}, again it assumption fails
by our assumption toward contradiction, and among its assumptions
clause (a) which says that ``$2^\lambda < 2^{\lambda^+} <
2^{\lambda^{++}}$"  holds by clause (a) of the present theorem, and clause
(c) which says ``${\frak s}$ has wnf-delayed uniqueness" has just been
proved.  So necessarily clause (b) of \scite{838-8h.37} fails which means
that ${\frak s}$ fails non-uniqueness for WNF, i.e. for some $M$.

Now we apply Observation \scite{838-8h.34}, noting that its assumption
``${\frak s}$ is categorical in $\lambda$" holds by clause (c) of the
present theorem,
so by the previous sentence one of the equivalent phrases the first
fails, hence all of them.  In particular ${\frak s}$ has uniqueness
for WNF.  \hfill$\square_{\scite{838-8h.61}}$
\enddemo
\goodbreak

\head {\S9 The combinatorial part} \endhead  \resetall \sectno=9
 \spuriousreset
\bigskip

We deal here with the ``relatively" pure-combinatorial parts.  We do
 just what is necessary.  We can get results on $\dot I \dot
E(\partial^+,{\frak K}_{\frak u})$, we can weaken the cardinal arithmetic
assumptions to $\emptyset \notin \text{ DfWD}_\partial$, see \cite{Sh:E45},
we can weaken the demands on ${\frak K}$; \ub{but} not here.
\bn
Recall the obvious by the definitions:
\proclaim{\stag{838-7f.7} Theorem}  If $2^\partial < 2^{\partial^+}$ 
\ub{then} $\{M_\eta/\cong :\eta \in {}^{\partial^+}(2^\partial)$ 
and $\|M_\eta\| = \partial^+\}$ has cardinality $\ge
\mu_{\text{unif}}(\partial^+,2^\partial)$ \ub{when} the
following conditions hold:
\mr
\item "{$\circledast$}"  $(a) \quad \bar M = \langle M_\eta:\eta \in
{}^{\partial^+ >}(2^\partial)\rangle$
\sn
\item "{${{}}$}"  $(b) \quad$ for $\eta \in {}^{\partial^+ >}
(2^\partial)$ the model $M_\eta$ has cardinality
$\le\partial$ and for notational
\nl

\hskip25pt   simplicity has universe an ordinal $< \partial^+$
\sn
\item "{${{}}$}"  $(c) \quad M_\eta \subseteq M_\nu$ if $\eta
\triangleleft \nu \in {}^{\partial^+ >}(2^\partial)$, so \ub{no}
a.e.c. appear here!
\sn
\item "{${{}}$}"  $(d) \quad \langle M_{\eta \restriction
\alpha}:\alpha < \ell g(\eta)\rangle$ is $\subseteq$-increasing
continuous for any $\eta \in {}^{\partial^+ >}(2^\partial)$
\sn
\item "{${{}}$}"  $(e) \quad M_\eta := \cup\{M_{\eta \restriction
\alpha}:\alpha < \partial^+\}$ for $\eta \in
{}^{\partial^+}(2^\partial)$
\sn
\item "{${{}}$}"  $(f) \quad$ if $\eta \in {}^{\partial^+
>}(2^\partial)$ and $\alpha_1 < \alpha_2 < 2^\partial$ and 
$\eta \char 94 \langle \alpha_\ell\rangle \trianglelefteq \nu_\ell \in
{}^\delta 2$ 
\nl

\hskip25pt for $\ell=1,2$ and $\delta < \partial^+$ \ub{then}
$M_{\nu_2},M_{\nu_1}$ are not isomorphic over $M_\eta$
\nl 

\hskip25pt or just
\sn
\item "{${{}}$}"  $(f)^-_1 \quad$ for $\eta \in {}^{\partial^+ >}2$,
there is ${\Cal U}_\eta \subseteq 2^\partial$ of cardinality
$2^\partial$ such that: if
\nl

\hskip25pt  $\alpha_0 \ne \alpha_1$ are from ${\Cal
U}_\eta$ and $\eta \char 94 \langle \alpha_\ell\rangle \trianglelefteq
\nu_\ell \in {}^\delta 2$ for $\ell < 2$ and $\delta < \partial^+$ 
\nl

\hskip25pt \ub{then} $M_{\nu_0},M_{\nu_1}$ are not isomorphic over $M_\eta$.
\endroster
\endproclaim
\bigskip

\demo{Proof}  Concerning clause (f)$^-_1$ we can by renaming get clause
(f), so in the rest of the proof of 
\scite{838-7f.7} we can ignore clause (f)$^-_1$.

Note that $\Xi_0 := \{\eta \in
{}^{\partial^+}(2^\partial):\|M_\eta\| < \partial^+\}$ has cardinality $\le
2^\partial$ (because for each $\eta \in \Xi_0$ there is
$\alpha_\eta < \partial^+$ such that $M_\eta = M_{\eta \restriction
\alpha_\eta}$; and note that by clause (f) we have 
$\eta_1 \in \Xi_0 \wedge \eta_2 \in \Xi_0
\wedge \alpha_{\eta_1} = \alpha_{\eta_2} \wedge \eta_1 \restriction
\alpha_{\eta_1} = \eta_2 \restriction \alpha_{\eta_2} \Rightarrow \eta_1 =
\eta_2)$.
So by clause (b) of $\circledast$ it follows that
$\eta \in {}^{\partial^+}(2^\partial) \backslash \Xi_0
\Rightarrow |M_\eta| = \partial^+$.

It suffices to assume that $\Xi \subseteq
{}^{\partial^+}(2^\partial)$ has cardinality $<
\mu_{\text{unif}}(\partial^+,2^\partial)$ and find $\eta \in
{}^{\partial^+}(2^\partial)$ such that $\nu \in \Xi \Rightarrow
M_\eta \ncong M_\nu$, because without loss of generality  $\Xi_0 \subseteq \Xi$.

Let $\langle \eta_\zeta:\zeta < |\Xi|\rangle$ list $\Xi$ and
let $N_\zeta := M_{\eta_\zeta}$ and
toward contradiction for every $\nu \in {}^{\partial^+}(2^\partial)$ we can
choose $\zeta_\nu = \zeta(\nu) < |\Xi|$ and an isomorphism $f_\nu$ from
$M_\nu$ onto $N_{\zeta(\nu)}$, so $f_\nu$ is a function from $M_\nu$ onto
$M_{\eta_{\zeta(\nu)}}$.

For $\zeta < |\Xi|$ let $W_\zeta = \{\nu \in
{}^{\partial^+}(2^\partial):\zeta_\nu = \zeta\}$, so clearly:
\mr
\item "{$(*)_1$}"  ${}^{\partial^+}(2^\partial)$ is equal to
$\cup\{W_\zeta:\zeta < |\Xi|\}$.
\ermn
[Why?  Obvious by our assumption toward contradiction.]
\mr
\item "{$(*)_2$}"  if $i < \partial^+$ and $\rho \in {}^i(2^\partial)$,
\ub{then} there are no $\varepsilon_1 \ne \varepsilon_2 < 2^\partial$
such that $\rho \char 94 \langle \varepsilon_\ell \rangle
\triangleleft \nu_\ell \in W_\zeta$ for $\ell = 1,2$ and $f_{\nu_1}
\rest M_\eta = f_{\nu_2} \rest M_\eta$.
\ermn
[Why?  By Clause (f) of the assumption.]

Together we get a contradiction to the definition of
$\mu_{\text{unif}}(\partial^+,2^\partial)$, see Definition \scite{838-0z.6}(7).  
\hfill$\square_{\scite{838-7f.7}}$ 
\enddemo
\bn
Similarly
\proclaim{\stag{838-7f.9} Claim}  1) In \scite{838-7f.7} we can replace
$2^\partial$ by $\langle \chi_i:i < \partial \rangle$ with $\chi_i \le
2^\partial$.
\nl
2) Also we can weaken clause (f) or (f)$^-_1$ there by demanding $\delta =
   \partial^+$.
\nl
3) Assume ${\frak K}$ is an a.e.c. and in \scite{838-7f.7} we demand
   $M_\nu \le_{\frak K} M_\eta \in {\frak K}$ for $\nu \triangleleft
\eta \in {}^{\partial^+ >}(2^\partial)$.  If we strengthen there clause 
(f)$^-_1$ by strengthening the conclusion to
``if $\eta \char 94 \langle \ell \rangle \triangleleft \eta_\ell \in
{}^{\partial^+}2$ for $\ell=1,2$ then $M_{\nu_0},M_{\nu_1}$ cannot be
$\le_{\frak K}$-amalgamated over $M_\eta"$ \ub{then}:
\mr
\item "{$(*)$}"  for every $\Xi \subseteq {}^{\partial^+}(2^\partial)$
of cardinality $< \mu_{\text{unif}}(\partial^+,2^\partial)$ for some
$\eta \in {}^{\partial^+}(2^\partial)$ the model $M_\eta$ has
cardinality $\partial^+$ and cannot be $\le_{\frak K}$-embedded in
$M_\nu$ for any $\nu \in \Xi$
\sn
\item "{$(**)$}"  if $2^{\partial^+} > (2^\partial)^+$ \ub{then} there
is $\Xi \subseteq {}^{\partial^+}2$ of cardinality $2^{\partial^+}$
such that if $\eta \ne \nu \in \Xi$ then $M_\eta$ cannot be
$\le_{\frak K}$-embedded into $M_\nu$.
\endroster
\endproclaim
\bigskip

\demo{Proof}  Left to the reader (easier than \scite{838-7f.21} below and
will not be used here).  \hfill$\square_{\scite{838-7f.9}}$
\enddemo
\bigskip

\remark{Remark}  Why do we prefer to state \scite{838-7f.7}?  As this is
how it is used.
\endremark
\bigskip

\proclaim{\stag{838-7f.10} Lemma}  Assuming $2^\theta = 2^{< \partial} <
2^\partial$ (and naturally but not used $2^\partial < 2^{\partial^+}$)
 and $\circledast(a)-(e)$ of \scite{838-7f.7}, a
sufficient condition for clause $(f)^-_1$ of \scite{838-7f.7} is:
\mr
\item "{$(a)^+$}"  $\bar M = \langle M_\eta:\eta \in
 {}^{\partial^+ >}(2^\partial)\rangle$ and $\langle M^*_{\eta,\zeta}:\zeta <
 \partial\rangle$ is $\subseteq$-increasing with union $M_\eta$ such
 that $\zeta < \partial \Rightarrow \|M_{\eta,\zeta}\| < \partial$
\sn
\item "{$(f)^-_2$}"  for each $\eta \in {}^{\partial^+ >}2$ we can find
$\langle M_{\eta,\rho}:\rho \in {}^{\partial \ge}2\rangle$ such that
{\roster
\itemitem{ $(\alpha)$ }  $\langle M_{\eta,\rho}:\rho \in {}^\partial 2
\rangle$ is a subsequence of $\langle M_{\eta \char 94 <\alpha>}:\alpha <
2^\partial\rangle$ with no repetitions so $M_{\eta,\rho} = 
M_{\eta \char 94 <\alpha(\rho)>}$ for some one-to-one function
$\rho \mapsto \alpha(\rho)$ from ${}^\partial 2$ to $2^\partial$
\sn
\itemitem{ $(\beta)$ }  if $\rho \in {}^{\partial >}2$ \ub{then}
$M_{\eta,\rho} \in K_{< \partial}$
\sn
\itemitem{ $(\gamma)$ }  if $\rho \in {}^{\partial >}2$ \ub{then} 
$\langle M_{\eta,\rho \restriction \alpha}:\alpha \le \ell
 g(\rho)\rangle$ is $\subseteq$-increasing continuous 
\sn
\itemitem{ $(\delta)$ }   $\cup\{M_{\eta,\rho
\restriction \varepsilon}:\varepsilon <\partial\}$ is equal to
$M_{\eta,\rho} = M_{\eta \char 94 <\alpha(\rho)>}$ for any $\rho 
\in {}^\partial 2$
\sn
\itemitem{ $(\varepsilon)$ }  $\partial$ is regular uncountable and
for some sequence $\langle
S_\varepsilon:\varepsilon < \partial\rangle$ of pairwise disjoint
non-small stationary subsets of $\partial$ (i.e. $\varepsilon <
 \partial \Rightarrow S_\varepsilon \in (\text{\rm WDmId}_\partial)^+)$
we have
\sn
\itemitem{ ${{}}$ }  $(*) \quad$ for every $\varepsilon < \partial$,
there is a pair $(\bar g,\bold c) = (\bar g^\varepsilon,\bold
c^\varepsilon)$, may not depend 
\nl

\hskip30pt on $\varepsilon$ such that:
\sn
\itemitem{ ${{}}$ }  $\quad \bullet_1 \quad \bar g = \langle
g_{\eta,\rho}:\rho \in {}^{\partial} 2\rangle$
\sn
\itemitem{ ${{}}$ }  $\quad \bullet_2 \quad g_{\eta,\rho}$ is a function from
$\partial$ to ${\Cal H}_{< \partial}(\partial^+)$
\sn
\itemitem{ ${{}}$ }  $\quad \bullet_3 \quad$ if $2^\partial >
\partial^+$ and $\rho_0,\rho_1 \in {}^\partial 2,
\rho_1 \restriction S_\varepsilon$ is constantly zero,
\nl

\hskip55pt $\delta < \partial^+,\eta \char 94 \langle \alpha(\rho_\ell)\rangle
\trianglelefteq \nu_\ell \in {}^\delta(2^\partial)$ for $\ell=0,1$ 
and $f$ is an 
\nl

\hskip55pt isomorphism from $M_{\nu_0}$ onto $M_{\nu_1}$ \ub{then} 
for some club $E$ of 
\nl

\hskip55pt $\partial$, if $\zeta \in E \cap S_\varepsilon$ we have
\nl

\hskip55pt $\rho_0(\zeta) = \bold c^\varepsilon(\rho_0 \restriction
\zeta,M_{\eta,\rho_0 \restriction \zeta},\rho_1 \restriction \zeta,
M_{\eta,\rho_1 \restriction\zeta},g_{\eta,\rho_0} 
\restriction \zeta$,
\nl

\hskip55pt $g_{\eta,\rho_1} \restriction \zeta,M_{\nu_0,\zeta},M_{\nu_\zeta},
f \restriction M_{\eta,\rho_0 \restriction \zeta})$
\sn
\itemitem{ ${{}}$ }  $\quad \bullet_4 \quad$ if 
$2^\partial = \partial^+$: as above but
$\bold c$ is preserved by any partial
\nl

\hskip55pt  order preserving function from
$\partial^+$ to $\partial^+$ extending {\rm id}$_{M^\eta_\partial}$.
\endroster}
\endroster
\endproclaim
\bigskip

\remark{Remark}  1) We can immitate \scite{838-7f.21}.
\nl
2) If $2^\partial = \partial^+$ then it follows that $\partial =
   \partial^{< \partial}$, so they give us stronger ways to construct.
\endremark
\bigskip

\demo{Proof} First
\mr
\item "{$\boxtimes$}"  for $\eta \in {}^{\partial^+ >}(2^\partial)$
and $\varepsilon < \partial$ there is $\varrho_{\eta,\varepsilon} \in
{}^{(S_\varepsilon)}2$ such that:
{\roster
\itemitem{ $(*)$ }   if $\rho_0 \ne \rho_1$ are from ${}^\partial 2,\eta
\char 94 \langle \alpha(\rho_\ell)\rangle \triangleleft \nu_\ell \in
{}^\delta 2,\delta < \partial^+$ and $f$ is an isomorphism from
$M_{\nu_0}$ onto $M_{\nu_1}$ \ub{then} for stationary many $\zeta \in
S_\varepsilon$ we have:

$$
\align
\varrho_\eta(\zeta) = &\bold c^\varepsilon(\rho_0 \rest \zeta,
M_{\eta,\rho_0 \rest \zeta},\rho_1 \rest \zeta,
M_{\eta,\rho_1 \rest \zeta},g_{\eta,\rho_0} \rest
\zeta,g_{\eta,\rho_1} \rest \zeta, \\
  &M^*_{\nu_0,\zeta},M^*_{\nu_1,\zeta},f \rest M_{\eta,\rho_0 \rest \zeta})
\endalign
$$
\endroster}
\ermn
[Why?  First if $2^\theta \ge \partial^+$, use the definition of
$S_\varepsilon \notin \text{ WDmId}(\partial)$, (see more in the proof
of \scite{838-7f.15}).  If $2^\theta = \partial \wedge 2^\partial =
\partial^+$, the proof is similar using the invariance of $\bold
c^\varepsilon$, i.e. $\bullet_4$.

Lastly, if $2^\theta = \partial \wedge 2^\partial > \partial^+$, use
$\mu_{\text{wd}}(\partial) > \partial^+$, see \scite{838-0z.7}(1A).]

Let $\eta \in {}^{\partial^+ >}(2^\partial)$.  For any $w \subseteq
\partial$ we define $\rho_{\eta,w} \in {}^\partial 2$ as follows:
$\rho_{\eta,w}(i)$ is $\varrho_{\eta,\varepsilon}$ \ub{if} for some
$\varepsilon <  \partial$ and $\ell < 2$ we have $i \in S_{2
\varepsilon + \ell} \wedge [\varepsilon \in w \equiv \ell =1]$ and is
zero otherwise.  So $\{\alpha_\eta(\rho_{\eta,w}):w \subset
\partial\}$ is as required in (f)$^-_1$.
  \hfill$\square_{\scite{838-7f.10}}$
\enddemo
\bigskip

\proclaim{\stag{838-7f.14} Theorem}  If 
$2^\partial < 2^{\partial^+}$ and $\mu =
\mu_{\text{unif}}(\partial^+,2^\partial)$ \ub{then}:
\mr
\item "{$(A)$}"  $\emptyset \notin \text{\rm UnfTId}_\mu(\partial^+)$
\sn
\item "{$(B)$}"  {\rm UnfTId}$_{\mu_1}(\partial^+)$ is
$\mu_1$-complete when $\aleph_0 \le \mu_1 = \text{\rm cf}(\mu_1) <
\mu$; see \scite{838-0z.6}(4),(5)
\sn
\item "{$(C)$}"   $\mu =
2^{\partial^+}$ except maybe \ub{when} (all the conditions below hold):
{\roster
\itemitem{ $\circledast$ }  $(a) \quad \mu < \beth_\omega$
\sn
\itemitem{ ${{}}$ }    $(b) \quad \mu^{\aleph_0} = 2^{\partial^+}$
\sn
\itemitem{ ${{}}$ }   $(c) \quad$ there is a family ${\Cal A} \subseteq
[\mu]^{\partial^+}$ of cardinality $\ge 2^{\partial^+}$ such that
\nl

\hskip30pt  the intersection of
any two distinct members of ${\Cal A}$ is finite.
\endroster}
\endroster 
\endproclaim
\bigskip

\remark{Remark}   So in the aleph sequence 
$\mu$ is much larger than $2^\partial$, when $\mu \ne 2^{\partial^+}$.
\endremark
\bigskip

\demo{Proof}  By \cite[AP,1.16]{Sh:f} we have clauses (b) + (c) of
$\circledast$ and they imply clause (a) by \cite{Sh:460} (or see
\cite{Sh:829}).   \hfill$\square_{\scite{838-7f.14}}$
\enddemo
\bigskip

\proclaim{\stag{838-7f.12} Claim}  Assume $\partial > \theta \ge \aleph_0$ 
is regular and $2^\theta = 2^{< \partial} < 2^\partial$.  
\ub{Then} $\{M_\eta/\cong:\eta \in
{}^\partial 2$ and $M_\eta$ has cardinality $\partial\}$ 
has cardinality $2^\partial$ \ub{when} the following condition holds:
\mr
\item "{$\circledast$}"  $(a) \quad \bar M = \langle M_\eta:\eta \in
{}^{\partial \ge} 2\rangle$ with $M_\eta$ a $\tau$-model
\sn
\item "{${{}}$}"  $(b) \quad$ for $\eta \in {}^\partial 2,\langle
M_{\eta \restriction \alpha}:\alpha \le \partial\rangle$ is
$\subseteq$-increasing continuous
\sn
\item "{${{}}$}"  $(c) \quad$ if $\eta \in {}^{\partial >}2$ and
$\eta \char 94 \langle \ell\rangle
\trianglelefteq \nu_\ell \in {}^\alpha 2$ for $\ell=0,1$ and 
$\alpha < \partial$ \ub{then} $M_{\nu_0},M_{\nu_1}$ 
\nl

\hskip25pt are not isomorphic over $M_{<>}$ or just for $\alpha = \partial$
\sn
\item "{${{}}$}"  $(d) \quad M_{<>}$ has cardinality $< \partial$
\sn
\item "{${{}}$}"  $(e) \quad M_\eta$ has cardinality $\le \partial$
for $\eta \in {}^{\partial >}2$.
\endroster 
\endproclaim
\bigskip

\demo{Proof}  As in the proof of \scite{838-7f.7} we can ignore the $\eta
\in {}^\partial 2$ for which $M_\eta$ has cardinality $< \partial$.

As $\partial^{\|M_{<>}\|} \le 2^{<\partial} < 2^\partial$
this is obvious, see \yCITE[0.9]{88r} or see Case 1 in the
proof of \scite{838-7f.21} below.  \hfill$\square_{\scite{838-7f.12}}$ 
\enddemo
\bn
The following is used in \yCITE[3c.22]{E46} (and can be used in
\yCITE[2b.32]{E46}; but compare with \scite{838-7f.12}!) 
\proclaim{\stag{838-7f.15} Claim}  The set $\{M_\eta /\cong :\eta\in
{}^\partial 2$ and $M_\eta$ has cardinality $\partial\}$ has cardinality
$\ge \mu$ when
\mr
\item "{$\boxtimes_1$}"  $\partial = { \text{\rm cf\/}}(\partial) 
> \aleph_0$ and 
{\roster
\itemitem{ (a) }  $M_\eta$ is a $\tau$-model of cardinality 
$< \partial$ for $\eta \in {}^{\partial >} 2$ 
\sn
\itemitem{ (b) }  for each $\eta \in {}^\partial 2,
\langle M_{\eta \restriction \alpha}:\alpha < \partial \rangle$ 
is $\subseteq$-increasing continuous with union, called $M_\eta$
\sn
\itemitem { (c) }  if $\eta \in {}^{\partial >}2,\eta \char 94 \langle
\ell \rangle \triangleleft \rho_\ell
\in {}^\partial 2$ for $\ell = 1,2$ \ub{then} $M_{\rho_1},M_{\rho_2}$ are not
isomorphic over $M_\eta$
\endroster}
\item "{$\boxtimes_2$}"  $\partial \notin 
{ \text{\rm WDmId\/}}_{< \mu}(\partial)$, e.g. $\mu =
\mu_{\text{wd}}(\partial)$. 
\endroster
\endproclaim
\bigskip

\demo{Proof}  Let $\Xi = \{\eta \in {}^\partial 2:M_\eta$ has
cardinality $< \partial\}$ and for $\eta \in \Xi$ let $\alpha_\eta
= \text{ min}\{\alpha \le \partial:M_\eta = M_{\eta \rest \alpha}\}$,
clearly 
\mr
\item "{$\boxdot_1$}"  $(a) \quad \eta \in \Xi$ implies $\alpha_\eta
< \partial$
\sn
\item "{${{}}$}"  $(b) \quad$ if $\eta \in \Xi$ and $\eta \rest
\alpha \triangleleft \nu \in \Xi \backslash \{\eta\}$ 
then $\alpha_\nu > \alpha_\eta$.
\ermn
For each $\varrho \in {}^\partial 2$ we define $F_\varrho:{}^{\partial
\ge}2 \rightarrow {}^{\partial \ge}2$ by: for $\eta \in {}^\alpha 2$
let $F_\varrho(\eta) \in {}^{2 \ell g(\eta)}2$ be defined by
$(F_\varrho(\eta))(2i) = \eta(i),(F_\varrho(\eta))(2i+1) = \varrho(i)$
for $i < \alpha$.  Easily $\langle \text{Rang}(F_\varrho):\varrho \in
{}^\partial 2\rangle$ are pairwise disjoint, hence for some $\varrho
\in {}^\gamma 2$, the sets Rang$(F_\varrho)$ is disjoint to $\Xi$
so without loss of generality  (by renaming):
\mr
\item "{$\boxdot_2$}"  $\Xi = \emptyset$.
\ermn
Let $\{N_\varepsilon:\varepsilon < \varepsilon_*\}$ be a maximal subset 
of $\{ M_\rho:\rho \in {}^\partial 2\}$ consisting of pairwise
non-isomorphic models.

Without loss of generality the universe of each $M_\eta,\eta \in
{}^{\partial >}2$ is an ordinal $\gamma_\eta < \partial$ and so the
universe of each $M_\eta,\eta \in {}^\partial 2$ is $\gamma_\eta :=
\cup\{\gamma_{\eta \rest i}:i < \partial\} = \partial$, in
particular the universe of $N_\varepsilon$ is $\partial$ and 
$\eta \in {}^\partial 2 \Rightarrow \gamma_\eta = \partial$.
For $\alpha < \partial$ and $\eta \in
{}^\alpha 2$ let the function $h_\eta$ be $h_\eta(i) = M_{\eta \restriction
(i+1)}$ for $i < \ell g(\eta)$.  For each $\varepsilon <
\varepsilon_*$ we define $\Xi_\varepsilon \subseteq {}^\partial 2$ by 
$\Xi_\varepsilon = \{\eta \in {}^\partial 2:M_\eta$ is 
isomorphic to $N_\varepsilon\}$.  

For $\eta \in \Xi_\varepsilon$ choose $f^\varepsilon_\eta:
M_\eta \rightarrow N_\varepsilon$, an isomorphism, hence 
$f^\varepsilon_\eta \in {}^\partial \partial$.  

By the assumption
\mr
\item "{$\boxdot_3$}"  if $\varepsilon < \varepsilon_*$ and $\eta \in
{}^{\partial >}2$ and $\eta \char 94 \langle \ell \rangle
\triangleleft \nu_\ell \in \Xi_\varepsilon$ for $\ell=0,1$ 
then $f^\varepsilon_{\nu_0} \rest \gamma_\eta \ne f^\varepsilon_{\nu_1} \rest 
\gamma_\eta$.
\ermn
We also for each $\varepsilon < \varepsilon_*$ define a function 
(= colouring) $\bold c_\varepsilon$ 
from $\dsize \bigcup_{\alpha < \partial}({}^\alpha 2 \times
{}^\alpha \partial$) to $\{0,1\}$ by:
\mr
\item "{$\boxdot_4$}"  $\bold c_\varepsilon(\eta,f) \text{ is }: 0$ \ub{if} 
there is $\nu$ such that $\eta \triangleleft \nu \in
 \Xi_\varepsilon$ and $f \subseteq f^\varepsilon_\nu$ and 
$\nu(\ell g(\eta))=0$
\nl
$\bold c_\varepsilon(\eta,f) \text{ is: } \, 1 \quad \text{
\ub{if} otherwise}$.
\ermn
Now for any $\eta \in \Xi_\varepsilon$, the set

$$
E^\varepsilon_\eta = \{\delta < \partial:
\gamma_{\eta \restriction \delta} = \delta
\text{ and } f^\varepsilon_\eta \restriction 
\delta \text{ is a function from } \delta \text{ to } \delta\}
$$
\mn
is clearly a club of $\partial$.  

Now
\mr
\item "{$\boxdot_5$}"  if $\varepsilon < \varepsilon_*,
\eta \in \Xi_\varepsilon$ and $\delta \in E^\varepsilon_\eta$ then 
$\bold c_\varepsilon(\eta \restriction \delta,f^\varepsilon_\eta 
\restriction \delta) =  \eta(\delta)$.
\ermn
[Why? If $\eta(\delta) = 0$ then $\eta \rest \delta$
witness that $\bold c_\varepsilon(\eta \rest \delta,f^\varepsilon_\eta
\rest \delta) = 0$.
If $\eta(\delta)=1$ just recall $\boxdot_3$.]

Hence we have $\Xi_\varepsilon \in \text{ WDmTId}(\partial)$.  To get a
contradiction it is enough to prove $\cup\{\Xi_\varepsilon:\varepsilon <
\varepsilon_*\} \ne {}^\partial 2$, but as $\varepsilon_* < \mu$ 
clearly $\dsize \bigcup_{\varepsilon < \varepsilon_*} \Xi_\varepsilon$
belongs to $\text{WDmId}_{< \mu}(\partial)$ hence is
not ${}^\partial 2$, so we are done.  \hfill$\square_{\scite{838-7f.15}}$
\enddemo
\bn
The following is used in \yCITE[2b.32]{E46}, 
\yCITE[3c.28]{E46}, \yCITE[3c.22]{E46} which
repeat the division to cases.
\proclaim{\stag{838-7f.21} Claim}  The set $\{M_\eta/\cong :\eta \in
{}^{\partial}2$ and $\|M_\eta\| = \partial\}$ has cardinality
$2^{\partial}$ \ub{when}:
\mr
\item "{$\boxtimes_1$}"   $M_\eta$ is a $\tau$-model of cardinality $<
\partial$ for $\eta \in {}^{\partial >}2,
\langle M_{\eta \restriction \alpha}:\alpha \le \ell g(\eta) \rangle$ is
$\subseteq$-increasing continuous, and: if $\delta <
\delta(1) < \partial$ are limit ordinals, $\eta_0,\eta_1 \in {}^\delta 2$ and
$\eta_0 \char 94 \langle \ell \rangle \triangleleft \nu_\ell \in
{}^{\delta(1)} 2$ and $\eta_1 \char 94 \langle 0 \rangle 
\triangleleft \nu'_\ell \in 
{}^{\delta(1)} 2$ for $\ell =0,1$ \ub{then} there are no $f_0,f_1$ such that
{\roster
\itemitem{ $(\alpha)$ }   $f_\ell$ is an
isomorphism from $M_{\nu_\ell}$ onto $M_{\nu'_\ell}$ for $\ell =0,1$
\sn
\itemitem{ $(\beta)$ }  $f_0 \restriction M_{\eta_0} = f_1 \restriction
M_{\eta_0}$ and $M_{\eta_1} = f_0(M_{\eta_0})$ 
\sn
\itemitem{ $(\gamma)$ }   for some
$\rho_0,\rho_1 \in {}^\partial 2$ we have $\nu'_\ell
\triangleleft \rho_\ell$ for $\ell=0,1$ and $M_{\rho_0},M_{\rho_1}$ are
isomorphic over $M_{\eta_1}$
\endroster} 
\item "{$\boxtimes_2$}"   $\partial = 
{ \text{\rm cf\/}}(\partial) > \aleph_0$, 
$\partial \notin { \text{\rm WDmId\/}}_{< \mu}(\partial)$ (hence 
$2^{< \partial} < 2^\partial$) and moreover
\sn
\item "{$\boxtimes_3$}"   $\partial$ is a successor 
cardinal, or at least there is no $\partial$-saturated normal ideal 
on $\partial$, or at least ${\text{\rm WDmId\/}}(\partial)$ 
is not $\partial$-saturated (which holds if for some 
$\theta < \partial,\{ \delta < \partial:
{\text{\rm cf\/}}(\delta) = \theta\} \notin { \text{\rm WDmId\/}}
(\partial)$ because the ideal is normal).
\endroster
\endproclaim
\bigskip

\remark{\stag{838-7f.24} Remark}  1) Compare with 
\yCITE[3.5]{88r} - which is
quite closed but speak on ${\frak K}$ rather than on a specific $\langle
M_\eta:\eta \in {}^{\partial >}2 \rangle$.  Can we get 
$\dot I \dot E_{\frak K}(\partial,{\frak K}) = 2^\partial$?
Also $\lambda^+$ there corresponds to $\partial$ here, a minor change.   
\nl
2) The parallel claim was inaccurate in the \cite[\S3]{Sh:576}.
\nl
3) Used in \yCITE[2b.32]{E46}.
\endremark
\bigskip

\demo{Proof of \scite{838-7f.21}}  Easily, as in the proof of \scite{838-7f.15}
\wilog
\mr
\item "{$\boxdot_1$}"  $\eta \in {}^\partial 2 \Rightarrow
\|M_\eta\| = \partial$ while, of course, preserving $\boxtimes_1$.
\ermn
We divide the
proof into cases according to the answer to the following:
\mn
\underbar{Question}:  Is there $\eta^* \in {}^{\partial >}2$ such that for
every $\nu$ satisfying $\eta^* \trianglelefteq \nu \in {}^{\partial >}2$ 
there are $\rho_0,\rho_1 \in 
{}^{\partial >}2$ such that: $\nu \triangleleft \rho_0,
\nu \trianglelefteq \rho_1$, and for any $\nu_0,\nu_1 \in {}^\partial 2$ 
satisfying $\rho_\ell \triangleleft \nu_\ell$, (for $\ell =0,1$)  the
models $M_{\nu_0},M_{\nu_1}$ are not isomorphic over $M_{\eta^*}$?
\mn
But first we can find a function 
$h:{}^{\partial >} 2 \rightarrow {}^{\partial >}2$, such that:
\mr
\widestnumber\item{$(b)_{\text{yes}}$}
\item "{$(*)$}"  the function $h$ is one-to-one, mapping ${}^{\partial
>}2$ to ${}^{\partial >}2$, preserving $\triangleleft$, satisfying 
$(h(\nu)) \char 94
\langle \ell \rangle \trianglelefteq h(\nu \char 94 \langle \ell
\rangle)$ and $h$ is continuous, for $\nu \in {}^\partial 2$ we 
let $h(\nu) := \dbcu_{\alpha < \partial}
h(\nu \restriction \alpha)$, so $\ell g(\eta) < \partial
\Leftrightarrow \ell g(h(\eta)) < \partial$ and:
\sn
\item "{$(b)_{\text{yes}}$}"  \ub{when} the answer to the question 
above is yes, it is exemplified by $\eta^* = h(\langle \rangle)$ and
$M_{h(\rho_0)},M_{h(\rho_1)}$ are not isomorphic over $M_{h(<>)}$ whenever $\nu
\in {}^{\partial >} 2$ and $h(\nu \char 94 <\ell>) \triangleleft
\rho_\ell \in {}^\partial 2$ for $\ell = 0,1$
\sn
\item "{$(b)_{\text{no}}$}"  when the answer to the question above is no,
$h(\langle \rangle) = \langle \rangle$ and if 
$\alpha +1 < \beta < \partial,\eta \in {}^{\alpha +1} 2$ and $h(\eta)
\triangleleft \rho_\ell \in {}^\beta 2$ for $\ell = 1,2$ \ub{then} we can
find $\nu_1,\nu_2$ and $g^*$ such that $\rho_\ell \triangleleft
\nu_\ell \in {}^\partial 2$ and $g^*$ is an isomorphism from
$M_{\nu_1}$ onto $M_{\nu_2}$ over $M_{h(\eta \restriction \alpha)}$.
\ermn
[Why can we get $(b)_{\text{no}}$?  We choose $h(\eta)$ for $\eta \in
{}^\alpha 2$ by induction on $\alpha$ such that $h(\eta) = \eta$ for
$\alpha = 0,h(\eta) = \cup\{h(\eta \restriction \beta):\beta <
\alpha\}$ when $\alpha$ is a limit ordinal, and if $\alpha = \beta
+1,\ell < 2$ apply the assumption (``the answer is no") with $h(\eta)
\char 94 <\ell>$ standing for $\eta^*$ and let $h(\eta \char 94
<\ell>)$ be a counterexample to ``for every $\nu$"; so we get even
more then the 
promise; the isomorphism is over $M_{h(\eta \rest \alpha) \char 94
\langle \ell\rangle}$ rather than $M_{h(\eta \rest \alpha)}$, and 
note that $h(\eta) \char 94 \langle \ell 
\rangle \triangleleft h(\eta \char 94 \langle
\ell \rangle)$.]
\enddemo
\bn
\underbar{Case 1}:  The answer is yes.   

We do not use the non-$\partial$-saturation of WDmId$(\partial)$ in this case. 
Without loss of generality $h$ is the identity, by renaming.
\newline
For any $\eta \in {}^\partial 2$ and $\subseteq$-embedding $g$ of 
$M_{\langle \rangle}$ into $M_\eta := 
\dsize \bigcup_{\alpha < \partial} M_{\eta \restriction \alpha}$, let

$$
\Xi_{\eta,g} := \{ \nu \in {}^\partial 2:
\text{there is an isomorphism from }  
M_\nu \text{ onto } M_\eta \text{ extending } g\}
$$

$$
\Xi_\eta := \{ \nu \in {}^\partial 2:\text{there is an isomorphism from }
M_\nu \text{ onto } M_\eta \}.
$$
\mn
So: 
\mr
\item "{$\boxdot_2$}"  $|\Xi_{\eta,g}| \le 1$ for any $g$ and 
$\eta \in {}^\partial 2$.
\ermn
[Why?  As if $\nu_0,\nu_1 \in \Xi_{\eta,g}$
are distinct then for some ordinal $\alpha < \partial$ and $\nu \in {}^\alpha
2$ we have $\nu := \nu_0 \restriction \alpha = \nu_1 \restriction
\alpha,\nu_0(\alpha) \ne \nu_1(\alpha)$ and use the choice of
$h(\nu \char 94 \langle \ell \rangle)$), see $(b)_{\text{yes}}$ above.]
\sn
Since $\Xi_\eta = \cup \{\Xi_{\eta,g}:g \text{ is a } \le_{\frak K}
\text{-embedding of } M_{\langle \rangle} \text{ into } M_\eta \}$, we have
\mr
\item "{$\boxdot_3$}"  $|\Xi_\eta| \le \partial^{\|M_{\eta^*}\|} 
\le 2^{< \partial}$. 
\ermn
Hence we can by induction on $\zeta < 2^\partial$ choose $\eta_\zeta \in
{}^\partial 2 \backslash \dsize \bigcup_{\xi < \zeta} 
\Xi_{\eta_\xi}$, (exist
by cardinality considerations as $2^{< \partial} < 2^\partial$).  Then 
$\xi < \zeta \Rightarrow M_{\eta_\xi} \ncong M_{\eta_\zeta}$ so we have
proved the desired conclusion. 
\bn
\underbar{Case 2}:  The answer is no.

Without loss of generality $M_\eta$ has as universe the ordinal
$\gamma_\eta < \partial$ for $\eta \in {}^{\partial >} 2$.
 \newline
Let $\langle S_i:i < \partial \rangle$ be a partition of $\partial$ to sets,
none of which is in WDmId$(\partial)$, possible by the assumption
$\boxtimes_3$.  For each $i < \partial$ 
we define a function $\bold c_i$ as follows:
\mr
\item "{$\boxdot_4$}"  if $\delta \in S_i$ and $\eta,\nu \in {}^\delta
2$ and $\gamma_\eta = \gamma_\nu = \delta = \gamma_{h(\eta)} 
= \gamma_{h(\nu)}$, and $f:\delta \rightarrow \delta$ then
{\roster
\itemitem{ $(a)$ }  $\bold c_i(\eta,\nu,f) = 0$ \ub{if} we can find  
$\eta_1,\eta_2 \in {}^\partial 2 \text{ satisfying } h(\eta) 
\char 94 \langle 0 \rangle \triangleleft \eta_1 \text{ and }
h(\nu) \char 94 \langle 0 \rangle \triangleleft \eta_2 \text{ such that } 
  f \text{ can be extended to an isomoprhism from } 
M_{\eta_1} \text{ onto } M_{\eta_2}$
\sn
\itemitem{ $(b)$ }  $\bold c_i(\eta,\nu,f) = 1$ otherwise.
\endroster}
\ermn
So for $i < \partial$, as $S_i \notin \text{ WDmId}(\partial)$, for some 
$\varrho^*_i \in {}^\partial 2$ we have:
\mr
\item "{$(*)_i$}"  for every $\eta \in {}^\partial 2,\nu \in {}^\partial 2$ and
$f \in {}^\partial \partial$ the following set of ordinals is 
stationary:
$$
\{\delta \in S_i:\bold c_i(\eta \restriction \delta,\nu \restriction \delta,
f \restriction \delta) = \varrho^*_i(\delta)\}.
$$
\ermn
Now for any $X \subseteq \partial$ let $\eta_X \in {}^\partial 2$ be 
defined by:
\mr
\item "{$\boxdot_5$}"  if $\alpha \in S_i$ then 
$i \in X \Rightarrow \eta_X(\alpha) = 1 - \varrho^*_i (\alpha)$ and
$i \notin X \Rightarrow \eta_X(\alpha) = 0$. 
\ermn
For $X \subseteq \partial$ let $\rho_X := \eta_{\{2i:i \in X\} 
\cup \{2i+1:i \in \partial \backslash X\}} \in {}^\partial 2$. Now we
shall show
\mr
\item "{$\oplus$}"  if $X,Y \subseteq \partial$, and $X \ne Y$ then 
$M_{h(\rho_X)}$ is not isomorphic to $M_{h(\rho_Y)}$.
\ermn
Clearly $\oplus$ will suffice for finishing the proof.

Assume toward a contradiction that $f$ is an isomorphism 
of $M_{h(\rho_X)}$ onto $M_{h(\rho_Y)}$; as $X \ne Y$ there 
is $i$ such that $i \in X \Leftrightarrow
i \notin Y$ so there is $j \in \{2i,2i+1\}$ such that
\mr
\item "{$\boxdot_6$}"   $\rho_X \restriction S_j = \langle 1 -
\varrho^*_j(\alpha):\alpha \in S_j \rangle$ and 
$\rho_Y \restriction S_j$ is identically zero. 
\ermn
Clearly the set $E = \{ \delta:f \text{ maps } \delta 
\text{ onto } \delta$ and $h(\rho_X \restriction
\delta),h(\rho_Y \restriction \delta) \in {}^\delta 2$ and the
universes of $M_{h(\rho_X \restriction \delta)},M_{h(\rho_Y
\restriction \delta)}$ are $\delta\}$ is a club of $\partial$ 
and hence $S_j \cap E \ne \emptyset$.

So if $\delta \in S_j \cap E$ then $f$ extends 
$f \restriction M_{h(\rho_X) \restriction \delta}$ and $f$ is 
an isomorphism from $M_{h(\rho_X)}$ onto $M_{h(\rho_Y)}$;
by the choice of $\varrho^*_j$ we can choose $\delta \in S_j \cap E$ such
that:
\mr
\item "{$\boxdot_7$}"  $\bold c_j(\rho_X \restriction \delta,
\rho_Y \restriction \delta,f \restriction \delta) = \varrho^*_j(\delta)$.
\ermn
Also by the choice of $j$, i.e. $\boxdot_6$ we have
\mr
\item "{$\boxdot_8$}"  $\rho_X(\delta) = 1 - \varrho^*_i(\delta)$ and
$\rho_Y(\delta) = 0$.
\endroster
\bn
\ub{Subcase 2A}:  $\rho_X(\delta)=0$.

Now $\rho_X \restriction \delta \triangleleft (\rho_X \restriction
\delta) \char 94 \langle \rho_X(\delta) \rangle = (\rho_X \restriction
\delta) \char 94 <0> \triangleleft \rho_X \in
{}^\partial 2$ and $(\rho_Y \restriction \delta) \triangleleft 
(\rho_Y \restriction \delta) \char 94 \langle
\rho_Y(\delta) \rangle = \rho_Y \char 94 <0> \triangleleft \rho_Y \in
{}^\partial 2$ (as $\rho_X(\delta)=0$ by the case and $\rho_Y(\delta) = 0$ as
$\delta \in S_j$ and the choice of $j$, i.e. by $\boxdot_6$).  
Hence $f,\rho_X,\rho_Y$ witness that by the definition of $\bold c_j$ we get 
\mr
\item "{$\otimes_1$}"  $\bold c_j(\rho_X \restriction \delta,\rho_Y 
\restriction \delta,f \restriction \delta) = 0$.
\ermn
Also, by $\boxdot_8$
\mr
\item "{$\otimes_2$}"  $0  = \rho_X(\delta) = 1 - \varrho^*_j(\delta)$
so $\varrho^*_j(\delta) = 1$.
\ermn
But $\otimes_1 + \otimes_2$ contradict the choice of $\delta$,
(indirectly the choice of $\varrho^*_j$), i.e., contradicts $\boxdot_7$. 
\bn
\ub{Subcase 2B}:  $\rho_X(\delta) =1$.

By $\boxdot_7$ and $\boxdot_6$ and the case assumption we have
$\bold c_j(\rho_X \restriction 
\delta,\rho_Y \restriction \delta,f
\restriction \delta) = \varrho^*_j(\delta) = 1 - \rho_X(\delta) = 0$ hence
by the definition of $\bold c_j$ there are $\eta_1,\eta_2 \in {}^\partial 2$
such that $h(\rho_X \restriction \delta) \char 94 <0> 
\triangleleft \eta_1,h(\rho_Y \restriction \delta) \char 94
<0> \triangleleft \eta_2$, and there is an isomorphism $g$ from
$M_{\eta_1}$ onto $M_{\eta_2}$ extending $f \rest \delta$.  
There is $\delta_1 \in
(\delta,\partial)$ such that: $f$ maps $M_{h(\rho_X) \restriction
\delta_1}$ onto $M_{h(\rho_Y) \restriction \delta_1}$ and $g$ maps
$M_{\eta_1 \restriction \delta_1}$ onto $M_{\eta_2 \restriction
\delta_2}$.  Now by the choice of $h$, i.e., clause $(b)_{\text{no}}$
above, with $h(\rho_Y \restriction \delta) \char 94
<0>,\eta_2 \restriction \delta_1,h(\rho_Y) \restriction
\delta_1$ here standing for $\nu,\rho_1,\rho_2$ there and get
$\nu_1,\nu_2,g^*$ as there so $\eta_2 \restriction \delta_1
\triangleleft \nu_1 \in {}^\partial 2,h(\rho_Y) \restriction \delta_1
\triangleleft \nu_2 \in {}^\partial 2$ and $g^*$ is an isomorphism
form $M_{\nu_1}$ onto $M_{\nu_2}$ over $M_{h(\rho_Y) \restriction
\delta \char 94 <0>}$.  So this contradicts $\boxtimes_1$ in the assumption
of the claim with $\delta,\delta_1,h(\rho_X \restriction
\delta),h(\rho_Y \restriction \delta),\eta_1 \restriction \delta_1,
h(\rho_X \restriction \delta_1),\eta_2 \rest \delta_1,
h(\rho_Y) \rest \delta_1,f \rest M_{\eta_1 \rest \delta_1},g \rest
M_{h(\rho_X \rest \delta_1)},\nu_1,\nu_2$ here standing
for $\delta,\delta(1),\eta_0,\eta_1,\nu_0,\nu_1,\nu'_0,\nu'_1,f_0,f_1,
\rho_0,\rho_1$ there.  \hfill$\square_{\scite{838-7f.21}}$
\goodbreak

\head {\S10 Proof of the non-structure theorems with choice functions}
\endhead  \resetall 
 \spuriousreset
\bigskip

When we try to apply several of the coding properties, we have to use
the weak diamond (as e.g. in \scite{838-7f.21}), but in order to use it we have
to fix some quite arbitrary choices; this is the role of the $\bar{\Bbb
F}$'s here.  Of course, we can weaken \scite{838-10k.1}, but no need here.

\demo{\stag{838-10k.1} Hypothesis}  ${\frak u}$ is a nice
$\partial$-construction framework (so $\partial$ is regular
uncountable) and $\tau$ is a ${\frak u}$-sub-vocabulary.
\enddemo
\bigskip

\definition{\stag{838-10k.3} Definition}  We call a model $M \in 
{\frak K}_{\frak u}$ standard if $M \in K^\circ_{\frak u} := 
\{M \in K_{\frak u}$: every member of $M$ is an ordinal $<
\partial^+\}$ and ${\frak K}^\circ_{\frak u} = (K^\circ_{\frak u},
\le_{\frak K} \restriction K^o_{\frak u})$.
\enddefinition
\bn
\ub{Convention}:  Models will be standard in this section if not said 
otherwise.
\bigskip

\definition{\stag{838-10k.5} Definition}  1) Let $K^{\text{rt}}_\partial
=K^{\text{rt}}_{\frak u}$ be the class of quadruples $(\bar M,\bar{\bold
J},\bold f,\bar{\Bbb F})$ such that:
\mr
\item "{$(A)$}"  $(\bar M,\bar{\bold J},\bold f) 
\in K^{\text{qt}}_{\frak u}$ recalling \scite{838-1a.29}, and $M_\partial =
\cup\{M_\alpha:\alpha < \partial\}$ has universe some ordinal $<
\partial^+$ divisible by $\partial$ hence $M_\alpha$ is standard for
$\alpha < \partial$
\sn
\item "{$(B)$}"  $\bar{\Bbb F} = \langle \Bbb F_\alpha:\alpha < 
\partial \rangle$ where $\Bbb F_\alpha$ is a
${\frak u}$-amalgamation choice function, see part (2) below
and\footnote{can demand this always} if $2^\partial = \partial^+$ then
each $\Bbb F_\alpha$ has strong uniqueness, see Definition
\scite{838-10k.7}(2) below.
\ermn
2)  We say that $\Bbb F$ is a ${\frak u}$-amalgamation function 
\ub{when}:
\mr
\item "{$(a)$}"  Dom$(\Bbb F) \subseteq \{(M_0,M_1,M_2,\bold J_1,\bold
J_2,A):M_\ell \in {\frak K}^\circ_{\frak u}$ for $\ell \le 2,M_0 \le_{\frak K}
M_\ell,(M_0,M_\ell,\bold J_\ell) \in \text{ FR}_\ell$ for $\ell=1,2$ 
and $M_1 \cap M_2 = M_0$ and $M_1 \cup M_2 \subseteq A \subseteq
\partial^+$, and $|A \backslash M_1 \backslash M_2| < \partial\}$ 
\sn
\item "{$(b)$}"  if $\Bbb F(M_0,M_1,M_2,\bold J_1,\bold J_2,A)$ is well
defined \ub{then} it has the form $(M_3,\bold J^+_1,\bold J^+_2)$ such that
{\roster
\itemitem{ $(\alpha)$ }  $M_\ell \le_{\frak K} M_3 \in 
K^\circ_{< \partial}$ for $\ell=1,2$
\sn
\itemitem{ $(\beta)$ }  $|M_3| = A$ 
\sn
\itemitem{ $(\gamma)$ }  $(M_0,M_\ell,\bold J_\ell) \le^\ell_{\frak u} 
(M_{3 - \ell},M_3,\bold J^+_\ell) \text{ for } \ell =1,2$.
\endroster}
\item "{$(c)$}"    if $(M_0,M_1,M_2,\bold J_1,\bold J_2)$ are as in
clause (a) then for some $A$ we have:\nl $\Bbb F(M_0,M_1,M_2,\bold
J_1,\bold J_2,A)$ is well defined and for any such $A$, 
the set $A \backslash M_1
\backslash M_2$ is disjoint to sup$\{\gamma +1:
\gamma \in M_1$ or $\gamma \in M_2\}$ 
\sn
\item "{$(d)$}"    if\footnote{Dropping clause (d) causes little change}
$\Bbb F(M_0,M_1,M_2,\bold J_1,\bold J_2,A^k)$ is
well defined for $k=1,2$ then $|A^1 \backslash M_1 \backslash M_2|
= |A^2 \backslash M_1 \backslash M_2|$
\nl
moreover
\sn
\item "{$(e)$}"    if $\Bbb F(M_0,M_1,M_2,\bold J_1,\bold J_2,A^1)$ is
well defined and $M_1 \cup M_2 \subseteq A^2 \subseteq \partial^+$ and
otp$(A^2 \backslash M_1 \backslash M_2) = \text{ otp}
(A^1 \backslash M_1 \backslash M_2)$ then also $\Bbb
F(M_0,M_1,M_2,\bold J_1, \bold J_2,A^2)$ is well defined.
\ermn
3) Let $(\bar M^1,\bar{\bold J}^1,\bold f^1,\bar{\Bbb F}^1)
<^{\text{at}}_{\frak u} (\bar M^2,\bar{\bold J}^2,
\bold f^2,\bar{\Bbb F}^2)$ with at standing for atomic, hold 
\ub{when} both quadruples are from $K^{\text{rt}}_{\frak u}$ and 
there are a club $E$ of $\partial$ and sequence 
$\bar{\bold I} = \langle \bold I_\alpha:\alpha < \partial
\rangle$ witnessing it which means that we have
\mr
\item "{$(a)$}"  $\delta \in E \Rightarrow \bold f^1(\delta) \le \bold
f^2(\delta) \and \text{ Min}(E \backslash (\delta +1)) > \bold
f^2(\delta)$
\sn
\item "{$(b)$}"  for $\delta \in E$, if $i \le \bold f^1(\delta)$ then
$M^1_{\delta +i} \le_{{\frak K}_{< \partial}} M^2_{\delta +i}$ and 
if $i < \bold f^i(\delta)$ then $\Bbb F^1_{\delta +i} = \Bbb F^2_{\delta +i}$
\sn
\item "{$(c)$}"   $\langle (M^0_\alpha,M^1_\alpha,\bold I_\alpha):
\alpha \in \cup\{[\delta,\delta +
\bold f^1(\delta)]:\delta \in E\} \rangle$ is $\le^1_{\frak u}$-increasing
continuous.
\sn
\item "{$(d)$}"  if $\delta \in E$ and $i < \bold f^1(\delta)$ and $A$
is the universe of $M^2_{\delta + i+1}$ \ub{then}
\nl
$(M^2_{\delta +i+1},\bold I^2_{\delta+i+1},\bold J^2_{\delta +i}) =
\Bbb F^1_{\delta +i}(M^1_{\delta+i},M^1_{\delta+i+1},
M^2_{\delta+i},\bold I_{\delta+i},\bold J^1_{\delta+i},A)$.
\ermn
4) We say that $(\bar M^\delta,\bar{\bold J}^\delta,\bold
f^\delta,\bar{\Bbb F}^\delta)$ is a canonical upper bound of 
$\langle (\bar M^\alpha,\bar{\bold J}^\alpha,\bold f^\alpha,\bar{\Bbb F}):
\alpha < \delta \rangle$ as in Definition \scite{838-1a.29}(4) adding: in
clause (c), case 1 subclause $(\gamma)$ to the conclusion $\Bbb
F^\delta_{\zeta +i} = \Bbb F^{\alpha_\xi}_{\zeta +i}$ (and similarly in
case 2).
\nl
4A) We say $(\langle \bar M^\alpha,\bar{\bold J}^\alpha,
\bold f^\alpha,\bar{\Bbb F}^\alpha):\alpha < \alpha(*) \rangle$ is a   
$<^{\text{at}}_{\frak u}$-tower \ub{if}:
\mr
\item "{$(a)$}"  $(\bar M^\alpha,\bar{\bold J}^\alpha,\bold
f^\alpha,\bar{\Bbb F}^\alpha) <^{\text{at}}_{\frak u} (\bar M^{\alpha
+1},\bold J^{\alpha +1},\bar{\Bbb F}^{\alpha +1})$ for when $\alpha +1
< \alpha(*)$
\sn
\item "{$(b)$}"  if $\delta < \alpha(*)$ is a limit ordinal, \ub{then}
$(\bar M^\delta,\bar{\bold J}^\delta,\bold f^\delta,
\bar{\Bbb F}^\delta) \in K^{\text{rt}}_{\frak u}$ is a
canonical upper bound of the
sequence $\langle(\bar M^\alpha,\bar{\bold J}^\alpha,\bold
f^\alpha,\Bbb F^\alpha):\alpha < \delta\rangle$. 
\ermn
5) Let $(\bar M,\bar{\bold J},\bold f,\bar{\Bbb F}) \le^{\text{rs}}_{\frak u} 
(\bar M',\bar{\bold J}',\bold f',\bar{\Bbb F}')$ means  that for some \nl
$<^{\text{at}}_{\frak u}$-tower $\langle (\bar M^\alpha,\bar{\bold J}^\alpha,
\bold f^\alpha,\bar{\Bbb F}^\alpha):\alpha \le \alpha(*) \rangle$
we have $(\bar M,\bar{\bold J},\bold f,\bar{\Bbb F}) = 
(\bar M^0_{\alpha(*)},\bar{\bold J}^0,\bold f^0,\bar{\Bbb F}^0)$ 
and $(\bar M',\bar{\bold J}',\bold f',\bar{\Bbb F}') =
(\bar M^{\alpha(*)},\bar{\bold J}^{\alpha(*)},\bold f,
\bar{\Bbb F}^{\alpha(*)})$. 
\nl
6) We say that the sequence $\langle (\bar M^\alpha,\bar{\bold
J}^\alpha,\bold f^\alpha,\bar{\Bbb F}^\alpha):\alpha <
\alpha(*))\rangle$ is $\le^{\text{rs}}_{\frak u}$-increasing
continuous if it is $\le^{\text{rs}}_{\frak u}$-increasing and for any
limit $\delta < \alpha(*)$ the tuple $(\bar M^\delta,\bar{\bold
J}^\delta,\bold f^\delta,\bar{\Bbb F}^\delta)$ is a canonical upper
bound of the sequence of $\langle(\bar M^\alpha,\bar{\bold
J}^\alpha,\bold f^\alpha,\bar{\Bbb F}^\alpha):\alpha < \delta\rangle$.
\enddefinition
\bigskip

\definition{\stag{838-10k.7} Definition}  For $\Bbb F$ a ${\frak
u}$-amalgamation choice function, see Definition \scite{838-10k.5}(2):
\nl
1) $\Bbb F$ has uniqueness \ub{when}:
\mr
\item "{$\circledast$}"  if $\Bbb F(M^\ell_0,M^\ell_1,M^\ell_2,\bold
J^\ell_1,\bold I^\ell_1,A^\ell) = (M^\ell_3,\bold J^\ell_2,\bold I^\ell_2)$
for $\ell=1,2$ (all models standard) and $f$ is a one to one function
from $A^1$ onto $A^2$ preserving the order of the ordinals, $f$ maps
$M^1_i$ to $M^2_i$ for $i=0,1,2$ (i.e. $f \restriction M^1_i$ is an
isomorphism from $M^1_i$ onto $M^2_i$) and $\bold J^1_1,\bold I^1_1$
onto $\bold J^2_1,\bold I^2_1$, respectively, \ub{then} $f$ is an
isomorphism from $M^1_3$ onto $M^2_3$ mapping $\bold J^1_2,\bold
I^1_2$ onto $\bold J^2_2,\bold I^2_2$ respectively.
\ermn
2) $\Bbb F$ has strong uniqueness \ub{when} $\Bbb F$ has uniqueness and
\mr
\item "{$\circledast$}"  if $\Bbb F(M^\ell_0,M^\ell_1,M^\ell_2,\bold
J^\ell_1,\bold I^\ell_1,A^\ell) = (M^\ell_3,\bold J^\ell_2,\bold I^\ell_2)$
for $\ell=1,2$ and $f$ is a one to one mapping from $M^1_1  \cup
M^1_2$ onto $M^2_1 \cup M^2_2$ such that:
$f \restriction M^1_i$ is an isomorphism from $M^1_i$ onto $M^2_i$ for
$i=0,1,2$ and it maps $\bold J^1_1,\bold I^1_1$ onto 
$\bold J^2_1,\bold I^2_1$, respectively, \ub{then}
$|A^1 \backslash M^1_1 \backslash M^1_2| = |A^2 \backslash M^2_1 \backslash
M^2_2|$, and there is an isomorphism $g$ from $M^1_3$ onto
$M^2_3$ extending $f$ and mapping $\bold J^1_2,\bold I^1_2$ onto 
$\bold J^2_2,\bold I^2_2$ respectively; and moreover, otp$(A^1
\backslash M^1_1,M^1_2) = \text{ otp}(A^2 \backslash M^2_1 \backslash
M^2_2)$ and $f \rest (A  \backslash M^1_1 \backslash M^1_2)$ is order
preserving. 
\endroster
\enddefinition
\bigskip

\remark{\stag{838-10k.8} Remark}  1) In Definition \scite{838-10k.5}, in part (2)
we can replace clause (e) by demanding 
$A$ is an interval of the form $[\gamma_*,\gamma_* + \theta]$
where $i < \partial,\gamma_* = \cup\{\gamma +1:\gamma \in M_1$ or
$\gamma \in M_2\}$.  Then in part (3) we have $M^1_\partial$ has
universe $\delta$ for some $\delta < \partial^+$ and $M^2_\gamma$ has universe
$\delta + \partial$.  Also the results in \scite{838-10k.11},
\scite{838-10k.13} becomes somewhat more explicit.
\nl
2) We can fix $\Bbb F_*$, i.e. demand $\Bbb F_\alpha = \Bbb F_*$ in
   Definition \scite{838-10k.5}.
\endremark
\bigskip

\proclaim{\stag{838-10k.9} Claim}  1) There is a ${\frak u}$-amalgamation
function with strong uniqueness.
\nl
2) $K^{\text{rt}}_{\frak u}$ is non-empty, moreover for any stationary $S
\subseteq \partial$ and triple $(\bar M,\bar{\bold J},\bold f) \in
K^{\text{qt}}_{\frak u}$, there is $\bar{\Bbb F}$ such that
$(\bar M,\bar{\bold J},\bold f,\bar{\Bbb F}) \in K^{\text{rt}}_{\frak
u}$ with $S = \{\delta < \partial:\bold f(\delta) > 0\}$.
\nl
3) If $(\bar M,\bar{\bold J}^1,\bold f^1,\bar{\Bbb F}^1) 
\in K^{\text{rt}}_{\frak u}$
\ub{then} for some $(\bar M^2,\bar{\bold J}^2,\bold f^2,\bar{\Bbb F}^2) \in
K^{\text{rt}}_{\frak u}$ we have $(\bar M^1,\bar{\bold J}^1,\bold
f^1,\bar{\Bbb F}^1)$ $
<^{\text{at}}_{\frak u} (\bar M^2,\bar{\bold J}^2,\bold f^2,\bar{\Bbb F}^2)$;
moreover, if $\delta$ is the universe of $M^1_\partial$ then $\alpha <
\partial \Rightarrow M^2_\alpha \backslash \delta \in \{[\delta,\delta
+i):i < \partial\}$. 
\nl
4) Canonical upper bound as in \scite{838-10k.5}(4) exists.
\endproclaim
\bigskip

\demo{Proof}  1) Let $\bold X$ be the set of quintuples $\bold x =
(M^{\bold x}_0,M^{\bold x}_1,M^{\bold x}_2,\bold J^{\bold x}_1,\bold
J^{\bold x}_2)$ as in clause (a) of Definition  \scite{838-10k.5}(2).
We define a two-place relation $E$ on $\bold X:\bold x E \bold y$ iff
$\bold x,\bold y \in \bold X$ and there is a one-to-one function $f$ from
$M^{\bold y}_1 \cup M^{\bold y}_2$ onto $M^{\bold x}_1 \cup M^{\bold
x}_2$ such that $f \restriction M^{\bold x}_\ell$ is an isomorphism
from $M^{\bold x}_\ell$ onto $M^{\bold y}_\ell$ for $\ell=0,1,2$ and
$f$ maps $\bold J^{\bold y}_\ell$ onto $\bold J^{\bold x}_\ell$ for
$\ell=1,2$. 
Clearly $E$ is an equivalence relation, and let $\bold Y \subseteq
\bold X$ be a set of representatives and for every $\bold x \in \bold
X$ let $\bold y(\bold x)$ be the unique $\bold y \in \bold Y$ which is
$E$-equivalent to it and let $f = f_{\bold x}$ be a one-to-one
function from $M^{\bold y}_1 \cup
M^{\bold y}_2$ onto $M^{\bold x}_1 \cup M^{\bold x}_2$ witnessing the
equivalence. 

For each $\bold y \in \bold Y$ by clause (F) of Definitin \scite{838-1a.3}
there is a triple $(M^{\bold y}_3,\bold J^{+,\bold y}_1,\bold
J^{+,\bold y}_2)$ such that
\mr
\item "{$(*)$}"  $(M^{\bold x}_0,M^{\bold x}_\ell,\bold J^{\bold
x}_\ell) \le^\ell_{\frak u} (M^{\bold y}_{3 - \ell},M^{\bold y}_3,
\bold J^{+,\bold y}_\ell)$ for $\ell=1,2$.
\ermn
Without loss of generality $M^{\bold y}_3$ has universe $\subseteq
\partial^+$, we can add has universe $M^{\bold y}_1 \cup M^{\bold
y}_2\cup [\gamma,\gamma + \theta)$ where $\gamma = \sup\{\alpha
+1:\alpha \in M^{\bold y}_1 \cup M^{\bold y}_2\}$ where $\theta$ is
the cardinality of 
$M^{\bold y}_3 \backslash M^{\bold y}_1 \backslash M^{\bold y}_2$.
 
Lastly, let us define $\Bbb F$ as follows:
\mr
\item "{$(a)$}"  $\zeta_{\bold x} = \text{otp}(M^{\bold y(\bold x)}_3
\backslash M^{\bold y(\bold x)}_1 \backslash M^{\bold y(\bold x)}_2) 
< \partial$
\nl
and
\item "{$(b)$}"  Dom$(\Bbb F) = \{(\bold x,A):\bold x \in \bold
X,M^{\bold x}_1 \cup M^{\bold x}_2 \subseteq A \subseteq \partial$ and
otp$(A \backslash M^{\bold x}_1 \backslash M^{\bold x}_2) =
\zeta_{\bold x}\}$ where $(\bold x,A)$ means $(M^{\bold x}_0,M^{\bold
x}_1,M^{\bold x}_2,\bold J^{\bold x}_2,\bold J^{\bold x}_2,A)$
\sn
\item "{$(c)$}"  for $(\bold x,A) \in \text{ Dom}(\Bbb F)$ let
$f_{\bold x,A}$ be the unique one to one function from $M^{\bold
y(\bold x)}_3$ onto $A$ which extends $f_{\bold x}$ and is order
preserving mapping from $M^{\bold y(\bold x)}_3 \backslash M^{\bold
y(\bold x)}_1 \backslash M^{\bold y(\bold x)}_2$ onto $A \backslash M^{\bold
x}_1 \backslash M^{\bold x}_2$
\sn
\item "{$(d)$}"  for $(\bold x,A) \in \text{ Dom}(\Bbb F)$ let $\Bbb
F(\bold x,A)$ be the image under $f_{\bold x,A}$ of $(M^{\bold y(\bold
x)}_3,\bold J^{+,\bold y(\bold x)}_1,\bold J^{+,\bold y(\bold x)}_2)$.
\ermn
Now check.
\nl
2) By part (1) and 
``$K^{\text{qt}}_{\frak u} \ne \emptyset$", see \scite{838-1a.37}(1).
\nl
3) Put together the proof of \scite{838-1a.37}(3) and part (2).
\nl
4) As in the proof of \scite{838-1a.37}(4).  \hfill$\square_{\scite{838-10k.9}}$.
\enddemo
\bigskip

\proclaim{\stag{838-10k.11} Claim}  1) There is a function $\bold m$
satisfying:
\mr
\item "{$\circledast_1$}"  if $(\bar M,\bar{\bold J},\bold f,\bar{\Bbb
F}) \le^{\text{rs}}_{\frak u} (\bar M',\bar{\bold J}',\bold
f',\bar{\Bbb F}')$, recalling \scite{838-10k.5}(5)
 \ub{then} for some function $h:\partial \rightarrow
{\Cal H}_{< \partial}(\partial^+)$, but if $2^\partial = \partial^+$
then $h:\partial \rightarrow {\Cal H}_{< \partial}(\partial)$ we have:
\sn
\item "{$\odot$}"  for a club of $\delta < \partial$ the object
$\bold m(h \restriction \delta,\bar M,\bar{\bold J},
\bold f,\bar{\Bbb F},M'_\delta)$ is a model $N \in K_{\frak u}$ such that
{\roster
\itemitem{ $(a)$ }  $M_{\delta + \bold f(\delta)} \le_{\frak u} N$
\sn
\itemitem{ $(b)$ }  $M'_\delta \le_{\frak u} N$
\sn
\itemitem{ $(c)$ }  there is $N' \le_u M'_\partial$ isomorphic to $N$
over $M_{\delta + \bold f(\delta)} + M'_\delta$ in fact $N' = M'_{\delta +
\bold f(\delta)}$
\endroster}
\item "{$\circledast_2$}"  $\bold m$ is preserved by partial, order
preserving functions from $\partial^+$ to $\partial^+$ compatible with
{\rm id}$_{M_\partial}$
\sn
\item "{$\circledast_3$}"  in fact in $\odot$ above, $\bold m(h \rest
\delta,\bar M,\bar{\bold J},\bold f,\bar{\Bbb F},M'_\delta)$ is
actually $\bold m(h \rest \delta,\bar M \rest [\delta,\delta + 
\bold f(\delta)+1,\bar{\bold J} \rest [\delta,\delta + \bold
f(\delta)),\bar{\Bbb F} \rest [\delta,\delta + \bold
f(\delta)],M'_\delta)$.
\endroster
\endproclaim
\bigskip

\demo{Proof}  By \scite{838-10k.13} we can define $\bold m$ explicitly.
\hfill$\square_{\scite{838-10k.11}}$ 
\enddemo
\bigskip

\proclaim{\stag{838-10k.13} Claim}  $(\bar M',\bar{\bold J}',\bold f',
\bar{\Bbb F}') \le^{\text{rs}}_{\frak u} (\bar M'',\bar{\bold J}'',\bold
f'',\bar{\Bbb F}'')$ \ub{iff} both are from $K^{\text{rt}}_{\frak
u}$ and we can find a $E,\alpha(*) = \alpha_*,\bar u,\bar{\bold d}$ such that:
\mr
\item "{$\circledast$}"   $(a) \quad E$ is a club of $\partial$
\sn
\item "{${{}}$}"  $(b) \quad \alpha_*$ an ordinal $< \partial^+$
\sn
\item "{${{}}$}"  $(c) \quad \bar u$ is a $\subseteq$-increasing 
continuous sequence $\langle u_i:i < \partial \rangle$
\sn
\item "{${{}}$}"  $(d) \quad i < \partial  \Rightarrow |u_i| <
\partial$ and $\alpha_* +1 = \cup\{\alpha_i:i < \partial\}$ and $\alpha
\in u_0,0 \in u_0$ and 
\nl

\hskip25pt $(\forall \beta < \alpha)(\beta \in u_i \equiv
\beta +1 \in u_i)$
\sn
\item "{${{}}$}"  $(e) \quad \bar{\bold d} = \langle \bold d_\delta:\delta \in
E\rangle$
\sn
\item "{${{}}$}"  $(f) \quad \bold d_\delta$ is a ${\frak u}$-free
$(\bold f(\delta),\text{\rm otp}(u_\delta))$-rectangle, see 
Definition \scite{838-1a.7}
\sn
\item "{${{}}$}"   $(g) \quad$ there is 
a $\le^{\text{at}}_{\frak u}$-tower
$\langle (\bar M^\alpha,\bar{\bold J}^\alpha,\bold
f^\alpha,\bar{\Bbb F}^\alpha):\alpha \le \alpha_*\rangle$ 
as in \scite{838-10k.5}(6) 
\nl

\hskip25pt witnessing
the assumption and $\langle \bar{\bold I}^\alpha:\alpha < \alpha_*
\rangle$ with $\bold I^\alpha$ witnessing
\nl

\hskip25pt $(\bar M^\alpha,\bar{\bold J}^\alpha,
\bold f^\alpha,\bar{\Bbb F}^\alpha)
<^{\text{at}}_{{\frak u}_*} (\bar M^{\alpha +1},\bar{\bold J}^{\alpha
+1},\bold f^{\alpha +1},\bar{\Bbb F}^{\alpha +1})$
\nl

\hskip25pt (so $(\bar M^0,\bar{\bold J}^0,\bold f^0,
\bar{\Bbb F}^0) = (\bar M',\bar{\bold J}',\bold f',\bar{\Bbb F}'),
(\bar M^{\alpha(*)},\bar{\bold J}^{\alpha(*)},\bold
f^{\alpha(*)},\bar{\Bbb F}^{\alpha(*)}) =$
\nl

\hskip25pt $(\bar M'',\bar{\bold J}'',\bold f'',\bar{\Bbb F}''))$ such that
{\roster
\itemitem{ ${{}}$ }  $\boxdot \quad$ if $\beta \in [0,\alpha_* +1),\delta
\in E$ and $\beta \in u_\delta$ and $j = \text{\rm otp}(\beta \cap
u_\delta)$ \ub{then}
\nl

\hskip50pt  for every $i \le \bold f(\delta)$
\sn
\itemitem{ ${{}}$ }  \hskip20pt $(\alpha) \quad 
M^{\bold d_\delta}_{i,j} = M^\beta_{\delta +i}$
\sn
\itemitem{ ${{}}$ }  \hskip20pt $(\beta) \quad \bold J^{\bold
d_\delta}_{i,0} = \bold J'_{\delta +i}$ when $i < \bold f(\delta)$
\sn
\itemitem{ ${{}}$ }  \hskip20pt $(\gamma) \quad \bold I^{\bold
d_\delta}_{0,j} = \bold I^\beta_\delta$ when $\beta < \alpha_*$
\sn
\itemitem{ ${{}}$ }  \hskip20pt $(\delta) \quad$ if $i < \bold
f(\delta)$ and $\beta < \alpha_*$ (so {\rm otp}$(\beta \cap u_\delta) <
\text{\rm otp}(u_\delta))$ then
\nl

\hskip50pt  $(M^{\bold d_\delta}_{i+1,j+1},\bold
I^{\bold d_\delta}_{i,j+1},\bold J^{\bold d_\delta}_{i+1,j}) = \Bbb
F_{\delta +i}(M^{\bold d_\delta}_{i,j},M^{\bold
d_\delta}_{i+1,j},M^{\bold d_\delta}_{i,j+1}$,
\nl

\hskip50pt $\bold I^{\bold d_\delta}_{i,j},\bold J^{\bold d_\delta}_{i,j},
|M^{\bold d_\delta}_{i+1,j+1}|)$.
\endroster}
\endroster
\endproclaim
\bigskip

\demo{Proof}  Straight.  \hfill$\square_{\scite{838-10k.13}}$
\enddemo
\bigskip

\remark{\stag{838-10k.14} Remark}  If we define the version of \scite{838-10k.5}(3) with
$|M^2| = |M^1| + \partial$ \ub{then} $\bold m(-)$ is O.K. not only up
to isomorphism but really given the value.
\endremark
\bigskip

\proclaim{\stag{838-10k.17} Claim}  Theorem \scite{838-2b.3} holds.

That is, $\dot I_\tau(\partial^+,K^{{\frak u},{\frak h}}_{\partial^+}) \ge
\mu_{\text{unif}}(\partial^+,2^\partial)$ \ub{when}:
\mr
\item "{$\circledast$}"  $(a) \quad 
2^\theta = 2^{< \partial} < 2^\partial$
\sn
\item "{${{}}$}"  $(b) \quad 2^\partial < 2^{\partial^+}$
\sn
\item "{${{}}$}"  $(c) \quad$ the ideal
{\rm WDmId}$(\partial)$ is not $\partial^+$-saturated
\sn
\item "{${{}}$}"  $(d) \quad {\frak u}$ has the weak $\tau$-coding,
see Definition \scite{838-2b.1}(5) (or just above
\nl

\hskip25pt  some triple $(\bar M^*,\bar{\bold J}^*,\bold f^*) \in
K^{\text{qt}}_{\frak u}$ with {\rm WDmId}$(\partial) 
\restriction (\bold f^*)^{-1}\{0\}$ not 
\nl

\hskip25pt $\partial^+$-saturated)
\sn
\item "{${{}}$}"  $(e) \quad {\frak h}$ is ${\frak u}-\{0,2\}$-appropriate.
\endroster
\endproclaim
\bigskip

\remark{\stag{838-7f.56R} Remark}  1) Similarly \scite{838-2b.7} holds.
\nl
2) We can below (in $\boxtimes$) immitate the proof of \scite{838-3r.3}.
\endremark
\bigskip

\demo{Proof}  Clearly when $(\bar M^*,\bar{\bold J}^*,\bold f^*)$ as
in clause (d) is not given, by \scite{838-10k.9}(2) we can choose it, even
with $\bold f^*$ constantly zero, so without loss of generality  such a triple is given.
By \scite{838-1a.49}(4) and clauses (d) + (e), without loss of generality :
\mr
\item "{$\otimes$}"  ${\frak h} = ({\frak h}_0,{\frak h}_2)$ witness
that $\{0,2\}$-almost every triple $(\bar M,\bar{\bold J},\bold f) \in
K^{\text{qt}}_{\frak u}$ above $(\bar M^*,\bar{\bold J}^*,\bold f^*)$
has the weak coding property.
\ermn
Let $\bar S$ be such that
\mr
\item "{$\odot$}"  $(a) \quad \bar S = \langle S^*_\zeta:\zeta <
\partial^+\rangle$
\sn
\item "{${{}}$}"  $(b) \quad S^*_\zeta \subseteq \partial$
\sn
\item "{${{}}$}"  $(c) \quad S^*_\zeta$ is increasing modulo
$[\partial]^{< \partial}$
\sn
\item "{${{}}$}"  $(d) \quad S^*_0$ and $S^*_{\zeta +1} \backslash S^*_\zeta
\notin \text{ WDmId}(\partial)$
\sn
\item "{${{}}$}"  $(e) \quad \bold f^* \rest (\partial \backslash
S^*_0)$ is constantly zero.
\ermn
Such sequence exists by clause (c) of the hypothesis.  It suffices to
deal with the case ${\frak h}_2$ is a ${\frak u}-2-S^*_0$-appropriate
function, see Definition \scite{838-1a.47}(2A).

We choose $\langle (\bar M^\eta,\bar{\bold J}^\eta,\bold f^\eta,
\bar{\Bbb F}^\eta):\eta \in {}^\gamma(2^\partial)\rangle$ 
by induction on $\gamma < \partial^+$ such that:
\mr
\item "{$\boxplus$}"  $(a) \quad (\bar M^\eta,\bar{\bold
J}^\eta,\bold f^\eta,\bar{\Bbb F}^\eta) \in K^{\text{rt}}_{\frak u}$
and if $\gamma=0$ then
\nl

\hskip25pt $(\bar M^\eta,\bar{\bold J}^\eta,\bold f^\eta) = (\bar
M^*,\bar{\bold J}^*,\bold f^*)$ 
\sn
\item "{${{}}$}"  $(b) \quad$ if $i \notin 
S^*_{\ell g(\eta)}$ then $\bold f_\eta(i)=0$
\sn
\item "{${{}}$}"  $(c) \quad \langle(\bar M^{\eta \restriction
\beta},\bar{\bold J}^{\eta \restriction \beta},\bold f^{\eta
\restriction \beta},\bar{\Bbb F}^{\eta \restriction \beta}:\beta \le
\gamma\rangle$ is $\le^{\text{rs}}_{\frak u}$-increasing continuous,
\sn
\item "{${{}}$}"  $(d) \quad$ if $\eta \in {}^{\gamma +1}(2^\partial)$ and
$\gamma$ is a non-limit ordinal and $\alpha < 2^\partial$ \ub{then}
the pair 
\nl

\hskip25pt $((\bar M^\eta,\bar{\bold J}^\eta,\bold f^\eta),(\bar M^{\eta
\char 94 <\alpha>},\bar{\bold J}^{\eta \char 94 <\gamma>},\bold
f^{\eta \char 94 <\alpha>}))$ strictly $1-S^*_0$-obey ${\frak h}_2$
\nl

\hskip25pt and $0$-obey ${\frak h}_0$ (see Definition
\scite{838-1a.43}(1), \scite{838-1a.47}(2)
\nl

\hskip25pt so without loss of generality  for all 
$\eta \char 94 \langle \alpha \rangle$ we choose the same value
\sn
\item "{${{}}$}"  $(e) \quad \eta \in {}^\gamma(2^\partial)$ and
$\alpha_1 \ne \alpha_2 < 2^\partial,\gamma$ a limit ordinal (even
$\partial$ divides $\gamma$) and
\nl

\hskip25pt  $(\bar M^{\eta \char 94 <\alpha_2>},\bar{\bold
J}^{\eta \char 94 <\alpha_2>},\bar{\Bbb F}^{\eta \char 94 <\alpha_2>})
\le^{\text{rt}}_{\frak u} (M',\bar{\bold J}',\bold f',\bar{\Bbb F}')$
\ub{then} 
\nl

\hskip25pt $M^{\eta \char 94 <\alpha_1>}_\partial$ cannot be $\le_{\frak
K}$-embedded into $M'_\partial$ over $M_\eta$.
\ermn
Why this is enough?  By \scite{838-7f.7}, noting that
\mr
\item "{$(*)$}"  if $\gamma(*) < \partial^+$ and
$\eta \char 94 \langle \alpha_i\rangle
\trianglelefteq \nu_i \in {}^{\gamma(*)}(2^\partial)$ for $i=0,1$ and
$\alpha_0 < \alpha_1 < 2^\partial$
\nl

\hskip25pt  and $f$ is an isomorphism from $M^{\nu_0}_\partial$
onto $M^{\nu_1}_\partial$ over $M^\eta_\partial$ \ub{then}
\nl

\hskip25pt  $\eta,f \restriction M^{\eta \char 94 <\alpha_0>}_\partial,
(\bar M_{\nu_1},\bar{\bold J}_{\nu_1},
\bold f_{\nu_2},\bar{\Bbb F}_{\nu_1})$
form a counterexample 
\nl

\hskip25pt  to clause (e) of $\boxplus$.
\ermn
For $\gamma=0$ clause (a) of $\boxplus$, i.e. choose $(\bar
M^*,\bar{\bold J}^*,\bold f^*)$ recallng our use of \scite{838-10k.9}(2).

For $\gamma$ limit use \scite{838-10k.9}(4).

So assume $\eta \in {}^{\gamma(*)}(2^\partial)$ and $(\bar
M^\eta,\bar{\bold J}^\eta,\bold f^\eta,\bar{\Bbb F}^\eta)$ has been
defined and we should deal with $\Xi_\eta := \{\eta \char 94 \langle
\alpha\rangle:\alpha < 2^\partial\}$.

We choose $(\alpha(0),N_0,\bold I_0)$ such that
\mr
\item "{$\oplus$}"  $(a) \quad \alpha(0) < \partial$
\sn
\item "{${{}}$}"  $(b) \quad (M^\eta_{\alpha(0)},N_0,\bold I_0) \in
\text{ FR}_1$
\sn
\item "{${{}}$}"  $(c) \quad N_0 \cap M^\eta_\partial =
M^\eta_{\alpha(0)}$
\sn
\item "{${{}}$}"  $(d) \quad$ if $\gamma(*)$ is a non-limit ordinal
then $(\alpha(0),N,\bold I)$ is as dictated by ${\frak h}$,
\nl

\hskip25pt  i.e. ${\frak h}_0$, see Definition \scite{838-1a.47}(1)(c)
\sn 
\item "{${{}}$}"  $(e) \quad$ if $\gamma(*)$ is a limit ordinal and
$(\bar M^\eta,\bar{\bold J}^\eta,\bold f^\eta)$ has the weak
coding$_1$
\nl

\hskip25pt  property (see Definition \scite{838-2b.1}(3)) \ub{then} for a
club of $\alpha(1) \in (\alpha(0),\gamma)$
\nl

\hskip25pt  we have:
{\roster
\itemitem{ ${{}}$ }  $\quad (*) \quad$  if
$(M_{\alpha(0)},N_0,\bold I_0) \le^1_{\frak u} 
(M_{\alpha(1)},N^*_1,\bold I^*_1)$
and $M^\eta_\partial \cap N^*_1 = M^\eta_{\alpha(1)}$
\nl

\hskip45pt  \ub{then}
there are $\alpha(2) \in (\alpha(1),\gamma)$ and $N^\ell_2,\bold
I^\ell_2$ for $\ell=1,2$
\nl

\hskip45pt  such that: $(M^\eta_{\alpha(1)},N_1,\bold
I_1) \le^1_{\frak u} (M^\eta_{\alpha(2)},N^\ell_2,\bold I^\ell_2)$ for
$\ell=1,2$ 
\nl

\hskip45pt and $N^1_2,N^2_2$ are $\tau$-incompatible
amalgamations of 
\nl

\hskip45pt $M^\eta_{\alpha(2)},N_1$ over $M^\eta_{\alpha(1)}$.
\endroster}
\ermn
Without loss of generality
\mr
\item "{${{}}$}"   $(f) \quad |N_0| \backslash M^\eta_{\alpha(0)}$ is
an initial segment of $\partial^+ \backslash |M^\eta_\partial|$.
\ermn
We shall use \scite{838-7f.10}.  Toward this we choose $u_i$ and if $i \in
 u_i$ also $E^\eta_\rho,M^\eta_\rho,
\bold I^\eta_\rho,\bold J^\eta_\rho$ (but $\bold J^\eta_\rho$ is 
chosen in the $(i+1)$-th step) for $\rho \in {}^i 2$ induction on $i \in
[\alpha(0),\partial)$ such that:
\mr
\item "{$\boxtimes_1$}"  $(a) \quad u_i \subseteq [\alpha(0),i]$ is
closed
\sn
\item "{${{}}$}"  $(b) \quad$ if $j < i$ then $u_i \cap (j+1) = u_j$
\sn
\item "{${{}}$}"  $(c) \quad E^\eta_\rho$ is a closed subset of $i
\cap u_i$
\sn
\item "{${{}}$}"  $(d) \quad j < \ell g(\rho) = E^\eta_{\rho \restriction
j} = E^\eta_\rho \cap j$
\sn
\item "{${{}}$}"  $(e) \quad$ if $j \in E^\eta_\rho$ then $\bold f_\eta(j)
< \text{ min}\{(E^\eta_\rho \backslash j)$ or $E^\eta_\rho \subseteq
j\}$ and $(j,j + f^\eta(j)] \subseteq u_i$
\sn
\item "{${{}}$}"  $(f) \quad M^\eta_\rho \in K^{\frak u}_{< \partial}$ and
$M^\eta_\rho \cap M^\eta = M^\eta_{\ell g(\rho)}$
\sn
\item "{${{}}$}"  $(g) \quad \langle M^\eta_{\rho \restriction j}:j
\in u_i\rangle$ is $\le_{{\frak K}_{< \partial}}$-increasing
continuous
\sn
\item "{${{}}$}"  $(h) \quad \langle (M^\eta_j,M^\eta_{\rho
\restriction j},\bold I^\eta_\rho):j \in u_i\rangle$ is
$\le^1_{\frak u}$-increasing continuous
\sn
\item "{${{}}$}"  $(i) \quad (\alpha) \quad 
(M^\eta_j,M^\eta_{j+1},\bold J^\eta_j)
\le^2_{\frak u} (M^\eta_{\rho \restriction j},M^\eta_{\rho
\restriction (j+1)},\bold J^\eta_{\rho \restriction j})$ when $j \in i
\cap u_i$ 
\sn
\item "{${{}}$}"  $\qquad (\beta) \quad$ if $j \in 
\cup\{[\zeta,\zeta + \bold f^\eta(\zeta)):
\zeta \in E^\eta_\rho\}$ then moreover we get
\nl

\hskip35pt  $(M^\eta_{\rho \rest
(j+1)},\bold I^\eta_{\rho \rest j},\bold J^\eta_{\rho \rest j})$ by
applying $\Bbb F^\eta_j$
\sn
\item "{${{}}$}"  $(j) \quad$ if $i=\ell g(\rho) = j+1,j$ is limit
$\in S^*_{\gamma(*)+1} \backslash S^*_{\gamma(*)}$ and
$\cup\{E_{\rho \restriction j}:j < \ell g(\rho)\}$ 
\nl

\hskip25pt is unbounded in $\ell g(\rho)$ and  
$\bold f_\eta(\ell g(\rho)) = 0$ and if we can \ub{then}
\nl

\hskip25pt  $u_i = u_j \cup \{i\}$ and 
$(M^\eta_{\rho \char 94 <\ell>},\bold I^\eta_{\rho \char 94 <\ell>},
\bold J^\eta_{\rho \char 94 <\ell>})$ 
\nl

\hskip25pt for $\ell=0,1$ are gotten as in $\oplus$ above so in
particular $M^\eta_{\rho \char 94 <0>},M^\eta_{\rho \char 94 <\ell>}$
\nl

\hskip25pt are $\tau$-incompatible amalgamations of $M^\eta_i,M^\eta_\rho$ over
$M^\eta_j$
\sn
\item "{${{}}$}"  $(k) \quad$ if $i=j+1,\delta = \text{
max}(E^\eta_{\rho \rest j}) \le i,\delta \in S^*_0,j = \delta +
\bold f^\eta(\delta)$, \ub{then} we act
\nl

\hskip25pt  as dictated by ${\frak h}$,
i.e. ${\frak h}_2$; moreover this holds for all the interval
\nl

\hskip25pt $[\delta + \bold f^\eta(\delta),\delta + \bold f^\eta(\delta)
+ i']$ for an appropriate $i' < \partial$
\nl

\hskip25pt  by the ``dictation" of ${\frak h}_2$ 
(see Definition \scite{838-1a.43}).
\ermn
In clause (j) this is possible for enough times, if $(\bar
M^\eta,\bar{\bold J}^\eta,\bold f^\eta)$ has the weak coding$_1$ property,
i.e. for $\rho \in {}^\partial 2$,
for a club of $i \in \bold f^{-1}_\eta\{0\}$ by the choice 
of ${\frak h}$.  Trace the Definitions.
\nl
Also $(\bar M^\eta,\bar{\bold J}^\eta,\bold f^\eta)$ has the weak
coding property by the choice of ${\frak h}$ and the induction hypothesis.

Clearly we can carry the induction on $i < \lambda$ and by
\scite{838-7f.10} carrying the induction to $\gamma(*)+1$, so we have
finished carrying the induction.  So by \scite{838-7f.7} we are done.
\nl
${{}}$   \hfill$\square_{\scite{838-10k.17}}$
\enddemo
\bigskip

\proclaim{\stag{838-10k.19} Claim}  Theorem \scite{838-2b.13} holds.

That is, $\dot I_\tau(\partial^+,K^{\frak u}_{\partial^+}) \ge
\mu_{\text{unif}}(\partial^+,2^\partial)$, \ub{when}:
\mr
\item "{$\circledast$}"  $(a) \quad
2^\theta = 2^{< \partial} < 2^\partial$
\sn
\item "{${{}}$}"  $(b) \quad 2^\partial < 2^{\partial^+}$
\sn
\item "{${{}}$}"  $(c) \quad {\frak u}$ has the vertical
$\tau$-coding$_1$ property above
some triple from $K^{\text{qt}}_{\frak u}$.
\endroster
\endproclaim
\bigskip

\demo{Proof}  Like the proof of \scite{838-10k.17} but:
\mn
\ub{Change (A)}:  We omit $\bar S$, i.e. $\odot$, the choice of $\langle
S^*_\varepsilon:\varepsilon < \partial^+\rangle$, can use $S^*_0 =
\partial$, or more transparently, use $S^*_\zeta = S^*_0$ stationary,
$\partial \backslash S^*_\zeta \in (\text{WDmId}(\partial))^+$, on
$S^*_0$ act as before on $\partial \backslash S^*_0$ and act as on
$S^*_{\zeta +1} \backslash S^*_\zeta$ before.
\mn
\ub{Change (B)}:  We change clause $\boxtimes(j)$ to fit the present
coding so any limit ordinal $j$.  \hfill$\square_{\scite{838-10k.19}}$
\enddemo
\bigskip

\proclaim{\stag{838-10k.21} Claim}  Theorem \scite{838-2b.19} holds.

That is $\dot I_\tau(\partial^+,K^{\frak u}_{\partial^+}) \ge
2^{\partial^+}$ \ub{when}:
\mr
\item "{$(a)$}"  $2^\theta = 2^{< \partial} < 2^\partial$ 
\sn
\item "{$(b)$}"  $2^\partial < 2^{\partial^+}$ 
\sn
\item "{$(c)$}"  ${\frak u}$ has horizontal $\tau$-coding property,
say just above $(\bar M^*,\bar{\bold J}^*,\bold f^*)$
\sn
\item "{$(d)$}"  the ideal {\rm WDmId}$(\partial)$ is not
$\partial^+$-saturated.
\endroster 
\endproclaim
\bigskip

\demo{Proof}  Similar to the proof of \scite{838-10k.17} but:
\mn
\ub{Change (A)}:  We omit $\odot$, i.e. $\bar S^*$ and use $S^*_0 = \partial$
\mn
\ub{Change (B)}:  In $\boxplus$ we use only $\eta \in
{}^{\partial^+} 2$ and clause (e) is changed to:
\mr
\item "{$(e)''$}"  if $(\bar M^{\eta \char 94 <1>},
\bar{\bold J}^{\eta \char 94 <1>},\bold f^{\eta \char 94 <1>},
\bar{\Bbb F}^{\eta \char 94 <1>}) \le^{\text{rt}}_{\frak u} (\bar M',\bar{\bold
J}',\bold f',\bar{\Bbb F}')$ \ub{then} $M^{\eta \char 94
<0>}_\partial$ cannot be $\le_{\frak K}$-embedded into $M'_\partial$
over $M^{<>}_\partial$.
\ermn 
\ub{Change (C)}:  We change $\boxtimes$ clause (j) to deal with the
present coding.
\mn
\ub{Change (D)}:  We use \scite{838-7f.12} rather than \scite{838-7f.7}.
\hfill$\square_{\scite{838-10k.21}}$ 
\enddemo
\bn
\centerline {$* \qquad * \qquad *$}
\bn
\margintag{10k.23}\ub{\stag{838-10k.23} Discussion}:  1) Instead constructing
$\le^{\text{at}}_{\frak u}$-successors $(\bar M^{\eta \char 94 \langle
\alpha\rangle},\bold J^{\eta \char 94 \langle \alpha\rangle},\bold
f^{\eta \char 94 \langle \alpha\rangle})$ of 
$(\bar M^\eta,\bar{\bold J}^\eta,\bold f^\eta)$, we may like
to build, for each $\alpha < 2^\partial$ 
an increasing sequence of length $\zeta$, first
with $\zeta < \partial$ then even $\zeta <\partial^+$ but a sequence of
approximations of height $\partial$.

We would like to have in quite many limit $\delta < \partial$ a ``real
choice" as the various coding properties says.  How
does this help?  If arriving to $\eta \in {}^\delta(2^\partial),\delta
< \partial^+,\eta \char 94 \langle \alpha\rangle$, the model
$M^{\eta \char 94 <\alpha>}_\delta$ is brimmed 
over $M_\delta$; this is certainly
beneficial and having a tower arriving to $\delta$ help toward this.  
But it has a price - we have to preserve it.  In case we
have existence for $K^{3,\text{up}}_{\frak s}$ this occurs, see the
proof of Theorem \scite{838-8h.37} (but was proved in an ad-hoc way).
\nl
2) So we have a function $\iota$ such that; 
so during the construction, for $\eta
   \in {}^{\partial^+}(2^\partial)$ letting $S_\eta := \{\xi < \ell
   g(\eta):\langle(\bar M^{(\eta \rest \xi) \char 94
   <\alpha>},\bar{\bold J}^{(\eta \rest \xi) \char 94 <\alpha>},\bold
   f^{(\eta \rest \xi) \char 94 <\alpha>}):\alpha < 2^\partial\rangle$
   is not constant$\}$ we have:
\mr
\item "{$(a)$}"  if $\xi = \ell g(\eta) = \sup(S_\eta) + 1$ and $\iota
   = \iota(\bar M^\eta \rest (\xi+1),\bar{\bold J}^{\eta \rest
   (\xi+1)},f^{\eta \rest (\xi +1)}),\bar{\Bbb F}^{\eta \rest
   (\xi+1)})$ and $\eta \triangleleft \nu_\ell \in {}^{\ell
   g(\eta)+\iota}(2^\partial)$ for $\ell=1,2$ then 
$(\bar M^{\nu_1},\bar{\bold J}^{\nu_1},\bold f^{\nu_1},
\bar{\Bbb F}^{\nu_1}) = (\bar M^{\nu_2},\bar{\bold J}^{\nu_2},\bold
   f^{\nu_2},\Bbb F^{\nu_2})$.
\endroster
\bn
To formalize this we can use (see a concrete example in the proof in
\scite{838-p.4.13}). 
\definition{\stag{838-10k.25} Definition}  1) We say $\bold i$ is a
$\partial$-parameter \ub{when}:
\mr
\item "{$(a)$}"  $\bold i = (\iota,\bar u)$
\sn
\item "{$(b)$}"  $\iota$ is an ordinal $\ge 1$ but $<
\partial^+$
\sn
\item "{$(c)$}"  $\bar u = \langle u_\varepsilon:\varepsilon < \partial\rangle$
is a $\subseteq$-increasing sequence of subsets of $\iota$ of
cardinality $< \partial$ with union $\iota$
\sn
\item "{$(d)$}"  if $\delta$ is a limit ordinal $< \partial$ then
$u_\delta$ is the closure of $\cup\{u_\varepsilon:\varepsilon < \delta\}$.
\ermn
2) We say $\bold d$ is a ${\frak u}$-free $([\varepsilon_1,\varepsilon_2],\bold
i)$-rectangle \ub{when}: 
\mr
\item "{$(a)$}"  $\varepsilon_1 \le \varepsilon_2 < \partial$
\sn
\item "{$(b)$}"  ${\bold d}$ is a ${\frak u}$-free
$([\varepsilon_1,\varepsilon_2],\iota)$-rectangle 
\ermn
(so $\iota$ may be $\ge \partial$, but then not serious; in fact, it
is an ${\frak u}$-free $([\varepsilon_1,\varepsilon_2],
u_\delta)$-rectangle but we complete it in the obvious way.)
\nl
3) The short case is when $\bold i$ is short, i.e. $\iota = 1,u_\alpha
= 1$.

The long case is $\iota = \partial,u_\varepsilon = \varepsilon +1$.
\enddefinition
\bn
\margintag{10k.27}\ub{\stag{838-10k.27} Discusssion}:  1) Above we have concentrated on what
we may call the ``short" case, the ``long" case as described in
 \scite{838-10k.23}, \scite{838-10k.25}  allows more
constructions by ``consuming" more levels.
\nl
2) Above we can restrict ourselves to the case $\partial = \lambda^+$ so
in \scite{838-10k.5} we then demand on $(A \backslash M_1 \backslash M_2)$
is just ``equal to $\lambda$" and
the possible variants of \scite{838-10k.5}(2),(a) + (b)(B) are irrelevant.
\goodbreak

\head {\S11 Remarks on pcf} \endhead  \resetall \sectno=11
 \spuriousreset
\bn

This section will provide us two pcf claims we use.  One is
\scite{838-pcf.1}, a set-theoretic division into 
cases when $2^\lambda < 2^{\lambda^+}$ (it is from pcf theory; note
that the definition of WDmId$(\lambda)$ is recalled in
\scite{838-0z.5}(4)(b) = \scite{838-0z.5}(4)(b) and of 
$\mu_{\text{wd}}(\lambda)$ is recalled in \scite{838-0z.5}(8) =
\scite{838-0z.5}(8)), we can replace 
$\lambda^+$ by regular $\lambda$ such that $2^\theta = 
2^{< \lambda} < 2^\lambda$ for some $\theta$).  The second deals with
the existence of large independent subfamilies of sets, \scite{838-pcf.11}.
This is a revised version of a part of \cite{Sh:603}.  See on history
related to \scite{838-pcf.1} in \cite{Sh:g} particularly in
\cite[II,5.11]{Sh:g} and \cite{Sh:430}.
\bigskip

\remark{Remark}  Recall that

$$
\align
\text{cov}(\chi,\mu,\theta,\sigma) = \chi + \text{ Min} \bigl\{
|{\Cal P}|:&{\Cal P} \subseteq [\chi]^{< \mu} \text{ and every member of} \\
  &[\chi]^{< \theta} \text{ is included in the union of } < \sigma
\text{ members of } {\Cal P} \bigr\}.
\endalign
$$
\endremark
\bigskip

\proclaim{\stag{838-pcf.1} Claim}   Assume $2^\lambda < 2^{\lambda^+}$.

\ub{Then} one of the following cases occurs: (clauses
  $(\alpha)-(\lambda)$ appear later)
\mr
\item "{$(A)_\lambda$}"  $\chi^* = 2^{\lambda^+}$ and for some $\mu$
clauses $(\alpha) - (\varepsilon)$ hold
\sn
\item "{$(B)_\lambda$}"  for some $\chi^* > 2^\lambda$ and $\mu$
clauses $(\alpha) - (\kappa)$ hold (note: $\mu$ appear only in
$(\alpha) - (\varepsilon))$
\sn
\item "{$(C)_\lambda$}"  $\chi^* = 2^\lambda$ and
clauses $(\eta) - (\mu)$ hold
\nl
\ub{where}
{\roster
\itemitem{ $(\alpha)$ }  $\lambda^+ < \mu \le 2^\lambda$ and {\rm cf}$(\mu) 
= \lambda^+$
\sn
\itemitem{ $(\beta)$ }   {\rm pp}$(\mu) = \chi^*$, moreover 
{\rm pp}$(\mu) =^+ \chi^*$ 
\sn
\itemitem{ $(\gamma)$ }  $(\forall \mu')(\text{\rm cf}(\mu') 
\le \lambda^+ < \mu' < \mu \Rightarrow 
\text{\rm pp}(\mu') < \mu)$ hence \newline
{\rm cf}$(\mu') \le \lambda^+ < \mu' < \mu \Rightarrow 
\text{\rm pp}_{\lambda^+}(\mu') < \mu$
\sn
\itemitem{ $(\delta)$ }  for every regular cardinal $\chi$ in the interval
$(\mu,\chi^*]$ there is an increasing sequence $\langle \lambda_i:
i < \lambda^+ \rangle$ of regular cardinals $> \lambda^+$ with limit $\mu$ 
such that $\chi = \text{\rm tcf } \left( \dsize
\prod_{i < \lambda^+} \lambda_i/J^{\text{bd}}_{\lambda^+} \right)$,
and $i < \lambda^+ \Rightarrow \text{\rm max pcf}\{\lambda_j:j < i\} <
\lambda_i < \mu$
\sn
\itemitem{ $(\varepsilon)$ }  for some regular $\kappa \le \lambda$, for any
$\mu' < \mu$ there is a tree ${\Cal T}$ with $\le \lambda$ 
nodes, $\kappa$ levels and $|\text{\rm lim}_\kappa({\Cal T})| 
\ge \mu'$ (in fact e.g. $\kappa = \text{\rm Min}\{\theta:2^\theta 
\ge \mu\}$ is appropriate; without loss of generality
${\Cal T} \subseteq {}^{\kappa >} \lambda)$
\sn
\itemitem{ $(\zeta)$ }   there is no normal $\lambda^{++}$-saturated ideal
on $\lambda^+$
\sn
\itemitem{ $(\eta)$ }   there is $\langle {\Cal T}_\zeta:
\zeta < \chi^* \rangle$ such that: ${\Cal T}_\zeta \subseteq 
{}^{\lambda^+ >}2$, a subtree of cardinality $\lambda^+$ 
and ${}^{\lambda^+}2 = \{\text{\rm lim}_{\lambda^+}
({\Cal T}_\zeta):\zeta < \chi^*\}$
\sn
\itemitem{ $(\theta)$ }  $\chi^* < 2^{\lambda^+}$ moreover
$\chi^* < \mu_{\text{unif}}(\lambda^+,2^\lambda)$, 
but $< \mu_{\text{unif}}(\lambda^+,2^\lambda)$ is not used here, 
\sn
\itemitem{ $(\iota)$ }  for some 
$\zeta < \chi^*$ we have {\rm lim}$_{\lambda^+}
({\Cal T}_\zeta) \notin$ {\rm UnfmTId}$_{(\chi^*)^+}(\lambda^+)$, not used here
\sn
\itemitem{ $(\kappa)$ }  {\rm cov}$(\chi^*,\lambda^{++},
\lambda^{++},\aleph_1) = \chi^*$ or 
$\chi^* = \lambda^+$, equivalently $\chi^* =  
\sup[\{\text{\rm pp}(\chi):\chi \le 2^\lambda,
\aleph_1 \le \text{\rm cf}(\chi) \le \lambda^+ < \chi\} 
\cup \{\lambda^+\}]$ by \cite[Ch.II,5.4]{Sh:g};
\ub{note} that clause 
$(\kappa)$ trivially follows from $\chi^* = 2^{\lambda^+}$
\sn
\itemitem{ $(\lambda)$ } for no $\mu \in (\lambda^+,2^\lambda]$ do we have 
{\rm cf}$(\mu) \le \lambda^+$, {\rm pp}$(\mu) > 2^\lambda$; equivalently 
$2^\lambda > \lambda^+ \Rightarrow \text{\rm cf}([2^\lambda]^{\lambda^+},
\subseteq) = 2^\lambda$ 
\sn
\itemitem{ $(\mu)$ }  if there is a normal $\lambda^{++}$-saturated
ideal on $\lambda^+$, moreover 
the ideal {\rm WDmId}$(\lambda^+)$ is, \ub{then} 
$2^{\lambda^+} = \lambda^{++}$ (so as $2^\lambda < 2^{\lambda^+}$
clearly $2^\lambda = \lambda^+$).
\endroster}
\endroster
\endproclaim
\bigskip

\demo{Proof}  This is related to \cite[II,5.11]{Sh:g}; we assume basic 
knowledge of pcf (or a readiness to believe quotations).
Note that by their definitions
\mr
\item "{$\circledast_1$}"  if $2^\lambda > \lambda^+$ then 
for any $\theta \in [\aleph_0,\lambda^+]$ we have
cf$([2^\lambda]^{\le \lambda^+},\subseteq) = 2^\lambda
\Leftrightarrow \text{ cov}(2^\lambda,\lambda^{++},\lambda^{++},2) =
2^\lambda \Leftrightarrow \text{
cov}(2^\lambda,\lambda^{++},\lambda^{++},\theta) = 2^\lambda$ .
\ermn
[Why?  Because $(2^\lambda)^{< \lambda^+} =
2^\lambda$ and cf$([2^\lambda]^{\le \lambda^+},\subseteq) = \text{
cov}(2^\lambda,\lambda^{++},\lambda^{++},2)$
\nl
$\ge \text{ cov}(2^\lambda,\lambda^{++},\lambda^{++},\aleph_0) \ge
\text{ cov}(2^\lambda,\lambda^{++},\lambda^{++},\theta) \ge
2^\lambda$ for $\theta \in [\aleph_0,\lambda]$.]

Note also that
\mr
\item "{$\circledast_2$}"  $\lambda^+ \notin 
\text{ WDmId}(\lambda^+)$ and ${}^{\lambda^+}2
\notin \text{ WDmTId}(\lambda^+)$.
\ermn
[Why?  Theorem \scite{838-0z.7}(2) with $\theta,\partial$ there standing
for $\lambda,\lambda^+$ here.] 
\enddemo
\bn
\underbar{Possibility 1}:  $2^\lambda > \lambda^+$ and 
$\text{cov}(2^\lambda,\lambda^{++},
\lambda^{++},\aleph_1) = 2^\lambda$ \ub{or} $2^\lambda = \lambda^+$;
and let $\chi^* := 2^\lambda$.

We shall show that case $(C)$ holds (for the cardinal $\lambda$),
the first assertion ``$\chi^* = 2^\lambda$" holds by our choice.
\newline
Now clause $(\kappa)$ is obvious.
 As for clause $(\eta)$, we have $\chi^* = 2^\lambda < 2^{\lambda^+}$.  Now
if $2^\lambda = \lambda^+$ we let 
${\Cal T}_\zeta = {}^{\lambda^+ >}2$, for $\zeta < \chi^*$ so clause
 $(\eta)$ holds, otherwise as ${}^{\lambda^+>}2$
 has cardinality $2^\lambda$, by the definitions of cov$(2^\lambda,\lambda^{++},\lambda^{++},\aleph_1)$ and the
 possibility assumption (and obvious equivalence) there is ${\Cal P} \subseteq
 [{}^{\lambda^+>} 2]^{\lambda^+}$ of cardinality $\chi^*$ such that
 any $A \in [{}^{\lambda^+>}2]^{\lambda^+}$ is included in the union
 of $\le \aleph_0$ of them.  So $\eta \in {}^{\lambda^+} 2 \Rightarrow
 (\exists A \in {\Cal P})(\exists^{\lambda^+} \alpha < \lambda^+)(\eta
 \rest \alpha \in A)$ so let $\langle A_\zeta:\zeta < \chi^*\rangle$
 list ${\Cal P}$ and let ${\Cal T}_\zeta = \{\eta \rest \alpha:\eta
 \in A_\zeta$ and $\alpha \le \ell g(\eta)\}$, now check that they are
 as required in clause  $(\eta)$.

As on the one hand by \cite[AP,1.16 + 1.19]{Sh:f} or see \scite{838-7f.14} we have
$\left( \mu_{\text{unif}}(\lambda^+,2^\lambda) \right)^{\aleph_0} = 
2^{\lambda^+} > 2^\lambda =
\chi^*$ and on the other hand $(\chi^*)^{\aleph_0} = (2^\lambda)^{\aleph_0} 
= 2^\lambda = \chi^*$ necessarily $\chi^* < \mu_{\text{unif}}
(\lambda^+,2^\lambda)$ 
so clause $(\theta)$ follows; next clause $(\iota)$ follows from clause 
$(\eta)$ by the definition of UnfTId$_{(\chi^*)^+}(\lambda^+)$.  
In fact in our
possibility for some $\zeta$, lim$_{\lambda^+}({\Cal T}_\zeta) \notin 
\text{\rm WDmTId}(\lambda^+)$ because
 WDmTId$(\lambda^+)$ is $(2^\lambda)^+$-complete by 
\scite{838-0z.7}(2),(4) recalling $\circledast_2$ and having 
chosen $\chi^* = 2^\lambda$.  

Now if $2^{\lambda^+} > \lambda^{++}$, (so $2^{\lambda^+}
\ge \lambda^{+3}$), then for some $\zeta < \chi^*,{\Cal T}_\zeta$ is 
(a tree with $\le \lambda^+$ nodes, $\lambda^+$ levels and) at 
least $\lambda^{+3} \,\, \lambda^+$-branches which is well 
known (see e.g. [J]) to imply ``no normal ideal on 
$\lambda^+$ is $\lambda^{++}$-saturated"; so we got clause $(\mu)$.
Also if $2^{\lambda^+} \le \lambda^{++}$ then $2^\lambda =
\lambda^+,2^{\lambda^+} = \lambda^{++}$.

As for clause $(\lambda)$, by the 
definition of $\chi^*$ and the assumption $\chi^* =
2^\lambda$ we have the first two phrases.
The ``equivalently" holds as $(2^\lambda)^{\aleph_0} = 2^\lambda$.
\bn
\underbar{Possibility 2}:  $\chi^* := \text{ cov}(2^\lambda,\lambda^{++},
\lambda^{++},\aleph_1) > 2^\lambda > \lambda^+$.

So (C)$_\lambda$ fails, and we have to show that $(A)_\lambda$ or
$(B)_\lambda$ holds. 

Let
\mr
\item "{$(*)_0$}"   $\mu := \text{ Min}\{\mu:\text{cf}(\mu) 
\le \lambda^+,\lambda^+ < \mu 
\le 2^\lambda$ and pp$(\mu) = \chi^*\}$.
\ermn
We know by \cite[II,5.4]{Sh:g} that 
$\mu$ exists and (by \cite[II,2.3]{Sh:g}(2)) clause $(\gamma)$ holds, also 
$2^\lambda < \text{ pp}(\mu) \le \mu^{\text{cf}(\mu)} \le
(2^\lambda)^{\text{cf}(\mu)} = 2^{\lambda + \text{cf}(\mu)}
$ hence 
cf$(\mu) = \lambda^+$.  So clauses $(\alpha),(\beta),(\gamma)$ hold 
(of course, for clause $(\beta)$ use \cite[Ch.II,5.4]{Sh:g}(2)), 
and by $(\gamma)$ + \cite[VIII,\S1]{Sh:g} also clause $(\delta)$ holds.

Toward trying to prove clause $(\varepsilon)$ let
\mr
\item "{$(*)_1$}"   $\Upsilon := \text{ Min}\{\theta:2^\theta \ge \mu\}$,
\nl
clearly\footnote{Below we show that $\Upsilon > \text{ cf}(\Upsilon)
\Rightarrow \text{ cf}(\Upsilon) > \aleph_0$.}
\sn
\item "{$(*)_2$}"   $\alpha < \Upsilon \Rightarrow 2^{|\alpha|} < \mu$ and
$\Upsilon \le \lambda$ (as $2^\lambda \ge \mu$) hence cf$(\Upsilon) \le
\Upsilon \le \lambda < \lambda^+ = \text{ cf}(\mu)$ hence 
$2^{<\Upsilon} < \mu$.
\ermn
Let
\mr
\item "{$(*)_3$}"  $(a) \quad u$ be a closed unbounded subset of
$\Upsilon$ of order type cf$(\Upsilon)$
\sn
\item "{${{}}$}"  $(b) \quad {\Cal T}^* = (\dbcu_{\alpha \in u}
{}^\alpha 2,\triangleleft)$ is a tree with cf$(\Upsilon)$ levels and
$\le 2^{< \Upsilon}$ nodes.
\endroster
\bn
Now we shall prove clause $(\varepsilon)$, i.e.
\mr
\item "{$(*)_4$}"  there 
is a tree with $\lambda$ nodes, cf$(\Upsilon)$ levels and 
$\ge \mu \,\Upsilon$-branches.
\endroster
\bn
\ub{Case A}:  $\Upsilon$ has cofinality $\aleph_0$.
\bn
In the case $\Upsilon = \aleph_0$ or just $2^{< \Upsilon} \le \lambda$
 clearly there is a tree as required, i.e. ${\Cal T}^*$ is a tree having $\le
 2^{< \Upsilon} \le \lambda$ nodes.
So we can assume $2^{< \Upsilon} > \lambda$ and
$\Upsilon > \text{ cf}(\Upsilon) = \aleph_0$ hence $\langle 2^\theta:\theta <
\Upsilon\rangle$ is not eventually constant. 

So necessarily $(\exists \theta < \Upsilon)(2^\theta \ge
\lambda)$hence $2^{< \Upsilon} > \lambda^+$ (and even $2^{< \Upsilon}
\ge \lambda^{+ \omega})$ and for some $\theta < \Upsilon$ we have
$\lambda^{++} < 2^\theta
< 2^{< \Upsilon} < \mu$.  Let $\chi' = \text{ cov}(2^{<
\Upsilon},\lambda^{++},\lambda^{++},\aleph_1)$, so $\chi' \ge 2^{<
\Upsilon}$ and $\chi' < \mu$ by \cite[II,5.4]{Sh:g} and
clause $(\gamma)$ of \scite{838-pcf.1}
which have been proved (in our present possibility).

We try to apply claim \scite{838-pcf.8} below with
$\aleph_0,\lambda^+,2^{< \Upsilon},\chi'$ here standing for $\theta,\kappa,
\mu,\chi$ there; we have to check the assumptions of
\scite{838-pcf.8} which means $2^{< \Upsilon} > \lambda^+ > \aleph_0$ and $\chi' =
\text{ cov}(2^{< \Upsilon},\lambda^{++},\lambda^{++},\aleph_1)$, both clearly
hold.  So the conclusion of \scite{838-pcf.8} holds which means that
$(\chi')^{\aleph_0} \ge \text{ cov}((2^{<
\Upsilon})^{\aleph_0},(\lambda^{++})^{\aleph_0},\lambda^{++},2)$, now
$(\chi')^{\aleph_0} \le \mu^{\aleph_0} \le (2^\lambda)^{\aleph_0} =
2^\lambda$ and $(2^{< \Upsilon})^{\aleph_0} = 2^\Upsilon \ge \mu$
because presently cf$(\Upsilon) = \aleph_0$ and the choice of
$\Upsilon$ and $(\lambda^{++})^{\aleph_0} \le (2^\theta)^{\aleph_0} = 
2^\theta$.  So by
monotonicity $2^\lambda \ge \text{ cov}
(\mu,2^\theta,\lambda^{++},\aleph_1)$.  But
cov$(2^\theta,\lambda^{++},\lambda^{++},\aleph_1) \le \chi' := \text{
cov}(2^{< \Upsilon},\lambda^{++},\lambda^{++},\aleph_1) < \mu \le
2^\lambda$ by clause $(\gamma)$ which we have proved (and
\cite[ChII,5.4]{Sh:g}) so by transitivity of cov, see \scite{838-pcf.9}(4), 
also $2^\lambda \ge
\text{ cov}(\mu,\lambda^{++},\lambda^{++},\aleph_1)$ contradicting the
present possibility.
\bn
\ub{Case B}:  cf$(\Upsilon) > \aleph_0$.

Let $h:{\Cal T}^* \rightarrow 2^{<
\Upsilon}$ be one-to-one, see $(*)_3(b)$.  Let ${\Cal P} \subseteq
[2^{< \Upsilon}]^{\le \lambda}$ be such that every $X \in 
[2^{< \Upsilon}]^{\le \lambda^+}$ is included in 
the union of countably members of ${\Cal P}$,
exists by clause $(\gamma)$ of \scite{838-pcf.1} by \cite[II.5.4]{Sh:g}.

Now for every $\bar \nu = \langle \nu_\varepsilon:\varepsilon <
\Upsilon\rangle \in \text{ lim}_{\text{cf}(\Upsilon)}({\Cal T}^*)$, for
some $A_{\bar \nu} \in {\Cal P}$ we have $\Upsilon = \sup\{\varepsilon
\in u:h(\nu_\varepsilon) \in A_{\bar \nu}\}$, so for every $\mu' \in
(2^{< \Upsilon},\mu)$ for some $A \in {\Cal P}$ we have $\mu' \le
|\{\bar\nu:\bar\nu \in \text{ lim}_{\text{cf}(\Upsilon)}({\Cal
T}^*)$ and $A_{\bar\nu} = A\}|$.

Now let ${\Cal T}'$ be the closure of $A$ to initial segments of
length $\in u$, easily ${\Cal T}'$ is as required.

So we have proved $(*)_4$ so in possibility (2) the demand $(\alpha) -
(\varepsilon)$ in $(A)_\lambda$ holds.
\bn
\underbar{Sub-possibility 2$\alpha$}:  $\chi^* < 2^{\lambda^+}$.

We shall prove $(B)_\lambda$, so by the above we are left with 
proving clauses $(\zeta)-(\kappa)$ when $\chi^* < 2^{\lambda^+}$.
By the choice of $\chi^*$, easily the demand in
clause $(\zeta)$ (in Case B of
\scite{838-pcf.1}) holds; that is let $\{u_\zeta:\zeta < \chi^*\}$ be a
family of subsets of ${}^{\lambda^+ >} 2$, a set of cardinality
$2^\lambda$, each of cardinality $\lambda^+$ such that any other such
subset is included in the union of $\le \aleph_0 < \aleph_1$ of them,
exist by the choice of $\chi^*$.

Let ${\Cal T}_\zeta = \{\nu \rest i:\nu \in u_\zeta$ and $i \le \ell
g(\nu)\}$.  Now $\langle {\Cal T}^*_\zeta:\zeta < \chi^*\rangle$ is as
required.

In clause $(\eta),``2^\lambda < \chi^*"$ holds 
as we are in possibility 2$\alpha$.
\medskip

Also as pp$(\mu) = \chi^*$ and cf$(\mu) = \lambda^+$  by 
the choice of $\mu$ necessarily 
(by transitivity of pcf, i.e., \cite[Ch.II,2.3]{Sh:g}(2)) we have 
cf$(\chi^*) > 
\lambda^+$ but $\mu > \lambda^+$.  
Easily $\lambda^+ < \chi \le \chi^* \wedge$ cf$(\chi) \le \lambda^+
\Rightarrow$ pp$(\chi) \le \chi^*$ hence cov$(\chi^*,\lambda^{++},
\lambda^{++},\aleph_1) = \chi^*$ by \cite[Ch.II,5.4]{Sh:g}, which gives 
clause $(\lambda)$.
Now let ${\Cal A} \subseteq [\chi^*]^{\lambda^+}$ exemplify
cov$(\chi^*,\lambda^{++},\lambda^{++},\aleph_1) = \chi^*$ and let
${\Cal A}' = \{B:B$ is an infinite countable subset of some $A \in
{\Cal A}\}$.  So ${\Cal A}' \subseteq [\chi^*]^{\aleph_0}$ and easily
$A \in [\chi^*]^{\lambda^+} \Rightarrow (\exists B \in {\Cal A}')(B
\subseteq A)$ and $|{\Cal A}'| \le \chi^*$  as 
$(\lambda^+)^{\aleph_0} \le 2^\lambda < \chi^*$ 
certainly there is no family of $>
\chi^*$ subsets of $\chi^*$ each of cardinality $\lambda^+$ with pairwise
finite intersections.  But by \scite{838-7f.14} there is ${\Cal A}'
\subseteq [\mu_{\text{unif}}(\lambda^{++},2^{\lambda^+})]^{\lambda^+}$
of cardinality $2^{\lambda^{++}}$ such that $A \ne B \in {\Cal A}'
\Rightarrow |A \cap B| < \aleph_0$, hence 
we have $\chi^* < \mu_{\text{unif}}(\lambda^+,2^\lambda)$ thus 
completing the proof of $(\theta)$.

Now clause $(\iota)$ follows by clauses $(\eta) + (\theta) + (\kappa)$
as $\emptyset \notin \text{ UnfTId}_{(\chi^*)^+}
\in (\lambda^+)$ which is $(\chi^*)^+$-complete ideals, see
\scite{838-7f.14}.  Note also that $\emptyset \in 
\text{ WDmTId}_{\chi^*}(\lambda^+)$ by \scite{838-0z.7}(2) and
it is a $(\chi^*)^+$ complete ideal by \scite{838-0z.7}(4).
Also by clause $(\alpha)$ which we have proved 
$2^{\lambda^+} \ne \lambda^{++}$ hence $2^{\lambda^+} \ge \lambda^{+3}$ so by
clause $(\eta)$ (as $\chi^* < 2^{\lambda^+}$), we have
$|\text{lim}_{\lambda^+}({\Cal T}_\zeta)| \ge \lambda^{+3}$ for some $\zeta$ 
which is well known (see \cite{J})
to imply no normal ideal on $\lambda^+$ is $\lambda^{++}$-saturated; 
i.e., clause $(\mu)$.  So we have proved clauses $(\alpha)-(\lambda)$
holds, i.e. that case $(B)_\lambda$ 
holds.
\bn
\underbar{Sub-possibility 2$\beta$}:  $\chi^* = 2^{\lambda^+}$ (and
$\chi^* > 2^\lambda > \lambda^+$).

We have proved that case $(A)_\lambda$ holds, as we already defined $\mu$ 
and $\chi^*$ and proved clauses 
$(\alpha),(\beta),(\gamma),(\delta),(\varepsilon)$ so we are done.
\hfill$\square_{\scite{838-pcf.1}}$ 
\bn
It may be useful to recall (actually $\lambda^{<\kappa>_{\text{tr}}} 
= \lambda$ suffice)
\demo{\stag{838-pcf.9} Fact}  1) Assume $\lambda > \theta \ge \kappa =
\text{ cf}(\kappa) \ge \kappa_1$.  Then
$\lambda^{<\kappa>_{\text{tr}}} \le 
\text{ cov}(\lambda,\theta^+,\kappa^+,\kappa_1)$ 
recalling $\mu^{<\kappa>_{\text{tr}}} =
\sup\{\text{lim}_\kappa({\Cal T}):{\Cal T}$ is a tree with $\le \mu$ nodes
and $\kappa$ levels, e.g. ${\Cal T}$ a subtree of ${}^{\kappa
>}\mu\}$.
\nl
2) If $\mu > \aleph_0$ is strong limit and $\lambda > \mu$ \ub{then} for
   some $\kappa < \mu$ we have cov$(\lambda,\mu,\mu,\kappa)$.
\nl
3) If ${\Cal T} \subseteq {}^{\lambda^+ >}2$ is a tree, 
$|{\Cal T}| \le \lambda^+$ and $\lambda \ge 
\beth_\omega$ \underbar{then} for every
regular $\kappa < \beth_\omega$ large enough, we can find $\langle Y_\delta:
\delta < \lambda^+,\text{cf}(\delta) = \kappa \rangle,|Y_\delta| \le
\lambda$ such that: 
\newline
for every $\eta \in \lim_{\lambda^+}({\Cal T})$ for a club
of $\delta < \lambda^+$ we have 
cf$(\delta) = \kappa \Rightarrow \eta \restriction \delta \in
Y_\delta$.
\nl
4) [Transitivity of cov]  If $\mu_3 \ge \mu_2 \ge \mu_1 \ge \theta \ge
   \sigma = \text{ cf}(\sigma)$ and $\lambda_2 = \text{
   cov}(\mu_3,\mu_2,\theta,\sigma)$ and $\mu < \mu_2 \Rightarrow
   \lambda_1 \ge \text{ cov}(\mu,\mu_1,\theta,\sigma)$ \ub{then}
   $\lambda_1 + \lambda_2 \ge \text{ cov}(\mu_3,\mu_1,\theta,\sigma)$.
\enddemo
\bigskip

\demo{Proof}  1) E.g. proved inside \scite{838-pcf.1}.
\nl
2) By \cite{Sh:460} or see \cite{Sh:829}.
\nl
3) Should be clear from part (2).
\nl
4) Let ${\Cal P}_2 \subseteq [\mu_3]^{< \mu_2}$ exemplify $\lambda_2 =
   \text{ cov}(\mu_3,\mu_2,\theta,\sigma)$ and for each $A \in P_2$
   let $f_A$ be one to one function from $|A|$ onto $A$ and ${\Cal
   P}_{1,A} \subseteq [|A|]^{< \mu_1}$ exemplify $\lambda_1 \ge \text{
   cov}(|A|,\mu_1,\theta,\sigma)$.  Lastly, let ${\bold P} =
   \{f_A(\alpha):\alpha \in u\}:u \in \bold P_{1,A},A \in \bold P_2\}$
   it exemplify $\lambda_1 + \lambda_2 \ge \text{
   cov}(\mu_3,\mu_1,\theta,\sigma)$ as required. 
 \hfill$\square_{\scite{838-pcf.8}}$
\enddemo
\bn
We have used in proving \scite{838-pcf.1} also
\demo{\stag{838-pcf.8} Observation}  Assume $\mu > \kappa > \theta$.

If $\chi = \text{\rm cov}(\mu,\kappa^+,\kappa^+,\theta^+)$ \ub{then}
$\chi^\theta \ge \text{\rm cov}(\mu^\theta,(\kappa^\theta)^+,\kappa^+,2)$.
\enddemo
\bigskip

\demo{Proof}  Let ${\Cal P} \subseteq [\mu]^\kappa$ exemplify $\chi =
\text{ cov}(\mu,\kappa^+,\kappa^+,\theta^+)$.

Let $\langle \eta_\alpha:\alpha < \mu^\theta\rangle$ list ${}^\theta
\mu$.  Now for ${\Cal U} \in [\mu]^\kappa$ define ${\Cal U}^{[*]} = \{\alpha <
\mu^\theta$: if $i < \theta$ then $\eta_\alpha(i) \in {\Cal U}\}$, and let
${\Cal P}_1 = {\Cal P},{\Cal P}_2 = \{\dbcu_{i < \theta} A_i:A_i \in
{\Cal P}$ for $i < \theta\}$ and ${\Cal P}_3 = \{{\Cal U}^{[*]}:{\Cal
U} \in {\Cal P}_2\}$.

So (of course $\chi \ge \mu$ as $\mu > \kappa$, hence $\chi > \kappa$)
\mr
\item "{$(*)_1$}"  $(a) \quad {\Cal P}_1 \subseteq [\mu]^\kappa$ has
cardinality $\chi$
\sn
\item "{${{}}$}"  $(b) \quad {\Cal P}_2 \subseteq
[\mu]^{\kappa^\theta}$ has cardinality $\le \chi^\theta$
\sn
\item "{${{}}$}"  $(c) \quad {\Cal P}_3 \subseteq
[\mu^\theta]^{\kappa^\theta}$ has cardinality $\le \chi^\theta +
\kappa^\theta = \chi^\theta$
\ermn
and
\mr
\item "{$(*)_2$}"  if ${\Cal U} \in [\mu^\theta]^{\le \kappa}$ then
{\roster
\itemitem{ $(a)$ }  ${\Cal U}' := \{\eta_\alpha(i):\alpha \in {\Cal
U}$ and $i < \theta\} \in [\mu]^{\le \kappa}$
\sn
\itemitem{ $(b)$ }  there are $A_i \in {\Cal P}_1$ for $i < \theta$ 
such that ${\Cal U}' \subseteq \dbcu_{i < \theta} A_i$
\sn
\itemitem{ $(c)$ }  $\dbcu_{i < \theta} A_i \in {\Cal P}_2$
\sn
\itemitem{ $(d)$ }  ${\Cal U} \subseteq (\dbcu_{i < \theta} A_i)^{[*]} 
\in {\Cal P}_3 \subseteq [\mu^\theta]^{\kappa^\theta}$.
\ermn
[Why?  Clause (a) holds by cardinal arithmetic, clause (b) holds by
the choice of ${\Cal P} = {\Cal P}_1$, clause (c) holds by the
definition of ${\Cal P}_2$ and clause (d) holds by the definition of
$(-)^{[*]}$ and of ${\Cal P}_3$.]
\endroster}
\ermn
Together we are done.  \hfill$\square_{\scite{838-pcf.8}}$ 
\enddemo
\bn
The following is needed when we like to get in the model theory not just
many models but many models no one $\le_{\frak K}$-embeddable into
another and even just for $\dot I$, (see \yCITE[4d.11]{E46}).
\bigskip

\proclaim{\stag{838-pcf.11} Claim}  Assume:
\mr
\item "{$(a)$}"  {\rm cf}$(\mu) \le \kappa < \mu,
\kappa^+ < \theta < \chi^*$ and
{\rm pp}$_\kappa(\mu) = \chi^*$, moreover {\rm pp}$_\kappa(\mu) =^+ \chi^*$
\sn
\item "{$(b)$}"  $\bold F$ is a function, with domain $[\mu]^\kappa$, such
that: 
for $a \in [\mu]^\kappa,\bold F(a)$ is a family of $< \theta$ members of
$[\mu]^\kappa$
\sn
\item "{$(c)$}"  $F$ is a function with domain $[\mu]^\kappa$ such that
$$
a \in [\mu]^\kappa \Rightarrow a \subseteq F(a) \in \bold F(a).
$$
\ermn
\underbar{Then} we can find pairwise distinct $a_i \in [\mu]^\kappa$ for
$i < \chi^*$ such that ${\Cal I} = \{a_i:i < \chi^*\}$ is 
$(F,\bold F)$-independent which means
\mr
\item "{$(*)_{F,\bold F,{\Cal I}}$}"  $\qquad a \ne b 
\and a \in {\Cal I} \and
b \in {\Cal I} \and c \in \bold F(a) \Rightarrow \neg(F(b) \subseteq c)$.
\endroster
\endproclaim
\bigskip

\remark{\stag{838-pcf.17} Remark}  1) Clearly this is a relative 
to Hajnal's free subset theorem \cite{Ha61}.
\nl
2) Note that we can choose $F(a) = a$. 
\nl
3) Also if $\mu_1 \le \mu$, cf$(\mu_1) \le \kappa \le \kappa + 
\theta < \mu_1$ and pp$_\kappa(\mu_1) \ge \mu$ then by 
\cite[Ch.II,2.3]{Sh:g} the Fact for $\mu_1$ implies the one for $\mu$. 
\nl
4) Note that if $\lambda = \text{ cf}([\mu]^\kappa,\subseteq)$
\ub{then} for some $\bold F,F$ as in the Fact we have
\mr
\item "{$\circledast$}"  if $a_i \in [\mu]^\kappa$ for 
$i < \lambda^+$ are pairwise
distinct then not every pair $\{a_i,a_j\}$ is $(\bold F,F)$-independent \nl
[why?  let ${\Cal P} \subseteq [\mu]^\kappa$ be cofinal (under $\subseteq$)
of cardinality $\lambda$, and let $F,\bold F$ be such that \nl
$\bold F(a) \subseteq \{b \in [\mu]^\kappa:a \subseteq b \text{ and }
b \in {\Cal P}\}$ has a $\subseteq$-maximal member $F(a)$; \nl
obviously there are such $F,\bold F$.
\ermn
 Now clearly
\mr
\item "{$(*)_1$}"  if $a \ne b$ are from $[\mu]^\kappa$ and 
$F(a) = F(b)$ then $\{a,b\}$ is not $(\bold F,F)$-independent.
\ermn
[Why?  Just look at the definition of $(\bold F,F)$-independent.]
\mr
\item "{$(*)_2$}"  if ${\Cal I} \subseteq [\mu]^\kappa$
is of cardinality $> \lambda$ (e.g. $\lambda^+$)
\ub{then} ${\Cal I}$ is not $(\bold F,F)$-independent.
\ermn
[Why?  As Rang$(F \rest {\Cal I}) \subseteq \text{ Rang}(F)
\subseteq {\Cal P}$ and ${\Cal P}$ has cardinality $\lambda$
necessarily there are $a \ne b$ from ${\Cal I}$ such that $F(a) =
F(b)$ and use $(*)_1$.]
\endremark
\bigskip

\demo{Proof}
\mr
\item "{$\boxtimes_1$}"  it suffices to prove the variant with 
$[\mu]^\kappa$ replaced by $[\mu]^{\le \kappa}$.
\ermn
[Why?    So we are given $\bold F,F$ as in the claim.
We define $g:[\mu]^{\le \kappa} \rightarrow
[\mu]^\kappa$ and functions $F',\bold F'$ with domain $[\mu]^{\le
\kappa}$ as follows:

$$
g(a) = \{\kappa + \alpha:\alpha \in a\} \cup \{\alpha:\alpha < \kappa\}
$$

$$
\bold F'(a) = \{\{\alpha:\kappa + \alpha \in b\}:b \in \bold F(g(a))\}
$$

$$
F'(a) = \{\alpha:\kappa + \alpha \in F(g(a))\}.
$$
\mn
Now $\bold F',F'$ are as in the claim only 
replacing everywhere $[\mu]^\kappa$ by
$[\mu]^{\le \kappa}$, and if ${\Cal I}' = \{a_i:i < \chi\} \subseteq
[\mu]^{\le \kappa}$ with no repetitions satisfying 
$(*)_{F',\bold F',{\Cal I}'}$ then we shall show that
${\Cal I} := \{g(a_i):i < \chi\}$ is with no repetitions and
$(*)_{F,\bold F,{\Cal I}}$ holds.
\sn
This clearly suffices, but why it holds?  
Clearly $g$ is a one-to-one function so $i \ne j < \chi \Rightarrow g(a_i)
\ne g(a_j)$ and Rang$(g) \subseteq [\mu]^\kappa$ so $g(a_i) \in
[\mu]^\kappa$.  Let $i \ne j$ and we should check that $[c' \in \bold
F(g(a_i)) \Rightarrow F(g(a_j)) \nsubseteq c']$, so fix $c'$ such that
$c' \in \bold F(g(a_i))$.

By the definition of $\bold F'(a_i)$ clearly $c := \{\alpha:\kappa +
\alpha \in c'\}$ belongs to $\bold F'(a_i)$.  By the choice of ${\Cal
I}' = \{a_i:i < \chi\}$ we know that $c \in \bold F'(a_i) \Rightarrow
F'(a_j) \nsubseteq c$, but by the previous sentence the antecedent
hold hence $F'(a_j) \nsubseteq c$ hence we can choose $\alpha \in
F'(a_j) \backslash c$.  By the choice of $F'(a_j)$ we have $\kappa +
\alpha \in F(g(a_j)$) and by the choice of $c$ we have $\kappa + \alpha
\notin c'$, so $\alpha$ witness $F(g(a_j)) \nsubseteq c'$ as required.]
\sn
So we conclude that we can replace $[\mu]^\kappa$ by
$[\mu]^{\le \kappa}$.  In fact we shall find the $a_i$ in
$[\mu]^{|{\frak a}|}$ where ${\frak a}$ chosen below.
\nl
As $\mu$ is a limit cardinal $\in (\kappa,\chi^*)$, if $\theta < \mu$
then we can replace
$\theta$ by $\theta^+$ but $\kappa^{++} < \mu$ 
so without loss of generality  $\kappa^{++} < \theta$.  

Now we prove
\mr
\item "{$\boxtimes_2$}"  for some unbounded subset $w$ of $\chi$
we have $\langle\text{Rang}(f_\alpha):\alpha \in w\rangle$ is $(\bold
F,F)$-independent \ub{when}:
{\roster
\itemitem{ $\oplus_{\chi,{\frak a},\bar f}$ }  $\theta < \chi = 
\text{ cf}(\Pi {\frak a}/J)$ where ${\frak a}
\subseteq \mu \cap \text{ Reg} \backslash \kappa^+,|{\frak a}| \le \kappa,
\sup({\frak a}) = \mu,J^{\text{bd}}_{\frak a} \subseteq J$ and for simplicity
$\chi = \text{ max pcf}({\frak a})$ and $\bar f = \langle f_\alpha:\alpha <
\chi \rangle$ is a sequence of members of $\Pi{\frak a},<_J$-increasing, and
cofinal in $(\Pi{\frak a},<_J)$, so, of course, $\chi \le \chi^*$.  
\endroster}
\ermn
Without loss of generality $f_\alpha(\lambda)
> \sup({\frak a} \cap \lambda)$ for $\lambda \in {\frak a}$. 

Also for every
$a \in [\mu]^\kappa$, define ch$_a \in \Pi{\frak a}$ by
ch$_a(\lambda) = \sup(a \cap \lambda)$ for $\lambda \in {\frak a}$ 
so for some $\zeta(a) < \chi$ we have ch$_a <_J f_{\zeta(a)}$ 
(as $\langle f_\alpha:\alpha < \chi \rangle$ is
cofinal in $(\Pi{\frak a},<_J)$).  So for each $a \in [\mu]^\kappa$, as
$|\bold F(a)| < \theta < \chi = \text{ cf}(\chi)$ clearly 
$\xi(a) := \sup\{\zeta(b):b \in
\bold F(a)\}$ is $< \chi$, and clearly $(\forall b \in \bold F(a))
[\text{ch}_b <_J f_{\xi(a)}]$.  So 
$C := \{\gamma < \chi:\text{for every } \beta < \gamma,
\xi(\text{Rang}(f_\beta)) < \gamma\}$ is a club of $\chi$.

For each 
$\alpha < \chi$, Rang$(f_\alpha) \in [\mu]^\kappa$, hence
$\bold F(\text{Rang}(f_\alpha))$ has cardinality $< \theta$, but 
$\theta <  \chi = \text{ cf}(\chi)$ hence for some $\theta_1 < \theta$ we 
have $\theta_1 > \kappa^+$ and $\chi = \sup\{\alpha < \chi:|\bold F
(\text{Rang}(f_\alpha)| \le \theta_1\}$, so without loss 
of generality $\alpha < \chi \Rightarrow \theta_1 \ge |\bold F(\text{Rang}
(f_\alpha))|$.
\medskip

As $\kappa^+ < \theta_1$, by \cite[\S1]{Sh:420} there are $\bar d,S$
such that
\mr
\item "{$(*)_1$}"  $(a) \quad S \subseteq \theta^+_1$ is a stationary 
\sn
\item "{${{}}$}"  $(b) \quad S \subseteq \{\delta <
\theta^+_1:\text{cf}(\delta) = \kappa^+\}$
\sn
\item "{${{}}$}"  $(c) \quad S$ belongs to $\check I[\theta^+_1]$, 
\sn
\item "{${{}}$}"  $(d) \quad \langle d_i:i < \theta^+_1 \rangle$ 
witness it, so otp$(d_i) \le \kappa^+,d_i \subseteq i,
[j \in d_i \Rightarrow d_j = d_i \cap i]$
\nl

\hskip25pt  and $i \in S \Rightarrow i = \sup(d_i)$, 
\nl

\hskip25pt and for simplicity (see \cite[III]{Sh:g})
\sn
\item "{${{}}$}"  $(e) \quad$ for every club $E$ of
$\theta^+_1$ for stationarily many $\delta \in S$ we have
\nl

\hskip25pt  $(\forall \alpha
\in d_\delta)[(\exists \beta \in E)(\sup(\alpha \cap d_\delta) < \beta <
\alpha)]$.
\ermn
Now try to choose by induction on $i < \theta^+_1$, a triple
$(g_i,\alpha_i,w_i)$ such that:
\mr
\item "{$(*)_2$}"   $(a) \quad g_i \in \Pi{\frak a}$
\sn
\item "{${{}}$}"  $(b) \quad$ if $j < i$ then\footnote{in fact, without loss
of generality
 min$({\frak a}) > \theta^+_1$, so we can demand $g_j < g_i$ so
clause (c) is redundant}
$g_j <_J g_i$
\sn
\item "{${{}}$}"  $(c) \quad (\forall \lambda \in {\frak a})
(\underset{j \in d_i} {}\to \sup \,g_j(\lambda) < g_i(\lambda))$
\sn
\item "{${{}}$}"  $(d) \quad \alpha_i < \chi$ and $\alpha_i > 
\sup(\dsize \bigcup_{j<i} w_j)$
\sn
\item "{${{}}$}"  $(e) \quad j < i \Rightarrow \alpha_j < \alpha_i$
\sn
\item "{${{}}$}"  $(f) \quad g_i <_J f_{\alpha_i}$
\sn
\item "{${{}}$}"  $(g) \quad \beta \in \dsize 
\bigcup_{j<i} w_j \Rightarrow \xi(\text{Rang}(f_\beta))
< \alpha_i \and f_\beta <_J g_i$
\sn
\item "{${{}}$}"  $(h) \quad w_i$ is a maximal subset of 
$(\alpha_i,\chi)$ satisfying
{\roster
\itemitem{ ${{}}$ }  $(*) \quad \beta \in 
w_i \and \gamma \in w_i \and \beta \ne \gamma
\and a \in \bold F(\text{Rang}(f_\beta)) \Rightarrow$
\nl

\hskip30pt $\neg(F(\text{Rang}(f_\gamma)) \subseteq a)$
\nl

and moreover
\sn
\itemitem{ ${{}}$ }   $(*)^+ \quad \beta \in w_i 
\and \gamma \in w_i \and \beta \ne
\gamma \and a \in \bold F(\text{Rang}(f_\beta)) \Rightarrow$ 
\nl

\hskip30pt $\{\lambda \in {\frak a}:f_\gamma (\lambda) \in a\} \in J$.
\endroster}
\ermn
Note that really (as indicated by the notation)
\mr 
\item "{$\otimes$}"  if $w \subseteq (\alpha_i,\chi)$ satisfies $(*)^+$ then
it satisfies $(*)$.
\ermn
[Why?  let us check $(*)$, so let $\beta \in w,\gamma \in w,\beta \ne \gamma$
and $a \in \bold F(\text{Rang}(f_\beta))$; by $(*)^+$ we know that
${\frak a}' = \{\lambda \in {\frak a}:f_\gamma(\lambda) \in a\} \in J$.
Now as $J$ is a
proper ideal on ${\frak a}$ clearly for some $\lambda \in {\frak a}$ we have
$\lambda \notin {\frak a}'$, hence $f_\gamma(\lambda) 
\notin a$ but $f_\gamma(\lambda) \in \text{ Rang}(f_\gamma)$ and by the
assumption on $(\bold F,F)$ we have Rang$(f_\gamma) 
\subseteq F(\text{Rang}(f_\gamma))$ hence
$f_\gamma(\lambda) \in F(\text{Rang}(f_\gamma)) \backslash a$ so $\neg(F
(\text{Rang}(f_\gamma)) \subseteq a)$, as required.]

We claim that we cannot carry the induction because if we succeed, then as
cf$(\chi) = \chi > \theta \ge \theta^+_1$ there is $\alpha$ such that
$\dsize \bigcup_{i < \theta^+_1} \alpha_i < \alpha < \chi$ and let
$\bold F(\text{Rang}(f_\alpha)) = 
\{a^\alpha_\zeta:\zeta < \theta_1\}$ (possible as
$1 \le |\bold F(\text{Rang}(f_\alpha))| \le \theta_1)$.  
Now for each $i < \theta^+_1$, by the choice of $w_i$ clearly $w_i \cup 
\{\alpha\}$ does not
satisfy the demand in clause $(h)$ and let it be exemplified by some pair
$(\beta_i,\gamma_i)$.  Now $\{\beta_i,\gamma_i\} \subseteq w_i$ is
impossible by the choice of $w_i$, i.e. as $w_i$ satisfies clause
(h).  Also $\beta_i \in w_i \wedge \gamma_i = \alpha$ is impossible 
as $\beta_i \in w_i \Rightarrow \xi(\text{Rang}(f_{\beta_i})) < 
\alpha_{i+1} < \alpha$, so necessarily $\gamma_i \in w_i$ and 
$\beta_i = \alpha$, so for some $a' \in \bold
F(\text{Rang}(f_{\beta_i})) = \bold F(\text{Rang}(f_\alpha))$ the
conclusion of $(*)^+$ fails, so as $\langle a^\alpha_\zeta:\zeta <
\theta_1\rangle$ list $\bold F(\text{Rang}(f_\alpha))$ it follows that for some
$\zeta_i < \theta_1$ we have

$$
{\frak a}_i = \{\lambda \in {\frak a}:f_{\gamma_i}(\lambda) \in 
a^\alpha_{\zeta_i}\} \notin J.
$$
\mn
[why use the ideal?  In order to show below that ${\frak b}_\varepsilon \ne
\emptyset$.]
But cf$(\theta^+_1) = \theta^+_1 > \theta_1$, so for some $\zeta(*) <
\theta^+_1$ we have $A := \{i:\zeta_i = \zeta(*)\}$ is unbounded in 
$\theta^+_1$.  Hence $E = \{\alpha < \theta^+_1:\alpha$ a limit ordinal and
$A \cap \alpha$ is unbounded in $\alpha\}$ is a club of $\theta^+_1$.  So for
some $\delta \in S$ we have $\delta = \sup(A \cap \delta)$, moreover
letting $\{\alpha_\varepsilon:\varepsilon < \kappa^+\}$ list
$d_\delta$ in increasing order, we have
$(\forall \varepsilon)[E \cap (\underset{\zeta < \varepsilon} {}\to \sup
\,\alpha_\zeta,\alpha_\varepsilon) \ne \emptyset]$ hence we can find
$i(\delta,\varepsilon) \in (\underset{\zeta < \varepsilon} {}\to \sup \,
\alpha_\zeta,\alpha_\varepsilon) \cap A$ for each $\varepsilon < \kappa^+$.
\mn
Clearly for each $\varepsilon < \kappa^+$

$$
\align
{\frak b}_\varepsilon = \bigl\{ \lambda \in {\frak a}:
g_{i(\delta,\varepsilon)} (\lambda) &<
f_{\alpha_{i(\delta,\varepsilon)}}
(\lambda) <
f_{\gamma_{i(\delta,\varepsilon)}}(\lambda) \\
  &< g_{i(\delta,\varepsilon)+1}(\lambda) <
f_{\alpha_{i(\delta,\varepsilon)+1}}(\lambda) < f_\alpha(\lambda)\} 
= {\frak a} \text{ mod } J
\endalign 
$$
\mn
hence ${\frak b}_\varepsilon \cap {\frak a}_{i(\delta,\varepsilon)} \notin
\emptyset$.  Moreover, ${\frak b}_\varepsilon \cap 
{\frak a}_{i(\delta,\varepsilon)}
\notin J$.  Now for each $\lambda \in {\frak a}$ let $\varepsilon
(\lambda)$ be
$\sup\{\varepsilon < \kappa^+:\lambda \in {\frak b}_\varepsilon \cap 
{\frak a}_{i(\delta,\varepsilon)}\}$ and let $\varepsilon(*) = 
\sup\{\varepsilon(\lambda):\lambda \in
{\frak a} \text{ and } \varepsilon(\lambda) < \kappa^+\}$ so as $|{\frak a}|
\le \kappa$ clearly $\varepsilon(*) < \kappa^+$.  Let $\lambda^* \in
{\frak b}_{\varepsilon(*)+1} \cap {\frak a}_{i(\delta,\varepsilon(*)+1)}$, so
$B := \{\varepsilon < \kappa^+:\lambda^* \in {\frak b}_\varepsilon \cap
{\frak a}_{i(\delta,\varepsilon)}\}$ is unbounded in $\kappa^+,\langle
f_{\beta_{i(\delta,\varepsilon)}}(\lambda^*):\varepsilon \in B \rangle$ is
strictly increasing (see clause $(c)$ above and the choice of
${\frak b}_\varepsilon$) and $\varepsilon \in B
\Rightarrow f_{\beta_{i(\delta,\varepsilon)}}(\lambda^*) \in 
a^\alpha_{\zeta(*)}$
(by the definition of ${\frak a}_{i(\delta,\varepsilon)}$, and $\zeta(*)$ as
$\zeta_{i(\delta,\varepsilon)} = \zeta(*)$).  We get contradiction to
$a \in \bold F(\text{Rang}(f_\alpha)) 
\Rightarrow |a| \le \kappa$.

So really we cannot carry the induction in $(*)_2$ so we are stuck 
at some $i < \theta^+_1$.  If
$i=0$, or $i$ limit, or $i = j+1 \and \sup(w_j) < \chi$ we can find $g_i$
and then $\alpha_i$ and then $w_i$ as required.  So necessarily
$i =j+1,\sup(w_j) = \chi$.  So we have finished proving $\boxtimes_2$.
\mr
\item "{$\boxtimes_3$}"  there is ${\Cal I} \subseteq [\mu]^{\le
\kappa}$ as required.
\ermn
Now if $\chi^*$ is regular, recalling that we assume pp$_\kappa(\mu)
=^+ \chi^*$ there are ${\frak a},J$ as required in $\oplus$ above for
$\chi = \chi^*$, hence also such $\bar f$.  Applying $\boxtimes_2$ to
$(\chi,{\frak a},J,\bar f)$ we get $w$ as there.  Now
$\langle\text{Rang}(f_\alpha):\alpha \in w\rangle$
is as
required in the fact.  So the only case left is when 
$\chi^*$ is singular.  Let $\chi^* = \underset{\varepsilon <
\text{ cf}(\chi^*)} {}\to \sup \chi_\varepsilon$ and 
$\chi_\varepsilon \in (\mu,\chi^*) \cap \text{ Reg}$ is 
(strictly) increasing with $\varepsilon$.  By \cite[Ch.II,\S3]{Sh:g} we 
can find, for each $\varepsilon < \text{ cf}(\chi^*),
{\frak a}_\varepsilon,J_\varepsilon,\bar f^\varepsilon = \langle
f^\varepsilon_\alpha:\alpha < \chi_\varepsilon \rangle$ satisfying the
demands in $\oplus$ above, but in addition
\mr
\item "{$\odot$}"  $\bar f^\varepsilon$ is $\mu^+$-free i.e. for every $u \in
[\chi_\varepsilon]^\mu$ there is a sequence 
$\langle {\frak b}_\alpha:\alpha \in u
\rangle$ such that ${\frak b}_\alpha \in J_\varepsilon$ and for each 
$\lambda \in {\frak a}_\varepsilon,\langle
f^\varepsilon_\alpha(\lambda):\alpha \text{ satisfies } \lambda 
\notin {\frak b}_\alpha \rangle$ is strictly 
increasing.
\ermn
So for every $a \in [\mu]^{\le\kappa}$ and $\varepsilon 
< \text{ cf}(\chi^*)$ we have

$$
\bigl\{\alpha < \chi_\varepsilon:\{\lambda \in {\frak a}_\varepsilon:f_\alpha
(\lambda) \in a\} \notin J_\varepsilon \bigr\} \text{ has cardinality } \le 
\kappa.
$$
\mn
Hence for each $a \in [\mu]^{\le\kappa}$

$$
\bigl\{(\varepsilon,\alpha):\varepsilon < \text{ cf}(\chi^*) \text{ and }
\alpha < \chi_\varepsilon \text{ and }\{\lambda \in {\frak a}_\varepsilon:
f_\alpha(\lambda) \in a\} \notin J_\varepsilon \bigr\}
$$
\mn
has cardinality $\le \kappa + \text{ cf}(\chi^*) = \text{ cf}(\chi^*)$
as for singular $\mu > \kappa \ge \text{ cf}(\mu)$ we have 
cf(pp$_\kappa(\mu)) > \kappa$.
\mn
Define: $X = \{(\varepsilon,\alpha):\varepsilon < \text{ cf}(\chi^*),\alpha
< \chi_\varepsilon\}$

$$
\align
F'\bigl((\varepsilon,\alpha)\bigr) = \bigl\{ (\varepsilon',\alpha'):
&(\varepsilon',\alpha') \in X \backslash \{(\varepsilon,\alpha)\} \text{ and
for some} \\
  &d \in \bold F(\text{Rang}(f^\varepsilon_\alpha))
\text{ we have }
\{\lambda \in {\frak a}_\varepsilon:f^{\varepsilon'}_{\alpha'}(\lambda) 
\in d\} \notin J_{\varepsilon'} \bigr\}
\endalign
$$
\mn
so $F'\bigl((\varepsilon,\alpha)\bigr)$ is a subset of $X$ of cardinality
$< \text{ cf}(\chi^*)^+ + \theta < \chi^*$.

So by Hajnal's free subset theorem \cite{Ha61} we finish proving
$\boxtimes_3$ (we could alternatively, for
$\chi^*$ singular, have imitated his proof). 

Recalling $\boxtimes_1$ we are done.  \hfill$\square_{\scite{838-pcf.11}}$
\enddemo

\nocite{ignore-this-bibtex-warning}

%
 REFERENCES.  
 \bibliographystyle{lit-plain}
 \bibliography{lista,listb,listx,listf,liste}

\enddocument